\documentclass[a4,10pt,english]{smfbook}
\usepackage{amsmath,amscd,amssymb}
\makeindex

\newcommand{\nbiga}{\mathcal{A}}
\newcommand{\nbigb}{\mathcal{B}}
\newcommand{\nbigc}{\mathcal{C}}

\newcommand{\nbige}{\mathcal{E}}
\newcommand{\nbigf}{\mathcal{F}}
\newcommand{\nbigg}{\mathcal{G}}

\newcommand{\nbigi}{\mathcal{I}}
\newcommand{\nbigj}{\mathcal{J}}
\newcommand{\nbigk}{\mathcal{K}}
\newcommand{\nbigl}{\mathcal{L}}
\newcommand{\nbigm}{\mathcal{M}}
\newcommand{\nbign}{\mathcal{N}}
\newcommand{\nbigo}{\mathcal{O}}
\newcommand{\nbigp}{\mathcal{P}}

\newcommand{\nbigr}{\mathcal{R}}
\newcommand{\nbigs}{\mathcal{S}}
\newcommand{\nbigt}{\mathcal{T}}
\newcommand{\nbigu}{\mathcal{U}}
\newcommand{\nbigv}{\mathcal{V}}
\newcommand{\nbigw}{\mathcal{W}}

\newcommand{\proj}{\mathbb{P}}
\newcommand{\seisuu}{{\mathbb Z}}
\newcommand{\rnum}{{\mathbb Q}}

\newcommand{\cnum}{{\mathbb C}}
\newcommand{\real}{{\mathbb R}}

\newcommand{\DD}{\mathbb{D}}
\newcommand{\EE}{\mathbb{E}}

\newcommand{\LL}{\mathbb{L}}

\newcommand{\gbiga}{\mathfrak A}
\newcommand{\gbigb}{\mathfrak B}
\newcommand{\gbigc}{\mathfrak C}

\newcommand{\gbigf}{\mathfrak F}

\newcommand{\gbigl}{\mathfrak L}
\newcommand{\gbigm}{\mathfrak M}
\newcommand{\gbign}{\mathfrak N}

\newcommand{\gbigp}{\mathfrak P}

\newcommand{\gbigs}{\mathfrak S}

\newcommand{\gminia}{\mathfrak a}
\newcommand{\gminib}{\mathfrak b}
\newcommand{\gminic}{\mathfrak c}

\newcommand{\gminik}{\mathfrak k}

\newcommand{\gminim}{\mathfrak m}

\newcommand{\gminiq}{\mathfrak q}

\newcommand{\gminiv}{\mathfrak v}

\newcommand{\vecxi}{{\boldsymbol \xi}}

\newcommand{\vece}{{\boldsymbol e}}

\newcommand{\vecv}{{\boldsymbol v}}
\newcommand{\vecu}{{\boldsymbol u}}
\newcommand{\vecw}{{\boldsymbol w}}

\newcommand{\veczero}{{\boldsymbol 0}}

\newcommand{\veca}{{\boldsymbol a}}
\newcommand{\vecb}{{\boldsymbol b}}

\newcommand{\vectau}{{\boldsymbol \tau}}

\newcommand{\vecm}{{\boldsymbol m}}

\newcommand{\vecx}{{\boldsymbol x}}

\newcommand{\vecV}{{\boldsymbol V}}

\newcommand{\vecL}{{\boldsymbol L}}

\newcommand{\lrarr}{\longrightarrow}

\newcommand{\pf}{{\bf Proof}\hspace{.1in}}

\def\Hom{\mathop{\rm Hom}\nolimits}

\def\End{\mathop{\rm End}\nolimits}

\def\Image{\mathop{\rm Im}\nolimits}

\def\Re{\mathop{\rm Re}\nolimits}

\def\Gr{\mathop{\rm Gr}\nolimits}

\def\GL{\mathop{\rm GL}\nolimits}

\def\rank{\mathop{\rm rank}\nolimits}

\def\Ker{\mathop{\rm Ker}\nolimits}

\def\length{\mathop{\rm length}\nolimits}
\def\Gr{\mathop{\rm Gr}\nolimits}
\def\Sym{\mathop{\rm Sym}\nolimits}

\def\ad{\mathop{\rm ad}\nolimits}
\def\Res{\mathop{\rm Res}\nolimits}

\def\ord{\mathop{\rm ord}\nolimits}
\def\degpar{\mathop{\rm par\textrm{-}deg}\nolimits}

\def\tr{\mathop{\rm tr}\nolimits}
\def\Tr{\mathop{\rm Tr}\nolimits}
\def\vol{\mathop{\rm dvol}\nolimits}
\def\dvol{\mathop{\rm dvol}\nolimits}

\def\id{\mathop{\rm id}\nolimits}

\def\Irr{\mathop{\rm Irr}\nolimits}

\newcommand{\del}{\partial}
\newcommand{\delbar}{\overline{\del}}

\newcommand{\pardeg}{\degpar}

\newcommand{\nhom}{{\mathcal Hom}}

\newcommand{\vecP}{\boldsymbol P}

\newcommand{\shikaku}{\sharp}

\newcommand{\barz}{\overline{z}}
\newcommand{\zbar}{\barz}
\newcommand{\zetabar}{\overline{\zeta}}
\newcommand{\baralpha}{\overline{\alpha}}
\newcommand{\alphabar}{\baralpha}
\newcommand{\barlambda}{\overline{\lambda}}
\newcommand{\lambdabar}{\barlambda}

\newcommand{\varphibar}{\overline{\varphi}}
\newcommand{\etabar}{\overline{\eta}}
\newcommand{\xibar}{\overline{\xi}}

\newcommand{\xbar}{\overline{x}}
\newcommand{\sbar}{\overline{s}}

\newcommand{\Abar}{\overline{A}}
\newcommand{\abar}{\overline{a}}

\newcommand{\Poin}{{\bf p}}

\newcommand{\DDlambda}{\DD^{\lambda}}

\newcommand{\DDlambdahat}{\widehat{\DD}^{\lambda}}
\newcommand{\DDhatlambda}{\DDlambdahat}

\newcommand{\nbigelambda}{\nbige^{\lambda}}

\newcommand{\KMS}{{\mathcal{KMS}}}

\newcommand{\Par}{{\mathcal Par}}
\newcommand{\Sp}{{\mathcal Sp}}

\newcommand{\kmsmap}{\gminik}

\newcommand{\lefttop}[1]{{}^{#1}\!}

\newcommand{\lamda}{\lambda}

\newcommand{\vecnbigl}{{\boldsymbol{\mathcal L}}}

\newcommand{\closedopen}[2]{[#1,#2[}

\newcommand{\openopen}[2]{]#1,#2[}
\newcommand{\closedclosed}[2]{[#1,#2]}

\newcommand{\Etilde}{\widetilde{E}}

\newcommand{\thetatilde}{\widetilde{\theta}}
\newcommand{\Vhat}{\widehat{V}}
\newcommand{\nablahat}{\widehat{\nabla}}
\newcommand{\Vtilde}{\widetilde{V}}
\newcommand{\nablatilde}{\widetilde{\nabla}}
\newcommand{\vecvtilde}{\widetilde{\vecv}}

\newcommand{\vecutilde}{\widetilde{\vecu}}

\newcommand{\vecVhat}{\widehat{\vecV}}

\newcommand{\wbar}{\overline{w}}

\newcommand{\htilde}{\widetilde{h}}

\newcommand{\ptilde}{\widetilde{p}}
\newcommand{\pitilde}{\widetilde{\pi}}

\newcommand{\vtilde}{\widetilde{v}}

\newcommand{\ftilde}{\widetilde{f}}
\newcommand{\utilde}{\widetilde{u}}
\newcommand{\Ftilde}{\widetilde{F}}

\newcommand{\gminiabar}{\overline{\gminia}}

\newcommand{\stilde}{\widetilde{s}}
\newcommand{\mubar}{\overline{\mu}}

\newcommand{\nbiglhat}{\widehat{\nbigl}}

\newcommand{\nbigltilde}{\widetilde{\nbigl}}

\newcommand{\DDtilde}{\widetilde{\DD}}

\newcommand{\nbigbtilde}{\widetilde{\nbigb}}

\newcommand{\nbigvtilde}{\widetilde{\nbigv}}
\newcommand{\nbigvhat}{\widehat{\nbigv}}

\newcommand{\Utilde}{\widetilde{U}}

\newcommand{\vecy}{\boldsymbol y}

\newcommand{\DDhat}{\widehat{\DD}}

\newcommand{\Ltilde}{\widetilde{L}}

\def\ord{\mathop{\rm ord}\nolimits}

\def\Gal{\mathop{\rm Gal}\nolimits}

\def\Loc{\mathop{\rm Loc}\nolimits}

\def\Ad{\mathop{\rm Ad}\nolimits}
\def\ad{\mathop{\rm ad}\nolimits}
\def\RFM{\mathop{\rm RFM}\nolimits}

\def\Hit{\mathop{\rm Hit}\nolimits}

\def\cov{\mathop{\rm cov}\nolimits}
\def\SU{\mathop{\rm SU}\nolimits}

\newcommand{\Dbar}{\overline{D}}

\newcommand{\Ptilde}{\widetilde{P}}

\newcommand{\ubar}{\overline{u}}

\newcommand{\gtilde}{\widetilde{g}}
\newcommand{\Ztilde}{\widetilde{Z}}

\newcommand{\Ctilde}{\widetilde{C}}

\newcommand{\Atilde}{\widetilde{A}}

\newcommand{\Ytilde}{\widetilde{Y}}

\newcommand{\nbigvlambda}{\nbigv^{\lambda}}

\newcommand{\gbigntilde}{\widetilde{\gbign}}

\newcommand{\betabar}{\overline{\beta}}

\newcommand{\nbigmlambda}{\nbigm^{\lambda}}

\newcommand{\taubar}{\overline{\tau}}

\newcommand{\veckappa}{{\boldsymbol \kappa}}

\newcommand{\Hhat}{\widehat{H}}

\newcommand{\nbigmbar}{\overline{\nbigm}}

\newcommand{\gammabar}{\overline{\gamma}}

\newcommand{\Lbar}{\overline{L}}
\newcommand{\Mbar}{\overline{M}}
\newcommand{\Bbar}{\overline{B}}

\newcommand{\Sptilde}{\widetilde{\Sp}}

\newcommand{\Gammabar}{\overline{\Gamma}}

\newcommand{\inftyhat}{\widehat{\infty}}

\newcommand{\bbar}{\overline{b}}

\newcommand{\ttilde}{\widetilde{t}}

\newcommand{\chitilde}{\widetilde{\chi}}

\newcommand{\LLhat}{\widehat{\LL}}
\newcommand{\Ihat}{\widehat{I}}
\newcommand{\vecG}{\boldsymbol G}

\newcommand{\vecell}{\boldsymbol \ell}
\newcommand{\su}{\mathfrak{su}}

\newtheorem{thm}{Theorem}[section]
\newtheorem{cor}[thm]{Corollary}

\newtheorem{rem}[thm]{Remark}
\newtheorem{lem}[thm]{Lemma}
\newtheorem{prop}[thm]{Proposition}
\newtheorem{df}[thm]{Definition}

\newtheorem{example}[thm]{Example}
\newtheorem{condition}[thm]{Condition}

\author{Takuro Mochizuki}
\address{Research Institute for Mathematical Sciences,
Kyoto University,
Kyoto 606-8502, Japan}
\email{takuro@kurims.kyoto-u.ac.jp}

\title{Periodic monopoles and difference modules}

\begin{document}

\frontmatter

\begin{abstract}
We study periodic monopoles satisfying some mild conditions,
called of GCK type.
Particularly, we give a classification of
periodic monopoles of GCK type 
in terms of difference modules
with parabolic structure,
which is a kind of Kobayashi-Hitchin correspondence
between differential geometric objects
and algebraic objects.
We also clarify the asymptotic behaviour of 
periodic monopoles of GCK type around infinity.
\end{abstract}


\maketitle

\tableofcontents

\mainmatter

\chapter{Introduction}

\section{Monopoles of GCK-type}

For $T>0$, we set $S^1_T:=\real/T\seisuu$.
\index{space $S^1_T$}
Let $(t,w)$ denote the standard local coordinate system
on $S^1_T\times\cnum$.
\index{coordinate system $(t,w)$}
We regard $S^1_T\times\cnum$
as a Riemannian manifold
with the metric $dt\,dt+dw\,d\wbar$.

A periodic monopole is a monopole
$(E,h,\nabla,\phi)$ on $S^1_T\times\cnum$.
\index{periodic monopole}
Namely, $E$ is a vector bundle on $S^1_T\times\cnum$
with a Hermitian metric $h$,
a unitary connection $\nabla$
and an anti-Hermitian endomorphism $\phi$
satisfying the Bogomolny equation
\begin{equation}
\label{eq;21.7.28.1}
 F(\nabla)-\ast\nabla\phi=0,
\end{equation}
where $F(\nabla)$ denote the curvature of $\nabla$,
and $\ast$ denotes the Hodge star operator.
\index{monopole}
\index{Bogomolny equation}
More precisely,
we admit that the monopole may have isolated singularities
at a finite subset $Z\subset S^1_T\times\cnum$,
i.e.,
the monopole
$(E,h,\nabla,\phi)$ is defined on $(S^1_T\times\cnum)\setminus Z$.
We impose the following conditions 
on the behaviour of $(E,h,\nabla,\phi)$ 
around $Z$ and $S^1_T\times\{\infty\}$.
\begin{itemize}
\item
 Each point of $Z$ is Dirac type singularity
 of $(E,h,\nabla,\phi)$.
\item
 $F(\nabla)\to 0$
 and $\phi=O(\log|w|)$ as $|w|\to\infty$.
\end{itemize}
In this paper,
such monopoles are called
of GCK-type (generalized Cherkis-Kapustin type).
\index{GCK-type}

Cherkis and Kapustin \cite{Cherkis-Kapustin1, Cherkis-Kapustin2}
studied such monopoles
under some more additional assumptions of genericity
on the behaviour around $Z$ and $S^1_T\times\{\infty\}$.
In particular, they studied the Nahm transforms between
such periodic monopoles and 
harmonic bundles on $(\proj^1,\{0,\infty,p_1,\ldots,p_m\})$.
Foscolo studied the deformation theory of such periodic monopoles
in \cite{Foscolo-deformation},
and the gluing construction in \cite{Foscolo-construction}.
More recently, Harland \cite{Harland}
classified $\SU(2)$-monopoles in the case $Z=\emptyset$
in terms of line bundles with parabolic structure
on spectral curves.
See also \cite{Maldonado, Maldonado-Ward, Harland-Ward}
on the geometry of the moduli space of some type of
periodic monopoles.

In this monograph,
we shall study periodic monopoles of GCK-type.
The first goal is to clarify their asymptotic behaviour 
around $S^1_T\times\{\infty\}$.
Then, we shall establish that
periodic monopoles are equivalent to difference modules.

There are two origins of this study.
One is the classifications of monopoles
in terms of holomorphic objects.
The other is the non-abelian Hodge theory
for harmonic bundles on compact Riemann surfaces.
Let us briefly recall them
in \S\ref{section;20.8.9.1}
and \S\ref{section;20.8.9.2}.

\section{Previous works on monopoles
and algebraic objects}
\label{section;20.8.9.1}

There are several interesting and deep studies
on the classification of monopoles in terms of algebraic objects.

The most classical and pioneering work is due to Donaldson and Hitchin.
Roughly speaking,
they obtained an equivalence between
$\SU(2)$-monopoles
$(E,h,\nabla,\phi)$ on $\real^3$
with finite energy
$\int_{\real^3}|F(\nabla)|^2<\infty$,
and based holomorphic maps $\varphi:\proj^1\lrarr\proj^1$,
where $\varphi$ is called based if $\varphi(0)=\infty$.
More precisely, 
Hitchin \cite{Hitchin-construction-monopole}
established that
the Nahm transform induces an equivalence
between monopoles and solutions of the Nahm equation,
and Donaldson \cite{Donaldson-Nahm} obtained 
an equivalence between
the solutions of the Nahm equation
and based holomorphic maps $\varphi:\proj^1\lrarr\proj^1$.
Hurtubise \cite{Hurtubise} clarified 
a more direct approach to obtain
based holomorphic maps from monopoles.
(Atiyah \cite{Atiyah-hyperbolic-monopoles}
gave an analogue construction
for monopoles on hyperbolic spaces.)
We shall review a construction of a based holomorphic map
from an $\SU(2)$-monopole in \S\ref{subsection;21.8.2.1}.
We remark that the periodic monopoles are excluded
by the finite energy condition.

The result was generalized by Hurtubise and Murray 
\cite{Hurtubise-classification-1989, Hurtubise-Murray1, Hurtubise-Murray2}
and Jarvis \cite{Jarvis1, Jarvis2} to the context of
a more general compact Lie group $G$ rather than $\SU(2)$.
Namely, they obtained a classification of
$G$-monopoles in terms of holomorphic maps from $\proj^1$
to flag varieties associated to $G$.

The interest in monopoles was renewed 
by the work of Kapustin and Witten \cite{Kapustin-Witten}
on the geometric Langlands theory from a physics viewpoint.
Inspired by their work,
Norbury \cite{Norbury} established that 
any Hecke transform of a holomorphic bundle
on a Riemann surface $\Sigma$
is represented by a singular monopole
on $\openopen{0}{1}\times\Sigma$
satisfying the Dirichlet condition at $t=0$
and the Neumann condition at $t=1$.
Here, for $a<b$, we set $\openopen{a}{b}:=\{a<t<b\}$.
\index{interval $\openopen{a}{b}$}
Charbonneau and Hurtubise \cite{Charbonneau-Hurtubise}
studied singular monopoles on 
$S^1\times\Sigma$,
and they established an equivalence
between such singular monopoles
and holomorphic bundles with a meromorphic automorphism
satisfying a stability condition on $\Sigma$.
Because our study is directly influenced by
the work of Charbonneau and Hurtubise,
we shall briefly review it in \S\ref{subsection;21.8.2.2}.

As noted above,
we shall briefly recall the constructions of
algebraic objects
in the correspondences due to Donaldson-Hitchin,
and Charbonneau-Hurtubise
in \S\ref{subsection;21.8.2.1}
and \S\ref{subsection;21.8.2.2}, respectively.
Though scattering maps are applied in the both cases,
there are also different flavors.
See also a remark in \S\ref{subsection;21.8.2.3}.

\subsection{$\SU(2)$-monopoles with finite energy on $\real^3$}
\label{subsection;21.8.2.1}

We recall an outline of the construction of
based holomorphic maps from $\SU(2)$-monopoles 
with finite energy on $\real^3$,
where some of fundamental concepts have already appeared.
We follow the explanation in
\cite[\S2 and \S16]{Atiyah-Hitchin-book}
and \cite{Hurtubise,Hurtubise-classification-1989}
though we omit the detail.
(See also \cite{Jarvis1}.)
Let $(E,h,\nabla,\phi)$ be a monopole
with finite energy
on $\real^3=\{(x_1,x_2,x_3)\in\real^3\}$
with the Euclidean metric $\sum dx_i^2$.
Following \cite{Atiyah-Hitchin-book},
on the basis of the results in \cite{Jaffe-Taubes-book},
we restrict ourselves to the case where
the following asymptotic conditions are satisfied.
\begin{itemize}
 \item $|\phi|_h=1-\frac{k}{2r}+O(r^{-2})$ as $r\to\infty$,
       where $r=\sqrt{\sum x_i^2}$,
       $|\phi|_h=-\frac{1}{2}\Tr(\phi^2)$
       and $k$ denotes a positive integer.
 \item $\frac{\del|\phi|_h}{\del\Omega}=O(r^{-2})$ as $r\to\infty$,
       where $\frac{\del}{\del\Omega}$ denotes
       the angular derivative.
 \item $|\nabla\phi|_h=O(r^{-2})$ as $r\to\infty$.
\end{itemize}

We identify $\real^3$ with $\real_t\times\cnum_z$
by setting $t=x_1$ and $z=x_2+\sqrt{-1}x_3$.
Note that $dt\,dt+dz\,d\zbar=\sum dx_i^2$.
We set $E^{t}:=E_{|\{t\}\times\cnum}$ for $t\in\real$.
They are enhanced to holomorphic vector bundles
by the differential operator
$\nabla_{\zbar}=\frac{1}{2}(\nabla_{x_2}+\sqrt{-1}\nabla_{x_3})$.
For any $t,t'\in\real$,
we have the isomorphism
$E^{t}\simeq E^{t'}$
obtained as the parallel transport
with respect to the differential operator
$\del_t:=\nabla_{x_1}-\sqrt{-1}\phi$.
The isomorphisms are called the {\em scattering maps}.
\index{scattering map}
Because the Bogomolny equation implies
the commutativity
\begin{equation}
\label{eq;21.7.28.2}
 \bigl[\del_t,\nabla_{\zbar}\bigr]=0,
\end{equation}
the scattering maps are holomorphic.
This is one of the key facts
in the study of the relationship between monopoles
and holomorphic or algebraic objects,
which goes back to \cite{Hitchin-monopoles-geodesics}.

\begin{rem}
In our study of periodic monopoles,
it is useful to formulate
the integrability condition {\rm(\ref{eq;21.7.28.2})}
as a mini-complex structure
which is an analogue of a complex structure.
(See {\rm\S\ref{subsection;17.10.28.1}}
for mini-complex structure.)
\hfill\qed
\end{rem}

There exist unitary frames
$u^{\pm}_1,u^{\pm}_2$ of $E$
such that the following holds.
\begin{itemize}
 \item $\nabla(u^{\pm}_1,u^{\pm}_2)=
       (u^{\pm}_1,u^{\pm}_2)(A^{\pm}_2\,dx_2+A^{\pm}_3\,dx_3)$,
       i.e., $\nabla_{x_1}u^{\pm}_i=0$.
       Moreover, $A^{\pm}_2,A^{\pm}_3\to 0$ as $x_1\to\pm\infty$.
 \item Let $\Phi^{\pm}$ be the $\su(2)$-valued functions
       determined by
       $\phi(u^{\pm}_1,u^{\pm}_2)=(u^{\pm}_1,u^{\pm}_2)\Phi^{\pm}$.
       Then, as $x_1\to\pm\infty$,
       $\Phi^{\pm}$ converge to
       a diagonal matrix 
       whose $(1,1)$-entry is $\sqrt{-1}$
       and $(2,2)$-entry is $-\sqrt{-1}$.
\end{itemize}
Let $(e^{\pm}_1,e^{\pm}_2)$ denote the asymptotic gauge
of $E$ at $x_1=\pm\infty$ induced by $(u^{\pm}_1,u^{\pm}_2)$.
They are well defined up to multiplications of
complex numbers with absolute value $1$ to $e^{\pm}_i$.

As explained in \cite{Hitchin-monopoles-geodesics},
there exists a frame $(v_1,v_2)$ of $E$
such that
(i) $\del_tv_i=0$,
(ii) $e^{t}t^{-k/2} v_1\to e^{+}_1$
and $e^{-t}t^{k/2} v_2\to e^+_2$ as $t\to\infty$.
It turns out that
$\nabla_{\zbar}v_1=0$.
Indeed,
by the commutativity (\ref{eq;21.7.28.2}),
we obtain $\del_t(\nabla_{\zbar}v_1)=0$.
Hence, we have the expression
$\nabla_{\zbar}v_1=\alpha_1(z)v_1+\beta_1(z)v_2$,
which implies
\[
 \nabla_{\zbar}(e^{t}t^{-k/2}v_1)
=\alpha_1(z)(e^{t}t^{-k/2}v_1)
+\beta_1(z)e^{2t}t^{-k}\cdot(e^{-t}t^{k/2}v_2).
\]
Because $\nabla_{\zbar}(e^{t}t^{-k/2}v_1)\to 0$ as $t\to\infty$,
we obtain $\alpha_1=\beta_1=0$.
Similarly,
$\nabla_{\zbar}v_2=\alpha_2(z) v_1$ holds
for a function $\alpha_2(z)$.
For any $t\in\real$,
let $L^{t}\subset E^t$ denote the subbundle
generated by
$f^t_{+,1}:=v_{1|\{t\}\times\cnum}$.
It is characterized by the following property.
\begin{itemize}
 \item For $(t_1,z_1)\in\real\times\cnum$
       and $s\in L^{t_1}_{+|z_1}\subset E^{t_1}_{|z_1}$,
       let $\stilde$ denote the section of
       $E_{|\real\times\{z_1\}}$ determined by the conditions
       $\del_t\stilde=0$
       and $\stilde_{|(t_1,z_1)}=s$.
       Then, $|\stilde|_h\to 0$ as $t\to\infty$.       
\end{itemize}
Because $\nabla_{\zbar}f^t_{+,1}=0$,
$L^{t}_+$ is a holomorphic subbundle of $E^t$,
and $f_{+,1}^t$ is a global holomorphic frame of $L^t_+$.
Moreover, $v_2$ induces a global holomorphic frame $f_{+,2}^t$
of $E^t/L^t_+$.
The scattering map $E^t\simeq E^{t'}$ induces
the isomorphisms $L^t_+\simeq L^{t'}_+$
and $E^t/L^t_+\simeq E^{t'}/L^{t'}_+$,
which preserve the distinguished frames.

By applying a similar consideration to the behaviour
as $t\to -\infty$,
we obtain a holomorphic subbundle $L^t_-\subset E^t$
characterized by the following condition.
\begin{itemize}
 \item For $(t_1,z_1)\in\real\times\cnum$
       and $s\in L^{t_1}_{-|z_1}\subset E^{t_1}_{|z_1}$,
       let $\stilde$ denote the section of
       $E_{|\real\times\{z_1\}}$ determined by the conditions
       $\del_t\stilde=0$
       and $\stilde_{|(t_1,z_1)}=s$.
       Then, $|\stilde|_h\to 0$ as $t\to-\infty$.
\end{itemize}
Moreover, the line bundles
$L^t_-$ and $E^t/L^t_-$ are equipped with
the distinguished frames
$f^t_{-,1}$ and $f^t_{-,2}$, respectively.
The subbundles and the distinguished frames
are preserved by the scattering maps.

There exists the naturally defined morphism
$\gamma:L^t_-\lrarr E^t/L^t_+$.
By the frames $f^t_{-,1}$ and $f^t_{+,2}$,
we regard $\gamma$ as an entire holomorphic function.
It turns out that $\gamma$ is a non-zero polynomial of degree $k$.
As remarked in \cite[\S16]{Atiyah-Hitchin-book},
there uniquely exists a holomorphic section $\ftilde^t_{+,2}$ of $E^t$
such that
(i) $\ftilde^t_{+,2}$ induces $f^t_{+,2}$,
(ii) for the expression
$f^t_{-,1}=\gamma \ftilde^t_{+,2}+\delta f^t_{+,1}$,
$\delta$ is a polynomial with $\deg\delta<k$.
Thus, we obtain the based holomorphic map
$\varphi=\delta/\gamma:\proj^1\lrarr\proj^1$.
Indeed, because of the ambiguity of
the choice of an asymptotic gauge $(e_1^{+},e_2^+)$,
we obtain a family of based holomorphic maps
parameterized by $S^1$.
Donaldson and Hitchin proved that this procedure induces
an equivalence between $\SU(2)$-monopoles with finite energy
and based holomorphic maps.

\begin{rem}
The construction also depends
on the choice of
an $\real$-linear isomorphism
$\real^3\simeq \real\times\cnum$
such that $\sum dx_i^2=dt\,dt+dw\,d\wbar$.
The choice of such an isomorphism corresponds
to a twistor parameter.
Hence, for each twistor parameter,
we obtain a family of based holomorphic maps
parameterized by $S^1$. 
\hfill\qed
\end{rem}

We may also formulate the associated algebraic object
in terms of $\cnum[z]$-modules,
where $\cnum[z]$ denotes the polynomial ring of the variable $z$.
Let $M$ be a free $\cnum[z]$-module of rank $2$.
Let $L_{\pm}\subset M$ be free $\cnum[z]$-submodules of rank $1$
such that $M/L_{\pm}$ are also free $\cnum[z]$-modules of rank $1$.
Let $g_{\pm,1}$ be frames of $L_{\pm}$,
and let $g_{2}$ be a frame of $M/L_{+}$.
Assume that the induced map
$L_-\lrarr M/L_+$ is expressed by a non-zero polynomial of degree $k$.

From the above monopole $(E,h,\nabla,\phi)$,
we obtain such $(M,L_{\pm},g_{\pm,1},g_2)$ as follows.
Let $\nbigo_{\proj^1}(\ast\infty)$ denote the sheaf of
meromorphic functions on $\proj^1$
which may allow poles along $\infty$.
We set
$\nbige=
 \nbigo_{\proj^1}(\ast\infty)f^0_{+,1}\oplus
 \nbigo_{\proj^1}(\ast\infty)\ftilde^0_{+,2}$.
It is equipped with the filtrations
$\nbigl_+=\nbigo_{\proj^1}(\ast\infty)f^0_{+,1}
\subset\nbige$
and
$\nbigl_-=\nbigo_{\proj^1}(\ast\infty)f^0_{-,1}
\subset\nbige$.
We set $M:=H^0(\proj^1,\nbige)$,
$L_{\pm}:=H^0(\proj^1,\nbigl_{\pm})$,
$g_{\pm,1}:=f^0_{\pm,1}$
and $g_2:=f^0_{+,2}$.
It is easy to see that
the associated based holomorphic map $\varphi$
is obtained from the associated $(M,L_{\pm},g_{\pm,1},g_2)$,
and vice versa.

\begin{rem}
In our study of periodic monopoles,
we prefer $\cnum[z]$-modules or $\nbigo$-modules
as algebraic objects
because we can easily add more structures
such as parabolic structures.
\hfill\qed
\end{rem}

\subsection{The correspondence due to Charbonneau and Hurtubise}
\label{subsection;21.8.2.2}

We also recall
the study of Charbonneau and Hurtubise \cite{Charbonneau-Hurtubise},
by omitting the details.
Let $\Sigma$ be a compact Riemann surface
equipped with a K\"ahler metric $g_{\Sigma}$.
We set $S^1=\real/T\seisuu$ for some $T>0$,
which is equipped with the standard metric $dt\,dt$.
Let $g$ denote the induced Riemannian metric of
$S^1\times\Sigma$.
Let $Z=\{(t_i,p_i)\,|\,i=1,\ldots,N\}
\subset S^1\times\Sigma$ be a finite subset.
For simplicity, we assume that
$t_i\neq t_j$ and $p_i\neq p_j$ if $i\neq j$.
We also assume $t_i\neq 0$.

Let $(E,h,\nabla,\phi)$ be a monopole on $(S^1\times\Sigma)\setminus Z$
such that each point of $Z$ is a Dirac type singularity
of $(E,h,\nabla,\phi)$.
(We shall review the Dirac type singularity of monopoles
in the Euclidean case in \S\ref{subsection;21.8.3.3}.)
We recall the way to obtain an algebraic object in this context.
Let $p_{\Sigma}:S^1\times \Sigma\lrarr \Sigma$ denote the projection.
Let $\nabla^{0,1}_{|\Sigma}:E\lrarr E\otimes p_{\Sigma}^{-1}\Omega_{\Sigma}^{0,1}$
denote the induced differential operator.
We set $\del_t:=\nabla_t-\sqrt{-1}\phi$.
The Bogomolny equation implies
the commutativity
$[\del_t,\nabla^{0,1}_{|\Sigma}]=0$.

Let $Z_1\subset S^1$ and $Z_2\subset \Sigma$
denote the image by the projections of $Z$.
For any $t\in S^1\setminus Z_1$,
we obtain the vector bundle
$E^t=E_{|\{t\}\times \Sigma}$ on $\Sigma$,
which is enhanced to a holomorphic bundle
by $\nabla^{0,1}_{|\Sigma}$.
As the parallel transport by $\del_t$,
we obtain the isomorphisms
$\Phi_{t',t}:
E^{t}_{|\Sigma\setminus Z_2}
\simeq
E^{t'}_{|\Sigma\setminus Z_2}$
for $t,t'\in S^1\setminus Z_1$,
called {\em the scattering maps}.
\index{scattering map}
By the commutativity
$[\del_t,\nabla^{0,1}_{|\Sigma}]=0$,
the maps are holomorphic.
Under the assumption that each point of $Z$
is a Dirac type singularity of $(E,h,\nabla,\phi)$,
it turns out that the scattering maps are at most meromorphic
along $Z_2$.
Thus, we obtain a holomorphic vector bundle $E^0$
with a meromorphic automorphism
$\Phi_{T,0}:E^0_{|\Sigma\setminus Z_2}\simeq
E^T_{|\Sigma\setminus Z_2}=E^0_{|\Sigma\setminus Z_2}$.
Such a pair is called a bundle pair in \cite{Charbonneau-Hurtubise}.

Let $\nbige$ be a holomorphic vector bundle on $\Sigma$ of rank $r$
equipped with a meromorphic automorphism
$\rho:\nbige_{|\Sigma\setminus Z_2}\simeq\nbige_{|\Sigma\setminus Z_2}$.
For each $p_i\in Z_2$,
there exist frames $v^{\pm}_1,\ldots,v^{\pm}_r$ of
$\nbige$ around $p_i$,
and a decreasing sequence of integers
$k_{i,1},\ldots,k_{i,r}$ 
such that
$\rho(v_j^-)=z^{k_{i,j}}v_j^+$,
where $z$ denotes a holomorphic coordinate around $p_i$
such that $z(p_i)=0$.
The tuple $(k_{i,j}\,|\,j=1,\ldots,r)$ is well defined
for each $p_i\in Z_2$.
We note that
$\sum_{j=1}^r\sum_{i=1}^N k_{i,j}=0$.
Then, Charbonneau and Hurtubise
introduced the degree for such a bundle pair $(\nbige,\rho)$,
depending on $Z$:
\[
 \deg(\nbige,\rho)=
 \deg(\nbige)
-\sum_{i=1}^N
 \frac{t_i}{T}\sum_{j=1}^r k_{i,j}.
\]
Let $\nbige'\subset\nbige$ be a subbundle
such that $\rho(\nbige'_{|\Sigma\setminus Z_2})
=\nbige'_{|\Sigma\setminus Z_2}$.
Then, we obtain the induced meromorphic automorphism
$\rho'$ of $\nbige'$,
and $\deg(\nbige',\rho')$ is defined similarly.
We say that $(\nbige,\rho)$ is stable
if
\[
 \frac{\deg(\nbige',\rho')}{\rank\nbige'}
<\frac{\deg(\nbige,\rho)}{\rank\nbige}
\]
for any proper subbundle $\nbige'\subset\nbige$
such that $\rho(\nbige'_{|\Sigma\setminus Z_2})=\nbige'_{|\Sigma\setminus Z_2}$.
We say that $(\nbige,\rho)$ is poly-stable
if it is a direct sum
$(\nbige,\rho)=\bigoplus(\nbige_{\ell},\rho_{\ell})$
such that (i) each $(\nbige_{\ell},\rho_{\ell})$ is stable,
(ii) $\deg(\nbige_{\ell},\rho_{\ell})/\rank(\nbige_{\ell})
=\deg(\nbige,\rho)/\rank(\nbige)$ for any $\ell$.
Charbonneau and Hurtubise proved that
if $(\nbige,\rho)$ is induced by
a monopole $(E,\nabla,h,\phi)$ with Dirac type singularity,
then $(\nbige,\rho)$ is polystable of degree $0$.
By using the fundamental result of Simpson in \cite{Simpson88},
they proved that this procedure induces an equivalence
between monopoles and polystable bundle pairs of degree $0$.

\begin{rem}
Let $(E,h,\nabla,\phi)$ be a vector bundle
$E$ on $(S^1\times\Sigma)\setminus Z$
with a Hermitian metric $h$, a unitary connection $\nabla$
and an anti-Hermitian endomorphism $\phi$.
As a generalization of the Bogomolny equation, 
Charbonneau and Hurtubise
studied the Hermitian-Einstein-Bogomolny equation
\begin{equation}
\label{eq;21.8.2.3}
 F(\nabla)-\sqrt{-1}c\omega_{\Sigma}\id_E=\ast\nabla\phi,
\end{equation}
where $\omega_{\Sigma}$ denotes the K\"ahler form
associated with $g_{\Sigma}$.
If $c=0$, it is the Bogomolny equation.
\hfill\qed
\end{rem}

\subsection{Remark}
\label{subsection;21.8.2.3}

Though the scattering maps are efficiently used in the both constructions
of algebraic objects
in {\rm\S\ref{subsection;21.8.2.1}}
and {\rm\S\ref{subsection;21.8.2.2}},
there are also different flavors.
Roughly speaking,
in \S\ref{subsection;21.8.2.1},
we obtain a module with filtrations on $Y$ 
from a monopole on $\real\times Y$,
and in \S\ref{subsection;21.8.2.2},
we obtain a module with an automorphism on $Y$
from a monopole on $S^1\times Y$.

In this monograph,
we construct a difference module with parabolic structure from a monopole
depending on the twistor parameter $\lambda$.
In the case $\lambda=0$, our construction is exactly an analogue of
that in \S\ref{subsection;21.8.2.2}
in the context of $\Sigma=\cnum$
though we need to study the additional issues
caused by the non-compactness of $\Sigma$.
In the case $\lambda\neq 0$, our construction is again
similar to that in {\rm\S\ref{subsection;21.8.2.2}},
but it is more complicated.

We also note that in the case $\lambda\neq 0$
it is also natural to expect to obtain
another algebraic object as the associated ``Betti'' object,
and the construction should be in some sense similar
to the construction in {\rm\S\ref{subsection;21.8.2.1}}
though we shall not discuss it in this monograph.
Indeed, in our study of doubly periodic monopoles \cite{Mochizuki-q-difference},
two kinds of the constructions appear.
From a monopole on the product of $\real$ and an elliptic curve,
we construct a $\gminiq^{\lambda}$-difference module with parabolic structure
depending on a twistor parameter $\lambda$
\cite[Theorem 1.5]{Mochizuki-q-difference}.
If $|\gminiq^{\lambda}|\neq 1$, or equivalently $|\lambda|\neq 1$,
by the Riemann-Hilbert correspondence
for $\gminiq^{\lambda}$-difference modules
in the local case due to
van der Put and Reversat \cite{van-der-Put-Reversat}
and  
Ramis, Sauloy and Zhang \cite{Ramis-Sauloy-Zhang}
and in the global case due to Kontsevich and Soibelman,
the monopole also induces a locally free sheaf
with some filtrations on the elliptic curve
$T^{\lambda}$
obtained as the quotient of $\cnum^{\ast}$
by the action of $(\gminiq^{\lambda})^{\seisuu}$
\cite[Theorem 1.7]{Mochizuki-q-difference}.
The construction of the $\gminiq^{\lambda}$-difference module
from the monopole
is similar to the construction in {\rm\S\ref{subsection;21.8.2.2}},
and the construction of the locally free module with filtrations on
$T^{\lambda}$ is similar to the construction in
{\rm \S\ref{subsection;21.8.2.1}}.

\section[Review of KH-correspondences for $\lambda$-flat bundles]{Review of the Kobayashi-Hitchin correspondences for
$\lambda$-flat bundles}
\label{section;20.8.9.2}

Recall that Simpson established the Kobayashi-Hitchin correspondence
for tame harmonic bundles on non-compact curves
\cite{Simpson90}.
It is our main goal to develop an analogue theory 
in the context of periodic monopoles.
Hence, let us begin with a brief review on the theory of
harmonic bundles on curves.

\subsection{Harmonic bundles and their underlying
 $\lambda$-flat bundles}

Let $C$ be any complex curve.
Let $(E,\delbar_E)$ be a holomorphic vector bundle on $C$.
Let $\theta$ be a Higgs field,
i.e., a holomorphic section of $\End(E)\otimes\Omega^1_C$.
Let $h$ be a Hermitian metric of $E$.
We obtain the Chern connection 
$\nabla_h=\delbar_E+\del_{E,h}$ determined by
$\delbar_E$ and $h$,
and the adjoint $\theta^{\dagger}_h$ of $\theta$
with respect to $h$.
\index{harmonic metric}
The metric $h$ is called harmonic
if the Hitchin equation
\begin{equation}
\label{eq;21.8.11.1}
 [\delbar_E,\del_{E,h}]+[\theta,\theta^{\dagger}_h]=0
\end{equation}
is satisfied.
\index{Hitchin equation}
A Higgs bundle with a harmonic metric
is called a harmonic bundle.
\index{harmonic bundle}
The Hitchin equation implies that
the connection 
$\delbar_E+\theta_h^{\dagger}+\del_{E,h}+\theta$
is flat.
More generally, for any complex number $\lambda$,
there exists the associated flat $\lambda$-connection.
Indeed, we set
$\DDlambda:=
 \delbar_E+\lambda\theta^{\dagger}
+\lambda\del_E+\theta$.
As explained in \cite{s3,s4},
it is a $\lambda$-connection in the sense of Deligne,
i.e.,
it satisfies a twisted Leibniz rule
$\DDlambda(fs)=
 (\lambda\del_C+\delbar_C)f\cdot s
+f\DDlambda(s)$
for any $f\in C^{\infty}(C,\cnum)$
and $s\in C^{\infty}(C,E)$.
\index{$\lambda$-connection}
The Hitchin equation means that
$\DDlambda$ is flat,
i.e.,
$\DDlambda\circ\DDlambda=0$.
\index{flat $\lambda$-connection}
Hence, harmonic bundles have the underlying
$\lambda$-flat bundles.
We may define the concept of harmonic metrics
for $\lambda$-flat bundles.
A $\lambda$-flat bundle with a harmonic metric
is equivalent to a Higgs bundle with a harmonic metric.
See \cite{Mochizuki-KHII},
for example.

\subsection{Kobayashi-Hitchin correspondences
in the smooth case}

Suppose that  $C$ is projective and connected.
We set $\deg(F)=\int_Cc_1(F)$
for any vector bundle $F$ on $C$.
\index{degree $\deg(F)$}
A $\lambda$-flat bundle $(V,\nabla^{\lambda})$
is called stable (semistable)
if we obtain
\[
 \deg(V')/\rank V'<(\leq)\deg(V)/\rank V
\]
for any $\lambda$-flat subbundle 
$(V',\nabla^{\lambda})
\subset (V,\nabla^{\lambda})$
with $0<\rank(V')<\rank(V)$.
A $\lambda$-flat bundle $(V,\nabla^{\lambda})$
is called polystable
if it is a direct sum
of stable $\lambda$-flat bundles
$\bigoplus (V_i,\nabla^{\lambda})$
such that
$\deg(V)/\rank V=\deg(V_i)/\rank V_i$
for any $i$.
Note that $\deg(V)=0$ always holds if $\lambda\neq 0$.
\index{stable ($\lambda$-flat bundle)}
\index{semistable ($\lambda$-flat bundle)}
\index{polystable ($\lambda$-flat bundle)}

The following is a fundamental theorem
in the study of harmonic bundles on projective curves,
due to 
Diederich-Ohsawa \cite{Diederich-Ohsawa},
Donaldson \cite{don2} and
Hitchin \cite{Hitchin-self-duality},
in the rank $2$ case,
and Corlette \cite{corlette} and 
Simpson \cite{Simpson88, s4,s5}
in the higher rank case,
or even in the higher dimensional case.
(See also \cite{Mochizuki-KHII} for the case of
general flat $\lambda$-connections.)

\begin{thm}
\label{thm;17.10.27.1}
Suppose that $C$ is projective and connected.
A $\lambda$-flat bundle $(V,\nabla^{\lambda})$
is polystable of degree $0$
if and only if $(V,\nabla^{\lambda})$ 
has a harmonic metric.
If $(V,\nabla^{\lambda})$ has two harmonic metrics
$h_j$ $(j=1,2)$,
there exists a decomposition
$(V,\nabla^{\lambda})=
 \bigoplus (V_i,\nabla^{\lambda})$
which is orthogonal with respect to both $h_1$ and $h_2$,
such that $h_{1|V_i}=a_i\cdot h_{2|V_i}$
for some constants $a_i>0$.
In particular, if $(V,\nabla^{\lambda})$ is stable,
it has a unique harmonic metric
up to the multiplication of positive constants.
\hfill\qed
\end{thm}

\begin{cor}
For each $\lambda$,
there exists a natural bijective correspondence
between the equivalence classes of
polystable $\lambda$-flat bundles of degree $0$
and the equivalence classes of
polystable Higgs bundles of degree $0$,
through harmonic bundles.
\hfill\qed
\end{cor}

The Kobayashi-Hitchin correspondence in the theorem
provides us with an interesting equivalence
between objects in differential geometry
and objects in algebraic geometry.
It is an origin of the hyperk\"ahler property of
the moduli spaces.
It is a starting point of the non-abelian Hodge theory
of Simpson. (See \cite{s4}.)

\subsection{Tame harmonic bundles and
regular filtered $\lambda$-flat bundles}
\label{subsection;17.10.27.10}

Simpson studied a generalization of 
Theorem \ref{thm;17.10.27.1}
to the context of harmonic bundles on quasi-projective curves
in \cite{Simpson90}.
To state his result,
let us recall the concepts
of tame harmonic bundles
and regular filtered $\lambda$-flat bundles.

Let $Y$ be a neighbourhood of $0$ in $\cnum$.
Let $(E,\delbar_E,\theta,h)$ be a harmonic bundle
on $Y\setminus\{0\}$.
We obtain the holomorphic endomorphism $f$ of $E$
such that $\theta=f\,dz/z$.
The harmonic bundle is called tame on $(Y,\{0\})$
if the following holds.
\index{tame harmonic bundle}
\begin{itemize}
\item 
 Let $\det(t\id_E-f)=
 t^{\rank E}+\sum_{j=0}^{\rank E-1} a_j(z)t^j$
 be the characteristic polynomial of $f$.
 Then, $a_j$ are holomorphic at $z=0$.
\end{itemize}

Let us explain the notion of
regular filtered $\lambda$-flat bundles on $(Y,\{0\})$.
Let $\nbigo_Y(\ast 0)$ denote the sheaf of
meromorphic functions on $Y$ which may have poles along $0$.
\index{sheaf $\nbigo_Y(\ast 0)$}
Recall that a filtered bundle on $(Y,\{0\})$
is a locally free $\nbigo_Y(\ast 0)$-module
$\nbigv$ of finite rank
with an increasing sequence of 
locally free $\nbigo_Y$-submodules
$\nbigp_a\nbigv\subset\nbigv$
$(a\in\real)$
satisfying the following conditions.
\index{filtered bundle}
\begin{itemize}
 \item
      $\nbigp_a\nbigv$ $(a\in\real)$ are lattices of
      $\nbigv$, i.e.,
      $\nbigp_a\nbigv(\ast \{0\})=\nbigv$.
\item
 $\nbigp_{a+n}\nbigv=\nbigp_a\nbigv(n\{0\})$
 for any $a\in\real$ and $n\in\seisuu$.
\item
 For any $a\in\real$,
 there exists  $\epsilon>0$
 such that
 $\nbigp_{a+\epsilon}\nbigv=\nbigp_a\nbigv$.
\end{itemize}
A regular filtered $\lambda$-flat bundle
is a filtered bundle $\nbigp_{\ast}\nbigv$
with a $\lambda$-connection
$\nabla^{\lambda}:
 \nbigv\lrarr\nbigv\otimes\Omega^1$
such that
$\nabla^{\lambda}$ is logarithmic
with respect to $\nbigp_{\ast}\nbigv$
in the sense of 
$\nabla^{\lambda}\nbigp_a\nbigv
\subset
\nbigp_a\nbigv\otimes\Omega^1(\log \{0\})$
for any $a\in\real$.
\index{regular filtered $\lambda$-flat bundle}
Note that
we obtain the finite dimensional $\cnum$-vector spaces
$\Gr^{\nbigp}_a(\nbigv):=
 \nbigp_{a}(\nbigv)\big/
 \sum_{b<a}\nbigp_b(\nbigv)$. 
\index{vector space $\Gr^{\nbigp}_a(\nbigv)$}
 
Let $(E,\delbar_E,\theta,h)$
be a tame harmonic bundle
on $Y\setminus\{0\}$.
For any complex number $\lambda$,
we obtain the holomorphic vector bundle
$(E,\delbar_E+\lambda\theta^{\dagger})$
on $Y\setminus\{0\}$.
Let $E^{\lambda}$ denote 
the locally free $\nbigo_{Y\setminus\{0\}}$-module
obtained as the sheaf of holomorphic sections of
$(E,\delbar_E+\lambda\theta^{\dagger})$.
It is equipped with a flat $\lambda$-connection
$\DDlambda$.
For any open subset $U\ni 0$,
let $\nbigp^hE^{\lambda}(U)$ denote the space of
holomorphic sections $s$ of $E^{\lambda}$ on $U\setminus\{0\}$
such that
$|s|_h=O(|z|^{-N})$ for some $N$,
where $|s|_h$ denotes the norm of $s$ with respect to $h$.
Thus, we obtain an $\nbigo_{Y}(\ast 0)$-module
$\nbigp^hE^{\lambda}$
as a meromorphic prolongation of $E^{\lambda}$.
\index{sheaf $\nbigp^hE^{\lambda}$}
For any $a\in\real$,
let $\nbigp^h_aE^{\lambda}(U)$ denote the space of
sections $s$ of $E^{\lambda}$ on $U\setminus\{0\}$
such that
$|s|_h=O(|z|^{-a-\epsilon})$ for any $\epsilon>0$.
Thus, we obtain an increasing sequence of
$\nbigo_Y$-submodules
$\nbigp^h_aE^{\lambda}\subset\nbigp^hE^{\lambda}$ $(a\in\real)$.
\index{sheaf $\nbigp^h_aE^{\lambda}$}
Simpson proved that
$(\nbigp^h_{\ast}E^{\lambda},\DDlambda)$ is
a regular filtered $\lambda$-flat bundle
in \cite{Simpson90}.

\begin{rem}
\label{rem;21.9.12.2}
We prefer to consider
filtered $\lambda$-flat bundle
$(\nbigp^h_{\ast}E^{\lambda},\DDlambda)$ on $(Y,0)$
rather than $(E^{\lambda},\DDlambda)$ on $Y\setminus \{0\}$
to keep the information of the behaviour of $h$ around $0$.
\hfill\qed
\end{rem}

Let $C$ be a smooth connected projective curve
with a finite subset $D\subset C$.
The concept of tame harmonic bundles on $(C,D)$
is defined in a natural way.
Let $\nbigo_C(\ast D)$ denote the sheaf of
meromorphic functions on $C$ which may have poles along $D$.
\index{sheaf $\nbigo_C(\ast D)$}
A filtered bundle on $(C,D)$
is a locally free $\nbigo_C(\ast D)$-module
$\nbigv$ of finite rank
with an increasing sequence of
$\nbigo_C$-locally free submodules
$\nbigp_{\veca}\nbigv$ $(\veca=(a_P\,|\,P\in D)\in\real^D)$
such that the following holds.
\index{filtered bundle}
\begin{itemize} 
\item For any $P\in D$,
 take a small neighbourhood $U_P$ of $P$ in $C$.
 Then, 
 $\nbigp_{\veca}\nbigv_{|U_P}$ depends only on  $a_P$,
 which we denote by 
 $\nbigp^{(P)}_{a_P}(\nbigv_{|U_P})$.
\item
 $\nbigp^{(P)}_{\ast}(\nbigv_{|U_P})$
 is a filtered bundle on $(U_P,P)$
 in the above sense.
\end{itemize}
A regular filtered $\lambda$-flat bundle
is a filtered bundle $\nbigp_{\ast}\nbigv$
with a $\lambda$-connection $\nabla^{\lambda}$
such that 
$\nabla^{\lambda}\nbigp_{\veca}\nbigv
 \subset
 \nbigp_{\veca}\nbigv\otimes\Omega^1_C(\log D)$
 for any $\veca\in\real^D$.
\index{regular filtered $\lambda$-flat bundle}
Moreover, for any regular filtered $\lambda$-flat bundle
$(\nbigp_{\ast}\nbigv,\nabla^{\lambda})$ on $(C,D)$,
we define
\index{degree $\deg(\nbigp_{\ast}\nbigv)$}
\[
 \deg(\nbigp_{\ast}\nbigv):=
 \deg(\nbigp_{\veczero}\nbigv)
-\sum_{P\in D}
 \sum_{-1<b\leq 0}
 b\dim_{\cnum}
 \Gr^{\nbigp^{(P)}}_b(\nbigv_{|U_P}).
\]
Here, $\veczero=(0,\ldots,0)\in\real^D$.
By using the degrees, we can define 
the stability condition and the polystability condition
in the context of regular filtered $\lambda$-flat bundles
in the natural way.

Any tame harmonic bundle on $(C,D)$
naturally induces a regular filtered $\lambda$-flat bundle
on $(C,D)$
by the procedure explained above.
Simpson established the following theorem 
in \cite{Simpson90},
which is the ideal generalization of
Theorem \ref{thm;17.10.27.1}
to the context of regular singular case.
(See also \cite{Mochizuki-KHI, Mochizuki-KHII}
for the case of $\lambda$-connections
and the higher dimensional case.)

\begin{thm}[Simpson]
\label{thm;17.10.27.2}
Suppose that $C$ is projective and connected.
The procedure induces a bijection
between the equivalence classes of
tame harmonic bundles on $(C,D)$
and the equivalence classes of
polystable regular filtered $\lambda$-flat bundles
with degree $0$.
\hfill\qed
\end{thm}

\subsection{Wild harmonic bundles
and good filtered $\lambda$-flat bundles}
\label{subsection;17.10.27.11}

Let us give a complement on
the generalization of Theorem \ref{thm;17.10.27.2}
to the context of non-regular case,
which was studied by Biquard-Boalch \cite{biquard-boalch}.
(See also \cite{Mochizuki-wild, Mochizuki-KH-Higgs}.)

Let $Y$ be a neighbourhood of $0$ in $\cnum$.
A harmonic bundle $(E,\delbar_E,\theta,h)$ on $Y\setminus\{0\}$
is called wild on $(Y,0)$ if the following holds.
\index{wild harmonic bundle}
\begin{itemize}
\item
     $a_j(z)$ are meromorphic at $z=0$.
     (See \S\ref{subsection;17.10.27.10}
     for the construction of $a_j$ from $(E,\delbar_E,\theta,h)$.)
\end{itemize}

To explain the concept of good filtered $\lambda$-flat bundles,
we make several preliminaries.
First, we recall the concept of unramifiedly good filtered 
$\lambda$-flat bundle.
\index{unramifiedly good filtered $\lambda$-flat bundle}
A filtered bundle $\nbigp_{\ast}\nbigv$ 
with a $\lambda$-connection $\nabla^{\lambda}$
on $(Y,\{0\})$ is called unramifiedly good
if there exist a subset $S\subset z^{-1}\cnum[z^{-1}]$
and the Hukuhara-Levelt-Turrittin decomposition 
\index{Hukuhara-Levelt-Turrittin decomposition}
\begin{equation}
\label{eq;21.8.6.1}
 (\nbigp_{\ast}\nbigv,\nabla^{\lambda})_{|\widehat{0}}
=\bigoplus_{\gminia\in S} 
 (\nbigp_{\ast}\nbigvhat_{\gminia},\nablahat^{\lambda}_{\gminia})
\end{equation}
such that 
$\nablahat^{\lambda}_{\gminia}-d\gminia\cdot\id_{\nbigvhat_{\gminia}}$
are logarithmic with respect to
$\nbigp_{\ast}\nbigvhat_{\gminia}$.
Here, for any $\nbigo_Y$-module $\nbigf$,
let $\nbigf_{|\widehat{0}}$
denote the formal completion of the stalk $\nbigf_0$ at $0$,
i.e., 
$\nbigf_0\otimes_{\nbigo_{Y,0}}\cnum[\![z]\!]$.
\index{formal completion $\nbigf_{|\widehat{0}}$}

Second,
let us recall the pull back of filtered bundles.
Take $q\in\seisuu_{\geq 1}$.
Let $\varphi_q:\cnum\lrarr \cnum$ be the ramified covering
given by $\varphi_q(\zeta)=\zeta^q$.
We set $Y_q=\varphi_q^{-1}(Y)$.
Let $\nbigp_{\ast}\nbigv$ be a filtered bundle on $(Y,\{0\})$.
We obtain the locally free $\nbigo_{Y_q}(\ast \{0\})$-module
$\nbigv':=\varphi_q^{\ast}\nbigv$.
We define 
\[
 \nbigp_a\nbigv':=
 \sum_{n+qb\leq a} \zeta^{-n}\varphi_q^{\ast}\nbigp_b\nbigv
\subset\nbigv'.
\]
Thus, we obtain a filtered bundle
$\nbigp_{\ast}\nbigv'$
over $(Y_q,\{0\})$
which is 
called the pull back of $\nbigp_{\ast}\nbigv$,
and denoted by
$\varphi_q^{\ast}\bigl(
\nbigp_{\ast}\nbigv\bigr)$.
\index{pull back $\varphi_q^{\ast}(\nbigp_{\ast}\nbigv)$}

A filtered bundle $\nbigp_{\ast}\nbigv$
with a $\lambda$-connection $\nabla^{\lambda}$
is called good
if there exists a ramified covering $\varphi_q:Y_q\lrarr Y$ as above
for an appropriate $q$
such that 
$\varphi_q^{\ast}(\nbigp_{\ast}\nbigv,\nabla^{\lambda})$
is unramifiedly good.
\index{good filtered $\lambda$-flat bundle}

Let $(E,\delbar_E,\theta,h)$
be a wild harmonic bundle on $(Y,0)$.
For any $\lambda$,
we obtain the $\lambda$-flat bundle
$(E^{\lambda},\DDlambda)$ on $Y\setminus\{0\}$
as above.
We obtain the $\nbigo_{Y}(\ast \{0\})$-module
$\nbigp^hE^{\lambda}$
and the increasing sequence of
$\nbigo_Y$-submodules
$\nbigp^h_aE^{\lambda}$ 
by the procedure explained in \S\ref{subsection;17.10.27.10}.
As proved in \cite{Mochizuki-wild},
$(\nbigp^h_{\ast}E^{\lambda},\DDlambda)$
is a good filtered $\lambda$-flat bundle.

Although we explained the concepts of 
wild harmonic bundles
and good filtered $\lambda$-flat bundles
on a neighbourhood of $0$ in $\cnum$,
they are naturally generalized
to the context of 
any complex curve $C$ with a discrete subset $D$,
and any wild harmonic bundle $(E,\delbar_E,\theta,h)$
on $(C,D)$
induces a good filtered $\lambda$-flat bundle
$(\nbigp_{\ast}\nbigelambda,\DDlambda)$.
The following is essentially proved in
\cite{biquard-boalch}.
(See also \cite{Mochizuki-wild, Mochizuki-KH-Higgs}.)

\begin{thm}
\label{thm;17.10.28.30}
Suppose that $C$ is projective and connected.
Then, the above procedure induces
a bijection between the equivalence classes of
wild harmonic bundles on $(C,D)$
and the equivalence classes of
polystable good filtered $\lambda$-flat bundles
with degree $0$.
\hfill\qed
\end{thm}

\section{Equivariant instantons 
and the underlying holomorphic objects}
\label{section;21.8.13.20}

In \S\ref{subsection;21.8.11.4}--\S\ref{subsection;21.8.11.5},
we explain our motivation to study
the Kobayashi-Hitchin correspondences for periodic monopoles
as analogue of the Kobayashi-Hitchin correspondences for
harmonic bundles (see \S\ref{section;20.8.9.2}).
The precise statements are postponed 
to \S\ref{section;21.8.11.2}--\S\ref{section;21.8.11.3}.
In \S\ref{subsection;21.8.4.22},
we recall that monopoles are regarded
as harmonic bundles of infinite rank,
which is a basic idea behind our study.
In particular, it is significant
when we study the asymptotic behaviour of periodic monopoles
(see \S\ref{section;21.8.4.21}).

\subsection{Instantons and the underlying holomorphic bundles}
\label{subsection;21.8.11.4}

Set $X:=\real^4$ with the standard Euclidean metric
$\sum_{i=1}^4 dx_i^2$.
An instanton on an open subset $U\subset X$
is a vector bundle $E$ with a Hermitian metric $h$
and a unitary connection $\nabla$
satisfying the anti-self duality (ASD) equation:
\index{instanton}
\index{ASD equation}
\begin{equation}
\label{eq;21.8.4.10}
 F(\nabla)+\ast F(\nabla)=0.
\end{equation}

We take an $\real$-linear isomorphism
$X\simeq\cnum^2=\{(z,w)\}$
such that the metric is $dz\,d\zbar+dw\,d\wbar$,
with which we regard $X$ as a K\"ahler surface.
The ASD equation (\ref{eq;21.8.4.10}) holds
if and only if
(i) the curvature is a $(1,1)$-form with respect to the complex structure,
and (ii) $\Lambda F(\nabla)=0$.
(See \cite{Kobayashi-vector-bundle} for the operator $\Lambda$.)
The $(0,1)$-part of $\nabla$ induces
a holomorphic structure of $E$.
Hence, from an instanton,
we obtain a holomorphic vector bundle
$(E,\delbar_E)$ with a Hermitian metric $h$
satisfying $\Lambda F(h)=0$,
where $F(h)$ denote the curvature of the Chern connection
determined by $h$ and $\delbar_E$.

Recall that such complex structures on $X$
are parameterized by $\proj^1$.
For each $\lambda\in\proj^1\setminus\{\infty\}$,
we obtain the complex coordinate system
$(\xi,\eta)=(z+\lambda\wbar,w-\lambda\zbar)$
which induces the complex structure corresponding to $\lambda$.
We obtain the complex manifold
$X^{\lambda}$ and the open subset $U^{\lambda}$.
An instanton $(E,h,\nabla)$ on $U$
induces a holomorphic vector bundle
$(E^{\lambda},\delbar_{E^{\lambda}})$
on $U^{\lambda}$.
\index{underlying holomorphic bundle $(E^{\lambda},\delbar_{E^{\lambda}})$}

Let $\Gamma$ be a closed subgroup of $\real^4$,
which naturally acts on $\real^4$ by the translation.
Suppose that $U$ is invariant under the $\Gamma$-action.
If $(E,h,\nabla)$ is $\Gamma$-equivariant,
then the underlying holomorphic bundles
$(E^{\lambda},\delbar_{E^{\lambda}})$ $(\lambda\in\cnum)$
are also $\Gamma$-equivariant.

\subsection{Instantons and harmonic bundles}
\label{subsection;21.8.4.20}

Recall that the concept of harmonic bundles
was discovered by Hitchin \cite{Hitchin-self-duality}
as the $2$-dimensional reduction of instantons.
Namely, 
in the case $\Gamma=\real^2$,
$\Gamma$-equivariant instantons on $U$
are equivalent to harmonic bundles on $U/\Gamma$.
Indeed, we may choose 
$\Gamma=\cnum_z\times\{0\}\subset \cnum_z\times\cnum_w$.
We may regard $U_0=U/\cnum_z$
as an open subset in $\cnum_w$.
Let $p:U\lrarr U_0$ denote the projection.
We obtain the induced vector bundle $E_0$ on $U_0$
with a $\cnum_z$-equivariant isomorphism $p^{\ast}E_0\simeq E$.
It is equipped with the induced metric $h_0$.
The $\cnum_z$-equivariant differential operator $\nabla_{\wbar}$ on $E$
naturally induces a holomorphic structure $\del_{E_0,\wbar}$
on $E_0$.
We obtain a holomorphic endomorphism $f_0$ on $E_0$
such that
$\nabla_{\zbar}(p^{\ast}(s))=p^{\ast}(f_0s)$
for any $s\in C^{\infty}(U_0,E_0)$.
Then, the equation $\Lambda F(h)=0$
is reduced to the Hitchin equation (\ref{eq;21.8.11.1})
for the metric $h_0$
of the Higgs bundle
$(E_0,\delbar_{E_0},f_0dw)$.

Note that $\cnum_z$-equivariant
holomorphic vector bundles on $U^{\lambda}$
are equivalent to 
$\lambda$-flat bundles on $U_0$.
We can easily check that 
the $\cnum_z$-equivariant holomorphic bundle
$(E^{\lambda},\delbar_{E^{\lambda}})$ on $U^{\lambda}$
corresponds to
the $\lambda$-flat bundle
$(E^{\lambda}_0,\DD^{\lambda})$ on $U_0$
associated with $(E,\delbar_E,\theta,h)$:
\[
\begin{CD}
\left(
 \mbox{
 \begin{tabular}{c}
  $\cnum_z$-equivariant\\
  instantons on $U$
 \end{tabular}}
 \right)
 @>>>
\left(
 \mbox{
 \begin{tabular}{c}
  $\cnum_z$-equivariant holomorphic\\
  vector bundles on $U^{\lambda}$
 \end{tabular}}
 \right)
 \\
 @V{=}VV @V{=}VV \\
\left(
 \mbox{\begin{tabular}{c}
	harmonic bundles \\
	on $U_0$
       \end{tabular}}
 \right)
 @>>>
\left(
 \mbox{\begin{tabular}{c}
	$\lambda$-flat bundles \\
	on $U_0$
	\end{tabular}}
\right).
\end{CD}
\]

Hence, the Kobayashi-Hitchin correspondence
between harmonic bundles and $\lambda$-flat bundles
in the case $X=\proj^1$
can be regarded as a correspondence
between $\cnum_z$-equivariant instantons 
and $\cnum_z$-equivariant holomorphic vector bundles
on $p^{-1}(\cnum_w\setminus D)$,
under some natural assumptions on the boundary behaviour
of the objects.

This kind of problems have been studied in various cases
of equivariant instantons,
sometimes in relation with the Nahm transforms.
For example, see
\cite{Biquard-Jardim},
\cite{Donaldson-GIT, Donaldson-Nahm},
\cite{Hurtubise-Murray1, Hurtubise-Murray2},
\cite{Jarvis1, Jarvis2},
\cite{Mochizuki-doubly-periodic}, \cite{Takayama},
etc.
It is useful for the classification of
equivariant instantons
because holomorphic objects are easier to study.
Conversely, it is useful for the study of deeper aspects
of such holomorphic objects,
for example the hyperk\"ahler property of the moduli spaces.
To establish the correspondence,
we need a careful study on the asymptotic behaviour of 
equivariant instantons depending on $\Gamma$,
and hence there seem to remain many things to be clarified
mathematically.

\subsection{Instantons and monopoles}
\label{subsection;21.8.11.5}

It is our purpose in this monograph is to study the issue
in the case $\Gamma\simeq\real\times\seisuu$.
We recall that monopoles are the $1$-dimensional reduction of
instantons.
Indeed, if $\Gamma\simeq \real\times\seisuu$,
then $\Gamma$-equivariant instantons on $U$
are equivalent to monopoles on $U/\Gamma$.
We can choose the coordinate system $(x_1,x_2,x_3,x_4)$ as
$\Gamma=\{(x_1,nT,0,0)\,|\,x_1\in\real,n\in\seisuu\}$ for some $T>0$.
We may regard
$U_1=U/\Gamma$ as an open subset of
$(\real/T\seisuu)\times\real^2_{x_3,x_4}$.
Let $p_1:U\lrarr U_1$ denote the projection.
If $(E,h,\nabla)$ is a $\Gamma$-equivariant tuple of
a vector bundle $E$ with a Hermitian metric
and a unitary connection $\nabla$ on $U$,
we obtain the vector bundle $E_1$ on $U_1$
with a $\Gamma$-equivariant isomorphism
$p_1^{\ast}(E_1)\simeq E$.
We obtain the induced Hermitian metric $h_1$ of $E_1$.
There exist the unitary connection $\nabla_1$ 
and the anti-Hermitian endomorphism $\phi$ of $(E_1,h_1)$
such that
$\nabla=p_1^{\ast}(\nabla)+p_1^{\ast}(\phi)\,dx_1$.
Then, the ASD equation (\ref{eq;21.8.4.10}) holds
if and only if $(E,h,\nabla,\phi)$ is a monopole.

We shall introduce the notion of mini-complex structure
of a $3$-dimensional manifold
(see \S\ref{subsection;17.10.28.1} for mini-complex structure).
We also introduce the notion of mini-holomorphic bundle
on a mini-complex manifold,
which is essentially a vector bundle
equipped with two commuting differential operators
as in \S\ref{subsection;21.8.2.1}
and \S\ref{subsection;21.8.2.2}.
(See \S\ref{subsection;17.10.28.42} for mini-holomorphic bundles.)
In our situation,
we obtain the mini-complex manifold $U_1^{\lambda}$
as $U_1$ equipped with the mini-complex structure
induced by the complex structure of $U^{\lambda}$.
(See \S\ref{subsection;17.10.3.1}.)
A $\Gamma$-equivariant holomorphic vector bundle on $U^{\lambda}$
is equivalent to a mini-holomorphic bundle on $U_1^{\lambda}$.
Moreover, the $\Gamma$-equivariant holomorphic bundle
$(E^{\lambda},\delbar_{E^{\lambda}})$
underlying $(E,h,\nabla)$
is equivalent to the mini-holomorphic bundle 
$(E_1^{\lambda},\delbar_{E_1^{\lambda}})$
underlying the monopole $(E_1,h_1,\nabla_1,\phi)$:
\[
\begin{CD}
\left(
 \mbox{
 \begin{tabular}{c}
  $\Gamma$-equivariant\\
  instantons on $U$
 \end{tabular}}
 \right)
 @>>>
\left(
 \mbox{
 \begin{tabular}{c}
  $\Gamma$-equivariant holomorphic\\
  vector bundles on $U^{\lambda}$
 \end{tabular}}
 \right)
 \\
 @V{=}VV @V{=}VV \\
\left(
 \mbox{\begin{tabular}{c}
	monopoles
	on $U_1$
       \end{tabular}}
 \right)
 @>>>
\left(
 \mbox{\begin{tabular}{c}
	mini-holomorphic\\
	bundles on $U_1^{\lambda}$
	\end{tabular}}
\right).
\end{CD}
\]
Very roughly,
when $U_1$ is the complement of a finite subset in
$(\real/T\seisuu)\times\real^2$,
our Kobayashi-Hitchin correspondence for periodic monopoles
are equivalences between
monopoles on $U_1$ and mini-holomorphic bundles
on $U_1^{\lambda}$ satisfying some conditions,
and the latter objects are translated to difference modules.
We shall explain our correspondences for more details
in \S\ref{section;21.8.11.2}--\S\ref{section;21.8.11.3}.

\subsection{Instantons and monopoles as harmonic bundles of infinite rank}
\label{subsection;21.8.4.22}

\subsubsection{Instantons as harmonic bundles of infinite rank}

If $U=\cnum_z\times U_0\subset \cnum_z\times\cnum_w$,
$\cnum_z$-equivariant harmonic bundles
on $U$ are equivalent to harmonic bundles on $U_0$,
as mentioned in \S\ref{subsection;21.8.4.20}.
It also implies that we may regard instantons on $U$
as harmonic bundles of infinite rank on $U_0$,
as we shall explain below.
Let $(E,\nabla,h)$ be an instanton on $U$.
Let $\nbige^{\infty}$ denote the sheaf of $C^{\infty}$-sections of
$E$ on $U$.
Let $\nbigc_{U_0}^{\infty}$ denote the sheaf of
$C^{\infty}$-functions on $U_0$.
\index{sheaf $\nbigc_{U_0}^{\infty}$}
Let $p:U\lrarr U_0$ denote the projection,
and let $p_!$ denote the proper push-forward functor for sheaves
with respect to $p$.
\index{functor \mbox{$p_!$}}
We obtain the $\nbigc_{U_0}^{\infty}$-module
$p_{!}\nbige^{\infty}$ on $U_0$.
We regard
$p_{!}\nbige^{\infty}$
as a $C^{\infty}$-vector bundle on $U_0$
of infinite rank.
The metric $h$ and the integration
induces a Hermitian metric $p_!h$ on $p_!\nbige^{\infty}$.
We may regard 
$\nabla_{\wbar}$ as a holomorphic structure
$\del_{p_!\nbige^{\infty},\wbar}$
of $p_{!}\nbige^{\infty}$.
We also have the differential operators
$\del_{p_!\nbige^{\infty},w}$ on
$p_!\nbige^{\infty}$ induced by $\nabla_{w}$.
Then,
$\del_{p_!\nbige^{\infty},\wbar}\,d\wbar
+\del_{p_!\nbige^{\infty},w}\,dw$
is a connection of $p_!\nbige^{\infty}$
which is unitary with respect to $p_!h$.
The differential operator
$f=\nabla_{\zbar}$ on $p_{!}\nbige^{\infty}$
is $\nbigc^{\infty}_{U_0}$-linear.
Because $[\nabla_{\zbar},\nabla_{\wbar}]=0$,
we may regard $f$ as a holomorphic endomorphism of
$(p_!\nbige^{\infty},\del_{p_!\nbige^{\infty},\wbar})$,
and we obtain the Higgs field $\theta=f\,dw$.
We also have the $\nbigc^{\infty}_{U_0}$-linear endomorphism
$f^{\dagger}=-\nabla_{\zbar}$.
We may regard $f^{\dagger}$ as the adjoint of $f$
with respect to $p_!h$.
We set $\theta^{\dagger}=f^{\dagger}\,d\wbar$.
Let $F$ denote the curvature of the connection
$\del_{p_!\nbige^{\infty},\wbar}\,d\wbar
+\del_{p_!\nbige^{\infty},w}\,dw$.
Then,
the equation
$[\nabla_{\wbar},\nabla_w]
+[\nabla_{\zbar},\nabla_z]=0$
is exactly the equation
$F+[\theta,\theta^{\dagger}]=0$.
In this sense,
$(p_!\nbige^{\infty},\theta,p_!h)$ is a harmonic bundle
of infinite rank on $U_0$.

\subsubsection{The underlying $\lambda$-flat bundles of infinite rank}
Let us pursue this analogy further.
For a complex number $\lambda$,
if $(\xi,\eta)=(z+\lambda\wbar,w-\lambda\zbar)$,
the following holds:
\[
 (1+|\lambda|^2)\del_{\xibar}=\lambda\del_w+\del_{\zbar},\quad
 (1+|\lambda|^2)\del_{\etabar}=\del_{\wbar}-\lambda\del_z.
\]
Hence, for any section $s$ of $p_!\nbige^{\infty}$
and $g\in\nbigc^{\infty}_{U_0}$,
we have the following equalities:
\[
(1+|\lambda|^2)\nabla_{\xibar}(gs)
=g(1+|\lambda|^2)\nabla_{\xibar}(s)
+\lambda\del_w(g)\cdot s,
\]
\[
(1+|\lambda|^2)\nabla_{\etabar}(gs)
=g(1+|\lambda|^2)\nabla_{\etabar}(s)
+\del_{\wbar}(g)\cdot s.
\]
We consider the differential operator
$\nabla^{\lambda}=(1+|\lambda|^2)\nabla_{\xibar}\,dw
+(1+|\lambda|^2)\nabla_{\etabar}\,d\wbar$
on $p_!\nbige^{\infty}$,
which is naturally regarded as a $\lambda$-connection.
The commutativity
of $\nabla_{\xibar}$ and $\nabla_{\etabar}$
is equivalent to
the flatness of the $\lambda$-connection
$\nabla^{\lambda}$,
i.e.,
$\nabla^{\lambda}\circ\nabla^{\lambda}=0$.
Moreover, 
we have
$\nabla^{\lambda}=
\del_{p_!\nbige^{\infty},\wbar}\,d\wbar
+\lambda \theta^{\dagger}
+\lambda\del_{p_!\nbige^{\infty},w}\,dw
+\theta$.
Hence, we may regard
the induced holomorphic bundle
$(E^{\lambda},\delbar_{E^{\lambda}})$ on $U^{\lambda}$
as the $\lambda$-flat bundle of infinite rank on $U_0$
associated with the harmonic bundle
$(p_!\nbige^{\infty},\theta,p_!h)$ of infinite rank.

\subsubsection{Monopoles as harmonic bundles of infinite rank}

Let us consider the case where
$\Gamma=\{(x_1+\sqrt{-1}nT,0)\,|\,x_1\in\real,n\in\seisuu\}
\subset\cnum_z\times\cnum_w$.
Suppose that $U=\cnum_z\times U_0$ for an open subset $U_0\subset\cnum_w$.
Set $U_1=U/\Gamma\subset (\real/T\seisuu)\times\cnum_w$.
Because $(\real/T\seisuu)$-equivariant monopoles on $U_1$
are equivalent to $\cnum_z$-equivariant instantons on $U$,
they are equivalent to harmonic bundles on $U_0$.
It means that we may regard monopoles on $U_1$
as harmonic bundles on $U_0$ of infinite rank.
We shall explain
the more precise relation between monopoles on $U_1$
and harmonic bundles on $U_0$
in \S\ref{subsection;20.7.30.42},
the relation of the mini-holomorphic bundles
on $U_1^{\lambda}$ and $\lambda$-flat bundles on $U_0$
in \S\ref{subsection;21.8.11.21}
(see also \S\ref{subsection;21.8.13.21},
 \S\ref{subsection;17.10.28.20} and \S\ref{subsection;20.7.31.22}),
and the compatibility of the relations in \S\ref{subsection;21.8.12.100}.
We shall give a complement on the compatibility
in \S\ref{subsection;17.10.12.1}.

\section{Difference modules with parabolic structure}
\label{section;21.8.11.2}

As an algebraic objects corresponding to
periodic monopoles,
we introduce the concept of parabolic structure
for difference modules.
Let $\cnum(y)$ denote the field of rational functions
of the variable $y$.
\index{field $\cnum(y)$}
Let $\cnum[y]$ denote the ring of polynomials.
\index{ring $\cnum[y]$}
For $y_0\in \cnum$,
$\cnum[\![y-y_0]\!]$
denote the ring of the formal power series
of $y-y_0$,
i.e.,
$\cnum[\![y-y_0]\!]
=\bigl\{
\sum_{j=0}^{\infty}a_j(y-y_0)^j
\,\,\big|\,a_j\in\cnum
\bigr\}$.
\index{\mbox{ring $\cnum[\![y-y_0]\!]$}}
Let $\cnum(\!(y-y_0)\!)$
denote the fraction field of
$\cnum[\![y-y_0]\!]$,
i.e.,
$\cnum(\!(y-y_0)\!)$
is the field of formal Laurent power series
of $y-y_0$.
\index{\mbox{field $\cnum(\!(y-y_0)\!)$}}
Similarly,
let $\cnum[\![y^{-1}]\!]$
denote the ring of the formal power series
of $y^{-1}$,
and
let $\cnum(\!(y^{-1})\!)$
denote the ring of the formal Laurent power series
of $y^{-1}$.
\index{\mbox{ring $\cnum[\![y^{-1}]\!]$}}
\index{\mbox{field $\cnum(\!(y^{-1})\!)$}}

\subsection{Difference modules}

Take $\varrho\in\cnum$.
Let $\Phi^{\ast}$ be the automorphism of
the field $\cnum(y)$ induced by
$\Phi^{\ast}(y)=y+\varrho$.
We prefer to regard it as the pull back
of functions by the automorphism $\Phi$
of $\cnum$ or $\proj^1$.
\begin{df}
\label{df;20.7.29.1}
\index{difference module}
\index{$\varrho$-difference module}
In this monograph,
a difference module is
a finite dimensional $\cnum(y)$-vector space $\vecV$
with a $\cnum$-linear automorphism $\Phi^{\ast}$
such that 
\[
 \Phi^{\ast}(g\,s)=
 \Phi^{\ast}(g)
 \cdot\Phi^{\ast}(s)
\]
for any $g(y)\in \cnum(y)$
 and $s\in\vecV$.
 When we emphasize $\varrho$,
it is called a $\varrho$-difference module.
 \hfill\qed
\end{df}

\begin{rem}
We set
$\nbiga:=\bigoplus_{n\in\seisuu}\cnum[y](\Phi^{\ast})^n$.
It is naturally a $\cnum$-vector space.
We define the multiplication $\nbiga\times\nbiga\lrarr\nbiga$
by
$\sum a_n(y)(\Phi^{\ast})^n\cdot
 \sum b_m(y)(\Phi^{\ast})^m
=\sum a_n(y)b_m(y+n\varrho)
 (\Phi^{\ast})^{m+n}$.
Thus, $\nbiga$ is an associative $\cnum$-algebra.
As a difference module,
it would be more natural to consider
a finitely generated left $\nbiga$-module $M$
such that
(i) $M$ is a torsion-free $\cnum[y]$-module,
(ii) $\dim_{\cnum(y)}\cnum(y)\otimes_{\cnum[y]}M<\infty$.
The data of such finitely generated $\nbiga$-modules
 are contained in the parabolic structure at finite place,
 which will be explained below.
Indeed,
let $\vecV$ be a difference module in the sense of
Definition {\rm\ref{df;20.7.29.1}}.
For any $\cnum[y]$-free submodule $V\subset\vecV$
such that $\cnum(y)\otimes_{\cnum[y]}V=\vecV$,
 we obtain the $\nbiga$-module
 $\nbiga\cdot V\subset\vecV$
 which satisfies the above finiteness conditions. 
 \hfill\qed
\end{rem}

\begin{rem}
If $\varrho=0$,
the difference module is just
a finite dimensional $\cnum(y)$-vector space
with a $\cnum(y)$-linear automorphism. 
\hfill\qed
\end{rem}

\subsection{Parabolic structure of difference modules
 at finite place}

Let $\vecV$ be a difference module.

\begin{df}[Definition
 \ref{df;17.12.1.10}]
\index{parabolic structure at finite place}
A parabolic structure of $\vecV$ at finite place
is a tuple as follows.
\begin{itemize}
\item 
      A $\cnum[y]$-free submodule $V\subset\vecV$
      such that $V\otimes_{\cnum[y]}\cnum(y)=\vecV$.
\item
 A function $m:\cnum\lrarr \seisuu_{\geq 0}$
     such that $\sum_{x\in\cnum} m(x)<\infty$.
     We assume
\[
 V\otimes_{\cnum[y]}\cnum[y]_D
=(\Phi^{\ast})^{-1}(V)\otimes_{\cnum[y]}\cnum[y]_D,
\]
where we set $D:=\{x\in\cnum\,|\,m(x)>0\}$,
and let $\cnum[y]_D$ denote the localization of
$\cnum[y]$ with respect to $y-x$ $(x\in D)$.
 \item
 A sequence of real numbers
 $0\leq \tau_{x}^{(1)} <\cdots<\tau_{x}^{(m(x))}<1$
 for each $x\in \cnum$.
 If $m(x)=0$, the sequence is assumed to be empty.
 The sequence is denoted by  $\vectau_{x}$.
 \index{tuple $\vectau_x$}
\item
 Lattices 
 $L_{x,i}\subset
 V\otimes_{\cnum[y]}\cnum(\!(y-x)\!)$
 for $x\in\cnum$ and $i=1,\ldots,m(x)-1$.
 We formally set
 $L_{x,0}:=
 V\otimes_{\cnum[y]}\cnum[\![y-x]\!]$
 and 
 $L_{x,m(x)}:=
 (\Phi^{\ast})^{-1}V\otimes\cnum[\![y-x]\!]$.
 The tuple of lattices is denoted by
 $\vecL_{x}$. \index{tuple $\vecL_x$}
Note that
$V\otimes\cnum[\![y-x]\!]
=(\Phi^{\ast})^{-1}(V)\otimes\cnum[\![y-x]\!]$
if $m(x)=0$.
\hfill\qed
\end{itemize}
\end{df}

\begin{rem}
We may regard a parabolic structure of finite place of
a difference module
as a reincarnation of a part of
an ordinary parabolic structure of
$\lambda$-flat bundles.
See \S{\rm\ref{subsection;21.8.6.22}}.
\hfill\qed
\end{rem}

\subsection{Good parabolic structure at $\infty$}
\label{subsection;18.2.1.1}

Note that
the automorphism $\Phi^{\ast}$ of $\cnum(y)$
induces an automorphism of $\cnum(\!(y^{-1})\!)$,
which is also denoted by $\Phi^{\ast}$.
We set
$\vecVhat:=
\vecV\otimes_{\cnum(y)}
 \cnum(\!(y^{-1})\!)$,
 which is a finite dimensional
 $\cnum(\!(y^{-1})\!)$-vector space
with the induced $\cnum$-linear automorphism $\Phi^{\ast}$
such that
$\Phi^{\ast}(fs)=\Phi^{\ast}(f)\Phi^{\ast}(s)$
for $f\in\cnum(\!(y^{-1})\!)$
and $s\in\vecVhat$.
Such $(\vecVhat,\Phi^{\ast})$ is called
a formal difference module.
\index{formal difference module}

For any $p\in\seisuu_{\geq 1}$,
we set
$S(p):=\bigl\{
 \sum_{j=1}^{p-1}
 \gminia_jy^{-j/p}\,\big|\,
 \gminia_j\in\cnum
 \bigr\}$.
\index{set $S(p)$}
According to the classification of formal difference modules 
\cite{Chen-Fahim, Praagman, Turrittin},
there exist $p\in\seisuu_{\geq 1}$
and a decomposition of the formal difference module
\begin{equation}
\label{eq;17.12.17.1}
 \vecVhat\otimes_{\cnum(\!(y^{-1})\!)}\cnum(\!(y^{-1/p})\!)
=\bigoplus_{\omega\in p^{-1}\seisuu}
 \bigoplus_{\alpha\in\cnum^{\ast}}
 \bigoplus_{\gminia\in S(p)}
 \bigl(
 \vecVhat_{p,\omega,\alpha,\gminia},
 \Phi^{\ast}
 \bigr),
\end{equation}
such that each $\vecVhat_{p,\omega,\alpha,\gminia}$
has a $\cnum[\![y^{-1/p}]\!]$-lattice
$L_{p,\omega,\alpha,\gminia}$
satisfying
\[
 \bigl(
 \alpha^{-1}y^{\omega}
 \Phi^{\ast}
-(1+\gminia)\id
 \bigr)
L_{p,\omega,\alpha,\gminia}
\subset
 y^{-1}
 L_{p,\omega,\alpha,\gminia}.
\]
We note that in the case $\varrho=0$
we can obtain the decomposition (\ref{eq;17.12.17.1})
from the generalized eigen decomposition of a linear map.

\begin{df}[Definition 
 \ref{df;20.7.28.10}, Definition \ref{df;20.7.28.11}]
 \label{df;20.7.30.1}
 \index{good parabolic structure at infinity}
 A good parabolic structure of $\vecV$ at $\infty$
is a filtered bundle $\nbigp_{\ast}\vecVhat$ over $\vecVhat$
 which is good with respect to $\Phi^{\ast}$
 in the following sense.
\begin{itemize}
\item
There exists the decomposition
 $\varphi_p^{\ast}(\nbigp_{\ast}\vecVhat)
 =\bigoplus \nbigp_{\ast}\vecVhat_{p,\omega,\alpha,\gminia}$,
 where
 $\varphi_p^{\ast}(\nbigp_{\ast}\vecVhat)$
 is the filtered bundle
 over $\vecVhat\otimes\cnum(\!(y^{-1/p})\!)$
 obtained as the pull back of
 $\nbigp_{\ast}\vecVhat$.
     (See {\rm\S\ref{subsection;17.10.27.11}} or
     {\rm \S\ref{subsection;20.7.20.20}}
     for the pull back of filtered bundles
     via a ramified covering.)
\item
$\bigl(
 \alpha^{-1}y^{\omega}
 \Phi^{\ast}
-(1+\gminia)\id
 \bigr)
 \nbigp_a\vecVhat_{p,\omega,\alpha,\gminia}
\subset
 y^{-1}
 \nbigp_a\vecVhat_{p,\omega,\alpha,\gminia}$
     for any $a\in\real$
     and for any $(\omega,\alpha,\gminia)$.
     \hfill\qed
\end{itemize}
\end{df}

\subsection{Parabolic difference modules}

\begin{df}
\index{parabolic difference module}
\label{df;21.9.17.1}
 A parabolic difference module
means a difference module
equipped with a parabolic structure at finite place
and a good parabolic structure at infinity. 
\hfill\qed
 \end{df}

Note that in the decomposition (\ref{eq;17.12.17.1}),
the numbers
\begin{equation}
\label{eq;17.12.18.1}
 r(\omega):=
 \sum_{\alpha,\gminia}
 \rank \vecVhat_{p,\omega,\alpha,\gminia}
\end{equation}
are independent of $p$,
and well defined for $\omega\in\rnum$.
Indeed, there exists the slope decomposition
$\vecVhat=\bigoplus_{\omega\in\rnum}\nbigs_{\omega}\vecVhat$
as explained in \S\ref{subsection;20.7.20.10}
for which we obtain
$\nbigs_{\omega}(\vecVhat)\otimes_{\cnum(\!(y^{-1})\!)}
 \cnum(\!(y_p^{-1})\!)
 =\bigoplus_{\alpha,\gminia}\vecVhat_{p,\omega,\alpha,\gminia}$
for any $p$ as above.
Hence, $r(\omega)=\dim_{\cnum}\nbigs_{\omega}\vecVhat$.
For $\nbigp_{\ast}\vecVhat$ as in Definition \ref{df;20.7.30.1},
we obtain the decomposition
$\nbigp_{\ast}\vecVhat
=\bigoplus\nbigp_{\ast}\nbigs_{\omega}(\vecVhat)$.
\index{slope decomposition}
\index{number $r(\omega)$}

\subsection{Degree and stability condition}

Consider a parabolic difference module
\[
\vecV_{\ast}=(\vecV,(V,m,(\vectau_{x},\vecL_x)_{x\in\cnum}),
 \nbigp_{\ast}\vecVhat).
\]
Let $\nbigf_V$ be the $\nbigo_{\proj^1}(\ast\infty)$-module
associated with $V$.
We obtain the filtered bundle
$\nbigp_{\ast}\nbigf_V$ 
over $\nbigf_V$
induced by $\nbigp_{\ast}\vecVhat$.
If $m(x)>0$,
for each $i=1,\ldots,m(x)$,
we set 
\index{relative degree $\deg(L_{i+1,x},L_{i,x})$}
\[
 \deg(L_{i,x},L_{i-1,x}):=
 \dim_{\cnum}\bigl(L_{i,x}/(L_{i,x}\cap L_{i-1,x})\bigr)
-\dim_{\cnum}\bigl(L_{i-1,x}/(L_{i,x}\cap L_{i-1,x})\bigr).
\]
Then, we define
\index{degree $\deg(\vecV_{\ast})$}
\begin{equation}
\label{eq;17.12.4.4}
 \deg(\vecV_{\ast}):=
 \deg(\nbigp_{\ast}\nbigf_V)
+\sum_{x\in\cnum}
 \sum_{i=1}^{m(x)}
 \bigl(1-\tau_x^{(i)}\bigr)
 \deg\bigl(L_{i,x},L_{i-1,x}\bigr)
+\sum_{\omega\in\rnum}
 \frac{\omega}{2} r(\omega).
\end{equation}
Here, $r(\omega)$ are defined in 
(\ref{eq;17.12.18.1}).
We also define
\begin{equation}
\label{eq;21.9.17.2}
\mu(\vecV_{\ast}):=
 \deg(\vecV_{\ast})\big/\dim_{\cnum(y)}\vecV.
\end{equation}
\index{slope $\mu(\vecV_{\ast})$}

Let $\vecV'$ be any difference submodule of $\vecV$,
i.e.,
$\vecV'$ is a $\cnum(y)$-subspace of $\vecV$
such that $\Phi^{\ast}(\vecV')=\vecV'$.
By setting $V':=\vecV'\cap V$
and $L_{i,x}':=L_{i,x}\cap (V'\otimes_{\cnum(y)}\cnum(\!(y-x)\!))$,
we obtain 
a parabolic structure 
$\bigl(V',m,(\vectau_{x},\vecL_x')_{x\in\cnum}\bigr)$
at finite place of $V'$.
By setting $\vecVhat':=\cnum(\!(y^{-1})\!)\otimes_{\cnum(y)}\vecV'$
and $\nbigp_{a}\vecVhat':=\vecVhat'\cap\nbigp_a\vecVhat$,
we obtain a good parabolic structure 
$\nbigp_{\ast}\vecVhat'$ at $\infty$.
Thus, we obtain the following induced parabolic difference module
from $\vecV'$:
\[
\vecV'_{\ast}=
\Bigl(\vecV',
\bigl(V',m,(\vectau_{x},\vecL_x')_{x\in\cnum}\bigr),
\nbigp_{\ast}\vecVhat'
\Bigr).
\]

\begin{df}
\label{df;21.9.17.3}
\index{stable (difference module)}
\index{semistable (difference module)}
\index{polystable (difference module)}
A parabolic difference module
$\vecV_{\ast}$
is called stable (resp. semistable) if 
$\mu(\vecV'_{\ast})
<\mu(\vecV_{\ast})$
(resp.
$\mu(\vecV'_{\ast})
\leq\mu(\vecV_{\ast})$)
for any 
difference submodule $\vecV'$ of $\vecV$
 with 
$0<\dim_{\cnum(y)}\vecV'
 <\dim_{\cnum(y)}\vecV$.
It is called polystable
if it is semistable
and a direct sum of stable ones.
\hfill\qed
\end{df}

\begin{rem}
 There exists the relation
\begin{equation}
\label{eq;20.7.30.2}
 \sum_{x\in\cnum}
 \sum_{i=1}^{m(x)}\deg(L_{i,x},L_{i-1,x})
 +\sum_{\omega\in\rnum} \omega\cdot r(\omega)=0.
\end{equation}
 Indeed,
 by the comparison of
 $V$ and $(\Phi^{\ast})^{-1}(V)$ in $\vecV$,
 we obtain
\[
 \deg(\nbigp_{\ast}(\nbigf_{(\Phi^{\ast})^{-1}(V)}))
 -\deg(\nbigp_{\ast}\nbigf_V)=\sum_{x,i}\deg(L_{i,x},L_{i-1,x}).
\]
There exists the isomorphism of
 $\nbigo_{\proj^1}(\ast\infty)$-modules
 $\Phi^{\ast}(\nbigf_{(\Phi^{\ast})^{-1}(V)})
 \simeq
 \nbigf_{V}$
 induced by $\Phi^{\ast}$.
 For the slope decomposition
 $\vecVhat=\bigoplus\nbigs_{\omega}(\vecVhat)$,
it induces
 $\Phi^{\ast}(\nbigp_a\nbigs_{\omega}(\vecVhat))
 \simeq
  \nbigp_{a-\omega}\nbigs_{\omega}(\vecVhat)$.
 Hence, we obtain
\[
 \deg(\nbigp_{\ast}\nbigf_V)
 =\deg(\nbigp_{\ast}\nbigf_{(\Phi^{\ast})^{-1}(V)})
 +\sum_{\omega\in\rnum} r(\omega)\cdot\omega.
\]
Thus, we obtain {\rm(\ref{eq;20.7.30.2})}. 
\hfill\qed
\end{rem}

\subsection{Easy examples of 
stable parabolic difference modules (1)}
\label{subsection;18.1.15.3}

Let us mention some easy examples of
parabolic difference modules.
Needless to say,
we can easily construct many more interesting examples.
For simplicity, we consider the case $\varrho=0$.

We take a non-empty finite subset $S\subset\cnum$
and a function $\ell:S\lrarr\seisuu$
such that
$\sum_{a\in S}|\ell(a)|\neq 0$.
We also assume that one of
$\ell(a)$ is an odd integer.
We set
$P(y):=\prod_{a\in S}(y-a)^{\ell(a)}\in\cnum(y)$.
We set
$\vecV:=\cnum(y)e_1\oplus\cnum(y)e_2$
and $V:=\cnum[y]e_1\oplus\cnum[y]e_2$.
We define the $\cnum(y)$-automorphism $\Phi^{\ast}$ of $\vecV$
by
\[
 \Phi^{\ast}(e_1,e_2)
=(e_1,e_2)
 \cdot
 \left(
 \begin{array}{cc}
 0 & P(y)\\
 1 & 0
 \end{array}
 \right).
\]
Let $m:\cnum\lrarr \seisuu_{\geq 0}$
be given by
$m(a)=1$ $(a\in S)$
and $m(a)=0$ $(a\not\in S)$.
We take $0\leq\tau_{a}<1$ for $a\in S$.
Thus, we obtain a parabolic structure
of $\vecV$ at finite place.
Note the following equality for any $a\in S$:
\[
 \deg(L_{a,1},L_{a,0})=-\ell(a).
\]

As for good parabolic structures at infinity,
there are two cases
where $\ell(\infty):=-\sum_{a\in S}\ell(a)$ is even or odd.

\subsubsection{The case where $\ell(\infty)$ is even}

Suppose $\ell(\infty)=-2k$ for $k\in\seisuu$.
There exists $\tau\in\cnum(\!(y^{-1})\!)$
such that $\tau^{2k}=P(y)$.
Note that $\tau/y\in \cnum[\![y^{-1}]\!]$
and that it is invertible.
We set $v_1:=\tau^ke_1+e_2$
and $v_2:=\tau^ke_1-e_2$.
Because
$\Phi^{\ast}(v_1)=\tau^kv_1$
and $\Phi^{\ast}(v_2)=-\tau^kv_2$,
we obtain the decomposition
$\vecVhat=
 \cnum(\!(y^{-1})\!)v_1
\oplus 
 \cnum(\!(y^{-1})\!)v_2$
compatible with the action of $\Phi^{\ast}$.
A filtered bundle $\nbigp_{\ast}\vecVhat$ over $\vecVhat$
is good with respect to $\Phi^{\ast}$
if and only if 
it satisfies
$\nbigp_{\ast}\vecVhat=
 \nbigp_{\ast}\bigl(
 \cnum(\!(y^{-1})\!)v_1
 \bigr)
\oplus
 \nbigp_{\ast}\bigl(
 \cnum(\!(y^{-1})\!)v_2
 \bigr)$.
Hence, such good $\nbigp_{\ast}\vecVhat$
is determined by the numbers
\[
 d_i= \deg^{\nbigp}(v_i):=
\min\{c\in\real\,|\,v_i\in\nbigp_{c}\vecVhat\}.
\]
We can easily see that
$\vecV$ has no non-trivial
difference submodule,
under the assumption that one of $\ell(a)$ is odd.
Hence, 
$\vecV_{\ast}=
(\vecV,(V,m,\vectau_{x},\vecL_x)_{x\in\cnum},\nbigp_{\ast}\vecVhat)$
is a stable parabolic difference module.
It is easy to see
$\deg(\nbigp_{\ast}\nbigf_{V})=k-d_1-d_2$,
and hence 
\begin{equation}
\label{eq;18.1.15.1}
\deg(\vecV_{\ast})
=-d_1-d_2-\sum_{a\in S}(1-\tau_{a})\ell(a).
\end{equation}

\subsubsection{The case where $\ell(\infty)$ is odd}
Let us consider the case where $\ell(\infty)$ is odd.
There exists $\tau\in\cnum(\!(y^{-1/2})\!)$
such that $\tau^{-2\ell(\infty)}=P(y)$
and $\tau/y^{1/2}\in\cnum[\![y^{-1/2}]\!]$.
We set
$v_1:=\tau^{-\ell(\infty)}e_1+e_2$
and $v_2:=\tau^{-\ell(\infty)}e_1-e_2$.
Because
$\Phi^{\ast}(v_1)=\tau^{-\ell(\infty)}v_1$
and 
$\Phi^{\ast}(v_2)=-\tau^{-\ell(\infty)}v_2$,
we obtain the decomposition
$\vecVhat\otimes\cnum(\!(y^{-1/2})\!)
=\cnum(\!(y^{-1/2})\!)v_1
\oplus
 \cnum(\!(y^{-1/2})\!)v_2$
 compatible with the action of $\Phi^{\ast}$.
A filtered bundle $\nbigp_{\ast}\vecVhat$ over $\vecVhat$
is good with respect to $\Phi^{\ast}$
 if and only if
 the induced filtered bundle
 over $\vecVhat\otimes\cnum(\!(y^{-1/2})\!)$
 satisfies
$\nbigp_{\ast}\bigl(
 \vecVhat\otimes\cnum(\!(y^{-1/2})\!)
 \bigr)
=\nbigp_{\ast}\bigl(
 \cnum(\!(y^{-1/2})\!)v_1
 \bigr)
\oplus
 \nbigp_{\ast}\bigl(
 \cnum(\!(y^{-1/2})\!)v_2
 \bigr)$.
Hence, it is determined by
the numbers
$d_i:=\min\bigl\{
 c\in\real\,\big|\,
 v_i\in\nbigp_{\ast}(\vecVhat\otimes\cnum(\!(y^{-1/2})\!))
 \bigr\}$.
As in the previous case,
$\vecV_{\ast}=
(\vecV,(V,m,(\vectau_{x},\vecL_x)),\nbigp_{\ast}\vecVhat)$
is stable.

Because $\nbigp_{\ast}(\vecVhat\otimes\cnum(\!(y^{-1/2})\!))$
is preserved by the natural action of the Galois group of
$\cnum(\!(y^{-1/2})\!)/\cnum(\!(y^{-1})\!)$,
we obtain $d_1=d_2=:d$.
It is easy to see
$\deg(\nbigp_{\ast}\nbigf_V)
=-d-\ell(\infty)/2$,
and hence
\begin{equation}
\label{eq;18.1.15.2}
\deg(\vecV_{\ast})
=-d-\sum_{a\in S}(1-\tau_{a})\ell(a).
\end{equation}

\subsection{Easy examples of 
stable parabolic difference modules (2)}
\label{subsection;18.1.31.10}

Take polynomials $P(y)$ and $Q(y)$
such that the following holds.
\begin{itemize}
\item
$2\deg(Q)\geq\deg(P)$.
If $2\deg(Q)=\deg(P)$,
we assume $q_{\deg(Q)}^2-4p_{\deg(P)}\neq 0$
for the expressions
$Q=\sum_{j=0}^{\deg(Q)}q_jy^j$
and $P=\sum_{j=0}^{\deg(P)}p_jy^j$.
We assume that
$(\deg(P),\deg(Q))\neq (0,0)$.
 \item
$P$ has simple zeroes.      
\end{itemize}

We set 
$\vecV:=
 \cnum(y)e_1
\oplus
 \cnum(y)e_2$,
 and we define
the $\cnum(y)$-automorphism $\Phi^{\ast}$
by
\[
 \Phi^{\ast}(e_1,e_2)
=(e_1,e_2)
 \left(
\begin{array}{cc}
 0 & P \\
 -1 & Q
\end{array}
 \right).
\]
Because $P$ has simple zeroes,
it is easy to see that
there is no difference submodule
$\vecV'\subset\vecV$
such that $\dim_{\cnum(y)}\vecV'=1$.

We put
$\alpha_1=
 2^{-1}\bigl(Q+Q(1-4PQ^{-2})^{1/2}\bigr)
 \in\cnum(\!(y^{-1})\!)$
 and $\alpha_2=\alpha_1^{-1}P$.
If $2\deg(Q)>\deg(P)$,
we choose $(1-4PQ^{-1})^{1/2}\in\cnum[\![y^{-1}]\!]$
such that
$(1-4PQ^{-1})^{1/2}-1\in y^{-1}\cnum[\![y^{-1}]\!]$.
Note that $\alpha_1\neq\alpha_2$.
Then, $\alpha_i$ $(i=1,2)$ are the roots of
the characteristic polynomial $T^2-QT+P$
of the automorphism.
Note that
$y^{-\deg Q}\alpha_1$
and 
$y^{-\deg P+\deg Q}\alpha_2$
are invertible in 
$\cnum[\![y^{-1}]\!]$.

We set $V:=\cnum[y]e_1\oplus\cnum[y]e_2$.
Set $S:=\{a\in\cnum\,|\,P(a)=0\}$.
Let $m:\cnum\lrarr\seisuu_{\geq 0}$
be the function determined by
$m(a)=1$ $(a\in S)$ and
$m(a)=0$ $(a\not\in S)$.
We take $0\leq \tau_{a}<1$ $(a\in S)$.
Then, the tuple $(V,m,(\tau_{a})_{a\in S})$ determines
a parabolic structure at finite place of $\vecV$.

Set $v_i:=Pe_1+\alpha_i e_2$,
for which $\Phi^{\ast}(v_i)=\alpha_iv_i$ holds.
Then, a good filtered bundle $\nbigp_{\ast}\vecVhat$
over $\vecVhat$ with respect to $\Phi^{\ast}$
is determined by $d_i:=\deg^{\nbigp}(v_i)$.
We obtain
$\deg(\nbigp_{\ast}\nbigf_V):=
 \deg(P)+\deg(Q)-d_1-d_2$.
The degree of 
$\vecV_{\ast}=
(\vecV,(V,m,\{\tau_{a}\}_{a\in S},\nbigp_{\ast}\vecVhat))$
is as follows:
\begin{multline}
\label{eq;18.1.31.2}
\deg(\vecV_{\ast})=
 \deg(P)+\deg(Q)-d_1-d_2
-\sum_{a\in S}(1-\tau_{a})\ell(a)
-\frac{\deg(Q)}{2}-\frac{\deg(P)-\deg(Q)}{2}
 \\
=\frac{1}{2}\deg(P)+\deg(Q)
-d_1-d_2-\sum_{a\in S}(1-\tau_{a}).
\end{multline}
For any given $\{\tau_{a}\}_{a\in S}$,
we may choose $(d_1,d_2)\in\real^2$
such that (\ref{eq;18.1.31.2}) vanishes.

Note that the numbers
$n:=\deg(P)$
and $k:=2^{-1}\bigl(\deg(Q)-(\deg(P)-\deg(Q))+n\bigr)=\deg(Q)$
correspond to 
the number of singularity
and the non-abelian charge in \cite[Page 5]{Cherkis-Kapustin2}
in the context of periodic monopoles.

\section{Kobayashi-Hitchin correspondences
for periodic monopoles}
\label{section;21.8.11.3}

Kobayashi-Hitchin correspondences
between periodic monopoles and parabolic difference modules
are stated as follows.
\begin{thm}[Corollary
\ref{cor;21.9.17.60}]
 \label{thm;17.12.4.10}
For any $T>0$ and $\lambda\in\cnum$,
there exists a natural bijective correspondence
between the isomorphism classes of the following objects.
\begin{itemize}
\item
 Periodic monopoles of GCK type
 on $S^1_T\times\cnum_w$.
\item
     Polystable parabolic $(2\sqrt{-1}\lambda T)$-difference modules
     of degree $0$.
\end{itemize}
 \end{thm}

We shall explain an outline
to obtain a parabolic $2\sqrt{-1}\lambda T$-difference module
from a periodic monopole
for $\lambda=0$ in \S\ref{section;21.8.3.30},
and for general $\lambda$ in \S\ref{section;21.8.3.31}.
The constructions are given by way of mini-holomorphic bundles
on a $3$-dimensional manifolds,
inspired by the work of
Charbonneau-Hurtubise \cite{Charbonneau-Hurtubise},
which we briefly recalled in \S\ref{subsection;21.8.2.2}.

\begin{cor}
In {\rm\S\ref{subsection;18.1.15.3}},
if we choose $(d_1,d_2)\in\real^2$ (resp. $d\in\real$)
such that {\rm(\ref{eq;18.1.15.1})}  
(resp. {\rm (\ref{eq;18.1.15.2})})
vanishes,
we obtain the periodic monopoles of GCK-type
 on $S^1_T\times\cnum$ for any $T>0$,
corresponding to the stable parabolic difference modules.
Similarly,
in {\rm\S\ref{subsection;18.1.31.10}},
if we choose $(d_1,d_2)\in\real^2$
such that {\rm(\ref{eq;18.1.31.2})} vanishes,
 we obtain the periodic monopoles of GCK-type
 on $S^1_T\times\cnum$ for any $T>0$,
corresponding to the stable parabolic difference modules.
\hfill\qed
\end{cor}
Note that even the existence of periodic monopoles was non-trivial,
which was studied by Foscolo in \cite{Foscolo-construction}.

\begin{rem}
A special case of this type of correspondences
is independently studied by Elliott and Pestun
in {\rm\cite{Elliott-Pestun}},
where rich studies of the related subjects
are described.
\hfill\qed
\end{rem}

\begin{rem}
In {\rm\cite{Harland}},
Harland studied a classification of periodic monopoles of rank $2$
without any singularity,
by using the Nahm transform
between harmonic bundles on $\cnum^{\ast}$ and periodic monopoles.
It is given in terms of spectral curves with a parabolic line bundle,
which corresponds to the classification at $\lambda=0$
in Theorem {\rm\ref{thm;17.12.4.10}}.
\hfill\qed
\end{rem}

\begin{rem}
In {\rm\cite{Mochizuki-q-difference}},
we studied similar correspondences between
$q$-difference modules and 
 doubly periodic monopoles,
 i.e.,
 monopoles on the product of
$\real$ and a $2$-dimensional torus.
In {\rm\cite{Mochizuki-triply-periodic-monopoles}},
we studied similar correspondences between
difference modules on elliptic curves
 and
 triply periodic monopoles,
 i.e.,
monopoles on a $3$-dimensional torus.
The triply periodic case is much easier
because we do not have to study
the asymptotic behaviour of monopoles
around infinity.
\hfill\qed
\end{rem}

\begin{rem}
Kontsevich and Soibelman proposed
a non-abelian Hodge theory
in the context of their holomorphic Floer theory
{\rm\cite{Kontsevich-Soibelman}}.
Correspondences for monopoles may be
regarded as $1$-dimensional examples of their theory.
\hfill\qed
\end{rem}

\subsection{The correspondence in the case $\lambda=0$}
\label{section;21.8.3.30}

We first explain our construction in the case $\lambda=0$,
which is conceptually and technically
easier than the case $\lambda\neq 0$.

\subsubsection{Mini-complex structure}

Let $\kappa$ denote the $\seisuu$-action on $\real_t\times\cnum_w$
defined by
$\kappa_n(t,w)=(t+nT,w)$ $(n\in\seisuu)$.
Let $\nbigm$ denote the quotient space of
$\real_t\times\cnum_w$
by the $\seisuu$-action $\kappa$,
i.e.,
$\nbigm=S^1_T\times\cnum_w$.
\index{space $\nbigm$}
\index{action $\kappa$}
We note that the complex vector fields
$\del_t$ and $\del_{\wbar}$
are naturally defined on $\nbigm$.

We have the naturally defined mini-complex structure
on $\real_t\times\cnum_w$.
(See {\rm\S\ref{subsection;17.10.28.1}}
for mini-complex structure.)
When we emphasize to consider this mini-complex structure,
we use the notation $\nbigm^0$,
i.e.,
let $\nbigm^0$ denote the $3$-dimensional manifold $\nbigm$
equipped with the mini-complex structure
induced by $(t,w)$.
\index{mini-complex manifold $\nbigm^{0}$}
A $C^{\infty}$-function $f$ on an open subset 
$U\subset\nbigm^0$ is called mini-holomorphic
if $\del_{t}f=\del_{\wbar}f=0$.
Let $\nbigo_{\nbigm^{0}}$
denote the sheaf of mini-holomorphic functions 
on $\nbigm^0$.
\index{sheaf $\nbigo_{\nbigm^0}$}
\index{mini-complex coordinate system $(t,w)$}

\subsubsection{Mini-holomorphic bundles associated with monopoles}

Let $Z$ be a finite subset of $\nbigm$.
Let $(E,h,\nabla,\phi)$ be a monopole on $\nbigm\setminus Z$.
Recall that the Bogomolny equation implies that
\[
 [\nabla_{t}-\sqrt{-1}\phi,\nabla_{\wbar}]=0.
\]
Hence, by considering the sheaf of local sections $s$ of $E$
such that 
$(\nabla_{t}-\sqrt{-1}\phi)s=
 \nabla_{\wbar}s=0$,
we obtain a locally free
$\nbigo_{\nbigm^{0}\setminus Z}$-module
$\nbige^{0}$.
\index{$\nbigo_{\nbigm^{0}\setminus Z}$-module $\nbige^{0}$}

\subsubsection{Dirac type singularity}
\label{subsection;21.8.3.40}

Let $P\in Z$.
We take a lift 
$\Ptilde=(t^0,w^0)
\in\real_{t}\times\cnum_{w}$
of $P$.
Let $\epsilon>0$ and $\delta>0$ be small,
and consider
$U_{-\epsilon}:=
\bigl\{(t^0-\epsilon,w)\,\big|\,|w|<\delta\bigr\}$
and
$U_{\epsilon}:=
\bigl\{(t^0+\epsilon,w)\,\big|\,|w|<\delta
\bigr\}$.
By taking the restriction of $\nbige^0$
to $U_{\pm\epsilon}$,
we obtain locally free $\nbigo_{U_{\pm\epsilon}}$-modules
$\nbige^{0}_{|U_{\pm\epsilon}}$.
Set 
$U^{\ast}_{\pm\epsilon}:=
U_{\pm\epsilon}\setminus\{(t^0\pm\epsilon,w^0)\}$.
As in the previous studies recalled
in \S\ref{subsection;21.8.2.1} and \S\ref{subsection;21.8.2.2},
by the parallel transport with respect to 
$\nabla_{t}-\sqrt{-1}\phi$,
we obtain the isomorphism 
$\Pi_P:\nbige^{0}_{|U^{\ast}_{-\epsilon}}
\simeq
\nbige^{0}_{|U^{\ast}_{\epsilon}}$,
called the scattering map.
\index{scattering map}
If $P\in Z$ is Dirac type singularity of
the monopole $(E,h,\nabla,\phi)$,
it is easy to see that
$\nbige^0$ is of Dirac type at $P$
in the sense that
$\Pi_P$ is meromorphic.
(See \S\ref{subsection;21.8.3.3} for
Dirac type singularity of monopoles.)

\subsubsection{Meromorphic extension and filtered extension at infinity}
\label{subsection;21.8.4.4}

We have the partial compactification
$\real_t\times\proj^1_w$ of $\real_t\times\cnum_w$.
The $\seisuu$-action $\kappa$ on $\real_t\times\cnum_w$
naturally extends to a $\seisuu$-action on $\real_t\times\proj^1_w$,
which is also denoted by $\kappa$.
The quotient space $\nbigmbar^0$ is naturally 
equipped with the mini-complex structure,
and it is a compactification of $\nbigm^0$.
\index{mini-complex manifold $\nbigmbar^0$}
Let $H^{0}_{\infty}:=
 \nbigmbar^{0}\setminus \nbigm^{0}\simeq S^1_T$.
Let $\nbigo_{\nbigmbar^{0}}(\ast H^{0}_{\infty})$
denote the sheaf of meromorphic functions on $\nbigmbar^{0}$
whose poles are contained in $H^{0}_{\infty}$.
\index{sheaf $\nbigo_{\nbigmbar^{0}}(\ast H_{\infty}^{0})$}
\index{space $H^{0}_{\infty}$}
(See \S\ref{subsection;21.8.3.10} for meromorphic functions
in the context of mini-complex structures.)

\begin{rem}
\label{rem;21.9.16.1}
In the study of tame or wild harmonic bundles
we prefer to consider meromorphic or filtered Higgs bundles
on $(C,D)$ rather than Higgs bundles on $C\setminus D$.
(See Remark {\rm\ref{rem;21.9.12.2}}.)
Similarly, we would like to consider meromorphic objects
or filtered objects on $(\nbigmbar^{0},H^0_{\infty})$
to keep the information of the growth orders with respect to $h$,
rather than the transcendental object on $\nbigm^0$.
It is the reason why we consider the compactification
$\nbigmbar^0$.
\hfill\qed
\end{rem}

Let $U$ be any open subset in $\nbigmbar^{0}\setminus Z$.
Let $\nbigp^h\nbige^0(U)$
denote the space of sections of $\nbige^0$
on $U\setminus H^{0}_{\infty}$
satisfying the following condition.
\begin{itemize}
\item For any $P\in U\cap H^{0}_{\infty}$,
there exists a neighbourhood $U_P$ of $P$ in $U$
such that 
$|s_{|U_P\setminus H^{0}_{\infty}}|_h=O\bigl(|w|^{N}\bigr)$
for some $N$.
\end{itemize}
Thus, we obtain the 
$\nbigo_{\nbigmbar^{0}\setminus Z}
 (\ast H^{0}_{\infty})$-module
 $\nbigp^h\nbige^{0}$.
\index{sheaf $\nbigp^h\nbige^{0}$}
We shall prove the following.
\begin{prop}[Proposition 
\ref{prop;17.9.17.50}]
\label{prop;21.8.5.11}
$\nbigp^h\nbige^0$ is a locally free
$\nbigo_{\nbigmbar^{0}\setminus Z}
 (\ast H^0_{\infty})$-module.
\end{prop}

Let $\pi^0:\nbigmbar^{0}\lrarr S^1_T$
denote the map induced by $(t,w)\longmapsto t$.
\index{projection $\pi^0$}
We set $\nbigmbar^0\langle t\rangle=(\pi^0)^{-1}(t)$
and
$\nbigm^0\langle t\rangle=\nbigm^0\cap\nbigmbar^0\langle t\rangle$.
\index{space $\nbigmbar^0\langle t\rangle$}
\index{space $\nbigm^0\langle t\rangle$}
We also set $Z_{t}:=(\pi^0)^{-1}(t)\cap Z$.
By taking the restriction,
we obtain the locally free 
$\nbigo_{\cnum_w\setminus Z_{t}}$-module
$\nbige^0_{|\nbigm^0\langle t\rangle\setminus Z}$.
As in \S\ref{subsection;17.10.27.10} and \S\ref{subsection;17.10.27.11},
we obtain an increasing sequence 
$\nbigp^h_a\bigl(
 \nbige^0_{|\nbigm^0\langle 0\rangle\setminus Z}
 \bigr)$
$(a\in\real)$
of $\nbigo_{\proj^1\setminus Z_{t}}$-modules
by considering local sections $s$ satisfying
$|s|_h=O(|w|^{a+\epsilon})$ for any $\epsilon>0$.
\index{filtration \mbox{$\nbigp^h_a\bigl(
 \nbige^0_{|\nbigm^0\langle 0\rangle\setminus Z}
 \bigr)$}}

\begin{thm}[Theorem \ref{thm;17.10.5.130}]
\label{thm;21.8.5.12}
 The tuple
$\bigl(
 \nbigp^h_{\ast}\bigl(\nbige^0_{|\nbigm^0\langle t\rangle\setminus Z}\bigr)
 \,\big|\,
 t\in S^1_T
\bigr)$
is a good filtered bundle
over $\nbigp^h\nbige^0$
in the sense of 
Definition {\rm\ref{df;20.7.28.20}}
and Definition {\rm\ref{df;17.10.28.10}}.
(See Theorem {\rm\ref{thm;17.10.7.10}}
for a more detailed description.)
\index{good filtered bundle}
\end{thm}

As we shall see in Proposition \ref{prop;17.10.14.21},
the tuple
$\bigl(
 \nbigp^h_{\ast}\bigl(
  \nbige^0_{|\nbigm^0\langle t\rangle\setminus Z}\bigr)
 \,\big|\,
 t\in S^1_T
 \bigr)$
 determines the behaviour of the metric $h$
around $H^{0}_{\infty}$ up to boundedness.
We shall also show that
the compatibility with such filtrations
implies the GCK condition around infinity.
(See \S\ref{subsection;17.10.13.500}.)

\begin{rem}
Although we shall often use the abbreviation 
$\nbigp^h_{\ast}\nbige^0$
to denote the tuple
$\bigl(
 \nbigp^h_{\ast}\bigl(\nbige^0_{|\nbigm^0\langle t\rangle\setminus Z}
 \bigr)
 \,\big|\,
 t\in S^1_T
\bigr)$,
the filtrations depend on $t\in S^1_T$
as basic examples show.
(See {\rm\S\ref{subsection;17.9.26.1}}.)
\hfill\qed
\end{rem}

\subsubsection{Kobayashi-Hitchin correspondence in the case $\lambda=0$}
\label{subsection;21.8.4.5}

Let $\nbigv$ be a locally free
$\nbigo_{\nbigmbar^{0}\setminus Z}
 (\ast H^{0}_{\infty})$-module.
Suppose that it is of Dirac type at each point of $Z$,
and that it is equipped with a tuple of filtered bundles
$\nbigp_{\ast}\nbigv=
\bigl(\nbigp_{\ast}(\nbigv_{|\nbigmbar^0\langle t\rangle\setminus Z})
 \,|\,t\in S^1_T\bigr)$
 over $\nbigv$
 which is good in the sense of
Definition {\rm\ref{df;20.7.28.20}}
and Definition {\rm\ref{df;17.10.28.10}}.
It is denoted by $\nbigp_{\ast}\nbigv$,
and called a good filtered bundle of Dirac type
over $(\nbigmbar^{0};Z,H^{0}_{\infty})$.
\index{good filtered bundle of Dirac type}

For any good filtered bundle of Dirac type 
$\nbigp_{\ast}\nbigv$
over $(\nbigmbar^{0};Z,H^{0}_{\infty})$,
we note that the numbers
$\deg\bigl(\nbigp_{\ast}(\nbigv_{|\nbigmbar^0\langle t\rangle\setminus Z})
\bigr)$
are well defined for $t\in S^1_T\setminus\pi(Z)$,
which induces affine functions
on the connected components of
$S^1_T\setminus \pi(Z)$.
We define
\begin{equation}
\label{eq;21.8.3.20}
 \deg(\nbigp_{\ast}\nbigv):=
 \frac{1}{T}
 \int_{0}^T
 \deg\bigl(\nbigp_{\ast}(\nbigv_{|\nbigmbar^0\langle t\rangle\setminus Z})
 \bigr)\,dt.
\end{equation}
\index{degree $\deg(\nbigp_{\ast}\nbigv)$}
We define the stability condition 
for good filtered bundles of Dirac type 
over $(\nbigmbar^{0};Z,H^{0}_{\infty})$
by using the degree as in the standard way.
\index{stability condition}

The following theorem is a natural analogue of 
Theorem \ref{thm;17.10.28.30} in the case $\lambda=0$,
in the context of periodic monopoles,
for which we apply the Kobayashi-Hitchin correspondence
for analytic stable bundles studied in \cite{Mochizuki-KH-infinite}.

\begin{thm}[Theorem
\ref{thm;17.9.30.20}]
\label{thm;21.8.3.22}
The construction from $(E,h,\nabla,\phi)$
to $\nbigp_{\ast}\nbige^0$
induces a bijective correspondence between
the equivalence classes of
monopoles of GCK-type
on $\nbigm\setminus Z$
and 
the equivalence classes of
polystable good filtered bundles of Dirac type
with degree $0$
on $(\nbigmbar^{0};Z,H^{0}_{\infty})$.
\end{thm}

\begin{rem}
Let $\Sigma$ be a compact Riemann surface,
and let $D\subset\Sigma$ denote a finite subset.
Let $g_{\Sigma\setminus D}$ be a K\"ahler metric
of $\Sigma\setminus D$.
Suppose that for each $P\in D$
there exists a neighbourhood $\Sigma_P$
such that
$\Sigma_P\setminus\{P\}$
with $g_{\Sigma\setminus D}$
is isomorphic to
$\{|z|>R\}$ with the metric $dz\,d\zbar$.
Let $Z$ denotes a finite subset of
$S^1_T\times(\Sigma\setminus D)$.
Theorem {\rm\ref{thm;17.10.28.31}} can be generalized
to correspondences
for monopoles on $(S^1_T\times(\Sigma\setminus D))\setminus Z$
with the metric $dt\,dt+g_{\Sigma\setminus D}$
satisfying similar conditions
around $(S^1_T\times D)\cup Z$.
\hfill\qed
\end{rem}

\subsubsection{$\nbigo_{\nbigmbar^{0}\setminus Z}(\ast H^0_{\infty})$-Modules
and $\cnum(w)$-modules with an automorphism}

Let $\nbigv$ be a locally free 
$\nbigo_{\nbigmbar^{0}\setminus Z}
(\ast H^0_{\infty})$-module
of Dirac type at $Z$.
Let $D$ denote the image of $Z$
by the projection $p_w:\nbigm^0=S^1_T\times\cnum_w\lrarr\cnum_w$.
Take a sufficiently small positive number $\epsilon$.
By the scattering map along the loop
$-\epsilon+s$ $(0\leq s\leq T)$,
we obtain the automorphism $\rho$
of the $\nbigo_{\proj^1}(\ast (D\cup\{\infty\}))$-module
$\nbigv_{|\nbigmbar^0\langle -\epsilon\rangle}$.

We set
$V:=H^0\bigl(
\proj^1,
\nbigv_{|\nbigmbar^0\langle -\epsilon\rangle}
\bigr)$,
and $\vecV:=\cnum(w)\otimes_{\cnum[w]} V$.
It is equipped with a $\cnum(w)$-linear automorphism $\rho$,
i.e.,
and $(\vecV,\rho)$ is a $0$-difference module.

For any $x\in \cnum$,
we set $m(x):=\bigl|p_w^{-1}(x)\cap Z\bigr|$.
If $m(x)>0$,
we obtain
$0\leq t^{(1)}_{x}<
t^{(2)}_{x}<\cdots<t^{m(x)}_{x}<T$
defined by
$p_w^{-1}(x)\cap Z
=\{(t^{(i)}_x,x)\}$.
We set
$\tau^{(i)}_x:=t^{(i)}_{x}/T$.
For $i=1,\ldots,m(x)-1$,
by choosing $t^{(i)}_{x}<s_i<t^{(i+1)}_{x}$,
we set
$L_{x,i}$
as $\nbigv_{(s_i,x)}\otimes \cnum[\![w-x]\!]$,
where $\nbigv_{(s_i,x)}$
denotes the stalk of $\nbigv$
at $(s_i,x)$.
Thus, we obtain a parabolic structure
at finite place
$(V,\{\vectau_x,\vecL_x\})$
of $\vecV$.

Let 
$(\nbigp_{\ast}(\nbigv_{|\nbigmbar^0\langle t\rangle })\,|\,t\in S^1_T)$
be a good filtered bundle
over $\nbigv$.
By the definition of good filtered bundles
over $\nbigv$ (Definition \ref{df;20.7.28.20}),
$\nbigp_{\ast}\nbigv_{|\nbigmbar^0\langle 0\rangle}$
induces a good parabolic structure of
$\vecV$ at $\infty$.

By this correspondence,
the degree (\ref{eq;21.8.3.20})
for $\nbigp_{\ast}\nbigv$
is translated to the degree
(\ref{eq;17.12.4.4})
for the parabolic difference module
$\vecV_{\ast}=
\bigl(\vecV,V,m,(\vectau_{x},\vecL_x)_{x\in\cnum},
\nbigp_{\ast}\vecVhat\bigr)$.
(See Lemma \ref{lem;21.9.6.5}.)
Hence, the stability condition for 
$\nbigp_{\ast}\nbigv$
is equivalent to
the stability condition for
$\vecV_{\ast}$.
Thus, we obtain the following proposition.

\begin{prop}[Proposition
\ref{prop;21.9.17.50}, Proposition \ref{prop;21.9.17.51}]
\label{prop;21.8.3.21}
The above construction induces an equivalence
between (stable, polystable)
good filtered bundle of Dirac type of degree $0$
and 
(stable, polystable)
parabolic $0$-difference modules of degree $0$.
\hfill\qed
 \end{prop}

Theorem \ref{thm;17.12.4.10} in the case $\lambda=0$
follows from Theorem \ref{thm;21.8.3.22}
and Proposition \ref{prop;21.8.3.21}.

\subsection{The correspondences in the general case}
\label{section;21.8.3.31}

Let us outline the construction for general $\lambda$
in a way as parallel to the case $\lambda=0$ as possible.

\subsubsection{Preliminary consideration}

In \S\ref{section;21.8.3.30},
we explained the construction of a $0$-difference module
from a monopole of GCK-type on $\nbigm\setminus Z$,
where one of the keys is the integrability
$[\nabla_t-\sqrt{-1}\phi,\nabla_{\wbar}]=0$
as in the previous works in \S\ref{section;20.8.9.1}.
If we choose
another $\real$-linear isomorphism
\begin{equation}
\label{eq;21.8.3.30}
 \real_{t}\times\cnum_w\simeq\real_{t_0}\times\cnum_{\beta_0}
\end{equation}
such that
$dt\,dt+dw\,d\wbar=dt_0\,dt_0+d\beta_0\,d\betabar_0$,
we also obtain
$[\nabla_{t_0}-\sqrt{-1}\phi,\nabla_{\betabar_0}]=0$.
It is natural to ask what objects we could obtain by using
the mini-complex structure induced by $(t_0,\beta_0)$ instead of $(t,w)$.
(See \S\ref{subsection;17.10.28.1} 
for mini-complex structure.)

\subsubsection{Mini-complex structure corresponding to 
the twistor parameter $\lambda$}
\label{subsection;21.8.3.1}

We can easily observe that isomorphisms (\ref{eq;21.8.3.30}) satisfying 
$dt\,dt+dw\,d\wbar=dt_0\,dt_0+d\beta_0\,d\betabar_0$
are naturally parameterized by $S^2=\proj^1$.
Conceptually, such a decomposition
is induced by a complex structure underlying
the hyperk\"ahler manifold $\real^4$ with
the standard Euclidean metric $\sum_{i=1}^4 dx_i^2$.
We choose a complex structure $(z,w)=(x_1+\sqrt{-1}x_2,x_3+\sqrt{-1}x_4)$
for which we have $dz\,d\zbar+dw\,d\wbar=\sum dx_i^2$.
We can regard the mini-complex manifold $\real_t\times\cnum_w$
as the quotient of $\cnum^2$ by the $\real_s$-action
defined as $s\bullet (z,w)=(z+s,w)$.
(See \S\ref{subsection;17.10.3.1}.)
In other words, we may regard that the coordinate system
$(t,w)$ on $\real_t\times\cnum_w$ is induced by
$(t,w)=(\Image z,w)$.
For the twistor parameter $\lambda\in\cnum$,
we have the corresponding complex structure of
$(\real^4,\sum dx_i^2)$,
for which a complex coordinate system is given by
$(\xi,\eta)=(z+\lambda\wbar,w-\lambda\zbar)$.
It induces a mini-complex structure on
$\real_t\times\cnum_w$,
which is different from the mini-complex structure
induced by $(t,w)$ unless $\lambda=0$.
To emphasize the mini-complex structure
depending on $\lambda$,
we use the notation $(\real_t\times\cnum_w)^{\lambda}$.

Note that the $\real_s$-action is expressed as
$s\bullet(\xi,\eta)=(\xi+s,\eta-\lambda s)$
in terms of $(\xi,\eta)$.
Let $(\alpha_0,\beta_0)$ be the complex coordinate system
defined as
\[
 (\alpha_0,\beta_0)
 =\frac{1}{1+|\lambda|^2}
 (\xi-\lambdabar\eta,\eta+\lambda\xi).
\]
Because the $\real_s$-action is expressed as
$s\bullet(\alpha_0,\beta_0)=(\alpha_0+s,\beta_0)$,
we obtain the induced mini-complex coordinate system
$(t_0,\beta_0)=(\Image\alpha_0,\beta_0)$
on $(\real_t\times\cnum_{w})^{\lambda}$,
i.e.,
\index{mini-complex coordinate system $(t_0,\beta_0)$}
\begin{equation}
\label{eq;20.7.30.3}
 t_0=
 \frac{1-|\lambda|^2}{1+|\lambda|^2}t
+\frac{2}{1+|\lambda|^2}\Image(\lambda\wbar),
\quad\quad
 \beta_0=\frac{1}{1+|\lambda|^2}
 \bigl(w+2\sqrt{-1}\lambda t+\lambda^2\wbar\bigr).
\end{equation}
We note $dt\,dt+dw\,d\wbar=dt_0\,dt_0+d\beta_0\,d\betabar_0$.
The $\seisuu$-action $\kappa$
is described as follows
in terms of $(t_0,\beta_0)$:
\begin{equation}
\label{eq;21.7.9.1}
 \kappa_n(t_0,\beta_0)=
 (t_0,\beta_0)
+nT\cdot
\Bigl(
 \frac{1-|\lambda|^2}{1+|\lambda|^2},\,\,
 \frac{2\sqrt{-1}\lambda}{1+|\lambda|^2}
\Bigr).
\end{equation}
Because
the complex vector fields
$\del_{t_0}$
and $\del_{\betabar_0}$
are $\kappa$-invariant,
we obtain the induced complex vector fields
on $\nbigm$,
which are also denoted by
$\del_{t_0}$ and $\del_{\betabar_0}$.

Let $\nbigm^{\lambda}$ denote the $3$-dimensional manifold
obtained as the $(\real_t\times\cnum_w)^{\lambda}$
by $\kappa$.
\index{mini-complex manifold $\nbigm^{\lambda}$}
Though $\nbigm^{\lambda}=\nbigm=\nbigm^{\lambda'}$
as $C^{\infty}$-manifolds,
the mini-complex structures of
$\nbigm^{\lambda}$ and $\nbigm^{\lambda'}$ are different
if $\lambda\neq \lambda'$.
A $C^{\infty}$-function $f$ on an open subset 
$U\subset\nbigmlambda$ is called mini-holomorphic
if $\del_{t_0}f=\del_{\betabar_0}f=0$.
Let $\nbigo_{\nbigm^{\lambda}}$
denote the sheaf of mini-holomorphic functions 
on $\nbigmlambda$.
\index{mini-complex manifold $\nbigm^{\lambda}$}
\index{sheaf $\nbigo_{\nbigmlambda}$}
\index{mini-complex coordinate system $(t_0,\beta_0)$}

\begin{rem}
Once the twistor parameter $\lambda$ is fixed,
the mini-complex coordinate system $(t_0,\beta_0)$
is uniquely determined by
$dt_0\,dt_0+d\beta_0\,d\betabar_0=dt\,dt+dw\,d\wbar$
up to the multiplication 
of a complex number $e^{\sqrt{-1}\varphi}$ $(\varphi\in\real)$
to $\beta_0$.
\hfill\qed 
\end{rem}

\subsubsection{Another coordinate system and
the compactification of $\nbigm^{\lambda}$}
\label{subsection;21.8.5.1}

As recalled in \S\ref{section;20.8.9.1},
it has been standard to use efficiently
the mini-complex coordinate system like $(t_0,\beta_0)$
in the study of monopoles,
which induces the orthogonal decomposition
of the $3$-dimensional Euclidean space.
However, in our study of periodic monopoles,
we also use another convenient coordinate system
$(t_1,\beta_1)$
given as follows:
\index{mini-complex coordinate system $(t_1,\beta_1)$}
\begin{equation}
\label{eq;21.7.9.2}
 (t_1,\beta_1)=
 \bigl(
 t_0+\Image(\lambdabar\beta_0),
 (1+|\lambda|^2)\beta_0
 \bigr)
=\bigl(
 t+\Image(\lambda\wbar),
 w+2\lambda\sqrt{-1}t+\lambda^2\wbar
 \bigr).
\end{equation}
We note that
a $C^{\infty}$-function $f$ on $U\subset\nbigm^{\lambda}$
is mini-holomorphic if and only if
$\del_{t_1}f=\del_{\betabar_1}f=0$.
The $\seisuu$-action $\kappa$ is described as
\begin{equation}
\label{eq;21.8.3.50}
 \kappa_n(t_1,\beta_1)
=(t_1,\beta_1)+n\cdot (T,2\sqrt{-1}\lambda T).
\end{equation}
We have the partial compactification
$(\real_t\times\cnum_w)^{\lambda}
=\real_{t_1}\times\cnum_{\beta_1}
\subset
\real_{t_1}\times\proj^1_{\beta_1}$.
By the same formula,
we define the $\seisuu$-action $\kappa$ on 
$\real_{t_1}\times\proj^1_{\beta_1}$.
The quotient space
$\nbigmbar^{\lambda}$
is naturally equipped with
the induced mini-complex structure,
and it is a compactification of $\nbigmlambda$.
\index{mini-complex manifold $\nbigmbar^{\lambda}$}
Let $H^{\lambda}_{\infty}:=
 \nbigmbar^{\lambda}\setminus \nbigm^{\lambda}
 \simeq S^1_T$.
Let $\nbigo_{\nbigmbar^{\lambda}}(\ast H^{\lambda}_{\infty})$
denote the sheaf of meromorphic functions on $\nbigmbar^{\lambda}$
whose poles are contained in $H^{\lambda}_{\infty}$.
\index{sheaf $\nbigo_{\nbigmbar^{\lambda}}(\ast H_{\infty}^{\lambda})$}
\index{space $H^{\lambda}_{\infty}$}

\begin{rem}
As in the case of $\lambda=0$
or the case of harmonic bundles
(see Remark {\rm\ref{rem;21.9.16.1}}),
we would like to consider meromorphic objects
or filtered objects on $(\nbigmbar^{\lambda},H^{\lambda}_{\infty})$
to keep the information of the growth orders with respect to $h$,
rather than the transcendental object on $\nbigm^{\lambda}$. 
It is the reason why we consider the compactification
$\nbigmbar^{\lambda}$.
\hfill\qed
\end{rem}

\begin{rem}
When we regard
the coordinate change from
$(t_1,\beta_1)$ to $(t_0,\beta_0)$
as a diffeomorphism
$\real_{t_1}\times\cnum_{\beta_1}
\simeq
\real_{t_0}\times\cnum_{\beta_0}$,
it does not extend to a diffeomorphism
$\real_{t_1}\times\proj^1_{\beta_1}\simeq
\real_{t_0}\times\proj^1_{\beta_0}$.
Similarly, the diffeomorphism
$\real_{t_1}\times\cnum_{\beta_1}\simeq
\real_{t}\times\cnum_{w}$
does not extend to a diffeomorphism
$\real_{t_1}\times\proj^1_{\beta_1}\simeq
\real_{t}\times\proj^1_{w}$.
Note that we may regard $\real_{t_1}\times\proj^1_{\beta_1}$
as the quotient of $\cnum_{\xi}\times\proj^1_{\eta}$ 
by the naturally induced $\real_s$-action.
\hfill\qed
\end{rem}

There are two reasons to use
$(t_1,\beta_1)$ instead of $(t_0,\beta_0)$.
It is one reason that
for the coordinate system $(t_0,\beta_0)$
the $\seisuu$-action in the $\real_{t_0}$-direction is trivial
if $|\lambda|=1$.
(Compare {\rm(\ref{eq;21.7.9.1})}
with {\rm(\ref{eq;21.8.3.50})}.)
See Remark \ref{rem;21.8.4.3} for another reason.
Clearly, we have $(t_0,\beta_0)=(t_1,\beta_1)$
in the case $\lambda=0$.

\begin{rem}
As explained above, we obtain
$(t_i,\beta_i)$ $(i=0,1)$ from $(\xi,\eta)$.
If $\lambda\neq 0$,
there exists another pair of mini-complex coordinate systems
$(t_i^{\dagger},\beta_i^{\dagger})$ $(i=0,1)$
obtained from
$(\xi^{\dagger},\eta^{\dagger})
=(\lambda^{-1}(w-\lambda\zbar),\lambda^{-1}(z+\lambda\wbar))
=(\lambda^{-1}\eta,\lambda^{-1}\xi)$
which we may use to obtain difference modules
from monopoles.
See Remark {\rm\ref{rem;21.8.10.3}}. 
\hfill\qed
\end{rem}

\subsubsection{Mini-holomorphic bundles associated with monopoles}

Let $Z$ be a finite subset of $\nbigm$.
Let $(E,h,\nabla,\phi)$ be a monopole on
$\nbigm\setminus Z$.
Because the Bogomolny equation implies that
$[\nabla_{t_0}-\sqrt{-1}\phi,\nabla_{\betabar_0}]=0$,
we obtain the locally free
$\nbigo_{\nbigm^{\lambda}\setminus Z}$-module
$\nbige^{\lambda}$
as the sheaf of local sections $s$ of $E$
such that 
$(\nabla_{t_0}-\sqrt{-1}\phi)s=
 \nabla_{\betabar_0}s=0$.
\index{$\nbigo_{\nbigm^{\lambda}\setminus Z}$-module $\nbige^{\lambda}$}

If $P\in Z$ is a Dirac type singularity of $(E,h,\nabla,\phi)$,
then $\nbigelambda$ is of Dirac type at $P$
as in \S\ref{subsection;21.8.3.40}.

\subsubsection{Meromorphic extension and filtered extension at
infinity}
\label{subsection;17.10.28.40}

Let $U$ be any open subset in $\nbigmbar^{\lambda}\setminus Z$.
As in \S\ref{subsection;21.8.4.4},
let $\nbigp^h\nbigelambda(U)$
denote the space of sections of $\nbigelambda$
on $U\setminus H^{\lambda}_{\infty}$
satisfying the following.
\begin{itemize}
\item For any $P\in U\cap H^{\lambda}_{\infty}$,
 there exists a neighbourhood $U_P$ of $P$ in $U$
 such that 
 $|s_{|U_P\setminus H^{\lambda}_{\infty}}|_h=O\bigl(|w|^{N}\bigr)$
 for some $N$.
\end{itemize}
Thus, we obtain the 
$\nbigo_{\nbigmbar^{\lambda}\setminus Z}
 (\ast H^{\lambda}_{\infty})$-module
 $\nbigp^h\nbige^{\lambda}$.
\index{sheaf $\nbigp^h\nbige^{\lambda}$}
We shall prove the following.
\begin{prop}[Proposition 
\ref{prop;17.9.17.50}]
\label{prop;21.8.4.1}
 $\nbigp^h\nbigelambda$ is a locally free
$\nbigo_{\nbigmbar^{\lambda}\setminus Z}
 (\ast H^{\lambda}_{\infty})$-module.
\end{prop}

Let $\pi^{\lambda}:\nbigmbar^{\lambda}\lrarr S^1_T$
denote the map induced by
$(t_1,\beta_1)\longmapsto t_1$.
\index{projection $\pi^{\lambda}$}
We set
$\nbigmbar^{\lambda}\langle t_1\rangle=(\pi^{\lambda})^{-1}(t_1)$
and
$\nbigm^{\lambda}\langle t_1\rangle=
\nbigm^{\lambda}\cap\nbigmbar^{\lambda}\langle t_1\rangle$
for $t_1\in S^1_T$.
\index{space $\nbigmbar^{\lambda}\langle t_1\rangle$}
\index{space $\nbigm^{\lambda}\langle t_1\rangle$}
By taking the restriction,
we obtain the locally free 
$\nbigo_{\nbigm^{\lambda}\langle t_1\rangle\setminus Z}(\ast \infty)$-module
$\nbigp^h(\nbigelambda_{|\nbigm^{\lambda}\langle t_1\rangle \setminus Z})$.
As in \S\ref{subsection;17.10.27.10} and \S\ref{subsection;17.10.27.11},
we obtain an increasing sequence 
$\nbigp^h_a\bigl(
\nbigelambda_{|\nbigm^{\lambda}\langle t_1\rangle\setminus Z}
 \bigr)$ $(a\in\real)$
of $\nbigo_{\nbigmbar^{\lambda}\langle t_1\rangle\setminus Z}$-modules
by considering local sections $s$
satisfying
$|s|_h=O(|w|^{a+\epsilon})$ for any $\epsilon>0$.
\index{sheaf $\nbigp^h_a\bigl(
\nbigelambda_{|\nbigm^{\lambda}\langle t_1\rangle\setminus Z}
 \bigr)$}

\begin{thm}[Theorem
\ref{thm;17.10.5.130}]
\label{thm;21.8.4.2}
The tuple
$\bigl(
 \nbigp^h_{\ast}\bigl(
 \nbigelambda_{|\nbigm^{\lambda}\langle t_1\rangle\setminus Z}
 \bigr)
 \,\big|\,
 t_1\in S^1_T
\bigr)$
is a good filtered bundle
over $\nbigp^h\nbigelambda$
in the sense of 
Definition {\rm\ref{df;20.7.28.20}}
and Definition {\rm\ref{df;17.10.28.10}}.
(See Theorem {\rm\ref{thm;17.10.7.10}}
for a more detailed description.)
\index{good filtered bundle}
\end{thm}

As we shall see in Proposition \ref{prop;17.10.14.21},
the tuple
$\bigl(
 \nbigp^h_{\ast}\bigl(
 \nbigelambda_{|\nbigm^{\lambda}\langle t_1\rangle\setminus Z}
 \bigr)
 \,\big|\,
 t_1\in S^1_T
 \bigr)$
 determines the behaviour of the metric $h$
around $H^{\lambda}_{\infty}$ up to boundedness.
We shall also show that
the compatibility with such filtrations
implies the GCK condition around infinity.
(See \S\ref{subsection;17.10.13.500}.)

We shall often use the abbreviation 
$\nbigp^h_{\ast}\nbige^{\lambda}$
to denote
$\bigl(
\nbigp^h_{\ast}\bigl(
\nbigelambda_{|\nbigm^{\lambda}\langle t_1\rangle\setminus Z}
\bigr)
\,\big|\,
t_1\in S^1_T
\bigr)$,
but we remark that 
the filtrations depend on $t_1\in S^1_T$
as basic examples show
(see \S\ref{subsection;17.9.26.1}).

\begin{rem}
\label{rem;21.8.4.3}
In the proof of 
 Proposition {\rm\ref{prop;21.8.4.1}}
and Theorem {\rm\ref{thm;21.8.4.2}},
one of the key facts is that
the holomorphic vector bundle with a Hermitian metric
$(\nbigelambda,h)_{|\nbigm^{\lambda}\langle t_1\rangle\setminus Z}$
is acceptable (Lemma {\rm\ref{lem;21.7.9.3}}).
(See {\rm\S\ref{subsection;20.8.8.31}}
for the acceptability condition.)
Hence, we may apply the general result
for acceptable bundles in {\rm\cite{Simpson90}}
to extend
$(\nbigelambda,h)_{|\nbigm^{\lambda}\langle t_1\rangle\setminus Z}$
 across $\infty\in\nbigmbar^{\lambda}\langle t_1\rangle\simeq\proj^1_{\beta_1}$
as a filtered bundle.
See {\rm\S\ref{subsection;21.8.10.4}}
for a more detailed explanation.

Let $\pi^{\lambda}_0:\nbigm^{\lambda}\lrarr S^1_T$
denote the projection induced by
$(t_0,\beta_0)\longmapsto t_0$.
We note that
$(\nbigelambda,h)_{|(\pi^{\lambda}_0)^{-1}(t_0)\setminus Z}$
is not acceptable, in general.

We may use $(t_1,\beta_1)$ 
not only $(t_0,\beta_0)$
thanks to the notion of mini-complex structure
instead of a pair of differential operators
$\nabla_{t_0}-\sqrt{-1}\phi$ and $\nabla_{\betabar_0}$
satisfying the integrability condition
$[\nabla_{t_0}-\sqrt{-1}\phi,\nabla_{\betabar_0}]$. 
This is a merit to consider mini-complex structure.
\hfill\qed
\end{rem}

\begin{rem}
To define the concept of good filtered bundle
in Theorem {\rm\ref{thm;21.8.4.2}},
we need to know the classification of 
locally free
$\nbigo_{\Hhat_{\infty}^{\lambda}}(\ast H_{\infty}^{\lambda})$-modules
of finite rank,
where $\Hhat_{\infty}^{\lambda}$ is the formal space
obtained as the completion of 
$\nbigmbar^{\lambda}$ along $H^{\lambda}_{\infty}$.
\index{ringed space $\Hhat_{\infty}^{\lambda}$}
Because such modules are naturally equivalent to
formal difference modules (see {\rm\S\ref{subsection;20.7.18.4}}),
we may apply the classical results
on the classification of formal difference modules
mentioned in {\rm\S\ref{subsection;18.2.1.1}}.
For our purpose,
it is also convenient to use
another equivalence of formal difference modules of level $\leq 1$
and formal differential modules whose Poincar\'e rank is strictly
smaller than $1$,
which will be explained in 
{\rm\S\ref{subsection;17.10.28.20}}.
\hfill\qed
\end{rem}

\subsubsection{Kobayashi-Hitchin correspondence
 of periodic monopoles of GCK type}

\label{subsection;17.10.28.41}

Let $\nbigv$ be a locally free
$\nbigo_{\nbigmbar^{\lambda}\setminus Z}
 (\ast H^{\lambda}_{\infty})$-module.
Suppose that it is of Dirac type at each point of $Z$,
and that it is equipped with a tuple of filtered bundles
$\nbigp_{\ast}\nbigv=
\bigl(\nbigp_{\ast}(\nbigv_{|\nbigmbar^{\lambda}\langle t_1\rangle})
 \,|\,t_1\in S^1_T\bigr)$
 over $\nbigv$
 which is good in the sense of
 Definition {\rm\ref{df;20.7.28.20}}
and Definition {\rm\ref{df;17.10.28.10}}.
We denote such object by
$\nbigp_{\ast}\nbigv$,
and call a good filtered bundle of Dirac type
over $(\nbigmbar^{\lambda};Z,H^{\lambda}_{\infty})$.
\index{good filtered bundle of Dirac type}

For any good filtered bundle of Dirac type 
$\nbigp_{\ast}\nbigv$
over $(\nbigmbar^{\lambda};Z,H^{\lambda}_{\infty})$,
we note that the numbers
$\deg\bigl(\nbigp_{\ast}(\nbigv_{|\nbigmbar^{\lambda}\langle t_1\rangle})\bigr)$
are well defined for $t_1\in S^1_T\setminus\pi^{\lambda}(Z)$,
which induces affine functions
on the connected components of
$S^1_T\setminus \pi^{\lambda}(Z)$.
We define
\index{degree $\deg(\nbigp_{\ast}\nbigv)$}
\begin{equation}
\label{eq;17.12.4.3}
 \deg(\nbigp_{\ast}\nbigv):=
 \frac{1}{T}
 \int_{0}^T
 \deg\bigl(\nbigp_{\ast}(\nbigv_{|\nbigmbar^{\lambda}\langle t_1\rangle})\bigr)\,dt_1.
\end{equation}
We define the stability condition 
for good filtered bundles of Dirac type 
over $(\nbigmbar^{\lambda};Z,H^{\lambda}_{\infty})$
by using the degree as in the standard way.

The following theorem is
a generalization of Theorem \ref{thm;21.8.3.22}
to the case of general $\lambda$,
and it is a natural analogue of 
Theorem \ref{thm;17.10.28.30}
in the context of periodic monopoles.
We again apply the Kobayashi-Hitchin correspondence
for analytic stable bundles studied in \cite{Mochizuki-KH-infinite}.

\begin{thm}[Theorem \ref{thm;17.9.30.20}]
\label{thm;17.10.28.31}
The construction from $(E,h,\nabla,\phi)$
to $\nbigp_{\ast}\nbigelambda$
induces a bijective correspondence between
the equivalence classes of
 monopoles of GCK-type
 on $\nbigm\setminus Z$
and 
the equivalence classes of
polystable good filtered bundles of Dirac type
 with degree $0$
 on $(\nbigmbar^{\lambda};Z,H^{\lambda}_{\infty})$.
\end{thm}

The following is an analogue of 
Corlette-Simpson correspondence
between flat bundles and Higgs bundles,
which is an immediate consequence of 
Theorem \ref{thm;17.10.28.31}.
\begin{cor}[Corollary \ref{cor;17.10.28.32}]
For any $\lambda\in\cnum$,
there exists the natural bijective correspondence
of the following objects
through periodic monopoles of GCK-type:
\begin{itemize}
\item
 Polystable good filtered bundles of Dirac type
 with degree $0$
 on $(\nbigmbar^{0};Z,H^{0}_{\infty})$.
\item
 Polystable good filtered bundles of Dirac type
 with degree $0$
 on $(\nbigmbar^{\lambda};Z,H^{\lambda}_{\infty})$.
\hfill\qed
\end{itemize}
\end{cor}

\subsubsection{Difference modules and 
 $\nbigo_{\nbigmbar^{\lambda}\setminus Z}
 (\ast H^{\lambda}_{\infty})$-modules}
\label{subsection;21.8.3.2}

Let $\nbigv$ be a locally free 
$\nbigo_{\nbigmbar^{\lambda}\setminus Z}
 (\ast H^{\lambda}_{\infty})$-module
 of Dirac type at $Z$.
Let $p_1:\real_{t_1}\times\proj^1_{\beta_1}\lrarr\real_{t_1}$
and $p_2:\real_{t_1}\times\proj^1_{\beta_1}\lrarr\proj^1_{\beta_1}$
denote the projections.
Let $\Ztilde\subset\real_{t_1}\times\proj^1_{\beta_1}$ 
be the pull back of $Z$ by the projection
$\varpi^{\lambda}:
 \real_{t_1}\times\proj^1_{\beta_1}\lrarr\nbigmbar^{\lambda}$.
\index{projection $\varpi^{\lambda}$}
We put
\[
 D:=p_2\Bigl(
 \Ztilde\cap p_1^{-1}\bigl(\{0\leq t_1<T\}\bigr)
 \Bigr).
\]
Take a sufficiently small positive number $\epsilon$.
By the scattering map,
we obtain the following isomorphism 
of $\nbigo_{\proj^1}(\ast (D\cup\{\infty\}))$-modules:
\begin{equation}
 \label{eq;17.12.3.1}
 (\varpi^{\lambda})^{\ast}
 (\nbigv)_{|p_1^{-1}(-\epsilon)}(\ast D)
\simeq
 (\varpi^{\lambda})^{\ast}
 (\nbigv)_{|p_1^{-1}(T-\epsilon)}(\ast D).
\end{equation}
We also have the following natural isomorphism
\begin{equation}
\label{eq;17.12.3.2}
 \Phi^{\ast}
 (\varpi^{\lambda})^{\ast}
 (\nbigv)_{|p_1^{-1}(T-\epsilon)}
\simeq
 (\varpi^{\lambda})^{\ast}
 (\nbigv)_{|p_1^{-1}(-\epsilon)}.
\end{equation}
We set
$V:=H^0\bigl(
 \proj^1,
 (\varpi^{\lambda})^{\ast}
 (\nbigv)_{|p_1^{-1}(-\epsilon)}
 \bigr)$.
Let $\cnum[\beta_1]_D$
denote the localization of $\cnum[\beta_1]$
by $\prod_{x\in D}(\beta_1-x)$.
By (\ref{eq;17.12.3.1}) and (\ref{eq;17.12.3.2}),
we obtain the $\cnum$-linear isomorphism
\begin{equation}
\label{eq;20.7.29.10}
 \Phi^{\ast}:V\otimes_{\cnum[\beta_1]}\cnum[\beta_1]_D
 \lrarr V\otimes_{\cnum[\beta_1]}
 \cnum[\beta_1]_{\Phi^{-1}(D)}.
\end{equation}
We set $\vecV:=\cnum(\beta_1)\otimes V$.
It is equipped with a $\cnum$-linear automorphism
$\Phi^{\ast}$ induced by (\ref{eq;20.7.29.10}),
and $(\vecV,\Phi^{\ast})$ is a difference module.

For any $x\in \cnum$,
we set $m(x):=\bigl|p_2^{-1}(x)\cap\Ztilde\bigr|$.
If $m(x)>0$,
we obtain
$0\leq t^{(1)}_{1,x}<
t^{(2)}_{1,x}<\cdots<t^{m(x)}_{1,x}<T$
as $p_2^{-1}(x)\cap \Ztilde$.
We set
$\tau^{(i)}_x:=t^{(i)}_{1,x}/T$.
For $i=1,\ldots,m(x)-1$,
by choosing $t^{(i)}_{1,x}<s_i<t^{(i+1)}_{1,x}$,
we set
$L_{x,i}$
as $\nbigv_{(s_i,x)}\otimes \cnum[\![\beta_1-x]\!]$,
where $\nbigv_{(s_i,x)}$
denotes the stalk of $\nbigv$
at $(s_i,x)$.
Thus, we obtain a parabolic structure
at finite place
$(V,\{\vectau_x,\vecL_x\})$
of $\vecV$.

Let 
$\bigl(\nbigp_{\ast}(
 \nbigv_{|\nbigmbar^{\lambda}\langle t_1\rangle\setminus Z})
  \,|\,t_1\in S^1_T\bigr)$
be a good filtered bundle
over $\nbigv$.
We obtain the filtered bundle
$\nbigp_{\ast}\bigl(
(\varpi^{\lambda})^{-1}(\nbigv)_{|p_1^{-1}(0)}
\bigr)$
over $(\varpi^{\lambda})^{-1}(\nbigv)_{|p_1^{-1}(0)}$
induced by the filtered bundle
$\nbigp_{\ast}(\nbigv_{|\nbigmbar^{\lambda}\langle 0\rangle})$.
By the definition of good filtered bundles
over $\nbigv$ (Definition \ref{df;20.7.28.20}),
it induces a good parabolic structure of
$\vecV$ at $\infty$.

By this correspondence,
the degree in (\ref{eq;17.12.4.3})
for $\nbigp_{\ast}\nbigv$
is translated to the degree
(\ref{eq;17.12.4.4})
for the parabolic difference module
$\vecV_{\ast}=
\bigl(\vecV,V,m,(\vectau_{x},\vecL_x)_{x\in\cnum},
\nbigp_{\ast}\vecVhat\bigr)$.
(See Lemma \ref{lem;21.9.6.5}.)
Hence, the stability condition for 
$\nbigp_{\ast}\nbigv$
is equivalent to
the stability condition for
$\vecV_{\ast}$.
Thus, we obtain the following proposition.

\begin{prop}[Proposition
\ref{prop;21.9.17.50}, Proposition \ref{prop;21.9.17.51}]
\label{prop;17.12.4.11}
The above construction induces an equivalence
between (stable, polystable)
good filtered bundle of Dirac type of degree $0$
and 
(stable, polystable)
difference modules of degree $0$.
\hfill\qed
\end{prop}

Theorem \ref{thm;17.12.4.10} for general $\lambda$
follows from Theorem \ref{thm;17.10.28.31}
and Proposition \ref{prop;17.12.4.11}.

\section{Asymptotic behaviour of
 periodic monopoles of GCK-type}
\label{section;21.8.4.21}

In the theory of wild harmonic bundles $(E,\delbar_E,\theta,h)$,
the first task is to study the asymptotic behaviour of
$h$ and $\theta$ around the singularity,
called Simpson's main estimate.
See \cite[Theorem 1]{Simpson90} for the tame case,
and \cite[\S7.2]{Mochizuki-wild} for the wild case.
The estimate is fundamental not only for
the study of
the filtered extension $(\nbigp^h_{\ast}E,\theta)$
of the Higgs bundle $(E,\delbar_E,\theta)$
but also for the study of the filtered extensions
$(\nbigp^h_{\ast}E^{\lambda},\DDlambda)$ ($\lambda\in\cnum$)
of the $\lambda$-flat bundles
$(E^{\lambda},\DDlambda)$ underlying the harmonic bundle.

We pursue a similar route in our study of periodic monopoles
by following the idea that
periodic monopoles are regarded as
harmonic bundles of infinite rank.
(See \S\ref{subsection;21.8.4.22}.)
Indeed, for the results in
\S\ref{subsection;21.8.4.4}--\S\ref{subsection;21.8.4.5}
and \S\ref{subsection;17.10.28.40}--\S\ref{subsection;17.10.28.41},
it is fundamental to understand the asymptotic behaviour of
monopoles $(E,h,\nabla,\phi)$ on
$\nbigb^{\ast}(R)=S^1_T\times\{w\in\cnum\,|\,|w|>R\}$
such that
$F(\nabla)\to 0$ and $|\phi|_h=O(\log|w|)$ as $|w|\to\infty$,
which we call GCK-condition.
We shall briefly describe the results in this section.

\subsection{Setting}
\label{subsection;21.8.29.1}
It is more convenient to study monopoles
on a ramified covering.
Namely, 
we set $U^{\ast}_{w,q}(R_1):=\{w_q\in\cnum\,|\,|w^q_q|>R_1\}$
for $q\in\seisuu_{\geq 1}$
and for $R_1>0$.
\index{space $U^{\ast}_{w,q}(R)$}
There exists the natural map
$U^{\ast}_{w,q}(R_1)\lrarr \cnum_w$
defined by $w_q\longmapsto w_q^q$.
We consider 
monopoles $(E,h,\nabla,\phi)$
on $\nbigb_q^{\ast}(R_1)=S^1_T\times U^{\ast}_{w,q}(R_1)$
with respect to
the Riemannian metric 
$dt\,dt+dw\,d\wbar=dt\,dt+q^2|w_q|^{2(q-1)}dw_q\,d\wbar_q$.
\index{space $\nbigb_q^{\ast}(R)$}
In the rest of this section,
we impose the GCK-condition,
i.e.,
$F(\nabla)\to 0$
and $|\phi|_h=O(\log|w_q|)$ as $|w_q|\to\infty$.

For each $\lambda\in\cnum$,
$\nbigb^{\ast}_q(R_1)$ has the mini-complex structure
induced by the covering map
$\nbigb^{\ast}_q(R_1)\lrarr \nbigb^{\ast}(R)$
and the mini-complex structure of $\nbigb^{\ast}(R)$
as an open subset of $\nbigm^{\lambda}$.
(See \S\ref{subsection;21.8.3.1}
for $\nbigm^{\lambda}$.)
When we emphasize the mini-complex structure depending on $\lambda$,
we use the notation $\nbigb^{\lambda\ast}_q(R_1)$.
Let $(E^{\lambda},\delbar_{E^{\lambda}})$ denote the mini-holomorphic vector bundle
on $\nbigb_{q}^{\lambda\ast}(R_1)$
underlying the monopole $(E,h,\nabla,\phi)$.
(See \S\ref{subsection;17.10.28.42} for mini-holomorphic bundles.)

As in the case of harmonic bundles,
we begin with the analysis of $(E^0,\delbar_{E^0})$ with $h$,
and we prove Theorem \ref{thm;17.10.28.50} below
which is an analogue of the asymptotic orthogonality for wild harmonic bundles
contained in Simpson's main estimate
(see \cite[Theorem 8.2.1]{Mochizuki-wild}).
It is fundamental even for the study of
$(E^{\lambda},\delbar_{E^{\lambda}},h)$
as in the case of harmonic bundles.
To show it, we outline our proof of
Proposition \ref{prop;21.8.4.1}
and Theorem \ref{thm;21.8.4.2}.

\subsection{Decomposition of mini-holomorphic bundles}

By considering the monodromy along $S^1_T\times\{w_q\}$,
we obtain the automorphism $M(w_q)$
of $E^0_{|\{0\}\times U^{\ast}_{w,q}(R_1)}$.
The eigenvalues of $M(w_q)$ $(w_q\in U^{\ast}_{w,q})$
determine a complex curve $\Sp(E^0)$
in $U^{\ast}_{w,q}(R_1)\times\cnum^{\ast}$,
called the spectral curve of the periodic monopole.
\index{spectral curve $\Sp(E^0)$}
(See \cite{Cherkis-Kapustin1, Cherkis-Kapustin2}.)
Set $U_{w,q}(R_1):=U^{\ast}_{w,q}(R_1)\cup\{\infty\}$
in $\proj^1_{w_q}$.
Under the GCK-condition,
the closure 
$\overline{\Sp(E^0)}$
of the spectral curve $\Sp(E^0)$
in $U_{w,q}(R_1)\times\proj^1$ is also
a complex analytic curve.
(See Proposition \ref{prop;20.7.29.20}.)
Hence, after taking a ramified covering,
we may assume that there exists
the following decomposition:
\[
 \overline{\Sp(E^0)}
=\coprod_{i\in\Lambda} S_i.
\]
Here,
each $S_i$ is a graph of 
a meromorphic function $g_i$ 
on $(U_{w,q},\infty)$.
There exist $(\ell_i,\alpha_i)\in\seisuu\times\cnum^{\ast}$
such that $g_i\sim \alpha_iw_q^{-\ell_i}$.
We obtain 
\[
 \overline{\Sp(E^0)}
=\coprod_{(\ell,\alpha)\in \seisuu\times\cnum^{\ast}}
  \Bigl(
 \coprod_{\substack{i\in\Lambda\\ (\ell_i,\alpha_i)=(\ell,\alpha)}}
 S_i
 \Bigr).
\]
We obtain the corresponding decomposition
of mini-holomorphic bundles:
\begin{equation}
\label{eq;20.8.7.1}
 (E^0,\delbar_{E^0})=\bigoplus_{(\ell,\alpha)}
 (E_{\ell,\alpha},\delbar_{E_{\ell,\alpha}}).
\end{equation}

\subsection{The induced Higgs bundles}
\label{subsection;21.8.5.2}

\subsubsection{Preliminary (1)}
\label{subsection;21.8.29.2}
Let $V$ be a locally free $\nbigo_{U_{w,q}^{\ast}(R_1)}$-module
with an endomorphism $g$.
Let $\Psi_q:\nbigb_q^{\ast}(R_1)\lrarr U_{w,q}^{\ast}(R_1)$
denote the projection.
Let $\Vtilde$ be the $C^{\infty}$-bundle
on $\nbigb_q^{\ast}(R_1)$
obtained as the pull back of $V$ by $\Psi_q$.
We obtain the naturally defined operator
$\del_{\Vtilde,\wbar}$ on $\Vtilde$
determined by
$\del_{\Vtilde,\wbar}(f \Psi_q^{-1}(s))=
 \del_{\wbar}(f)\Psi_q^{-1}(s)+f\Psi_q^{-1}(\del_{V,\wbar}s)$
for any $f\in C^{\infty}(\nbigb_q^{\ast}(R_1))$
and $C^{\infty}$ section $s$ of $V$ on $U_{w,q}^{\ast}(R_1)$.
We also have the operator
$\del_{\Vtilde,t}$ on $\Vtilde$
determined by the following condition
for $f$ and $s$ as above:
\[
 \del_{\Vtilde,t}(f\Psi_q^{-1}(s))
=\del_t(f)\Psi_q^{-1}(s)
-2\sqrt{-1}f\Psi_q^{-1}(gs).
\]
Thus, we obtain a mini-holomorphic structure
$\delbar_{\Vtilde}$ on $\Vtilde$ on $\nbigb^{0\ast}_q(R_1)$.
We set $\Psi_q^{\ast}(V,g):=(\Vtilde,\delbar_{\Vtilde})$.
\index{mini-holomorphic bundle $\Psi_q^{\ast}(V,g)$}
\begin{rem}
The monodromy along $S^1_T\times\{w_q\}$ is
$\exp(2\sqrt{-1}Tg(w_q))$.
\hfill\qed
\end{rem}
 
\subsubsection{Preliminary (2)}

For each $(\ell,\alpha)\in\seisuu\times\cnum^{\ast}$,
there exists a basic example of monopole
\index{monopole $\LL_q^{\ast}(\ell,\alpha)$}
\begin{equation}
\label{eq;21.8.5.3}
 (\LL^{\ast}_q(\ell,\alpha),h_{\LL,q,\ell,\alpha},
 \nabla_{\LL,q,\ell,\alpha},\phi_{\LL,q,\ell,\alpha})
\end{equation}
on $\nbigb_q^{\ast}(R_1)$,
for which the monodromy along $S^1_T\times\{w_q\}$
with respect to $\nabla_t-\sqrt{-1}\phi$
is given as the multiplication of $\alpha w_q^{-\ell}$.
(See \S\ref{subsection;21.8.14.2}.)
Let $\LL^{0\ast}_q(\ell,\alpha)$ denote the underlying
mini-holomorphic bundle on $\nbigb_q^{0\ast}(R_1)$.
\index{mini-holomorphic bundle $\LL^{0\ast}_q(\ell,\alpha)$}

\subsubsection{The induced Higgs bundles}

It is not difficult to see that,
for each $(E_{\ell,\alpha},\delbar_{E_{\ell,\alpha}})$
in (\ref{eq;20.8.7.1}),
there exist a holomorphic vector bundle with an endomorphism
$(V_{\ell,\alpha},\delbar_{V_{\ell,\alpha}},g_{\ell,\alpha})$
on $U_{w,q}^{\ast}(R_1)$
and an isomorphism
\[
 (E_{\ell,\alpha},\delbar_{E_{\ell,\alpha}})
\simeq
 \LL_q^{0\ast}(\ell,\alpha)
 \otimes
 \Psi_q^{\ast}(V_{\ell,\alpha},\delbar_{V_{\ell,\alpha}},g_{\ell,\alpha}).
\]
We may choose $g_{\ell,\alpha}$ such that
the eigenvalues of $g_{\ell,\alpha|w_q}$
goes to $0$ as $|w_{q}|\to\infty$.
By setting $\theta_{\ell,\alpha}:=g_{\ell,\alpha}\,dw$,
we obtain Higgs bundles
$(V_{\ell,\alpha},\delbar_{V_{\ell,\alpha}},\theta_{\ell,\alpha})$.

\subsection{Asymptotic orthogonality}

Let $h_{\ell,\alpha}$ be the restriction of $h$
to $E_{\ell,\alpha}$.
We obtain the metric
$h_{\ell,\alpha}\otimes h^{-1}_{\LL,q,\ell,\alpha}$
of $\Psi_q^{\ast}(V_{\ell,\alpha},\delbar_{V_{\ell,\alpha}},g_{\ell,\alpha})$.
There exists the Fourier expansion of
$h_{\ell,\alpha}\otimes h^{-1}_{\LL,q,\ell,\alpha}$
along the fibers $S^1_T\times\{w_q\}$,
and it turns out that 
the invariant part induces a metric
$h_{V,\ell,\alpha}$ of $V_{\ell,\alpha}$.
We obtain the following metric of $E$:
\[
 h^{\shikaku}:=
 \bigoplus_{(\ell,\alpha)}
 h_{\LL,q,\ell,\alpha}
 \otimes
 \Psi_q^{-1}\bigl(
 h_{V,\ell,\alpha}
 \bigr).
\]
Let $s$ be the automorphism of $E$
determined by $h=h^{\shikaku}\cdot s$.
The following theorem implies
that the difference of $h$ and $h^{\shikaku}$ decays rapidly.
\begin{thm}[Theorem
\ref{thm;16.9.20.22}, Proposition \ref{prop;17.10.5.3}]
\label{thm;17.10.28.50}
For any $m\in\seisuu_{\geq 0}$,
there exist positive constants $C_i(m)$ $(i=1,2)$
such that
\[
 \Bigl|
 \nabla_{\kappa_1}\circ
 \cdots
 \circ\nabla_{\kappa_m}
 (s-\id_E)
 \Bigr|_h
\leq
 C_1(m)\exp\bigl(-C_2(m)|w_q^q|\bigr)
\]
for any
$(\kappa_1,\ldots,\kappa_m)
\in \{t,w,\wbar\}^m$.
 \end{thm}

As a consequence of Theorem \ref{thm;17.10.28.50},
$(V_{\ell,\alpha},\delbar_{V_{\ell,\alpha}},
 \theta_{\ell,\alpha},h_{V,\ell,\alpha})$
asymptotically satisfies the Hitchin equation.
Namely,
let $F(h_{V,\ell,\alpha})$ denote the curvature
of the Chern connection determined by
$\delbar_{V_{\ell,\alpha}}$ and $h_{V,\ell,\alpha}$,
and let $\theta^{\dagger}_{\ell,\alpha}$
be the adjoint of $\theta_{\ell,\alpha}$
with respect to $h_{V,\ell,\alpha}$.
Then, we obtain the following decay
for some $\epsilon>0$:
\[
 F(h_{V,\ell,\alpha})
+\bigl[
 \theta_{\ell,\alpha},\,
 \theta^{\dagger}_{\ell,\alpha}
 \bigr]
=O\Bigl(
 \exp(-\epsilon|w_q^q|)
 \Bigr).
\]
We also obtain similar estimates
for any derivatives of the left hand side.
In this sense,
$(V_{\ell,\alpha},\delbar_{V_{\ell,\alpha}},
 \theta_{\ell,\alpha},h_{V,\ell,\alpha})$
is an asymptotic harmonic bundle.

\subsection{Curvature decay}

Many of the estimates for harmonic bundles
can be generalized to estimates for asymptotic harmonic bundles.
(See \S\ref{subsection;17.10.21.10}
and \cite[\S5.5]{Mochizuki-doubly-periodic}.)
It allows us to obtain estimates for 
periodic monopoles of GCK-type.
For example, we obtain the estimate
for the curvature $F(\nabla)$ of the monopole $(E,h,\nabla,\phi)$
of GCK-type.
\begin{cor}[Corollary
\ref{cor;17.10.6.1}]
\label{cor;21.8.5.5}
For the expression
\[
F(\nabla)=
F(\nabla)_{w,\wbar}dw\,d\wbar
+F(\nabla)_{w,t}dw\,dt
+F(\nabla)_{\wbar,t}d\wbar\,dt,
\]
the following estimates hold:
\[
 \bigl|
 F(\nabla)_{w\wbar}
 \bigr|_h
=O\Bigl(
 |w_q^q|^{-2}(\log|w_q|)^{-2}
 \Bigr),
\]
\[
 \bigl|
 F(\nabla)_{wt}
 \bigr|_h
=O(|w_q^q|^{-1}),
\quad
 \bigl|
 F(\nabla)_{\wbar t}
 \bigr|_h
=O(|w_q^q|^{-1}).
\]
\end{cor}

The estimates in Corollary \ref{cor;21.8.5.5}
are useful in the study of
$(E^{\lambda},\delbar_{E^{\lambda}})$ for any $\lambda$.
Let $(t_1,\beta_1)$ be the mini-complex coordinate system
as in \S\ref{subsection;21.8.5.1}.
Note that the complex vector fields
$\del_{\betabar_1}$ are defined on $\nbigb^{\ast}_q(R_1)$.
We obtain the differential operator $\del_{E^{\lambda},\betabar_1}$
on $E^{\lambda}$ induced by
$\del_{\betabar_1}$ and the mini-holomorphic structure
$\delbar_{E^{\lambda}}$.
We also obtain $\del_{E^{\lambda},h,\beta_1}$
induced by $h$ and $\del_{E^{\lambda},\betabar_1}$ in the standard way.
According to Proposition \ref{prop;17.9.29.101},
there exists the following relation:
\[
 \bigl[\del_{E^{\lambda},h,\beta_1},\del_{E^{\lambda},\betabar_1}\bigr]
=\frac{1}{1+|\lambda|^2}F(\nabla)_{w,\wbar}.
\]
We also note that $\beta_1=(1+|\lambda|^2)w+2\lambda\sqrt{-1}t_1$.
Hence, we obtain the following estimate
as $|\beta_1|\to\infty$ where $t_1$ varies in a compact set:
\begin{equation}
\label{eq;21.8.5.10}
 \Bigl|
 \bigl[
\del_{E^{\lambda},h,\beta_1},\del_{E^{\lambda},\betabar_1}
 \bigr]
 \Bigr|_h
 =O\Bigl(|\beta_1|^{-2}(\log|\beta_1|)^{-2}
 \Bigr).
\end{equation}
This means the acceptability mentioned
in Remark \ref{rem;21.8.4.3},
which is useful in the proof of Proposition {\rm\ref{prop;21.8.4.1}}
and Theorem {\rm\ref{thm;21.8.4.2}}.

\subsection{The filtered extension in the case $\lambda=0$}
\label{subsection;21.8.29.10}

Let us explain an outline of the proof of
Proposition \ref{prop;21.8.5.11}
and Theorem \ref{thm;21.8.5.12} in the case $\lambda=0$,
which are easier than
the claims for general $\lambda$
(Proposition \ref{prop;21.8.4.1} and
Theorem \ref{thm;21.8.4.2}).
We discuss it in the ramified case
under the setting in \S\ref{subsection;21.8.29.1}.

Let $U_{w,q}$ be a neighbourhood of $\infty$ in $\proj^1_{w_q}$,
and $U_{w,q}^{\ast}=U_{w,q}\setminus\{\infty\}$.
\index{space $U_{w,q}$}
\index{space $U_{w,q}^{\ast}$}
We assume that $U_{w,q}^{\ast}\subset U_{w,q}^{\ast}(R_1)$.
We have the natural partial compactification
$\nbigb^0_{q}=S^1_T\times U_{w,q}$
of the mini-complex manifold
$\nbigb^{0\ast}_q=S^1_T\times U_{w,q}^{\ast}$.
\index{mini-complex manifold $\nbigb^0_q$}
Let $H^0_{\infty,q}:=\nbigb^0_q\setminus\nbigb^{0\ast}_q
=S^1_T\times\{\infty\}$.
\index{space $H^0_{\infty,q}$}
Let $\pi_q^0:\nbigb^0_q\lrarr S^1_T$ denote the projection
induced by $(t,w_q)\longmapsto t$.
\index{projection $\pi_q^0$}
For $t\in S^1_T$, we set
$\nbigb^0_q\langle t\rangle=
(\pi_q^0)^{-1}(t)$
and 
$\nbigb^{0\ast}_q\langle t\rangle=
\nbigb^{0\ast}\cap \nbigb^0_q\langle t\rangle$.
\index{space $\nbigb^0_q\langle t\rangle$}
\index{space $\nbigb^{0\ast}_q\langle t\rangle$}

We define the $\nbigo_{\nbigb^0_q}(\ast H^0_{\infty,q})$-module
$\nbigp^hE^0$
by the procedure in \S\ref{subsection;21.8.4.4}
from $E^0$ with $h$.
It is easy to prove that $\nbigp^hE^0$
is a locally free $\nbigo_{\nbigb^0_q}(\ast H^0_{\infty,q})$-module
by the acceptability (\ref{eq;21.8.5.10})
as in Proposition \ref{prop;21.8.5.11}.
Moreover, we obtain a tuple of
filtered bundles
$\{\nbigp^h_{\ast}(E^0_{|\nbigb^{0\ast}_q\langle t\rangle})\,
|\,t\in S^1_T\}$.
To obtain Theorem \ref{thm;21.8.5.12},
we need to prove that the tuple is good.

By using the estimates for the asymptotic harmonic bundles,
we obtain the filtered bundle
$\nbigp_{\ast}V_{\ell,\alpha}$
on $(U_{w,q},\infty)$
from the holomorphic bundle
$(V_{\ell,\alpha},\delbar_{V_{\ell,\alpha}})$
with the metric $h_{V,\ell,\alpha}$.
(See \S\ref{section;21.7.7.20}.)
Note that the projection
$\Psi_q:\nbigb_q^{0\ast}\lrarr U_{w,q}^{\ast}$
given by $(t,w_q)\longmapsto w_q$
naturally extends to
$\Psi_q^0:\nbigb_q^{0}\lrarr U_{w,q}$.
\index{projection $\Psi_q^0$}
We can naturally generalize the construction
in \S\ref{subsection;21.8.29.2}
to the construction of
filtered bundle on $(\nbigb^0_{q},H^0_{\infty,q})$
from a filtered bundle with endomorphism on $(U_{w,q},\infty)$.
Hence, 
we obtain a filtered extension
$(\Psi^0_q)^{\ast}(\nbigp_{\ast}V_{\ell,\alpha},g_{\ell,\alpha})$
of $\Psi_q^{\ast}(V_{\ell,\alpha},g_{\ell,\alpha})$,
which is good
(see Proposition \ref{prop;20.7.24.10}
which goes back to Proposition \ref{prop;20.7.20.131}).
\index{filtered bundle
$(\Psi^0_q)^{\ast}(\nbigp_{\ast}V_{\ell,\alpha},g_{\ell,\alpha})$}
It is easy to see that
$\Psi_q^{\ast}(\nbigp_{\ast}V_{\ell,\alpha},g_{\ell,\alpha})$
is equal to
the filtered extension of
$\Psi_q^{\ast}(V_{\ell,\alpha},g_{\ell,\alpha})$
with respect to $\Psi_q^{-1}(h_{V,\ell,\alpha})$.

By explicit computations in \S\ref{section;17.10.2.20},
we obtain that
$\bigl\{
\nbigp_{\ast}\bigl(
\LL^{0\ast}_q(\ell,\alpha)_{|\nbigb^{0\ast}_q\langle t\rangle}
\bigr)\,\big|\,t\in S^1_T
\bigr\}$
is good.
Because
\[
 \nbigp^h_{\ast}(E^0_{|\nbigb^{0\ast}_q\langle t\rangle})
=\nbigp_{\ast}\bigl(
\LL^{0\ast}_q(\ell,\alpha)_{|\nbigb^{0\ast}_q\langle t\rangle}
\bigr)
\otimes
\Psi_q^{\ast}(\nbigp_{\ast}V_{\ell,\alpha},g_{\ell,\alpha})
_{|\nbigb^{0\ast}_q\langle t\rangle}
\quad(t\in S^1_T)
\]
we obtain that 
$\{\nbigp^h_{\ast}(E^0_{|\nbigb^{0\ast}_q\langle t\rangle})
\,|\,t\in S^1_T\}$
is good.
This is an outline of the proof of Theorem \ref{thm;21.8.5.12}
in the case $\lambda=0$.

\subsection{The filtered extension for general $\lambda$}
\label{subsection;21.9.16.2}

\subsubsection{Ramified covering space}

We continue to use the notation in \S\ref{subsection;21.8.29.10}.
There exists the map
$\nbigmbar^{\lambda}\lrarr \proj^1_w$
induced by
$(t_1,\beta_1)\longmapsto
(1+|\lambda|^2)^{-1}(\beta_1-2\sqrt{-1}t_1)$.
Let $\nbigb^{\lambda}_q$ denote the fiber product of
$\nbigmbar^{\lambda}$ and $U_{w,q}$ over $\proj^1_w$.
\index{space $\nbigb^{\lambda}_q$}
Let $\Psi_q^{\lambda}:\nbigb^{\lambda}_q\lrarr U_{w,q}$
denote the projection.
\index{projection $\Psi_q^{\lambda}$}
We set $H^{\lambda}_{\infty,q}:=(\Psi_q^{\lambda})^{-1}(\infty)$
and $\nbigb^{\lambda\ast}_q:=(\Psi_q^{\lambda})^{-1}(U_{w,q}^{\ast})$.
\index{space $H^{\lambda}_{\infty,q}$}
\index{space $\nbigb^{\lambda\ast}_q$}
Because the induced map
$\nbigb^{\lambda\ast}_q\lrarr\nbigm^{\lambda}$
is a local diffeomorphism,
$\nbigb^{\lambda\ast}_q$ inherits
the locally Euclidean metric and the mini-complex structure.
Note that
$\nbigb^{\lambda\ast}_q=\nbigb^{0\ast}_q$
as Riemannian manifolds
because $\nbigm^{\lambda}=\nbigm^0$ as Riemannian manifolds.
We can observe that
the mini-complex structure of
$\nbigb^{\lambda\ast}_q$ uniquely extends
to a mini-complex structure of
$\nbigb^{\lambda}_q$.
(See \S\ref{subsection;20.7.22.100}.)
We also have the naturally defined map
$\pi_q^{\lambda}:\nbigb^{\lambda}_q\lrarr S^1_T$
which is induced by the projection
$(t_1,\beta_1)\longmapsto t_1$
where $(t_1,\beta_1)$ denotes the mini-complex coordinate system
in \S\ref{subsection;21.8.5.1}.
\index{projection $\pi_q^{\lambda}$}
We set
$\nbigb^{\lambda}_q\langle t_1\rangle
=(\pi_q^{\lambda})^{-1}(t_1)$
and
$\nbigb^{\lambda\ast}_q\langle t_1\rangle
=\nbigb^{\lambda}_q\langle t_1\rangle
\cap
\nbigb^{\lambda\ast}_q$.
\index{space $\nbigb^{\lambda}_q\langle t_1\rangle$}
\index{space $\nbigb^{\lambda\ast}_q\langle t_1\rangle$}

We obtain the mini-holomorphic bundle
$(E^{\lambda},\delbar_{E^{\lambda}})$
on $\nbigb^{\lambda\ast}_q$
as in \S\ref{subsection;21.8.29.1}.
We obtain
the $\nbigo_{\nbigb^{\lambda}_q}(\ast H^{\lambda}_{\infty,q})$-module
$\nbigp^h(E^{\lambda})$
from $(E^{\lambda},\delbar_{E^{\lambda}},h)$
by the procedure in \S\ref{subsection;17.10.28.40}.
\index{sheaf $\nbigp^h(E^{\lambda})$}
We can prove that it is a locally free
$\nbigo_{\nbigb^{\lambda}_q}(\ast H^{\lambda}_{\infty,q})$-module
by the acceptability (\ref{eq;21.8.5.10}),
as in Proposition \ref{prop;21.8.4.1}.
Moreover, we obtain the tuple of filtered bundles
$\bigl\{
 \nbigp^h_{\ast}\bigl(
 E^{\lambda}_{|\nbigb^{\lambda\ast}_q\langle t_1\rangle}
 \bigr)\,\big|\,t_1\in S^1_T
\bigr\}$.
To obtain Theorem \ref{thm;21.8.4.2},
we need to prove that the tuple is good.

\subsubsection{Approximation}

As an analogue of the Chern connection
for holomorphic vector bundles with a Hermitian metric,
we construct
$\del^{\shikaku}_{E^0}$ and $\phi^{\shikaku}$
from $(E^0,\delbar_{E^0})$ with the metric $h^{\shikaku}$.
By Theorem \ref{thm;17.10.28.50},
we obtain
$\nabla-(\delbar_{E^0}+\del^{\shikaku}_{E^0})=
O\bigl(\exp(-\epsilon|w_q^q|)\bigr)$
and
$\phi-\phi^{\shikaku}=O\bigl(\exp(-\epsilon|w_q^q|)\bigr)$
for some $\epsilon>0$.
We construct the differential operator
$\delbar^{\shikaku}_{E^{\lambda}}$ of $E$
by using the connection
$\delbar_{E^0}+\del_{E^0}^{\shikaku}$ and $\phi^{\shikaku}$
as in \S\ref{subsection;17.10.14.1}.
Then, $\delbar^{\shikaku}_{E^{\lambda}}$ is asymptotically 
a mini-holomorphic structure,
i.e.,
$\delbar^{\shikaku}_{E^{\lambda}}
\circ\delbar^{\shikaku}_{E^{\lambda}}
=O\bigl(\exp(-\epsilon|w_q|^q)\bigr)$.
We also have
$\delbar_{E^{\lambda}}-\delbar^{\shikaku}_{E^{\lambda}}
=O\bigl(\exp(-\epsilon|w_q|^q)\bigr)$.
Let $\nbigc^{\infty}_{\nbigb^{\lambda}_q}$ denote
the sheaf of $C^{\infty}$-functions on $\nbigb^{\lambda}_q$.
\index{sheaf $\nbigc^{\infty}_{\nbigb^{\lambda}_q}$}
It turns out that
$\delbar^{\shikaku}_{E^{\lambda}}$
induces a $C^{\infty}$-differential operator
on $\nbigc^{\infty}_{\nbigb^{\lambda}_q}
\otimes_{\nbigo_{\nbigb^{\lambda}_q}}
\nbigp^h(E^{\lambda})$,
and $\delbar^{\shikaku}_{E^{\lambda}}$
and $\delbar_{E^{\lambda}}$
induce the same operator
on the formal completion of $\nbigp^h(E^{\lambda})$
along $H^{\lambda}_{\infty,q}$.

\subsubsection{Formal completion of asymptotic harmonic bundles at infinity}
\label{subsection;21.8.5.10}

Let $(V_{\ell,\alpha},\delbar_{V_{\ell,\alpha}},
\theta_{\ell,\alpha},h_{V,\ell,\alpha})$
be the asymptotic harmonic bundle in \S\ref{subsection;21.8.5.2}.
Let $\delbar_{V_{\ell,\alpha}}+\del_{V_{\ell,\alpha}}$
denote the Chern connection
associated with $h_{V,\ell,\alpha}$.
We obtain the holomorphic vector bundle
$V^{\lambda}_{\ell,\alpha}=
(V_{\ell,\alpha},\delbar_{V_{\ell,\alpha}}+\lambda\theta^{\dagger}_{\ell,\alpha})$.
It turns out that
$(V^{\lambda}_{\ell,\alpha},h_{V,\ell,\alpha})$ is acceptable,
and hence we obtain the filtered extension
$\nbigp_{\ast}V^{\lambda}_{\ell,\alpha}$
on $(U_{w,q},\infty)$.
We have the $\lambda$-connection
$\DDlambda_{\ell,\alpha}=
\delbar_{V_{\ell,\alpha}}+\lambda\theta^{\dagger}_{\ell,\alpha}
+\lambda\del_{V_{\ell,\alpha}}
+\theta_{\ell,\alpha}$ of $V_{\ell,\alpha}$.
Though it is not necessarily flat,
we have
$\DDlambda_{\ell,\alpha}\circ\DDlambda_{\ell,\alpha}
=O\bigl(
\exp(-\epsilon|w_q^q|)\bigr)$.
It turns out that
as the formal completion at infinity,
we obtain a formal good filtered $\lambda$-flat bundle
$(\nbigp_{\ast}\Vhat^{\lambda}_{\ell,\alpha},
\DDhat^{\lambda}_{\ell,\alpha})$
as in the case of wild harmonic bundles.
(See \S\ref{subsection;17.12.16.4}--\S\ref{subsection;20.7.27.3}.)
\index{formal good filtered $\lambda$-flat bundle
$(\nbigp_{\ast}\Vhat^{\lambda}_{\ell,\alpha},
\DDhat^{\lambda}_{\ell,\alpha})$}

\subsubsection{The formal structure of
$\nbigp^hE^{\lambda}$ at infinity}

We can consider the formal space
$\Hhat^{\lambda}_{\infty,q}$
as the formal completion of $\nbigb^{\lambda}_q$
along $H^{\lambda}_{\infty,q}$.
\index{ringed space $\Hhat^{\lambda}_{\infty,q}$}
(See \S\ref{subsection;20.7.30.22}.)
We set
$H^{\lambda}_{\infty,q}\langle t_1\rangle
=\nbigb^{\lambda}_{q}\cap H^{\lambda}_{\infty,q}$,
and let $\Hhat^{\lambda}_{\infty,q}\langle t_1\rangle$
denote the formal completion of
$\nbigb^{\lambda}_q\langle t_1\rangle$
along $H^{\lambda}_{\infty,q}\langle t_1\rangle$.
\index{set $H^{\lambda}_{\infty,q}\langle t_1\rangle$}
\index{ringed space $\Hhat^{\lambda}_{\infty,q}\langle t_1\rangle$}

We obtain the
$\nbigo_{\Hhat^{\lambda}_{\infty,q}}(\ast H^{\lambda}_{\infty,q})$-module
$\nbigp(E^{\lambda})_{|\Hhat^{\lambda}_{\infty,q}}$
induced by $\nbigp(E^{\lambda})$.
We also have the tuple of filtered bundles
$\nbigp_{\ast}(E^{\lambda})
 _{|\Hhat^{\lambda}_{\infty,q}\langle t_1\rangle}$
$(t_1\in S^1_T)$.

We have the mini-holomorphic bundles
$\LL^{\lambda\ast}_q(\ell,\alpha)$ on $\nbigb^{\lambda\ast}_q$
underlying the monopoles (\ref{eq;21.8.5.3}).
\index{mini-holomorphic bundle
$\LL^{\lambda\ast}_q(\ell,\alpha)$}
As their filtered extension,
we obtain 
$\nbigo_{\nbigb^{\lambda}_{q}}(\ast H^{\lambda}_{\infty,q})$-modules
$\LL^{\lambda}_q(\ell,\alpha)$
and the tuples of filtered bundles
$\bigl\{
 \nbigp_{\ast}
  \LL^{\lambda}_q(\ell,\alpha)
   _{|\nbigb^{\lambda}_q\langle t_1\rangle}
\,\big|\,t_1\in S^1_T
\bigr\}$.
\index{$\nbigo_{\nbigb^{\lambda}_{q}}(\ast H^{\lambda}_{\infty,q})$-module
$\LL^{\lambda}_q(\ell,\alpha)$}
\index{filtered bundle $\nbigp_{\ast}
  \LL^{\lambda}_q(\ell,\alpha)
   _{|\nbigb^{\lambda}_q\langle t_1\rangle}$}
As the formal completion along $H^{\lambda}_{\infty,q}$,
we obtain 
$\nbigo_{\Hhat^{\lambda}_{\infty,q}}(\ast H^{\lambda}_{\infty,q})$-modules
$\LL^{\lambda}_q(\ell,\alpha)_{|\Hhat^{\lambda}_{\infty,q}}$
and the good tuple of the filtered bundles
$\bigl\{
 \nbigp_{\ast}\LL^{\lambda}_q(\ell,\alpha)
  _{|\Hhat^{\lambda}_{\infty,q}\langle t_1\rangle}
\,\big|\,
t_1\in S^1_T
\bigr\}$.

As we shall explain in \S\ref{subsection;17.10.28.20},
from the formal good filtered $\lambda$-flat bundles
$(\nbigp_{\ast}\Vhat_{\ell,\alpha},\DDhat^{\lambda}_{\ell,\alpha})$,
we construct good filtered 
$\nbigo_{\Hhat^{\lambda}_{\infty,q}}(\ast H^{\lambda}_{\infty,q})$-modules
$(\Psi^{\lambda}_q)^{\ast}(
\nbigp_{\ast}\Vhat_{\ell,\alpha},\DDhat^{\lambda}_{\ell,\alpha})$.
(See Proposition \ref{prop;20.7.20.131}.)
\index{filtered bundle $(\Psi^{\lambda}_q)^{\ast}(
\nbigp_{\ast}\Vhat_{\ell,\alpha},\DDhat^{\lambda}_{\ell,\alpha})$}

By the approximation in \S\ref{subsection;21.8.5.10},
and by the compatibility as in
Proposition \ref{prop;21.8.12.5}
and Lemma \ref{lem;21.8.12.60},
there exists an isomorphism of
$\nbigo_{\Hhat^{\lambda}_{\infty,q}}(\ast H^{\lambda}_{\infty,q})$-modules
\[
 \nbigp(E^{\lambda})_{|\Hhat^{\lambda}_{\infty,q}}
 \simeq
 \bigoplus
 \LL^{\lambda}_q(\ell,\alpha)_{|\Hhat^{\lambda}_{\infty,q}}
 \otimes
 (\Psi^{\lambda}_q)^{\ast}(\Vhat_{\ell,\alpha},
 \DDhat^{\lambda}_{\ell,\alpha}),
\]
which induces an isomorphism of the tuples of filtered bundles
\[
 \nbigp(E^{\lambda})_{|\Hhat^{\lambda}_{\infty,q}\langle t_1\rangle}
 \simeq
 \bigoplus
  \nbigp_{\ast}\LL^{\lambda}_q(\ell,\alpha)
   _{|\Hhat^{\lambda}_{\infty,q}\langle t_1\rangle}
  \otimes
  (\Psi^{\lambda}_q)^{\ast}(
  \nbigp_{\ast}\Vhat_{\ell,\alpha},\DDhat^{\lambda}_{\ell,\alpha})
  _{|\Hhat^{\lambda}_{\infty,q}\langle t_1\rangle}.
\]
Thus, we obtain that
the tuple
$\bigl\{
\nbigp_{\ast}(E^{\lambda})_{|\Hhat^{\lambda}_{\infty,q}\langle t_1\rangle}
\,\big|\,t_1\in S^1_T
\bigr\}$
is good.
This is an outline of our proof of Theorem \ref{thm;21.8.4.2},

\section{Acknowledgement}

I thank Hiraku Nakajima 
whose lectures attracted me to the study of monopoles.
A part of this study was done during my stay in 
the University of Melbourne,
and I am grateful to Kari Vilonen and Ting Xue
for their excellent hospitality and their support.
I am grateful to Carlos Simpson 
whose works provide the most important foundation
with this study.
I am directly influenced by the interesting works of 
Benoit Charbonneau and Jacques Hurtubise
\cite{Charbonneau-Hurtubise},
and Sergey Cherkis and Anton Kapustin
\cite{Cherkis-Kapustin1, Cherkis-Kapustin2}.
I was inspired by 
a talk of Hurtubise in Tata Institute of Fundamental Research.
I thank Maxim Kontsevich  and Yan Soibelman
for their comments and discussions.
I am grateful to Masaki Yoshino for our discussions.
I am grateful to Claude Sabbah
for his kindness and discussions on many occasions.
I thank Yoshifumi Tsuchimoto and Akira Ishii
for their constant encouragement.
I am heartily grateful to the reviewers
for their careful and patient readings and for their constructive
comments to improve this manuscript.

My interest in ``Kobayashi-Hitchin correspondence''
was renewed when I made a preparation for
a talk in the 16th Oka Symposium,
which drove me to this study. I thank the organizers,
particularly Junichi Matsuzawa and Ken-ichi Yoshikawa.

I am partially supported by
the Grant-in-Aid for Scientific Research (S) (No. 17H06127),
the Grant-in-Aid for Scientific Research (S) (No. 16H06335),
the Grant-in-Aid for Scientific Research (A) (No. 21H04429),
the Grant-in-Aid for Scientific Research (C) (No. 15K04843),
and the Grant-in-Aid for Scientific Research (C) (No. 20K03609),
Japan Society for the Promotion of Science.

\chapter{Preliminaries}

In \S\ref{subsection;17.10.28.1},
we introduce the notion of mini-complex structure
of $3$-dimensional manifolds.
In \S\ref{subsection;17.10.28.42},
we introduce the notion of mini-holomorphic vector bundles
on mini-complex $3$-dimensional manifolds.
The notions of mini-complex structure
and mini-holomorphic bundle
have been implicitly used in the study of monopoles.
In \S\ref{subsection;20.7.30.30},
we recall the notion of monopole
as mini-holomorphic bundles equipped with
a Hermitian metric satisfying (\ref{eq;20.7.30.31}).
We also recall the notion of Dirac type singularity of monopoles
\cite{Kronheimer-Master-Thesis},
and a characterization which easily follows from
\cite{Mochizuki-Yoshino}.

In \S\ref{subsection;20.7.30.41},
we recall the dimensional reduction of instantons
to monopoles.
We also explain the underlying dimensional reduction
of holomorphic bundles to mini-holomorphic bundles.
In \S\ref{subsection;20.7.30.42},
we explain the dimensional reduction of monopoles
to harmonic bundles.
We also explain the underlying dimensional reduction
of mini-holomorphic bundles to Higgs bundles
in \S\ref{subsection;17.9.29.2}.

In \S\ref{subsection;17.10.2.1},
we study the twistor family of mini-complex structures
on $\nbigm=(\real/T\seisuu)\times \cnum$ for $T>0$.
For $\lambda\in\cnum$,
we obtain the mini-complex manifold $\nbigm^{\lambda}$
as $\nbigm$
equipped with a mini-complex structure corresponding to
$\lambda$.
We introduce two convenient mini-complex local coordinate systems
$(t_0,\beta_0)$ and $(t_1,\beta_1)$ for $\nbigm^{\lambda}$.
We introduce a compactification $\nbigmbar^{\lambda}$
of $\nbigm^{\lambda}$
by using the coordinate system $(t_1,\beta_1)$.
Set $H^{\lambda}_{\infty}:=\nbigmbar^{\lambda}\setminus \nbigmlambda$.

In \S\ref{section;21.8.12.101},
we study the dimensional reduction from
$\nbigo_{\nbigm^{\lambda}}$-modules
to $\lambda$-flat bundles on $\cnum_w$ for general $\lambda$,
and its compatibility with the dimensional reduction from
monopoles to harmonic bundles.
The case of $\lambda=0$ is already essentially explained
in \S\ref{subsection;20.7.30.42}.
The case of $\lambda\neq 0$ is also useful for our study.
We also explain the dimensional reduction from
$\nbigo_{\nbigmbar^{\lambda}}$-modules
to an $\nbigo_{\proj^1_w}$-modules equipped with
a meromorphic $\lambda$-connection.

In \S\ref{subsection;20.7.30.60},
we shall introduce a section $G(h)$
of $\End(E)$
for a mini-holomorphic bundle $(E,\delbar_E)$
with a Hermitian metric $h$ on
an open subset of $\nbigm^{\lambda}$.
It is an analogue of the contraction of
the curvature of the Chern connection
for a holomorphic vector bundle with a Hermitian metric
on a K\"ahler manifold.
In \S\ref{subsection;20.7.30.61}--\S\ref{subsection;20.7.30.62},
we explain some general formulas for $G(h)$
which are analogues of standard and useful formulas
in the study of Hermitian-Einstein metrics \cite{Simpson88}.
We give a complement to \S\ref{section;21.8.12.101}
in \S\ref{subsection;17.10.12.1},
which is mainly a preliminary
for the proof of Proposition \ref{prop;17.10.12.6}.
In \S\ref{subsection;17.10.24.32},
we recall the dimensional reduction of instantons
with respect to the Hopf fibration due to Kronheimer,
and prove a formula which will be useful
in the proof of Proposition \ref{prop;17.9.10.10}.

In \S\ref{subsection;17.10.28.100},
we explain an equivalence between
difference modules with parabolic structure at finite place
and locally free
$\nbigo_{\nbigmlambda\setminus Z}(\ast H_{\lambda}^{\infty})$-modules
for finite subsets $Z$.

In \S\ref{subsection;17.10.5.122},
we review filtered bundles on punctured curves,
and we recall the general construction
of a filtered bundle 
from a holomorphic bundle with a Hermitian metric
satisfying the acceptability condition.

\section{Mini-complex structures on $3$-manifolds}
\label{subsection;17.10.28.1}

\subsection{Mini-holomorphic functions on
$\real\times\cnum$}

Let $t$ and $w$ be the standard coordinates of $\real$ and $\cnum$,
respectively.
The orientation of $\real\times\cnum$ is given
as the product of the orientation of $\real$ and $\cnum$.
Let $U$ be an open subset in $\real\times\cnum$.
A $C^{\infty}$-function $f$ on $U$ is called
mini-holomorphic if $\del_tf=0$ and $\del_{\wbar}f=0$.
\index{mini-holomorphic function}
\index{mini-holomorphic diffeomorphism}

Let $U_i$ $(i=1,2)$ be open subsets in $\real\times\cnum$.
Let $F:U_1\lrarr U_2$ be a diffeomorphism.
It is called mini-holomorphic if 
(i) $F$ preserves the orientations,
(ii) $F^{\ast}(f)$ is mini-holomorphic 
for any mini-holomorphic function $f$ on $U_2$.

\begin{lem}
Let $F=(F_t(t,w),F_w(t,w)):U_1\lrarr U_2$
be a diffeomorphism.
Then, $F$ is mini-holomorphic if and only if
we have $\del_tF_t>0$,
$\del_tF_w=0$ and $\del_{\wbar}F_w=0$.
\end{lem}
\pf
Suppose that $F$ is mini-holomorphic.
Because $F^{\ast}(w)=F_w$ is mini-holomorphic,
we obtain $\del_tF_w=\del_{\wbar}F_w=0$.
Because $F$ is orientation preserving,
we obtain  $\del_t F_t\cdot|\del_{w}F_w|^2>0$,
which implies $\del_tF_t>0$.
The converse is also easily proved.
\hfill\qed

\subsection{Mini-complex structure on $3$-dimensional manifolds}

Let us define the notion of mini-complex structure
for $3$-dimensional manifolds
as in the case of smooth structure.
(For example, see \cite{Lee-book}).
\index{mini-complex structure}
Let $M$ be an oriented $3$-dimensional $C^{\infty}$-manifold.
A mini-complex atlas on $M$ 
is a family of open subsets
$U_{\lambda}$ $(\lambda\in\Lambda)$
with an orientation-preserving embedding
$\varphi_{\lambda}:
 U_{\lambda}\lrarr \real\times\cnum$
satisfying the following conditions.
\index{mini-complex atlas}
\begin{itemize}
\item
 $M=\bigcup_{\lambda\in\Lambda}U_{\lambda}$.
\item
 The coordinate change
 $\varphi_{\lambda}(U_{\lambda}\cap U_{\mu})
\lrarr
 \varphi_{\mu}(U_{\lambda}\cap U_{\mu})$
is mini-holomorphic.
\end{itemize}
Such $(U_{\lambda},\varphi_{\lambda})$
is called a mini-complex chart.
\index{mini-complex chart}
We shall often use $(t,w)$ instead of $\varphi_{\lambda}$.
Two mini-complex atlas
$\{(U_{\lambda},\varphi_{\lambda})\,|\,\lambda\in\Lambda\}$
and 
$\{(V_{\gamma},\psi_{\gamma})\,|\,\gamma\in \Gamma\}$
are defined to be equivalent
if the union is also a mini-complex atlas.
There exists the partial order
on the family of mini-complex atlases
defined by inclusions.
Each equivalence class of mini-complex atlases
has a unique maximal mini-complex atlas.
A mini-complex structure on $M$
is defined to be a maximal mini-complex atlas on $M$.

When $M$ is equipped with a mini-complex structure,
a $C^{\infty}$-function $f$ on $M$ is called mini-holomorphic
if its restriction to any mini-complex chart
is mini-holomorphic.
Let $\nbigo_M$ denote the sheaf of mini-holomorphic functions.

\index{mini-holomorphic function}
\index{sheaf $\nbigo_M$}

\subsection{Tangent bundles}

Suppose that $M$ is equipped with a mini-complex structure.
Let $(U;t,w)$ be a mini-complex chart.
The real vector field $\del_t$ determines
an oriented subbundle $T_SU$ of $TU$.
The quotient bundle $T_QU$ is equipped
with the complex structure $J$,
where $J$ is an automorphism of $T_QU$
such that $J^2=-1$.
Because the mini-complex coordinate change
preserves the subbundles,
we obtain a globally defined subbundle
$T_SM$ of $TM$ of rank $1$.
We also obtain the quotient bundle $T_QM$,
which is equipped with the complex structure.
\index{$T_SM$}
\index{$T_QM$}

Let us consider 
$T^{\cnum}M:=TM\otimes_{\real}\cnum$. \index{$T^{\cnum}M$}
There exists the decomposition of complex vector bundles
$T^{\cnum}_QM
=T_Q^{1,0}M\oplus T_Q^{0,1}M$,
where
$T_Q^{1,0}M$ is the $\sqrt{-1}$-eigen bundle of $J$,
and $T_Q^{0,1}M$ is the $-\sqrt{-1}$-eigen bundle of $J$.
\index{$T_Q^{\cnum}M$, $T_Q^{1,0}M$, $T_Q^{0,1}M$} 
We obtain the exact sequence of complex vector bundles on $M$:
\[
\begin{CD}
 0 @>>>
 T^{\cnum}_SM
 @>>>
 T^{\cnum}M
 @>{a_1}>>
 T_Q^{1,0}M\oplus T_Q^{0,1}M
 @>>> 0.
\end{CD}
\] 
Let $\Theta_M^{1,0}$ and $\Theta^{0,1}_M$
denote the inverse image of 
$T_Q^{1,0}M$ and $T_Q^{0,1}M$ by $a_1$,
respectively. \index{$\Theta_M^{1,0}$, $\Theta_M^{0,1}$}
Note that
$\Theta_M^{1,0}\cap \Theta_M^{0,1}
=T^{\cnum}_SM$.
The complex conjugate on $\cnum$
induces the complex conjugate on $TM\otimes\cnum$
for which we obtain
$\overline{\Theta^{1,0}_M}=\Theta^{0,1}_M$.
It also induces the complex conjugate on
$T_Q^{1,0}M\oplus T_Q^{0,1}M$
for which
$\overline{T_Q^{1,0}M}=T_Q^{0,1}M$ holds.

\begin{lem}
\label{lem;17.10.3.2}
Let $F:M_1\lrarr M_2$ be a diffeomorphism
of oriented $3$-dimensional manifolds.
Suppose that $M_i$ are equipped with
 mini-complex structure,
and that the following holds:
\begin{itemize}
\item
 $F$ is orientation preserving.
\item
 $dF(T_SM_1)=T_SM_2$.
\item
 The induced isomorphism
 $T_QM_1\simeq T_QM_2$
 is $\cnum$-linear.
\end{itemize}
Then, $F$ preserves
the mini-complex structures.
\end{lem}
\pf
It is enough to study the claim locally around any point of $M_i$.
Let $U_i$ $(i=1,2)$ be open subsets of $\real\times\cnum$
which are naturally equipped with mini-complex structures.
Let $F:U_1\lrarr U_2$ be an orientation preserving diffeomorphism
such that (i) $dF(T_SU_1)=T_SU_2$,
(ii) the induced isomorphism
$T_QU_1\simeq T_QU_2$ is $\cnum$-linear.
We have the expression
$F=(F_t(t,w),F_w(t,w))$.
By the condition (i),
we obtain $\del_tF_w(t,w)=0$.
By the condition (ii),
we obtain $\del_{\wbar}F_w(t,w)=0$.
Hence, $F$ is mini-holomorphic.
\hfill\qed

\subsection{Cotangent bundles}

Let $\Omega_Q^{1,0}$ and $\Omega_Q^{0,1}$
denote the complex dual bundles of 
$T_Q^{1,0}M$ and $T_Q^{0,1}M$,
respectively. \index{$\Omega_Q^{1,0}$, $\Omega_Q^{0,1}$}
Let $\Omega_M^{1,0}$ and $\Omega_M^{0,1}$
denote the complex dual bundle of
$\Theta^{1,0}_M$ and $\Theta_M^{0,1}$.
\index{$\Omega_M^{1,0}$, $\Omega_M^{0,1}$}
We obtain the following exact sequences:
\[
 0\lrarr
 \Omega_Q^{1,0}
 \oplus
 \Omega_Q^{0,1}
\lrarr
 T^{\ast}M\otimes\cnum
\lrarr
  (T^{\cnum}_SM)^{\lor}
\lrarr 0,
\]
\[
 0\lrarr
 \Omega_Q^{1,0}
\lrarr
 \Omega_M^{1,0}
\lrarr
  (T^{\cnum}_SM)^{\lor}
\lrarr 0,
\quad\quad
 0\lrarr
 \Omega_Q^{0,1}
\lrarr
 \Omega_M^{0,1}
\lrarr
  (T^{\cnum}_SM)^{\lor}
\lrarr 0,
\]
\[
 0\lrarr
 \Omega_Q^{1,0}
\lrarr
 T^{\ast}M\otimes\cnum
\lrarr
 \Omega_M^{0,1}
\lrarr 0,
\quad\quad
 0\lrarr
 \Omega_Q^{0,1}
\lrarr
 T^{\ast}M\otimes\cnum
\lrarr
 \Omega_M^{1,0}
\lrarr 0.
\]
Note 
$\rank_{\cnum}\Omega^{0,1}_M
=\rank_{\cnum}\Omega^{1,0}_M=2$.
We set
$\Omega^{0,2}_M:=\bigwedge^2\Omega^{0,1}_M$
and 
$\Omega^{2,0}_M:=\bigwedge^2\Omega^{1,0}_M$.

For any $C^{\infty}$-vector bundle $E$ on $M$,
let $C^{\infty}(M,E)$ denote the space of
$C^{\infty}$-sections of $E$.
\index{$C^{\infty}(M,E)$}

The exterior derivative
$d:C^{\infty}(M,\cnum)\lrarr C^{\infty}(M,T^{\ast}M\otimes\cnum)$
induces the following differential operators:
\[
 \delbar_M:
 C^{\infty}(M,\cnum)
\lrarr
 C^{\infty}(M,\Omega^{0,1}_M),
\quad\quad
 \del_M:
 C^{\infty}(M,\cnum)
\lrarr
 C^{\infty}(M,\Omega^{1,0}_M).
\]
\index{operators $\delbar_M$, $\del_M$}
We also obtain the following operators induced by
$d:C^{\infty}(M,T^{\ast}M\otimes\cnum)
\lrarr
  C^{\infty}(M,\bigwedge^2T^{\ast}M\otimes\cnum)$:
\[
 \delbar_M:
 C^{\infty}(M,\Omega^{0,1}_M)
\lrarr
 C^{\infty}(M,\Omega^{0,2}_M),
\quad\quad
 \del_M:
 C^{\infty}(M,\Omega^{1,0}_M)
\lrarr
 C^{\infty}(M,\Omega^{2,0}_M).
\]
We have $\delbar_M\circ\delbar_M=0$
and $\del_M\circ\del_M=0$.
It is easy to see that a $C^{\infty}$-function $f$ on $M$ 
is mini-holomorphic
if and only if $\delbar_Mf=0$.

\subsection{Meromorphic functions}
\label{subsection;21.8.3.10}

Let $M$ be a $3$-dimensional manifold
with a mini-complex structure.
Let $H\subset M$ be a $1$-dimensional submanifold 
such that the tangent bundle $TH$ is contained in $(T_SM)_{|H}$.
Let $U$ be an open subset of $M$.
A holomorphic function $f$ on $U\setminus H$
is called meromorphic of at most order $N$ along $H$
if the following holds.
\begin{itemize}
\item
 Let $P$ be any point of $H\cap U$.
 We take a mini-complex coordinate neighbourhood
 $(U_P,t,w)$ around $P$.
 Note that $H$ is described as $\{w=w_0\}$.
 Then, $(w-w_0)^Nf_{|U_P}$ gives a mini-holomorphic 
 function on $U_P$.
\end{itemize}
Let $\nbigo_M(NH)$ denote the sheaf of 
meromorphic functions of order $N$ along $H$.
We obtain the sheaf
$\nbigo_M(\ast H):=
 \varinjlim_N \nbigo_M(NH)$.
A local section of $\nbigo_M(\ast H)$
is called a meromorphic function on $M$
whose poles are contained in $H$.
\index{meromorphic functions}
\index{sheaf $\nbigo_M(\ast H)$}

\section{Mini-holomorphic bundles}
\label{subsection;17.10.28.42}

\subsection{Mini-holomorphic bundles}
\label{subsection;17.9.29.200}

Let $M$ be a $3$-dimensional manifold
with a mini-complex structure.
Let $E$ be a $C^{\infty}$-bundle on $M$ of finite rank.
Let us consider a differential operator
\[
 \delbar_E:
 C^{\infty}(M,E)\lrarr
 C^{\infty}(M,E\otimes\Omega^{0,1}_M)
\]
satisfying
$\delbar_E(fs)=\delbar_M(f)\,s+f\delbar_E(s)$
for any $f\in C^{\infty}(M,\cnum)$
and $s\in C^{\infty}(M,E)$.
As in the ordinary case of vector bundles on complex manifolds,
it induces 
$\delbar_E:
 C^{\infty}(M,E\otimes\Omega^{0,1}_M)
\lrarr
 C^{\infty}(M,E\otimes\Omega^{0,2}_M)$.
Such a differential operator $\delbar_E$
is called a mini-holomorphic structure of $E$
if $\delbar_E\circ\delbar_E=0$.
\index{mini-holomorphic structure}
\index{mini-holomorphic vector bundle}
\index{operator $\delbar_E$}

In terms of a mini-complex coordinate system $(t,w)$,
a mini-holomorphic structure is equivalent
to a pair of differential operators
$\del_{E,t}$ and $\del_{E,\wbar}$
on $C^{\infty}(E)$
such that 
$\del_{E,t}(fs)=\del_t(f)s+f\del_{E,t}(s)$
and 
$\del_{E,\wbar}(fs)=\del_{\wbar}(f)s+f\del_{E,\wbar}(s)$
for $f\in C^{\infty}(M,\cnum)$
and $s\in C^{\infty}(M,E)$,
satisfying the commutativity condition
$[\del_{E,t},\del_{E,\wbar}]=0$.
Indeed, the operators $\del_{E,\wbar}$ and $\del_{E,t}$
are induced by the inner product of
$\del_{\wbar}$ and $\del_t$
with $\delbar_E(s)$ $(s\in C^{\infty}(U,E))$,
respectively.
\index{operators $\del_{E,t}$, $\del_{E,\wbar}$}

A local section $s$ of $E$ is called mini-holomorphic 
if $\delbar_E(s)=0$. \index{mini-holomorphic section}
By considering the sheaf of mini-holomorphic sections 
of $E$,
we obtain a locally free $\nbigo_U$-module,
which is often denoted by the same notation $E$.
The following is standard.
\begin{lem}
The above procedure induces 
an equivalence of mini-holomorphic bundles of finite rank
and locally free $\nbigo_M$-modules of finite rank.
\hfill\qed
\end{lem}

\subsection{Metrics and the induced operators}

Let $(E,\delbar_E)$ be a mini-holomorphic bundle on $M$.
Let $h$ be a $C^{\infty}$-metric of 
a mini-holomorphic bundle $E$.
As in the case of complex differential geometry
(see \cite{Kobayashi-vector-bundle}),
we obtain the induced differential operator
$\del_{E,h}:C^{\infty}(M,E)\lrarr
 C^{\infty}(M,E\otimes\Omega_M^{1,0})$
such that 
$\del_{E,h}(fs)=\del_M(f)s+f\del_{E,h}(s)$
for $f\in C^{\infty}(M,\cnum)$
and $s\in C^{\infty}(M,E)$,
determined by the condition
$\delbar_Mh(u,v)=
h(\delbar_Eu,v)+
 h(u,\del_{E,h}v)$
for any $u,v\in C^{\infty}(M,E)$.
We have $\del_{E,h}\circ\del_{E,h}=0$.
\index{operator $\del_{E,h}$}

\subsection{Splittings}

Let $T_QM\lrarr TM$ be a splitting
of the exact sequence
$0\lrarr T_SM\lrarr TM\lrarr T_QM\lrarr 0$.
We obtain the induced splitting
$(T_S^{\cnum}M)^{\lor}
\lrarr
 (T^{\cnum}M)^{\lor}$,
and the decomposition:
\[
 (T^{\cnum}M)^{\lor}
=(T^{\cnum}_SM)^{\lor}
\oplus
 \Omega_Q^{1,0}
\oplus
 \Omega_Q^{0,1}.
\]
We also obtain the decompositions:
\begin{equation}
\label{eq;20.7.15.1}
 \Omega_M^{1,0}
=\Omega_Q^{1,0}\oplus
 (T^{\cnum}_SM)^{\lor},
\quad\quad
 \Omega_M^{0,1}
=\Omega_Q^{0,1}\oplus
 (T^{\cnum}_SM)^{\lor}.
\end{equation}

Let $(E,\delbar_E)$ be a mini-holomorphic bundle
with a metric $h$.
We obtain the decompositions
$\delbar_E=\delbar_E^Q+\delbar_E^S$
and 
$\del_{E,h}=\del_{E,h}^Q+\del_{E,h}^S$
according to (\ref{eq;20.7.15.1}).
\index{Operator $\delbar_E^Q$, $\delbar_E^S$, $\del_{E,h}^Q$, $\del_{E,h}^S$}

We obtain the unitary connection
$\nabla_h$
and the anti-self-adjoint section
$\phi_h$
of $\End(E)\otimes(T_SM\otimes\cnum)^{\lor}$:
\[
\nabla_h:=\delbar^Q_{E}+\del^Q_{E,h}+
 \frac{1}{2}(\delbar_E^S+\del^S_{E,h}),
\quad\quad
 \phi_h=
 \frac{\sqrt{-1}}{2}
 (\delbar_E^S-\del^S_{E,h}).
\]
They are called the Chern connection
and the Higgs field.
Note that the construction of the Chern connection
and the Higgs field depend on the choice of 
the splitting $T_QM\lrarr TM$.
\index{Chern connection $\nabla_h$}
\index{Higgs field $\phi_h$}

\begin{rem}
If $M$ is equipped with a Riemannian metric,
the orthogonal complement of
$T_SM$ in $TM$ is naturally isomorphic to
$T_QM$.
Hence, we obtain a splitting,
which we shall use without mention.
Moreover,
by the Riemannian metric,
$T_SM$ is identified with
the product bundle $\real\times M$.
Hence, we regard $\phi_h$
as an anti-Hermitian section of $\End(E)$.
\hfill\qed
\end{rem}

\subsection{Scattering maps}

Let $(E,\delbar_E)$ be a mini-holomorphic bundle on $M$.
Let $\gamma:[0,1]\lrarr M$ be a path
such that $T\gamma(\del_s)$ is contained in $T_SM$.
Then, the mini-holomorphic structure of $E$
induces a connection of $\gamma^{\ast}E$.
Hence, we obtain the isomorphism
$E_{\gamma(0)}\lrarr E_{\gamma(1)}$
as the parallel transport of the connection,
called the scattering map \cite{Charbonneau-Hurtubise}.
\index{scattering map}

Let $U$ be an open subset of $M$
with a mini-complex coordinate $\varphi$
such that 
$\varphi(U)\supset\openopen{-2\epsilon}{2\epsilon}\times U_w$,
where $U_w$ is an open subset in $\cnum$.
Here, for $a<b$, we set $\openopen{a}{b}:=\{t\in\real\,|\,a<t<b\}$.
\index{interval \mbox{$\openopen{a}{b}$}}
We regard 
$\closedclosed{-\epsilon}{\epsilon}\times U_w$
as a subset of $M$.
By the mini-holomorphic structure,
$(E^t,\delbar_{E^t}):=
 (E,\delbar_E)_{|\{t\}\times U_w}$ is naturally a holomorphic
vector bundle on $U_w$ for any 
$t\in\closedclosed{-\epsilon}{\epsilon}$.
By considering the scattering map along the path 
$\gamma_w:[0,1]\lrarr \bigl(t_1+s(t_2-t_1),w\bigr)$ $(w\in U_w)$,
we obtain the isomorphism
$(E^{t_1},\delbar_{E^{t_1}})
 \simeq 
 (E^{t_2},\delbar_{E^{t_2}})$
because of the commutativity
$[\del_{E,t},\del_{E,\wbar}]=0$.

\subsection{Dirac type singularity of mini-holomorphic bundles}

Let $U$ be a neighbourhood of $(0,0)$ in $\real\times\cnum$.
We take $\epsilon>0$ and $\delta>0$
such that 
$\openopen{-2\epsilon}{2\epsilon}\times \{|w|<\delta\}$
is contained in $U$.
Set $U_w:=\{|w|<\delta\}$
and $U_w^{\ast}:=U_w\setminus\{0\}$.
Let $(E,\delbar_E)$ be a mini-holomorphic bundle
on $U\setminus \{(0,0)\}$.
For any $t\in\real$ with $0<|t|<2\epsilon$,
we obtain the holomorphic vector bundle
$(E^t,\delbar_{E^t})$ on $U_w$
as the restriction of $(E,\delbar_E)$.
We have the scattering map
$\Psi:E^{-\epsilon}_{|U_w^{\ast}}
\simeq E^{\epsilon}_{|U_w^{\ast}}$.
Recall that $(0,0)$ is called Dirac type singularity of
$(E,\delbar_E)$ if $\Psi$ is meromorphic at $0$.
(See \cite{Mochizuki-Yoshino}.)
In that case,
we say that $(E,\delbar_E)$ is of Dirac type 
on $(U,(0,0))$.
\index{Dirac type singularity (mini-holomorphic bundle)}

Let $M$ be a $3$-dimensional manifold
with a mini-complex manifold.
Let $Z$ be a discrete subset in $M$.
Let $(E,\delbar_E)$ be a mini-holomorphic bundle
on $M\setminus Z$.
We say that $(E,\delbar_E)$ is of Dirac type on $(M,Z)$
if the following holds.
\begin{itemize}
\item
 For any $P\in Z$,
 we take a mini-complex chart
 $(M_P,\varphi_P)$ 
 such that (i) $\varphi_P(P)=(0,0)$,
 (ii) the induced mini-holomorphic bundle
 on $\varphi_P(M_P\setminus\{P\})$
 is of Dirac type on $(\varphi_P(M_P),(0,0))$.
\end{itemize}
We can easily see that the condition is independent
of the choice of a mini-complex chart
$(M_P,\varphi_P)$.

\subsection{Kronheimer resolution of Dirac type singularity}
\label{subsection;17.9.19.2}

Let $U$ be a neighbourhood of $(0,0)$ in $\real\times\cnum$.
Let $(E,\delbar_E)$ be a mini-holomorphic bundle
on $U^{\ast}:=U\setminus \{(0,0)\}$.
Let $\varphi:\cnum^2\lrarr \real\times\cnum$
be the map given by
\begin{equation}
\label{eq;20.7.30.40}
 \varphi(u_1,u_2)=(|u_1|^2-|u_2|^2,2u_1u_2).
\end{equation}
We set $\Utilde:=\varphi^{-1}(U)$
and $\Utilde^{\ast}:=\Utilde\setminus\{(0,0)\}$.
The induced map
$\Utilde^{\ast}\lrarr U^{\ast}$ is 
also denoted by $\varphi$.
We set $\Etilde:=\varphi^{-1}(E)$.
We obtain the holomorphic structure
$\delbar_{\Etilde}$ on $\Etilde$
determined by the following condition:
\[
 \del_{\Etilde,\ubar_1}\varphi^{-1}(s)
=u_1\varphi^{-1}(\del_{E,t}s)
+2\ubar_2\varphi^{-1}(\del_{E,\wbar}s),
\]
\[
 \del_{\Etilde,\ubar_2}\varphi^{-1}(s)
=-u_2\varphi^{-1}(\del_{E,t}s)
+2\ubar_1\varphi^{-1}(\del_{E,\wbar}s).
\]
Here, $s$ denotes a $C^{\infty}$-section of $E$.
(See \cite{Charbonneau-Hurtubise, Mochizuki-Yoshino}.)
\begin{lem}[\cite{Charbonneau-Hurtubise, Mochizuki-Yoshino}]
$(0,0)$ is Dirac type singularity of $(E,\delbar_E)$
if and only if $(\Etilde,\delbar_{\Etilde})$
extends to a holomorphic vector bundle
$(\Etilde_0,\delbar_{\Etilde_0})$ on $\Utilde$.
\hfill\qed
\end{lem}
We call such $(\Etilde_0,\delbar_{\Etilde_0})$
the Kronheimer resolution of $(E,\delbar_E)$
at $(0,0)$.
\index{Kronheimer resolution}

\begin{rem}
Note that there exists a natural morphism
$\varphi^{-1}(\nbigo_{U^{\ast}})\lrarr\nbigo_{\Utilde^{\ast}}$.
Let $\nbige$ denote the $\nbigo_{U^{\ast}}$-module
obtained as the sheaf of mini-holomorphic sections of $E$.
It is easy to see that
the sheaf of holomorphic sections of $(\Etilde,\delbar_{\Etilde})$
is naturally isomorphic to
$\nbigo_{\Utilde^{\ast}}\otimes_{\varphi^{-1}(\nbigo_U)}
\varphi^{-1}(\nbige)$.
\hfill\qed
\end{rem}

\subsection{Precise description of Dirac type singularities}
\label{subsection;17.12.5.1}

Let us describe mini-holomorphic bundles of Dirac type 
more precisely.
For simplicity, we assume that
$U=\openopen{-2\epsilon}{2\epsilon}\times U_w$,
where $U_w$ is as above.
We put $U^{\ast}:=U\setminus\{(0,0)\}$.
We set $H_{>0}:=\openopen{0}{2\epsilon}\times \{0\}$.

\begin{lem}
\label{lem;21.9.6.2}
For any mini-holomorphic bundle $\nbige$ of rank $r$,
there exist a tuple of integers
$\ell_1\leq \ell_2\leq\cdots\leq\ell_r$
and an isomorphism
$\nbige\simeq
 \bigoplus_{i=1}^r
 \nbigo_{U^{\ast}}(\ell_i H_{>0})$.
The tuple of the integers $(\ell_1,\ldots,\ell_r)$ is called
the weight of  $\nbige$ at the Dirac type singularity $(0,0)$.
\end{lem}
\pf
Let $\nbige_{(\pm\epsilon,0)}$ 
denote the stalks of $\nbige$ at $(\pm\epsilon,0)$,
which are free $\nbigo_{\cnum,0}$-modules.
We have the isomorphism
$\nbige_{(-\epsilon,0)}(\ast 0)
\simeq
 \nbige_{(\epsilon,0)}(\ast 0)$.
It is a standard fact that
there exists a frame
$v_1,\ldots,v_r$ of 
$\nbige_{(-\epsilon,0)}$
such that
$\nbige_{(\epsilon,0)}(\ast 0)
=\bigoplus_{i=1}^r
 \nbigo_{\cnum,0} w^{-\ell_i}v_i$.
Then, the claim follows.
\hfill\qed

\begin{lem}
Let $\nbige$ be as in Lemma {\rm\ref{lem;21.9.6.2}}.
We set
$S^2_{\epsilon}=\{(t,w)\,|\,t^2+|w|^2=\epsilon^2\}\subset
U\setminus\{(0,0)\}$
with the orientation as the boundary of
$\{(t,w)\,|\,t^2+|w|^2\leq \epsilon^2\}$.
Then, we have
$\int_{S^2_{\epsilon}}c_1(\nbige)=\sum_{i=1}^r\ell_i$. 
\end{lem}
\pf
It is enough to check the claim in the case $\rank\nbige=1$.
There exists a frame $v$ of $\nbige_{(-\epsilon,0)}$
such that $w^{-\ell}v$ is a frame of $\nbige_{(\epsilon,0)}$
under the isomorphism
$\nbige_{(-\epsilon,0)}(\ast 0)\simeq
\nbige_{(\epsilon,0)}$.
In the case $\nbige'=\nbigo_{U}\cdot u$,
we clearly have
$\int_{S^2_{\epsilon}}c_1(\nbige')=0$.
There exists a complex structure of $S^2_{\epsilon}$
such that (i) it is compatible with the orientation,
(ii) $w$ induces a holomorphic coordinate around $(\epsilon,0)$.
Then, by the correspondence $v\mapsto u$,
we can identify $\nbige_{|S^2_{\epsilon}}$
with $\nbige'_{|S^2_{\epsilon}}\bigl(\ell(\epsilon,0) \bigr)$.
Hence, we obtain $\int_{S^2_{\epsilon}}c_1(\nbige)=\ell$.
\hfill\qed

\vspace{.1in}

We set
$H_{<0}:=
 \openopen{-2\epsilon}{0}\times \{0\}$.
Note that for any integers $\ell_i$ $(i=1,2)$
there exists the isomorphism
\[
 \nbigo_{U^{\ast}}((\ell_1+\ell_2) H_{>0})
\simeq
 \nbigo_{U^{\ast}}(\ell_1H_{>0}-\ell_2H_{<0})
\]
induced by the multiplication of $w^{\ell_2}$.

Let $j:U^{\ast}\lrarr U$ denote the inclusion.
We set $H_{\geq 0}:=\closedopen{0}{2\epsilon}\times \{0\}$.
We obtain
\[
 j_{\ast}\nbigo_{U^{\ast}}(\ell H_{>0})
\simeq
 \left\{
 \begin{array}{ll}
 \nbigo_U(\ell H_{\geq 0}) & (\ell<0),\\
 \nbigo_U(\ell H_{>0}) & (\ell\geq 0).
 \end{array}
 \right.
\]
Let $\iota_0$ denote the inclusion
$\{0\}\times U_w\lrarr U$.
We obtain
\[
 \iota_0^{-1}\bigl(j_{\ast}\nbigo_{U^{\ast}}(\ell H_{>0})\bigr)
\simeq
 \left\{
 \begin{array}{ll}
 \nbigo_{U_w}(\ell\{0\}) & (\ell<0),\\
 \nbigo_{U_w} & (\ell\geq 0).
 \end{array}
 \right.
\]

\subsection{Subbundles and quotient bundles}
\label{subsection;17.10.24.100}

Let $U$ be as in \S\ref{subsection;17.12.5.1}.
Let $\nbige$ be a mini-holomorphic bundle of Dirac type
on $(U,0)$.
\begin{lem}\mbox{{}}
\begin{itemize}
\item
Let $\nbige'$ be an $\nbigo_{U^{\ast}}$-submodule
of $\nbige$ which is also locally free.
Then, $(0,0)$ is Dirac type singularity of $\nbige'$.
\item
Let $\nbige''$ be a quotient $\nbigo_{U^{\ast}}$-module
of $\nbige$ which is also locally free.
Then, $(0,0)$ is Dirac type singularity of $\nbige''$.
\end{itemize}
\end{lem}
\pf
It is easy to see that the scattering maps of 
$\nbige'$  and $\nbige''$
are meromorphic at $0$.
Hence, $(0,0)$ is Dirac type singularity.
\hfill\qed

\begin{rem}
Let $0\lrarr \nbige'\lrarr \nbige\lrarr\nbige''\lrarr 0$
be an exact sequence of mini-holomorphic bundles
of Dirac type on $U^{\ast}$.
We have the Kronheimer resolutions
$\widetilde{\nbige}'$,
$\widetilde{\nbige}$
and
$\widetilde{\nbige}''$
at $(0,0)$,
as in {\rm\S\ref{subsection;17.9.19.2}}.
Note that 
$0\lrarr \widetilde{\nbige}'\lrarr\widetilde{\nbige}
 \lrarr\widetilde{\nbige}''\lrarr 0$
is not necessarily exact at $(0,0)$.
\hfill\qed
\end{rem}

\subsection{Basic functoriality}

Let $(E_i,\delbar_{E_i})$ $(i=1,2)$
be a $3$-dimensional mini-complex manifolds.
Let
$(E_1,\delbar_{E_1})
\oplus
(E_2,\delbar_{E_2})$
denote the mini-holomorphic bundle
$(E_1\oplus E_2,\delbar_{E_1\oplus E_2})$,
where 
the operator
$\delbar_{E_1\oplus E_2}:
C^{\infty}(M,E_1\oplus E_2)
\lrarr
C^{\infty}(M,(E_1\oplus E_2)\otimes\Omega^{0,1}_M)$
is determined by
$\delbar_{E_1\oplus E_2}(s_1\oplus s_2)
=\delbar_{E_1}(s_1)\oplus
\delbar_{E_2}(s_2)$
for $s_i\in C^{\infty}(M,E_i)$.
Let
$(E_1,\delbar_{E_1})
\otimes
(E_2,\delbar_{E_2})$
denote the mini-holomorphic bundle
$(E_1\otimes E_2,\delbar_{E_1\otimes E_2})$,
where the operator
$\delbar_{E_1\otimes E_2}:
C^{\infty}(M,E_1\otimes E_2)
\lrarr
C^{\infty}(M,(E_1\otimes E_2)\otimes\Omega^{0,1}_M)$
satisfies
\[
\delbar_{E_1\otimes E_2}(s_1\otimes s_2)
=\delbar_{E_1}(s_1)\otimes s_2 
+s_1\otimes \delbar_{E_2}(s_2)
\]
for $s_i\in C^{\infty}(M,E_i)$.
Let
$\Hom\bigl(
(E_1,\delbar_{E_1}),
(E_2,\delbar_{E_2})
\bigr)$
denote the mini-holomorphic bundle
\[
 (\Hom(E_1,E_2),\delbar_{\Hom(E_1,E_2)}),
\]
where the operator
$\delbar_{\Hom(E_1,E_2)}:
C^{\infty}(M,\Hom(E_1,E_2))
\lrarr
C^{\infty}(M,\Hom(E_1,E_2)\otimes\Omega^{0,1}_M)$
satisfies
\[
\delbar_{E_2}(f(s_1))
=\delbar_{\Hom(E_1,E_2)}(f)(s_1)
+f(\delbar_{E_1}s_1)
\]
for $f\in C^{\infty}(M,\Hom(E_1,E_2))$,
and $s_1\in C^{\infty}(M,E_1)$.

For a mini-holomorphic bundle $(E,\delbar_E)$ on $M$,
the dual $(E,\delbar_E)^{\lor}$
is defined as $\Hom((E,\delbar_E),(\cnum\times M,\delbar_M))$.
\index{direct sum (mini-holomorphic bundle)}
\index{tensor product (mini-holomorphic bundle)}
\index{inner homomorphism (mini-holomorphic bundle)}
\index{dual (mini-holomorphic bundle)}

\section{Monopoles}
\label{subsection;20.7.30.30}

\subsection{Monopoles and mini-holomorphic bundles}
\label{subsection;17.9.29.3}

Let $M$ be an oriented $3$-dimensional manifold
with a Riemannian metric $g_M$.
\begin{df}
 Let $E$ be a complex vector bundle with a Hermitian metric $h$,
 a unitary connection $\nabla$
 and an anti-self-adjoint endomorphism $\phi$.
 The tuple $(E,h,\nabla,\phi)$ is called a monopole on $M$
 if the Bogomolny equation
 $F(\nabla)-\ast\nabla\phi=0$,
 where
 $F(\nabla)$ denotes the curvature of $\nabla$,
 and
 $\ast$ denotes the Hodge star operator with respect to $g_M$.
\hfill\qed
\end{df}
\index{monopole $(E,h,\nabla,\phi)$}

Suppose that $M$ is equipped with a mini-complex structure
such that
the orthogonal decomposition
$TM=T_SM\oplus (T_SM)^{\bot}$
induces a splitting of $TM\lrarr T_QM$,
and that
the multiplication of $\sqrt{-1}$ on $T_QM$
is an isometry with respect to the induced metric of $T_QM$.
We identify $T_SM$ with the product bundle $\real\times M$.

Let $(E,\delbar_E)$ be a mini-holomorphic bundle
on $M$.
For any Hermitian metric $h$ of $E$,
the Chern connection $\nabla_h$
and the Higgs field $\phi_h$ are associated.
\begin{df}
We say that $(E,\delbar_E,h)$ is a monopole
if $(E,h,\nabla_h,\phi_h)$ is a monopole.
In the case,
we say that $(E,\delbar_E)$ is the mini-holomorphic bundle
underlying the monopole.
\hfill\qed
\end{df}
\index{monopole $(E,\delbar_E,h)$}
 
\subsection{Euclidean monopoles}
\label{subsection;16.9.20.20}

In this paper,
we study monopoles
on spaces which are locally isomorphic to
$\real\times\cnum$ with 
the natural metric and the natural mini-complex structure.
Let us look at the condition more explicitly 
in terms of local coordinate systems.

Let $U$ be an open subset of $\real\times\cnum$.
It is equipped with the metric $dt\,dt+dw\,d\wbar$.
It is also equipped with the mini-complex structure.
Let $(E,\delbar_E)$ be a mini-holomorphic bundle on $U$.
Let $h$ be a Hermitian metric of $E$.
We obtain the differential operator
$\del_{E,h,w}:E\lrarr E$ determined by the condition
\index{operator $\del_{E,h,w}$}
\[
 \del_{\wbar}
 h\bigl(u,v\bigr)
=h\bigl(\del_{E,\wbar}u,v\bigr)
+h\bigl(u,\del_{E,h,w}v\bigr).
\]
We also obtain the differential operator
$\del_{E,h,t}':E\lrarr E$ determined by the condition
\index{operator $\del_{E,h,t}'$}
\[
\del_th(u,v)
=h(\del_{E,t}u,v)+h(u,\del_{E,h,t}'v). 
\]
We set
\[
 \nabla_{h,x}:=
\del_{E,\wbar}+\del_{E,w},
\quad
 \nabla_{h,y}:=
 -\sqrt{-1}\bigl(
 \del_{E,\wbar}
-\del_{E,w}
 \bigr),
\]
\[
 \nabla_{h,t}=
 \frac{1}{2}
 \bigl(
 \del_{E,t}+\del_{E,h,t}'
 \bigr),
\quad
 \phi_h=\frac{\sqrt{-1}}{2}\bigl(
 \del_{E,t}-\del_{E,h,t}'
 \bigr).
\]
Thus, from the metric $h$,
we obtain the unitary connection
$\nabla_h(s)=
\nabla_{h,x}(s)\,dx+\nabla_{h,y}(s)\,dy+\nabla_{h,y}(s)\,dt$,
and the anti-self-adjoint section $\phi_h$ of $\End(E)$.
Note that the mini-holomorphic structure is recovered by
\[
\del_{E,\wbar}=\frac{1}{2}(\nabla_{h,x}+\sqrt{-1}\nabla_{h,y}),
\quad\quad
\del_{E,t}=\nabla_{h,t}-\sqrt{-1}\phi_h.
\]
We also write
$\del_{E,\wbar}$ and $\del_{E,h,w}$
as $\nabla_{h,\wbar}$ and $\nabla_{h,w}$.
Let $F(h)$ denote the curvature of $\nabla_h$,
which is expressed as
$F(h)=F(h)_{\wbar,t}d\wbar\,dt
+F(h)_{w,t}dw\,dt+F(h)_{\wbar,w}d\wbar\,dw$.
\index{curvature $F(h)$}

\begin{lem}
\label{lem;17.9.29.40}
We obtain
$F(h)_{\wbar,t}-\sqrt{-1}\nabla_{h,\wbar}\phi_h=0$
and 
$F(h)_{w,t}+\sqrt{-1}\nabla_{h,w}\phi_h=0$.
\end{lem}
\pf
We have the commutativity
$[\del_{E,\wbar},\del_{E,t}]=0$,
which implies
\[
 \bigl[
 \nabla_{h,\wbar},
 \nabla_{h,t}-\sqrt{-1}\phi_h
 \bigr]
=F(h)_{\wbar,t}-\sqrt{-1}\nabla_{h,\wbar}\phi_h=0.
\]
As the adjoint, we obtain 
$F(h)_{w,t}+\sqrt{-1}\nabla_{h,w}\phi_h=0$.
\hfill\qed

\begin{cor}
\label{cor;20.7.16.1}
A mini-holomorphic bundle $(E,\delbar_E)$
with a Hermitian metric $h$
is a monopole,
i.e.,
the associated tuple
 $(E,h,\nabla_h,\phi_h)$ is a monopole
if and only if the following equation is satisfied.
\begin{equation}
\label{eq;20.7.30.31}
\frac{\sqrt{-1}}{2}\nabla_{h,t}\phi_h
+F(h)_{\wbar,w}=0.
\end{equation}
\end{cor}
\pf
Because
$\ast(dw)=-\sqrt{-1}dw\,dt$,
$\ast(d\wbar)=\sqrt{-1}d\wbar\,dt$,
and $\ast(dt)=-\frac{\sqrt{-1}}{2}d\wbar\,dw$,
the claim follows from Lemma \ref{lem;17.9.29.40}.
\hfill\qed

\vspace{.1in}

Conversely,
let $E$ be a $C^{\infty}$-vector bundle on $U$
with a Hermitian metric $h$,
a unitary connection $\nabla$
and an anti-Hermitian endomorphism $\phi$.
We set
$\del_{E,\wbar}:=\nabla_{\wbar}$
and $\del_{E,t}:=\nabla_t-\sqrt{-1}\phi$.
We obtain the differential operator
$\delbar_E:C^{\infty}(U,E)\lrarr C^{\infty}(U,E\otimes\Omega^{0,1}_U)$
by $\delbar_{E}(s)=\del_{E,\wbar}(s)\,d\wbar+\del_{E,t}(s)\,dt$
for any $s\in C^{\infty}(U,E)$.
In general, $(E,\delbar_E)$ is not necessarily mini-holomorphic.
If $(E,h,\nabla,\phi)$ is a monopole,
the Bogomolny equation implies
the commutativity
$\bigl[\del_{E,\wbar},\del_{E,t}\bigr]=0$,
and hence
$(E,\delbar_E)$ is a mini-holomorphic bundle,
which is called the underlying mini-holomorphic bundle
of the monopole $(E,h,\nabla,\phi)$.
\index{the underlying mini-holomorphic bundle}
Note that $(E,h,\nabla,\phi)$ is recovered from
$(E,\delbar_E)$ with $h$.

\begin{rem}
\label{rem;21.8.22.11}
Let $(E,\delbar_E)$ be a mini-holomorphic bundle
with a Hermitian metric $h$.
Let $\nabla_h$ and $\phi_h$ be
the Chern connection and the Higgs field.
We also have the operator
$\del'_{E,h,t}=\nabla_{h,t}+\sqrt{-1}\phi_h$.
By the construction,
we have
\[
 [\nabla_{h,\wbar},\del'_{E,h,t}]
=2[\nabla_{h,\wbar},\nabla_{h,t}],\quad
[\nabla_{h,w},\del_{E,t}] 
=2[\nabla_{h,w},\nabla_{h,t}].
\]
We also have
\[
\nabla_{h,\wbar}\phi_h
=-\frac{\sqrt{-1}}{2}[\nabla_{h,\wbar},\del'_{E,h,t}],\quad
\nabla_{h,w}\phi_h
 =\frac{\sqrt{-1}}{2}[\nabla_{h,w},\del_{E,t}],
\]
\[ \nabla_{h,t}(\phi_h)
=-\frac{\sqrt{-1}}{2}[\del_{E,t},\del'_{E,h,t}].
\]
\hfill\qed
\end{rem}

\subsection{Dirac type singularity}
\label{subsection;21.8.3.3}

We recall the notion of Dirac type singularity of monopoles.
It is originally introduced by Kronheimer \cite{Kronheimer-Master-Thesis},
and it was later generalized by Pauly \cite{Pauly}
to the context of general $3$-dimensional Riemannian manifolds,
and by Charbonneau and Hurtubise
\cite{Charbonneau-Hurtubise} in the higher rank case.

Let $\varphi:\cnum^2\lrarr \real\times\cnum$
be the map (\ref{eq;20.7.30.40}).
Let $U$ be a neighbourhood of $(0,0)$ in $\real\times\cnum$.
We set $U^{\ast}:=U\setminus\{(0,0)\}$.

Let $(E,h,\nabla,\phi)$ be a monopole on $U^{\ast}$.
We set
$\xi:=
 -u_1d\ubar_1+\ubar_1du_1
 -\ubar_2du_2+u_2d\ubar_2$.
We put
$(\Etilde,\htilde):=
 \varphi^{-1}(E,h)$.
We obtain the unitary connection
$\nablatilde:=
 \varphi^{\ast}(\nabla)
+\sqrt{-1}\varphi^{\ast}(\phi)\otimes\xi$.
As proved in \cite{Kronheimer-Master-Thesis},
$(\Etilde,\htilde,\nablatilde)$ is an instanton
on $\varphi^{-1}(U^{\ast})$.
Recall that 
$(0,0)$ is called Dirac type singularity of 
$(E,h,\nabla,\phi)$ 
if $(\Etilde,\htilde,\nablatilde)$
extends to an instanton on $\varphi^{-1}(U)$.
It particularly implies that
the underlying mini-holomorphic bundle $(E,\delbar_E)$ of
$(E,h,\nabla,\phi)$ is of Dirac type 
on $(U,(0,0))$.
\index{Dirac type singularity (monopole)}

We have the following simple characterization of
Dirac type singularity
\begin{thm}[\mbox{\cite[Theorem 1]{Mochizuki-Yoshino}}]
$(0,0)$ is a Dirac type singularity of $(E,h,\nabla,\phi)$
if and only if
$|\phi|_h=O\bigl((t^2+|w|^2)^{-1/2}\bigr)$.
\hfill\qed
\end{thm}

We also have the following characterization
which easily follows from \cite[Theorem 2]{Mochizuki-Yoshino}.
\begin{prop}
\label{prop;17.10.13.151}
Let $(E,\delbar_E,h)$ be a monopole on $U^{\ast}$.
It has a Dirac type singularity at $(0,0)$
if and only if the following holds.
\begin{itemize}
\item
 $(E,\delbar_E)$ is a mini-holomorphic bundle 
 of Dirac type on $(U,(0,0))$.
 In particular,
 we obtain the Kronheimer resolution
 $(\Etilde_0,\delbar_{\Etilde_0})$ on $\varphi^{-1}(U)$.
 (See {\rm\S\ref{subsection;17.9.19.2}}.)
\item
 Let $h_1$ be a $C^{\infty}$-metric of $E$
 such that $\varphi^{-1}(h_1)$ induces
 a $C^{\infty}$-metric on $\Etilde_0$.
 Then, there exist positive constants
 $C_1>1$ and $N_1$ such that
$C_1^{-1}(|t|^2+|w|^2)^{-N_1}\cdot h_1
\leq h
\leq
 C_1(|t|^2+|w|^2)^{N_1}\cdot h_1$.
\hfill\qed
\end{itemize}
\end{prop}

\subsubsection{Dirac monopoles (examples)}

The fundamental examples are Dirac monopoles.
Let us recall the description from
\cite[\S5.2]{Mochizuki-Yoshino}
with a minor correction.
We set $U=\real\times\cnum$,
$U^{\ast}=U\setminus\{(0,0)\}$,
and $A_{\pm}:=U\setminus \{(t,0)\,|\,\pm t\leq 0\}$.
We have $U^{\ast}=A_+\cup A_-$.
We regard $U$ as a mini-complex manifold
in a natural way.
For any integer $m$,
let $L^{(m)}$ be the mini-holomorphic line bundle on $U^{\ast}$
equipped with mini-holomorphic frames
$\sigma^{(m)}_{\pm}$ of $L^{(m)}_{|A_{\pm}}$
such that
$\sigma^{(m)}_-=(w/2)^m\sigma^{(m)}_+$
on $A_+\cap A_-$.
Let $h^{(m)}$ be the Hermitian metric of $L^{(m)}$
determined by
\[
h^{(m)}(\sigma_+^{(m)},\sigma_+^{(m)})
=2^{m}(t+R)^{-m},
\]
\[
 h^{(m)}(\sigma_-^{(m)},\sigma_-^{(m)})
 =2^{-m}(-t+R)^{-m},
\]
where $R=\sqrt{t^2+|w|^2}$.
We obtain the Chern connection $\nabla^{(m)}$
and the Higgs field $\phi^{(m)}$
associated with $(L^{(m)},\delbar_{L^{(m)}},h^{(m)})$.

We obtain the holomorphic line bundle
$\Ltilde^{(m)}$ on
$\cnum^2\setminus\{(0,0)\}=\varphi^{-1}(U^{\ast})$
as in \S\ref{subsection;17.9.19.2}.
It is equipped with a holomorphic frame $e^{(m)}$
such that
$e^{(m)}=u_1^m\varphi^{-1}(\sigma_+^{(m)})$
on $\varphi^{-1}(A_+)$,
and
$e^{(m)}=u_2^{-m}\varphi^{-1}(\sigma_-^{(m)})$
on $\varphi^{-1}(A_-)$.
Let $\htilde^{(m)}_0$ be the Hermitian metric
of $\Ltilde^{(m)}$
defined as
$\htilde^{(m)}_0(e^{(m)},e^{(m)})=1$.
It is easy to check that $\htilde^{(m)}_0$
is the pull back of $h^{(m)}$.
Because $(\Ltilde^{(m)},\htilde^{(m)})$ is an instanton,
we obtain that
$(L^{(m)},h^{(m)},\nabla^{(m)},\phi^{(m)})$
is a monopole.

Let us compute $\nabla^{(m)}$ and $\phi^{(m)}$ explicitly.
$\del_{L^{(m)},\wbar}\sigma^{(m)}_+=\del_{L(m),t}\sigma^{(m)}_+=0$.
We have
\[
\del_{L^{(m)},w}\sigma^{(m)}_+
=\sigma^{(m)}_+\del_w\log h^{(m)}(\sigma^{(m)}_+,\sigma^{(m)}_+)
=\sigma^{(m)}_+\frac{-m\wbar}{2R(t+R)}.
\]
We also have
\[
 \del'_{L^{(m)},t}\sigma^{(m)}_+
=\sigma^{(m)}_+\del_t\log h^{(m)}(\sigma^{(m)}_+,\sigma^{(m)}_+)
=\sigma^{(m)}_+\frac{-m}{R}.
\]
Therefore, we obtain
\[
 \nabla^{(m)}\sigma^{(m)}_+
 =\sigma^{(m)}_+
 \Bigl(
 -\frac{m\wbar}{2R(t+R)}\,dw
-\frac{m}{2R}\,dt
 \Bigr),
 \quad\quad
 \phi^{(m)}=\frac{\sqrt{-1}m}{2R}.
\]
\begin{rem}
In {\rm\cite[\S5.2]{Mochizuki-Yoshino}},
$\phi^{(m)}=\frac{\sqrt{-1}m}{R}$ should be corrected to
$\phi^{(m)}=\frac{\sqrt{-1}m}{2R}$.
\hfill\qed 
\end{rem}

\subsection{Basic functoriality}

Let $M$ be a $3$-dimensional mini-complex manifold
with a Riemannian metric.
Assume that it is locally isomorphic to $\real\times\cnum$
with the canonical mini-complex structure
and the Euclidean metric.
Let $(E_i,h_i,\nabla_i,\phi_i)$ $(i=1,2)$
be monopoles on $M$.
The underlying mini-holomorphic bundle
is denoted by $(E_i,\delbar_{E_i})$.

The vector bundle
$E_1\oplus E_2$
is equipped with the naturally induced Hermitian metrics
$h_{E_1\oplus E_2}$,
the unitary connection $\nabla_{E_1\oplus E_2}$,
and the Higgs field 
$\phi_{E_1\oplus E_2}
=\phi_{E_1}\oplus \phi_{E_2}$.
The tuple 
\[
(E_1,h_1,\nabla_1,\phi_1)
\oplus
(E_2,h_2,\nabla_2,\phi_2)
=
(E_1\oplus E_2,h_{E_1\oplus E_2},
 \nabla_{E_1\oplus E_2},
 \phi_{E_1\oplus E_2}) 
\]
is a monopole.
The underlying mini-holomorphic bundle
is naturally isomorphic to
$(E_1,\delbar_{E_1})\oplus
(E_2,\delbar_{E_2})$.

The vector bundle
$E_1\otimes E_2$
is equipped with the naturally induced Hermitian metrics
$h_{E_1\otimes E_2}$,
the unitary connection $\nabla_{E_1\otimes E_2}$,
and the Higgs field 
$\phi_{E_1\otimes E_2}
=\phi_{E_1}\otimes\id_{E_2}+\id_{E_1}\otimes\phi_{E_2}$.
The tuple 
\[
(E_1,h_1,\nabla_1,\phi_1)
\otimes
(E_2,h_2,\nabla_2,\phi_2)
=
(E_1\otimes E_2,h_{E_1\otimes E_2},
 \nabla_{E_1\otimes E_2},
 \phi_{E_1\otimes E_2}) 
\]
is a monopole.
The underlying mini-holomorphic bundle
is naturally isomorphic to
$(E_1,\delbar_{E_1})\otimes
(E_2,\delbar_{E_2})$.

The vector bundle
$\Hom(E_1,E_2)$
is equipped with
the naturally induced Hermitian metrics
$h_{\Hom(E_1,E_2)}$,
the unitary connection $\nabla_{\Hom(E_1,E_2)}$,
and the Higgs field 
$\phi_{\Hom(E_1,E_2)}$
given as
$\phi_{\Hom(E_1,E_2)}(s)=\phi_2\circ s-s\circ\phi_1$.
The tuple 
\begin{multline}
\Hom\Bigl(
(E_1,h_1,\nabla_1,\phi_1),
(E_2,h_2,\nabla_2,\phi_2)
\Bigr)
= \\
 (\Hom(E_1,E_2),h_{\Hom(E_1,E_2)},
 \nabla_{\Hom(E_1,E_2)},
 \phi_{\Hom(E_1,E_2)}) 
\end{multline}
is a monopole.
The underlying mini-holomorphic bundle
is naturally isomorphic to
$\Hom\bigl(
(E_1,\delbar_{E_1}),
(E_2,\delbar_{E_2})
\bigr)$.
\index{direct sum (monopole)}
\index{tensor product (monopole)}
\index{inner homomorphism (monopole)}
\index{dual (monopole)}

The product line bundle
$\cnum\times M$ on $M$
is equipped with the metric $h_0$
defined by $h_0(1,1)=1$,
and the trivial unitary connection $\nabla_0$
and the Higgs field $0$,
and the tuple
$(\cnum\times M,h_0,\nabla_0,0)$ is a monopole.
For any monopole $(E,h,\nabla,\phi)$ on $M$,
we obtain the induced monopole
$(E,h,\nabla,\phi)^{\lor}
=\Hom\Bigl(
(E,h,\nabla,\phi),
(\cnum\times M,h_0,\nabla_0,0)
\Bigr)$
as the dual.
The underlying mini-holomorphic bundle
of $(E,h,\nabla,\phi)^{\lor}$
is naturally isomorphic to
the dual of the mini-holomorphic bundle
underlying $(E,h,\nabla,\phi)$.

\section{Dimensional reduction from $4D$ to $3D$}
\label{subsection;20.7.30.41}

\subsection{Instantons induced by monopoles}
\label{subsection;17.10.3.10}

Monopoles are regarded as the $1$-dimensional reduction
of instantons.
Namely, monopoles are equivalent to
$\real$-equivariant instantons.
Let us recall the construction explicitly
in terms of coordinate systems.

Let $U$ be an open subset of $\real\times\cnum$.
We use real and complex coordinates
$t$ and $w=x+\sqrt{-1}y$
for $\real$ and $\cnum$,
respectively.
Let $p:\cnum\times\cnum\lrarr \real\times\cnum$
be the map given by
$(s+\sqrt{-1}t,w)\longmapsto (t,w)$.
We use the standard Euclidean metrics
$dt\,dt+dw\,d\wbar$ on $\real\times\cnum$,
and 
$dz\,d\zbar+dw\,d\wbar$ on $\cnum\times\cnum$.

Let $E$ be a $C^{\infty}$-bundle on $U$
with a Hermitian metric $h$,
a unitary connection $\nabla_E$
and an anti-self-adjoint endomorphism $\phi$
of $E$.
We set
$(\Etilde,\htilde):=
 p^{\ast}(E,h)$
and 
$\nabla_{\Etilde}:=p^{\ast}\nabla_E+p^{\ast}\phi\,ds$
on $p^{-1}(U)$.

\begin{prop}[Hitchin]
$(E,h,\nabla_E,\phi)$ is a monopole
if and only if 
$(\Etilde,\htilde,\nabla_{\Etilde})$
is an instanton.
\end{prop}
\pf
Let $\ast_3$ denote the Hodge star operator
of $\real^3$ with the orientation
$dt\wedge dx\wedge dy$.
For example,
we have
$\ast_3(dx\wedge dy)=dt$
and 
$\ast_3(dx)=dy\wedge dt$.
Note that $\ast_3\circ\ast_3=\id$.
Let $\ast_4$ denote the Hodge star operator on
$\real^4$
with the orientation
$ds\wedge dt\wedge dx\wedge dy$.
For example,
we have
$\ast_4(dx\wedge dy)
=ds\wedge dt$
and
$\ast_4 (dx\wedge ds)=-dy\wedge dt$.

By the construction of $\nabla_{\Etilde}$,
we obtain
$F(\nabla_{\Etilde})
=p^{\ast}F(\nabla_E)
+p^{\ast}(\nabla_E \phi)\,ds$,
which implies
\[
 \ast_4F(\nabla_{\Etilde})
=-p^{\ast}(\ast_3F(\nabla_E))\,ds
-p^{\ast}(\ast_3\nabla_E\phi).
\]
Hence, we obtain
\[
 \ast_4F(\nabla_{\Etilde})
+F(\nabla_{\Etilde})
=p^{\ast}\bigl(F(\nabla_E)-\ast_3\nabla_E\phi\bigr)
+p^{\ast}\bigl(\nabla_E\phi-\ast_3F(\nabla_E)\bigr)\,ds.
\]
Therefore,
$\ast_4F(\nabla_{\Etilde})+F(\nabla_{\Etilde})=0$
holds
if and only if
$F(\nabla_E)=\ast_3\nabla_E\phi$
holds.
\hfill\qed

\begin{rem}
Let us consider an $\real$-action
on $p^{-1}(U)$
by $T\cdot (z,w)=(z+T,w)$.
Then,
the proposition means  that
$\real$-equivariant instantons on $p^{-1}(U)$
correspond to
monopoles on $U$.
\hfill\qed
\end{rem}

The following lemma is clear by the construction.

 \begin{lem}
Let $f$ be any section of $E$.
Then,
$(\nabla_{\Etilde,s}+\sqrt{-1}\nabla_{\Etilde,t})p^{\ast}f=0$
holds
if and only if
$(\nabla_{E,t}-\sqrt{-1}\phi)f=0$ holds.
In other words,
$\real$-invariant $\nablatilde_{\zbar}$-holomorphic
sections on $p^{-1}(U)$
correspond to
$(\nabla_{E,t}-\sqrt{-1}\phi)$-flat sections
on $U$.
\hfill\qed
 \end{lem}

The following lemma follows from
the estimate for instantons
due to Uhlenbeck \cite{Uhlenbeck1}.

\begin{lem}
\label{lem;17.10.5.10}
Let $(E,h,\nabla_E,\phi)$ be a monopole on $U$.
Suppose $\bigl|\nabla_E\phi\bigr|_h\leq \epsilon$,
or equivalently 
$\bigl|F(\nabla_E)\bigr|_h\leq \epsilon$
for a positive number $\epsilon$.
Let $K$ be any relatively compact subset of $U$.
If $\epsilon$ is sufficiently small,
the higher derivatives of
$\nabla_E\phi$ and $F(\nabla_E)$
are dominated by $C_K\epsilon$ on $K$
for some $C_K>0$,
where $C_K$ depends on the order of derivatives and $K$,
but independent of $(E,h,\nabla_E,\phi)$.
\hfill\qed
\end{lem}

We write some formulas.
We clearly have
$\nabla_{\Etilde,\wbar}(p^{-1}(u))
=p^{-1}(\nabla_{E,\wbar}u)$
and 
$\nabla_{\Etilde,w}(p^{-1}(u))
=p^{-1}(\nabla_{E,w}u)$
for any $u\in C^{\infty}(U,E)$.
As for the derivative in the $z$-direction,
we have 
\[
 \nabla_{\Etilde,\zbar}(p^{-1}(u))
=\frac{\sqrt{-1}}{2}p^{-1}\bigl((\nabla_t-\sqrt{-1}\phi)u\bigr),
\]
\[
 \nabla_{\Etilde,z}(p^{-1}(u))
=-\frac{\sqrt{-1}}{2}p^{-1}\bigl(
 (\nabla_t+\sqrt{-1}\phi) u\bigr).
\]
For the real coordinate system $(s,t)$
induced by $z=s+\sqrt{-1}t$,
we obtain
$\nabla_{\Etilde,s}p^{-1}(u)
 =p^{-1}(\phi u)$
and
$\nabla_{\Etilde,t}p^{-1}(u)
=p^{-1}(\nabla_tu)$.

\subsection{Holomorphic bundles and mini-holomorphic bundles}
\label{subsection;13.11.29.2}

We have the corresponding procedure 
to construct a holomorphic bundle
from a mini-holomorphic bundle.
Let $E$ be a $C^{\infty}$-vector bundle on $U$
with a mini-holomorphic structure
$\delbar_E:C^{\infty}(U,E)\lrarr C^{\infty}(U,E\otimes\Omega^{0,1})$.
We have the operators
$\del_{E,\wbar}$
and $\del_{E,t}$ on $C^{\infty}(U,E)$
satisfying 
$[\del_{E,\wbar},\del_{E,t}]=0$.
We set $\Etilde:=p^{-1}(E)$.
We have the holomorphic structure $\delbar_{\Etilde}$
satisfying the following formula
for any $C^{\infty}$-section $v$ of $E$:
\begin{equation}
\label{eq;17.9.29.21}
 \delbar_{\Etilde}\bigl(p^{-1}(v)\bigr)
=d\wbar\wedge 
 p^{-1}\bigl(\del_{\wbar}v\bigr)
+d\zbar\wedge \frac{\sqrt{-1}}{2}
 p^{-1}\bigl(\del_tv\bigr).
\end{equation}

Let $h$ be a Hermitian metric of $E$.
It induces a Hermitian metric $\htilde$ of $\Etilde$.
We can easily check that
$(E,\delbar_E,h)$ is a monopole
on $U$ with respect to
the Euclidean metric $dt\,dt+dw\,d\wbar$
if and only if 
$(\Etilde,\delbar_{\Etilde},\htilde)$
is an instanton on $p^{-1}(U)$
with respect to
the Euclidean metric $dz\,d\zbar+dw\,d\wbar$.
This is compatible
with the construction in \S\ref{subsection;17.9.29.3}.

\section{Dimensional reduction from $3D$ to $2D$}
\label{subsection;20.7.30.42}

\subsection{Monopoles induced by harmonic bundles}
\label{subsection;17.9.29.1}

Let us recall that the concept of harmonic bundle
was discovered by Hitchin \cite{Hitchin-self-duality}
as the $2$-dimensional reduction of instantons.
Because monopoles are regraded as the $1$-dimensional reduction
of instantons,
harmonic bundles are the $1$-dimensional reduction of
monopoles.
We recall the construction in an explicit way.

Let $U$ be a Riemann surface 
equipped with a holomorphic $1$-form $\varphi$
which is nowhere vanishing on $U$.
Let $(V,\delbar_V,h_V,\theta_V)$ be a harmonic bundle
on $U$.
We obtain the holomorphic endomorphism $f$ of $V$
determined by $\theta_V=f\,\varphi$.
Let $Y$ be a real $1$-dimensional manifold
equipped with a closed real $1$-form $\tau$
which is nowhere vanishing on $Y$.
The product $Y\times U$ is naturally equipped with
a mini-holomorphic structure
and the metric
$\tau\tau+\varphi\,\varphibar$.

Let $p_2:Y\times U\lrarr U$ be the projection.
We set 
\begin{equation}
 \label{eq;17.10.6.21}
 (E,h_E):=p_2^{\ast}(V,h_V),
\quad
 \nabla_E=p_2^{\ast}\nabla_V-\sqrt{-1}p_2^{\ast}(f+f^{\dagger})\,\tau,
\quad
 \phi_E=p_2^{\ast}(f-f^{\dagger}).
\end{equation}
Here, $p_2^{\ast}\nabla_V$
denotes the connection of $E$
induced as the pull back of $\nabla_V$.
Then, 
$(E,h_E,\nabla_E,\phi_E)$ is a monopole
on $Y\times U$
with the metric $\tau\tau+\varphi\varphibar$.
We set
$\Hit_2^3(V,\delbar_V,\theta_V,h_V):=
(E,h_E,\nabla_E,\phi_E)$.

Let $\gminiv$ be the vector field of $Y$
such that $\langle\tau,\gminiv\rangle=1$.
It naturally defines a vector field on $Y\times U$,
which is also denoted by $\gminiv$.
Let $L_{\gminiv}$ be the differential operator of $E$
defined by
$L_{\gminiv}(g\,p_2^{\ast}v)
=\gminiv(g)p_2^{\ast} v$
for any $g\in C^{\infty}(Y\times U)$
and any $C^{\infty}$-section $v$ of $V$.
Then, we have
\[
 \nabla_{E,\gminiv}-\sqrt{-1}\phi_E
=L_{\gminiv}
-\sqrt{-1}p_2^{\ast}(f+f^{\dagger})
-\sqrt{-1}p_2^{\ast}(f-f^{\dagger})
=L_{\gminiv}
-2\sqrt{-1}p_2^{\ast}f.
\]

\subsection{Mini-holomorphic bundles
 induced by holomorphic bundles with a Higgs field}
\label{subsection;17.9.29.2}

We have the corresponding procedure 
to construct a mini-holomorphic bundle
on $Y\times U$
from a holomorphic bundle $V$ with an endomorphism $f$
on $U$.
We use the natural splitting
$\Omega^{0,1}_{Y\times U}
=(T^{\cnum}Y)^{\lor}\oplus \Omega^{0,1}_U$.

Let $(V,\delbar_V,\theta_V)$ be a Higgs bundle on $U$.
We have the expression $\theta_V=f\,\varphi$.
We set
\[
 E:=p_2^{-1}(V),
\quad
 \delbar_E:=
 p_2^{\ast}(\delbar_V)-2\sqrt{-1}p_2^{\ast}f\,\tau.
\]
Here, $p_2^{\ast}(\delbar_V)$ is the mini-holomorphic structure
of $E$ obtained as the pull back of $\delbar_V$.
We obtain a mini-holomorphic bundle
$(E,\delbar_E)$ on $Y\times U$.
We denote it by
$p_2^{\ast}(V,\delbar_V,\theta_V)$
or 
$p_2^{\ast}(V,\delbar_V,f)$.
\index{Mini-holomorphic bundle
$p_2^{\ast}(V,\delbar_V,\theta_V)$,
$p_2^{\ast}(V,\delbar_V,f)$}

If $V$ is equipped with a Hermitian metric $h_V$,
we have the induced metric $h=p_2^{\ast}(h_V)$ 
on $p_2^{\ast}(V,\delbar_V,\theta_V)$.
The Chern connection and the Higgs fields
are induced by the formula (\ref{eq;17.10.6.21}).
If $h_V$ is a harmonic metric 
of the Higgs bundle $(V,\delbar_V,\theta_V)$,
then $h$ satisfies the Bogomolny equation
for $p_2^{\ast}(V,\delbar_V,\theta_V)$
with respect to the metric
$\tau\tau+\varphi\varphibar$.
This is compatible
with the construction in \S\ref{subsection;17.9.29.1}.

\begin{rem}
The above procedure describes 
the dimensional reduction 
of the underlying objects
of monopoles and harmonic bundles
at the twistor parameter $0$.
There exists a similar procedure for each twistor parameter $\lambda$
which we shall describe in
{\rm\S\ref{section;21.8.12.101}}.
\hfill\qed
\end{rem}

\subsection{Mini-holomorphic sections and monodromy}
\label{subsection;17.10.5.1}

Let us consider $Y_1=\real$ equipped with
the nowhere vanishing closed $1$-form $\tau_1=dt$,
where $t$ is the standard coordinate of $\real$.
Let $p_{2,1}:Y_1\times U\lrarr U$ denote the projection.
For any Higgs bundle $(V,\delbar_V,\theta_V)$ on $U$,
we obtain the mini-holomorphic bundle
$p_{2,1}^{\ast}(V,\delbar_V,\theta_V)$
on $Y_1\times U$
as explained in \S\ref{subsection;17.9.29.2}.
The following lemma is clear by the construction.
 \begin{lem}
For any holomorphic section $v$ of $V$ on $U$,
\[
   \exp\bigl(2\sqrt{-1}t p_{2,1}^{\ast}(f)\bigr)p_{2,1}^{\ast}(v)
\]
is a mini-holomorphic section of
$p_{2,1}^{\ast}(V,\delbar_V,\theta_V)$.
\hfill\qed
 \end{lem}

Let $Y_2=\real/T\seisuu$ equipped with $\tau_2=dt$.
Let $p_{2,2}:Y_2\times U\lrarr U$ denote the projection.
For any Higgs bundle $(V,\delbar_V,\theta_V)$ on $U$,
we obtain the mini-holomorphic bundle
$p_{2,2}^{\ast}(V,\delbar_V,\theta_V)$
on $Y_2\times U$.

For $P\in U$,
let $\gamma_P$ denote the loop 
$[0,T]\lrarr Y\times U$
given by
$s\longmapsto (s,P)$.
We identify $V_{|P}$
and $E_{|(0,P)}$.
By the mini-holomorphic structure,
we obtain the parallel transport along $\gamma_P$
which induces the monodromy automorphism
$M(\gamma_P)$ of $V_{|P}$.
We obtain the following lemma
from the previous one.
\begin{lem}
$M(\gamma_P)=\exp(2\sqrt{-1}Tf_{|P})$.
\hfill\qed
\end{lem}

\subsection{Appendix: monopoles as harmonic bundles of infinite rank}

Let $Y\times U$ with $\tau\tau+\varphi\varphibar$ be
as in \S\ref{subsection;17.9.29.1}.
Let $E$ be a $C^{\infty}$-vector bundle
on $Y\times U$
equipped with a Hermitian metric $h$,
a unitary connection $\nabla$
and an anti-Hermitian endomorphism $\phi$.
Let $\nbige^{\infty}$ denote
the sheaf of $C^{\infty}$-sections of $E$.
Let $\nbigc^{\infty}_U$ denote the sheaf of $C^{\infty}$-functions
on $U$.
\index{sheaf $\nbigc^{\infty}_U$}
We obtain the $\nbigc^{\infty}_U$-module $p_{2!}\nbige^{\infty}$,
where $p_{2!}$ denote the proper push-forward with respect to $p_2$.
It is equipped with the Hermitian metric $p_{2!}h$
defined as
\[
p_{2!}h(s_1,s_2)=\int_{Y}h(s_1,s_2)\tau
\]
for sections $s_i$ of $p_{2!}\nbige^{\infty}$,
which we can naturally regard
as $C^{\infty}$-sections of $E$ whose supports
are proper over $U$.
By using the $\Omega^{0,1}(U)$-part
$\delbar_{E}^Q$ of $\nabla$
in the decomposition
$T^{\ast}(Y\times U)\otimes\cnum
=T^{\ast}(Y)\otimes\cnum
\oplus \Omega^{1,0}(U)\oplus\Omega^{0,1}(U)$,
we obtain the holomorphic structure
$\delbar_{p_{2!}\nbige^{\infty}}$ of $p_{2!}\nbige^{\infty}$
by
$\delbar_{p_{2!}\nbige^{\infty}}(s)=
\delbar_{E}^Q(s)$.
Similarly, by using the $\Omega^{1,0}(U)$-part
$\del_{E}^Q$ of $\nabla$,
we obtain the differential operator
$\del_{p_{2!}\nbige^{\infty}}$ of $p_{2!}\nbige^{\infty}$
by 
$\del_{p_{2!}\nbige^{\infty}}(s)=
\del_{E}^Q(s)$.
It is easy to check that
$\delbar_{p_{2!}\nbige^{\infty}}
+\del_{p_{2!}\nbige^{\infty}}$
is unitary with respect to $p_{2!}h$.
We define the $\nbigc^{\infty}_U$-endomorphisms $f$ and $f^{\dagger}$ of
$p_{2!}\nbige^{\infty}$ by
\[
 f(s)=\frac{\sqrt{-1}}{2}(\nabla_{\gminiv}-\sqrt{-1}\phi)s,
 \quad
 f^{\dagger}(s)=\frac{\sqrt{-1}}{2}(\nabla_{\gminiv}+\sqrt{-1}\phi) s.
\]
We set $\theta=f\,\varphi$ and $\theta^{\dagger}=f^{\dagger}\varphibar$,
which are mutually adjoint with respect to $p_{2!}h$.
It is easy to see that
$(E,h,\nabla,\phi)$ is a monopole
if and only if
$(p_{2!}\nbige^{\infty},\delbar_{p_{2!}\nbige^{\infty}},
\theta,p_{2!}h)$ is a harmonic bundle of infinite rank
in the sense that
$\delbar_{p_{2!}\nbige^{\infty}}(f)=0$
and
$[\delbar_{p_{2!}\nbige^{\infty}},\del_{p_{2!}\nbige^{\infty}}]
+[\theta,\theta^{\dagger}]=0$
are satisfied.

Let $(E_1,\delbar_{E_1})$ be a mini-holomorphic bundle 
on $Y\times U$ with the natural mini-complex structure.
Similarly, let $\nbige^{\infty}_1$ denote the sheaf of
$C^{\infty}$-sections of $E_1$,
and we obtain the $\nbigc_U^{\infty}$-module
$p_{2!}\nbige^{\infty}_1$.
It is equipped with the holomorphic structure
$\delbar_{p_{2!}\nbige^{\infty}_1}$
induced by
the $\Omega^{0,1}(U)$-component of $\delbar_{E_1}$.
Let $\del_{E_1,\gminiv}$ be the differential operator of $E_1$
induced by $\delbar_{E_1}$ and $\gminiv$.
We obtain the endomorphism $f_1$ of $p_{2!}\nbige_1^{\infty}$
by $f_1(s)=\del_{E_1,\gminiv}(s)$,
which is holomorphic with respect to $\delbar_{p_{2!}\nbige^{\infty}_1}$.
We obtain the Higgs field $\theta_1=f_1\,\varphi$.
Thus, we obtain a Higgs bundle
$(p_{2!}\nbige_1^{\infty},\delbar_{p_2!\nbige_1^{\infty}},\theta_1)$
of infinite rank on $U$.
Clearly, if the mini-holomorphic bundle
$(E_1,\delbar_{E_1})$ underlies a monopole,
the Higgs bundle 
$(p_{2!}\nbige_1^{\infty},\delbar_{p_2!\nbige_1^{\infty}},\theta_1)$
underlies the corresponding harmonic bundle of infinite rank.

\section{Twistor families of mini-complex structures
on $\real\times\cnum$ and $(\real/T\seisuu)\times\cnum$}
\label{subsection;17.10.2.1}

\subsection{Preliminary}
\label{subsection;17.10.3.1}

Let $X$ be a $2$-dimensional complex vector space.
Let $L\subset X$ be an oriented $1$-dimensional real vector space.
Let $M:=X/L$ be the quotient $3$-dimensional real vector space.
Note that $M$ is equipped with a naturally induced mini-complex
structure.
Indeed, there exists a $\cnum$-linear isomorphism
$\varphi:X\simeq \cnum^2=\{(z,w)\}$
such that $\varphi(L)=\{(s,0)\,|\,s\in\real\}$
in a way compatible with the orientations.
Let $t$ be the imaginary part of $z$.
Then, by the mini-complex coordinate system $(t,w)$,
we obtain the isomorphism
$M\simeq\real\times\cnum$
which induces a mini-complex structure on $M$.
We can check the following claim
by a direct computation.
\begin{lem}
The mini-complex structure
is independent of a choice of $\varphi$.
\hfill\qed
\end{lem}

\begin{rem}
\label{rem;21.8.10.2}
Let $L_{\cnum}\subset X$ denote the $1$-dimensional complex vector space
generated by $L$.
There exists the natural exact sequence
of $\real$-vector spaces
\begin{equation}
\label{eq;21.8.10.1}
 0\lrarr L_{\cnum}/L\lrarr M\lrarr X/L_{\cnum}\lrarr 0.
\end{equation}
Once we fix isomorphisms $L_{\cnum}/L\simeq \real$
and $X/L_{\cnum}\simeq\cnum$,
a choice of linear mini-complex coordinate system $(t,w)$ as above
induces a splitting of the exact sequence
{\rm(\ref{eq;21.8.10.1})}.
\hfill\qed
\end{rem}

\subsection{Spaces}

Let $X$ be a $2$-dimensional $\cnum$-vector space.
\index{space $X$}
We take a $\cnum$-linear isomorphism
$X\simeq\cnum^2=\{(z,w)\}$,
and consider the hyperk\"ahler metric
$g_X:=dz\,d\zbar+dw\,d\wbar$.
We consider the $\real$-subspace
$L=\{(s,0)\,|\,s\in\real\}\subset X$
with the natural orientation.
We set $M:=X/L$.
Let $t:=\Image(z)$.
We obtain the mini-complex coordinate system $(t,w)$ on $M$.
\index{space $M$}

We take $T>0$.
We consider the $\seisuu$-action $\kappa$ on $M$
given by $\kappa_n(t,w)=(t+Tn,w)$ $(n\in\seisuu)$.
The quotient space is denoted by $\nbigm$.
\index{space $\nbigm$}
\index{action $\kappa$}

\subsection{Twistor family of complex structures}
\label{subsection;17.10.3.11}

For any complex number $\lambda$,
we have the complex structure on $X$
whose twistor parameter is $\lambda$.
We use the notation $X^{\lambda}$
to denote the complex manifold $X$
equipped with the complex structure corresponding to $\lambda$.
\index{complex manifold $X^{\lambda}$}
We define the complex coordinate system
$(\xi,\eta)$ on $X^{\lambda}$
as follows:
\begin{equation}
\label{eq;17.10.3.12}
 \xi=z+\lambda\wbar,
\quad
 \eta=w-\lambda\zbar.
\end{equation}
\index{complex coordinate system $(\xi,\eta)$}
The inverse transform is described as follows:
\[
 z=\frac{1}{1+|\lambda|^2}
 (\xi-\lambda\etabar),
\quad
 w=\frac{1}{1+|\lambda|^2}
 (\eta+\lambda\xibar).
\]
The metric $dz\,d\zbar+dw\,d\wbar$
is equal to
$(1+|\lambda|^2)^{-2}
 (d\xi\,d\xibar+d\eta\,d\etabar)$.
We obtain the following relations
of complex vector fields:
\[
 \del_{\xibar}
=\frac{1}{1+|\lambda|^2}
 (\del_{\zbar}+\lambda\del_w),
\quad
 \del_{\etabar}
=\frac{1}{1+|\lambda|^2}
 (\del_{\wbar}-\lambda\del_z),
\]
\[
 \del_{\xi}
=\frac{1}{1+|\lambda|^2}
 (\del_z+\lambdabar\del_{\wbar}),
\quad
 \del_{\eta}
=\frac{1}{1+|\lambda|^2}
 (\del_w-\lambdabar\del_{\zbar}).
\]

The $\real$-subspace $L$ and the $\seisuu$-action $\kappa$
are described as follows 
in terms of the coordinate system $(\xi,\eta)$:
\[
 L=\bigl\{s(1,-\lambda)\,\big|\,s\in\real\bigr\},
\]
\[
\kappa_n(\xi,\eta)=
 (\xi+\sqrt{-1}Tn,\eta+\lambda\sqrt{-1}Tn)
=(\xi,\eta)
+nT\cdot(\sqrt{-1},\lambda\sqrt{-1}).
\]

\begin{rem}
\label{rem;21.8.10.3}
Recall that the space of the twistor parameters is $\proj^1$.
Let $\proj^1=\cnum_{\lambda}\cup\cnum_{\mu}$ be the covering
where $\lambda$  and $\mu$ are related by $\lambda\mu=1$.
Let $X^{\dagger\mu}$ denote the complex manifold $X$
equipped with the complex structure corresponding to $\mu$.
\index{complex manifold $X^{\dagger\mu}$}
Let $M^{\dagger\mu}$ denote the mini-complex manifold
obtained as $X^{\dagger\,\mu}/L$. 
\index{mini-complex manifold $M^{\dagger\mu}$}
We have $X^{\dagger\mu}=X^{\lambda}$ 
and $M^{\dagger\mu}=M^{\lambda}$ if $\lambda\mu=1$.
At $\mu=0$, i.e., $\lambda=\infty$,
the complex structure of $X^{\dagger\,0}$ is
induced by $(z^{\dagger},w^{\dagger})=(-\zbar,\wbar)$.
\index{complex coordinate system $(z^{\dagger},w^{\dagger})$}
We note that $L=\{(s,0)\}$ with respect to $(z^{\dagger},w^{\dagger})$,
but that the orientations of $L$ induced by $z$ and $z^{\dagger}=-\zbar$
are mutually reversed.
The complex structure of $X^{\dagger\,\mu}$
is induced by
\index{complex coordinate system $(\xi^{\dagger},\eta^{\dagger})$}
\[
(\xi^{\dagger},\eta^{\dagger})=
(z^{\dagger}+\mu\wbar^{\dagger},w^{\dagger}-\mu\zbar^{\dagger})
=(\lambda^{-1}(w-\lambda\zbar),\lambda^{-1}(z+\lambda\wbar))
=(\lambda^{-1}\eta,\lambda^{-1}\xi).
\]
We shall explain how to obtain 
difference modules from monopoles
behind which the complex coordinate system
$(\xi,\eta)$ is implicitly used.
For $\lambda=\mu^{-1}\neq 0$,
we may also obtain another difference module
from the monopole 
similarly by using $(\xi^{\dagger},\eta^{\dagger})$.
It is analogue to the situation that
a harmonic bundle $(E,\delbar_E,\theta,h)$ on a complex manifold $Y$ induces
a $\lambda$-flat bundle
$(E,\delbar_E+\lambda\theta^{\dagger}+\lambda\del_E+\theta)$ on $Y$
and a $\mu$-flat bundle
$(E,\del_E+\mu\theta+\mu\delbar_E+\theta^{\dagger})$
on $Y^{\dagger}$ which is the conjugate of $Y$.
The analogy can be enhanced to the more precise correspondence
by the Nahm transforms between
wild harmonic bundles on $\proj^1$ and periodic monopoles,
which will be explained elsewhere. 
\hfill\qed
\end{rem}

\subsection{Family of mini-complex structures}

Corresponding to the complex structure of $X^{\lambda}$,
we obtain the mini-complex structures on $M$ and $\nbigm$.
(See \S\ref{subsection;17.10.3.1}.)
The $3$-dimensional manifolds with mini-complex structure
are denoted by $M^{\lambda}$ and $\nbigm^{\lambda}$.
\index{mini-complex manifold $M^{\lambda}$}
\index{mini-complex manifold $\nbigm^{\lambda}$}

In the case $\lambda=0$,
we use the mini-complex coordinate system
$(t,w)=(\Image(z),w)$ on $M^0$.
It induces local mini-complex coordinate
systems on $\nbigm^0$.
\index{mini-complex coordinate system $(t,w)$}

We shall introduce two mini-complex coordinate systems
$(t_i,\beta_i)$ $(i=0,1)$ on $M^{\lambda}$,
which induce mini-complex local coordinate systems
on $\nbigm^{\lambda}$.
For that purpose, we shall introduce two complex coordinate systems 
$(\alpha_i,\beta_i)$ $(i=0,1)$ on $X^{\lambda}$
from which $(t_i,\beta_i)$ are induced.
\begin{rem}
The coordinate systems
$(t_i,\beta_i)$ depend on $\lambda$,
but we omit to denote the dependence
to simplify the description.
\hfill\qed
\end{rem}

\subsection{The mini-complex coordinate system $(t_0,\beta_0)$}
\label{subsection;17.10.4.1}

Let $(\alpha_0,\beta_0)$ be the complex coordinate system
of $X^{\lambda}$ given by the following relation:
\index{complex coordinate system $(\alpha_0,\beta_0)$}
\[
 (\xi,\eta)
=\alpha_0(1,-\lambda)
+\beta_0(\lambdabar,1)
=(\alpha_0+\lambdabar\beta_0,\beta_0-\lambda\alpha_0).
\]
The inverse transform is described as follows:
\[
 \alpha_0=
 \frac{\xi-\lambdabar\eta}{1+|\lambda|^2},
\quad\quad
 \beta_0=
 \frac{\eta+\lambda\xi}{1+|\lambda|^2}.
\]
It is easy to check
$d\alpha_0\,d\alphabar_0
+d\beta_0\,d\betabar_0
=dz\,d\zbar+dw\,d\wbar$.

The $\real$-subspace $L$ 
and the $\seisuu$-action $\kappa$
are described as follows  
in terms of $(\alpha_0,\beta_0)$:
\begin{equation}
\label{eq;17.9.29.10}
 L=\{(s,0)\,|\,s\in\real\},
\end{equation}
\begin{equation}
 \kappa_n(\alpha_0,\beta_0)
=(\alpha_0,\beta_0)
+nT\cdot
 \Bigl(
 \frac{1-|\lambda|^2}{1+|\lambda|^2}
 \sqrt{-1},
 \frac{2\lambda\sqrt{-1}}{1+|\lambda|^2}
 \Bigr). 
\end{equation}

We set $t_0:=\Image(\alpha_0)$.
Because $L$ is described as (\ref{eq;17.9.29.10}),
$(t_0,\beta_0)$
is a mini-complex coordinate system of $M^{\lambda}$.
\index{mini-complex coordinate system $(t_0,\beta_0)$}
The induced $\seisuu$-action $\kappa$ on $M^{\lambda}$
is described as follows in terms of $(t_0,\beta_0)$:
\[
 \kappa_n(t_0,\beta_0)
=(t_0,\beta_0)
+nT\cdot
 \Bigl(
 \frac{1-|\lambda|^2}{1+|\lambda|^2},
 \frac{2\lambda\sqrt{-1}}{1+|\lambda|^2}
 \Bigr). 
\]
Clearly, $\kappa_n$
is given along the integral curve of 
the real vector field:
\begin{equation}
\label{eq;17.9.29.30}
\del_t=
 \frac{1-|\lambda|^2}{1+|\lambda|^2}
 \del_{t_0}
+\frac{2\lambda}{1+|\lambda|^2}
 \sqrt{-1}\del_{\beta_0}
-\frac{2\lambdabar}{1+|\lambda|^2}
 \sqrt{-1}\del_{\betabar_0}.
\end{equation}

\begin{rem}
If $(t'_0,\beta'_0)$ is another mini-complex coordinate system 
of $M^{\lambda}$
such that $dt_0'\,dt_0'+d\beta'_0\,d\betabar'_0$,
we obtain $(t_0',\beta'_0)=(t_0,a\beta_0)$,
where $a$ is a complex number such that $|a|=1$.
\hfill\qed
\end{rem}

\begin{rem}
For the twistor parameter $\mu\in\cnum_{\mu}\subset\proj^1$,
we obtain $(\alpha_0^{\dagger},\beta_0^{\dagger})$
from $(\xi^{\dagger},\mu^{\dagger})$ as above:
\index{complex coordinate system $(\alpha_0^{\dagger},\beta_0^{\dagger})$}
\[
 (\alpha_0^{\dagger},\beta_0^{\dagger})
 =\frac{1}{1+|\lambda|^2}(\xi^{\dagger}-\mubar\eta^{\dagger},
 \eta^{\dagger}+\mu\xi^{\dagger})
\]
Because the orientations of $L$ induced by $z$ and $z^{\dagger}=-\zbar$
are mutually reversed,
the induced mini-complex coordinate system of $M^{\dagger\mu}$ is
$(-\Image(\alpha^{\dagger}_0),\beta^{\dagger}_0)$.
If $\mu=\lambda^{-1}\neq 0$,
we have
$(\alpha_0^{\dagger},\beta_0^{\dagger})
=\bigl(-\alpha_0,(\lambda^{-1}\lambdabar)\beta_0\bigr)$,
and hence
$(t_0^{\dagger},\beta_0^{\dagger})=
 (-\Image(\alpha^{\dagger}_0),\beta^{\dagger}_0)=
\bigl(t_0,(\lambda^{-1}\lambdabar)\beta_0\bigr)$.
\index{mini-complex coordinate system $(t_0^{\dagger},\beta_0^{\dagger})$}
\hfill\qed
\end{rem}

\subsection{The mini-complex coordinate system
   $(t_1,\beta_1)$}
\label{subsection;17.10.3.12}

Let $(\alpha_1,\beta_1)$ be 
the complex coordinate system
of $X^{\lambda}$
determined by the relation
$(\xi,\eta)=
 \alpha_1(1,-\lambda)+\beta_1(0,1)$.
 \index{complex coordinate system $(\alpha_1,\beta_1)$}
The transformations are described as follows:
\[
\left\{
\begin{array}{l}
 \xi=\alpha_1,
\\
 \eta=\beta_1-\lambda\alpha_1,
\end{array}
\right.
\quad\quad
\left\{
 \begin{array}{l}
 \alpha_1=\xi,\quad
  \\
 \beta_1=\eta+\lambda\xi.
 \end{array}
\right.
\]
The $\real$-subspace $L$ and the $\seisuu$-action $\kappa$
are described as follows in terms of
$(\alpha_1,\beta_1)$:
\begin{equation}
\label{eq;17.9.29.11}
 L=\{(s,0)\,|\,s\in\real\},
\end{equation}
\begin{equation}
\kappa_n(\alpha_1,\beta_1)
=(\alpha_1,\beta_1)
+nT(\sqrt{-1},2\lambda\sqrt{-1}).
\end{equation}

We set $t_1:=\Image\alpha_1$.
Because $L$ is described as (\ref{eq;17.9.29.11}),
we obtain the mini-complex coordinate system
$(t_1,\beta_1)$ on $M^{\lambda}$.
\index{mini-complex coordinate system $(t_1,\beta_1)$}
The induced $\seisuu$-action $\kappa$
on $M^{\lambda}$
is described as follows:
\[
 \kappa_n(t_1,\beta_1)=(t_1,\beta_1)
+nT(1,2\lambda\sqrt{-1}).
\]

\begin{rem}
For the twistor parameter $\mu\in\cnum_{\mu}\subset\proj^1$,
we obtain the complex coordinate system
$(\alpha_1^{\dagger},\beta_1^{\dagger})
=(\xi^{\dagger},\eta^{\dagger}+\mu\xi^{\dagger})$
of $X^{\dagger\,\mu}$,
and the mini-complex coordinate system
 $(t_1^{\dagger},\beta_1^{\dagger})
 =(-\Image(\alpha_1)^{\dagger},\beta_1^{\dagger})$
of $M^{\dagger\,\mu}$.
\index{complex coordinate system $(\alpha_1^{\dagger},\beta_1^{\dagger})$}
\index{mini-complex coordinate system $(t_1^{\dagger},\beta_1^{\dagger})$}
If $\mu=\lambda^{-1}$,
we have
$(\alpha_1^{\dagger},\beta_1^{\dagger})
=(\lambda^{-1}\eta,\lambda^{-2}(\eta+\lambda\xi))$,
and hence
$(t_1^{\dagger},\beta_1^{\dagger})
=(t_1-\Image(\lambda^{-1}\beta_1),\lambda^{-2}\beta_1)$.
We obtain the splittings of $M^{\lambda}=M^{\dagger\mu}$
induced by $(t_1,\beta_1)$ and $(t_1^{\dagger},\beta_1^{\dagger})$
as in Remark {\rm\ref{rem;21.8.10.2}},
and they are different.
Therefore, we obtain two difference modules
depending on the choice of $(\xi,\eta)$ or $(\xi^{\dagger},\eta^{\dagger})$
as mentioned in Remark {\rm\ref{rem;21.8.10.3}}. 
There are essentially only these two choices.
See {\rm\S\ref{subsection;21.8.10.4}}.
\hfill\qed
\end{rem}

\begin{rem}
If $\lambda=0$,
we have $(t_0,\beta_0)=(t_1,\beta_1)=(t,w)$.
\hfill\qed
\end{rem}

\subsection{Coordinate change}

We have the following relation:
\[
 \left\{
 \begin{array}{l}
 \alpha_1=\alpha_0+\lambdabar\beta_0\\
 \beta_1=(1+|\lambda|^2)\beta_0,
 \end{array}
 \right.
\quad\quad
 \left\{
 \begin{array}{l}
 \alpha_0=\alpha_1-(1+|\lambda|^2)\lambdabar\beta_1\\
 \beta_0=(1+|\lambda|^2)^{-1}\beta_1.
 \end{array}
 \right.
\]
Hence, we have the following relation:
\[
 \left\{
 \begin{array}{l}
 t_1=t_0+\Image(\lambdabar\beta_0)\\
 \beta_1=(1+|\lambda|^2)\beta_0,
 \end{array}
 \right.
\quad\quad
 \left\{
 \begin{array}{l}
 t_0=t_1-(1+|\lambda|^2)^{-1}\Image(\lambdabar\beta_1)\\
 \beta_0=(1+|\lambda|^2)^{-1}\beta_1.
 \end{array}
 \right.
\]
We obtain the following relation of vector fields:
\begin{equation}
\label{eq;17.9.29.31}
 \del_{t_1}=\del_{t_0},
\quad
 \del_{\betabar_1}
=\frac{\lambda}{1+|\lambda|^2}
 \frac{1}{2\sqrt{-1}}
 \del_{t_0}
+\frac{1}{1+|\lambda|^2}
 \del_{\betabar_0}.
\end{equation}

\subsection{Compactification}
\label{subsection;21.8.21.10}

Set $\Mbar^{\lambda}:=\real_{t_1}\times\proj^1_{\beta_1}$,
which is equipped with the natural mini-complex structure.
\index{mini-complex manifold $\Mbar^{\lambda}$}
We have the $\seisuu$-action $\kappa$
on $\Mbar^{\lambda}$
given by
$\kappa_n(t_1,\beta_1)=(t_1+Tn,\beta_1+2\sqrt{-1}\lambda Tn)$.
The quotient space is denoted by $\nbigmbar^{\lamda}$.
\index{mini-complex manifold $\nbigmbar^{\lambda}$}
It is a compactification of $\nbigmlambda$,
and equipped with the naturally induced mini-complex structure.
We set 
$H^{\lambda\cov}_{\infty}:=
 \Mbar^{\lambda}\setminus M^{\lambda}$
and
$H^{\lambda}_{\infty}:=
 \nbigmbar^{\lambda}\setminus\nbigm^{\lambda}$.
Let $\varpi^{\lambda}:\Mbar^{\lambda}\lrarr\nbigmbar^{\lambda}$
denote the projection.
\index{space $H^{\lambda}_{\infty}$}
\index{space $H^{\lambda\cov}_{\infty}$}
\index{projection $\varpi^{\lambda}$}
For $\lambda_1\neq \lambda_2$,
we have $\nbigmbar^{\lambda_1}\neq\nbigmbar^{\lambda_2}$
though $\nbigm^{\lambda_1}=\nbigm^{\lambda_2}$
as $C^{\infty}$-manifolds.

\begin{rem}
If $\lambda=\mu^{-1}\neq 0$,
we may obtain another natural compactification
induced by $(t_1^{\dagger},\beta_1^{\dagger})$. 
The compactifications depend on the splittings
induced by $(t_1,\beta_1)$ and $(t_1^{\dagger},\beta_1^{\dagger})$
(see also Remark {\rm\ref{rem;21.8.10.2}}).
\hfill\qed
\end{rem}

\subsection{Mini-holomorphic bundles associated with monopoles}
\label{subsection;21.8.23.2}

Let $\varphi:U\lrarr \nbigm$ be a local diffeomorphism.
We regard $U$ as the Riemannian manifold
whose metric is obtained as the pull back of
$dt\,dt+dw\,d\wbar$.
For $\lambda\in\cnum$,
let $U^{\lambda}$ denote the $3$-dimensional manifold
with the mini-complex structure
obtained as the pull back of the mini-complex structure of
$\nbigm^{\lambda}$.

For any monopole  $(E,h,\nabla,\phi)$ on $U$,
we obtain the mini-holomorphic bundle
on $U^{\lambda}$ underlying $(E,h,\nabla,\phi)$
for each $\lambda$
(see \S\ref{subsection;16.9.20.20}),
which we denote by $(E^{\lambda},\delbar_{E^{\lambda}})$.
\index{mini-holomorphic bundle $(E^{\lambda},\delbar_{E^{\lambda}})$}

\subsubsection{Compatibility with the dimensional reduction from 4D to 3D}
\label{subsection;20.7.30.51}

Let $U$ be as above.
We set $\Utilde:=\real_s\times U$
on which we consider the Riemannian metric
$ds\,ds+dt\,dt+dw\,d\wbar$.
Let $\Utilde^{0}$
denote the complex manifold
whose complex structure is given by
local coordinate systems $(z,w)=(s+\sqrt{-1}t,w)$.
For any $\lambda\in\cnum$,
let $\Utilde^{\lambda}$
denote the complex manifold
whose complex structure
is given by the local complex coordinate systems
$(\xi,\eta)$ as in \S\ref{subsection;17.10.3.11}.

Let $(E,h,\nabla,\phi)$ be a monopole on $U$.
We obtain the induced instanton
$(\Etilde,\htilde,\nablatilde)$ on $\Utilde$
as in \S\ref{subsection;17.10.3.10}.
Let $(\Etilde^{\lambda},\delbar_{\Etilde^{\lambda}})$
denote the holomorphic bundle on
$\Utilde^{\lambda}$
underlying the instanton.
The following lemma is clear by the constructions,
or can be checked by a direct computation.
\begin{lem}
$(\Etilde^{\lambda},\delbar_{\Etilde^{\lambda}})$
is equal to the holomorphic bundle
induced by 
the mini-holomorphic bundle $(E^{\lambda},\delbar_{E^{\lambda}})$
as in {\rm \S\ref{subsection;13.11.29.2}}.
\hfill\qed
\end{lem}

\section{$\nbigo_{\nbigm^{\lambda}}$-modules and $\lambda$-connections}
\label{section;21.8.12.101}

\subsection{Dimensional reduction from
  $\nbigo_{\nbigm^{\lambda}}$-modules to
$\lambda$-flat bundles}
\label{subsection;21.8.11.21}

\subsubsection{Setting}
\label{subsection;21.8.11.20}

Let $\Psi:\nbigm=S^1_T\times\cnum_w\lrarr\cnum_w$ denote the projection.
Let $U_w\lrarr \cnum_w$ be a local diffeomorphism.
Let $U$ denote the fiber product of $U_w$ and $\nbigm$ over $\cnum_w$.
The induced morphism $U\lrarr U_w$ is also denoted by $\Psi$.
There exists a natural isomorphism $U\simeq S^1_T\times U_w$
which is equipped with the naturally defined $S^1_T$-action.
We obtain the naturally induced local diffeomorphism $U\lrarr\nbigm$.
We regard $U$ as a Riemannian manifold
whose metric is obtained as the pull back of $dt\,dt+dw\,d\wbar$.
For $\lambda\in\cnum$,
let $U^{\lambda}$ denote the $3$-dimensional manifold
equipped with the mini-complex structure
obtained as the pull back of the mini-complex structure of
$\nbigm^{\lambda}$.

\vspace{.1in}
We shall explain that $\lambda$-flat bundles on $U_w$
are the dimensional reduction of mini-holomorphic bundles on $U^{\lambda}$
(Corollary \ref{cor;21.8.11.22}),
as the special case of a more general equivalence
including the non-integrable case (Lemma \ref{lem;21.8.12.2}).

\subsubsection{Some vector fields and forms}
We have the local complex coordinates on $U_w$ induced by $w$,
which induce the differential forms
$dw$ and $d\wbar$,
and the vector fields $\del_w$ and $\del_{\wbar}$
globally defined on $U_w$.
We have the local mini-complex coordinate systems on $U^{\lambda}$
induced by $(t_1,\beta_1)$,
and
we obtain the globally defined differential forms
$dt_1$, $d\beta_1$ and $d\betabar_1$,
and the globally defined complex vector fields
$\del_{t_1}$, $\del_{\beta_1}$ and $\del_{\betabar_1}$ on $U$.
Note that $\Psi$ is described as
$\Psi(t_1,\beta_1)=(1+|\lambda|^2)^{-1}(\beta_1-2\sqrt{-1}\lambda t_1)$
in terms of the coordinate systems $(t_1,\beta_1)$
and $w$.
The tangent map of $\Psi$ is described as follows:
\begin{equation}
\label{eq;17.10.1.1}
 T\Psi(\del_{t_1})
=\frac{-2\sqrt{-1}\lambda }{1+|\lambda|^2}\del_w
+\frac{2\sqrt{-1}\cdot \lambdabar}{1+|\lambda|^2}\del_{\wbar},
\quad\quad
 T\Psi(\del_{\betabar_1})
=\frac{1}{1+|\lambda|^2}
 \del_{\wbar}.
\end{equation}
Similarly, 
we have the local mini-complex coordinate systems on $U^{\lambda}$
induced by $(t_0,\beta_0)$,
the globally defined differential forms
$dt_0$, $d\beta_0$ and $d\betabar_0$,
and the globally defined complex vector fields
$\del_{t_0}$, $\del_{\beta_0}$ and $\del_{\betabar_0}$ on $U$.
We also have the local mini-complex coordinate systems on $U^0$
induced by $(t,w)$,
the globally defined differential forms
$dt$, $dw$ and $d\wbar$,
and the globally defined complex vector fields
$\del_t$, $\del_{w}$ and $\del_{\wbar}$ on $U$.

\subsubsection{A general equivalence}
\label{subsection;21.8.12.3}
Let $V$ be a $C^{\infty}$-vector bundle on $U_w$.
Recall that
a $\lambda$-connection of $V$ in the $C^{\infty}$-sense
is a linear differential operator
$\DDlambda:C^{\infty}(U_w,V)\lrarr
C^{\infty}(U_w,V\otimes\Omega_{U_w}^1)$
such that
$\DDlambda(gs)=(\lambda\del+\delbar)g\cdot s
+f\DDlambda(s)$
for any $g\in C^{\infty}(U_w)$
and $s\in C^{\infty}(U_w,V)$.
\index{$\lambda$-connection ($C^{\infty}$)}
(See \cite[\S2.2]{Mochizuki-KHII}.)
It is called flat if $\DDlambda\circ\DDlambda=0$.

Let $\del_{V,\wbar}$ (resp. $\DDlambda_w$)
denote the differential operator of $V$
induced by $\del_{\wbar}$ (resp. $\del_w$)
and $\DDlambda$.
\index{operators $\del_{V,\wbar}$, $\DDlambda_w$}
The flatness of $\DDlambda$ is equivalent to
$[\del_{V,\wbar},\DDlambda_w]=0$.

We set $\Vtilde=\Psi^{-1}(V)$.
We say that a linear differential operator
$\delbar_{\Vtilde}:
C^{\infty}(U^{\lambda},\Vtilde)
\lrarr
C^{\infty}(U^{\lambda},\Vtilde\otimes\Omega^{0,1}_{U^{\lambda}})$
satisfies a mini-complex Leibniz rule
if
$\delbar_{\Vtilde}(fu)
=\delbar_{U^{\lambda}}(f)u+f\delbar_{\Vtilde}(u)$
for any $f\in C^{\infty}(U)$
and $u\in C^{\infty}(U,\Vtilde)$.
\index{mini-complex Leibniz rule}

Let $\DDlambda$ be a $\lambda$-connection of $V$
in the $C^{\infty}$-sense.
Noting the relation (\ref{eq;17.10.1.1}),
we obtain the linear differential operators
$\del_{\Vtilde,\betabar_1}$
and $\del_{\Vtilde,t_1}$
of $\Vtilde$
such that the following holds
for any $s\in C^{\infty}(U_w,V)$ and $f\in C^{\infty}(U)$:
\begin{equation}
\label{eq;21.8.12.30}
 \del_{\Vtilde,\betabar_1}(f \Psi^{-1}(s))
=\del_{\betabar_1}f\cdot \Psi^{-1}(s)
+\frac{1}{1+|\lambda|^2}
 f \Psi^{-1}(\del_{V,\wbar}s).
\end{equation}
\begin{equation}
\label{eq;21.8.12.31}
 \del_{\Vtilde,t_1}(f \Psi^{-1}(s))
=\del_{t_1}f\cdot \Psi^{-1}(s)
+
\frac{1}{1+|\lambda|^2}
 f \Psi^{-1}\bigl(
 -2\sqrt{-1}\DDlambda_w s
+2\sqrt{-1}\cdot\lambdabar\del_{V,\wbar}s
 \bigr).
\end{equation}
We define
$\delbar_{\Vtilde}:C^{\infty}(U,\Vtilde)\lrarr
C^{\infty}(U,\Vtilde\otimes\Omega^{0,1}_{U^{\lambda}})$
by
$\delbar_{\Vtilde}(u)=
\del_{\Vtilde,\betabar_1}(u)d\betabar_1
+\del_{\Vtilde,t_1}(u)\,dt_1$
for any $u\in C^{\infty}(U,\Vtilde)$.
By the relation (\ref{eq;17.10.1.1}),
$\delbar_{\Vtilde}$ satisfies
the mini-complex Leibniz rule.
Note that
$(\Vtilde,\delbar_{\Vtilde})$
is $S^1_T$-equivariant.
By the construction, the following holds:
\begin{equation}
\label{eq;21.8.12.1}
 \bigl[\del_{\Vtilde,\betabar_1},\del_{\Vtilde,t_1} \bigr]
 =\frac{-2\sqrt{-1}}{(1+|\lambda|^2)^2}
 \Psi^{-1}\bigl(
 \bigl[\del_{\wbar},\DD^{\lambda}_w\bigr]
 \bigr).
\end{equation}

\begin{lem}
\label{lem;21.8.12.2}
 The above procedure induces an equivalence between
the following objects.
\begin{itemize}
 \item Vector bundles $V$ on $U_w$
       equipped with a $\lambda$-connection
       in the $C^{\infty}$-sense.
 \item $S^1_T$-equivariant vector bundles
       $\Vtilde$ on $U$
       equipped with
       an $S^1_T$-equivariant linear differential operator
       $\delbar_{\Vtilde}:
       C^{\infty}(U^{\lambda},\Vtilde)\lrarr
       C^{\infty}(U^{\lambda},\Vtilde\otimes\Omega^{0,1}_{U^{\lambda}})$
       satisfying the mini-complex Leibniz rule.
\end{itemize}
\end{lem}
\pf
Let us indicate the inverse construction.
Let $(\Vtilde,\delbar_{\Vtilde})$ be
as in the statement of Lemma \ref{lem;21.8.12.2}.
Let $\del_{\Vtilde,\betabar_1}$ (resp. $\del_{\Vtilde,t_1}$)
denote the $S^1_T$-equivariant differential operator
on $\Vtilde$ induced by
$\delbar_{\Vtilde}$ and $\del_{\betabar_1}$
(resp. $\del_{t_1}$).
There exists a $C^{\infty}$-bundle $V$ on $U_w$
with an $S^1_T$-equivariant isomorphism
$\Psi^{-1}(V)\simeq \Vtilde$.
We obtain the differential operators $\del_{V,\wbar}$
and $\DDlambda_w$ of $V$
by the following condition for any $s\in C^{\infty}(U_w,V)$:
\[
 \Psi^{-1}(\del_{V,\wbar}(s))=
 (1+|\lambda|^2)\del_{\Vtilde,\betabar_1}\Psi^{-1}(s),
\]
\[
 \Psi^{-1}(\DDlambda_ws)
 =\frac{\sqrt{-1}}{2}(1+|\lambda|^2)
 \del_{\Vtilde,t_1}\Psi^{-1}(s)
 +(1+|\lambda|^2)\lambdabar\del_{\Vtilde,\betabar_1}
 \Psi^{-1}(s).
\]
Then, by the relation (\ref{eq;17.10.1.1}),
we have
$\del_{V,\wbar}(gs)=\del_{\wbar}(g)s+g\del_{V,\wbar}(s)$
and
$\DDlambda_w(gs)=\lambda\del_w(g)s+g\DDlambda_w(s)$
for any $g\in C^{\infty}(U_w)$ and $s\in C^{\infty}(U_w,V)$.
Hence, $\del_{V,\wbar}$ and $\DDlambda_w$
induces a $\lambda$-connection of $V$.
Two constructions are mutually inverse.
\hfill\qed

\subsubsection{Mini-holomorphic bundles and flat $\lambda$-connections}

As a consequence of Lemma \ref{lem;21.8.12.2}
and (\ref{eq;21.8.12.1}),
we obtain the following.
\begin{cor}
 \label{cor;21.8.11.22}
The construction in {\rm\S\ref{subsection;21.8.12.3}}
induces an equivalence between
$\lambda$-flat bundles on $U_w$
and 
$S^1_T$-equivariant mini-holomorphic vector bundles
on $U$. 
\hfill\qed 
\end{cor}

\subsubsection{$\lambda$-flat bundles of infinite rank}
\label{subsection;21.8.13.31}

Let $E$ be a $C^{\infty}$-bundle on $U^{\lambda}$
equipped with a differential operator
$\delbar_E:C^{\infty}(U^{\lambda},E)
\lrarr C^{\infty}(U^{\lambda},E\otimes\Omega^{0,1}_{U^{\lambda}})$
satisfying the mini-complex Leibniz rule.
Let $E_{C^{\infty}}$ denote the sheaf of $C^{\infty}$-sections of
$E$.
Let $\nbigc^{\infty}_{U_w}$ denote the sheaf of
$C^{\infty}$-functions on $U_w$.
We obtain the $\nbigc^{\infty}_{U_w}$-module
$V_1:=\Psi_!(E_{C^{\infty}})$,
where $\Psi_!$ denotes the proper push-forward with respect to $\Psi$.
We define 
$\del_{V_1,\wbar}$
and $\DDlambda_{V_1,w}$ on the sheaf $V_1$
by
\[
\del_{V_1,\wbar}(u)=(1+|\lambda|^2)\del_{E,\betabar_1}(u),
\quad
\DDlambda_{V_1,w}(u)=
\frac{\sqrt{-1}}{2}(1+|\lambda|^2)\del_{E,t_1}(u)
+\lambdabar(1+|\lambda|^2)\del_{E,\betabar_1}(u)
\]
for sections $u$ of $V_1$.
We can easily observe that
$\del_{V_1,\wbar}$ and $\DDlambda_{V_1,w}$
determines a $\lambda$-connection $\DDlambda_{V_1}$ of $V_1$
in the $C^{\infty}$-sense.
If $(E,\delbar_E)$ is a mini-holomorphic bundle,
then $\DDlambda_{V_1}$ is a flat $\lambda$-connection.

\subsubsection{Remark}

This analogy between
mini-holomorphic bundles on $U^{\lambda}$
and $\lambda$-flat bundles on $U_w$
is one of the motivations of this study,
as explained in {\rm\S\ref{section;21.8.13.20}}.
(See also \S\ref{subsection;21.8.13.101}.)
It is also useful when we study the property of
the mini-holomorphic bundles
underlying monopoles induced by harmonic bundles.
(See {\rm Corollary \ref{cor;21.8.12.6}}
and the computations in \S\ref{subsection;17.10.25.130}
and \S\ref{subsection;21.8.13.22} below.)
We shall generalize this analogy to the case
where $U_w$ is a ramified covering
of a neighbourhood of $\infty\in\proj^1$
(Corollary \ref{cor;21.8.12.20} and Corollary \ref{cor;21.8.13.24}).
We study the formal and ramified version in \S\ref{subsection;17.10.7.2}.
Note that in the formal version
we can remove the $S^1_T$-equivariance condition
(see Proposition \ref{prop;20.7.18.1}).
We obtain the filtered version
(Proposition \ref{prop;20.7.20.131} and
Proposition \ref{prop;20.7.24.10}) as a refinement.
It is also useful to study the non-integrable case
(Lemma \ref{lem;21.8.12.2}, Proposition \ref{prop;21.8.12.5},
Lemma \ref{lem;21.8.12.10} and Lemma \ref{lem;21.8.12.34}).

We recall that 
in the study of harmonic bundles
we understand the asymptotic behaviour of
harmonic bundles
through the filtered extension of the underlying $\lambda$-flat bundles,
pioneered by Simpson \cite{Simpson90},
and further studied in 
\cite{mochi2, Mochizuki-wild}.
Similarly, it is our viewpoint 
to understand the asymptotic behaviour of monopoles
through the filtered extension of the underlying mini-holomorphic bundles,
studied in \S\ref{section;20.7.25.10}.
Moreover, we would like to apply the results for
filtered extension of $\lambda$-flat bundles underlying
harmonic bundles
to study the filtered extension of the mini-holomorphic bundles.
That is the main reason to pursue this analogy.
(See \S\ref{subsection;21.9.16.2} for an outline of the argument.)

\subsection{Comparison of some induced operators}
\label{subsection;21.8.12.100}

Let $U_w$ and $U$ be as in \S\ref{subsection;21.8.11.20}.
Let $V$ be a $C^{\infty}$-bundle on $U_w$
with a Hermitian metric $h$.
Let $\delbar_V$ be a holomorphic structure of $V$,
and let $\theta\in C^{\infty}(U_w,\End(V)\otimes\Omega^{1,0})$
which does not necessarily satisfy
$\delbar_V\theta=0$.
Let $h$ be a Hermitian metric of $V$.
We obtain the Chern connection $\delbar_V+\del_V$,
and the adjoint
$\theta^{\dagger}\in C^{\infty}(U_w,\End(V)\otimes\Omega^{0,1})$
of $\theta$.
Let $\del_{V,\wbar}$ (resp. $\del_{V,w}$)
denote the differential operator of $V$
induced by $\delbar_V$ and $\del_{\wbar}$
(resp. $\del_V$ and $\del_w$).
Let $f$ and $f^{\dagger}$ be the endomorphisms of $V$
determined by
$\theta=f\,dw$ and $\theta^{\dagger}=f^{\dagger}d\wbar$,
respectively.
We set $\Vtilde=\Psi^{-1}(V)$ on $U$
with the metric $\htilde=\Psi^{-1}(h)$.

On one hand,
we obtain the unitary connection $\nabla$
and an anti-Hermitian endomorphism $\phi$
of $\Vtilde$ satisfying the following condition
for $s\in C^{\infty}(U_w,V)$
as in (\ref{eq;17.10.6.21}):
\[
 \nabla_{\wbar}(\Psi^{-1}(s))=\Psi^{-1}(\del_{V,\wbar}s),
 \quad
  \nabla_{w}(\Psi^{-1}(s))=\Psi^{-1}(\del_{V,w}s),
\]
\[
 \nabla_t(\Psi^{-1}(s))
 =\Psi^{-1}\bigl(
 -\sqrt{-1}(f+f^{\dagger})(s)
 \bigr),
 \quad
 \phi(\Psi^{-1}(s))
 =\Psi^{-1}\bigl(
 (f-f^{\dagger})s
 \bigr).
\]
From $\Vtilde$ with $\nabla$ and $\phi$,
we obtain a linear differential operator
$\delbar^{(1)}_{\Vtilde}:
C^{\infty}(U^{\lambda},\Vtilde)\lrarr
C^{\infty}(U^{\lambda},\Vtilde\otimes\Omega^{0,1}_{U^{\lambda}})$
as in \S\ref{subsection;16.9.20.20}
which satisfies the mini-complex Leibniz rule.

On the other hand, we obtain
the $\lambda$-connection
$\DDlambda
=\delbar_V+\lambda\theta^{\dagger}
+\lambda\del_V+\theta$ of $V$ in the $C^{\infty}$-sense.
From $(V,\DDlambda)$,
we obtain
$\delbar^{(2)}_{\Vtilde}:
C^{\infty}(U^{\lambda},\Vtilde)\lrarr
C^{\infty}(U^{\lambda},\Vtilde\otimes\Omega^{0,1}_{U^{\lambda}})$
as in \S\ref{subsection;21.8.11.21}.

\begin{prop}
\label{prop;21.8.12.5}
 We have $\delbar^{(1)}_{\Vtilde}=\delbar^{(2)}_{\Vtilde}$.
\end{prop}
\pf
Note that
\[
 \nabla_{t_0}=
 \frac{1-|\lambda|^2}{1+|\lambda|^2}\nabla_t
-\frac{2\lambda\sqrt{-1}}{1+|\lambda|^2}\nabla_w
+\frac{2\lambdabar\sqrt{-1}}{1+|\lambda|^2}\nabla_{\wbar},
\]
\[
 \nabla_{\betabar_0}
=\frac{\sqrt{-1}\lambda}{1+|\lambda|^2}\nabla_t
+\frac{\lambda^2}{1+|\lambda|^2}\nabla_w
+\frac{1}{1+|\lambda|^2}\nabla_{\wbar}.
\]
We also obtain the following:
\begin{multline}
 \del^{(1)}_{\Vtilde,\betabar_1}
=\frac{\lambda}{1+|\lambda|^2}
 \frac{1}{2\sqrt{-1}}\del^{(1)}_{\Vtilde,t_0}
+\frac{1}{1+|\lambda|^2}
 \del^{(1)}_{\Vtilde,\betabar_0} \\
=\frac{1}{1+|\lambda|^2}\nabla_{\betabar_0}
-\frac{\lambda}{1+|\lambda|^2}\frac{1}{2}(\phi+\sqrt{-1}\nabla_{t_0}),
\end{multline}
\begin{equation}
 \del^{(1)}_{\Vtilde,t_1}
=\del^{(1)}_{\Vtilde,t_0}
=\nabla_{t_0}-\sqrt{-1}\phi.
\end{equation}
Hence, we obtain the following:
\begin{equation}
\label{eq;17.10.7.310}
 \del^{(1)}_{\Vtilde,\betabar_1}
=\frac{1}{1+|\lambda|^2}
 \Bigl(
 \frac{\lambda\sqrt{-1}}{2}\nabla_t
+\nabla_{\wbar}
-\frac{\lambda}{2}\phi
 \Bigr),
\end{equation}
\begin{equation}
\label{eq;17.10.7.311}
 \del^{(1)}_{\Vtilde,t_1}
=\frac{1-|\lambda|^2}{1+|\lambda|^2}\nabla_t
-\frac{2\lambda\sqrt{-1}}{1+|\lambda|^2}\nabla_w
+\frac{2\lambdabar\sqrt{-1}}{1+|\lambda|^2}\nabla_{\wbar}
-\sqrt{-1}\phi.
\end{equation}
By the construction of $\nabla$ and $\phi$,
we obtain the following:
\[
 \del^{(1)}_{\Vtilde,\betabar_1}\Psi^{-1}(s)
=\frac{1}{1+|\lambda|^2}
 \Psi^{-1}\bigl(
 (\del_{V,\wbar}+\lambda f^{\dagger})s
 \bigr),
\]
\[
 \del^{(1)}_{\Vtilde,t_1}
 \Psi^{-1}(s)
=\Psi^{-1}
 \Bigl(
 \frac{2\lambdabar\sqrt{-1}}{1+|\lambda|^2}
 (\del_{V,\wbar}+\lambda f^{\dagger})s
-\frac{2\sqrt{-1}}{1+|\lambda|^2}
 (\lambda\del_{V,w}+f)s
 \Bigr).
\]
This is equal to $\del^{(2)}_{\Vtilde}$
constructed from
$\DDlambda=\delbar_V+\lambda\theta^{\dagger}+
\lambda\del_V+\theta$
in \S\ref{subsection;21.8.11.21}.
\hfill\qed

\subsubsection{Comparison of mini-holomorphic bundles
induced by harmonic bundles}

Suppose that $(V,\delbar_V,\theta,h)$ is a harmonic bundle on $U_w$.
On one hand,
$(\Vtilde,\htilde,\nabla,\phi)$ is a monopole on $U$
associated with $(V,\delbar_V,\theta,h)$ as in \S\ref{subsection;17.9.29.1},
and $(\Vtilde,\delbar^{(1)}_{\Vtilde})$
is the mini-holomorphic bundle on $U^{\lambda}$
underlying $(\Vtilde,\htilde,\nabla,\phi)$.
On the other hand, $\DDlambda$
is a flat $\lambda$-connection associated with $(V,\delbar_V,\theta,h)$,
and $(\Vtilde,\delbar_{\Vtilde}^{(2)})$
is the mini-holomorphic bundle associated with
$(V,\DDlambda)$
in Corollary \ref{cor;21.8.11.22}.
We obtain the following from Proposition \ref{prop;21.8.12.5}.
\begin{cor}
\label{cor;21.8.12.6}
The associated mini-holomorphic bundles
$(\Vtilde,\delbar^{(1)}_{\Vtilde})$
and
$(\Vtilde,\delbar^{(2)}_{\Vtilde})$ on $U^{\lambda}$
are the same.
\hfill\qed
\end{cor}

\subsection{$\nbigo_{\nbigmbar^{\lambda}}$-modules
and $\lambda$-connections}
\label{subsection;21.8.13.21}

\subsubsection{Setting}
Let $\Psi^{\lambda}:\nbigmbar^{\lambda}\lrarr \proj^1_w$
be the map induced by 
$\Psi^{\lambda}(t_1,\beta_1)=
 (1+|\lambda|^2)^{-1}(\beta_1-2\sqrt{-1}\lambda t_1)$.
\index{map $\Psi^{\lambda}$}
Note that the restriction of $\Psi^{\lambda}$ to $\nbigm^{\lambda}$
is equal to $\Psi$.

Let $U_w$ be any open subset of $\proj^1_w$
such that $\infty\in U_w$.
We set $U^{\lambda}=(\Psi^{\lambda})^{-1}(U_w)\subset\nbigmbar^{\lambda}$,
which is equipped with the mini-complex structure
as an open subset of $\nbigmbar^{\lambda}$.
We shall generalize Lemma \ref{lem;21.8.12.2}
and Corollary \ref{cor;21.8.11.22} to this context.

\begin{rem}
We shall later generalize the construction to the case where
$U_w\lrarr \proj^1$ is a ramified covering around $\infty$.
(See {\rm\S\ref{subsection;21.8.13.1}}.) 
\hfill\qed
\end{rem}

\subsubsection{A general equivalence}
\label{subsection;21.8.12.11}

Let $(V,\delbar_V)$ be a holomorphic vector bundle on $U_w$.
Let $\del_{V,\wbar}$ denote the differential operator of $V$
induced by $\delbar_V$ and $\del_{\wbar}$.
Let $\DDlambda_w$ be a linear differential operator of $V$
satisfying
$\DDlambda_w(fs)=\lambda\del_w(f)\cdot s+f\DDlambda_w(s)$
for any $f\in C^{\infty}(U_w)$ and $s\in C^{\infty}(U_w,V)$.
\index{operators $\del_{\wbar}$, $\DDlambda_w$}
Note that $\DDlambda_w$ is equivalent to
a linear differential operator
$\DD^{\lambda\,(1,0)}:C^{\infty}(U_w,V)
\lrarr C^{\infty}(U_w,V\otimes\Omega^{1,0}_{U_w}(2\infty))$
such that
$\DD^{\lambda\,(1,0)}(fs)
=\lambda\del_{U_w}(f)s+f\DD^{\lambda\,(1,0)}(s)$
for any $f\in C^{\infty}(U_w)$
and $s\in C^{\infty}(U_w,V)$.
\index{operator $\DD^{\lambda\,(1,0)}$}
The restriction of $\delbar_V+\DD^{\lambda(1,0)}$
to $U_w\setminus\{\infty\}$
is a $\lambda$-connection of $V_{|U_w\setminus\{\infty\}}$
in the $C^{\infty}$-sense.

We obtain a vector bundle
$\Vtilde:=(\Psi^{\lambda})^{-1}(V)$ on $U^{\lambda}$.
We set
$U^{\lambda\ast}:=U^{\lambda}\setminus H^{\lambda}_{\infty}$
and
$\Vtilde^{\ast}:=\Vtilde_{|U^{\lambda\ast}}$.
We have already constructed
the $S^1_T$-equivariant differential operator
$\delbar_{\Vtilde^{\ast}}:
 C^{\infty}(U^{\lambda\ast},\Vtilde^{\ast})
 \lrarr
 C^{\infty}(U^{\lambda\ast},
 \Vtilde^{\ast}\otimes\Omega^{0,1}_{U^{\lambda\ast}})$
satisfying
$\delbar_{\Vtilde^{\ast}}(fu)=
\delbar_{U^{\lambda\ast}}(f)\cdot u
+f\delbar_{\Vtilde^{\ast}}(u)$
for any $f\in C^{\infty}(U^{\lambda\ast})$
and $u\in C^{\infty}(U^{\lambda\ast},\Vtilde^{\ast})$.

\begin{lem}
\label{lem;21.8.12.33}
The differential operator
$\delbar_{\Vtilde^{\ast}}$ uniquely extends to
a linear differential operator
$\delbar_{\Vtilde}:
 C^{\infty}(U^{\lambda},\Vtilde)
 \lrarr
 C^{\infty}(U^{\lambda},
 \Vtilde\otimes\Omega^{0,1}_{U^{\lambda}})$
satisfying the mini-complex Leibniz rule.
\end{lem}
\pf
We may assume that $U_w$ is an open disc around $\infty$.
Let $\vecv$ be a holomorphic frame of $(V,\delbar_V)$ on $U_w$.
We have $\del_{V,\wbar}\vecv=0$.
Let $A$ be the matrix valued $C^{\infty}$-function on
$U_{w}$ determined by
$\DDlambda_w\vecv=\vecv A$.
The pull back $\vecvtilde:=(\Psi^{\lambda})^{-1}\vecv$
is a $C^{\infty}$-frame of $\Vtilde$ on $U^{\lambda}$.

For any $P\in H^{\lambda}_{\infty}$,
we choose $\Ptilde\in H^{\lambda\cov}_{\infty}$
such that $\varpi^{\lambda}(\Ptilde)=P$.
Then, $(t_1,\tau_1)=(t_1,\beta_1^{-1})$ around $\Ptilde$
induces a local mini-complex coordinate system
on a neighbourhood $U^{\lambda}_1$ of $P$.
On $U^{\lambda}_1\setminus H^{\lambda}_{\infty}$,
the operators
$\del_{\Vtilde^{\ast},t_1}$
and
$\del_{\Vtilde^{\ast},\taubar_1}$
are described as follows
with respect to the frame $\vecvtilde$:
\[
 \del_{\Vtilde^{\ast},t_1}\vecvtilde
=\vecvtilde \frac{-2\sqrt{-1}}{1+|\lambda|^2}
 (\Psi^{\lambda})^{\ast}(A),
\quad
 \del_{\Vtilde^{\ast},\taubar_1}\vecvtilde
=-\betabar_1^{2}\del_{\Vtilde^{\ast},\betabar_1}\vecvtilde=0
\]
Note that $A$ is $C^{\infty}$ with respect to
$w^{-1}$ and $\wbar^{-1}$,
and that
$(\Psi^{\lambda})^{\ast}(w^{-1})
=(1+|\lambda|^2)\tau_1(1-2\sqrt{-1}\lambda t_1\tau_1)^{-1}$
is $C^{\infty}$ around $P$.
Hence, $\del_{\Vtilde^{\ast},t_1}$
and $\del_{\Vtilde^{\ast},\taubar_1}$
uniquely extend to operators 
$\del_{\Vtilde,t_1}$
and $\del_{\Vtilde,\taubar_1}$
on $C^{\infty}(U^{\lambda},\Vtilde)$.
Then, the claim of the lemma is clear.
\hfill\qed

\begin{lem}
\label{lem;21.8.12.10}
The above construction induces an equivalence between
the following objects.
\begin{itemize}
 \item Holomorphic vector bundles
       $(V,\delbar_V)$ on $U_w$ equipped with
       a linear differential operator
       $\DD^{\lambda}_w$ on $C^{\infty}(U_w,V)$
       such that
$\DD^{\lambda}_w(fs)
=\lambda\del_{w}(f)s+f\DD^{\lambda}_w(s)$
for any $f\in C^{\infty}(U_w)$
and $s\in C^{\infty}(U_w,V)$.
 \item $S^1_T$-equivariant vector bundles $\Vtilde$ on $U^{\lambda}$
       equipped with an $S^1_T$-equivariant differential linear operator
       $\delbar_{\Vtilde}:
       C^{\infty}(U^{\lambda},\Vtilde)
       \lrarr
       C^{\infty}(U^{\lambda},\Vtilde\otimes\Omega^{0,1}_{U^{\lambda}})$
       satisfying the mini-complex Leibniz rule.
\end{itemize}
\end{lem}
\pf
We indicate the inverse construction.
Let $(\Vtilde,\delbar_{\Vtilde})$ be as in the statement of
Lemma \ref{lem;21.8.12.10}.
There exists a $C^{\infty}$-vector bundle $V$ on $U_w$
equipped with an $S^1_T$-equivariant isomorphism
$\Vtilde\simeq (\Psi^{\lambda})^{-1}(V)$.
We set $V^{\ast}:=V_{|U_w\setminus\{\infty\}}$.
We obtain a $\lambda$-connection
$\DD^{\lambda}_{V^{\ast}}$
corresponding to $\delbar_{\Vtilde|U^{\lambda\ast}}$.

There exists the map $\pi^{\lambda}:U^{\lambda}\lrarr S^1_T$
induced by $(t_1,\beta_1)\longmapsto t_1$.
Fix $t_1^0\in S^1_T$.
The fiber $C:=(\pi^{\lambda})^{-1}(t_1^0)$
is naturally equipped with a complex structure,
and $(\Vtilde,\delbar_{\Vtilde})_{|C}$
is naturally a holomorphic vector bundle on $C$.
The induced morphism $\Psi^{\lambda}_C:C\lrarr U_w$
is holomorphic, and the image is an open neighbourhood of $\infty$.
We have the isomorphism of $C^{\infty}$-vector bundles
$\Vtilde_{|C}\simeq (\Psi^{\lambda})^{-1}(V)_{|C}$.
Note that
$(\Psi^{\lambda})^{-1}(V)_{|C\setminus H^{\lambda}_{\infty}}$
is equipped with the holomorphic structure induced by
$\delbar_{V^{\ast}}$,
and that
$\Vtilde_{|C\setminus H^{\lambda}_{\infty}}
\simeq
(\Psi^{\lambda})^{-1}(V)_{|C\setminus H^{\lambda}_{\infty}}$
is holomorphic.
It implies that
$\delbar_{V^{\ast}}$
uniquely extends to a holomorphic structure $\delbar_V$ of $V$.

Let $\vecv$ be a holomorphic frame of $(V,\delbar_V)$ on $U_w$.
We obtain a frame $\vecvtilde=(\Psi^{\lambda})^{-1}(\vecv)$
of $\Vtilde$.
There exists a matrix-valued $C^{\infty}$-function $A$ on $U_w$
determined by 
\[
 \del_{\Vtilde,t_1}\vecvtilde
 =\vecvtilde\cdot \frac{-2\sqrt{-1}}{1+|\lambda|^2}
 (\Psi^{\lambda})^{\ast}(A).
\]
Let $\DDlambda_{V^{\ast},w}$ denote the differential operator
of $V^{\ast}$
induced by $\DD^{\lambda}_{V^{\ast}}$ and $\del_w$.
By the construction of $\DD^{\lambda}_{V^{\ast}}$,
we have
$\DD^{\lambda}_{V^{\ast},w}\vecv_{|U_w\setminus\{\infty\}}
 =\vecv_{|U_w\setminus\{\infty\}}\cdot A$.
It implies that $\DD^{\lambda}_{V^{\ast},w}$
uniquely extends to a differential operator $\DD^{\lambda}_w$
of $V$.

In this way,
we obtain $(V,\delbar_V,\DD^{\lambda}_w)$
from $(\Vtilde,\delbar_{\Vtilde})$.
The two constructions are mutually inverse.
\hfill\qed

\subsubsection{Mini-holomorphic bundles and meromorphic flat $\lambda$-connections}

The following corollary is
an immediate consequence of Lemma \ref{lem;21.8.12.10}
and (\ref{eq;21.8.12.1}).

\begin{cor}
\label{cor;21.8.12.20}
 The construction in {\rm\S\ref{subsection;21.8.12.11}}
induces an equivalence of the following objects.
\begin{itemize}
 \item Locally free $\nbigo_{U_w}$-modules $\nbigv$
       equipped with
       a meromorphic $\lambda$-connection
       $\DD^{\lambda}:\nbigv\lrarr\nbigv\otimes\Omega^{1}_{U_w}(2\infty)$.
 \item $S^1_T$-equivariant locally free $\nbigo_{U^{\lambda}}$-modules
       $\nbigvtilde$.
       \hfill\qed
\end{itemize}
\end{cor}

\subsubsection{Another description of the construction}
\label{subsection;21.8.12.32}

Let us explain another more direct description of
the equivalence in Corollary \ref{cor;21.8.12.20}.

For any open subset $\nbigu\subset\nbigmbar^{\lambda}$,
let $\nbigk_{\nbigmbar^{\lambda}}(\nbigu)$
denote the space of $C^{\infty}$-functions $f$ on $\nbigu$
such that 
$\del_{\betabar_1}\bigl(
f_{|\nbigu\setminus H^{\lambda}_{\infty}}
\bigr)=0$.
For $P\in \nbigu\cap H^{\lambda}_{\infty}$,
a choice $\Ptilde\in(\varpi^{\lambda})^{-1}(P)$
induces a local mini-complex coordinate system
$(t_1,\tau_1)=(t_1,\beta_1^{-1})$ around $P$.
For any $f\in\nbigk_{\nbigmbar^{\lambda}}(\nbigu)$,
we obtain $\del_{\taubar_1}f=0$.
Because $[\del_{t_1},\del_{\betabar_1}]=0$,
we obtain the naturally defined action of $\del_{t_1}$
on $\nbigk_{\nbigmbar^{\lambda}}(\nbigu)$,
and $f\in\nbigk_{\nbigmbar^{\lambda}}(\nbigu)$ is
mini-holomorphic if and only if $\del_{t_1}f=0$.

Thus, we obtain a sheaf $\nbigk_{\nbigmbar^{\lambda}}$.
\index{sheaf $\nbigk_{\nbigmbar^{\lambda}}$}
We have the natural action $\del_{t_1}$ on
$\nbigk_{\nbigmbar^{\lambda}}$,
and the kernel is $\nbigo_{\nbigmbar^{\lambda}}$.
For any open subset $\nbigu\subset\nbigmbar^{\lambda}$,
let $\nbigk_{\nbigu}$ denote the restriction of
$\nbigk_{\nbigmbar^{\lambda}}$ to $\nbigu$.

Let $U_w$ be an open subset of $\proj^1_w$,
and $U^{\lambda}=(\Psi^{\lambda})^{-1}(U_w)
\subset\nbigmbar^{\lambda}$.
We have the naturally defined monomorphism of sheaves
$(\Psi^{\lambda})^{-1}(\nbigo_{U_w})
\lrarr \nbigk_{U^{\lambda}}$.

Let $\nbigv$ be a locally free $\nbigo_{U_w}$-module
equipped with a meromorphic $\lambda$-connection
$\DDlambda:\nbigv\lrarr\nbigv\otimes\Omega^1_{U_w}(2\infty)$.
From $\nbigo_{U_w}$-module $\nbigv$,
we obtain the locally free $\nbigk_{U^{\lambda}}$-module
\index{sheaf $\nbigvtilde^{\infty}$}
\[
 \nbigvtilde^{\infty}:=
 \nbigk_{U^{\lambda}}\otimes_{(\Psi^{\lambda})^{-1}(\nbigo_{U_w})}
 (\Psi^{\lambda})^{-1}(\nbigv). 
\]
We obtain the differential operator
$\del_{\nbigvtilde^{\infty},t_1}$ on $\nbigvtilde^{\infty}$
determined by the following condition
for local sections $f$ and $s$ of
$\nbigk_U$ and $\nbigv$
as follows:
\[
 \del_{\nbigvtilde^{\infty},t_1}\bigl(f(\Psi^{\lambda})^{-1}(s)\bigr)
=\del_{t_1}(f)(\Psi^{\lambda})^{-1}(s)
+\frac{1}{1+|\lambda|^2}f\cdot
(\Psi^{\lambda})^{-1}(-2\sqrt{-1}\DDlambda_ws).
\]
Note that $\nbigvtilde^{\infty}$
is equal to the sheaf of $C^{\infty}$-sections $u$ of
$\Vtilde$ such that $\del_{\Vtilde,\betabar_1}u=0$.
Hence, the kernel of $\del_{\nbigvtilde^{\infty},t_1}$
is naturally isomorphic to
the sheaf of mini-holomorphic sections of $\Vtilde$.

\begin{rem}
\label{rem;20.8.7.11}
 Although we use the map 
$\Psi^{\lambda}(t_1,\beta_1)=
 (1+|\lambda|^2)^{-1}(\beta_1-2\sqrt{-1}\lambda t_1)$,
it is more natural to consider
$(t_1,\beta_1)\longmapsto
 \beta_1-2\sqrt{-1}\lambda t_1$.
We adopt $\Psi^{\lambda}$
for the consistency with the dimensional reduction 
from monopoles to harmonic bundles
in {\rm\S\ref{subsection;17.9.29.1}}.
See Remark {\rm\ref{rem;20.8.7.10}}.
\hfill\qed
\end{rem}

\section{Curvatures of mini-holomorphic bundles 
with metric on $\nbigm^{\lambda}$}
\label{subsection;20.7.30.60}

\subsection{Contraction of curvature and analytic degree}

Let $U\lrarr \nbigm^{\lambda}$
be a local diffeomorphism.
We obtain the induced metric and 
the induced mini-complex structure on $U$.
We also obtain the complex vector fields
$\del_{\beta_i}$, $\del_{\betabar_i}$,
and $\del_{t_i}$ $(i=0,1)$ on $U$.

Let $(E,\delbar_E)$ be a mini-holomorphic bundle on $U$.
We have the operators
$\del_{E,\betabar_0}$ and $\del_{E,t_0}$.
Let $h$ be a Hermitian metric of $E$,
which induces
the Higgs field $\phi_h$
and the Chern connection $\nabla_h$.
Let $F(h)$ denote the curvature of $\nabla_h$.
We have the expression
$F(h)=F(h)_{\beta_0\betabar_0}d\beta_0\,d\betabar_0
+F(h)_{\beta_0,t_0}d\beta_0\,dt_0
+F(h)_{\betabar_0,t_0}d\betabar_0\,dt_0$.
We set 
\begin{equation}
\label{eq;20.7.30.10}
 G(h):=2F(h)_{\beta_0\betabar_0}
-\sqrt{-1}\nabla_{h,t_0}\phi_h.
\end{equation}
\index{endomorphism $G(h)$}
Note that the Bogomolny equation is equivalent to
$G(h)=0$
(see Corollary \ref{cor;20.7.16.1}).

\begin{df}
Suppose that $\Tr G(h)$ is 
described as a sum $g_1+g_2$,
where $g_1$ is $L^1$ on $U$
and $g_2$ is non-positive everywhere.
Then, we set 
$\deg(E,\delbar_E,h):=
 \int_U\Tr G(h)\dvol_U$,
where $\dvol_U$ is the volume form
induced by the Riemannian metric.
Note that $\deg(E,\delbar_E,h)\in\real\cup\{-\infty\}$.
\index{degree $\deg(E,\delbar_E,h)$}
\hfill\qed
\end{df}

\begin{rem}
\label{rem;17.9.29.20}
Let $\real_{s_0}\times U\lrarr \real_{s_0}\times\nbigm^{\lambda}$
be the induced local diffeomorphism.
We obtain the complex structure on $\real_{s_0}\times U$
such that
$(\alpha_0,\beta_0)=(s_0+\sqrt{-1}t_0,\beta_0)$
induces local complex coordinate systems.
Let $(\Etilde,\delbar_{\Etilde})$
be the holomorphic vector bundle 
with the metric $\htilde$
induced by $(E,\delbar_E)$ with $h$
as in \S{\rm\ref{subsection;13.11.29.2}}.
Let $\Ftilde$ denote the curvature of the Chern connection
associated with $(\Etilde,\delbar_{\Etilde},\htilde)$.
We have the expression 
$\Ftilde=
 \Ftilde_{\alpha_0,\alphabar_0}d\alpha_0\,d\alphabar_0
+\Ftilde_{\alpha_0,\betabar_0}d\alpha_0\,d\betabar_0
+\Ftilde_{\beta_0,\alphabar_0}d\beta_0\,d\alphabar_0
+\Ftilde_{\beta_0,\betabar_0}d\beta_0\,d\betabar_0$.
Then,
$\sqrt{-1}\Lambda\Ftilde
=2(\Ftilde_{\alpha_0,\alphabar_0}
+\Ftilde_{\beta_0,\betabar_0})$
is equal to the pull back of 
$G(h)$ by the projection
$\real_{s_0}\times U\lrarr U$,
where $\Lambda$ is the adjoint of 
the multiplication of the K\"ahler  form of
$\real_{s_0}\times U$.
(See {\rm\cite{Kobayashi-vector-bundle}}.)
Hence, $\deg(E,\delbar_E,h)$
is an analogue of the analytic degree
in {\rm\cite{Simpson88}}.
\hfill\qed
\end{rem}

\begin{rem}
In {\rm \cite{Mochizuki-triply-periodic-monopoles}}
and a previous version of this monograph,
we set
``$G(h)=
F(h)_{\beta_0\betabar_0}
-\frac{\sqrt{-1}}{2}\nabla_{h,t_0}\phi_{h}$'',
 which is the half of 
 {\rm(\ref{eq;20.7.30.10})}.
 There is no essential difference
 though constants
 in some formulas are changed.
\hfill\qed
\end{rem}

\subsection{Chern-Weil formula}
\label{subsection;20.7.30.61}

The standard Chern-Weil formula \cite{Simpson88, Simpson90}
is translated as follows,
which is implicitly contained in \cite{Charbonneau-Hurtubise}.
\index{Chern-Weil formula}

\begin{lem}
\label{lem;17.10.24.41}
Let $V$ be a mini-holomorphic subbundle of $E$.
Let $h_V$ be the induced metric of $V$.
Let $p_V$ denote the orthogonal projection of $E$
 onto $V$,
 which we regard as an endomorphism of $E$
in a natural way.
Then, we obtain the following formula:
\begin{equation}
\label{eq;21.8.22.40}
 \Tr G(h_V)=
 \Tr\bigl(G(h)p_V\bigr)
-2\bigl|
 \del_{E,\betabar_0}p_V
 \bigr|^2
-\frac{1}{2}\bigl|
 \del_{E,t_0}p_V
 \bigr|^2.
\end{equation}
\end{lem}
\pf
We use the notation in Remark \ref{rem;17.9.29.20}.
We obtain the induced holomorphic subbundle
$\Vtilde$ of $\Etilde$
with the metric $\htilde_{\Vtilde}$
induced by $V$.
We obtain the Chern connection
$\nabla_{\htilde}=\delbar_{\Etilde}+\del_{\Etilde,h}$.
Let $p_{\Vtilde}$ be the projection of $\Etilde$
onto $\Vtilde$,
which is the pull back of $p_V$.
We regard $p_{\Vtilde}$ as the endomorphism of
$\Etilde$ in a natural way.
Let $\iota_{\Vtilde}$ and $\iota_{\Vtilde^{\bot}}$
denote the inclusions of
$\Vtilde$ and $\Vtilde^{\bot}$
into $\Etilde$, respectively.
We set
$A:=p_{\Vtilde}\circ\delbar_{\Etilde}
 \circ\iota_{\Vtilde^{\bot}}$.
Because
$p_{\Vtilde}\circ F(\htilde)\circ\iota_{\Vtilde}
=F(\htilde_{\Vtilde})
-A\circ A^{\dagger}$,
we obtain
\[
 \sqrt{-1}\Tr\Lambda F(\htilde_{\Vtilde})
=\sqrt{-1}\Tr(p_{\Vtilde}\circ\Lambda F(\htilde)\circ\iota_{\Vtilde})
+\sqrt{-1}\Lambda\Tr(A\circ A^{\dagger}).
\]
Note that 
$\sqrt{-1}\Lambda F(\htilde)$
and 
$\sqrt{-1}\Lambda F(\htilde_{\Vtilde})$
are the pull back of
$G(h)$ and $G(h_V)$,
respectively.
Because of (\ref{eq;17.9.29.21}),
$\sqrt{-1}\Lambda\Tr(A\circ A^{\dagger})$
is the pull back of
$-2\bigl|
 \del_{E,\betabar_0}p_V
 \bigr|^2
-\frac{1}{2}\bigl|
 \del_{E,t_0}p_V
 \bigr|^2$.
Then, we obtain the claim of the lemma.
\hfill\qed

\vspace{.1in}
As a direct consequence of the lemma,
if $|G(h)|_h$ is $L^1$ on $U$,
then $\deg(V,h_V)$ makes sense in $\real\cup\{-\infty\}$
for any mini-holomorphic subbundle $V$ of $E$.

\subsection{Another description of $G(h)$}
\label{subsection;17.10.5.120}

We have the differential operator $\del_{E,\betabar_1}$
on $C^{\infty}(U,E)$,
which is given by the inner product of
$\del_{\betabar_1}$ and $\delbar_E$.
We have the differential operator 
$\del_{E,h,\beta_1}$ on $C^{\infty}(U,E)$
determined by the condition
$\del_{\betabar_1}h(u,v)
=h(\del_{E,\betabar_1}u,v)
+h(u,\del_{E,h,\beta_1}v)$.
\index{operator $\del_{E,h,\beta_1}$}

\begin{prop}
\label{prop;17.9.29.101}
The following formula holds:
\begin{multline}
\label{eq;17.9.29.100}
 G(h)= \\
 2(1+|\lambda|^2)^2
 \bigl[
 \del_{E,h,\beta_1},
 \del_{E,\betabar_1}
 \bigr] 
-\sqrt{-1}
 \Bigl(
 (1-|\lambda|^2)\nabla_{h,t_0}\phi
+2\lambda\sqrt{-1}\nabla_{h,\beta_0}\phi
-2\lambdabar\sqrt{-1}\nabla_{h,\betabar_0}\phi
 \Bigr)
\\
=2(1+|\lambda|^2)^2
 \bigl[
 \del_{E,h,\beta_1},
 \del_{E,\betabar_1}
 \bigr]
-\sqrt{-1}
 (1+|\lambda|^2)
 \nabla_{h,t}\phi.
\end{multline}
Here, $\nabla_{h,t}$
denote the inner product of 
$\nabla_h$
and $\del_t$.
(See {\rm(\ref{eq;17.9.29.30})}.)
\end{prop}
\pf
To simplify the description,
we omit to denote the dependence on $h$,
i.e.,
$\del_{E,h,\beta_1}$ is denoted by $\del_{E,\beta_1}$,
for example.
By the relation (\ref{eq;17.9.29.31}),
we obtain the following:
\begin{equation}
 \del_{E,\betabar_1}
=\frac{\lambda}{1+|\lambda|^2}
 \frac{1}{2\sqrt{-1}}\del_{E,t_0}
+\frac{1}{1+|\lambda|^2}
  \del_{E,\betabar_0}
=\frac{1}{1+|\lambda|^2}
 \nabla_{\betabar_0}
-\frac{\lambda}{1+|\lambda|^2}
 \frac{1}{2}(\phi+\sqrt{-1}\nabla_{t_0}).
\end{equation}
We obtain the following:
\[
 \del_{E,\beta_1}
=\frac{1}{1+|\lambda|^2}
 \nabla_{\beta_0}
-\frac{\lambdabar}{1+|\lambda|^2}
 \frac{1}{2}(\phi-\sqrt{-1}\nabla_{t_0}).
\]
We have the following:
\begin{multline}
 \bigl[
 \del_{E,\beta_1},
 \del_{E,\betabar_1}
 \bigr]
=\Bigl(
 \frac{1}{1+|\lambda|^2}
 \Bigr)^2
 \Bigl(
 \bigl[
 \nabla_{\beta_0},\nabla_{\betabar_0}
 \bigr]
-\frac{\lambda}{2}
 \bigl(
 \nabla_{\beta_0}\phi
+\sqrt{-1}[\nabla_{\beta_0},\nabla_{t_0}]
 \bigr)
 \Bigr.
 \\
 \Bigl.
+\frac{\lambdabar}{2}
 \bigl(
 \nabla_{\betabar_0}\phi
-\sqrt{-1}[\nabla_{\betabar_0},\nabla_{t_0}]
 \bigr)
+\frac{|\lambda|^2}{4}
 \bigl[\phi-\sqrt{-1}\nabla_{t_0},
 \phi+\sqrt{-1}\nabla_{t_0}\bigr]
 \Bigr).
\end{multline}
By Lemma \ref{lem;17.9.29.40},
we obtain
\[
 -\frac{\lambda}{2}
 \nabla_{\beta_0}\phi
-\frac{\lambda}{2}\sqrt{-1}
 \bigl[\nabla_{\beta_0},\nabla_{t_0}\bigr]
=-\lambda\nabla_{\beta_0}\phi,
\quad\quad
 \frac{\lambdabar}{2}
 \nabla_{\betabar_0}\phi
-\frac{\sqrt{-1}}{2}\lambdabar
 \bigl[
 \nabla_{\betabar_0},\nabla_{t_0}
 \bigr]
=\lambdabar\nabla_{\betabar_0}\phi.
\]
We also have
\[
 \frac{|\lambda|^2}{4}
 \bigl[
 \phi-\sqrt{-1}\nabla_{t_0},
 \phi+\sqrt{-1}\nabla_{t_0}
 \bigr]
=-\frac{|\lambda|^2}{2}
 \sqrt{-1}\nabla_{t_0}\phi.
\]
Hence, we obtain the following:
\begin{multline}
 \bigl[
 \del_{E,\beta_1},\del_{E,\betabar_1}
 \bigr]
=\left(\frac{1}{1+|\lambda|^2}\right)^2
 \Bigl(
 F_{\beta_0,\betabar_0}
-\lambda\nabla_{\beta_0}\phi
+\lambdabar\nabla_{\betabar_0}\phi
-\frac{|\lambda|^2}{2}
 \sqrt{-1}\nabla_{t_0}\phi
 \Bigr)
\\
=\left(\frac{1}{1+|\lambda|^2}\right)^2
 \Bigl(
\frac{1}{2}G(h)
-\lambda\nabla_{\beta_0}\phi
+\lambdabar\nabla_{\betabar_0}\phi
+\frac{1-|\lambda|^2}{2}
 \sqrt{-1}\nabla_{t_0}\phi
 \Bigr).
\end{multline}
Then, we obtain (\ref{eq;17.9.29.100}).
\hfill\qed

\vspace{.1in}
We state some consequences.

\begin{cor}
\label{cor;17.9.29.120}
Suppose that $U:=\nbigm^{\lambda}\setminus Z$,
where $Z$ is a finite subset of $\nbigm^{\lambda}$.
Suppose that 
$\Tr\bigl[\del_{E,h,\beta_1},\del_{E,\betabar_1}\bigr]$
and 
$\Tr\nabla_{h,t}\phi_h$ are $L^1$ on 
 $\nbigm^{\lambda}\setminus Z$.
Then, we obtain
\begin{equation}
 \label{eq;17.9.29.102}
 \deg(E,\delbar_E,h)
=\int_{0}^Tdt_1
 \int_{\cnum_{\beta_1}}
 \Tr\Bigl(
 \bigl[\del_{E,h,\beta_1},\del_{E,\betabar_1}\bigr]
 \Bigr)
 \sqrt{-1}d\beta_1\,d\betabar_1.
\end{equation}
\end{cor}
\pf
The following formula holds:
\[
 \dvol=\frac{\sqrt{-1}}{2}
 d\beta_0d\betabar_0\,dt_0
=\frac{\sqrt{-1}}{2}
 (1+|\lambda|^2)^{-2}
 d\beta_1\,d\betabar_1\,dt_1
=\frac{\sqrt{-1}}{2}
 dw\,d\wbar\,dt.
\]
By using Fubini theorem,
we obtain the following:
\[
 \int_{\nbigm^{\lambda}\setminus Z}
 \Tr\bigl(\nabla_{t}\phi\bigr)
 \dvol
=\frac{\sqrt{-1}}{2}
 \int_{\cnum_w}
 dw\,d\wbar
 \int_{\real/T\seisuu}
 \del_{t}\Tr(\phi_h)\,dt
=0.
\]
Hence, by Proposition \ref{prop;17.9.29.101},
we obtain the following.
\begin{multline}
 \int_{\nbigm^{\lambda}\setminus Z}
 \Tr G(h)\dvol
=
\int_{\nbigm^{\lambda}}
 \Tr\bigl(
 [\nablatilde_{\beta_1},
 \nablatilde_{\betabar_1}]
 \bigr)
 \sqrt{-1}
 d\beta_1\,d\betabar_1\,dt_1 \\
=
\int_{0}^T dt_1
 \int_{\cnum_{\eta_1}}
 \Tr\bigl(
 \bigl[
 \nablatilde_{\beta_1},
 \nablatilde_{\betabar_1}
\bigr]
 \bigr)
 \sqrt{-1}
 d\beta_1\,d\betabar_1.
\end{multline}
Thus, we obtain (\ref{eq;17.9.29.102}).
\hfill\qed

\vspace{.1in}
Note that for almost all $0\leq t_1\leq T$,
the analytic degree of
$(E,\delbar_E,h)_{|\{t_1\}\times\cnum_{\beta_1}}$
is defined as the integration of the first Chern form
(see \cite{Simpson88}):
\begin{equation}
\label{eq;21.8.19.10}
 \deg\bigl(
 (E,\delbar_E,h)_{|\{t_1\}\times\cnum_{\beta_1}}
 \bigr)
=\frac{\sqrt{-1}}{2\pi}
 \int_{\cnum_{\beta_1}} 
 \Tr\bigl([\del_{E,h,\beta_1},\del_{E,\betabar_1}]\,
 d\beta_1\,d\betabar_1
  \bigr).
\end{equation}
Then, (\ref{eq;17.9.29.102}) is rewritten as follows:
\begin{equation}
\label{eq;21.8.19.11}
 \deg(E,\delbar_E,h)
=\int_0^T
 2\pi\deg\bigl((E,\delbar_E,h)_{|\{t_1\}\times\cnum_{\beta_1}}\bigr)\,
 dt_1.
\end{equation}

We also obtain the following useful formula
to relate the curvatures of the mini-holomorphic bundles
on $\nbigm^{\lambda}$ and $\nbigm^0$,
underlying a monopole.

\begin{cor}
\label{cor;17.12.16.2}
Suppose that $(E,\delbar_E,h)$ is a monopole on $U$,
i.e., $G(h)=0$.
Then, we obtain
\[
 \bigl[
 \del_{E,h,\beta_1},
 \del_{E,\betabar_1}
 \bigr]
=\frac{\sqrt{-1}}{2}
 \frac{1}{1+|\lambda|^2}
 \nabla_{h,t}\phi
=\frac{1}{1+|\lambda|^2}F_{w,\wbar}.
\]
\hfill\qed
\end{cor}

\subsection{Change of metrics}
\label{subsection;20.7.30.62}

Let $h_1$ be another Hermitian metric of $E$.
Let $s$ be the automorphism of $E$
determined by $h_1=h\cdot s$,
which is self-adjoint with respect to both $h$ and $h_1$.
The following is a variant of \cite[Lemma 3.1]{Simpson88}.
\begin{lem}
\label{lem;17.10.13.110}
The following holds:
\begin{equation}
\label{eq;17.9.29.110}
 G(h_1)=G(h)
-2\del_{E,\betabar_0}
 \bigl(
 s^{-1}\del_{E,h,\beta_0}s
 \bigr)
-\frac{1}{2}
 \Bigl[
 \nabla_{h,t_0}-\sqrt{-1}\phi_h,\,
 s^{-1}\bigl[
 \nabla_{h,t_0}+\sqrt{-1}\phi_{h},s
 \bigr]
 \Bigr].
\end{equation}
\end{lem}
\pf
Because
$\del_{E,h_1}=\del_{E,h}+s^{-1}\del_{E,h}s$,
the following holds:
\[
 F(h_1)_{\beta_0,\betabar_0}
=\bigl[
 \del_{E,h_1,\beta_0},\del_{E,\betabar_0}
 \bigr]
=F(h)_{\beta_0,\betabar_0}
-\del_{E,\betabar_0}(s^{-1}\del_{E,h,\beta_0}s).
\]
We also obtain the following:
\[
 \nabla_{h_1,t_0}=
 \nabla_{h,t_0}
+\frac{1}{2}s^{-1}
 \bigl[
 \nabla_{h,t_0}+\sqrt{-1}\phi_h,s
 \bigr],
\quad\quad
 \phi_{h_1}
=\phi_{h}
-\frac{\sqrt{-1}}{2}
 s^{-1}\bigl[
 \nabla_{h,t_0}+\sqrt{-1}\phi_h,s
 \bigr].
\]
We obtain the following:
\[
 -\sqrt{-1}
 \nabla_{h_1,t_0}\phi_{h_1}
=-\sqrt{-1}
 \nabla_{h,t}\phi_{h}
-\frac{1}{2}
 \Bigl[
 \nabla_{h,t_0}-\sqrt{-1}\phi_{h},
 s^{-1}
 \bigl[
  \nabla_{h,t_0}+\sqrt{-1}\phi_{h},s
 \bigr]
 \Bigr].
\]
Hence, we obtain (\ref{eq;17.9.29.110}).
\hfill\qed

\vspace{.1in}
We obtain the following direct consequence,
which is also a variant of \cite[Lemma 3.1]{Simpson88}.
\begin{cor}
\label{cor;17.10.15.12}
The following equality holds:
\begin{multline}
 -\Bigl(\del_{\betabar_0}\del_{\beta_0}+\frac{1}{4}\del_{t_0}^2\Bigr)
 \Tr s
 =\\
 \frac{1}{2}\Tr\Bigl(
 s\bigl(G(h_1)-G(h)\bigr)
 \Bigr)
-\bigl|s^{-1/2}\del_{E,h,\beta_0}s\bigr|_h^2
-\frac{1}{4}
 \bigl|s^{-1/2}\del'_{E,h,t_0}s\bigr|_h^2.
\end{multline}
We also have the following inequality:
\[
 -\Bigl(
 \del_{\betabar_0}
 \del_{\beta_0}
+\frac{1}{4}\del_{t_0}^2
 \Bigr)
 \log\bigl(\Tr(s)\bigr)
 \leq
 \frac{1}{2}\Bigl(
 \bigl|G(h)\bigr|_h
 +\bigl|G(h_1)\bigr|_{h_1}
 \Bigr).
\]
\hfill\qed
\end{cor}

\begin{cor}
\label{cor;17.10.24.31}
If $\rank E=1$,
we have
$G(h_1)-G(h)
=2^{-1}\Delta\log s$ on $U$.
Here, $\Delta$
denote the Laplacian 
of the Riemannian manifold $\nbigm^{\lambda}$.
\hfill\qed
\end{cor}

\subsection{Relation with $\lambda$-connections}
\label{subsection;17.10.12.1}

This subsection is a complement to \S\ref{section;21.8.12.101}.
We use the notation there.

Let us recall the condition for harmonic bundles
given in terms of $\lambda$-connections
\cite[\S2.2]{Mochizuki-KHII}.
Let $(V,\DDlambda)$ be a $\lambda$-flat bundle
on $U_w$ with a Hermitian metric $h_V$.
Let $\DDlambda=d''_V+d'_V$ denote the decomposition
into the $(0,1)$-part and the $(1,0)$-part.
We have the $(1,0)$-operator $\delta'$
and the $(0,1)$-operator $\delta''$
determined by the conditions
\index{operators $\delta'$, $\delta''$}
\[
 \delbar h_V(u,v)
=h_V(d''_Vu,v)+h_V(u,\delta'v),
\quad\quad
 \lambda \del h_V(u,v)
=h_V(d'_Vu,v)+h_V(u,\delta''v)
\]
for $u,v\in C^{\infty}(U_w,V)$.
We have 
the $(0,1)$-operator $\delbar_V$,
the $(1,0)$-operator $\del_V$,
the section
  $\theta\in C^{\infty}(U_w,V\otimes\Omega^{1,0})$
and the section
  $\theta^{\dagger}\in C^{\infty}(U_w,V\otimes\Omega^{0,1})$
determined by 
\index{operators $\delbar_V$, $\del_V$, $\theta$, $\theta^{\dagger}$}
\[
d''_V=\delbar_V+\lambda\theta^{\dagger},
\quad
d'_V=\lambda\del_V+\theta,
\quad
\delta'=\del_V-\lambdabar\theta,
\quad
\delta''=\lambdabar\del_V-\theta^{\dagger}.
\]
Then, $(V,\DDlambda,h)$ is called a harmonic bundle
if $(V,\delbar_V,\theta,h)$ is a harmonic bundle.
If $\lambda\neq 0$,
it is equivalent to $\delbar_V\theta=0$.
We set 
$\DD^{\lambda\star}:=\delta'-\delta''$.
Then, 
$(V,\DDlambda,h_V)$ is a harmonic bundle
if and only if $[\DD^{\lambda},\DD^{\lambda\star}]=0$.
\index{operator $\DD^{\lambda\star}$}

\vspace{.1in}

Let $d''_{V,\wbar}$ denote the inner product
of $d''_V$ and $\del_{\wbar}$.
\index{operator $d''_{V,\wbar}$, $\delta''_{\wbar}$, etc.}
We use the notation
$d'_{V,w}$, $\delta''_{\wbar}$, etc., in similar meanings.
Note that
$d'_w-\lambda\delta'_w$ and $\delta''_{\wbar}-\lambdabar d'_{\wbar}$
are endomorphisms of $V$,
which follow from
$d'-\lambda\delta'=(1+|\lambda|^2)\theta$
and
$\delta''-\lambdabar d'=-(1+|\lambda|^2)\theta^{\dagger}$.

As explained in \S\ref{section;21.8.12.101},
we obtain a mini-holomorphic bundle
$(\Vtilde,\delbar_{\Vtilde})$
on $U^{\lambda}=\Psi^{-1}(U_w)
\subset\nbigm^{\lambda}$
from $(V,\DDlambda)$.
We also obtain the Hermitian metric
$\htilde:=(\Psi^{\lambda})^{-1}(h_V)$.
Let $\nabla_{\Vtilde}$ and $\phi$
denote the Chern connection and the Higgs field
of $(\Vtilde,\delbar_{\Vtilde},\htilde)$.
Let $\del'_{\Vtilde,t_0}$
be the differential operator
induced by $\del_{\Vtilde,t_0}$ and $h$ 
as in \S\ref{subsection;16.9.20.20}.

\begin{lem}
\label{lem;17.10.13.101}
We obtain the following formulas:
\[
 \bigl[
 \nabla_{\Vtilde,\betabar_0},\nabla_{\Vtilde,\beta_0}
 \bigr]
=\frac{1}{(1+|\lambda|^2)^2}
 (\Psi^{\lambda})^{-1}
 \Bigl(
 \bigl[
 d''_{\wbar},\delta'_w
 \bigr]
+\lambdabar
 \bigl[
 d''_{\wbar},\delta''_{\wbar}
 \bigr]
+\lambda[d'_w,\delta'_{w}
 \bigr]
+|\lambda|^2
 \bigl[
 d'_w,\delta''_{\wbar}
 \bigr]
 \Bigr),
\]
\[
 \bigl[
 \nabla_{\Vtilde,\betabar_0},
 \del'_{\Vtilde,t_0}
 \bigr]
=\frac{2\sqrt{-1}}{(1+|\lambda|^2)^2}
 (\Psi^{\lambda})^{-1}
 \Bigl(
 [d''_{\wbar},\delta''_{\wbar}]
-\lambda[d''_{\wbar},\delta'_w]
+\lambda[d'_w,\delta''_{\wbar}]
-\lambda^2[d'_w,\delta'_w]
 \Bigr),
\]
\[
 \bigl[
 \nabla_{\Vtilde,\beta_0},
 \del_{\Vtilde,t_0}
 \bigr]
=\frac{-2\sqrt{-1}}{(1+|\lambda|^2)^2}
 (\Psi^{\lambda})^{-1}
 \Bigl(
 \bigl[
 \delta'_w,d'_w
 \bigr]
-\lambdabar[\delta'_w,d''_{\wbar}]
+\lambdabar[\delta''_{\wbar},d'_w]
-\lambdabar^2[\delta''_{\wbar},d''_{\wbar}]
 \Bigr),
\]
\[
 \bigl[
 \del_{\Vtilde,t_0},
 \del'_{\Vtilde,t_0}
 \bigr]
=\frac{4}{(1+|\lambda|^2)^2}
 (\Psi^{\lambda})^{-1}
 \Bigl(
 [d'_w,\delta''_{\wbar}]
-\lambda[d'_w,\delta'_w]
-\lambdabar[d''_{\wbar},\delta''_{\wbar}]
+|\lambda|^2[d''_{\wbar},\delta'_w]
 \Bigr).
\]
We also have the following formula:
\[
 \phi=\frac{1}{1+|\lambda|^2}
 (\Psi^{\lambda})^{-1}
 \bigl(
 (d'_w-\lambda\delta'_w)+(\delta''_{\wbar}-\lambdabar d'_{\wbar})
 \bigr).
\]
\end{lem}
\pf
By the construction,
the following holds
for any $s\in C^{\infty}(U_w,V)$:
\[
 \del_{\Vtilde,t_1}(\Psi^{\lambda})^{-1}(s)
=\frac{-2\sqrt{-1}}{1+|\lambda|^2}(\Psi^{\lambda})^{-1}\Bigl(
 \bigl(
 d'_w-\lambdabar d''_{\wbar}
 \bigr)s
 \Bigr).
\]
We denote it as 
\[
 \del_{\Vtilde,t_1}
=\frac{-2\sqrt{-1}}{1+|\lambda|^2}
 (\Psi^{\lambda})^{-1}
\Bigl(
d'_w-\lambdabar d''_{\wbar}
 \Bigr).
\]
Similarly, we have 
\begin{eqnarray}
\label{eq;17.12.9.1}
\begin{cases}
  \del'_{\Vtilde,t_1}
=\frac{2\sqrt{-1}}{1+|\lambda|^2}
  (\Psi^{\lambda})^{-1}
 \Bigl(
 \delta''_{\wbar}
-\lambda \delta'_w
 \Bigr) &
 \\
 \del_{\Vtilde,\betabar_1}
=\frac{1}{1+|\lambda|^2}
 (\Psi^{\lambda})^{-1}\bigl(d''_{\wbar}\bigr)
 &
 \\
  \del_{\Vtilde,\beta_1}
=\frac{1}{1+|\lambda|^2}
 (\Psi^{\lambda})^{-1}(\delta'_w).
 &
\end{cases}
\end{eqnarray}
Note that $\del_{\Vtilde,t_0}=\del_{\Vtilde,t_1}$
and $\del'_{\Vtilde,t_0}=\del'_{\Vtilde,t_1}$.
We obtain
\[
 \nabla_{\Vtilde,\betabar_0}
=(1+|\lambda|^2)\del_{\Vtilde\betabar_1}
-\frac{\lambda}{2\sqrt{-1}}\del_{\Vtilde,t_1}
=
 \frac{1}{1+|\lambda|^2}
 (\Psi^{\lambda})^{-1}(d''_{\wbar})
+\frac{\lambda}{1+|\lambda|^2}(\Psi^{\lambda})^{-1}(d'_w).
\]
Similarly, we obtain
\[
 \nabla_{\Vtilde,\beta_0}
=\frac{1}{1+|\lambda|^2}
 (\Psi^{\lambda})^{-1}(\delta'_w)
+\frac{\lambdabar}{1+|\lambda|^2}
 (\Psi^{\lambda})^{-1}(\delta''_{\wbar}).
\]
Then, we obtain the desired equalities
by direct computations.
\hfill\qed

\begin{lem}
\label{lem;17.10.12.2}
We have the following formula:
\begin{equation}
\label{eq;17.10.11.20}
G(\htilde)=
\frac{1}{(1+|\lambda|^2)}
(\Psi^{\lambda})^{-1}
 \Bigl(
 \sqrt{-1}\Lambda_{U_w} [\DD^{\lambda},\DD^{\lambda\star}]
 \Bigr).
\end{equation}
Here, $\Lambda_{U_w}:\Omega^{1,1}_{U_w}\lrarr \Omega^{0,0}_{U_w}$
is determined by
$\Lambda_{U_w}(dw\,d\wbar)=-2\sqrt{-1}$.
\end{lem}
\pf
By definition, we have
\[
 -\sqrt{-1}
 \nabla_{\Vtilde,t_0}\phi
=-\frac{1}{2}
 \bigl[
 \del_{\Vtilde,t_0},\del'_{\Vtilde,t_0}
 \bigr].
\]
By a direct computation,
we obtain
\begin{multline}
2\bigl[
 \nabla_{\Vtilde,\beta_0},\nabla_{\Vtilde,\betabar_0}
 \bigr]
-\sqrt{-1}\nabla_{t_0}\phi
 = \\
 \frac{2}{(1+|\lambda|^2)^2}
 (\Psi^{\lambda})^{-1}
 \Bigl(
-\bigl[
  d''_{\wbar}+\lambda d'_{w},
 \delta'_{w}+\lambdabar\delta''_{\wbar}
 \bigr]
-\bigl[
 d'_w-\lambdabar d''_{\wbar},
 \delta''_{\wbar}
-\lambda \delta'_w
 \bigr]
\Bigr)
\\
=\frac{2}{1+|\lambda|^2}
 (\Psi^{\lambda})^{-1}
 \Bigl(
 \bigl[
 \delta'_w,d''_{\wbar}
\bigr]
+\bigl[
 \delta''_{\wbar},
 d'_{w} \bigr]
 \Bigr).
\end{multline}
We also have
\[
 \sqrt{-1}
 \Lambda_{U_w}
 \bigl[
 \DD^{\lambda},\DD^{\lambda\star}\bigr]
=\sqrt{-1}
 \Lambda_{U_w}
 \Bigl(
 [d'',\delta']
-[d',\delta'']
 \Bigr)
=2\Bigl(
 \bigl[
 \delta'_w,d''_{\wbar}
\bigr]
+\bigl[
 \delta''_{\wbar},
 d'_{w} \bigr]
 \Bigr).
\]
Hence, we obtain (\ref{eq;17.10.11.20}).
\hfill\qed

\begin{cor}
\label{cor;21.8.13.30}
$(\Vtilde,\delbar_{\Vtilde},\htilde)$ is a monopole
if and only if 
$(V,\DDlambda,h)$ is a harmonic bundle
\hfill\qed
\end{cor}

\begin{lem}
\label{lem;17.10.24.50}
Let $\del_{\Vtilde,\htilde,\beta_1}$
be the operator induced by
$\del_{\Vtilde,\betabar_1}$
and $\htilde$ as in {\rm\S\ref{subsection;17.10.5.120}}.
Then, we obtain the following equality:
\[
 \bigl[
 \del_{\Vtilde,\betabar_1},
 \del_{\Vtilde,\htilde,\beta_1}
 \bigr]
=
 \frac{1}{(1+|\lambda|^2)^2}
 \bigl( \Psi^{\lambda}\bigr)^{-1}
 \Bigl(
 [d''_{\wbar},\delta'_{h,w}]
 \Bigr).
\]
\end{lem}
\pf
It follows from (\ref{eq;17.12.9.1}).
\hfill\qed

\subsubsection{$\lambda$-flat bundles of infinite rank
with a harmonic metric}
\label{subsection;21.8.13.100}

Let $(E,\delbar_E)$ be a mini-holomorphic bundle on $U^{\lambda}$.
We obtain the Chern connection $\nabla$
and the Higgs field $\phi$.
Let $E_{C^{\infty}}$ denote the sheaf of $C^{\infty}$-sections of $E$.
We obtain
the $\nbigc^{\infty}_{U_w}$-module
$V_1=\Psi_!(E_{C^{\infty}})$
equipped with the induced $\lambda$-connection $\DDlambda_{V_1}$
as explained in \S\ref{subsection;21.8.13.31}.

Let $h$ be a Hermitian metric of $E$.
It induces a Hermitian metric $h_1=\Psi_!(h)$ of $V_1$.
We obtain the operators
\[
\delta_w'=(1+|\lambda|^2)\del_{E,h,\beta_1},
\quad
\delta''_{\wbar}=\frac{-\sqrt{-1}}{2}
(1+|\lambda|^2)
\del'_{E,h,t_0}
+\lambda(1+|\lambda|^2)\del_{E,h,\beta_1}
\]
on $V_1$.
We set 
$\delta'(s)=\delta'_w(s)\,dw$
and
$\delta''(s)=\delta''_{\wbar}(s)\,d\wbar$
for any local section $u$ of $V_1$.
We obtain $\DD^{\lambda\star}_{V_1}=\delta'-\delta''$ on $V_1$.
We have the decomposition $\DDlambda_{V_1}=d''+d'$
into the $(0,1)$-part and the $(1,0)$-part.
By the construction,
we have
$\delbar h_1(u_1,u_2)
=h_1(d''u_1,u_2)+h_1(u_1,\delta'u_2)$
and
$\del h_1(u_1,u_2)
=h_1(d'u_1,u_2)+h_1(u_1,\delta''u_2)$
for local sections $u_i$ of $V_1$.
By similar computations,
we obtain that
$\sqrt{-1}\Lambda_{U_w}[\DD^{\lambda}_{V_1},\DD^{\lambda\star}_{V_2}]$
is the multiplication of $(1+|\lambda|^2)G(h)$.
Hence, $[\DD^{\lambda}_{V_1},\DD^{\lambda\star}_{V_1}]=0$
if and only if $(E,\delbar_E,h)$ is a monopole.
In this sense,
we may regard a monopole on $U^{\lambda}$
as a $\lambda$-flat bundle with a harmonic metric of infinite rank.

\subsubsection{Remark}
\label{subsection;21.8.13.101}

The analogy between monopoles on $U^{\lambda}$
and $\lambda$-flat bundles with a harmonic metric on $U_w$
(Corollary \ref{cor;21.8.13.30} and \S\ref{subsection;21.8.13.100})
is one of the motivation of this study,
as explained in \S\ref{section;21.8.13.20}.
It is also useful for the construction of
a Hermitian metric $h_1$ of a given mini-holomorphic bundle
such that $G(h_1)$ is small
in Proposition \ref{prop;17.10.12.6}.

In the study of harmonic bundles,
we study the asymptotic behaviour of
a harmonic metric $h$ of $(V,\DDlambda)$
by constructing a Hermitian metric $h_0$ of $V$
such that $\nbigp^h_{\ast}V=\nbigp^{h_0}_{\ast}V$,
and that $[\DDlambda_{h_0},\DD^{\lambda\star}_{h_0}]$ is small.
We can explicitly construct such $h_0$,
and we can prove that $h$ and $h_0$ are mutually bounded.
This kind of argument was pioneered by Simpson \cite{Simpson90},
and further applied in \cite{mochi2,Mochizuki-wild}.
Such a Hermitian metric
is also useful in the proof of Kobayashi-Hitchin correspondence,
i.e., the existence of globally defined harmonic metrics.
By adopting a similar strategy,
in Proposition \ref{prop;17.10.12.6} below,
we shall construct a Hermitian metric $h_1$
of a mini-holomorphic bundle such that $G(h_1)$ is small.
For the construction of such $h_1$,
we use a metric $h_0$ for $\lambda$-flat bundles as above
through Lemma \ref{lem;17.10.13.101},
Lemma \ref{lem;17.10.12.2} and Lemma \ref{lem;17.10.24.50}.

\subsection{Dimensional reduction of Kronheimer}
\label{subsection;17.10.24.32}

Let us recall the dimensional reduction of Kronheimer
\cite{Kronheimer-Master-Thesis}.
Let $\varphi:\cnum^2\lrarr \real\times\cnum$
be the map (\ref{eq;20.7.30.40}).
Let $U$ be a neighbourhood of 
$(0,0)$ in $\real\times\cnum$.
We set
$\Utilde:=\varphi^{-1}(U)$.
We put
$U^{\ast}:=
  U\setminus\{(0,0)\}$
and 
$\Utilde^{\ast}:=
 \Utilde\setminus\{(0,0)\}$.

Let $(E,\delbar_E)$
be a mini-holomorphic bundle on $U^{\ast}$
with a Hermitian metric $h_E$.
We obtain the Chern connection $\nabla$ and 
the Higgs field $\phi$.
Let $F$ be the curvature of $\nabla$.

We put
$\Etilde:=\varphi^{-1}(E)$
on $\Utilde^{\ast}$.
It is equipped with the unitary connection
$\nablatilde:=
 \varphi^{\ast}(\nabla)
+\sqrt{-1}\varphi^{\ast}(\phi)\otimes \xi$,
where
\[
 \xi=-u_1\,d\ubar_1+\ubar_1\,du_1
-\ubar_2\,du_2+u_2d\ubar_2.
\]
The curvature $\Ftilde$
is equal to
$\varphi^{\ast}(F)
+\sqrt{-1}\varphi^{\ast}(\nabla\phi)
\wedge \xi$.

\begin{lem}
\label{lem;17.10.26.10}
We have the following equality
of currents on $U$:
\begin{equation}
\label{eq;17.9.5.2}
\varphi_{\ast}
 \Bigl(
 \Tr\bigl(
 \Lambda\Ftilde
\bigr)\dvol_{\Utilde}
\Bigr)
=
\pi\Tr G(h)\dvol_U.
 \end{equation}
Here, $\varphi_{\ast}$ denote the push-forward of
currents by the proper map $\varphi$.
\end{lem}
\pf
Note 
$(|u_1|^2+|u_2|^2)^2=\varphi^{\ast}(|w|^2+|t|^2)$.
We have $\varphi^{\ast}dw=2u_1du_2+u_2du_1$,
$\varphi^{\ast}d\wbar=2\ubar_1d\ubar_2+\ubar_2d\ubar_1$
and 
$\varphi^{\ast}dt=\ubar_1du_1+u_1d\ubar_1-\ubar_2du_2-u_2d\ubar_2$.
The forms
$\varphi^{\ast}dw$,
$\varphi^{\ast}d\wbar$,
$\varphi^{\ast}dt$
and $\xi$ are orthogonal.
Moreover, 
we have
$|\varphi^{\ast}dw|^2=|\varphi^{\ast}d\wbar|^2=8(|u_1|^2+|u_2|^2)$,
$|\varphi^{\ast}dt|^2=4(|u_1|^2+|u_2|^2)$
and $|\xi|^2=4(|u_1|^2+|u_2|^2)$.
Hence, we have
\[
 \bigl|
 \varphi^{\ast}(dw\,d\wbar\,dt)\,\xi
 \bigr|
=8(|u_1|^2+|u_2|^2)^2
 \bigl|
 du_1\,d\ubar_1\,du_2\,d\ubar_2
 \bigr|.
\]
Let $\psi_{(u_1,u_2)}:\real/2\pi\lrarr \cnum^2$
be given by
$\psi_{(u_1,u_2)}(e^{\sqrt{-1}\theta})
=\bigl(u_1e^{\sqrt{-1}\theta},u_2e^{-\sqrt{-1}\theta}\bigr)$.
We have
\[
 \psi_{(u_1,u_2)}^{\ast}\xi
=2\sqrt{-1}(|u_1|^2+|u_2|^2)d\theta.
\]
Hence, we obtain
\begin{multline}
 \varphi_{\ast}\bigl(|du_1\,d\ubar_1\,du_2\,d\ubar_2|\bigr)
=\frac{1}{8(|w|^2+|t|^2)}
 \bigl|dw\,d\wbar\,dt\bigr|\cdot 4\pi(|w|^2+|t|^2)^{1/2}
 \\
=\frac{\pi}{2(|w|^2+|t|^2)^{1/2}}
\bigl|dw\,d\wbar\,dt\bigr|. 
\end{multline}
We have
$|dw\,d\wbar\,dt|=2\dvol_U$
and 
$|du_1\,d\ubar_1\,du_2\,d\ubar_2|
=4\dvol_{\Utilde}$.
Hence, we obtain the following:
\[
 \varphi_{\ast}\dvol_{\Utilde}
=\frac{\pi}{4(|w|^2+|t|^2)^{1/2}}
 \dvol_{U}.
\]

By direct computations,
we have
\[
 \Ftilde_{u_1,\ubar_1}=
 4u_1u_2\varphi^{\ast}F_{wt}
-4\ubar_1\ubar_2\varphi^{\ast}F_{\wbar t}
+4|u_2|^2F_{w\wbar}
-2|u_1|^2\sqrt{-1}\varphi^{\ast}\nabla_t\phi,
\]
\[
 \Ftilde_{u_2\ubar_2}=
-4u_1u_2\varphi^{\ast}F_{wt}
+4\ubar_1\ubar_2\varphi^{\ast}F_{\wbar t}
+4|u_1|^2F_{w\wbar}
-2|u_2|^2\sqrt{-1}\varphi^{\ast}\nabla_t\phi.
\]
Hence, we obtain 
\[
 \Ftilde_{u_1\ubar_1}
+\Ftilde_{u_2\ubar_2}
=4(|u_1|^2+|u_2|^2)
 \varphi^{\ast}
 \bigl(
 F_{w\wbar}-\frac{\sqrt{-1}}{2}
 \nabla_t\phi
 \bigr)
=\varphi^{\ast}
 \Bigl(2(|w|^2+|t|^2)^{1/2}
 G(h)
 \Bigr).
\]
Thus, we obtain (\ref{eq;17.9.5.2}).
\hfill\qed

\subsection{Appendix: Ambiguity of the choice of a splitting}
\label{subsection;21.8.10.4}

We explain one of the main reasons why we use $(t_1,\beta_1)$
(or $(t_1^{\dagger},\beta_1^{\dagger})$).
Let $Z$ be a finite subset of $\nbigm$.
Let $(E,\nabla,h,\phi)$ be a monopole on $\nbigm\setminus Z$.
Let $F(\nabla)=F(\nabla)_{w\wbar}dw\,d\wbar
+F(\nabla)_{tw}dt\,dw+F(\nabla)_{t\wbar}\,dt\,d\wbar$
denote the curvature of $\nabla$.
Suppose the following conditions.
\begin{condition}
\label{condition;21.8.10.10}\mbox{{}}
\begin{itemize}
 \item  $|F(\nabla)_{w\wbar}|_h=O\bigl(|w|^{-2}(\log|w|)^{-2}\bigr)$
	as $|w|\to\infty$.
	It implies $|\nabla_t\phi|_h=O\bigl(|w|^{-2}(\log|w|)^{-2}\bigr)$
	as $|w|\to\infty$.
 \item There exist positive numbers $R,C,\epsilon>0$,
       a finite subset $S\subset\rnum$ and an orthogonal decomposition
\[
 E=\bigoplus_{\omega\in S}E_{\omega}
\]
on $S^1_T\times\{|w|>R\}$ such that
       $F(\nabla)_{tw}-\bigoplus_{\omega\in S}
       \omega w^{-1}\id_{E_{\omega}}
       =O\bigl(|w|^{-1-\epsilon}\bigr)$ as $|w|\to\infty$.
       Note that it implies that
$F(\nabla)_{t\wbar}+\bigoplus_{\omega\in S}\omega \wbar^{-1}\id_{E_{\omega}}
       =O\bigl(|w|^{-1-\epsilon}\bigr)$ as $|w|\to\infty$.
       It also implies that
       $\nabla_w\phi-\bigoplus (\sqrt{-1}\omega)w^{-1}\id_{E_{\omega}}
       =O\bigl(|w|^{-1-\epsilon}\bigr)$
       and 
       $\nabla_{\wbar}\phi-\bigoplus (\sqrt{-1}\omega)\wbar^{-1}\id_{E_{\omega}}
       =O\bigl(|w|^{-1-\epsilon}\bigr)$
       as $|w|\to\infty$.
\item We assume that $S\neq \{0\}$.
 \end{itemize}
\end{condition}
Let $p_1:M\lrarr \nbigm$ denote the projection.
We set $Z_1:=p_1^{-1}(Z)$.
We set
$(E_1,h_1,\nabla_1,\phi_1):=
p_1^{-1}(E,h,\nabla,\phi)$
on $M\setminus Z_1$.
Let $\lambda\neq 0$.
Let $(E_1^{\lambda},\delbar_{E^{\lambda}_1})$
denote the holomorphic vector bundle on
$M^{\lambda}\setminus Z_1$
underlying $(E_1,h_1,\nabla_1,\phi_1)$.

Let $(t_2,\beta_2)$ be a mini-complex coordinate system
of $M^{\lambda}$.
Let $F^{t_2}_{\beta_2,\betabar_2}d\beta_2\,d\betabar_2$
denote the curvature of the Chern connection of
the holomorphic bundle
$(E_1^{\lambda},\delbar_{E_1^{\lambda}})_{|\{t_2\}\times\cnum_{\beta_2}}$
with the metric induced by $h_1$.
Suppose the following conditions
\begin{condition}
\label{condition;21.8.10.11}\mbox{{}}
 \begin{itemize}
 \item $|F^{t_2}_{\beta_2,\betabar_2}|_h=O\bigl(
       |\beta_2|^{-2}(\log|\beta_2|)^{-2}
       \bigr)$ as $|\beta_2|\to\infty$.
       It is equivalent to
       $|F^{t_2}_{\beta_2,\betabar_2}|_h=O\bigl(
       |w|^{-2}(\log|w|)^{-2}
       \bigr)$
       as $|w|\to\infty$.
 \end{itemize}
\end{condition}
 
\begin{prop}
Either one of the following holds.
\begin{itemize}
 \item There exist a positive number $c_1$ and
       a non-zero complex number $c_2$
       such that
       $(t_2,\beta_2)=c_1(t_1,c_2\beta_1)$.
 \item There exist a positive number $c_1$
       a non-zero complex number $c_2$
       such that
       $(t_2,\beta_2)=c_1(t_1^{\dagger},c_2\beta_1^{\dagger})
       =c_1(t_1-\Image(\lambda^{-1}\beta_1),c_2\lambda^{-2}\beta_1)$.
\end{itemize}
\end{prop}
\pf
After changing $(t_2,\beta_2)$ to $c_1(t_2,\beta_2)$ for some $c_1>0$,
we may assume that
there uniquely exists a linear complex coordinate system
$(\alpha_2,\beta_2)$ of $X^{\lambda}$
such that
(i) the $\real$-action is described as
$s\bullet(\alpha_2,\beta_2)=(\alpha_2+s,\beta_2)$,
(ii) $(\alpha_2,\beta_2)$ induces $(t_2,\beta_2)$ on $M^{\lambda}$.
We have the instanton $(\Etilde,\nablatilde,\htilde)$
on $X\setminus\Ztilde_1$,
corresponding to the monopole $(E_1,h_1,\nabla_1,\phi_1)$,
where $\Ztilde$ denote the pull back of $Z$ by
the projection $X^{\lambda}\lrarr M^{\lambda}$.
Let $(\Etilde^{\lambda},\delbar_{\Etilde^{\lambda}})$
be the holomorphic vector bundle on 
$X^{\lambda}\setminus\Ztilde_1$
underlying $(\Etilde,\htilde,\nablatilde)$.
Note that
$(\Etilde^{\lambda},\delbar_{\Etilde^{\lambda}})
_{|\{\sqrt{-1}t_2\}\times\cnum_{\beta_2}}
\simeq
(E^{\lambda}_1,\delbar_{E^{\lambda}_1})_{|\{t_2\}\times\cnum_{\beta_2}}$.
Hence,
according to Condition \ref{condition;21.8.10.11},
for the expression
\[
       F(\nablatilde)=
       F(\nablatilde)_{\alpha_2\alphabar_2}d\alpha_2d\alphabar_2
       +F(\nablatilde)_{\alpha_2\betabar_2}d\alpha_2d\betabar_2
       +F(\nablatilde)_{\beta_2\alphabar_2}d\beta_2d\alphabar_2
       +F(\nablatilde)_{\beta_2\betabar_2}d\beta_2d\betabar_2,
\]
       we have
       $F(\nablatilde)_{\beta_2\betabar_2|\{\sqrt{-1}t_2\}\times\cnum_{\beta_2}}
       =O\bigl(|w|^{-2}(\log|w|)^{-2}\bigr)$
       as $|\beta_2|\to\infty$.
Note that for the expression
\[
       F(\nablatilde)=
       F(\nablatilde)_{z\zbar}dz\,d\zbar
       +F(\nablatilde)_{z\wbar_2}dz\,d\wbar
       +F(\nablatilde)_{w\zbar}dw\,dz
       +F(\nablatilde)_{w\wbar}dw\,d\wbar,
\]
Condition \ref{condition;21.8.10.10} implies that
$|F(\nablatilde)_{z\zbar}|_{\htilde}=|F(\nablatilde)_{w\wbar}|_{\htilde}
=O\bigl(|w|^{-2}(\log|w|)^{-2}\bigr)$.
Moreover, there exists a decomposition
$\Etilde=\bigoplus \Etilde_{\omega}$ such that
\[
|F(\nablatilde)_{z\wbar}
+\bigoplus C(\omega \wbar^{-1})\id_{E_{\omega}}|_{\htilde}
=O\bigl(|w|^{-1-\epsilon}\bigr).
\]

Because $\beta_2$ is $\real$-invariant,
there exists a non-zero complex number $a$
such that $\beta_2=a(\eta+\lambda\xi)=a\beta_1$.
There exists a complex number $b$ such that
$\alpha_2=\xi+b(\eta+\lambda\xi)$.
We obtain
\[
 z=\frac{1}{1+|\lambda|^2}
 \bigl(\alpha_2-a^{-1}b\beta_2+|\lambda|^2\alphabar_2
 -\lambda\abar^{-1}\betabar_2(1+\lambdabar\bbar)\bigr),
\]
\[
 w=\frac{1}{1+|\lambda|^2}
 \bigl(-\lambda\alpha_2+a^{-1}\beta_2(1+\lambda b)
 +\lambda\alphabar_2-\lambda\abar^{-1}\bbar\betabar_2\bigr).
\]
Hence, the restriction of $dz\,d\wbar$
to $\{\sqrt{-1}t_2\}\times\cnum_{\beta_2}$
is equal to
\[
 b(1+\lambdabar\bbar)\frac{-1}{(1+|\lambda|^2)}|a|^2
 d\beta_2\,d\betabar_2.
\]
Similarly, the restriction of $dw\,d\zbar$
to $\{\sqrt{-1}t_2\}\times\cnum_{\beta_2}$
is equal to
\[
  \bbar(1+\lambda b)\frac{-1}{(1+|\lambda|^2)}|a|^2
 d\beta_2\,d\betabar_2.
\]
We obtain $b(1+\lambdabar\bbar)=0$,
i.e., $b=0$ or $b=-\lambda^{-1}$.
If $b=0$, we obtain $\alpha_2=\xi=\alpha_1$,
and hence $t_2=t_1$.
If $b=-\lambda^{-1}$,
we obtain $\alpha_2=-\lambda^{-1}\eta=-\alpha_1^{\dagger}$,
and hence $t_2=t_1^{\dagger}$.
Thus, we are done.
\hfill\qed

\section{Difference modules and
$\nbigo_{\nbigmbar^{\lambda}\setminus Z}(\ast
  H^{\lambda}_{\infty})$-modules}
\label{subsection;17.10.28.100}

\subsection{Difference modules
   with parabolic structure at finite place}

Let $\varrho\in\cnum$.
Let $\Phi^{\ast}$ be an automorphism of
the rational function field $\cnum(y)$ induced by
$\Phi^{\ast}(y)=y+\varrho$.
\index{automorphism $\Phi^{\ast}$}
Let $\vecV$ be a difference module
as in Definition \ref{df;20.7.29.1},
i.e.,
$\vecV$ is a finite dimensional
$\cnum(y)$-module
equipped with a $\cnum$-linear automorphism
$\Phi^{\ast}$
such that
$\Phi^{\ast}(fs)=\Phi^{\ast}(f)\cdot\Phi^{\ast}(s)$
for $f\in \cnum(y)$ and $s\in \vecV$.
\index{difference module}

Let $V$ be a lattice of $\vecV$,
i.e.,
a $\cnum[y]$-free module $V\subset\vecV$
such that $\cnum(y)\otimes_{\cnum[y]}V=\vecV$.
\index{lattice}
We set
$V':=V+(\Phi^{\ast})^{-1}(V)$.
We obtain
the associated free $\nbigo_{\proj^1}(\ast\infty)$-module
$\nbigf_V$.
We also obtain
the associated free $\nbigo_{\proj^1}(\ast\infty)$-modules
$\nbigf_{(\Phi^{\ast})^{-1}V}$
and $\nbigf_{V'}$.
There exists the natural isomorphism
$\nbigf_V\simeq
 \Phi^{\ast}\nbigf_{(\Phi^{\ast})^{-1}(V)}$.
We may regard 
$\nbigf_V$ and $\nbigf_{(\Phi^{\ast})^{-1}(V)}$
as $\nbigo_{\proj^1}(\ast\infty)$-submodules of
$\nbigf_{V'}$,
and 
$\nbigf_{V'}=
 \nbigf_V+\nbigf_{(\Phi^{\ast})^{-1}(V)}$
holds.

\begin{df}
\index{parabolic structure at finite place}
\label{df;17.12.1.10}
A parabolic structure of $\vecV$ at finite place
is a tuple as follows.
\begin{itemize}
\item 
      A $\cnum[y]$-free submodule $V\subset\vecV$
      such that
      $\cnum(y)\otimes_{\cnum[y]}V=\vecV$.
\item
 A function $m:\cnum\lrarr \seisuu_{\geq 0}$
 such that $\sum_{x\in\cnum} m(x)<\infty$.
 We assume that
 $\nbigf_V(\ast D)=
   \nbigf_{(\Phi^{\ast})^{-1}(V)}(\ast D)$,
where $D:=\{x\in\cnum\,|\,m(x)>0\}$.
\item
 A sequence of real numbers
 $0\leq \tau_{x}^{(1)} <\cdots<\tau_{x}^{(m(x))}<T$
 for each $x\in \cnum$.
 If $m(x)=0$, the sequence is assumed to be empty.
 The sequence is denoted by  $\vectau_{x}$.
\item
 Lattices 
 $L_{x,i}\subset
 V\otimes\cnum(\!(y-x)\!)$
 for $x\in\cnum$ and $i=1,\ldots,m(x)-1$.
 We formally set
 $L_{x,0}:=
 V\otimes\cnum[\![y-x]\!]$
 and 
 $L_{x,m(x)}:=
 (\Phi^{\ast})^{-1}V\otimes\cnum[\![y-x]\!]$.
 The tuple of lattices is denoted by
 $\vecL_{x}$.
\end{itemize}
 Here, $\cnum[\![y-x]\!]$
 denotes the ring of formal power series
 with the variable $y-x$,
 and $\cnum(\!(y-x)\!)$
 denotes the ring of formal Laurent power series
 with the variable $y-x$.
\hfill\qed
\end{df}

\begin{rem}
The notion of parabolic structure
of a difference module at finite place
is apparently different from the ordinary notion of
parabolic structure in the context of harmonic bundles.
However, as we shall explain in {\rm\S\ref{subsection;21.8.6.22}},
a parabolic structure of $2\sqrt{-1}\lambda$-difference module
at finite place
is a reincarnation 
of an ordinary parabolic structure of
the regular part of 
$\lambda$-flat bundle at $\{0,\infty\}$
through the Mellin transform or the algebraic Nahm transform.
\hfill\qed
\end{rem}
 
\subsection{Construction of difference modules from
 $\nbigo_{\nbigmbar^{\lambda}\setminus Z}(\ast H^{\lambda}_{\infty})$-modules}
\label{subsection;17.11.18.2}

Let $Z$ be a finite subset in $\nbigm^{\lambda}$.
Let $\nbige$ be a locally free 
$\nbigo_{\nbigmbar^{\lambda}\setminus Z}
 (\ast H^{\lambda}_{\infty})$-module of Dirac type.
Let us observe that we obtain the associated
difference module with parabolic structure at finite place
with $\varrho=2\sqrt{-1}\lambda T$.

Let $\varpi^{\lambda}:\Mbar^{\lambda}\lrarr\nbigmbar^{\lambda}$
denote the projection.
\index{projection $\varpi^{\lambda}$}
Set $Z^{\cov}:=(\varpi^{\lambda})^{-1}(Z)$.
Let $p_1:\Mbar^{\lambda}\lrarr \real_{t_1}$
and $p_2:\Mbar^{\lambda}\lrarr \proj^1_{\beta_1}$
be the projections.
We obtain 
the $\nbigo_{\Mbar^{\lambda}\setminus Z}(\ast H^{\lambda\cov}_{\infty})$-module
$\nbige^{\cov}:=(\varpi^{\lambda})^{\ast}\nbige$.
Let $j:\Mbar^{\lambda}\setminus Z\lrarr\Mbar^{\lambda}$
denote the inclusion.
We obtain 
the $\nbigo_{\Mbar^{\lambda}}(\ast H^{\lambda\cov}_{\infty})$-module
$j_{\ast}\nbige^{\cov}$.
For any $a\in\real$,
let $\iota_a:p_1^{-1}(a)\lrarr \Mbar^{\lambda}$.
We obtain the locally free $\nbigo_{\proj^1}(\ast\infty)$-modules
$\iota_a^{-1}j_{\ast}\nbige^{\cov}$
$(a\in\real)$.

Because $j_{\ast}\nbige^{\cov}$
is naturally $\seisuu$-equivariant,
we have the following isomorphism for any $a\in\real$
and $n\in\seisuu$:
\begin{equation}
\label{eq;17.9.30.1}
 (\Phi^n)^{\ast}
 \iota_{a+Tn}^{-1}\bigl(
 j_{\ast}(\nbige^{\cov})\bigr)
\simeq
  \iota_{a}^{-1}\bigl(
 j_{\ast}(\nbige^{\cov})\bigr).
\end{equation}
Here, $\Phi:\proj^1\lrarr \proj^1$
is given by $\Phi(\beta_1)=\beta_1+2\sqrt{-1}\lambda T$.

We set
$D_{a,b}:=p_2\Bigl(
 Z^{\cov}\cap\bigl(\{a\leq t_1\leq b\}\times \proj^1_{\beta_1}\bigr)
 \Bigr)$
for any $a\leq b$.
By the scattering map,
we obtain the isomorphism
\begin{equation}
\label{eq;17.9.30.2}
 \iota_{a}^{-1}\bigl(
 j_{\ast}(\nbige^{\cov})
 \bigr)(\ast D_{a,b})
\simeq
 \iota_{b}^{-1}\bigl(
 j_{\ast}(\nbige^{\cov})
 \bigr)(\ast D_{a,b}).
\end{equation}

For any $a\in\real$,
we set
\[
 V_a:=H^0\Bigl(\proj^1,
 \iota_a^{-1}\bigl(j_{\ast}(\nbige^{\cov})\bigr)
 \Bigr).
\]
For any subset $S\subset \cnum$,
let $\cnum[\beta_1]_S$
denote the localization of $\cnum[\beta_1]$
with respect to $\beta_1-x$ $(x\in S)$.
By (\ref{eq;17.9.30.2}),
we obtain  the isomorphism of
$\cnum[\beta_1]_{D_{a,b}}$-modules:
\[
 V_{a}\otimes\cnum[\beta_1]_{D_{a,b}}
\simeq
 V_{b}\otimes\cnum[\beta_1]_{D_{a,b}}.
\]
Hence, we obtain the isomorphism of 
$\cnum(\beta_1)$-modules for any $a\leq b$.
\begin{equation}
\label{eq;17.9.30.10}
 V_{a}\otimes\cnum(\beta_1)
\simeq
 V_{b}\otimes\cnum(\beta_1).
\end{equation}

From (\ref{eq;17.9.30.1}),
we have the $\cnum$-linear isomorphism
\[
 (\Phi^{\ast})^n: V_{a+nT}\simeq V_{a}
\]
such that
$(\Phi^{\ast})^n(\beta_1^{\ell}s)
=(\beta_1+2\sqrt{-1}\lambda nT)^{\ell}(\Phi^{\ast})^n(s)$.
We obtain the $\cnum$-linear isomorphism
\begin{equation}
\label{eq;17.9.30.11}
 (\Phi^{\ast})^n:
 V_{a+nT}\otimes\cnum(\beta_1)
\simeq
 V_{a}\otimes\cnum(\beta_1)
\end{equation}
such that
$(\Phi^{\ast})^n(g(\beta_1)s)
=g(\beta_1+2\sqrt{-1}\lambda nT)\cdot (\Phi^{\ast})^n(s)$
for any $g(\beta_1)\in\cnum(\beta_1)$.

We take a small $\epsilon>0$
such that
$Z^{\cov}\cap
 \bigl(
 \{-\epsilon\leq t_1<0\}\times\proj^1_{\beta_1}
 \bigr)=\emptyset$.
Set $\vecV(\nbige):=V_{-\epsilon}\otimes\cnum(\beta_1)$
which is a finite dimensional
$\cnum(\beta_1)$-vector space.
It is identified with $V_{b}\otimes\cnum(\beta_1)$
for any $b\in\real$ by (\ref{eq;17.9.30.10}).
It is equipped with
$\cnum$-linear automorphism
$\Phi^{\ast}$
by (\ref{eq;17.9.30.11}).

Set $V(\nbige):=V_{-\epsilon}$
for which we obtain
$\cnum(\beta_1)\otimes V(\nbige)
=\vecV(\nbige)$.
Let $m_Z:\cnum\lrarr \seisuu_{\geq 0}$
be the function such that
\begin{equation}
\label{eq;21.9.17.10}
m_Z(x)=
 \Bigl|
 p_2^{-1}(x)\cap 
 p_1^{-1}(\{0\leq t_1<T\})
 \cap Z^{\cov}
 \Bigr|.
\end{equation}
For each $x\in \cnum$,
the sequence
$0\leq t_{1,Z,x}^{(1)}<\cdots <t_{1,Z,x}^{(m_Z(x))}<T$
determined by
\begin{equation}
\label{eq;21.9.17.11}
 \{(t_{1,Z,x}^{(i)},x)\,|\,i=1,\ldots,m_Z(x)\}
=p_2^{-1}(x)\cap 
 p_1^{-1}(\{0\leq t_1<T\})
 \cap Z^{\cov}.
\end{equation}
We set
\begin{equation}
\label{eq;21.9.17.12}
 \tau_{Z,x}^{(i)}:=t_{1,Z,x}^{(i)}/T.
\end{equation}
The tuple
$\tau_{Z,x}^{(i)}$ $(i=1,\ldots,m_Z(x))$
is denoted by $\vectau_{Z,x}$.
For $x\in \cnum$ and for $i=1,\ldots,m_Z(x)-1$,
we choose $t_{1,Z,x}^{(i)}<b(x,i)<t_{1,Z,x}^{(i+1)}$,
and we set
\[
 L_{Z,x,i}(\nbige):=
 V_{b(x,i)}\otimes\cnum[\![\beta_1-x]\!],
\]
which is independent of the choice of $b(x,i)$.
The tuple $L_{Z,x,i}$ $(i=1,\ldots,m_Z(x)-1)$
is denoted by $\vecL_{Z,x}$.
Thus, we obtain a parabolic structure at finite place
\[
 \bigl(V(\nbige),m_Z,
 (\vectau_{Z,x},\vecL_{Z,x}(\nbige))_{x\in\cnum}
 \bigr)
\]
of $\vecV(\nbige)$.

\subsection{Construction of 
$\nbigo_{\nbigmbar^{\lambda}\setminus Z}(\ast H^{\lambda})$-modules
from difference modules}
\label{subsection;17.10.5.100}

Let $\vecV$ be a $(2\sqrt{-1}\lambda T)$-difference module
with a parabolic structure at finite place
$(V,m,(\vectau_{x},\vecL_x)_{x\in\cnum})$.
We set 
$Z_1:=\coprod_{x\in\cnum}
 \bigl\{
 (T\tau_{x}^{(i)},x)\,\big|\,\tau_{x}^{(i)}\in\vectau_{x}
 \bigr\}\subset\real_{t_1}\times\cnum_{\beta_1}$.
 We set $Z:=\varpi^{\lambda}(Z_1)
 \subset\nbigmlambda$.
Let us construct
an $\nbigo_{\nbigmbar^{\lambda}\setminus Z}(\ast H^{\lambda})$-module
$\nbige(\vecV,V,m,(\vectau_{x},\vecL_x)_{x\in\cnum})$.

We set 
$U:=\openopen{-\epsilon}{T-\epsilon/2}\times \proj^1_{\beta_1}$
and 
$H^{\lambda}_{\infty,\epsilon}:=
\openopen{-\epsilon}{T-\epsilon/2}\times\{\infty\}$.
We set $D:=\{x\in\cnum\,|\,m(x)>0\}$.
For $x\in D$,
we set $t_x^{(i)}=T\tau_x^{(i)}$ $(i=1,\ldots,m(x))$,
$t_x^{(0)}=-\epsilon/2$
and $t_x^{(m(x)+1)}=T-\epsilon$.
The following conditions determine
a unique 
$\nbigo_{U\setminus Z_1}(\ast H^{\lambda}_{\infty,\epsilon})$-module
$\nbige_1$:
\begin{itemize}
\item
 $\nbige_1(\ast p_2^{-1}(D))$ is isomorphic to
 the pull back of $\nbigf_V(\ast D)$
 by the projection $U\lrarr \proj^1_{\beta_1}$.
\item
 Take $t_x^{(i)}<a_i<t_x^{(i+1)}$
 for $i=0,\ldots,m(x)$.
 Take a small neighbourhood $U_x$ of $x$ in $\cnum$.
 Let $\iota_{\{a_i\}\times U_x}$ denote the inclusion
 of $\{a_i\}\times U_x$ to $U\setminus Z_1$.
 Then, 
 $\iota_{\{a_i\}\times U_x}^{\ast}\nbige_1
 \subset
  \nbigf_V(\ast D)_{|U_x}$
 is equal to the $\nbigo_{U_x}$-submodule
 determined by $L_{x,i}$.
\end{itemize}

Let $\kappa_1\colon
 \openopen{-\epsilon}{-\epsilon/2}\times \proj^1
\lrarr 
 \openopen{T-\epsilon}{T-\epsilon/2}\times\proj^1$
be the isomorphism
defined by
$\kappa_1(t_1,\beta_1)=(t_1+T,\beta_1+2\sqrt{-1}\lambda T)$.
By the construction,
there exists the isomorphism
\[
 \kappa_1^{\ast}
 \Bigl(
 \nbige_{1|\openopen{T-\epsilon}{T-\epsilon/2}\times\proj^1}
 \Bigr)
\simeq
 \nbige_{1|\openopen{-\epsilon}{-\epsilon/2}\times\proj^1}.
\]
Hence, we obtain
a locally free 
$\nbigo_{\nbigmbar^{\lambda}\setminus Z}
 (\ast H^{\lambda}_{\infty})$-module
$\nbige(\vecV,V,m,(\tau_x,\vecL_x)_{x\in\cnum})$ of Dirac type.
By the construction,
the following is clear.
\begin{lem}
\label{lem;17.12.4.2}
The constructions in {\rm\S\ref{subsection;17.11.18.2}}
and {\rm\S\ref{subsection;17.10.5.100}}
are mutually inverse up to canonical isomorphisms.
\hfill\qed 
\end{lem}

\begin{rem}
The relation between 
finitely generated difference modules
and 
$\nbigo_{\nbigmbar^{\lambda}\setminus Z}(\ast
 H^{\lambda}_{\infty})$-modules
of Dirac type is analogue to
the relation between 
$\nbigo_C(\ast D)$-modules
and filtered bundles over $(C,D)$
for a compact Riemann surface $C$
with a finite subset $D\subset C$.
\hfill\qed
\end{rem}

\subsection{Appendix: Mellin transform and parabolic structures at finite place}
\label{subsection;21.8.6.22}

Let us explain that
the notion of parabolic structure at finite place of difference modules
is related with
the usual notion of parabolic structure of filtered
$\lambda$-flat bundles
through the algebraic Nahm transform,
which is an analogue of Mellin transform.

\subsubsection{Mellin transform}
\label{subsection;21.8.6.2}
\index{Mellin transform}

Let $\lambda$ be any complex number.
Let $M$ be a $\cnum[u,u^{-1}]$-module
equipped with a $\cnum$-linear endomorphism
$\DDlambda_u:M\lrarr M$
such that
$\DDlambda_u(fs)=f\DDlambda_u(s)+\lambda\del_u(f)s$.

We define the automorphism $\Phi^{\ast}$
of $\cnum$-algebra $\cnum[\beta_1]$
by $\Phi^{\ast}(\beta_1)=\beta_1+2\sqrt{-1}\lambda$.
We define
$\cnum[\beta_1]$-action on $M$
by
$\beta_1s=-2\sqrt{-1}u\DDlambda_u(s)$
for any $s\in M$.
We define the $\cnum$-linear automorphism $\Phi^{\ast}_M$ of $M$
by
$\Phi^{\ast}(s)=us$ for any $s\in M$.
Because $\Phi_M^{\ast}(\beta_1s)=\Phi^{\ast}(\beta_1)\Phi_M^{\ast}(s)$,
we obtain $2\sqrt{-1}\lambda$-difference module,
which is called the Mellin transform.

\vspace{.1in}
We have the natural geometrization of the transformation
of the Mellin transform.
Let $p_i$ $(i=1,2)$ denote the projections of
$\proj^1_u\times\cnum_{\beta_1}$
onto the $i$-th component.
Let $\nbigv$ be
an algebraic quasi-coherent $\nbigo_{\proj^1}(\ast\{0,\infty\})$-module
with a $\lambda$-connection $\DDlambda$.
We obtain the following complex $\nbigc^{\bullet}(\nbigv,\DDlambda)$
on $\proj^1_u\times\cnum_{\beta_1}$:
\[
\begin{CD}
 p_1^{\ast}\nbigv
 @>{u\DDlambda_u-\frac{\sqrt{-1}}{2}\beta_1}>>
 p_1^{\ast}\nbigv,
\end{CD}
\]
where the first term sits in the degree $0$.
We have $R^ip_{2\ast}(\nbigc^{\bullet}(\nbigv,\DDlambda))=0$
unless $i=1$.
We obtain the algebraic $\nbigo_{\cnum_{\beta_1}}$-module
$\gbigm(\nbigv,\DDlambda):=
R^1p_{2\ast}(\nbigc^{\bullet}(\nbigv,\DDlambda))$.

Let $f_u:\nbigv\lrarr\nbigv$
be the automorphism defined as the multiplication of $u$.
It induces the isomorphism
$(\nbigv,\DDlambda+\lambda du/u)\simeq (\nbigv,\DDlambda)$.

Let $\Phi:\cnum_{\beta_1}\lrarr\cnum_{\beta_1}$
be defined by $\Phi(\beta_1)=\beta_1+2\sqrt{-1}\lambda$.
We obtain the following morphisms:
\[
\begin{CD}
 (\id\times\Phi)^{\ast}\nbigc^{\bullet}(\nbigv,\DDlambda)
 \simeq
 \nbigc^{\bullet}(\nbigv,\DDlambda+\lambda du/u)
 @>{p_1^{\ast}f_u}>>
 \nbigc^{\bullet}(\nbigv,\DDlambda).
\end{CD}
\]
Hence, we obtain the isomorphism
$\Phi^{\ast}\gbigm(\nbigv,\DDlambda)
\simeq
\gbigm(\nbigv,\DDlambda)$.

If $M=H^0(\proj^1,\nbigv)$,
then $H^0(\proj^1,\gbigm(\nbigv,\DDlambda))$
is the Mellin transform.

\begin{rem}
The stationary phase formula for Mellin transform
has been studied in
{\rm\cite{Garcia-Lopez}} and {\rm\cite{Graham-Squire}}.
\hfill\qed
\end{rem}

\subsubsection{Algebraic Nahm transform
for filtered $\lambda$-flat bundles (special case)}

Let $D\subset \cnum^{\ast}$ be a finite subset.
We set $\Dbar=D\cup\{0,\infty\}$.
Let $\nbigv$ be a locally free
$\nbigo_{\proj^1}(\ast\Dbar)$-module of finite rank
equipped with a $\lambda$-connection $\DDlambda$.
Let $(\nbigp_{\ast}\nbigv,\DDlambda)$ be
a good filtered $\lambda$-flat bundle on $(\proj^1,\Dbar)$.
For simplicity, we assume the following.
\begin{condition}\mbox{}
\label{condition;21.8.6.30}
 \begin{itemize}
 \item $D$ is non-empty.
       Moreover, at each point $P$ of $D$,
       $(\nbigp_{\ast}\nbigv,\DDlambda)$ does not have the regular part,
       i.e.,
       in the Hukuhara-Levelt-Turrittin decomposition {\rm(\ref{eq;21.8.6.1})}
       of the pull back by an appropriate ramified covering of
       a neighborhood of $P$,
       we have $\nbigp_{\ast}\nbigvhat_{0}=0$.
  \item $\nbigp_{\ast}\nbigv$ is regular at $0$ and $\infty$,
	i.e.,
       $\DDlambda$ is logarithmic with respect to
	$\nbigp_{\veca}(\nbigv)$ for any $\veca\in\real^{\Dbar}$
	at $0$ and $\infty$.
       \hfill\qed
 \end{itemize} 
\end{condition}

We shall construct
an $\nbigo_{\nbigmbar^{\lambda}\setminus Z}(\ast H^{\lambda})$-module
from $(\nbigp_{\ast}\nbigv,\DDlambda)$,
where $Z$ denotes a finite subset explained below.
We assume $T=1$ in the construction of $\nbigm^{\lambda}$.

Let $\nbigv_0$ denote the stalk of $\nbigv$ at $0$.
We have the filtration $\nbigp_a(\nbigv_0)$ $(a\in\real)$.
We set
$\Gr^{\nbigp}_a(\nbigv_0)
=\nbigp_a(\nbigv_0)\big/\nbigp_{<a}(\nbigv_0)$.
Because $\DDlambda$ is logarithmic at $0$,
we obtain the endomorphism
$\Res_0(\DDlambda)$
of $\Gr^{\nbigp}_a(\nbigv_0)$
as the residue of $\DDlambda$.
Let 
$\EE_{\alpha} \Gr^{\nbigp}_a(\nbigv_0)$
denote the generalized eigen space corresponding to
$\alpha\in\cnum$.
Let $\KMS(\nbigp_{\ast}\nbigv,0)$
denote the set of $(a,\alpha)\in\real\times\cnum$
such that
$\EE_{\alpha} \Gr^{\nbigp}_a(\nbigv_0)\neq 0$.
Similarly,
let $\nbigv_{\infty}$ denote the stalk of $\nbigv$ at $\infty$,
which is equipped with the filtration $\nbigp_{\ast}(\nbigv_{\infty})$
such that $\DDlambda$ is logarithmic.
We obtain
$\Gr^{\nbigp}_a(\nbigv_{\infty}):=
\nbigp_a(\nbigv_{\infty})\big/\nbigp_{<a}(\nbigv_{\infty})$
which is equipped with the endomorphism
$\Res_{\infty}(\DDlambda)$
obtained as the residue of $\DDlambda$.
Let $\EE_{\alpha}\Gr^{\nbigp}_a(\nbigv_{\infty})$
denote the generalized eigen space corresponding to
$\alpha\in\cnum$.
We obtain the set
$\KMS(\nbigp_{\ast}\nbigv,\infty)$
of $(a,\alpha)\in\real\times\cnum$
such that
$\EE_{\alpha} \Gr^{\nbigp}_a(\nbigv_{\infty})\neq 0$.

Let $\Ztilde\subset\real\times\cnum$ denote the union of the following sets
\[
 \Bigl\{(t_1,\beta_1)\in\real\times\cnum\,\Big|\,
  \Bigl(t_1,\frac{\sqrt{-1}}{2}\beta_1\Bigr)
  \in \KMS(\nbigp_{\ast}\nbigv,0)
  \Bigr\}
\]
\[
 \Bigl\{
 (t_1,\beta_1)\in\real\times\cnum\,\Big|\,
  \Bigl(-t_1,-\frac{\sqrt{-1}}{2}\beta_1\Bigr)
  \in \KMS(\nbigp_{\ast}\nbigv,\infty)
 \Bigr\}.
\]
Note that $\Ztilde$ is preserved by the $\seisuu$-action $\kappa$
on $\real\times\cnum$
defined as $\kappa_n(t_1,\beta_1)=(t_1+n,\beta_1+2\sqrt{-1}\lambda n)$.
Let $Z\subset\nbigm^{\lambda}$ denote the quotient of $\Ztilde$.

Let $\pitilde:\real\times\proj^1\lrarr \real$ denote the projection.
We obtain the set $\pitilde(\Ztilde)\subset\real$.

For any $(a,b)\in\real^2$,
let $\nbigp_{a,b}(\nbigv)\subset\nbigv$ be
the $\nbigo_{\proj^1}$-module such that
(i) $\nbigp_{a,b}(\nbigv)(\ast\{0,\infty\})=\nbigv$,
(ii) the stalk of $\nbigp_{a,b}(\nbigv)$ at $0$
is $\nbigp_a(\nbigv_0)$,
(iii) the stalk of $\nbigp_{a,b}(\nbigv)$ at $\infty$
is $\nbigp_{b}(\nbigv_{\infty})$.

Let $p_i$ denote the projections of $\proj^1_u\times \proj^1_{\beta_1}$
onto the $i$-th components.
We obtain the following complex
$\nbigc^{\bullet}_{a,b}$
on $\proj^1_u\times \proj^1_{\beta_1}$:
\[
\begin{CD}
 p_1^{\ast}\nbigp_{a,b}\nbigv(\ast(\proj^1\times\{\infty\}))
 @>{u\DDlambda_u-\frac{\sqrt{-1}}{2}\beta_1}>>
 p_1^{\ast}\nbigp_{a,b}\nbigv(\ast(\proj^1\times\{\infty\}))
\end{CD}
\]
where the first term sits in the degree $0$.
We obtain
$R^ip_{2\ast}\nbigc^{\bullet}_{a,b}=0$
unless $i=1$,
and
$\gbigntilde_{a,b}(\nbigp_{\ast}\nbigv,\DDlambda):=
R^1p_{2\ast}\nbigc^{\bullet}_{a,b}$
is a locally free $\nbigo_{\proj^1}(\ast\infty)$-module.

Let $t_1\in\real$.
Suppose $t_1\not\in \pitilde(\Ztilde)$.
There exists $\epsilon>0$
such that $\{t_1-\epsilon\leq a\leq t_1+\epsilon\}\cap\pitilde(\Ztilde)=\emptyset$.
For any $t_1-\epsilon\leq t_1'\leq t_1+\epsilon$,
we have
$\nbigp_{t_1,-t_1}\nbigv=\nbigp_{t_1',-t_1'}\nbigv$,
and hence there exists a natural isomorphism
\begin{equation}
\label{eq;21.8.6.10}
 \gbigntilde_{t_1,-t_1}(\nbigp_{\ast}\nbigv,\DDlambda)
 \simeq
 \gbigntilde_{t_1',-t_1'}(\nbigp_{\ast}\nbigv,\DDlambda).
\end{equation}
Suppose $t_1\in\pitilde(\Ztilde)$.
There exists $\epsilon>0$
such that $\{t_1-\epsilon\leq a\leq t_1+\epsilon\}\cap\pitilde(\Ztilde)=\{t_1\}$.
For any $t_1-\epsilon\leq t_1'\leq t_1+\epsilon$,
we have
$\nbigp_{t_1-\epsilon,-t_1-\epsilon}\nbigv
\subset
 \nbigp_{t'_1,-t'_1}\nbigv$,
and hence there exists the following morphism
$\gbigntilde_{t_1-\epsilon,-t_1-\epsilon}(\nbigp_{\ast}\nbigv,\DDlambda)
 \lrarr
 \gbigntilde_{t'_1,-t_1'}(\nbigp_{\ast}\nbigv,\DDlambda)$.

\begin{lem}
\label{lem;21.8.6.12}
We set $\Ztilde_{t_1}:=\pitilde^{-1}(t_1)\cap \Ztilde$,
which we naturally regard as a subset of $\cnum$.
Then, the induced morphism
\begin{equation}
 \label{eq;21.8.6.11}
  \gbigntilde_{t_1-\epsilon,-t_1-\epsilon}(\nbigp_{\ast}\nbigv,\DDlambda)
  (\ast\Ztilde_{t_1})
 \lrarr
 \gbigntilde_{t'_1,-t_1'}(\nbigp_{\ast}\nbigv,\DDlambda)(\ast\Ztilde_{t_1})
\end{equation}
is an isomorphism.
\end{lem}
\pf
If $\alpha$ is an eigenvalues of $\Res_0(\DDlambda)$
on $\Gr^{\nbigp}_{t_1}(\nbigv_0)$,
then $\frac{\sqrt{-1}}{2}\alpha$
is contained in $\Ztilde_{t_1}$.
If $\alpha$ is an eigenvalues of $\Res_{\infty}(\DDlambda)$
on $\Gr^{\nbigp}_{-t_1}(\nbigv_{\infty})$,
then $-\frac{\sqrt{-1}}{2}\alpha$
is contained in $\Ztilde_{t_1}$.

Let $\beta_1\in\cnum\setminus\Ztilde_{t_1}$.
There exists a neighbourhood $U$ of $\beta_1$
in $\cnum\setminus\Ztilde_{t_1}$.
By the consideration in the previous paragraph,
$\Res_0(\DDlambda)-\frac{\sqrt{-1}}{2}\beta_1'$ $(\beta_1'\in U)$
are invertible
on
$\Gr^{\nbigp}_{t_1}(\nbigv_0)$,
and 
$\Res_{\infty}(\DDlambda)+\frac{\sqrt{-1}}{2}\beta_1'$
$(\beta_1'\in U)$
are invertible
on
$\Gr^{\nbigp}_{-t_1}(\nbigv_{\infty})$.
Hence, we can easily check the claim of the lemma.
\hfill\qed

\vspace{.1in}
By Lemma \ref{lem;21.8.6.12},
we obtain the following isomorphism
for any $t_1-\epsilon\leq t_1'\leq t_1+\epsilon$:
\begin{equation}
\label{eq;21.8.6.13}
 \gbigntilde_{t_1,-t_1}(\nbigp_{\ast}\nbigv,\DDlambda)_{|\proj^1\setminus\Ztilde_{t_1}}
 \simeq
 \gbigntilde_{t'_1,-t'_1}(\nbigp_{\ast}\nbigv,\DDlambda)_{|\proj^1\setminus\Ztilde_{t_1}}.
\end{equation}

By the isomorphisms (\ref{eq;21.8.6.10})
and (\ref{eq;21.8.6.13}),
we obtain
an $\nbigo_{(\real\times\proj^1)\setminus\Ztilde}
 (\ast(\real\times\{\infty\}))$-module
$\gbigntilde(\nbigp_{\ast}\nbigv)$
from
$\gbigntilde_{t_1,-t_1}(\nbigp_{\ast}\nbigv)_{|\proj^1\setminus\Ztilde_{t_1}}$
$(t_1\in\real)$.
The isomorphism $f_u$ on $\nbigv$
induces
$\kappa_1^{\ast}\gbigntilde(\nbigp_{\ast}\nbigv,\DDlambda)\simeq
 \gbigntilde(\nbigp_{\ast}\nbigv,\DDlambda)$.
Hence, we obtain
an $\nbigo_{\nbigmbar^{\lambda}\setminus Z}(\ast H^{\lambda}_{\infty})$-module
$\gbign(\nbigp_{\ast}\nbigv)$.

\vspace{.1in}

Let $(t_1^0,\beta_1^0)\in\Ztilde$.
There exist
integers $k_1\geq \ldots\geq k_r$
and a frame
$v_1,\ldots,v_r$ of
$\gbigntilde(\nbigp_{\ast}\nbigv,\DDlambda)_{(t_1^0-\epsilon,\beta_1)}$
such that
the tuple $(\beta_1-\beta_1^0)^{k_i}v_i$ $(i=1,\ldots,r)$
is a frame of
$\gbigntilde(\nbigp_{\ast}\nbigv,\DDlambda)_{(t_1^0+\epsilon,\beta_1)}$.

We set $\alpha=\frac{\sqrt{-1}}{2}\beta_1^0$.
We obtain the filtration $W$
on the spaces
$\EE_{\alpha}\Gr^{\nbigp}_{t_1^0}(\nbigv_0)$
and
$\EE_{-\alpha}\Gr^{\nbigp}_{-t_1^0}(\nbigv_{\infty})$
obtained as the weight filtration of the nilpotent parts of
$\Res_0(\DDlambda)$ and $\Res_{\infty}(\DDlambda)$,
respectively.
For $k\geq 0$,
let $P\Gr^W_k\EE_{\alpha}\Gr^{\nbigp}_{t_1^0}(\nbigv_0)$
and $P\Gr^W_k\EE_{-\alpha}\Gr^{\nbigp}_{-t_1^0}(\nbigv_{\infty})$
denote the primitive parts.

\begin{prop}
$\dim P\Gr^W_k\EE_{\alpha}\Gr^{\nbigp}_{t_1^0}(\nbigv_0)$
is equal to the number of
$v_i$ such that $k_i=k+1$,
and 
$\dim P\Gr^W_k\EE_{-\alpha}\Gr^{\nbigp}_{-t_1^0}(\nbigv_{\infty})$
is equal to the number of
$v_i$ such that $k_i=-k-1$. 
\end{prop}
\pf
It is enough to consider the case
$(t_1^0,\beta_1^0)=(0,0)$.
Let $c_1:\nbigp_{-\epsilon,-\epsilon}\nbigv
\lrarr\nbigp_{-\epsilon,\epsilon}\nbigv$
and
$c_2:\nbigp_{-\epsilon,-\epsilon}\nbigv
\lrarr\nbigp_{\epsilon,-\epsilon}\nbigv$
denote the natural inclusions.
Let $d_1:\nbigp_{-\epsilon,\epsilon}\nbigv
\lrarr\nbigp_{\epsilon,\epsilon}\nbigv$
and $d_2:\nbigp_{\epsilon,-\epsilon}\nbigv
\lrarr\nbigp_{\epsilon,\epsilon}\nbigv$
denote the natural inclusions.
We obtain the following exact sequence:
\[
0\lrarr
\nbigp_{-\epsilon,-\epsilon}\nbigv
\stackrel{c_1\oplus c_2}{\lrarr}
\nbigp_{\epsilon,-\epsilon}\nbigv
\oplus
\nbigp_{-\epsilon,\epsilon}\nbigv
\stackrel{d_1-d_2}{\lrarr} 
\nbigp_{\epsilon,\epsilon}\nbigv
\lrarr 0.
\]
It induces the following exact sequence
\begin{multline}
 0\lrarr
 \gbigntilde_{-\epsilon,-\epsilon}(\nbigp_{\ast}\nbigv,\DDlambda)
 \lrarr
 \gbigntilde_{\epsilon,-\epsilon}(\nbigp_{\ast}\nbigv,\DDlambda)
 \oplus
 \gbigntilde_{-\epsilon,\epsilon}(\nbigp_{\ast}\nbigv,\DDlambda)
 \lrarr \\
 \gbigntilde_{\epsilon,\epsilon}(\nbigp_{\ast}\nbigv,\DDlambda)
 \lrarr 0.
\end{multline}

On $\cnum$,
we consider the following complex:
\[
\begin{CD}
\Gr^{\nbigp}_{0}(\nbigv_0)
\otimes\nbigo_{\cnum}
@>{\Res_0(\DDlambda)-\frac{\sqrt{-1}}{2}\beta_1}>>
\Gr^{\nbigp}_{0}(\nbigv_0)
\otimes\nbigo_{\cnum}
\end{CD}
\]
The quotient is isomorphic to
the quotient of
$\gbigntilde_{-\epsilon,-\epsilon}(\nbigp_{\ast}\nbigv,\DDlambda)
\lrarr
\gbigntilde_{\epsilon,-\epsilon}(\nbigp_{\ast}\nbigv,\DDlambda)$.
We have a similar description of the quotient of
$\gbigntilde_{-\epsilon,-\epsilon}(\nbigp_{\ast}\nbigv,\DDlambda)
\lrarr
\gbigntilde_{-\epsilon,\epsilon}(\nbigp_{\ast}\nbigv,\DDlambda)$.
Then, we can easily deduce the claim of the proposition
from Lemma \ref{lem;21.8.6.20} below.
\hfill\qed

\begin{lem}
\label{lem;21.8.6.20}
Let $V$ be an $r$-dimensional vector space
with a base $e_1,\ldots,e_r$.
Let $f$ be the endomorphism of $V$
determined by $f(e_i)=e_{i+1}$ $(i=1,\ldots,r-1)$
and $f(e_r)=0$.
Then, the cokernel
$f-z\id_V:V\otimes\nbigo_{\cnum_z}\lrarr V\otimes\nbigo_{\cnum_z}$
is isomorphic to
$\nbigo_{\cnum}/z^r\nbigo_{\cnum}$.
\end{lem}
\pf
We set $F=f-z\id_V$.
Let $V'$ denote the vector space
$\bigoplus\cnum_{i=2}^{r} e_i$,
and let $\rho:V\lrarr V'$ denote the projection.
It is easy to see that $\rho\circ F$ is an epimorphism,
and the kernel of $\rho\circ F$
is generated by
$\sum_{j=1}^r z^{r-j}e_j$.
Then, the claim of the lemma easily follows.
\hfill\qed

\begin{rem}
We can consider this kind of transformation for
stable good filtered $\lambda$-flat bundles of degree $0$
which do not necessarily satisfy Condition {\rm\ref{condition;21.8.6.30}}. 
It will be explained elsewhere.
The transformation naturally appear as
the algebraic counterpart of 
the Nahm transform between
periodic monopoles
and wild harmonic bundles on $(\proj^1,\Dbar)$.
Some special cases were 
studied in \cite{Cherkis-Kapustin1, Cherkis-Kapustin2, Harland},
and a more general case is
studied in {\rm\cite{Mochizuki-periodic-Nahm}}.
\hfill\qed
\end{rem}

\section{Filtered prolongation of acceptable bundles}
\label{subsection;17.10.5.122}

\subsection{Filtered bundles on a neighbourhood of $0$ in $\cnum$}
\label{subsection;20.8.1.40}
\label{subsection;17.10.5.301}

\index{filtered bundle}
Let $Y$ be a neighbourhood of $0$ in $\cnum$.
Let $\nbige$ be a locally free $\nbigo_Y(\ast 0)$-module
of finite rank.
A filtered bundle over $\nbige$ is
an increasing sequence of 
locally free $\nbigo_Y$-submodules
$\nbigp_a\nbige\subset\nbige$
$(a\in\real)$
satisfying the following conditions.
\begin{itemize}
 \item
      $\nbigp_a\nbige$ $(a\in\real)$ are lattices of
      $\nbige$, i.e.,
      $\nbigp_a\nbige(\ast \{0\})=\nbige$.
\item
 $\nbigp_{a+n}\nbige=\nbigp_a\nbige(n\{0\})$
 for any $a\in\real$ and $n\in\seisuu$.
\item
 For any $a\in\real$,
 there exists $\epsilon>0$
 such that
 $\nbigp_{a+\epsilon}\nbige=\nbigp_a\nbige$.
\end{itemize}
In that case,
we also say that $\nbigp_{\ast}\nbige$ is a filtered bundle
on $(Y,0)$, for simplicity.

For any $a\in\real$,
we set
$\nbigp_{<a}\nbige:=
\sum_{b<a}\nbigp_b(\nbige)$,
and
\[
 \Gr^{\nbigp}_{a}(\nbige):=
 \nbigp_a(\nbige)\big/
 \nbigp_{<a}(\nbige).
\]
We may naturally regard
$\Gr^{\nbigp}_{a}(\nbige)$
as a finite dimensional $\cnum$-vector space.
\index{vector space $\Gr^{\nbigp^h}_a(\nbige)$}

A frame $\vecv=(v_1,\ldots,v_{\rank\nbige})$ of $\nbigp_a\nbige$
is called compatible with the parabolic structure
if there exists a decomposition
$\vecv=\coprod_{a-1<b\leq a}\vecv_b$
such that the following holds.
\index{compatible frame}
\begin{itemize}
 \item $\vecv_b$ is a tuple of sections of
       $\nbigp_b(\nbige)$,
       and $\vecv_b$ induces a base of
       $\Gr^{\nbigp}_b(\nbige)$.
\end{itemize}
For any non-zero section $s$,
the number
\index{degree $\deg^{\nbigp}(s)$}
\begin{equation}
\label{eq;20.8.8.50}
 \deg^{\nbigp}(s):=\min\bigl\{
 c\in\real\,\big|\,
 s\in\nbigp_c\nbige
 \bigr\}
\end{equation}
is called the parabolic degree of $s$.
If $s=0$, we set
$\deg^{\nbigp}(s)=-\infty$.

\subsubsection{$G$-equivariance}

Let $G$ be a finite group acting on $Y$
preserving $0$,
which is not necessarily effective.
We say that a filtered bundle $\nbigp_{\ast}\nbige$
on $(Y,0)$ is $G$-equivariant
if $\nbige$ is $G$-equivariant,
and each $\nbigp_a\nbige$ is preserved by the $G$-action.
\index{equivariant filtered bundle}

For a $G$-equivariant filtered bundle $\nbigp_{\ast}\nbige$
on $(Y,0)$,
each $H^0(Y,\nbigp_a\nbige)$ is naturally
a $G$-representation.
The following lemma is obvious because $G$ is finite.
\begin{lem}
\label{lem;21.8.25.1}
There exist $G$-invariant subspaces
$H_b\subset H^0(Y\,\nbigp_b\nbige)$ $(b\in\real)$
such that
the natural morphism
$H^0(Y,\nbigp_b\nbige)\lrarr \Gr^{\nbigp}_b(\nbige)$ 
induces an isomorphism
$H_b\simeq \Gr^{\nbigp}_b(\nbige)$
of $G$-representations.
\hfill\qed 
\end{lem}

\begin{cor}
\label{cor;21.8.25.10}
Let $\Gr^{\nbigp}_b(\nbige)=\bigoplus_{i=1}^{m(b)} V_{b,i}$
$(a-1<b\leq a)$
be decompositions of $G$-representations.
Then, there exist a $G$-invariant neighbourhood $Y'$ of $0$ in $Y$
and a $G$-equivariant decomposition
$\nbigp_{\ast}\nbige_{|Y'}=\bigoplus_{a-1<b\leq a}\bigoplus_{i=1}^{m(b)}
\nbigp_{\ast}\nbige_{b,i}$
such that the following holds.
\begin{itemize}
 \item We have $\Gr^{\nbigp}_{c}(\nbige_{b,i})=0$
       for $a-1<c\leq a$ unless $c=b$,
       and
       $\Gr^{\nbigp}_b(\nbige_{b,i})=V_{b,i}$.
\end{itemize} 
\end{cor}
\pf
Let $H_b$ be as in Lemma \ref{lem;21.8.25.1}.
We obtain the induced decompositions
of $G$-representations
$H_b=\bigoplus_{i=1}^{m(i)} H_{b,i}$.
Let $\nbige'_{b,i}$ denote the $\nbigo_Y(\ast 0)$-submodule
of $\nbige$ generated by $H_{b,i}$.
There exists a $G$-invariant neighbourhood $Y'$ of $0$ in $Y$
such that
$\nbige_{|Y'}=\bigoplus \nbige'_{b,i|Y'}$.
Then, the claim of the lemma is easy to check.
\hfill\qed

\begin{cor}
\label{cor;21.8.25.3}
 If $G$ is cyclic,
there exist a $G$-invariant neighbourhood $Y'$ of $0$ in $Y$
and a compatible frame $\vecv$ of $\nbigp_a\nbige_{|Y'}$
such that
$g^{\ast}(v_i)=g^{n(i)}v_i$ for some $n(i)\in\seisuu$.
\index{equivariant compatible frame}
\hfill\qed
\end{cor}

The following lemma is obvious.
\begin{lem}
Suppose that the action of $G$ on $Y$ is trivial.
Let $\Irr(G)$ denote the set of
the isomorphism classes of the irreducible representations of $G$.
There exists a unique $G$-equivariant decomposition
$\nbigp_{\ast}\nbige=\bigoplus_{\kappa\in\Irr(G)}
\nbigp_{\ast}\nbige_{\kappa}$
such that for any $y\in Y$ and for any $a\in\real$,
the $G$-representation
$\nbigp_a(\nbige_{\kappa})_{|y}$
is isomorphic to a direct sum of $\kappa$.
\hfill\qed
\end{lem}

\subsubsection{Subbundles, quotient and splitting}
\label{subsection;21.8.26.51}

Let $\nbigp_{\ast}\nbige$ be a $G$-equivariant filtered bundle
on $(Y,0)$.
Recall that 
a $G$-equivariant locally free $\nbigo_Y(\ast 0)$-submodule
$\nbige'\subset\nbige$ is called saturated
if $\nbige''=\nbige/\nbige'$ is also locally free.
\index{saturated}
In that case, we define
$\nbigp_a(\nbige'):=\nbigp_a(\nbige)\cap\nbige'$
for any $a$,
and we obtain the $G$-equivariant induced filtered bundle over $\nbige'$.
\index{induced filtered bundle (subbundle)}
Let $\nbigp_a(\nbige'')$ be the image of
$\nbigp_a(\nbige)\lrarr\nbige''$,
and we obtain the $G$-equivariant induced filtered bundle over $\nbige''$.
\index{induced filtered bundle (quotient bundle)}

\begin{lem}
\label{lem;21.8.24.1}
There exists a $G$-invariant neighbourhood $Y'$ of $0$ in $Y$
and 
a splitting $\nbige''_{|Y'}\lrarr\nbige_{|Y'}$
with which we obtain
$\nbigp_{\ast}(\nbige)_{|Y'}\simeq
\nbigp_{\ast}(\nbige')_{|Y'}\oplus\nbigp_{\ast}(\nbige'')_{|Y'}$.
\end{lem}
\pf
We may naturally regard
$\Gr^{\nbigp}_b(\nbige')$ as
a $G$-invariant subspace of $\Gr^{\nbigp}_b(\nbige)$.
Because $G$ is finite,
there exists a $G$-invariant decomposition
$\Gr^{\nbigp}_b(\nbige)
=\Gr^{\nbigp}_b(\nbige')\oplus V_b''$,
and there exists a $G$-invariant subspace
$H''_b\subset H^0(Y,\nbigp_b\nbige)$
such that the map
$H^0(Y,\nbigp_b\nbige)\lrarr\Gr^{\nbigp}_b(\nbige)$
induces an isomorphism $H''_b\simeq V_b''$
of $G$-representations.
The composite of the natural morphisms
$H^0(Y,\nbigp_b\nbige)
\lrarr
H^0(Y,\nbigp_b\nbige'')\lrarr
\Gr^{\nbigp}_b(\nbige'')$
induces an isomorphism
$H''_b\simeq \Gr^{\nbigp}_b(\nbige'')$
of $G$-representations.

Let $\nbige''_1$ be the $\nbigo_Y(\ast 0)$-submodule
of $\nbige$ generated by $\bigoplus_{a-1<b\leq a} H''_b$.
There exists a $G$-invariant neighbourhood $Y'$ of $0$ in $Y$
such that
$\nbige_{|Y'}=\nbige'_{|Y'}\oplus\nbige''_{1|Y'}$.
Then, it is easy to check the claims of the lemma.
\hfill\qed

\subsubsection{Basic functoriality}
\label{subsection;21.8.15.1}

Let $\nbigp_{\ast}\nbige_i$ $(i=1,2)$
be filtered bundles
over locally free $\nbigo_Y(\ast 0)$-modules $\nbige_i$.
Let $\nbigp_{\ast}\nbige_1\oplus\nbigp_{\ast}\nbige_2$
denote the filtered bundle
$\nbigp_{\ast}(\nbige_1\oplus\nbige_2)$
over $\nbige_1\oplus\nbige_2$
obtained as
\[
\nbigp_a(\nbige_1\oplus\nbige_2)
=\nbigp_a(\nbige_1)\oplus\nbigp_a(\nbige_2). 
\]
\index{direct sum (filtered bundle)}
Let $\nbigp_{\ast}\nbige_1\otimes\nbigp_{\ast}\nbige_2$
denote the filtered bundle
$\nbigp_{\ast}(\nbige_1\otimes\nbige_2)$
over $\nbige_1\otimes\nbige_2$
defined as
\[
 \nbigp_a(\nbige_1\otimes\nbige_2)
 =\sum_{b+c\leq a}
  \nbigp_b(\nbige_1)\otimes\nbigp_c(\nbige_2).
\]
\index{tensor product (filtered bundle)}
Let $\nhom(\nbige_1,\nbige_2)$
denote the sheaf of $\nbigo_Y(\ast 0)$-morphisms
from $\nbige_1$ to $\nbige_2$.
It is naturally a locally free
$\nbigo_Y(\ast 0)$-module.
Let $\nhom(\nbigp_{\ast}\nbige_1,\nbigp_{\ast}\nbige_2)$
denote the filtered bundle
$\nbigp_{\ast}\nhom(\nbige_1,\nbige_2)$
defined as follows
for any open subset $U\subset Y$:
\[
 \nbigp_a\nhom(\nbige_1,\nbige_2)(U)
 =\bigl\{
 f\in \nhom(\nbige_1,\nbige_2)(U)\,\big|\,
 f(\nbigp_b\nbige_{1|U})
 \subset\nbigp_{b+a}(\nbige_{2|U})\,\,
 (\forall b\in\real)
 \bigr\}.
\]
\index{inner homomorphism (filtered bundle)}
\begin{rem}
If $\nbige_i$ and $\nbigp_{\ast}\nbige_i$ are $G$-equivariant,
the induced filtered bundles are also $G$-equivariant
by the construction.
\hfill\qed
\end{rem}

Note that
$\nbigo_Y(\ast 0)$ is equipped with the canonical
filtered bundle
$\nbigp_{\ast}\bigl(\nbigo_Y(\ast 0)\bigr)$,
i.e.,
$\nbigp_a(\nbigo_Y(\ast 0))
=\nbigo_Y([a]0)$,
where $[a]:=\max\{n\in\seisuu\,|\,n\leq a\}$.
For a filtered bundle $\nbigp_{\ast}\nbige$,
let $\nbigp_{\ast}\nbige^{\lor}$ denote the 
filtered bundle obtained as
$\nhom(\nbigp_{\ast}\nbige,\nbigp_{\ast}\nbigo_Y(\ast 0))$.
\index{dual (filtered bundle)}
Note that $\nbigp_{\ast}\nbigo_Y(\ast 0)$
is naturally $G$-equivariant.
Hence, if $\nbigp_{\ast}\nbige$ is $G$-equivariant,
then $\nbigp_{\ast}\nbige^{\lor}$ is also $G$-equivariant.

\subsubsection{Pull back}

Let $q\in\seisuu_{\geq 1}$.
Let $\varphi_q:\cnum\lrarr \cnum$ be the morphism defined by
$\varphi_q(\zeta)=\zeta^q$.
Let $Y_q:=\varphi_q^{-1}(Y)$.

Let $\nbigp_{\ast}\nbige$ be a filtered bundle
over a locally free $\nbigo_Y(\ast 0)$-module $\nbige$
on $(Y,0)$.
We obtain the $\nbigo_{Y_q}(\ast 0)$-module
$\varphi_q^{\ast}(\nbige)$
and the filtered bundle
$\nbigp_{\ast}(\varphi_q^{\ast}\nbige)$ defined as follows:
\[
 \nbigp_a(\varphi_{q}^{\ast}\nbige)=
 \sum_{\substack{b\in\real,n\in\seisuu\\
 qb+n\leq a}}
 \zeta^{-n}
 \varphi_q^{\ast}(\nbigp_b\nbige)
\]
The filtered bundle $\nbigp_{\ast}(\varphi_q^{\ast}\nbige)$ is denoted by
$\varphi_q^{\ast}(\nbigp_{\ast}\nbige)$.
\index{pull back (filtered bundle)}
We can check the following lemma by using Lemma \ref{lem;21.8.24.1}.
\begin{lem}
Let $\nbige$ and $\nbigp_{\ast}\nbige$ be as above.
Let $\nbige'$ be a saturated locally free $\nbigo_Y(\ast 0)$-submodule
of $\nbige$.
We set $\nbige''=\nbige/\nbige'$. 
We have the induced filtered bundles $\nbigp_{\ast}\nbige'$ and
$\nbigp_{\ast}\nbige''$
over $\nbige'$ and $\nbige''$, respectively.
Then, the pull back of the filtered bundle
$\varphi_{q}^{\ast}(\nbigp_{\ast}\nbige')$
(resp. $\varphi_{q}^{\ast}(\nbigp_{\ast}\nbige'')$)
is equal to the filtered bundle over $\varphi_q^{\ast}(\nbige')$
(resp. $\varphi_q^{\ast}(\nbige'')$)
induced by $\varphi_q^{\ast}(\nbigp_{\ast}\nbige)$
and the inclusion
$\varphi_q^{\ast}(\nbige')\subset\varphi_q^{\ast}(\nbige)$
(resp. the projection
$\varphi_q^{\ast}(\nbige)\lrarr\varphi_q^{\ast}(\nbige'')$).
\hfill\qed
\end{lem}

\begin{lem}
\label{lem;21.8.25.21}
 Let $\nbigp_{\ast}\nbige_i$ $(i=1,2)$
be filtered bundles on $(Y,0)$.
Then, there exist natural isomorphisms
 $\varphi_q^{\ast}(\nbigp_{\ast}\nbige_1\oplus\nbigp_{\ast}\nbige_2)
 \simeq
 \varphi_q^{\ast}(\nbigp_{\ast}\nbige_1)
 \oplus
 \varphi_q^{\ast}(\nbigp_{\ast}\nbige_2)$,
 $\varphi_q^{\ast}(\nbigp_{\ast}\nbige_1\otimes\nbigp_{\ast}\nbige_2)
 \simeq
 \varphi_q^{\ast}(\nbigp_{\ast}\nbige_1)
 \otimes
 \varphi_q^{\ast}(\nbigp_{\ast}\nbige_2)$,
 and
 $\varphi_q^{\ast}\nhom\bigl(\nbigp_{\ast}\nbige_1,
 \nbigp_{\ast}\nbige_2\bigr)
 \simeq
 \nhom\bigl(\varphi_q^{\ast}(\nbigp_{\ast}\nbige_1),
 \varphi_q^{\ast}(\nbigp_{\ast}\nbige_2)\bigr)$.
\end{lem}
\pf
By using compatible frames,
we can observe that
it is enough to check the claims in the case $\rank(\nbige_i)=1$.
In that case, we can check the claims by direct computations.
(For example, see \S\ref{subsection;21.8.25.2}.)
\hfill\qed

\subsubsection{Push-forward}

Let $\nbigp_{\ast}\nbige_1$
be a filtered bundle over
a locally free $\nbigo_{Y_q}(\ast 0)$-module $\nbige_1$.
We obtain the $\nbigo_Y(\ast 0)$-module
$\varphi_{q\ast}(\nbige_1)$
over which we obtain the following filtered bundle:
\[
 \nbigp_a(\varphi_{q\ast}\nbige_1)=
 \varphi_{q\ast}
 \nbigp_{a/q}(\nbige_1).
\]
The filtered bundle
$\nbigp_{\ast}(\varphi_{q\ast}\nbige_1)$ is denoted by
$\varphi_{q\ast}(\nbigp_{\ast}\nbige_1)$.
\index{push-forward (filtered bundle)}
We can check the following lemma by using Lemma \ref{lem;21.8.24.1}.
\begin{lem}
Let $\nbige_1$ be as above.
Let $\nbige_1'$ be 
a saturated locally free $\nbigo_{Y_q}(\ast 0)$-submodule
of $\nbige_1$.
We set $\nbige_1''=\nbige_1/\nbige_1'$.
We obtain the induced filtered bundles
$\nbigp_{\ast}\nbige'_1$
and $\nbigp_{\ast}\nbige''_1$
over $\nbige'_1$ and $\nbige''_1$, respectively.
Then, $\varphi_{q\ast}(\nbigp_{\ast}\nbige'_1)$
(resp. $\varphi_{q\ast}\nbigp_{\ast}\nbige''_1$)
is equal to the filtered bundle over $\varphi_{q\ast}(\nbige'_1)$
(resp. $\varphi_{q\ast}(\nbige''_1)$)
induced by $\varphi_{q\ast}\nbigp_{\ast}\nbige_1$
and the inclusion $\varphi_{q\ast}\nbige_1'\subset\varphi_{q\ast}\nbige_1$
(resp. the projection
$\varphi_{q\ast}\nbige_1\lrarr\varphi_{q\ast}\nbige_1''$).
\hfill\qed
\end{lem}

\vspace{.1in}
Let $\nbigp_{\ast}\nbige$ be a filtered bundle on $(Y,0)$.
Let $\nbigp_{\ast}\nbige_1$ be a filtered bundle on $(Y_q,0)$.
There exist natural isomorphisms
$\varphi_{q\ast}(\nbige_1\otimes\varphi_q^{\ast}\nbige)
\simeq
\varphi_{q\ast}(\nbige_1)\otimes\nbige$
and
$\varphi_{q\ast}\nhom\bigl(
\varphi_q^{\ast}\nbige,\nbige_1
\bigr)
\simeq
\nhom(\nbige,\varphi_{q\ast}\nbige_1)$.

\begin{lem}
There exist natural isomorphisms
 $\varphi_{q\ast}(
 \nbigp_{\ast}\nbige_1\otimes\varphi_q^{\ast}(\nbigp_{\ast}\nbige))
 \simeq
 \varphi_{q\ast}\nbigp_{\ast}\nbige_1\otimes\nbigp_{\ast}\nbige$
and 
$\varphi_{q\ast}\nhom\bigl(
\varphi_q^{\ast}\nbigp_{\ast}\nbige,\nbigp_{\ast}\nbige_1
\bigr)
\simeq
\nhom(\nbigp_{\ast}\nbige,\varphi_{q\ast}\nbigp_{\ast}\nbige_1)$.
\end{lem}
\pf
By using compatible frames,
we can observe that it is enough to study the case
$\rank\nbige=1$.
Then, we can check the claims by direct computations.
\hfill\qed

\subsubsection{Descent}

We set $\Gal_q:=\{\alpha\in\cnum\,|\,\alpha^q=1\}$.
It acts on $Y_q$ by the multiplication of $c$ to $\zeta$
with which we may regard $\Gal_q$ as the Galois group
of the ramified covering $\varphi_q$.
Let $\nbigp_{\ast}\nbige_1$ be a $\Gal_q$-equivariant
filtered bundle on $(Y_q,0)$.
We have the naturally induced $\Gal_q$-action
on $\varphi_{q\ast}\nbige_1$.
Let $\nbige_2$ be the $\Gal_q$-invariant part of $\varphi_{q\ast}\nbige_1$
which is also a locally free $\nbigo_Y(\ast 0)$-module.
Let $\nbigp_{a}\nbige_2$ denote the $\Gal_q$-invariant part of
$\nbigp_a\varphi_{q\ast}\nbige_1$.
Thus, we obtain the filtered bundle
$\nbigp_{\ast}\nbige_2$,
which is called the descent of $\nbigp_{\ast}\nbige_1$.
\index{descent (filtered bundle)}

\begin{lem}
\label{lem;21.8.25.30}
$\varphi_q^{\ast}(\nbigp_{\ast}\nbige_2)$
is naturally isomorphic to $\nbigp_{\ast}\nbige_1$. 
\end{lem}
\pf
By the adjunction,
there exists a natural morphism
$\varphi_q^{\ast}(\varphi_{q\ast}\nbige_1)\lrarr\nbige_1$.
By using Corollary \ref{cor;21.8.25.3},
we can check that the induced morphism
$\varphi_q^{\ast}\nbigp_{\ast}\nbige_2\lrarr\nbigp_{\ast}\nbige_1$
is an isomorphism.
\hfill\qed

\vspace{.1in}
We can check the following lemma
by using Lemma \ref{lem;21.8.24.1}.
\begin{lem}
\label{lem;21.8.25.32}
Let $\nbige_1'$ be a $\Gal_{q}$-equivariant
saturated locally free $\nbigo_Y(\ast 0)$-submodule of $\nbige_1$.
We set $\nbige''_1=\nbige_1/\nbige_1'$.
\begin{itemize}
 \item 
 We obtain the $\nbigo_{Y}(\ast 0)$-submodule $\nbige_2'\subset\nbige_2$
as the descent of $\nbige'_1$.
Then, the decent $\nbigp_{\ast}\nbige'_2$ of $\nbigp_{\ast}\nbige'_1$
is equal to the filtered bundle over $\nbige_2'$
induced by $\nbigp_{\ast}\nbige_2$ with the inclusion
$\nbige'_2\subset\nbige_2$.
 \item
We obtain the $\nbigo_Y(\ast 0)$-quotient module
      $\nbige_2\lrarr\nbige_2''$
as the descent of $\nbige''_1$      .
Then, the decent $\nbigp_{\ast}\nbige''_2$ of $\nbigp_{\ast}\nbige''_1$
is equal to the filtered bundle over $\nbige_2''$
induced by $\nbigp_{\ast}\nbige_2$ with the projection
$\nbige_2\lrarr\nbige_2''$.
\hfill\qed
\end{itemize}
\end{lem}

We can easily check the following lemma
by using a compatible frame.

\begin{lem}
\label{lem;21.8.25.31}
For a filtered bundle $\nbigp_{\ast}\nbige$ on $(Y,0)$,
$\varphi_q^{\ast}(\nbigp_{\ast}\nbige)$ is $\Gal_q$-equivariant,
and the decent is naturally isomorphic to $\nbigp_{\ast}\nbige$.
\hfill\qed
\end{lem}

\subsubsection{Some examples}
\label{subsection;21.8.25.2}

Let $c\in\real$.
Let $\nbigo_Y(\ast 0)\cdot e$ denote
the locally free $\nbigo_Y(\ast 0)$-module of rank one
with a global frame $e$.
Let $\nbigp^{(c)}_{\ast}\bigl(\nbigo_Y(\ast 0)\cdot e\bigr)$
denote the filtered bundle over $\nbigo_Y(\ast 0)\cdot e$
defined as
$\nbigp^{(c)}_{a}(\nbigo_Y(\ast 0)\cdot e)
=\nbigo_{Y}([a-c]0)\cdot e$,
where $[d]:=\max\{n\in\seisuu\,|\,n\leq d\}$ for any $d\in\real$.

\begin{lem}
For any filtered bundle $\nbigp_{\ast}\nbige$ of rank one,
there exist a neighbourhood $Y'$ of $0$ in $Y$
and $-1<c\leq 0$ such that 
 $\nbigp_{\ast}\nbige_{|Y'}
 \simeq\nbigp^{(c)}_{\ast}(\nbigo_{Y'}(\ast 0)\cdot e)$.
\end{lem}
\pf
There exists a neighbourhood $Y'$
and a frame $v$ of $\nbigp_0\nbige_{|Y'}$.
There exists $-1<c\leq 0$
such that $\deg^{\nbigp}(v)=c$.
Then, we obtain
 $\nbigp_{\ast}\nbige_{|Y'}
 \simeq\nbigp^{(c)}_{\ast}(\nbigo_{Y'}(\ast 0)\cdot e)$.
\hfill\qed

\vspace{.1in}

We have the natural isomorphism
$\nbigo_Y(\ast 0)\cdot e_1\otimes\nbigo_Y(\ast 0)\cdot e_2
\simeq \nbigo_Y(\ast 0)\cdot(e_1\otimes e_2)$,
under which we have
\[
 \nbigp^{(c_1)}_{\ast}(\nbigo_Y(\ast0)\cdot e_1)
 \otimes
 \nbigp^{(c_2)}_{\ast}(\nbigo_Y(\ast 0)\cdot e_2)
 \simeq
 \nbigp^{(c_1+c_2)}_{\ast}(\nbigo_Y(\ast 0)\cdot(e_1\otimes e_2)).
\]
We have the natural isomorphism
$\nhom_{\nbigo_Y(\ast 0)}
\bigl(
\nbigo_Y(\ast 0)e_1,\nbigo_Y(\ast 0)e_2
\bigr)\simeq \nbigo_Y(\ast 0)\,e_1^{\lor}\otimes e_2$,
where $(e_1^{\lor}\otimes e_2)$ is determined by
$(e_1^{\lor}\otimes e_2)(e_1)=e_2$.
Under the isomorphism, we have
\[
 \nhom\bigl(
 \nbigp^{(c_1)}_{\ast}(\nbigo_Y(\ast0)\cdot e_1),
 \nbigp^{(c_2)}_{\ast}(\nbigo_Y(\ast 0)\cdot e_2)
 \bigr)
 \simeq
 \nbigp^{(-c_1+c_2)}_{\ast}(\nbigo_Y(\ast 0)\cdot (e_1^{\lor}\otimes e_2)).
\]
Under the natural isomorphism
$\varphi_q^{\ast}(\nbigo_Y(\ast 0)\cdot e)
\simeq
\nbigo_{Y_q}(\ast 0)\cdot \varphi_q^{\ast}(e)$,
we have
\[
 \varphi_q^{\ast}\bigl(\nbigp^{(c)}_{\ast}(\nbigo_{Y}(\ast 0)
  \cdot e)\bigr)
 \simeq
 \nbigp^{(qc)}_{\ast}\bigl(
 \nbigo_{Y_q}(\ast 0)\cdot \varphi_q^{\ast}(e)\bigr).
\]
We have the natural isomorphism
\[
 \varphi_{q\ast}(\nbigo_{Y_q}(\ast 0)\cdot e)
 \simeq
 \bigoplus_{i=0}^{q-1}
 \nbigo_{Y}(\ast 0)\cdot v_i
\]
given by $\zeta^i\cdot e\longmapsto v_i$ $(i=0,\ldots,q-1)$.
Under the isomorphism,
we obtain
\[
 \varphi_{q\ast}(\nbigp^{(c)}_{\ast}\nbigo_{Y_q}(\ast 0)\cdot e)
 \simeq
 \bigoplus_{i=0}^{q-1}
 \nbigp^{(c-i)/q}_{\ast}\bigl(
  \nbigo_Y(\ast 0)\cdot v_i
  \bigr).
\]
Let $m$ be an integer.
Let us consider the $\Gal_q$-action on
$\nbigp_{\ast}(\nbigo_{Y_q}(\ast 0)\,e)$
defined by 
$\alpha^{\ast}(f(\zeta)e)=\alpha^m\cdot f(\alpha\zeta)e$.
The naturally induced $\Gal_q$-action on
$\varphi_{q\ast}(\nbigo_{Y_q}(\ast 0)\cdot e)$
is described as
$\alpha^{\ast}(v_i)=\alpha^{m+i}v_i$.
Let $i_0$ be the unique integer satisfying
$0\leq i\leq q-1$ and $i_0+m\in q\seisuu$.
Then, the descent of
$\nbigo_{Y_q}(\ast 0)\cdot e$ is
$\nbigo_{Y}(\ast 0)\cdot v_{i_0}$,
and
the descent of $\nbigp^{(c)}_{\ast}\cdot e$ is
$\nbigp^{(c-i_0)/q}\nbigo_Y(\ast 0)\cdot v_{i_0}$.

\subsection{Acceptable bundles on a punctured disc}
\label{subsection;20.8.8.31}

Let $(E,\delbar_E)$ be a holomorphic vector bundle
on $Y\setminus\{0\}$
with a Hermitian metric $h$.
For any $a\in\real$ and open subset $U\subset Y$
with $0\in U$,
let $\nbigp^h_aE(U)$ denote the space of
holomorphic sections
$s$ of $E$ on $U\setminus\{0\}$
such that
$|s|_h=O\bigl(|z|^{-a-\epsilon}\bigr)$ for any $\epsilon>0$,
where $|s|_h$ denotes the function on $U\setminus\{0\}$
obtained as the norm of $s$ with respect to $h$.
\index{norm \mbox{$|s|_h$}}
For any open subset $U\subset Y$ with $0\not\in U$,
let $\nbigp^h_aE(U)$ denote the space of
holomorphic sections of $s$ of $E$ on $U$.
Thus, we obtain an $\nbigo_U$-module
$\nbigp^h_aE$.
Similarly,
for any open subset $U\subset Y$ with $0\in U$,
let $\nbigp^h E(U)$ denote the space of
holomorphic sections $s$ of $E$ on $U\setminus\{0\}$
such that
$|s|_h=O\bigl(|z|^{-N}\bigr)$ for some $N>0$.
For any open subset $U\subset Y$ with $0\not\in U$,
let $\nbigp^h E(U)$ denote the space of
holomorphic sections $s$ of $E$ on $U$.
Thus, we obtain an $\nbigo_U(\ast 0)$-module
$\nbigp^h E$.
\index{sheaf $\nbigp^hE$}
\index{sheaf $\nbigp^h_aE$}

Suppose that $(E,\delbar_E,h)$ is acceptable,
i.e.,
the curvature $F(h)$ of the Chern connection
satisfies the following estimate with respect to $h$:
\index{acceptable bundle}
\begin{equation}
\label{eq;20.7.30.20}
F(h)=O\bigl(|z|^{-2}(-\log|z|)^{-2}\bigr)\,dz\,d\zbar.
\end{equation}

\begin{prop}[\cite{cg,Simpson90}]
\label{prop;20.7.26.1}
Under the assumption of the acceptability,
$\nbigp^h_aE$ are locally free $\nbigo_U$-modules,
and $\nbigp^h E$ is a locally free $\nbigo_U(\ast 0)$-module.
\hfill\qed
\end{prop}

In this way, we obtain a filtered bundle 
$\nbigp^h_{\ast}E=
\bigl(
 \nbigp^h_aE\,\big|\,a\in\real
 \bigr)$
 over $\nbigp^h E$.

\begin{rem}
 We shall often use the notation $\nbigp_{\ast}E$
 instead of $\nbigp^h_{\ast}E$
 to simplify the description
 if there is no risk of confusion.
\hfill\qed
\end{rem}

Suppose that $\vecv$ is a frame of $\nbigp^h_a(E)$
compatible with the parabolic structure.
Let $h_0$ be the Hermitian metric
determined as follows:
\[
 h_0(v_i,v_j):=
\left\{
 \begin{array}{ll}
 |z|^{-2\deg^{\nbigp^h}(v_i)} & (i=j),\\
 0   & (i\neq j).
 \end{array}
\right.
\]
Recall the following lemma.
\begin{lem}[\cite{Mochizuki-wild}]
\label{lem;17.10.7.100}
If $(E,\delbar_E,h)$ is acceptable,
there exist $C\geq 1$ and $N>0$
such that
$C^{-1}(-\log |z|)^{-N} h_0
\leq
 h
\leq
 C(-\log|z|)^N h_0$.
\hfill\qed
\end{lem}

\subsubsection{Basic functoriality}

Let $(E_i,\delbar_{E_i})$ $(i=1,2)$ be holomorphic vector bundle
on $Y\setminus \{0\}$ with a Hermitian metric $h_i$.
Suppose that $h_i$ are acceptable.
We obtain the filtered bundles
$\nbigp^{h_i}_{\ast}E_i$.
The bundles
$E_1\oplus E_2$, $E_1\otimes E_2$
and $\nhom(E_1,E_2)$ are naturally equipped with
the holomorphic structure and the metric $\htilde$,
and they are acceptable.
According to \cite{Simpson90},
the following holds.
(It also follows from Lemma \ref{lem;17.10.7.100}.)
\begin{lem}
 \label{lem;17.10.7.200} 
There exist natural isomorphisms
$\nbigp^{\htilde}_{\ast}(E_1\oplus E_2)
\simeq
 \nbigp^{h_1}_{\ast}E_1\oplus\nbigp^{h_2}_{\ast}E_2$,
$\nbigp^{\htilde}_{\ast}(E_1\otimes E_2)
\simeq
 \nbigp^{h_1}_{\ast}E_1\otimes\nbigp^{h_2}_{\ast}E_2$,
and 
$\nbigp^{\htilde}_{\ast}\nhom(E_1,E_2)
\simeq
 \nhom(\nbigp^{h_1}_{\ast}E_1,\nbigp^{h_2}_{\ast}E_2)$.
\hfill\qed
\end{lem}

Let $(E^{\lor},\delbar_{E^{\lor}})$
denote the dual bundle of $(E,\delbar_E)$.
It is equipped with the induced Hermitian metric $h^{\lor}$.
It is easy to see that
$(E^{\lor},\delbar_{E^{\lor}},h^{\lor})$
is acceptable.
As a special case of Lemma \ref{lem;17.10.7.200},
we obtain
$\nbigp_{\ast}^{h^{\lor}}(E^{\lor})
=(\nbigp^h_{\ast}E)^{\lor}$.

\subsubsection{Pull back and descent}

Let $(E,\delbar_E)$ be a holomorphic bundle
on $Y\setminus\{0\}$ with a Hermitian metric $h$
such that $(E,\delbar_E,h)$ is acceptable on $(Y_q,0)$.
We obtain the holomorphic vector bundle
$\varphi_q^{-1}(E,\delbar_E)$
with the induced metric $\varphi_q^{-1}(h)$.
Obviously,
$\varphi_q^{-1}(E,\delbar_E,h)$ is acceptable on $(Y_q,0)$.
We can easily check the following by using Lemma \ref{lem;17.10.7.100}:
\[
 \nbigp^{\varphi_q^{-1}(h)}_{\ast}\bigl(
 \varphi_q^{-1}(E)
 \bigr)
 \simeq
 \varphi_q^{\ast}(\nbigp^h_{\ast}E).
\]
Let $(E_1,\delbar_{E_1})$ be a holomorphic vector bundle
on $Y_q\setminus\{0\}$
with a Hermitian metric $h_1$
such that
$(E_1,\delbar_{E_1},h_1)$ is acceptable.
We naturally obtain a holomorphic vector bundle
$\varphi_{q\ast}(E_1,\delbar_{E_1})$
with the induced metric $h_1$.
It is easy to see that
$\varphi_{q\ast}(E_1,\delbar_{E_1},h_1)$
is acceptable,
and
\[
\nbigp_{\ast}^{\varphi_{q\ast}h_1}\bigl(
 \varphi_{q}(E_1)\bigr)
 \simeq
 \varphi_{q\ast}(\nbigp^{h_1}_{\ast}E_1).
\]
Suppose that
$(E_1,\delbar_{E_1},h_1)$ is $\Gal_q$-equivariant.
We obtain
a holomorphic vector bundle
$(E_2,\delbar_{E_2})$ on $Y_q\setminus \{0\}$
with the induced metric $h_2$
as the descent of $(E_1,\delbar_{E_1},h_1)$.
It is easy to see that
$\nbigp^{h_2}_{\ast}(E_2)$
is the descent of $\nbigp^{h_1}_{\ast}E_1$.

\subsection{Global case}

\subsubsection{Filtered bundles}
\label{subsection;21.8.15.2}

Let $C$ be a Riemann surface with a discrete subset $D\subset Y$.
Let $\nbige$ be a locally free $\nbigo_Y(\ast D)$-module of finite rank.
A filtered bundle
$\nbigp_{\ast}\nbige$ over $\nbige$ is
a tuple
$\bigl(\nbigp_{\veca}\nbige\,\big|\,\veca\in\real^D\bigr)$
of locally free $\nbigo_Y$-submodules
such that the following holds.
\index{filtered bundle (global)}
\begin{itemize}
 \item For $P\in D$, let $(C_P,z_P)$ denote
       a holomorphic coordinate neighbourhood of $P$
       such that $C_P\cap D=\{P\}$ and $z_P(P)=0$.
       Then, for each $\veca=(a(Q)\,|\,Q\in D)\in\real^D$,
       $\nbigp_{\veca}(\nbige)_{|C_P}$ depends only on
       $a(P)$, denoted by
       $\nbigp_{a(P)}(\nbige_{|C_P})$,
       and $\nbigp_{\ast}(\nbige_{|C_P})$
       is a filtered bundle over
       the $\nbigo_{C_P}(\ast P)$-module $\nbige_{|C_P}$.
\end{itemize}
For any locally free $\nbigo_C(\ast D)$-submodule $\nbige'\subset\nbige$
such that $\nbige''=\nbige/\nbige'$ is also locally free,
we obtain the induced filtered bundles
$\nbigp_{\ast}(\nbige')$ and $\nbigp_{\ast}(\nbige'')$
over $\nbige'$ and $\nbige''$, respectively,
applying the constructions in \S\ref{subsection;17.10.5.301}.
\index{induced filtered bundle (subbundle)}
\index{induced filtered bundle (quotient-bundle)}

For filtered bundles
$\nbigp_{\ast}(\nbige_i)$ $(i=1,2)$
over locally free $\nbigo_Y(\ast D)$-modules $\nbige_i$,
we obtain the induced filtered bundles
$\nbigp_{\ast}(\nbige_1\oplus\nbige_2)$
over $\nbige_1\oplus\nbige_2$,
$\nbigp_{\ast}(\nbige_1\otimes\nbige_2)$
over $\nbige_1\otimes\nbige_2$,
and
$\nhom(\nbigp_{\ast}\nbige_1,\nbigp_{\ast}\nbige_2)$
over $\nhom(\nbige_1,\nbige_2)$
as in \S\ref{subsection;21.8.15.1}.
\index{direct sum (filtered bundle)}
\index{tensor product (filtered bundle)}
\index{inner homomorphism (filtered bundle)}
\index{dual (filtered bundle)}

If $C$ is a compact Riemann surface,
then $D$ is finite,
and the degree of a filtered bundle $\nbigp_{\ast}\nbige$
on $(C,D)$ is defined as 
\index{degree $\deg(\nbigp_{\ast}\nbige)$}
\begin{equation}
\label{eq;21.9.6.3}
\deg(\nbigp_{\ast}\nbige)
=\deg(\nbigp_{\veca}\nbige)
-\sum_{P\in D}
 \sum_{a(P)-1<b\leq a(P)}
 b\dim\Gr^{\nbigp}_b(\nbige_{|C_P}).
\end{equation}
The number is independent of the choice of $\veca$.

\subsubsection{Acceptable bundles}

Let $(E,\delbar_E)$ be a holomorphic vector bundle on 
$C\setminus D$ with a Hermitian metric $h$.
For each $P\in D$,
we take a holomorphic coordinate neighbourhood
$(C_P,z_P)$ around $P$ such that $C_P\cap D=\{P\}$
and $z_P(P)=0$.
We set $C_P^{\ast}=C_P\setminus\{P\}$.
By applying the procedure \ref{subsection;20.8.8.31}
to the restriction of $(E,\delbar_E,h)_{|C_P\setminus P}$,
we obtain a tuple of sheaves
$\nbigp^h_{\ast}(E_{|C_P\setminus P})
=\bigl(\nbigp^h_a(E_{|C_P^{\ast}})\,\big|\,a\in\real\bigr)$.
For $\veca=(a(P)\,|\,P\in D)\in\real^D$,
we obtain the $\nbigo_C$-module
$\nbigp^h_{\veca}(E)$
from $(E,\delbar_E)$
and $\nbigp^{h}_{\ast}(E_{|C_P^{\ast}})$ $(P\in D)$.
Similarly, we obtain the $\nbigo_C(\ast D)$-module
$\nbigp^h(E)$.
\index{sheaf $\nbigp^h_{\veca}(E)$}
\index{sheaf $\nbigp^h(E)$}

We say that $(E,\delbar_E,h)$ is acceptable
if the following holds.
\index{acceptable bundle}
\begin{itemize}
 \item Let $P\in D$.
       If $(C_P,z_P)$ is a holomorphic coordinate neighbourhood
       such that $C_P\cap D=\{P\}$ and $z_P(P)=0$,
       then $(E,\delbar_E,h)_{|C_P\setminus P}$
       is acceptable in the sense of \S\ref{subsection;20.8.8.31}.
\end{itemize}
If $(E,\delbar_E,h)$ is acceptable,
then $\nbigp^hE$ is a locally free $\nbigo_C(\ast D)$-module,
and $\nbigp^h_{\ast}(E)$ is a filtered bundle over $\nbigp^hE$.

\begin{prop}[\cite{Simpson88, Simpson90}]
If $C$ is compact,
the following holds.
\[
 \deg(\nbigp_{\ast}^hE)
 =\frac{\sqrt{-1}}{2\pi}
 \int_{C\setminus D}
 \Tr F(h).
\]
\hfill\qed
\end{prop}

\chapter[Formal good parabolic difference modules]{Formal difference modules and good parabolic structure}
\label{section;20.8.8.20}

We recall the classification of formal difference modules
\cite{Chen-Fahim, Duval, Praagman, Turrittin}
in \S\ref{subsection;20.7.20.10}.
We shall introduce the concept of good filtered bundles
in the context of formal difference modules
in \S\ref{subsection;20.7.30.21}.

We explain a geometrization of difference modules
in \S\ref{subsection;20.7.30.22}.
Namely, we introduce a ringed space
$\Hhat_{\infty,q}$,
and explain an equivalence between
formal difference modules
and  some sheaves on $\Hhat_{\infty,q}$.
It is the formal version of the equivalence
in \S\ref{subsection;17.10.28.100}.
Then, in \S\ref{subsection;20.7.30.23},
we define the notion of good filtered bundles 
on $\Hhat_{\infty,q}$
as the translation of the parabolic structure
of difference modules.

In \S\ref{subsection;17.10.28.20},
we explain another equivalence between
some sheaves on $\Hhat_{\infty,q}$
and formal $\lambda$-connections,
which is the ramified and formal version of the construction
in \S\ref{subsection;21.8.12.32}.

\section{Formal difference modules}
\label{subsection;20.7.20.10}

\index{\mbox{field \mbox{$\cnum(\!(y^{-1})\!)$}}}
Let $\cnum(\!(y^{-1})\!)$
denote the field of formal Laurent power series of $y^{-1}$.
We fix an algebraic closure of $\cnum(\!(y^{-1})\!)$.
We also fix a $q$-th root $y_q$ of $y$ 
for each $q\in\seisuu_{\geq 1}$
such that $y_p^{p/q}=y_q$
for any $p\in q\seisuu_{\geq 1}$.
\index{variable $y_q$}

Let $\varrho\in\cnum$.
Let $\Phi^{\ast}$ denote the automorphism of
the field $\cnum(y)$
defined by $\Phi^{\ast}(f)(y)=f(y+\varrho)$.
Note that
$\Phi^{\ast}(y^{-1})=(y+\varrho)^{-1}
=y^{-1}(1+\varrho y^{-1})^{-1}$.
There exists the expansion
$(1+\varrho y^{-1})^{-1}
=\sum_{j=0}^{\infty} (-\varrho y^{-1})^j$
in $\cnum(\!(y^{-1})\!)$.
We obtain
$\Phi^{\ast}\bigl((y^{-1})^k\bigr)
=(y^{-1})^k\bigl(
 \sum_{j=0}^{\infty}(-\varrho y^{-1})^j
 \bigr)^k$
 in $\cnum(\!(y^{-1})\!)$,
and
$\Phi^{\ast}\bigl((y^{-1})^k\bigr)
-(y^{-1})^k
\in (y^{-1})^{k+1}\cnum[\![y^{-1}]\!]$.
Hence, we obtain the induced automorphism
of $\cnum(\!(y^{-1})\!)$,
which is also denoted by $\Phi^{\ast}$.
\index{automorphism $\Phi^{\ast}$}
It induces automorphisms of
$\cnum(\!(y_q^{-1})\!)$
determined by
$\Phi^{\ast}(y_q^{-1})=y_q^{-1}(1+\varrho y^{-1})^{-1/q}$
for any $q\in\seisuu_{>0}$,
where $(1+\varrho y^{-1})^{-1/q}$ denotes the unique power series
$1+\sum_{j>0}a_{q,j}y^{-j}$ such that
$\bigl(
1+\sum_{j>0}a_{q,j}y^{-j}
\bigr)^q=\sum_{j=0}^{\infty}(-\varrho y)^j$.
We regard $y^{-j}=y_q^{-qj}$.
In this section, a difference $\cnum(\!(y_q^{-1})\!)$-module
means a finite dimensional
$\cnum(\!(y_q^{-1})\!)$-module $\nbign$
with a $\cnum$-linear automorphism $\Phi^{\ast}:\nbign\lrarr\nbign$
such that $\Phi^{\ast}(fs)=\Phi^{\ast}(f)\cdot\Phi^{\ast}(s)$.
\index{a difference \mbox{$\cnum(\!(y_q^{-1})\!)$-module}}
A morphism of $\cnum(\!(y_q^{-1})\!)$-difference modules
$g:\nbign_1\lrarr\nbign_2$
is a $\cnum(\!(y_q^{-1})\!)$-linear map such that
$\Phi^{\ast}(g(v))=g(\Phi^{\ast}(v))$
for any $v\in\nbign_1$.

\subsection{Formal difference modules of level $\leq 1$}
\label{subsection;20.7.20.30}

\index{level $\leq 1$ (formal difference module)}
Let $q\in\seisuu_{\geq 1}$.
Let $\nbign$ be a difference $\cnum(\!(y_q^{-1})\!)$-module.
A $\cnum[\![y_q^{-1}]\!]$-lattice $\nbigl$
of $\nbign$ is called $\Phi^{\ast}$-invariant
if $\Phi^{\ast}(\nbigl)=\nbigl$.
\index{\mbox{$\Phi^{\ast}$-invariant $\cnum[\![y_q^{-1}]\!]$-lattice}}
In that case,
we obtain the induced automorphism
$\Phi^{\ast}_{|\infty}$
of $\nbigl_{|\infty}=\nbigl/y_q^{-1}\nbigl$.
\index{vector space $\nbigl_{|\infty}$}

By following \cite[\S1]{Chen-Fahim}
(see also \cite[\S6.2]{Chen-Fahim}),
we say that the level of $\nbign$ is less than $1$
if there exists a $\Phi^{\ast}$-invariant lattice
$\nbigl$ of $\nbign$ such that
the induced automorphism $\Phi^{\ast}_{|\infty}$ of
$\nbigl_{|\infty}$ is unipotent.
\index{automorphism $\Phi^{\ast}_{|\infty}$}

\begin{prop}[\cite{Chen-Fahim, Duval, Praagman, Turrittin}]
\label{prop;20.7.19.40}
For any
difference  $\cnum(\!(y_q^{-1})\!)$-module $\nbign$
of level $\leq 1$,
there exist
$p\in q\seisuu_{>0}$,
a finite subset
$\nbigi(\nbign)\subset
\bigl\{
\gminib\in y_p^{-1}\cnum[y_p^{-1}]\,\big|\,
\deg_{y_p^{-1}}\gminib<p
\bigr\}$
and a decomposition of difference modules
\begin{equation}
\label{eq;20.7.19.1}
 \nbign\otimes_{\cnum(\!(y_q^{-1})\!)}
 \cnum(\!(y_p^{-1})\!)
 =\bigoplus_{\gminib\in \nbigi(\nbign)}
 \nbign_{\gminib},
 \quad (\nbign_{\gminib}\neq 0)
\end{equation}
 such that
each $\nbign_{\gminib}$ has a $\Phi^{\ast}$-invariant lattice
 $\nbigl_{\gminib}$ satisfying
 $
 (\Phi^{\ast}-(1+\gminib)\id_{\nbign_{\gminib}})\nbigl_{\gminib}
 \subset y_p^{-p}\nbigl_{\gminib}$.
 The set $\nbigi(\nbign)$
 and the decomposition {\rm(\ref{eq;20.7.19.1})}
 are uniquely determined.
\end{prop}
\pf
If $\varrho=0$, under the assumption,
the eigenvalues of the $\cnum(\!(y^{-1}_q)\!)$-automorphism
$\Phi^{\ast}$ are invertible elements in
$\cnum[[y_p^{-1}]]$ for some $p\in q\seisuu_{>0}$.
Hence, the claim is easily checked.
If $\varrho\neq 0$,
the claim is known as a part of
the classification of difference modules.
See \cite[Proposition 5]{Chen-Fahim}.
(It is not difficult to prove it directly
by a standard method in the study of formal connections.)
\hfill\qed

\vspace{.1in}

By following \cite{Chen-Fahim},
we say that the level of
a difference module $\nbign$ is $0$
if there exists a $\Phi^{\ast}$-invariant lattice
$\nbigl$ of $\nbign$ such that
$(\Phi^{\ast}-\id)(\nbigl)=y_q^{-q}\nbigl=y^{-1}\nbigl$.
\index{level $0$ (formal difference module)}
In Proposition \ref{prop;20.7.19.40},
each $\nbign_{\gminib}$ is isomorphic to
the tensor product of a difference module of the level $0$
and a difference module
$\cnum(\!(y_p^{-1})\!)\,e$
with the difference operator
$\Phi^{\ast}$ defined as
$\Phi^{\ast}(e)=(1+\gminib)e$.

\begin{rem}
In Proposition {\rm\ref{prop;20.7.19.40}},
$\Phi^{\ast}-(1+\gminib)\id$
is well defined as a $\cnum$-linear endomorphism of
$\nbign_{\gminib}$.
Note that
 $\bigl(\Phi^{\ast}-(1+\gminib)\id\bigr)\nbigl_{\gminib}
 \subset y_p^{-p}\nbigl_{\gminib}$
 holds
 if and only if
 there exists a frame $\vecv=(v_1,\ldots,v_r)$
 of $\nbigl_{\gminib}$
 such that
 $\bigl(\Phi^{\ast}-(1+\gminib)\id\bigr)v_i
 \in y_p^{-p}\nbigl_{\gminib}$.
 Indeed,
 the ``only if'' part of the claim is clear.
 Suppose that there exists such a frame $\vecv$
 of $\nbigl_{\gminib}$.
 For any $f\in \cnum[\![y_p^{-1}]\!]$,
 we obtain
\begin{equation}
\label{eq;20.7.22.1}
 \bigl(\Phi^{\ast}-(1+\gminib)\id\bigr)(fv_i)
 =\Phi^{\ast}(f)
 \bigl(\Phi^{\ast}-(1+\gminib)\id\bigr)(v_i)
 +(1+\gminib)(\Phi^{\ast}(f)-f)v_i.
\end{equation}
Because
 $\Phi^{\ast}(f)-f\in y_p^{-p-1}\cnum[\![y_p^{-1}]\!]$,
 {\rm(\ref{eq;20.7.22.1})} is contained in
 $y_p^{-p}\nbigl_{\gminib}$.
 \hfill\qed
 \end{rem}

The following lemma is standard.
\begin{lem}
\label{lem;21.8.25.40}
Let $g:\nbign_1\lrarr \nbign_2$ be a morphism of
$\cnum(\!(y_q^{-1})\!)$-difference modules of level $\leq 1$.
Choose $p\in q\seisuu_{\geq 0}$ such that
both $\nbign_i\otimes \cnum(\!(y_p^{-1})\!)$ 
have decompositions
of $\cnum(\!(y_p^{-1})\!)$-difference modules
as in {\rm(\ref{eq;20.7.19.1})}
\[
 \nbign_i\otimes\cnum(\!(y_p^{-1})\!)
 =\bigoplus_{\gminib\in y_p^{-1}\cnum[y_p^{-1}]}
 \nbign_{i,\gminib}.
\]
Then, we have
 $g\bigl(\nbign_{1,\gminib}\bigr)
 \subset
 \nbign_{2,\gminib}$,
where we set $\nbign_{i,\gminib}=0$
if $\gminib\not\in \nbigi(\nbign_i)$.
\end{lem}
\pf
We indicate only an outline.
We may assume $p=q$.
It is enough to prove that
$g=0$ if $\nbign_1=\nbign_{i,\gminib_i}$ with $\gminib_1\neq\gminib_2$.
There exist $\Phi^{\ast}$-invariant lattices
$\nbigl_i\subset\nbign_i$
such that 
$(\Phi^{\ast}-(1+\gminib_i)\nbign_{i})\nbigl_i\subset
 y_q^{-q}\nbigl_i$.
Let $\vecv_i$ be a frame of $\nbigl_i$.
Let $A$ be the matrix representing $g$
with respect to the frames $\vecv_i$,
i.e., $g(\vecv_1)=\vecv_2A$.
There exists an expansion
$A=\sum_{j\geq N} A_j (y_q^{-1})^{j}$.
Suppose $A\neq 0$.
We may assume that $A_N\neq 0$.
Let $B_i$ be the matrices determined by
$\Phi^{\ast}(\vecv_i)=\vecv_iB_i$.
We have the expansions
$B_i=(1+\gminib_i)I_{r(i)}+\sum_{j\geq q} B_{i,j}y_q^{-j}$,
where we set $r(i)=\rank\nbign_i$,
and $I_{r(i)}$ denote the $r(i)$-square identity matrices.
We also have the expansions
$\gminib_i=\sum_{j=1}^{q-1} \gminib_{i,j}y_q^{-j}$.
Let $j(0):=\min\{j\,|\,\gminib_{1,j}-\gminib_{2,j}\}$.
Because $B_2\Phi^{\ast}(A)=AB_1$
and $1\leq j_0$,
we obtain $\gminib_{2,j(0)}A_{N}=\gminib_{1,j(0)}A_N$
by taking the coefficients of $y_q^{-N-j_0}$.
It implies $A_N=0$, which contradicts with the assumption.
Hence, we obtain $A=0$.
\hfill\qed

\vspace{.1in}

For $p\in q\seisuu_{>0}$,
we set $\Gal_{q,p}:=\bigl\{\gamma\in\cnum\,|\,\gamma^{p/q}=1\bigr\}$
which acts on $\cnum(\!(y_p^{-1})\!)$
by $\gamma\bullet f(y_p^{-1})=f(\gamma^{-1}y_p^{-1})$.
\index{group $\Gal_{q,p}$}
It is identified with the Galois group of
the extension
$\cnum(\!(y_p^{-1})\!)$ over $\cnum(\!(y_q^{-1})\!)$.
Note that
the $\Gal_{q,p}$-action on $\cnum(\!(y_p^{-1})\!)$
induces a $\Gal_{q,p}$-action on 
$\nbign\otimes_{\cnum(\!(y_q^{-1})\!)}\cnum(\!(y_p^{-1})\!)$.
We obtain the following lemma from
the uniqueness of the index set $\nbigi(\nbign)$
and the decomposition (\ref{eq;20.7.19.1}).
\begin{lem}
\label{lem;20.7.22.4}
The $\Gal_{q,p}$-action on $\cnum(\!(y_p^{-1})\!)$
preserves $\nbigi(\nbign)\subset \cnum(\!(y_p^{-1})\!)$.
Moreover, the following holds.
 \begin{itemize}
 \item For $\gminib\in \nbigi(\nbign)$,
	let $\nbigl_{\gminib}$ be a lattice of $\nbign_{\gminib}$
	satisfying
	$(\Phi^{\ast}-(1+\gminib)\id)\nbigl_{\gminib}
	\subset y_p^{-p}\nbigl_{\gminib}$.
	Then, for any $\gamma\in \Gal_{q,p}$,
	$\gamma\bullet\nbigl_{\gminib}$
	is a lattice of $\nbign_{\gamma\bullet\gminib}$
	satisfying
	$(\Phi^{\ast}-(1+\gamma\bullet\gminib)\id)
	 \nbigl_{\gamma\bullet\gminib}
       \subset y_p^{-p}\nbigl_{\gamma\bullet\gminib}$.
       \hfill\qed
 \end{itemize}
\end{lem}

\begin{rem} 
If $\varrho\neq 0$,
for a difference module $\nbign$ of level $0$,
there exist frame $\vecu=(u_1,\ldots,u_r)$
and a matrix $A\in M_r(\cnum)$
such that
 $\Phi^{\ast}(\vecu)=\vecu (I_r+y_q^{-1}A)$,
 where $I_r\in M_r(\cnum)$ denotes
 the identity matrix.
 (See {\rm\cite{Chen-Fahim}}.)
 \hfill\qed
\end{rem}

\begin{rem}
In {\rm\S\ref{subsection;17.10.28.20}},
we shall explain that
formal difference modules of level $\leq 1$ on
$\cnum(\!(y_q^{-1})\!)$
are equivalent to
formal $\lambda$-connections on $\cnum(\!(w_q^{-1})\!)$
 whose Poincar\'{e} rank is strictly smaller than $q$.
\hfill\qed
\end{rem}

\subsection{Formal difference modules of pure slope}
\label{subsection;20.7.20.31}
\index{pure slope}

Let $\nbign$ be a 
difference $\cnum(\!(y_q^{-1})\!)$-module.
For integers $m$ and $\ell$,
a $\cnum[\![y_q^{-1}]\!]$-lattice of $\nbign$
is called
$y_q^m(\Phi^{\ast})^{\ell}$-invariant
if
$y_q^m(\Phi^{\ast})^{\ell}(\nbigl)=\nbigl$.
\index{$y_q^m(\Phi^{\ast})^{\ell}$-invariant
\mbox{$\cnum[\![y_q^{-1}]\!]$-lattice}}
Note that
the $\cnum$-linear automorphism
of $\nbigl_{|\infty}=\nbigl/y_q^{-1}\nbigl$
is induced by $y_q^m(\Phi^{\ast})^{\ell}$,
which is denoted by
$(y_q^m(\Phi^{\ast})^{\ell})_{|\infty}$.

For any $\omega\in\rnum$,
we set
$Z(q,\omega):=
\bigl\{p\in q\seisuu_{>0}\,\big|\,p\omega\in\seisuu\bigr\}$,
i.e.,
$Z(q,\omega)=
q\seisuu_{>0}\cap \omega^{-1}\seisuu$ if $\omega\neq 0$,
and 
$Z(q,0)=q\seisuu_{>0}$.

\begin{df}
We say that 
a difference $\cnum(\!(y_q^{-1})\!)$-module
$\nbign$ has pure slope $\omega\in\rnum$
if there exist $p\in Z(q,\omega)$
and a $y_q^{p\omega}(\Phi^{\ast})^{p/q}$-invariant lattice
 of $\nbign$.
 \index{slope (formal difference module)}
 \index{pure slope (formal difference module)}
\hfill\qed 
\end{df}

\begin{lem}
\mbox{{}}\label{lem;20.7.19.21}
\begin{itemize}
 \item Suppose that
       a difference $\cnum(\!(y_q^{-1})\!)$-module $\nbign$
       has a pure slope $\omega$.
       Then, for any $p\in Z(q,\omega)$,
       there exists a $y_q^{p\omega}(\Phi^{\ast})^{p/q}$-invariant
       lattice of $\nbign$.
 \item Take any $s\in q\seisuu_{>0}$.
       Then, a difference $\cnum(\!(y_q^{-1})\!)$-module $\nbign$
       has a pure slope $\omega$
       if and only if 
       the induced 
       difference $\cnum(\!(y_s^{-1})\!)$-module
       $\cnum(\!(y_s^{-1})\!)\otimes_{\cnum(\!(y_q^{-1})\!)} \nbign$
       has a pure slope $\omega$.
 \item Take any $s\in Z(q,\omega)$.
       If a difference $\cnum(\!(y_q^{-1})\!)$-module $\nbign$
       has a pure slope $\omega$,
       there exists a $\Gal_{q,s}$-invariant lattice $\nbigl$ of
       the induced 
       difference $\cnum(\!(y_s^{-1})\!)$-module
       $\cnum(\!(y_s^{-1})\!)\otimes_{\cnum(\!(y_q^{-1})\!)} \nbign$
       such that
       $y_s^{s\omega}\Phi^{\ast}(\nbigl)=\nbigl$.
\end{itemize}
\end{lem}
\pf
Let us study the first claim.
By the assumption,
there exists $p_1\in Z(q,\omega)$
such that
there exists a $y_q^{p_1\omega}(\Phi^{\ast})^{p_1/q}$-invariant
lattice $\nbigl_1$ of $\nbign$.
Let $p$ be any element of $Z(q,\omega)$.
If $p\in p_1\seisuu_{>0}$,
it is easy to see that
$\nbigl_1$ is
$y_q^{p\omega}(\Phi^{\ast})^{p/q}$-invariant.
Let us consider the case where
there exists $\ell\in\seisuu_{>0}$
such that $p_1=p\ell$.
We set
$\nbigl=\sum_{i=0}^{\ell-1}
y_q^{ip\omega}\cdot
(\Phi^{\ast})^{ip/q}(\nbigl_1)$.
Then, we obtain
\begin{multline}
 y_q^{p\omega}
 (\Phi^{\ast})^{p/q}(\nbigl_1)
 =\sum_{i=0}^{\ell-1}
 y_q^{(i+1)p\omega}\cdot
 (\Phi^{\ast})^{(i+1)p/q}(\nbigl_1)
  \\
  =\sum_{i=1}^{\ell-1}
  y_q^{ip\omega}
  (\Phi^{\ast})^{ip/q}(\nbigl_1)
  +y_q^{(\ell-1)p\omega}
  \cdot y_q^{p\omega}
  y_q^{-p\omega}\nbigl_1
 =\nbigl.
\end{multline}
The first claim in the general case immediately follows.

Let us study the second claim.
The ``only if'' part of the claim is obvious.
Let us study the ``if'' part of the claim.
It is enough to study the case
$s\in Z(q,\omega)$.
Suppose that
$\cnum(\!(y_s^{-1})\!)\otimes_{\cnum(\!(y_q^{-1})\!)} \nbign$
has a pure slope $\omega$.
Let $\nbigl$ be a $y_s^{s\omega}\Phi^{\ast}$-invariant lattice
of $\cnum(\!(y_{s}^{-1})\!)\otimes_{\cnum(\!(y_q^{-1})\!)}\nbign$.
It is also
a $(y_s^{s\omega})^{s/q}(\Phi^{\ast})^{s/q}$-invariant lattice
of $\cnum(\!(y_{s}^{-1})\!)\otimes_{\cnum(\!(y_q^{-1})\!)}\nbign$.
Note that
$(y_s^{s\omega})^{s/q}=y_q^{s\omega}$.
We set
$\nbigl_1:=\nbigl\cap\nbign$
for the natural inclusion
$\nbign\lrarr
\cnum(\!(y_{s}^{-1})\!)\otimes_{\cnum(\!(y_q^{-1})\!)}\nbign$.
Because
both $\nbigl$ and $\nbign$
are preserved by the action of 
$y_q^{s\omega}(\Phi^{\ast})^{s/q}$
on
$\cnum(\!(y_{s}^{-1})\!)\otimes_{\cnum(\!(y_q^{-1})\!)}\nbign$,
we obtain that
$\nbigl_1$ is 
$y_q^{s\omega}(\Phi^{\ast})^{s/q}$-invariant.
Hence,
$\nbign$ has pure slope $\omega$.

As for the third claim,
we have already proved that
there exists a lattice $\nbigl_1$ of 
$\cnum(\!(y_{s}^{-1})\!)\otimes_{\cnum(\!(y_q^{-1})\!)}\nbign$,
such that
$y_s^{s\omega}\Phi^{\ast}(\nbigl_1)=\nbigl_1$.
Then,
$\nbigl=\sum_{\gamma\in\Gal_{q,s}}
 \gamma\bullet\nbigl_1$
has the desired property.
\hfill\qed

\begin{example}
 A  difference $\cnum(\!(y_q^{-1})\!)$-module
 of level $\leq 1$
 has pure slope $0$.
\hfill\qed
\end{example}

\begin{example}
\index{difference module $\vecL_q(\ell,\alpha)$}
For any $\ell\in\seisuu$ and $\alpha\in\cnum^{\ast}$,
let $\vecL_q(\ell,\alpha)$ denote
the difference $\cnum(\!(y_q^{-1})\!)$-module
$\vecL_q(\ell,\alpha)
=\cnum(\!(y_q^{-1})\!)\,e$
with the difference operator $\Phi^{\ast}$ defined as
\[
 \Phi^{\ast}(e)=y_q^{-\ell}\alpha\cdot e.
\]
Then, $\vecL_q(\ell,\alpha)$
has pure slope $\ell/q$.
\hfill\qed
\end{example}

\begin{lem}
\label{lem;20.7.19.10}
 Let $\nbign_i$ $(i=1,2)$
 be difference $\cnum(\!(y_q^{-1})\!)$-modules
 with pure slopes $\omega_i$.
Any morphism of 
 difference $\cnum(\!(y_q^{-1})\!)$-modules
 $\nbign_1\lrarr\nbign_2$ is $0$
 unless the following condition is satisfied.
\begin{itemize}
 \item $\omega_1=\omega_2$.
 \item For any $p\in Z(q,\omega)$,
       and for any
       $y_q^{p\omega_1}(\Phi^{\ast})^{p/q}$-invariant
       lattices $\nbigl_i$ of $\nbign_i$,
       the induced automorphisms
       on $\nbigl_{i|\infty}$
       have a common eigenvalue.
\end{itemize}
\end{lem}
\pf
Let $f:\nbign_1\lrarr\nbign_2$
be a morphism of 
difference $\cnum(\!(y_q^{-1})\!)$-modules.
Take any $p\in Z(q,\omega_1)\cap Z(q,\omega_2)$.
There exist
$y_q^{p\omega_i}(\Phi^{\ast})^{p/q}$-invariant
lattices $\nbigl_i$ of $\nbign_i$.
Let $\vecv_i$ be frames of $\nbigl_i$.
Let $r_i:=\rank\nbign_i$.
We obtain the matrices
$A_i\in M_{r_i}(\cnum[\![y_q^{-1}]\!])$
such 
$y_q^{p\omega}(\Phi^{\ast})^{p/q}\vecv_i=\vecv_i A_i$.
For the expansion
$A_i=\sum_{j=0}^{\infty} A_{i,j}y_q^{-j}$,
$A_{i,0}$ are invertible.

Let $B$ be the matrix determined by
$f(\vecv_1)=\vecv_2\cdot B$.
If $B\neq 0$,
there exists the expansion
$B=\sum_{j\leq j_0} B_{j}y_q^{j}$
such that $B_{j_0}\neq 0$.
Because $f$ is compatible with
the difference operators,
we obtain
\begin{equation}
\label{eq;20.7.19.20}
 B(y_q^{-1})\cdot A_1\cdot y_q^{-p\omega_1}
 =A_2\cdot B\bigl((1+(p/q)\varrho y_q^{-q})^{-1/q}\bigr)
 \cdot y_q^{-p\omega_2}.
\end{equation}
By comparing the expansions
of the both sides of (\ref{eq;20.7.19.20}),
we obtain $B_{j_0}=0$
unless $\omega_1=\omega_2$ holds
and $A_{i,0}$ have a common eigenvalue.
Then, the claim of the lemma follows.
\hfill\qed

\vspace{.1in}
We can classify
difference modules with pure slope
by the following proposition
and Lemma \ref{lem;20.7.19.21}.

\begin{prop}[\cite{Chen-Fahim, Duval, Praagman, Turrittin}]
\label{prop;20.7.19.22}
 Let $\nbign$ be a difference $\cnum(\!(y_q^{-1})\!)$-module
with pure slope $\omega$.
Suppose that $q\omega\in\seisuu$.
Then, there exists a unique decomposition of
 difference  $\cnum(\!(y_q^{-1})\!)$-modules
\begin{equation}
\label{eq;20.7.19.30}
 \nbign=\bigoplus_{\alpha\in\cnum^{\ast}}
  \nbign_{\alpha}
\end{equation}
such that each  $\nbign_{\alpha}$ is isomorphic to
 the tensor product of
 $\vecL_q(q\omega,\alpha)$
 and a difference $\cnum(\!(y_q^{-1})\!)$-module
 of level $\leq 1$.
\end{prop}
\pf
Let $\nbigl$ be a $y_q^{q\omega}\Phi^{\ast}$-invariant
lattice of $\nbign$.
We obtain the generalized eigen decomposition
\begin{equation}
 \nbigl_{|\infty}
 =\bigoplus_{\alpha\in\cnum^{\ast}}
  \EE_{\alpha}(\nbigl_{|\infty})
\end{equation}
of $(y_{q}^{q\omega}\Phi^{\ast})_{|\infty}$,
i.e.,
the restriction of
$\alpha^{-1}(y_{q}^{q\omega}\Phi^{\ast})_{|\infty}$
to 
$\EE_{\alpha}(\nbigl_{|\infty})$
are unipotent.
It is a standard fact that
there exists a decomposition
 $\nbigl=\bigoplus\nbigl_{\alpha}$
 of $\cnum[\![y_q^{-1}]\!]$-lattice
 such that
 (i) $y_q^{q\omega}\Phi^{\ast}(\nbigl_{\alpha})
 =\nbigl_{\alpha}$,
 (ii) $\nbigl_{\alpha|\infty}=\EE_{\alpha}(\nbigl_{|\infty})$.
 (See \cite[\S3]{Chen-Fahim}, for example.)
Then, the induced decomposition
$\nbign=\bigoplus
\cnum(\!(y_q^{-1})\!)\otimes_{\cnum[\![y_q^{-1}]\!]}
\nbigl_{\alpha}$
has the desired property.
The uniqueness follows from
Lemma \ref{lem;20.7.19.10}.
\hfill\qed

\vspace{.1in}
Let $\nbign_i$ $(i=1,2)$
be difference $\cnum(\!(y_q^{-1})\!)$-modules
with pure slope $\omega$.
Suppose $q\omega\in\seisuu$.
Each $\nbign_i$ has a decomposition
$\nbign_i=\bigoplus_{\alpha\in\cnum^{\ast}}
\nbign_{i,\alpha}$
as in (\ref{eq;20.7.19.30}).
We obtain the following lemma
from Lemma \ref{lem;20.7.19.10}.
\begin{lem}
\label{lem;21.8.26.12}
 Any morphism
 $F:\nbign_1\lrarr\nbign_2$
 preserves the decompositions
 $\nbign_i=\bigoplus_{\alpha\in\cnum^{\ast}}
\nbign_{i,\alpha}$. 
\hfill\qed
\end{lem}

Let $\nbign$ be a difference $\cnum(\!(y_q^{-1})\!)$-difference
of pure slope $\omega\in\rnum$.
Take $p\in q\seisuu_{>0}$
such that $p\omega\in\seisuu$.
Then,
there exists a decomposition
\begin{equation}
\label{eq;20.7.22.2}
 \cnum(\!(y_p^{-1})\!)
 \otimes_{\cnum(\!(y_q^{-1})\!)} \nbign
 =\bigoplus
  \nbign_{\alpha}. 
\end{equation}
The following lemma is easy to see.
\begin{lem}
\label{lem;20.7.22.3}
 For any $\gamma\in\Gal_{q,p}$,
we obtain
$\gamma\bullet\nbign_{\alpha}
=\nbign_{\alpha\gamma^{-p\omega}}$.
 As a result, the following holds.
\begin{itemize}
 \item Let $F$ be the automorphism of
        $\cnum(\!(y_p^{-1})\!)
       \otimes_{\cnum(\!(y_q^{-1})\!)} \nbign$
       obtained by the multiplication of
       $\alpha y_p^{-p\omega}$ on
       $\nbign_{\alpha}$ in {\rm(\ref{eq;20.7.22.2})}.
       Then, $F$ is equivariant with respect to
       the $\Gal_{q,p}$-action
\hfill\qed
\end{itemize}
\end{lem}

\subsection{Slope decomposition of formal difference modules}
\label{subsection;20.7.18.30}

We can understand the structure of
general formal difference modules
by the slope decomposition
in the following theorem.

\begin{thm}[see \cite{Chen-Fahim,Duval, Praagman, Turrittin}]
\label{thm;20.7.19.41}
For any difference $\cnum(\!(y_q^{-1})\!)$-module $\nbign$,
there exists a unique decomposition of
difference $\cnum(\!(y_q^{-1})\!)$-modules
\[
 \nbign=\bigoplus_{\omega\in\rnum}
 \nbigs_{\omega}(\nbign)
\]
such that
each $\nbigs_{\omega}(\nbign)$ has pure slope $\omega$. 
It is called the slope decomposition.
\index{slope decomposition (formal difference module)}
\hfill\qed
\end{thm}

Lemma \ref{lem;20.7.19.10}
implies the following lemma.
\begin{lem}
\label{lem;21.8.26.13}
Any morphism  of difference $\cnum(\!(y_q^{-1})\!)$-modules
$F:\nbign_1\lrarr\nbign_2$
preserves the slope decompositions
$\nbign_i=
\bigoplus
 \nbigs_{\omega}(\nbign_i)$.
 \hfill\qed
\end{lem}

\begin{cor}
Let $\nbign$ be any difference $\nbigc(\!(y_q^{-1})\!)$-module.
For any difference $\nbigc(\!(y_q^{-1})\!)$-submodule
$\nbign_1\subset\nbign$,
we have
$\nbigs_{\omega}(\nbign_1)\subset\nbigs_{\omega}(\nbign)$.
For any difference $\nbigc(\!(y_q^{-1})\!)$-quotient module
$\nbign_2$ of $\nbign$,
$\nbigs_{\omega}(\nbign_2)$
is a quotient of  
$\nbigs_{\omega}(\nbign)$
\hfill\qed
\end{cor}

\section{Good filtered bundles of formal difference modules}
\label{subsection;20.7.30.21}

\subsection{Filtered bundles over $\cnum(\!(y_q^{-1})\!)$-modules}
\label{subsection;20.7.20.20}

Let $\nbige$ be any $\cnum(\!(y_q^{-1})\!)$-module
of finite rank.
A filtered bundle over $\nbige$
is defined to be an increasing sequence
$\nbigp_{\ast}\nbige
=\bigl(
\nbigp_a\nbige\,\big|\,a\in\real
\bigr)$
of $\cnum[\![y_q^{-1}]\!]$-lattices of
$\nbige$
such that
(i) for any $a\in\real$,
there exists $\epsilon>0$
such that $\nbigp_{a}\nbige=\nbigp_{a+\epsilon}\nbige$,
(ii) $\nbigp_{a+n}\nbige=y_q^{n}\nbigp_a\nbige$
for any $a\in\real$ and $n\in\seisuu$.
\index{filtered bundle (formal)}
This is the formal version of the notion of filtered bundle
in \S\ref{subsection;20.8.1.40}.
We set
\begin{equation}
\label{eq;21.9.17.80}
 \Gr^{\nbigp}_a(\nbige):=
  \nbigp_a(\nbige)\big/\nbigp_{<a}(\nbige),
\end{equation}
where $\nbigp_{<a}(\nbige)=\sum_{b<a}\nbigp_b(\nbige)$.
\index{vector space $\Gr^{\nbigp}_a(\nbige)$}
A frame $\vecv$ of $\nbigp_a\nbige$
is called compatible
if there exists a decomposition $\vecv=\coprod_{a-1<b\leq a} \vecv_b$
such that
$\vecv_b$ is a tuple of elements of $\nbigp_b\nbige$
and induces a base of $\Gr^{\nbigp}_b(\nbige)$.
\index{compatible frame}
The parabolic degrees $\deg^{\nbigp}(s)$ for $s\in\nbige$
are also defined
as in (\ref{eq;20.8.8.50}).
\index{degree $\deg^{\nbigp}(s)$}

\subsubsection{$G$-equivariance}

Let $G$ be a finite group acting on
the field $\cnum(\!(y_q^{-1})\!)$
such that
$g^{\ast}(y_q^{-1})\in y_q^{-1}\cnum[\![y_q^{-1}]\!]$
for any $g\in G$.
A $\cnum(\!(y_q^{-1})\!)$-module $\nbige$
is called $G$-equivariant
if $\nbige$ is a $G$-representation over $\cnum$
such that $g^{\ast}(fv)=g^{\ast}(f)\cdot g^{\ast}(v)$
for any $g\in G$.
\index{$G$-equivariant \mbox{$\cnum(\!(y_q^{-1})\!)$-module}}
We say that a filtered bundle $\nbigp_{\ast}\nbige$ over $\nbige$
is $G$-equivariant
if each $\nbigp_a\nbige$ $(a\in\real)$
is preserved by the $G$-action.
\index{$G$-equivariant filtered bundle}
Because $G$ is finite, the following lemma is obvious.
\begin{lem}
If $\nbige$ and $\nbigp_{\ast}\nbige$ are $G$-equivariant,
then there exist
$G$-invariant subspaces $H_b\subset\nbigp_b\nbige$ $(b\in\real)$
such that
the natural morphisms $\nbigp_b\nbige\lrarr \Gr^{\nbigp}_b(\nbige)$ 
induce isomorphisms of the $G$-representations
$H_b\simeq \Gr^{\nbigp}_b(\nbige)$.
\hfill\qed
\end{lem}

As in the case of Corollary \ref{cor;21.8.25.10},
we obtain the following.
\begin{cor}
\label{cor;21.8.25.11}
Let $\Gr^{\nbigp}_b(\nbige)=\bigoplus_{i=1}^{m(b)} V_{b,i}$
$(a-1<b\leq a)$
be decompositions of $G$-representations.
Then, there exists a $G$-equivariant decomposition
$\nbigp_{\ast}\nbige=\bigoplus_{a-1<b\leq a}\bigoplus_{i=1}^{m(b)}
\nbigp_{\ast}\nbige_{b,i}$
such that the following holds.
\begin{itemize}
 \item We have $\Gr^{\nbigp}_{c}(\nbige_{b,i})=0$
       for $a-1<c\leq a$ unless $c=b$,
       and
       $\Gr^{\nbigp}_b(\nbige_{b,i})=V_{b,i}$.
       \hfill\qed
\end{itemize} 
\end{cor}

\begin{cor}
\label{cor;21.8.25.12}
If $G$ is cyclic,
there exists a compatible frame $\vecv$ of $\nbigp_a\nbige$
such that
$g^{\ast}(v_i)=g^{n(i)}v_i$ for some $n(i)\in\seisuu$.
\index{equivariant compatible frame}
\hfill\qed
\end{cor}

\subsubsection{Submodules, quotient modules and splittings}

Let $\nbige$ be a $G$-equivariant
$\cnum(\!(y_q^{-1})\!)$-module
equipped with a $G$-equivariant filtered bundle
$\nbigp_{\ast}\nbige$.
For any $G$-equivariant $\cnum(\!(y_q^{-1})\!)$-submodule
$\nbige'\subset\nbige$,
we set $\nbigp_a(\nbige')=\nbigp_a(\nbige)\cap\nbige'$
$(a\in\real)$
and we obtain the $G$-equivariant induced filtered bundle over $\nbige'$.
\index{induced filtered bundle (subbundle)}
Put $\nbige''=\nbige/\nbige'$.
Let $\nbigp_a(\nbige'')$ $(a\in\real)$ denote the image of
$\nbigp_a(\nbige)\lrarr\nbige''$,
and we obtain the $G$-equivariant induced filtered bundle over $\nbige''$.
\index{induced filtered bundle (quotient-bundle)}
The following lemma is similar to Lemma \ref{lem;21.8.24.1}.
\begin{lem}
\label{lem;21.8.25.20}
There exists a $G$-equivariant splitting
$\nbige''\lrarr\nbige$
which induces
$\nbigp_{\ast}\nbige\simeq\nbigp_{\ast}\nbige'\oplus\nbigp_{\ast}\nbige''$.
\hfill\qed
 \end{lem}

\subsubsection{Basic functoriality}
Let $\nbige_i$ $(i=1,2)$ be $\cnum(\!(y_q^{-1})\!)$-modules
of finite rank.
Let $\nbigp_{\ast}\nbige_i$ be filtered bundles
over $\nbige_i$.
By setting
$\nbigp_{a}(\nbige_1\oplus\nbige_2)
=\nbigp_a\nbige_1\oplus\nbigp_a\nbige_2$,
we obtain the filtered bundle
$\nbigp_{\ast}(\nbige_1\oplus\nbige_2)$
over $\nbige_1\oplus\nbige_2$,
which is denoted by
$\nbigp_{\ast}(\nbige_1)\oplus\nbigp_{\ast}(\nbige_2)$.
By setting
$\nbigp_a(\nbige_1\otimes\nbige_2)
 =\sum_{b+c\leq a}
 \nbigp_b(\nbige_1)\otimes\nbigp_c(\nbige_2)$,
we obtain
the filtered bundle
$\nbigp_{\ast}(\nbige_1\otimes\nbige_2)$
over $\nbige_1\otimes\nbige_2$,
which is denoted by
$\nbigp_{\ast}(\nbige_1)
\otimes
\nbigp_{\ast}(\nbige_2)$.
\index{direct sum (filtered bundle)}
\index{tensor product (filtered bundle)}
By setting
\[
 \nbigp_{a}\Hom(\nbige_1,\nbige_2)
 =\bigl\{
  f\in \Hom(\nbige_1,\nbige_2)\,\big|\,
  f(\nbigp_b\nbige_1)
  \subset \nbigp_{b+a}\nbige_2\,\,(\forall b\in\real)
  \bigr\},
\]
we obtain the filtered bundle
$\nbigp_{\ast}\Hom(\nbige_1,\nbige_2)$
over $\Hom(\nbige_1,\nbige_2)$,
which is denoted by
$\Hom(\nbigp_{\ast}\nbige_1,\nbigp_{\ast}\nbige_2)$.
\index{inner homomorphism (filtered bundle)}

Note that $\cnum(\!(y_q^{-1})\!)$
is equipped with the canonical filtered bundle
$\nbigp_{\ast}\cnum(\!(y_q^{-1})\!)$
defined by
$\nbigp_a(\cnum(\!(y_q^{-1})\!))
=y_q^{[a]}\cnum[[y_q^{-1}]]$,
where $[a]:=\max\{n\in\seisuu\,|\,n\leq a\}$.
For a filtered bundle $\nbigp_{\ast}\nbige$,
the dual $(\nbigp_{\ast}\nbige)^{\lor}$
is defined as
$\Hom(\nbigp_{\ast}\nbige,\nbigp_{\ast}(\cnum(\!(y_q^{-1})\!)))$.
\index{dual (filtered bundle)}

\begin{rem}
If $\nbigp_{\ast}\nbige_i$ are $G$-equivariant,
the induced filtered bundles are also $G$-equivariant.
Because $\nbigp_{\ast}\cnum(\!(y_q^{-1})\!)$ is naturally
$G$-equivariant,
$\nbigp_{\ast}\nbige^{\lor}$ is $G$-equivariant
if $\nbigp_{\ast}\nbige$  is $G$-equivariant.
\hfill\qed
\end{rem}

\subsubsection{Pull back}
Let $p\in q\seisuu_{>0}$.
For any  $\cnum(\!(y_q^{-1})\!)$-module $\nbige$
of finite rank,
we set
$\varphi_{q,p}^{\ast}(\nbige):=
\cnum(\!(y_p^{-1})\!)
\otimes_{\cnum(\!(y_q^{-1})\!)}\nbige$.
(See also \S\ref{subsection;20.7.21.1} below
for this notation.)
\index{pull back}
Let $\nbigp_{\ast}\nbige$
be a filtered bundle over
a $\cnum(\!(y_q^{-1})\!)$-module $\nbige$
of finite rank.
For any $a\in\real$,
we define
\[
 \nbigp_a(\varphi_{q,p}^{\ast}\nbige)
 =\sum_{n+b(p/q)\leq a}
 y_p^{n}
 \nbigp_{b}(\nbige)
 \otimes_{\cnum[\![y_q^{-1}]\!]}\cnum[\![y_p^{-1}]\!].
\]
Thus, we obtain a filtered bundle
$\nbigp_{\ast}\varphi_{q,p}^{\ast}\nbige$
over $\varphi_{q,p}^{\ast}\nbige$,
which we denote by
$\varphi_{q,p}^{\ast}(\nbigp_{\ast}\nbige)$,
and called the pull back of $\nbigp_{\ast}\nbige$.
\index{pull back (filtered bundle)}
The filtered bundle
$\nbigp_{\ast}\nbige$ is naturally
$\Gal_{q,p}$-equivariant.
(See \S\ref{subsection;20.7.20.30}
for $\Gal_{q,p}$.)
By the induced $\Gal_{q,p}$-action
on $\Gr^{\nbigp}_a(\nbige)$,
we obtain the canonical decomposition
\begin{equation}
\label{eq;20.7.20.101}
 \Gr^{\nbigp}_a(\varphi_{q,p}^{\ast}\nbige)
 =\bigoplus_{m=0}^{(p/q)-1}
  \nbigg_m\Gr^{\nbigp}_a(\varphi_{q,p}^{\ast}\nbige),
\end{equation}
where $\alpha\in\Gal_{q,p}$ acts
on $\nbigg_m\Gr^{\nbigp}_a(\nbige)$
by the multiplication of $\alpha^m$.
There exists the natural isomorphism
\begin{equation}
\label{eq;20.7.20.102}
\nbigg_0\Gr^{\nbigp}_a(\varphi^{\ast}_{q,p}\nbige)
\simeq
\Gr^{\nbigp}_{aq/p}(\nbige).
\end{equation}

We can check the following lemma by using
Lemma \ref{lem;21.8.25.20}.
\begin{lem}
\label{lem;21.8.25.42}
Let $\nbige$ and $\nbigp_{\ast}\nbige$ be as above.
Let $\nbige'\subset\nbige$ be a $\cnum(\!(y_q^{-1})\!)$-submodule.
We set $\nbige''=\nbige/\nbige''$.
We have the induced filtered bundles
$\nbigp_{\ast}\nbige'$ and $\nbigp_{\ast}\nbige''$
over $\nbige'$ and $\nbige''$, respectively. 
 Then, $\varphi_{q,p}^{\ast}(\nbigp_{\ast}\nbige')$
 (resp. $\varphi_{q,p}^{\ast}(\nbigp_{\ast}\nbige'')$)
is equal to the filtered bundle
 over $\varphi_{q,p}^{\ast}\nbige'$
 (resp. $\varphi_{q,p}^{\ast}\nbige''$)
induced by $\varphi_{q,p}^{\ast}(\nbigp_{\ast}\nbige)$ 
and the inclusion
 $\varphi_{q,p}^{\ast}(\nbige')\subset\varphi_{q,p}^{\ast}(\nbige)$
 (resp. the projection
  $\varphi_{q,p}^{\ast}(\nbige)\lrarr\varphi_{q,p}^{\ast}(\nbige'')$). 
 \hfill\qed
\end{lem}

We can check the following lemma by using compatible frames
as in the case of Lemma \ref{lem;21.8.25.21}.
\begin{lem}
\label{lem;21.8.26.1}
 Let $\nbige_i$ be $\cnum(\!(y_q^{-1})\!)$-modules of finite rank
 equipped with a filtered bundle $\nbigp_{\ast}\nbige_i$.
Then, there exist natural isomorphisms
 $\varphi_{q,p}^{\ast}(\nbigp_{\ast}\nbige_1\oplus\nbigp_{\ast}\nbige_2)
 \simeq
 \varphi_{q,p}^{\ast}(\nbigp_{\ast}\nbige_1)
 \oplus
 \varphi_{q,p}^{\ast}(\nbigp_{\ast}\nbige_2)$,
 $\varphi_{q,p}^{\ast}(\nbigp_{\ast}\nbige_1\otimes\nbigp_{\ast}\nbige_2)
 \simeq
 \varphi_{q,p}^{\ast}(\nbigp_{\ast}\nbige_1)
 \otimes
 \varphi_{q,p}^{\ast}(\nbigp_{\ast}\nbige_2)$,
 and
 $\varphi_{q,p}^{\ast}\Hom\bigl(\nbigp_{\ast}\nbige_1,
 \nbigp_{\ast}\nbige_2\bigr)
\simeq
 \Hom\bigl(\varphi_{q,p}^{\ast}(\nbigp_{\ast}\nbige_1),
\varphi_{q,p}^{\ast}(\nbigp_{\ast}\nbige_2)\bigr)$.
\hfill\qed
\end{lem}

\subsubsection{Push-forward}

Let $p\in q\seisuu_{\geq 1}$.
Any $\cnum(\!(y_p^{-1})\!)$-module $\nbige_1$
of finite rank
is naturally a $\cnum(\!(y_q^{-1})\!)$-module
of finite rank,
which is denoted by
$\varphi_{q,p\ast}(\nbige_1)$.
\index{push-forward}
Let $\nbigp_{\ast}\nbige_1$ be a filtered bundle
over a $\cnum(\!(y_p^{-1})\!)$-module $\nbige_1$
of finite rank.
For any $a\in\real$,
we obtain the following
$\cnum[\![y_q^{-1}]\!]$-submodule
\[
 \nbigp_a\bigl(\varphi_{q,p\ast}(\nbige_1)\bigr)
 =\nbigp_{aq/p}\nbige_1
 \subset
 \varphi_{q,p\ast}(\nbige_1).
\]
Thus, we obtain a filtered bundle
$\nbigp_{\ast}(\varphi_{q,p\ast}\nbige_1)$
over $\varphi_{q,p\ast}\nbige_1$,
which is denoted by 
$\varphi_{q,p\ast}(\nbigp_{\ast}\nbige_1)$,
and called the push-forward of
$\nbigp_{\ast}\nbige_1$.
\index{push-forward (filtered bundle)}

\begin{lem}
\label{lem;21.8.26.2}
Let $\nbige_1$ and $\nbigp_{\ast}\nbige_1$ be as above.
Let $\nbige'_1\subset\nbige_1$ be a $\cnum(\!(y_q^{-1})\!)$-submodule.
We set $\nbige''_1:=\nbige_1/\nbige'_1$.
They are equipped with the induced filtered bundles
$\nbigp_{\ast}\nbige'_1$ and $\nbigp_{\ast}\nbige''_1$,
respectively. 
Then, $\varphi_{q,p\ast}(\nbigp_{\ast}\nbige'_1)$
(resp. $\varphi_{q,p\ast}(\nbigp_{\ast}\nbige''_1)$)
is equal to the filtered bundle induced by 
$\varphi_{q,p\ast}(\nbigp_{\ast}\nbige_1)$
and the inclusion 
$\varphi_{q,p\ast}(\nbige'_1)\subset\varphi_{q,p\ast}(\nbige_1)$
(resp. the projection
$\varphi_{q,p\ast}(\nbige_1)\lrarr\varphi_{q,p\ast}(\nbige''_1)$).
\hfill\qed
\end{lem}

Let $\nbige$ be a $\cnum(\!(y_q^{-1})\!)$-module of finite rank
equipped with a filtered bundle $\nbigp_{\ast}\nbige$.
Let $\nbige_1$ and $\nbigp_{\ast}\nbige_1$ be as above.
We obtain the following lemma by using a compatible frame of
$\nbigp_{\ast}\nbige$.
\begin{lem}
\label{lem;21.8.26.3}
The natural isomorphism
$\varphi_{q,p\ast}(\nbige_1\otimes\varphi_{q,p}^{\ast}\nbige)
\simeq
\varphi_{q,p\ast}(\nbige_1)\otimes\nbige$
induces 
$\varphi_{q,p\ast}(\nbigp_{\ast}\nbige_1
 \otimes\varphi_{q,p}^{\ast}\nbigp_{\ast}\nbige)
\simeq
 \varphi_{q,p\ast}(\nbigp_{\ast}\nbige_1)\otimes\nbigp_{\ast}\nbige$.
The natural isomorphism
 $\varphi_{q,p\ast}(
 \Hom(\varphi_{q,p}^{\ast}\nbige,\nbige_1))
\simeq
  \Hom(\nbige,\varphi_{q,p\ast}\nbige_1)$
induces 
 $\varphi_{q,p\ast}(
 \Hom(\varphi_{q,p}^{\ast}\nbigp_{\ast}\nbige,\nbigp_{\ast}\nbige_1))
\simeq
  \Hom(\nbigp_{\ast}\nbige,\varphi_{q,p\ast}\nbigp_{\ast}\nbige_1)$.
\hfill\qed
\end{lem}

\subsubsection{Descent}

Let $\nbige_1$ be a $\cnum(\!(y_p^{-1})\!)$-module
of finite rank equipped with a filtered bundle
$\nbigp_{\ast}\nbige_1$.
Suppose that $\nbige_1$ and $\nbigp_{\ast}\nbige_1$
are $\Gal_{q,p}$-equivariant.
We obtain the $\cnum(\!(y_q^{-1})\!)$-module of finite rank
$\nbige_2$ as the $\Gal_{q,p}$-invariant part of
$\varphi_{q,p\ast}(\nbige_1)$,
equipped with the induced filtered bundle
$\nbigp_{\ast}\nbige_2$.
It is called the descent of $\nbigp_{\ast}\nbige_1$.
\index{descent (filtered bundle)}
The following lemma is similar to Lemma \ref{lem;21.8.25.30}.
\begin{lem}
\label{lem;21.8.26.4}
 $\varphi_{q,p}^{\ast}\nbigp_{\ast}\nbige_2$ is naturally isomorphic to
$\nbigp_{\ast}\nbige_1$.
\hfill\qed
\end{lem}

The following lemma is similar to Lemma \ref{lem;21.8.25.32}.
\begin{lem}
\label{lem;21.8.25.41}
Let $\nbige_1'$ be a $\Gal_{q,p}$-equivariant
$\cnum(\!(y_p^{-1})\!)$-submodule of $\nbige_1$.
We set $\nbige''_1:=\nbige_1/\nbige'_1$,
which is naturally $\Gal_{q,p}$-equivariant.
We obtain the $\cnum(\!(y_q^{-1})\!)$-submodule
$\nbige_2'\subset\nbige_2$
as the descent of $\nbige'_1$.
We also obtain the $\cnum(\!(y_q^{-1})\!)$-quotient module
$\nbige_2\lrarr\nbige_2''$
as the descent of $\nbige''_1$.
Then, the decent $\nbigp_{\ast}\nbige'_2$ of $\nbigp_{\ast}\nbige'_1$
is equal to the filtered bundle over $\nbige_2'$
induced by $\nbigp_{\ast}\nbige_2$ with the inclusion
$\nbige'_2\subset\nbige_2$.
The descent of $\nbigp_{\ast}\nbige''_2$ of $\nbigp_{\ast}\nbige''_1$
is equal to the filtered bundle over $\nbige_2''$
induced by $\nbigp_{\ast}\nbige_2$ with the projection
$\nbige_2\lrarr\nbige_2''$.
\hfill\qed
\end{lem}

The following lemma is similar to Lemma \ref{lem;21.8.25.31},
and easy to check.
\begin{lem}
\label{lem;21.8.26.5}
For a filtered bundle
$\nbigp_{\ast}\nbige$ over
a $\cnum(\!(y_q^{-1})\!)$-module $\nbige$,
the pull back
$\varphi_{q,p}^{\ast}(\nbigp_{\ast}\nbige)$
is naturally equipped with
the $\Gal_{q,p}$-action,
and the descent of
$\varphi_{q,p}^{\ast}(\nbigp_{\ast}\nbige)$
is naturally isomorphic to
$\nbigp_{\ast}\nbige$.
\hfill\qed
\end{lem}

\subsection{Good filtered bundles over formal difference modules}
\label{subsection;20.8.8.10}

Let $(\nbign,\Phi^{\ast})$
be a difference $\cnum(\!(y_q^{-1})\!)$-module.
We set
$S(q):=\bigl\{
 \sum_{j=1}^{q-1}\gminib_jy_q^{-j}\,\big|\,
 \gminib_j\in\cnum
 \bigr\}$.
\index{set $S(q)$}
\begin{df}
\label{df;20.7.20.2}\index{unramified difference module}
 We say that
 $(\nbign,\Phi^{\ast})$ is unramified
 if there exists a decomposition of
difference $\cnum(\!(y_q^{-1})\!)$-modules
\begin{equation}
\label{eq;20.7.20.1}
 \nbign
=\bigoplus_{\ell\in\seisuu}
 \bigoplus_{\alpha\in\cnum^{\ast}}
 \bigoplus_{\gminib\in S(q)}
 \nbign_{\ell,\alpha,\gminib}
\end{equation}
 such that each $\nbign_{\ell,\alpha,\gminib}$
 has a lattice $\nbigl_{\ell,\alpha,\gminib}$
 satisfying
\[
 (y_q^{\ell}\alpha^{-1}\Phi-(1+\gminib)\id)\nbigl_{\ell,\alpha,\gminib}
 \subset y_q^{-q}\nbigl_{\ell,\alpha,\gminib}.
\]
 In the case,
 such a lattice
 $\bigoplus\nbigl_{\ell,\alpha,\gminib}$
 is called an unramifiedly good lattice of
 $(\nbign,\Phi^{\ast})$.
 \hfill\qed 
\end{df}
\index{unramifiedly good lattice (formal difference module)}

\begin{df}
\label{df;20.7.28.10}
\index{unramifiedly good filtered difference module}
If $(\nbign,\Phi^{\ast})$ is unramified,
a filtered bundle $\nbigp_{\ast}\nbign$
over $\nbign$ is called unramifiedly good
for $(\nbign,\Phi^{\ast})$
if each $\nbigp_a\nbign$ is unramifiedly good lattices
of $(\nbign,\Phi^{\ast})$.
Such $(\nbigp_{\ast}\nbign,\Phi^{\ast})$
is called an unramifiedly good filtered difference
$\cnum(\!(y_q^{-1})\!)$-module. 
\hfill\qed
\end{df}

By Proposition \ref{prop;20.7.19.40},
Proposition \ref{prop;20.7.19.22},
and Theorem \ref{thm;20.7.19.41},
the following holds.
\begin{prop}
For any $\cnum(\!(y_q^{-1})\!)$-difference module
$(\nbign,\Phi^{\ast})$,
there exist $p\in q\seisuu_{>0}$
such that
$\varphi_{q,p}^{\ast}(\nbign,\Phi)$
 is unramified.  
\hfill\qed
\end{prop}

\begin{df}
\label{df;20.7.28.11}
\index{good filtered difference module}
A filtered bundle $\nbigp_{\ast}\nbign$ over $\nbign$
is called good with respect to $\Phi^{\ast}$
if there exists $p\in q\seisuu_{>0}$
such that 
$\varphi_{q,p}^{\ast}(\nbigp_{\ast}\nbign)$
is unramifiedly good
for $(\nbign,\Phi^{\ast})$. 
 Such $(\nbigp_{\ast}\nbign,\Phi^{\ast})$
is called a good filtered difference $\cnum(\!(y_q^{-1})\!)$-module.
 \hfill\qed
\end{df}

The following lemma is obvious.
\begin{lem}
\label{lem;21.9.17.70}
Let $(\nbigp_{\ast}\nbign,\Phi^{\ast})$
 be a good filtered difference module.
 Then, the filtered bundle
 $\nbigp_{\ast}\nbign$ is compatible
 with the slope decomposition
 in Theorem {\rm\ref{thm;20.7.19.41}},
 i.e.,
 $\nbigp_{\ast}\nbign
 =\bigoplus_{\omega\in\rnum}
 \nbigp_{\ast}(\nbigs_{\omega}\nbign)$.
 Moreover,
 for any $p\in Z(q,\omega)$,
 each $\nbigp_a\nbigs_{\omega}\nbign$
 is a $y_q^{p\omega}(\Phi^{\ast})^{p/q}$-invariant
 lattice.
\hfill\qed
\end{lem}

\begin{lem}
\label{lem;21.9.17.100}
Let $(\nbigp_{\ast}\nbign,\Phi^{\ast})$
 be a good filtered difference $\cnum(\!(y_q^{-1})\!)$-module.
 If $(\nbign,\Phi^{\ast})$ is unramified,
 then
 $(\nbigp_{\ast}\nbign,\Phi^{\ast})$ is an unramifiedly good
 filtered difference $\cnum(\!(y_q^{-1})\!)$-module. 
\end{lem}
\pf
Because $(\nbign,\Phi^{\ast})$
is assumed to be unramified,
there exists a decomposition
(\ref{eq;20.7.20.1})
of $\nbign$ as in Definition \ref{df;20.7.20.2}.
There exists $p\in q\seisuu_{>0}$
such that
$\varphi_{q,p}^{\ast}(\nbigp_{\ast}\nbign,\Phi^{\ast})$
is unramifiedly good.
There exists a decomposition
\begin{equation}
\label{eq;20.7.20.3}
 \varphi_{q,p}^{\ast}(\nbign)
  =\bigoplus_{\ell\in\seisuu}
  \bigoplus_{\alpha\in\cnum^{\ast}}
  \bigoplus_{\gminib\in S(p)}
  (\varphi_{q,p}^{\ast}(\nbign))_{\ell,\alpha,\gminib}
\end{equation}
as in Definition \ref{df;20.7.20.2}.
It is easy to observe that
$\varphi_{q,p}^{\ast}\nbign_{\ell,\alpha,\gminib}
=\nbign_{\ell(p/q),\alpha,\gminib}$.
The decomposition (\ref{eq;20.7.20.3})
is preserved by the natural action of
$\Gal_{q,p}$.
Because the filtered bundle
$\varphi^{\ast}(\nbigp_{\ast}\nbign)$
is compatible with
the decomposition (\ref{eq;20.7.20.3}),
we obtain that
the filtered bundle
$\nbigp_{\ast}\nbign$
is compatible with the decomposition
(\ref{eq;20.7.20.1}),
i.e.,
$\nbigp_{\ast}\nbign
=\bigoplus \nbigp_{\ast}\nbign_{\ell,\alpha,\gminib}$.
Because
\[
y_p^p\bigl(y_p^{p\ell/q}\alpha^{-1}
\Phi^{\ast}-(1+\gminib)\id\bigr)
\nbigp_b (\varphi_{q,p}^{\ast}\nbign)_{p\ell/q,\alpha,\gminib}
\subset
\nbigp_{b}(\varphi_{q,p}^{\ast}\nbign)_{p\ell/q,\alpha,\gminib},
\]
for any $b\in\real$,
we easily obtain
$\bigl(y_q^{\ell}\alpha^{-1}
\Phi^{\ast}-(1+\gminib)\id\bigr)
\nbigp_b \nbign_{\ell,\alpha,\gminib}
\subset
y_q^{-q}
\nbigp_{b}\nbign_{\ell,\alpha,\gminib}$.
Hence, $\nbigp_{\ast}\nbign$ is good
with respect to $\Phi^{\ast}$.
\hfill\qed

\begin{lem}
\label{lem;21.8.26.14}
Let $(\nbigp_{\ast}\nbign,\Phi^{\ast})$ be a good filtered
difference $\cnum(\!(y_q^{-1})\!)$-module.
Let $\nbign'\subset\nbign$ be a difference $\cnum(\!(y_q^{-1})\!)$-submodule.
We set $\nbign''=\nbign/\nbign'$.
They are equipped with the induced filtered bundle
$\nbigp_{\ast}\nbign'$ and $\nbigp_{\ast}\nbign''$, respectively.
Then, $(\nbigp_{\ast}\nbign',\Phi^{\ast})$
and  $(\nbigp_{\ast}\nbign'',\Phi^{\ast})$
are good filtered difference $\cnum(\!(y_q^{-1})\!)$-modules.
\end{lem}
\pf
We explain the proof for $\nbign'$.
The other case can be argued similarly.
By Lemma \ref{lem;21.8.25.42} and 
Lemma \ref{lem;21.8.25.41},
it is enough to study the case where
both $(\nbign,\Phi^{\ast})$
and $(\nbign',\Phi^{\ast})$ are unramified,
i.e.,
there exist decompositions
\begin{equation}
\label{eq;21.8.25.43}
\nbign=
 \bigoplus_{\ell\in\seisuu}\bigoplus_{\alpha\in\cnum^{\ast}}
 \bigoplus_{\gminib\in S(q)}
 \nbign_{\ell,\alpha,\gminib},
\end{equation}
\begin{equation}
\label{eq;21.8.25.44}
 \nbign'=
 \bigoplus_{\ell\in\seisuu}\bigoplus_{\alpha\in\cnum^{\ast}}
 \bigoplus_{\gminib\in S(q)}
 \nbign'_{\ell,\alpha,\gminib}
\end{equation}
as in (\ref{eq;20.7.20.1}).
By Lemma \ref{lem;21.8.25.40} and Lemma \ref{lem;20.7.19.10},
the decompositions (\ref{eq;21.8.25.43}) and (\ref{eq;21.8.25.44})
are compatible with the inclusion $\nbign'\subset\nbign$.
Hence, the induced filtration $\nbigp_{\ast}\nbign'$
is compatible with the decomposition (\ref{eq;21.8.25.44}).
It is easy to check that
each $\nbigp_a(\nbign')$ is an unramifiedly good lattice.
\hfill\qed

\begin{lem}
\label{lem;21.8.26.31}
 If $(\nbigp_{\ast}\nbign_i,\Phi^{\ast})$ $(i=1,2)$
 be good filtered difference $\cnum(\!(y_q^{-1})\!)$-modules,
 $(\nbigp_{\ast}(\nbign_1)\oplus\nbigp_{\ast}(\nbign_2),\Phi^{\ast})$,
 $(\nbigp_{\ast}(\nbign_1)\otimes\nbigp_{\ast}(\nbign_2),\Phi^{\ast})$
 and $(\Hom(\nbigp_{\ast}\nbign_1,\nbigp_{\ast}\nbign_2),\Phi^{\ast})$
 are also good filtered difference
 $\cnum(\!(y_q^{-1})\!)$-modules.
\end{lem}
\pf
By Lemma \ref{lem;21.8.26.1} and Lemma \ref{lem;21.8.26.5},
it is enough to consider the case where
$(\nbigp_{\ast}\nbign_i,\Phi^{\ast})$ are unramified.
Then, it is easy to check.
\hfill\qed

\subsection{The induced endomorphisms
   on the graded pieces}
\label{subsection;20.7.20.112}
   
Let $(\nbigp_{\ast}\nbign,\Phi^{\ast})$
be a good filtered difference $\cnum(\!(y_q^{-1})\!)$-module.
If $\nbign$ is unramified,
there exists the decomposition (\ref{eq;20.7.20.1})
as in Definition \ref{df;20.7.20.2}.
By the condition,
$y_q^q(y_q^{\ell}\alpha^{-1}
 \Phi^{\ast}-(1+\gminib)\id)$
 induces the $\cnum$-linear endomorphism
 $\gbigf_{a,\ell,\alpha,\gminib}$ of
 $\Gr^{\nbigp}_a(\nbign_{\ell,\alpha,\gminib})$.
We obtain 
the endomorphism
$\bigoplus_{\ell,\alpha,\gminib}
 \gbigf_{a,\ell,\alpha,\gminib}$
 on $\Gr^{\nbigp}_a(\nbign)$
which we denote by $\Res(\Phi^{\ast})$.
We obtain the monodromy weight filtration $W$
on $\Gr^{\nbigp}_a(\nbign)$
with respect to the nilpotent part of
$\Res(\Phi^{\ast})$.
It is compatible with the decomposition
$\Gr^{\nbigp}_a(\nbign)
=\bigoplus_{\ell,\alpha,\gminib}
\Gr^{\nbigp}_a(\nbign_{\ell,\alpha,\gminib})$.
\index{endomorphism $\Res(\Phi^{\ast})$}

In general,
there exists $p\in q\seisuu_{>0}$
such that
$\varphi_{q,p}^{\ast}(\nbign,\Phi^{\ast})$
is unramified.
We obtain the endomorphism
$\Res(\Phi^{\ast})$
on
$\Gr^{\nbigp}_a(\varphi_{q,p}^{\ast}\nbign)$.
\begin{lem}
\label{lem;20.7.22.5}
 $\Res(\Phi^{\ast})$ is equivariant
 with respect to
 the natural $\Gal_{q,p}$-action
 on $\Gr^{\nbigp}_a(\varphi_{q,p}^{\ast}\nbign)$.
\end{lem}
\pf
It follows from
Lemma \ref{lem;20.7.22.4}
and Lemma \ref{lem;20.7.22.3}.
\hfill\qed

\vspace{.1in}

By Lemma \ref{lem;20.7.22.5},
$\Res(\Phi^{\ast})$ is compatible
with the canonical decomposition (\ref{eq;20.7.20.101}).
By using the identification
$\Gr^{\nbigp}_a(\nbign)
\simeq
\nbigg_0\Gr^{\nbigp}_{pa/q}(\varphi_{q,p}^{\ast}\nbign)$,
we obtain an endomorphism
on $\Gr^{\nbigp}_a(\nbign)$,
denoted by $\Res(\Phi^{\ast})$.
We also obtain the monodromy weight filtration $W$
on $\Gr^{\nbigp}_a(\nbign)$
with respect to the nilpotent part of
$\Res(\Phi^{\ast})$.
\index{monodromy weight filtration $W$}

\section{Geometrization of formal difference modules}
\label{subsection;20.7.30.22}

\subsection{Ringed spaces}
\label{subsection;20.7.31.1}
A ringed space is defined to be a topological space $X$
equipped with a sheaf of algebras $A_X$ on $X$,
called the structure sheaf.
\index{ringed space}
In this monograph, the structure sheaf of any ringed space
is assumed to be commutative.
A morphism of ringed spaces
$F:(X,A_X)\lrarr (Y,A_Y)$ is a continuous map
$F:X\lrarr Y$
with a morphism of sheaves of algebras
$F^{-1}(A_Y)\lrarr A_X$.
For an $A_Y$-module $\nbigf$,
we obtain an $A_X$-module
$F^{\ast}(\nbigf):=
A_X\otimes_{F^{-1}A_Y}F^{-1}(\nbigf)$.
If $X$ is a subspace of $Y$
and if $F^{-1}(A_Y)=A_X$,
$F^{\ast}(\nbigf)$ is denoted as $\nbigf_{|X}$.
Even if $F^{-1}(A_Y)\neq A_X$,
we shall often denote
$F^{\ast}(\nbigf)$ as $\nbigf_{|X}$
to simplify the description
if there is no risk of confusion.
\index{pull back}

For an $A_X$-module $\nbigf'$,
we obtain an $A_Y$-module $F_{\ast}(\nbigf')$
as the standard push-forward.
\index{push-forward}

\subsection{Some formal spaces}
\label{subsection;20.7.21.1}

For any $q\in\seisuu_{\geq 1}$,
let $\inftyhat_{y,q}$ denote
the ringed space which consists of
the point $\infty_{y,q}$ with the ring
$\nbigo_{\inftyhat_{y,q}}:=\cnum[\![y_q^{-1}]\!]$.
\index{ringed space $\inftyhat_{y,q}$}
\index{point $\infty_{y,q}$}
In other words,
$\inftyhat_{y,q}$ is the completion of
the projective line $\proj^1_{y_q}$ along $y_q=\infty$.
We also set
$\nbigo_{\inftyhat_{y,q}}(\ast \infty_{y,q})
:=\cnum(\!(y_q^{-1})\!)$.
\index{ring $\nbigo_{\inftyhat_{y,q}}(\ast \infty_{y,q})$}
If $p\in q\seisuu_{\geq 1}$,
we obtain the naturally defined morphism of the ringed spaces
$\varphi_{q,p}:\inftyhat_{y,p}\lrarr \inftyhat_{y,q}$
induced by
the inclusions
$\varphi_{q,p}^{\ast}:
\cnum[\![y_q^{-1}]\!]
\lrarr
\cnum[\![y_p^{-1}]\!]$.
\index{map $\varphi_{q,p}$}
There also exists the natural ring homomorphism
$\varphi_{q,p}^{\ast}:
\nbigo_{\inftyhat_{y,q}}(\ast\infty)
\lrarr
\nbigo_{\inftyhat_{y,p}}(\ast\infty)$.
The space $\inftyhat_{y,1}$
is also denoted by $\inftyhat_y$.
\index{ringed space $\inftyhat_y$}

We set
$H^{\cov}_{\infty,q}:=\real_t\times\{\infty_{y,q}\}$.
\index{space $H^{\cov}_{\infty,q}$}
Let $\nbigo_{\Hhat^{\cov}_{\infty,q}}$
(resp. $\nbigo_{\Hhat^{\cov}_{\infty,q}}(\ast H^{\cov}_{\infty,q})$)
denote the sheaf of
locally constant functions on $H^{\cov}_{\infty,q}$
to $\cnum[\![y_q^{-1}]\!]$
(resp. $\cnum(\!(y_q^{-1})\!)$).
Let $\Hhat^{\cov}_{\infty,q}$ denote the ringed space
which consists of the topological space
$H^{\cov}_{\infty,q}$
with the sheaf of algebras $\nbigo_{\Hhat^{\cov}_{\infty,q}}$.
\index{formal space $\Hhat^{\cov}_{\infty,q}$}
\index{sheaf $\nbigo_{\Hhat^{\cov}_{\infty,q}}$}
\index{sheaf $\nbigo_{\Hhat^{\cov}_{\infty,q}}(\ast H^{\cov}_{\infty,q})$}

Let $T\in\real_{>0}$.
Let $\kappa$ denote the natural $\seisuu$-action
on $H^{\cov}_{\infty,q}$
induced by
$\kappa_n(t)=t+nT$ $(n\in\seisuu)$.
\index{action $\kappa$}
We obtain the isomorphism
$\kappa_n^{\ast}\nbigo_{\Hhat^{\cov}_{\infty,q}}
\simeq
\nbigo_{\Hhat^{\cov}_{\infty,q}}$
induced by
$\kappa_n^{\ast}(y)=y+2\sqrt{-1}nT\lambda$.
In this sense,
the ringed space
$\Hhat^{\cov}_{\infty,q}$
is equipped with a $\seisuu$-action $\kappa$.
The sheaf
$\nbigo_{\Hhat^{\cov}_{\infty,q}}(\ast H^{\cov}_{\infty,q})$
is also naturally $\seisuu$-equivariant.

Let $H_{\infty,q}$ denote the quotient space
of $H^{\cov}_{\infty,q}$
by the action $\kappa$.
\index{space $H_{\infty,q}$}
We obtain the sheaves
$\nbigo_{\Hhat_{\infty,q}}$
and
$\nbigo_{\Hhat_{\infty,q}}(\ast H_{\infty,q})$
as the descents of
$\nbigo_{\Hhat^{\cov}_{\infty,q}}$
and
$\nbigo_{\Hhat^{\cov}_{\infty,q}}(\ast H^{\cov}_{\infty,q})$,
respectively.
Let $\Hhat_{\infty,q}$ denote the ringed space
obtained as
the topological space
$H_{\infty,q}$
with the sheaf of algebras
$\nbigo_{\Hhat_{\infty,q}}$.
\index{ringed space $\Hhat_{\infty,q}$}
\index{sheaf $\nbigo_{\Hhat_{\infty,q}}$}
\index{sheaf $\nbigo_{\Hhat_{\infty,q}}(\ast H_{\infty,q})$}

The projections 
$\Hhat^{\cov}_{\infty,q}\lrarr
 \Hhat_{\infty,q}$
and 
$H^{\cov}_{\infty,q}\lrarr
 H_{\infty,q}$
 are denoted by $\varpi_q$.
\index{morphism $\varpi_q$}
We also denote
$\Hhat^{\cov}_{\infty,q}$,
$\Hhat_{\infty,q}$, etc.,
by 
$\Hhat^{\cov}_{\infty,q,T}$,
$\Hhat_{\infty,q,T}$, etc.,
respectively,
when we emphasize the dependence on $T$.
\index{ringed space $\Hhat^{\cov}_{\infty,q,T}$}
\index{space $\Hhat^{\cov}_{\infty,q,T}$}

\begin{rem}
The space $\Hhat_{\infty}$
and $\Hhat^{\cov}_{\infty}$
are the formal completion of
the spaces
$\nbigmbar^{\lambda}$
and $\Mbar^{\lambda}$
along $H^{\lambda}_{\infty}$
and $H^{\lambda\,\cov}_{\infty}$.
(See {\rm\S\ref{subsection;17.10.2.1}}
for the notation.)
\hfill\qed
\end{rem}

The ringed space $\real$ with the sheaf of $\cnum$-valued
locally constant functions is also denoted by $\real$.
\index{ringed space $\real$}
We set $S^1_T:=\real/T\seisuu$.
\index{space $S^1_T$}
The projection 
$\Hhat^{\cov}_{\infty,q}\lrarr \real$
induces a morphism
$\pi_q:\Hhat_{\infty,q}\lrarr S^1_T$.
\index{projection $\pi_q$}
It also induces an isomorphism
$H_{\infty,q}\simeq S^1_T$.
For $t\in S^1_T$,
let $\Hhat_{\infty,q}\langle t\rangle$
denote the formal space
obtained as the point set $\{t\}\subset H_{\infty,q}\simeq S^1_T$
with the ring
$\nbigo_{\Hhat_{\infty,q},t}$.
\index{ringed space $\Hhat^t_{\infty,q}$}
For any $\nbigo_{\Hhat_{\infty,q}}(\ast H_{\infty,q})$-module $\nbige$,
let $\nbige_{|\Hhat_{\infty,q}\langle t\rangle}$ denote
the pull back of $\nbige$
by the natural morphism
$\Hhat_{\infty,q}\langle t\rangle\lrarr \Hhat_{\infty,q}$.
\index{restriction $\nbige_{|\Hhat^t_{\infty,q}}$}
Note that if we choose $\ttilde\in \real$ which is mapped to $t\in S^1_T$,
there exists the induced isomorphism
$\{\ttilde\}\times\inftyhat_{y,q}
\simeq
\Hhat_{\infty,q}\langle t\rangle$.
Hence, we may regard $\nbige_{|\Hhat_{\infty,q}\langle t\rangle}$
as $\cnum(\!(y_q^{-1})\!)$-module
in a way depending on the lift $\ttilde$.

Let $\nbige^{\cov}$ be 
a locally free
$\nbigo_{\Hhat^{\cov}_{\infty,q}}(\ast H^{\cov}_{\infty,q})$-module
of finite rank.
\index{sheaf $\nbige^{\cov}$}
For any $t\in\real$,
let $\nbige^{\cov}_{|\{t\}\times\inftyhat_{y,q}}$
denote the pull back of
$\nbige^{\cov}$ by
$\{t\}\times\inftyhat_{y,q}
\lrarr
 \Hhat^{\cov}_{\infty,q}$.
\index{restriction $\nbige^{\cov}_{|\{t\}\times\inftyhat_{y,q}}$}
For any $t_1,t_2\in\real$,
there exists the natural isomorphism
$\nbige^{\cov}_{|\{t_1\}\times\inftyhat_{y,q}}
\lrarr
\nbige^{\cov}_{|\{t_2\}\times\inftyhat_{y,q}}$.

\subsection{Difference modules
and $\nbigo_{\Hhat_{\infty,q}}(\ast H_{\infty,q})$-modules}
\label{subsection;20.7.18.4}

Let $\nbige$ be a locally free 
$\nbigo_{\Hhat_{\infty,q}}(\ast H_{\infty,q})$-module
of finite rank.
We obtain 
$\nbige^{\cov}:=\varpi_q^{-1}\nbige$ on 
$\Hhat^{\cov}_{\infty,q}$.
There exists the isomorphism
\[
 \Pi_{T,0}:
 \nbige^{\cov}_{|\{0\}\times\inftyhat_{y,q}}
\lrarr
 \nbige^{\cov}_{|\{T\}\times\inftyhat_{y,q}}
\]
induced by the parallel transport.
There also exists the natural identification
\[
 \nbige^{\cov}_{|\{T\}\times\inftyhat_{y,q}}
=
 \nbige_{|\Hhat_{\infty,q}\langle 0\rangle}
=
 \nbige^{\cov}_{|\{0\}\times\inftyhat_{y,q}}.
\]
We define the difference operator
$\Phi^{\ast}:
\nbige_{|\Hhat_{\infty,q}\langle 0\rangle}
\lrarr
\nbige_{|\Hhat_{\infty,q}\langle 0\rangle}$
as the composite of the following morphisms:
\[
 \nbige_{|\Hhat_{\infty,q}\langle 0\rangle}=
 \nbige^{\cov}_{|\{0\}\times\inftyhat_{y,q}}
 \stackrel{\Pi_{T,0}}{\lrarr}
  \nbige^{\cov}_{|\{T\}\times\inftyhat_{y,q}}
=\nbige_{|\Hhat_{\infty,q}\langle 0\rangle}.
\]
\index{formal difference module
$\nbige_{|\Hhat_{\infty,q}\langle 0\rangle}$}
It is also equivalent to the following constructions.
Let $H^0\bigl(H^{\cov}_{\infty,q},\nbige^{\cov}\bigr)$
denote the space of the sections of
$\nbige^{\cov}$ on $H^{\cov}_{\infty,q}$.
We define
$\Phi^{\ast}$ on 
$H^0\bigl(H^{\cov}_{\infty,q},\nbige^{\cov}\bigr)$
as the composite of the following maps:
\[
 H^0\bigl(H^{\cov}_{\infty,q},\nbige^{\cov}\bigr)
 \stackrel{\kappa_1^{\ast}}{\lrarr}
 H^0\bigl(H^{\cov}_{\infty,q},\kappa_1^{\ast}\nbige^{\cov}\bigr)
=H^0\bigl(H^{\cov}_{\infty,q},\nbige^{\cov}\bigr).
\]
Here,
we obtain the second equality
by the identifications
$\kappa_1^{\ast}\nbige^{\cov}
=(\varpi_q\circ\kappa_1)^{\ast}\nbige
=\varpi_q^{\ast}\nbige=\nbige^{\cov}$.
Under the natural isomorphism
$H^0\bigl(H^{\cov}_{\infty,q},\nbige^{\cov}\bigr)
\simeq
 \nbige^{\cov}_{|\{0\}\times\inftyhat_{y,q}}
 =\nbige_{|\Hhat_{\infty,q}\langle 0\rangle}$,
the two constructions are the same.
Note that the constructions are compatible with
direct sums, tensor products and inner homomorphisms. 

\begin{prop}
\label{prop;17.11.20.1}
The above construction induces an equivalence of
the categories of
locally free $\nbigo_{\Hhat_{\infty,q}}(\ast H_{\infty,q})$-modules
of finite rank
and $(2\sqrt{-1}\lambda T)$-difference $\cnum(\!(y_q^{-1})\!)$-modules.
\end{prop}
\pf
Let $(\nbign,\Phi^{\ast})$ be
a $(2\sqrt{-1}\lambda T)$-difference
$\cnum(\!(y_q^{-1})\!)$-module.
Let $p_2:\Hhat^{\cov}_{\infty,q}\lrarr \inftyhat_{y,q}$
denote the projection.
We obtain 
the $\nbigo_{\Hhat^{\cov}_{\infty,q}}(\ast H^{\cov}_{\infty,q})$-module
$p_2^{-1}\nbign$
which is naturally $\seisuu$-equivariant
by the action
$\kappa_n^{\ast}(p_2^{-1}(s))=p_2^{-1}\bigl((\Phi^{\ast})^n(s)\bigr)$.
We obtain 
an $\nbigo_{\Hhat_{\infty,q}}(\ast H_{\infty,q})$-module
as the descent.
This is a converse construction of the previous one.
\hfill\qed

\vspace{.1in}
The properties for difference $\cnum(\!(y_q^{-1})\!)$-modules
are translated to the properties of
$\nbigo_{\Hhat_{\infty,q}}(\ast H_{\infty,q})$-modules.

\begin{df}
\label{df;20.7.31.10}
Let $\nbige$ be a locally free
$\nbigo_{\Hhat_{\infty,q}}(\ast H_{\infty,q})$-module.
\begin{itemize}
 \item
\index{level $0$ ($\nbigo_{\Hhat_{\infty,q}}(\ast H_{\infty,q})$-module)}
\index{level $\leq 1$ ($\nbigo_{\Hhat_{\infty,q}}(\ast H_{\infty,q})$-module)}
The level of $\nbige$ is $0$
(resp. less than $1$)
if the level of
the difference $\cnum(\!(y_q^{-1})\!)$-module
$\nbige_{|\Hhat_{\infty,q}\langle 0\rangle}$ is $0$
(resp. less than $1$).
(See {\rm\S\ref{subsection;20.7.20.30}}.)
 \item
 \index{pure slope ($\nbigo_{\Hhat_{\infty,q}}(\ast H_{\infty,q})$-module)}
 $\nbige$ has pure slope $\omega\in\rnum$
 if the difference $\cnum(\!(y_q^{-1})\!)$-module
 $\nbige_{|\Hhat_{\infty,q}\langle 0\rangle}$
 has pure slope $\omega$.
(See {\rm\S\ref{subsection;20.7.20.31}})
 \item
 \index{unramified ($\nbigo_{\Hhat_{\infty,q}}(\ast H_{\infty,q})$-module)}
 $\nbige$ is unramified
 if 
 the difference $\cnum(\!(y_q^{-1})\!)$-module
 $\nbige_{|\Hhat_{\infty,q}\langle 0\rangle}$ is unramified. 
(See Definition {\rm\ref{df;20.7.20.2}}.)
\hfill\qed
\end{itemize}
\end{df}

Any $\nbigo_{\Hhat_{\infty,q}}(\ast H_{\infty,q})$-module
$\nbige$ of finite rank
has the unique slope decomposition 
\begin{equation}
\label{eq;20.7.20.100}
 \nbige=\bigoplus_{\omega\in\rnum}\nbigs_{\omega}\nbige,
\end{equation}
where each $\nbigs_{\omega}\nbige$
has pure slope $\omega$.
\index{slope decomposition 
($\nbigo_{\Hhat_{\infty,q}}(\ast H_{\infty,q})$-module)}

\subsection{Lattices and the induced local systems}

Let $\nbige$ be 
a locally free
$\nbigo_{\Hhat_{\infty,q}}(\ast H_{\infty,q})$-module
of finite rank.
A lattice of $\nbige$ means
a locally free $\nbigo_{\Hhat_{\infty,q}}$-submodule
$\nbige_0\subset\nbige$
such that
$\nbigo_{\Hhat_{\infty,q}}(\ast H_{\infty,q})\otimes
\nbige_0=\nbige$.
\index{lattice}
Similarly,
for a locally free
$\nbigo_{\Hhat^{\cov}_{\infty,q}}(\ast H^{\cov}_{\infty,q})$-module
$\nbige^{\cov}$ of finite rank,
a lattice of $\nbige^{\cov}$ means
a locally free $\nbigo_{\Hhat^{\cov}_{\infty,q}}$-submodule
$\nbige^{\cov}_0\subset\nbige^{\cov}$
such that
$\nbigo_{\Hhat^{\cov}_{\infty,q}}(\ast H^{\cov}_{\infty,q})\otimes
\nbige^{\cov}_0=\nbige^{\cov}$.

Let $\nbige$ be a locally free
$\nbigo_{\Hhat_{\infty,q}}(\ast H_{\infty,q})$-module
of finite rank.
We set
$\nbige^{\cov}:=\varpi_q^{-1}(\nbige)$.
Any lattice $\nbige_0$ of $\nbige$ induces
a lattice
$\nbige_0^{\cov}=\varpi_q^{-1}(\nbige_0)$
of $\nbige^{\cov}$.
For any $k\in\seisuu_{>0}$,
we obtain a $\seisuu$-equivariant
local system
$\nbige_0^{\cov}/y_q^{-k}\nbige_0^{\cov}$
on $H^{\cov}_{\infty,q}$.
It induces a local system
$\Loc_k(\nbige_0)$ on $H_{\infty,q}$.
\index{local system $\Loc_k(\nbige_0)$}
The following lemma is obvious.
\begin{lem}
Let $\nbige$ be 
 an $\nbigo_{\Hhat_{\infty,q}}(\ast H_{\infty,q})$-module.
\begin{itemize}
 \item There exists a lattice $\nbige_0$ of $\nbige$
      if and only if $\nbige$ has pure slope $0$.
 \item The level of $\nbige$ is less than $1$
 if and only if
 there exists
 a lattice $\nbige_0$ of $\nbige$
 such that the monodromy of
 $\Loc_1(\nbige_0)$ is unipotent.
 \item  The level of $\nbige$ is $0$
 if and only if
 there exists a
 lattice $\nbige_0$ of $\nbige$
 such that the monodromy of
 $\Loc_q(\nbige_0)$ is the identity.
\hfill\qed
\end{itemize}
\end{lem}

\section{Filtered bundles in the formal case}
\label{subsection;20.7.30.23}

\subsection{Pull back and descent of
$\nbigo_{\Hhat_{\infty,p}}(\ast H_{\infty,p})$-modules}

For any $p\in q\seisuu_{\geq 1}$,
there exists the naturally induced morphism
$\nbigr_{q,p}:\Hhat_{\infty,p}\lrarr \Hhat_{\infty,q}$,
whose underlying map
$\nbigr_{q,p}:
H_{\infty,p}\lrarr H_{\infty,q}$ is the identity,
and the morphism of sheaves
$\nbigr_{q,p}^{\ast}\nbigo_{\Hhat_{\infty,q}}
\lrarr \nbigo_{\Hhat_{\infty,p}}$
is induced by the extension
$\cnum[\![y_q^{-1}]\!]
\lrarr\cnum[\![y_p^{-1}]\!]$. \index{morphism $\nbigr_{q,p}$}
\index{morphism $\nbigr_{q,p}$}
The $\Gal_{q,p}$-action on $\cnum[\![y_p^{-1}]\!]$
induces a $\Gal_{q,p}$-action on
$\Hhat_{\infty,p}$,
which we may regard as
the Galois groups action
for the ramified covering
$\Hhat_{\infty,p}\lrarr\Hhat_{\infty,q}$.

For any $\nbigo_{\Hhat_{\infty,q}}$-module $\nbige$,
we obtain
the $\nbigr_{q,p}^{\ast}\nbige$
as the pull back
as in \S\ref{subsection;20.7.31.1}.
If $\nbige$ is 
an $\nbigo_{\Hhat_{\infty,q}}(\ast H_{\infty,q})$-module,
then 
$\nbigr_{q,p}^{\ast}\nbige$ 
is naturally an
$\nbigo_{\Hhat_{\infty,p}}(\ast H_{\infty,p})$-module.
\index{pull back}

Let $\nbige_1$ be 
an $\nbigo_{\Hhat_{\infty,p}}$-module.
We obtain the $\nbigo_{\Hhat_{\infty,q}}$-module
$\nbigr_{q,p\ast}\nbige_1$
obtained as the push-forward.
\index{push-forward}
If $\nbige_1$ is an $\nbigo_{\Hhat_{\infty,p}}(\ast H_{\infty,p})$-module,
then $\nbigr_{q,p\ast}\nbige_1$ is also
an $\nbigo_{\Hhat_{\infty,q}}(\ast H_{\infty,q})$-module.
If  $\nbige_1$ is a $\Gal_{q,p}$-equivariant
$\nbigo_{\Hhat_{\infty,p}}$-module,
$\nbigr_{q,p\ast}\nbige_1$ is 
equipped with the induced $\Gal_{q,p}$-action.
The invariant part of $\nbigr_{q,p\ast}\nbige_1$
is called the descent of $\nbige_1$
with respect to the $\Gal_{q,p}$-action.
\index{descent}

\subsection{Filtered bundles}

We consider the family version of
filtered bundles in \S\ref{subsection;20.7.20.20}
in a straightforward way.

\begin{df}
\label{df;20.7.21.2}
\index{filtered bundle}
A filtered bundle over
a locally free
$\nbigo_{\Hhat_{\infty,q}}(\ast H_{\infty,q})$-module
$\nbige$ of finite rank
is defined to be a family 
of filtered bundles
$\nbigp_{\ast}(\nbige_{|\Hhat_{\infty,q}\langle t\rangle})$
$(t\in S^1_T)$
over $\nbige_{|\Hhat_{\infty,q}\langle t\rangle}$.
\index{filtered bundle}
The family is also denoted as
$\nbigp_{\ast}\nbige$, for simplicity.
 We also say that
 $\nbigp_{\ast}\nbige$ is a filtered bundle
 on $(\Hhat_{\infty,q},H_{\infty,q})$.
 \index{filtered bundle on $(\Hhat_{\infty,q}H_{\infty,q})$}
 \hfill\qed
\end{df}

\begin{rem}
In Definition {\rm\ref{df;20.7.21.2}},
we do not impose any condition 
between 
$\nbigp_{\ast}(\nbige_{|\Hhat_{\infty,q}\langle t_1\rangle})$
and
$\nbigp_{\ast}(\nbige_{|\Hhat_{\infty,q}\langle t_2\rangle})$
under the natural isomorphisms
$\nbige_{|\Hhat_{\infty,q}\langle t_1\rangle}\simeq
\nbige_{|\Hhat_{\infty,q}\langle t_2\rangle}$.
We shall explain what condition should be imposed
in {\rm\S\ref{subsection;17.10.7.1}--\S\ref{subsection;21.8.16.1}},
and we shall eventually introduce the goodness condition
in Definition {\rm\ref{df;20.7.31.12}}.
\hfill\qed
\end{rem}

\subsubsection{Subbundles and quotient bundles}

Let $\nbigp_{\ast}\nbige$ be a filtered bundle
on $(\Hhat_{\infty,q},H_{\infty,q})$.
For any locally free $\nbigo_{\Hhat_{\infty,q}}(\ast H_{\infty,q})$-module
$\nbige'\subset\nbige$,
we obtain the filtered bundle $\nbigp_{\ast}\nbige'$
over $\nbige'$
by setting
$\nbigp_a(\nbige'_{|\Hhat_{\infty,q}\langle t\rangle})
:=\nbigp_a(\nbige_{|\Hhat_{\infty,q}\langle t\rangle})
\cap \nbige'_{|\Hhat_{\infty,q}\langle t\rangle}$.
\index{induced filtered bundle (subbundle)}
We also obtain the filtered bundle over
$\nbige''=\nbige/\nbige'$
by setting
$\nbigp_{a}(\nbige''_{|\Hhat_{\infty,q}\langle t\rangle})$
as the image of
$\nbigp_a(\nbige_{|\Hhat_{\infty,q}\langle t\rangle})
\lrarr
\nbige''_{|\Hhat_{\infty,q}\langle t\rangle}$.
\index{induced filtered bundle (quotient bundle)}

\subsubsection{Basic functoriality}

Let $\nbige_i$ $(i=1,2)$ be locally free
$\nbigo_{\Hhat_{\infty,q}}(\ast H_{\infty,q})$-modules.
Let
$\nbigp_{\ast}(\nbige_{i})
=\bigl(
\nbigp_{\ast}(\nbige_{i|\Hhat_{\infty,q}\langle t\rangle})
\,\big|\,t\in S^1_T
\bigr)$
be filtered bundles over $\nbige_i$.
The induced filtered bundles
$\Bigl(
 \nbigp_{\ast}
  (\nbige_{1|\Hhat_{\infty,q}\langle t\rangle})
  \oplus
  \nbigp_{\ast}
  (\nbige_{2|\Hhat_{\infty,q}\langle t\rangle})
  \Big|\,
  t\in S^1_T
  \Bigr)$ over $\nbige_1\oplus\nbige_2$
is denoted by
$\nbigp_{\ast}(\nbige_1\oplus\nbige_2)$.
The filtered bundle
$\Bigl(
 \nbigp_{\ast}
  (\nbige_{2|\Hhat_{\infty,q}\langle t\rangle})
   \otimes
   \nbigp_{\ast}
  (\nbige_{2|\Hhat_{\infty,q}\langle t\rangle})
  \,\Big|\,
  t\in S^1_T
\Bigr)$
over $\nbige_1\otimes\nbige_2$
is denoted by
$\nbigp_{\ast}(\nbige_1\otimes\nbige_2)$.
Let $\nhom(\nbige_1,\nbige_2)$
denote the sheaf of
$\nbigo_{\Hhat_{\infty,q}}(\ast H_{\infty,q})$-homomorphisms
from $\nbige_1$ to $\nbige_2$.
Similarly, we obtain a naturally defined filtered bundle
over $\nhom(\nbige_1,\nbige_2)$,
which is denoted by
$\nhom(\nbigp_{\ast}\nbige_1,\nbigp_{\ast}\nbige_2)$.
\index{direct sum (filtered bundle)}
\index{tensor product (filtered bundle)}
\index{inner homomorphism (filtered bundle)}
Let
$\nbigp_{\ast}\bigl(
\nbigo_{\Hhat_{\infty,q}}(\ast H_{\infty,q})\bigr)$
denote the filtered bundle over
$\nbigo_{\Hhat_{\infty,q}}(\ast H_{\infty,q})$
defined by
$\nbigp_a\bigl(
 \nbigo_{\Hhat_{\infty,q}}(\ast H_{\infty,q})
 \bigr)
 =y_q^{[a]}\nbigo_{\Hhat_{\infty,q}}$.
We define 
$\nbigp_{\ast}\nbige^{\lor}
=\nhom\bigl(
\nbigp_{\ast}\nbige,
\nbigp_{\ast}(\nbigo_{\Hhat_{\infty,q}}(\ast H_{\infty,q}))
\bigr)$.
\index{dual (filtered bundle)}

\subsubsection{Pull back}

Let $\nbigp_{\ast}\nbige$ be a filtered bundle
over a locally free
$\nbigo_{\Hhat_{\infty,q}}(\ast H_{\infty,q})$-module
$\nbige$.
Because 
$\nbigr_{q,p}^{\ast}(\nbige)_{|\Hhat_{\infty,p}\langle t\rangle}$
$(t\in S^1_T)$
are the pull back of
$\nbige_{|\Hhat_{\infty,q}\langle t\rangle}$
by $\varphi_{q,p}:\inftyhat_{y,p}\lrarr \inftyhat_{y,q}$,
we obtain the
filtered bundles
$\nbigp_{\ast}\bigl(
\nbigr_{q,p}^{\ast}(\nbige)_{|\Hhat_{\infty,p}\langle t\rangle}
\bigr)=
\varphi_{q,p}^{\ast}
\bigl(\nbigp_{\ast}(\nbige_{|\Hhat_{\infty,q}\langle t\rangle})
\bigr)$
as in \S\ref{subsection;20.7.20.20}.
 The tuple
$\bigl(
\nbigp_{\ast}(
 \nbigr_{q,p}^{\ast}(\nbige)_{|\Hhat_{\infty,p}\langle t\rangle})
 \,\big|\,
 t\in S^1_T
\bigr)$
is denoted by
$\nbigr_{q,p}^{\ast}(\nbigp_{\ast}\nbige)$,
and called the pull back of
$\nbigp_{\ast}\nbige$.
\index{pull back (filtered bundle)}
We obtain the following lemma from Lemma \ref{lem;21.8.25.42}.
\begin{lem}
\label{lem;21.8.26.10}
Let $\nbige$ and $\nbigp_{\ast}\nbige$ be as above.
Let $\nbige'$ be
a free $\nbigo_{\Hhat_{\infty,q}}(\ast H_{\infty,q})$-submodule
of $\nbige$.
We put $\nbige''=\nbige/\nbige'$.
They are equipped with the induced filtered bundle
$\nbigp_{\ast}\nbige'$ and $\nbigp_{\ast}\nbige''$,
respectively.
Then,
$\nbigr_{q,p}^{\ast}(\nbigp_{\ast}\nbige')$
(resp. $\nbigr_{q,p}^{\ast}(\nbigp_{\ast}\nbige'')$)
is equal to the filtered bundle over $\nbigr_{q,p}^{\ast}\nbige'$
(resp. $\nbigr_{q,p}^{\ast}\nbige'$)
induced by $\nbigr_{q,p}^{\ast}(\nbigp_{\ast}\nbige)$
and the inclusion 
$\nbigr_{q,p}^{\ast}\nbige'\subset\nbigr_{q,p}^{\ast}\nbige$
(the projection
$\nbigr_{q,p}^{\ast}\nbige\lrarr\nbigr_{q,p}^{\ast}\nbige''$).
\hfill\qed
\end{lem}

We obtain the following lemma from Lemma \ref{lem;21.8.26.1}.
\begin{lem}
\label{lem;21.8.26.30}
Let $\nbigp_{\ast}\nbige_i$ be filtered bundles over
$(\Hhat_{\infty,q},H_{\infty,q})$.
Then, there exist natural isomorphisms
 $\nbigr_{q,p}^{\ast}(\nbigp_{\ast}\nbige_1\oplus\nbigp_{\ast}\nbige_2)
 \simeq
 \nbigr_{q,p}^{\ast}(\nbigp_{\ast}\nbige_1)
 \oplus
 \nbigr_{q,p}^{\ast}(\nbigp_{\ast}\nbige_2)$,
 $\nbigr_{q,p}^{\ast}(\nbigp_{\ast}\nbige_1\otimes\nbigp_{\ast}\nbige_2)
 \simeq
 \nbigr_{q,p}^{\ast}(\nbigp_{\ast}\nbige_1)
 \otimes
 \nbigr_{q,p}^{\ast}(\nbigp_{\ast}\nbige_2)$,
 and
 $\nbigr_{q,p}^{\ast}\nhom\bigl(\nbigp_{\ast}\nbige_1,
 \nbigp_{\ast}\nbige_2\bigr)
\simeq
 \nhom\bigl(\nbigr_{q,p}^{\ast}(\nbigp_{\ast}\nbige_1),
 \nbigr_{q,p}^{\ast}(\nbigp_{\ast}\nbige_2)\bigr)$.
\hfill\qed
\end{lem}

\subsubsection{Push-forward}

Let $\nbigp_{\ast}\nbige_1$ be a filtered bundle
over a locally free 
$\nbigo_{\Hhat_{\infty,p}}(\ast H_{\infty,p})$-module
$\nbige_1$.
Because
$\nbigr_{q,p\ast}(\nbige_1)_{|\Hhat_{\infty,q}\langle t\rangle}
=\varphi_{q,p\ast}(\nbige_{1|\Hhat_{\infty,p}\langle t\rangle})$,
we obtain the induced filtered bundle
$\nbigp_{\ast}\bigl(
 \nbigr_{q,p\ast}(\nbige_1)
_{|\Hhat_{\infty,q}\langle t\rangle}\bigr)
=\varphi_{q,p\ast}\bigl(
\nbigp_{\ast}(\nbige_{1|\Hhat_{\infty,p}\langle t\rangle})
 \bigr)$.
The tuple
$\bigl(
\nbigp_{\ast}(\nbigr_{q,p\ast}(\nbige_1)_{|\Hhat_{\infty,q}\langle t\rangle})
\,\big|\,
t\in S^1_T
\bigr)$
is denoted by
$\nbigr_{q,p\ast}(\nbigp_{\ast}\nbige_1)$,
and called the push-forward of
$\nbigp_{\ast}\nbige_1$.
\index{push-forward (filtered bundle)}
We obtain the following lemma from Lemma \ref{lem;21.8.26.2}.
\begin{lem}
Let $\nbige_1$ and $\nbigp_{\ast}\nbige_1$ be as above.
Let $\nbige'_1\subset\nbige_1$ be a locally free
$\nbigo_{\Hhat_{\infty,p}}(\ast H_{\infty,p})$-submodule.
We set $\nbige''_1=\nbige_1/\nbige''_1$. 
They are equipped with the induced filtered bundles
$\nbigp_{\ast}(\nbige'_1)$
and  $\nbigp_{\ast}(\nbige''_1)$, respectively.
Then, $\nbigr_{q,p\ast}(\nbigp_{\ast}\nbige'_1)$
(resp. $\nbigr_{q,p\ast}(\nbigp_{\ast}\nbige''_1)$)
is equal to the filtered bundle induced by 
$\nbigr_{q,p\ast}(\nbigp_{\ast}\nbige_1)$
and the inclusion 
$\nbigr_{q,p\ast}(\nbige'_1)\subset\nbigr_{q,p\ast}(\nbige_1)$
(the projection
$\nbigr_{q,p\ast}(\nbige_1)\lrarr\nbigr_{q,p\ast}(\nbige''_1)$).
\hfill\qed
\end{lem}

Let $\nbigp_{\ast}\nbige$ be a filtered bundle
on $(\Hhat_{\infty,q},H_{\infty,q})$.
Let $\nbigp_{\ast}\nbige_1$ be a filtered bundle
on $(\Hhat_{\infty,p},H_{\infty,p})$.
We obtain the following lemma from Lemma \ref{lem;21.8.26.3}.
\begin{lem}
The natural isomorphism
$\nbigr_{q,p\ast}(\nbige_1\otimes\nbigr_{q,p}^{\ast}\nbige)
\simeq
\nbigr_{q,p\ast}(\nbige_1)\otimes\nbige$
induces 
$\nbigr_{q,p\ast}(\nbigp_{\ast}\nbige_1
 \otimes\nbigr_{q,p}^{\ast}\nbigp_{\ast}\nbige)
\simeq
 \nbigr_{q,p\ast}(\nbigp_{\ast}\nbige_1)\otimes\nbigp_{\ast}\nbige$.
The natural isomorphism
 $\nbigr_{q,p\ast}(
 \nhom(\nbigr_{q,p}^{\ast}\nbige,\nbige_1))
\simeq
  \nhom(\nbige,\nbigr_{q,p\ast}\nbige_1)$
induces 
 $\nbigr_{q,p\ast}(
 \nhom(\nbigr_{q,p}^{\ast}\nbigp_{\ast}\nbige,\nbigp_{\ast}\nbige_1))
\simeq
  \nhom(\nbigp_{\ast}\nbige,\nbigr_{q,p\ast}\nbigp_{\ast}\nbige_1)$.
\hfill\qed
\end{lem}

\subsubsection{Descent}
If $\nbigp_{\ast}\nbige_1$ is $\Gal_{q,p}$-equivariant,
then we obtain the induced filtered bundle $\nbigp_{\ast}\nbige_2$
over the descent $\nbige_2$ of $\nbige_1$
by taking the $\Gal_{q,p}$-invariant part of
$\nbigp_{\ast}\bigl(
\nbigr_{q,p\ast}(\nbige_1)_{|\Hhat_{q,\infty}\langle t\rangle}\bigr)$
$(t\in S^1_T)$
as in \S\ref{subsection;20.7.20.20}. \index{descent (filtered bundle)}
We obtain the following lemma from Lemma \ref{lem;21.8.26.4}.
\begin{lem}
$\nbigr_{q,p}^{\ast}(\nbigp_{\ast}\nbige_2)$
is naturally isomorphic to
$\nbigp_{\ast}\nbige_1$.
\hfill\qed
\end{lem}

We obtain the following lemma from Lemma \ref{lem;21.8.25.41}.
\begin{lem}
Let $\nbige_1'$ be a $\Gal_{q,p}$-equivariant locally free
$\nbigo_{\Hhat_{\infty,p}}(\ast H_{\infty,p})$-submodule
of $\nbige_1$.
We set $\nbige''_1:=\nbige_1/\nbige'_1$. 
We obtain the locally free
$\nbigo_{\Hhat_{\infty,q}}(\ast H_{\infty,q})$-submodule
$\nbige_2'\subset\nbige_2$
as the descent of $\nbige'_1$.
We also obtain the locally free 
$\nbigo_{\Hhat_{\infty,q}}(\ast H_{\infty,q})$-quotient module
$\nbige_2\lrarr\nbige_2''$
as the descent of $\nbige''_1$.
Then, the decent $\nbigp_{\ast}\nbige'_2$ of $\nbigp_{\ast}\nbige'_1$
is equal to the filtered bundle over $\nbige_2'$
induced by $\nbigp_{\ast}\nbige_2$ with the inclusion
$\nbige'_2\subset\nbige_2$.
Similarly, 
the decent $\nbigp_{\ast}\nbige''_2$ of $\nbigp_{\ast}\nbige''_1$
is equal to the filtered bundle over $\nbige_2''$
induced by $\nbigp_{\ast}\nbige_2$ with the projection
$\nbige_2\lrarr\nbige_2''$.
\hfill\qed
\end{lem}

We obtain the following lemma from Lemma \ref{lem;21.8.26.5}.
\begin{lem}
\label{lem;21.8.26.11}
For a filtered bundle $\nbigp_{\ast}\nbige$ on
$(\Hhat_{\infty,q},H_{\infty,q})$,
$\nbigr_{q,p}^{\ast}(\nbigp_{\ast}\nbige)$
is naturally $\Gal_{q,p}$-equivariant,
and
$\nbigp_{\ast}\nbige$ is isomorphic to the descent of
$\nbigr_{q,p}^{\ast}\nbigp_{\ast}\nbige$.
\hfill\qed 
\end{lem}

\subsection{Basic filtered objects with pure slope}
\label{subsection;17.10.7.1}

Let $q\in\seisuu_{>0}$,
$\ell\in\seisuu$ and $\alpha\in\cnum^{\ast}$.
We consider the $\seisuu$-action 
on $\LLhat^{\lambda\cov}_q(\ell,\alpha):=
 \nbigo_{\Hhat^{\cov}_{\infty,q}}(\ast H^{\cov}_{\infty,q})
 \,e_{q,\ell,\alpha}$
given by
\[
\kappa_1^{\ast}(e_{q,\ell,\alpha})=
 \alpha\,
 (1+|\lambda|^2)^{\ell/q}
 y_q^{-\ell}(1+2\sqrt{-1}\lambda T y^{-1})^{-\ell/q}
 \exp\Bigl(\frac{\ell}{q}G(2\sqrt{-1}T\lambda y^{-1})\Bigr)
 \,
 e_{q,\ell,\alpha}.
\]
\index{sheaf $\LLhat^{\lambda\cov}_q(\ell,\alpha)$}
Here, $G(x)=1-x^{-1}\log(1+x)$.
We obtain the induced 
$\nbigo_{\Hhat_{\infty,q}}(\ast H_{\infty,q})$-module,
which is denoted by
$\LLhat^{\lambda}_q(\ell,\alpha)$.
It has pure slope $\ell/q$.
\index{sheaf $\LLhat^{\lambda}_q(\ell,\alpha)$}

Let $a\in\real$.
Let $e^t_{q,\ell,\alpha}$ denote 
the restriction of $e_{q,\ell,\alpha}$
to $\{t\}\times\inftyhat_{y,q}$.
We define the filtrations 
\[
 \nbigp^{(a)}_{\ast}\Bigl(
 \LLhat^{\lambda\cov}_q(\ell,\alpha)
 _{|\{t\}\times\inftyhat_{y,q}}
 \Bigr)
\]
by setting 
\[
 \deg^{\nbigp}(e^t_{q,\ell,\alpha})=a-\frac{\ell t}{T}.
\]
Because the filtrations are preserved by
the $\seisuu$-action,
we obtain a tuple of the filtrations
$\nbigp^{(a)}_{\ast}\LLhat^{\lambda}_q(\ell,\alpha)
=\bigl(
\nbigp^{(a)}_{\ast}(\LLhat^{\lambda}_q(\ell,\alpha)
 _{|\Hhat_{\infty,q}\langle t\rangle})\,\big|\,
t\in S^1_T
\bigr)$.
\index{filtered bundle
$\nbigp^{(a)}_{\ast}\LLhat^{\lambda}_q(\ell,\alpha)$ }

\begin{rem}
 The filtered bundles in this section
 naturally appear in the context of periodic monopoles.
 See {\rm\S\ref{subsection;17.9.26.1}}.
\hfill\qed
\end{rem}

\subsection{Good filtered bundles
   over $\nbigo_{\Hhat_{\infty,q}}(\ast H_{\infty,q})$-modules
   with level $\leq 1$}

Let $\nbige$ be a locally free 
$\nbigo_{\Hhat_{\infty,q}}(\ast H_{\infty,q})$-module
of finite rank
with level $\leq 1$.
(See Definition \ref{df;20.7.31.10}.)
Note that any $\Phi^{\ast}$-invariant lattice $\nbigl$
of the difference $\cnum(\!(y_q^{-1})\!)$-module
$\nbign:=\nbige_{|\Hhat_{\infty,q}\langle 0\rangle}$
induces a lattice
$\Upsilon(\nbigl)$ of $\nbige$.
In particular,
it induces lattices
$\Upsilon(\nbigl)_{|\Hhat_{\infty,q}\langle t\rangle}$
of $\nbige_{|\Hhat_{\infty,q}\langle t\rangle}$ $(t\in S^1_T)$.

 \begin{df}
\label{df;20.7.31.11}
\index{good filtered bundle (level $\leq 1$)}
Suppose that the level of $\nbige$ is less than $1$.
A filtered bundle
$\nbigp_{\ast}(\nbige_{|\Hhat_{\infty,q}\langle t\rangle}\,|\,t\in S^1_T)$
is called good
if there exists a filtered bundle
 $\nbigp_{\ast}(\nbign)$ over $\nbign$
 which is good with respect to $\Phi^{\ast}$
 such that
 $\nbigp_{\ast}(\nbige_{|\Hhat_{\infty,q}\langle t\rangle})
 =\Upsilon(\nbigp_{\ast}\nbign)_{|\Hhat_{\infty,q}\langle t\rangle}$.
\hfill\qed
\end{df}

\begin{lem}
\label{lem;21.8.26.20}
Let $\nbigp_{\ast}\nbige$ be a good filtered bundle over $\nbige$.
Let $\nbige'$  be a locally free
$\nbigo_{\Hhat_{\infty,q}}(\ast H_{\infty,q})$-submodule
of $\nbige$.
We set $\nbige''=\nbige/\nbige'$.
They are equipped with the induced filtered bundles
$\nbigp_{\ast}\nbige'$ and $\nbigp_{\ast}\nbige''$,
respectively.
Then,  $\nbigp_{\ast}\nbige'$ and $\nbigp_{\ast}\nbige''$
are also good in the sense of
Definition {\rm\ref{df;20.7.31.11}}.
\end{lem}
\pf
We explain the proof for $\nbige'$.
The other case can be argued similarly.
Note that the level of $\nbige'$ is also less than $1$,
which follows from Lemma \ref{lem;21.8.26.12}
and Lemma \ref{lem;21.8.26.13}.
We set $\nbign'=\nbige'_{|\Hhat_{\infty,q}\langle 0\rangle}$,
which is a $\cnum(\!(y_q^{-1})\!)$-difference submodule
of $\nbign$,
which is equipped with the induced filtration $\nbigp_{\ast}\nbign'$.
By Lemma \ref{lem;21.8.26.14},
$\nbigp_{\ast}\nbign'$ is good.
It is easy to see that
$\nbigp_{\ast}(\nbige'_{|\Hhat_{\infty,q}\langle t\rangle})=
\Upsilon(\nbigp_{\ast}\nbign)_{|\Hhat_{\infty,q}\langle t\rangle}$.
Hence, $\nbigp_{\ast}\nbige'$ is good.
\hfill\qed

\vspace{.1in}

If moreover $\nbige$ is unramified,
there exists a decomposition
\begin{equation}
\label{eq;20.7.20.50}
 \nbige=\bigoplus_{\gminib\in S(q)}\nbige_{\gminib}
\end{equation}
such that
for a lattice $\nbige_{\gminib,0}$ of $\nbige_{\gminib}$,
the monodromy of $\Loc_q(\nbige_{\gminib,0})$
is equal to the multiplication of $1+\gminib$.
The following lemma is easy to see.
  \begin{lem}
   Suppose that the level of $\nbige$ is less than $1$
   and that $\nbige$ is unramified.
Then, a filtered bundle
 $\nbigp_{\ast}(\nbige)
 =\bigl(
  \nbigp_{\ast}(\nbige_{|\Hhat_{\infty,q}\langle t\rangle})\,\big|\,
 t\in S^1_T
 \bigr)$
 is good if and only if the following condition is satisfied.
\begin{itemize}
 \item For each $a\in\real$, there exists a lattice
       $\nbigp_a(\nbige)$ of $\nbige$
       such that
       $\nbigp_a(\nbige_{|\Hhat_{\infty,q}\langle t\rangle})
       =\nbigp_a(\nbige)_{|\Hhat_{\infty,q}\langle t\rangle}$.
 \item The lattice $\nbigp_a(\nbige)$ is compatible
       with the decomposition {\rm(\ref{eq;20.7.20.50})},
       i.e., by setting
       $\nbigp_a(\nbige_{\gminib})=\nbigp_a\nbige\cap\nbige_{\gminib}$,
       we obtain
       $\nbigp_a(\nbige)
       =\bigoplus\nbigp_a\nbige_{\gminib}$.
 \item The monodromy of $\Loc_q(\nbigp_a\nbige_{\gminib})$
       are equal to the multiplication of $1+\gminib$.
\hfill\qed
\end{itemize}
  \end{lem}

\subsection{Good filtered bundles
over $\nbigo_{\Hhat_{\infty,q}}(\ast H_{\infty,q})$-modules}
\label{subsection;21.8.16.1}

Let $\nbige$ be a locally free 
$\nbigo_{\Hhat_{\infty,q}}(\ast H_{\infty,q})$-module
of finite rank.
\begin{df}
\index{unramified modulo level $\leq 1$}
We say that $\nbige$ is unramified modulo level $\leq 1$
 if there exists a decomposition
\begin{equation}
\label{eq;20.7.20.60}
 \nbige=\bigoplus_{\ell\in\seisuu}
 \bigoplus_{\alpha\in\cnum^{\ast}}
 \nbige_{\ell,\alpha}
\end{equation}
such that
each $\nbige_{\ell,\alpha}$
is the tensor product of
$\LLhat^{\lambda}_q(\ell,\alpha)$
and
a locally free
$\nbigo_{\Hhat_{\infty,q}}(\ast H_{\infty,q})$-module
 of level $\leq 1$.
 Note that the decomposition
 {\rm(\ref{eq;20.7.20.60})}
 is uniquely determined.
\hfill\qed
 \end{df}

\begin{df}
 \label{df;20.7.31.12}
 \index{good filtered bundle (unramified modulo level $\leq 1$)}
 If $\nbige$ is unramified modulo level $\leq 1$,
 a filtered bundle $\nbigp_{\ast}\nbige$
 is called good if the following conditions are satisfied.
\begin{itemize}
 \item $\nbigp_{\ast}\nbige$ is compatible
       with the decomposition {\rm(\ref{eq;20.7.20.60})},
       i.e.,
       $\nbigp_{\ast}(\nbige)
       =\bigoplus_{\ell,\alpha}
       \nbigp_{\ast}(\nbige_{\ell,\alpha})$.
 \item There exists a good filtered bundle
       $\nbigp_{\ast}\bigl(
       \LLhat^{\lambda}_q(\ell,\alpha)^{-1}
       \otimes
       \nbige_{\ell,\alpha}
       \bigr)$
       over
       $\LLhat^{\lambda}_q(\ell,\alpha)^{-1}
       \otimes
       \nbige_{\ell,\alpha}$
       (Definition {\rm\ref{df;20.7.31.11}})
       such that
\[
\nbigp_{\ast}(\nbige_{\ell,\alpha|\Hhat_{\infty,q}\langle t\rangle})
=\nbigp^{(0)}_{\ast}\bigl
(\LLhat^{\lambda}_q(\ell,\alpha)_{|\Hhat_{\infty,q}\langle t\rangle}
\bigr)
\otimes\nbigp_{\ast}\bigl(
\LLhat^{\lambda}_q(\ell,\alpha)^{-1}
\otimes
\nbige_{\ell,\alpha}
\bigr)_{|\Hhat_{\infty,q}\langle t\rangle}.
\]
\end{itemize}
\hfill\qed
\end{df}

Recall that for any locally free
$\nbigo_{\Hhat_{\infty,q}}(\ast H_{\infty,q})$-module
$\nbige$ of finite rank
there exists $p\in q\seisuu_{>0}$
such that
$\nbigr_{q,p}^{\ast}\nbige$
is unramified modulo level $\leq 1$.
(See Proposition \ref{prop;20.7.19.22}
and Theorem \ref{thm;20.7.19.41}.)

\begin{df}
\label{df;20.7.28.20}
\index{good filtered bundle}
 For a locally free
 $\nbigo_{\Hhat_{\infty,q}}(\ast H_{\infty,q})$-module
 $\nbige$ of finite rank,
 a filtered bundle
 $\nbigp_{\ast}\nbige$ over $\nbige$ is good
 if there exist $p\in q\seisuu_{>0}$
 such that
 (i)  $\nbigr_{q,p}^{\ast}(\nbige)$
 is unramified modulo level $\leq 1$,
 (ii) $\nbigr_{q,p}^{\ast}(\nbigp_{\ast}\nbige)$
 is a good filtered bundle
 over $\nbigr_{q,p}^{\ast}(\nbige)$
 (Definition {\rm\ref{df;20.7.31.12}}).
\hfill\qed
\end{df}

\subsubsection{An equivalence}

For a good filtered bundle $\nbigp_{\ast}\nbige$
over 
a locally free $\nbigo_{\Hhat_{\infty,q}}(\ast H_{\infty,q})$-module
$\nbige$,
by taking the restriction to
$\Hhat_{\infty,q}\langle 0\rangle$,
we obtain a filtered bundle
$\nbigp_{\ast}(\nbige_{|\Hhat_{\infty,q}\langle 0\rangle})$
over $\nbige_{|\Hhat_{\infty,q}\langle 0\rangle}$.
By the construction, it is good for
the difference $\cnum(\!(y_q^{-1})\!)$-module
$\nbige_{|\Hhat_{\infty,q}\langle 0\rangle}$.
The following proposition is obvious
by the definitions.
 \begin{prop}
\label{prop;21.9.17.41}
  By the restriction from
  $\nbigp_{\ast}\nbige$
 to $\nbigp_{\ast}(\nbige_{|\Hhat_{\infty,q}\langle 0\rangle})$,
we obtain an equivalence between
 good filtered bundles on
 locally free $\nbigo_{\Hhat_{\infty,q}}(\ast H_{\infty,q})$-modules
 and good filtered
  difference $\cnum(\!(y_q^{-1})\!)$-modules.
 \hfill\qed
 \end{prop}

\subsubsection{Some properties}

The following lemma is obvious.
\begin{lem}
Let $\nbigp_{\ast}\nbige$
be a good filtered bundle
on $(\Hhat_{\infty,q},H_{\infty,q})$.
Then, it is compatible with the slope decomposition
{\rm(\ref{eq;20.7.20.100})},
i.e.,
$\nbigp_{\ast}\nbige
=\bigoplus \nbigp_{\ast}\nbigs_{\omega}\nbige$.
\hfill\qed
\end{lem}

\begin{lem}
\label{lem;21.8.26.40}
 Let $\nbigp_{\ast}\nbige$ be a good filtered bundle
on $(\Hhat_{\infty,q},H_{\infty,q})$.
Let $\nbige'$ be a locally free
$\nbigo_{\Hhat_{\infty,q}}(\ast H_{\infty,q})$-submodule
of $\nbige$. 
We set $\nbige''=\nbige/\nbige'$.
 Then, the induced filtered bundle $\nbigp_{\ast}\nbige'$
(resp. $\nbigp_{\ast}\nbige''$)
over $\nbige'$ (resp. $\nbige''$) is good. 
\end{lem}
\pf
By Lemma \ref{lem;21.8.26.10} and Lemma \ref{lem;21.8.26.11},
it is enough to study the case where
$\nbige$ and $\nbige'$ are unramified modulo $\leq 1$.
Because $\nbigp_{\ast}\nbige$ is compatible with
the decomposition (\ref{eq;20.7.20.60}),
it is enough to study the case where
the level of $\nbige$ is less than $1$,
which we have already studied in Lemma \ref{lem;21.8.26.20}.
\hfill\qed

\begin{lem}
\label{lem;21.8.26.50}
Let $\nbigp_{\ast}\nbige_i$ $(i=1,2)$
be good filtered bundles over
locally free
$\nbigo_{\Hhat_{\infty,q}}(\ast H_{\infty,q})$-modules $\nbige_i$.
Then,
$\nbigp_{\ast}(\nbige_1)\oplus\nbigp_{\ast}(\nbige_2)$,
$\nbigp_{\ast}(\nbige_1)\otimes\nbigp_{\ast}(\nbige_2)$,
and
$\nhom(\nbigp_{\ast}\nbige_1,\nbigp_{\ast}\nbige_2)$
are also good filtered bundles. 
In particular,
for a good filtered bundle $\nbigp_{\ast}\nbige$,
its dual
$(\nbigp_{\ast}\nbige)^{\lor}$
is also a good filtered bundle. 
\end{lem}
\pf
By Lemma \ref{lem;21.8.26.30} and Lemma \ref{lem;21.8.26.11},
it is enough to consider the case where
$\nbige_i$ are unramified modulo level $\leq 1$.
Because the good filtered bundles $\nbigp_{\ast}\nbige_i$
are compatible with the decompositions as in (\ref{eq;20.7.20.60}),
we have only to consider the case where
the levels of $\nbige_i$ are less than $1$.
Then, it follows from the claims in the case of
$\cnum(\!(y_q^{-1})\!)$-difference modules
(Lemma \ref{lem;21.8.26.31}).
\hfill\qed

\subsection{Global lattices on the covering space}
\label{subsection;20.7.22.31}

Let $\nbigp_{\ast}\nbige$ be a good filtered bundle
over a locally free
$\nbigo_{\Hhat_{\infty,q}}(\ast H_{\infty,q})$-module
$\nbige$.
We set $\nbige^{\cov}:=\varpi_q^{-1}(\nbige)$.
We obtain the induced filtrations
$\nbigp_{\ast}\bigl(
\nbige^{\cov}_{|\{t\}\times\inftyhat_{y,q}}
\bigr)$
$(t\in H_{\infty,q}^{\cov})$.
The slope decomposition
$\nbige=\bigoplus \nbigs_{\omega}(\nbige)$
induces a decomposition
$\nbige^{\cov}
=\bigoplus \nbigs_{\omega}(\nbige)^{\cov}$.
For each $t\in H^{\cov}_{\infty,q}$,
we obtain the decomposition
$\nbigp_{\ast}(\nbige^{\cov}_{|\{t\}\times\inftyhat_{y,q}})
=\bigoplus
\nbigp_{\ast}(\nbigs_{\omega}(\nbige)^{\cov}
 _{|\{t\}\times\inftyhat_{y,q}})$.

\begin{lem}
\label{lem;20.7.20.110}
 For any $t_1,t_2\in H^{\cov}_{\infty,q}$,
the isomorphism
 $\nbigs_{\omega}(\nbige)^{\cov}_{|\{t_1\}\times\inftyhat_{y,q}}
 \simeq
 \nbigs_{\omega}(\nbige)^{\cov}_{|\{t_2\}\times\inftyhat_{y,q}}$
 induces the following isomorphisms for any $a\in\real$:
\begin{equation}
\label{eq;20.7.20.103}
 \nbigp_{a-t_1q\omega/T}
  \nbigs_{\omega}(\nbige)^{\cov}_{|\{t_1\}\times\inftyhat_{y,q}}
 \simeq
 \nbigp_{a-t_2q\omega/T}
 \nbigs_{\omega}(\nbige)^{\cov}_{|\{t_2\}\times\inftyhat_{y,q}}.
\end{equation}
\end{lem}
\pf
If $\nbige$ is unramified,
the claim is clear by the definition of good filtrations.
In general,
there exists $p\in q\seisuu_{>0}$
such that
$\nbigr_{q,p}^{\ast}(\nbige)$ is unramified.
We obtain the isomorphism
\[
 \nbigp_{a-t_1p\omega/T}
  \nbigs_{\omega}(\nbigr_{q,p}^{\ast}\nbige)^{\cov}
   _{|\{t_1\}\times\inftyhat_{y,q}}
 \simeq
 \nbigp_{a-t_2p\omega/T}
 \nbigs_{\omega}(\nbigr_{q,p}^{\ast}\nbige)^{\cov}
  _{|\{t_2\}\times\inftyhat_{y,q}}.
\]
It is compatible with the natural $\Gal_{q,p}$-action.
By using the natural isomorphisms (\ref{eq;20.7.20.102}),
we obtain (\ref{eq;20.7.20.103}).
\hfill\qed

\vspace{.1in}

By Lemma \ref{lem;20.7.20.110},
for each $a\in\real$,
there uniquely exists
a locally free $\nbigo_{\Hhat^{\cov}_{\infty,q}}$-submodule
$\gbigp_a(\nbige^{\cov})$
of $\nbige^{\cov}$
satisfying 
\[
 \gbigp_a(\nbige^{\cov})_{|\{t\}\times\inftyhat_{y,q}}
 :=\bigoplus_{\omega\in\rnum}
 \nbigp_{a-tq\omega/T}
 (\nbigs_{\omega}\nbige^{\cov}_{|\{t\}\times\inftyhat_{y,q}}).
\]
Thus, we obtain
the global filtration $\gbigp_{\ast}$ of $\nbige^{\cov}$.
The following lemma is obvious by the construction.
\index{global filtration $\gbigp_{\ast}$}

\begin{lem}
\label{lem;20.7.20.111}
 Let $(\nbigp_{\ast}\nbign,\Phi^{\ast})$
 be a good filtered difference $\cnum(\!(y_q^{-1})\!)$-module
 obtained as the restriction of
 $\nbigp_{\ast}\nbige$
 to $\Hhat_{\infty,q}\langle 0\rangle$.
 Let $p_2:\Hhat^{\cov}_{\infty,q}\lrarr \inftyhat_{y,q}$
 denote the projection.
 Then,
 under the natural isomorphism
 $\nbige^{\cov}\simeq p_2^{-1}(\nbign)$,
 we obtain
 $\gbigp_{\ast}(\nbige^{\cov})
 =p_2^{-1}(\nbigp_{\ast}\nbign)$. 
\hfill\qed
\end{lem}

We set
$\Gr^{\gbigp}_a(\nbige^{\cov}):=
\gbigp_a(\nbige^{\cov})\big/
\gbigp_{<a}(\nbige^{\cov})$.
We obtain the decomposition
$\Gr^{\gbigp}_a(\nbige^{\cov})=
\bigoplus
\Gr^{\gbigp}_a(\nbigs_{\omega}\nbige^{\cov})$
induced by the slope decomposition.
\index{local system $\Gr^{\gbigp}_a(\nbige^{\cov})$}
By the construction, there exist the natural isomorphisms:
\[
 \Gr^{\gbigp}_a(\nbigs_{\omega}\nbige^{\cov})_{\{t\}\times\inftyhat_{y,q}}
 \simeq
 \Gr^{\nbigp}_{a-tq\omega/T}
 (\nbigs_{\omega}\nbige^{\cov}_{|\{t\}\times\inftyhat_{y,q}}).
\]
By Lemma \ref{lem;20.7.20.111}
and the construction in \S\ref{subsection;20.7.20.112},
there exist the endomorphism
$\gbigf_{\nbige}$
of $\Gr^{\gbigp}_a(\nbige^{\cov})$
induced by
$\Res(\Phi^{\ast})$
for the corresponding
good filtered difference $\cnum(\!(y_q^{-1})\!)$-module.
We also obtain 
the monodromy weight filtration $W$
of the nilpotent part of $\gbigf_{\nbige}$
on each $\Gr^{\gbigp}_{a}(\nbige^{\cov})$.
They are compatible with the slope decomposition.

\subsection{Local lattices}
\label{subsection;17.10.10.1}

Let $\nbigp_{\ast}\nbige=
 \bigl(\nbigp_{\ast}\nbige_{|\Hhat_{\infty,q}\langle t\rangle}
 \,\big|\,t\in S^1_T\bigr)$
denote a good filtered bundle
over a locally free 
$\nbigo_{\Hhat_{\infty,q}}(\ast H_{\infty,q})$-module
$\nbige$.
For any $t\in S^1_T$,
we set
\index{set $\Par(\nbigp_{\ast}\nbige,t)$}
\[
 \Par(\nbigp_{\ast}\nbige,t):=
 \bigl\{
 b\in\real\,\big|\,
 \Gr^{\nbigp}_b(\nbige_{|\Hhat_{\infty,q}\langle t\rangle})\neq 0
 \bigr\}.
\]
We clearly have
$\Par(\nbigp_{\ast}\nbige,t)
=\bigcup_{\omega\in\rnum}
 \Par\bigl(\nbigp_{\ast}(\nbigs_{\omega}\nbige),t\bigr)$.
The following lemma is clear by the construction.
\begin{lem}
\label{lem;17.10.8.20}
For $t_0,t_1\in S^1_T$,
there exists the following relation:
\[
 \Par(\nbigp_{\ast}\nbigs_{\omega}\nbige,t_1)
=\bigl\{
 b-q\omega T^{-1}(\ttilde_1-\ttilde_0)
 \,\big|\,
 b\in\Par(\nbigp_{\ast}\nbigs_{\omega}\nbige,t_0)
 \bigr\}.
\]
 Here, $\ttilde_i\in\real$ are any lifts of
 $t_i\in S^1_T=\real/T\seisuu$.
\hfill\qed
\end{lem}

Take $t_0\in S^1_T\simeq H_{\infty,q}$.
For $0<\epsilon<T/10$,
we set
\index{space $H_{\infty,q}\langle t_0,\epsilon\rangle$}
\[
H_{\infty,q}\langle t_0,\epsilon\rangle:=
\openopen{t_0-\epsilon}{t_0+\epsilon}
\subset H_{\infty,q},
\]
and let
$\Hhat_{\infty,q}\langle t_0,\epsilon\rangle$
denote the ringed space obtained as
$H_{\infty,q}\langle t_0,\epsilon\rangle$
with
$\nbigo_{\Hhat_{\infty,q}\langle t_0,\epsilon\rangle}:=
\nbigo_{\Hhat_{\infty,q}|\Hhat_{\infty,q}\langle t_0,\epsilon\rangle}$.
\index{ringed space $\\Hhat_{\infty,q}\langle t_0,\epsilon\rangle$}
We set
$\nbigo_{\Hhat_{\infty,q}\langle t_0,\epsilon\rangle}
 (\ast H_{\infty,q}\langle t_0,\epsilon\rangle)
 =\nbigo_{\Hhat_{\infty,q}}(\ast H_{\infty,q})
 _{|\Hhat_{\infty,q}\langle t_0,\epsilon\rangle}$.
The following lemma is clear
by the condition of good filtered bundles.
\begin{lem}
For each $a\in\real$,
there uniquely exists a locally free 
$\nbigo_{\Hhat_{\infty,q}\langle t_0,\epsilon\rangle}$-submodule
$\vecP^{(t_0)}_a(\nbige_{|\Hhat_{\infty,q}\langle t_0,\epsilon\rangle})
 \subset\nbige_{|\Hhat_{\infty,q}\langle t_0,\epsilon\rangle}$
 satisfying the following condition
 for any $t\in H_{\infty,q}\langle t_0,\epsilon\rangle$:
\index{lattice $\vecP^{(t_0)}_a
(\nbige_{|\Hhat_{\infty,q}\langle t_0,\epsilon\rangle})$}
\[
\vecP^{(t_0)}_a\bigl(
\nbige_{|\Hhat_{\infty,q}\langle t_0,\epsilon\rangle}
\bigr)_{|\Hhat_{\infty,q}\langle t\rangle}=
\bigoplus_{\omega}
\nbigp_{a-q\omega(\ttilde-\ttilde_0)/T}
(\nbigs_{\omega}\nbige_{|\Hhat_{\infty,q}\langle t\rangle}).
\]
Here, $\ttilde_0\in\real$ denotes any lift of $t_0$,
and $\ttilde$ denotes the lift of $t$
such that $|\ttilde-\ttilde_0|<\epsilon$. 
It satisfies the following.
\begin{itemize}
\item
$\vecP^{(t_0)}_a(\nbige_{|\Hhat_{\infty,q}\langle t_0,\epsilon\rangle})
=\bigoplus_{\omega}
     \vecP^{(t_0)}_a(\nbigs_{\omega}
      \nbige_{|\Hhat_{\infty,q}\langle t_0,\epsilon\rangle})$.
\item
     $\vecP^{(t_0)}_a(\nbige_{|\Hhat_{\infty,q}\langle
     t_0,\epsilon\rangle})
     \otimes_{\nbigo_{\Hhat_{\infty,q}\langle t_0,\epsilon\rangle}}
     \nbigo_{\Hhat_{\infty,q}\langle
     t_0,\epsilon\rangle}(\ast H_{\infty,q}\langle t_0,\epsilon\rangle)
=\nbige_{|\Hhat_{\infty,q}\langle t_0,\epsilon\rangle}$.
\hfill\qed
\end{itemize}
 \end{lem}

We set
$\vecP^{(t_0)}_{<a}(\nbige_{|\Hhat_{\infty,q}\langle t_0,\epsilon\rangle})
=\sum_{b<a}
\vecP^{(t_0)}_{b}(\nbige_{|\Hhat_{\infty,q}\langle t_0,\epsilon\rangle})$.
We obtain the following locally constant sheaves
on $H_{\infty,q}\langle t_0,\epsilon\rangle$:
\index{local system $ \vecG^{(t_0)}_a(\nbige_{|\Ihat_q(t_0,\epsilon)})$}
\[
 \vecG^{(t_0)}_a(\nbige_{|\Hhat_{\infty,q}\langle t_0,\epsilon\rangle}):=
 \vecP^{(t_0)}_{a}(\nbige_{|\Hhat_{\infty,q}\langle t_0,\epsilon\rangle})
 \Big/
 \vecP^{(t_0)}_{<a}(\nbige_{|\Hhat_{\infty,q}\langle t_0,\epsilon\rangle}).
\]
There exists the decomposition:
\[
 \vecG^{(t_0)}_{\ast}(\nbige_{|\Hhat_{\infty,q}\langle t_0,\epsilon\rangle})
=\bigoplus_{\omega\in\rnum}
 \vecG^{(t_0)}_{\ast}
 (\nbigs_{\omega}\nbige_{|\Hhat_{\infty,q}\langle t_0,\epsilon\rangle}).
\]
The following lemma is clear by the construction.
\begin{lem}
Let $\ttilde_0\in\real$ which is a lift of $t_0$.
Let $\Hhat_{\infty,q}^{\cov}\langle \ttilde_0,\epsilon\rangle:=
 \openopen{\ttilde_0-\epsilon}{\ttilde_0+\epsilon}\times\inftyhat_{y,q}$.
For $\omega\in\rnum$,
there exists a natural isomorphism
\[
 \vecG^{(t_0)}_b(
 \nbigs_{\omega}\nbige_{|\Hhat_{\infty,q}\langle t_0,\epsilon\rangle})
\simeq
 \Gr^{\gbigp}_{b+\omega q\ttilde_0/T}
 (\nbigs_{\omega}\nbige^{\cov})
  _{|\Hhat_{\infty,q}^{\cov}\langle\ttilde_0,\epsilon\rangle}.
\]
As a result, there exist
the nilpotent endomorphism
$\gbign$ and the monodromy weight filtration $W$
on 
$\vecG^{(t_0)}_b(\nbige_{|\Hhat_{\infty,q}\langle t_0,\epsilon\rangle})$. 
\hfill\qed
\end{lem}

We set
$\nu(\nbige):=\max\bigl\{|\omega_1|+|\omega_2|\,\big|\,
\nbigs_{\omega_i}(\nbige)\neq 0
\bigr\}$.\index{number $\nu(\nbige)$}
We assume that $\epsilon$ satisfies
\begin{equation}
\label{eq;20.7.21.11}
 10 q\nu(\nbige)\epsilon/T
 <
 \min\bigl\{|a-b|\,\big|\,
 a,b\in\Par(\nbigp_{\ast}\nbige,t_0),\,a\neq b
 \bigr\}.
\end{equation}
Let $v$ be a non-zero section of
$\vecP^{(t_0)}_a(\nbige_{|\Hhat_{\infty,q}\langle t_0,\epsilon\rangle})$.
According to the decomposition
\[
 \vecP^{(t_0)}_a(\nbige_{|\Hhat_{\infty,q}\langle t_0,\epsilon\rangle})
=\bigoplus
 \vecP^{(t_0)}_a\bigl(
 \nbigs_{\omega}\nbige_{|\Hhat_{\infty,q}\langle t_0,\epsilon\rangle}
 \bigr),
\]
we obtain the decomposition
$v=\sum v_{\omega}$.
For any $v_{\omega}$
and for any $t\in H_{\infty,q}\langle t_0,\epsilon\rangle$,
we set
\[
 \deg^{\nbigp}(v_{\omega|\Hhat_{\infty,q}\langle t\rangle})
 :=
 \left\{
\begin{array}{ll}
 \inf\bigl\{
 b\in\real\,\big|\,
 v_{\omega}\in \nbigp_b(\nbige_{|\Hhat_{\infty,q}\langle t\rangle})
 \bigr\}
 & (v_{\omega}\neq 0)\\
 -\infty & (v_{\omega}=0).
\end{array}
 \right.
\]

\begin{lem}
\label{lem;20.7.21.31}
Suppose that there exists $\omega_0$
such that
\[
  \deg^{\nbigp}(v_{\omega|\Hhat_{\infty,q}\langle t_0\rangle})
 <\deg^{\nbigp}(v_{\omega_0|\Hhat_{\infty,q}\langle t_0\rangle})
\]
for any $\omega\neq \omega_0$.
Then, we obtain
\[
 \deg^{\nbigp}(v_{\omega|\Hhat_{\infty,q}\langle t\rangle})
 <\deg^{\nbigp}(v_{\omega_0|\Hhat_{\infty,q}\langle t\rangle})
\]
for $t\in H_{\infty,q}\langle t_0,\epsilon\rangle$
and any $\omega\neq \omega_0$.
\end{lem}
\pf
It follows from Lemma \ref{lem;17.10.8.20}
and (\ref{eq;20.7.21.11}).
\hfill\qed

\subsection{Complement for good filtered bundles
with  level $\leq 1$}

Let $\nbige$ be a locally free
$\nbigo_{\Hhat_{\infty,q}}(\ast H_{\infty,q})$-module
of level $\leq 1$.
Let
$\nbigp_{\ast}\nbige=\bigl(
\nbigp_{\ast}(\nbige_{|\Hhat_{\infty,q}\langle t\rangle}) \,\big|\,
t\in S^1_T
\bigr)$
be a good filtered bundle over $\nbige$.
In this case,
for each $a\in\real$,
there exist $\nbigo_{\Hhat_{\infty,q}}$-lattices
$\nbigp_a\nbige\subset\nbige$
such that
$\nbigp_a(\nbige)_{|\Hhat_{\infty,q}\langle t\rangle}
=\nbigp_a\bigl(\nbige_{|\Hhat_{\infty,q}\langle t\rangle}
\bigr)$.
The other filtrations are also simply described
in terms of $\nbigp_a\nbige$.
Indeed,
we obtain
$\gbigp_a(\nbige^{\cov})
=\varpi_q^{-1}(\nbigp_a\nbige)$.
For any $t_0\in S^1_T$,
we also obtain
$\vecP^{(t_0)}_a(\nbige_{|\Hhat_{\infty,q}\langle t_0,\epsilon\rangle})
=\nbigp_a(\nbige)_{|\Hhat_{\infty,q}\langle t_0,\epsilon\rangle}$.

\section{Formal difference modules
of level $\leq 1$ and formal $\lambda$-connections}
\label{subsection;17.10.28.20}

Let $q\in\seisuu_{\geq 1}$.
Let us explain that
$(2\sqrt{-1}\lambda T)$-difference
$\cnum(\!(y_q^{-1})\!)$-modules
of level $\leq 1$
are equivalent to
formal $\lambda$-connections
whose Poincar\'{e} rank are strictly smaller than $q$.

\subsection{Formal $\lambda$-connections}

We introduce a formal variable $x$,
and we fix a $q$-th root $x_q$ of $x$.
\index{variable $x_q$}
Let $\nbigv$ be a $\cnum(\!(x_q^{-1})\!)$-module
with a $\lambda$-connection $\nabla^{\lambda}$.
A $\cnum[\![x_q^{-1}]\!]$-lattice $\nbigl$ of $\nbigv$
is called $\nabla^{\lambda}_x$-invariant
if $\nabla^{\lambda}_{x}\nbigl\subset\nbigl$
is satisfied,
 where
 $\nabla^{\lambda}_x
 =q^{-1}x^{-1}x_q\nabla^{\lambda}_{x_q}
 =q^{-1}x_q^{-q+1}\nabla^{\lambda}_{x_q}$.
\index{$\nabla^{\lambda}_x$-invariant lattice}
For any $\nabla_x^{\lambda}$-invariant lattice $\nbigl$,
we obtain the endomorphism
$F_{\nabla^{\lambda}}$ of
$\nbigl_{|\infty_{x,q}}=\nbigl/x_q^{-1}\nbigl$
as follows;
for any $v\in \nbigl_{|\infty_{x,q}}$,
we take $\vtilde\in\nbigl$ which induces $v$,
then we put
$F_{\nabla^{\lambda}}(v)=\nabla^{\lambda}(\vtilde)_{|\infty_{x,q}}$.
\index{endomorphism $F_{\nabla^{\lambda}}$}

Recall that there exists $p\in q\seisuu_{>0}$
and a decomposition
\begin{equation}
 \varphi_{q,p}^{\ast}(\nbigv,\nabla^{\lambda})
 =\bigoplus_{\gminia\in x_p\cnum[x_p]}
 (\nbigv_{\gminia},
 \nabla^{\lambda}_{\gminia})
\end{equation}
such that
there exist $\cnum[\![x_q^{-1}]\!]$-lattices
$\nbigl_{\gminia}$ of $\nbigv_{\gminia}$
such that
$\nabla^{\lambda}_{\gminia}-d\gminia\id$
are logarithmic with respect to
$\nbigl_{\gminia}$.
(See \S\ref{subsection;20.7.20.20}
or \S\ref{subsection;20.7.21.1}
for $\varphi_{q,p}^{\ast}$.)
The set
$\nbigi(\nbigv,\nabla^{\lambda})=\{\gminia\,|\,\nbigv_{\gminia}\neq 0\}$
is uniquely determined.
If $\nbigi(\nbigv,\nabla^{\lambda})$
is contained in $x_q\cnum[x_q]$,
then $\nbigi(\nbigv,\nabla^{\lambda})$ is called unramified.
In that case, $(\nbigv,\nabla^{\lambda})$
is called unramified.
The Poincar\'{e} rank of
$(\nbigv,\nabla^{\lambda})$
is defined to be
\[
 \max\Bigl\{
 \frac{q}{p}\deg_{x_p}\gminia(x_p)\,\big|\,
 \gminia\in\nbigi(\nbigv,\nabla^{\lambda})
 \Bigr\}.
\]
We say that the $(\nbigv,\nabla^{\lambda})$ is regular
if the Poincar\'{e} rank of $(\nbigv,\nabla^{\lambda})$ is $0$.
It is equivalent to the existence of
a $\cnum[\![x_q^{-1}]\!]$-lattice
$\nbigl$ of $\nbigv$
for which $\nabla^{\lambda}$ is logarithmic,
i.e.,
$x\nabla^{\lambda}_x\nbigl\subset\nbigl$.
\index{Poincar\'{e} rank}
\index{regular}

The following lemma is standard.
 \begin{lem}
  Poincar\'{e} rank of $(\nbigv,\nabla^{\lambda})$
  is less than $q$
  if and only if there exists
  a $\nabla^{\lambda}_x$-invariant lattice
  of $(\nbigv,\nabla^{\lambda})$.
 If Poincar\'{e} rank of $(\nbigv,\nabla^{\lambda})$
 is strictly less than $q$,
 for any $\nabla^{\lambda}_x$-invariant lattice $\nbigl$,
 the induced endomorphism
 $F_{\nabla^{\lambda}}$ of $\nbigl_{|\infty_{x,q}}$
 is nilpotent.
 \hfill\qed
 \end{lem}

\subsection{Some sheaves of algebras on $\Hhat_{\infty,q}$}
\label{subsection;21.8.12.54}

For any open subset $U\subset H^{\cov}_{\infty,q}$,
let
$\nbigk_{\Hhat^{\cov}_{\infty,q}}(U)$
denote the space of formal power series with the variable
$y_q^{-1}$
over $C^{\infty}(U)$,
i.e.,
$\nbigk_{\Hhat^{\cov}_{\infty,q}}(U):=
\Bigl\{
 \sum_{j=0}^{\infty}a_{j}(t)
y_q^{-j}\,\Big|\,
 a_{j}(t)\in C^{\infty}(U)
 \Bigr\}$.
\index{sheaf $\nbigk_{\Hhat^{\cov}_{\infty,q}}$}
Similarly,
for any open subset $U\subset H^{\cov}_{\infty,q}$,
let 
$\nbigk_{\Hhat^{\cov}_{\infty,q}}(\ast H^{\cov}_{\infty,q})(U)$
denote the space of
formal Laurent power series with the variable $y_q^{-1}$
over $C^{\infty}(U)$.
\index{sheaf $\nbigk_{\Hhat^{\cov}_{\infty,q}}(\ast H^{\cov}_{\infty,q})$}
Thus, we obtain sheaves of algebras
$\nbigk_{\Hhat^{\cov}_{\infty,q}}$
and
$\nbigk_{\Hhat^{\cov}_{\infty,q}}(\ast H^{\cov}_{\infty,q})$
on $H^{\cov}_{\infty,q}$.
Note that
$\nbigk_{\Hhat^{\cov}_{\infty,q}}(\ast H^{\cov}_{\infty,q})$
is naturally isomorphic to
$\nbigo_{\Hhat^{\cov}_{\infty,q}}(\ast H_{\infty,q})
\otimes_{\nbigo_{\Hhat^{\cov}_{\infty}}}
\nbigk_{\Hhat^{\cov}_{\infty,q}}$.
We obtain the operators
$\del_t:\nbigk_{\Hhat^{\cov}_{\infty,q}}
\lrarr\nbigk_{\Hhat^{\cov}_{\infty,q}}$
and
$\del_t:\nbigk_{\Hhat^{\cov}_{\infty,q}}(\ast H^{\cov}_{\infty,q})
\lrarr\nbigk_{\Hhat^{\cov}_{\infty,q}}(\ast H^{\cov}_{\infty,q})$
defined by
$\del_t\sum a_j(t)y_q^{-j}
=\sum \del_ta_j(t)y_q^{-j}$.
The kernels are
$\nbigo_{\Hhat^{\cov}_{\infty,q}}$
and
$\nbigo_{\Hhat^{\cov}_{\infty,q}}(\ast H^{\cov}_{\infty,q})$,
respectively.

The sheaves
$\nbigk_{\Hhat^{\cov}_{\infty,q}}$
and
$\nbigk_{\Hhat^{\cov}_{\infty,q}}(\ast H^{\cov}_{\infty,q})$
are equivariant with respect to the $\seisuu$-action $\kappa$.
The operators
$\del_t$ on
$\nbigk_{\Hhat^{\cov}_{\infty,q}}$
and
$\nbigk_{\Hhat^{\cov}_{\infty,q}}(\ast H^{\cov}_{\infty,q})$
are also equivariant.

We obtain the sheaves
$\nbigk_{\Hhat_{\infty,q}}$
and
$\nbigk_{\Hhat_{\infty,q}}(\ast H_{\infty,q})$
on $H_{\infty,q}$
as the descents of
$\nbigk_{\Hhat^{\cov}_{\infty,q}}$ and
$\nbigk_{\Hhat^{\cov}_{\infty,q}}(\ast H^{\cov}_{\infty,q})$,
respectively.
\index{sheaf $\nbigk_{\Hhat_{\infty,q}}$}
\index{sheaf $\nbigk_{\Hhat_{\infty,q}}(\ast H_{\infty,q})$}
We obtain the induced operators $\del_t$
on 
$\nbigk_{\Hhat_{\infty,q}}$
and
$\nbigk_{\Hhat_{\infty,q}}(\ast H_{\infty,q})$.
The kernels are 
$\nbigo_{\Hhat_{\infty,q}}$ and
$\nbigo_{\Hhat_{\infty,q}}(\ast H_{\infty,q})$,
respectively.

\subsection{From formal $\lambda$-connections
  to formal difference modules}
\label{subsection;17.10.7.2}

Let $\Psi$ denote the map
$H_{\infty}\lrarr \{\infty\}$.
\index{map $\Psi$}
There exists the morphism of sheaves
$\Psi^{\ast}:
\Psi^{-1}\nbigo_{\inftyhat_{x}}
\lrarr
\nbigk_{\Hhat_{\infty}}$
induced by
$\Psi^{\ast}(x^{-1})=y^{-1}(1-2\sqrt{-1}\lambda t/y)^{-1}$.
We regard $(\Psi,\Psi^{\ast})$
as a ``$C^{\infty}$''-map
$\Psi:\Hhat_{\infty}\lrarr\inftyhat_x$.
\index{morphism $(\Psi,\Psi^{\ast})$}
Similarly,
we obtain the map
$\Psi_q:H_{\infty,q}\lrarr \{\infty_{x,q}\}$
and the morphism of sheaves
$\Psi_q^{\ast}:
 \Psi_q^{-1}\nbigo_{\inftyhat_{x,q}}
 \lrarr
 \nbigk_{\Hhat_{\infty,q}}$
 induced by
 $\Psi_q^{\ast}(x_q^{-1})
 =y_q^{-1}(1-2\sqrt{-1}\lambda t/y_q^q)^{-1/q}$.
We regard the pair
$(\Psi_q,\Psi_q^{\ast})$
as a $C^{\infty}$-map
$\Psi_q:\Hhat_{\infty,q}\lrarr \inftyhat_{x,q}$.
\index{morphism $(\Psi_q,\Psi_q^{\ast})$}
Let $\inftyhat_{x,q}\lrarr \inftyhat_x$
be the ramified covering given by $x_q\longmapsto x=x_q^q$.
Then, the following diagram is commutative:
\[
 \begin{CD}
  \Hhat_{\infty,q}
 @>>>
 \Hhat_{\infty}\\
 @V{\Psi_q}VV @V{\Psi}VV \\
 \inftyhat_{x,q}
 @>>>
 \inftyhat_{x}.
 \end{CD}
\]
\index{map $\Psi$}
\index{map $\Psi_q$}

Let $(\nbigv,\nabla^{\lambda})$
be a $\cnum(\!(x_q^{-1})\!)$-module
with a $\lambda$-connection.
We obtain a locally free
$\nbigk_{\Hhat_{\infty,q}}(\ast H_{\infty,q})$-module
\index{sheaf $\nbigvtilde^{\infty}$}
\begin{equation}
\label{eq;20.7.18.2}
 \nbigvtilde^{\infty}:=
\nbigk_{\Hhat_{\infty,q}}(\ast H_{\infty,q})
\otimes_{\Psi_q^{-1}\nbigo_{\inftyhat_{x,q}}(\ast \infty_{x,q})}
\Psi_q^{-1}(\nbigv)
\simeq
\nbigk_{\Hhat_{\infty,q}}
\otimes_{\Psi_q^{-1}\nbigo_{\inftyhat_{x,q}}}
\Psi_q^{-1}(\nbigv).
\end{equation}
\index{sheaf $\nbigvtilde^{\infty}$}
There uniquely exists an operator
$\del_{\nbigvtilde^{\infty},t}:
 \nbigvtilde^{\infty}\lrarr\nbigvtilde^{\infty}$
determined by the following condition
for any local sections
$f$ and $s$ of
$\nbigk_{\Hhat_{\infty,q}}$
and $\nbigv$:
\index{operator $\del_{\nbigvtilde^{\infty},t}$}
\begin{equation}
\label{eq;20.7.18.3}
\del_{\nbigvtilde^{\infty},t}\bigl(
f\cdot\Psi_q^{-1}(s)\bigr)
=\del_t(f)\Psi_q^{-1}(s)
+f\Psi_q^{-1}(-2\sqrt{-1}\nabla^{\lambda}_xs).
\end{equation}
Let $\Psi_q^{\ast}(\nbigv,\nabla^{\lambda})$
denote the $\nbigo_{\Hhat_{\infty,q}}(\ast H_{\infty,q})$-module
obtained as the kernel of
$\del_{\nbigvtilde^{\infty},t}:
\nbigvtilde^{\infty}\lrarr\nbigvtilde^{\infty}$.
\index{sheaf $\Psi_q^{\ast}(\nbigv,\nabla^{\lambda})$}

If $\nabla^{\lambda}_x$-invariant lattice $\nbigl$
of $\nbigv$ exists,
we obtain a $\nbigk_{\Hhat_{\infty,q}}$-submodule
\index{sheaf $\nbigltilde^{\infty}$}
\begin{equation}
\label{eq;21.8.20.1}
\nbigltilde^{\infty}:=
\nbigk_{\Hhat_{\infty,q}}
\otimes_{\Psi_q^{-1}\nbigo_{\inftyhat_{x,q}}}
\Psi_q^{-1}(\nbigl)
\subset\nbigvtilde^{\infty}
\end{equation} \index{sheaf $\nbigltilde^{\infty}$}
such that
$\del_{\nbigvtilde^{\infty},t}
 \nbigltilde^{\infty}
 \subset\nbigltilde^{\infty}$.
We set
$\Psi_q^{\ast}(\nbigl,\nabla^{\lambda}):=
\nbigltilde^{\infty}\cap
\Psi_q^{\ast}(\nbigv,\nabla^{\lambda})$
which is an $\nbigo_{\Hhat_{\infty,q}}$-module.
\index{sheaf $\Psi_q^{\ast}(\nbigl,\nabla^{\lambda})$}
 
\begin{prop}
Suppose that the Poincar\'{e} rank
of $(\nbigv,\nabla^{\lambda})$
 is less than $q$.
\begin{itemize}
 \item
$\Psi_q^{\ast}(\nbigv,\nabla^{\lambda})$
 is a locally free $\nbigo_{\Hhat_{\infty,q}}(\ast H_{\infty,q})$-module
 such  that
 (i) $\rank\Psi_q^{\ast}(\nbigv,\nabla^{\lambda})=\rank\nbigv$,
 (ii) 
 $\Psi_q^{\ast}(\nbigv,\nabla^{\lambda})$ has pure slope $0$.
      If the Poincar\'{e} rank of
      $(\nbigv,\nabla^{\lambda})$ is strictly less than $q$,
      then the level of 
      $\Psi_q^{\ast}(\nbigv,\nabla^{\lambda})$ is
      less than $1$.
     If $(\nbigv,\nabla^{\lambda})$ is regular,
      the level of 
      $\Psi_q^{\ast}(\nbigv,\nabla^{\lambda})$ is $0$.
 \item
 For any $\nabla^{\lambda}_x$-invariant lattice $\nbigl$ of $\nbigv$,
 $\Psi_q^{\ast}(\nbigl,\nabla^{\lambda})$
 is an $\nbigo_{\Hhat_{\infty,q}}$-lattice of
 $\Psi_q^{\ast}(\nbigv,\nabla^{\lambda})$.
      If the induced endomorphism
      $F_{\nabla^{\lambda}}$ of $\nbigl_{|\infty_{x,q}}$ is nilpotent,
      then the monodromy of
      $\Loc_1\bigl(\Psi_q^{\ast}(\nbigl,\nabla^{\lambda})\bigr)$
      is unipotent.      
      If $\nabla^{\lambda}$ is logarithmic with respect to $\nbigl$,
      then the monodromy of
      $\Loc_q\bigl(\Psi_q^{\ast}(\nbigl,\nabla^{\lambda})\bigr)$
      is the identity.
\end{itemize} 
\end{prop}
\pf
Let $\nbigl$ be a $\nabla^{\lambda}_x$-invariant lattice
of $(\nbigv,\nabla^{\lambda})$.
Set $r:=\rank(\nbigv)$.
Let $\vecv=(v_1,\ldots,v_{r})$
be a frame of $\nbigl$.
We obtain $A\in M_r(\cnum[\![x_q^{-1}]\!])$
determined by 
$\nabla^{\lambda}_x\vecv
=\vecv A$.
We obtain the induced frame
$\vecvtilde=(\vtilde_1,\ldots,\vtilde_r)
=\Psi_q^{-1}(\vecv)$
of $\nbigltilde^{\infty}$.
The action of $\del_{\nbigvtilde^{\infty},t}$
is expressed as follows
with respect to $\vecvtilde$:
\[
 \del_{\nbigvtilde^{\infty},t}\vecvtilde
=\vecvtilde\cdot\Bigl(
 -2\sqrt{-1}\Psi_q^{\ast}A
 \Bigr).
\]
Note that
for the expansion
$A=\sum_{j=0}^{\infty} A_jx_q^{-j}$,
where $A_j\in M_r(\cnum)$,
we obtain
\[
 \Psi_q^{\ast}(A)(t,y_q)
=\sum_{j=0}^{\infty} A_jy_q^{-j}(1-2\sqrt{-1}\lambda t/y_q^q)^{-j/q}.
\]

Let $P\in H_{\infty,q}$ be any point,
and let $U$ be a simply connected neighbourhood of $P$
in $H_{\infty,q}$.
We put
$\vecvtilde^{(1)}_P:=\vecv_{|U}\exp(2\sqrt{-1}A_0t)$.
Let
$\Atilde^{(1)}\in M_r(\nbigk_{\Hhat_{\infty,q}}(U))$
be determined by
$\del_{\nbigvtilde^{\infty},t}\vecvtilde^{(1)}_P
 =\vecvtilde_P^{(1)}\Atilde^{(1)}$.
It is easy to check that
$\Atilde^{(1)}\in y_q^{-1}
M_r(\nbigk_{\Hhat_{\infty,q}}(U))$.
\begin{lem}
\label{lem;20.7.17.2}
 There exist matrices
 $G^{(i)}\in M_r(C^{\infty}(U))$ $(i=2,3,\ldots)$
 such that the following condition is satisfied.
 \begin{itemize}
  \item Let $I_r\in M_r(\cnum)$ denote the identity matrix.
	We set
	$\vecvtilde_P^{(i)}=
	\vecvtilde_P^{(1)}(I_r-y_q^{-1}G^{(2)})
	\cdots (I_r-y_q^{-i+1}G^{(i)})$.
	Let $\Atilde^{(i)}\in M_r(\nbigk_{\Hhat_{\infty,q}}(U))$
	be determined by
	$\del_{\nbigvtilde^{\infty},t}\vecvtilde_P^{(i)}=
	\vecvtilde_P^{(i)}\Atilde^{(i)}$.
	Then, $\Atilde^{(i)}$ is contained in
	$y_q^{-i}M_r(\nbigk_{\Hhat_{\infty,q}}(U))$.
 \end{itemize}
\end{lem}
\pf
We construct $G^{(i)}$ inductively.
Suppose that $G^{(i)}$ $(i\leq \ell)$ are constructed.
For the expansion
$A^{(\ell)}=\sum_{j\geq \ell}A_j^{(\ell)}y_q^{-j}$,
there exists $G^{(\ell+1)}\in M_r(C^{\infty}(U))$
such that
$\del_tG^{(\ell+1)}=A^{(\ell)}_{\ell}$.
We can check that the condition is satisfied for
$i\leq \ell+1$.
\hfill\qed

\vspace{.1in}

By Lemma \ref{lem;20.7.17.2},
there exists a frame
$\vecvtilde^{(\infty)}_P$
of $\nbigltilde^{\infty}_{|U}$
such that $\del_{\nbigvtilde^{\infty},t}\vecvtilde^{(\infty)}_P=0$.
Hence,
$\Psi^{\ast}(\nbigl,\nabla^{\lambda})$ is a locally free
$\nbigo_{\Hhat_{\infty,q}}$-module
whose rank is equal to $\rank\nbigv$.
Note that
the conjugacy class of the monodromy of the local system
$\Loc_1(\Psi^{\ast}(\nbigl,\nabla^{\lambda}))$
is equal to the conjugacy class of
$\exp(-2\sqrt{-1}\lambda A_0T)$.
Hence, if $A_0$ is nilpotent,
the monodromy of
$\Loc_1(\Psi^{\ast}(\nbigl,\nabla^{\lambda}))$
is unipotent.
If $\nabla^{\lambda}$ is logarithmic
with respect to $\nbigl$,
we obtain $A_j=0$ for $j=0,\ldots,q-1$,
and hence the monodromy of
$\Loc_q(\Psi^{\ast}(\nbigl,\nabla^{\lambda}))$
is the identity.
Thus, we obtain the claims for 
$\Psi^{\ast}(\nbigl,\nabla^{\lambda})$.
The claims for
$\Psi^{\ast}(\nbigv,\nabla^{\lambda})$
immediately follow.
\hfill\qed

\begin{rem}
Suppose that
the Poincar\'{e} rank of $(\nbigv,\nabla^{\lambda})$
is less than $q$.
Let $\nbigl$ be a lattice such that
 $\nabla^{\lambda}_x\nbigl\subset\nbigl$.
 If any eigenvalues $\alpha$ of $F_{\nabla}$ on
 $\nbigl_{|\infty_{x,q}}$
 satisfy
 $\exp(2\sqrt{-1}T\alpha)=1$,
 then the level of $\Psi^{\ast}(\nbigv,\nabla^{\lambda})$
 is strictly less than $1$.
\hfill\qed
\end{rem}

Let $P\in H_{\infty,q}$.
By evaluating $C^{\infty}$-functions at $t$,
we obtain the morphisms of algebras
$\nbigk_{\Hhat_{\infty,q},P}
\lrarr
\nbigo_{\Hhat_{\infty,q},P}$
and
$\nbigk_{\Hhat_{\infty,q}}(\ast H_{\infty,q})_P
\lrarr
\nbigo_{\Hhat_{\infty,q}}(\ast H_{\infty,q})_P$.
The natural morphism
$\Psi_q^{\ast}(\nbigv,\nabla^{\lambda})
\lrarr
 \nbigvtilde^{\infty}$
induces
\begin{equation}
\label{eq;20.7.22.20}
 \Psi_q^{\ast}(\nbigv,\nabla^{\lambda})_P
\stackrel{c_{\nbigv,1}}{\lrarr}
 \nbigvtilde^{\infty}_P
\stackrel{c_{\nbigv,2}}{\lrarr}
  \nbigvtilde^{\infty}_P
   \otimes_{\nbigk_{\Hhat_{\infty,q}}(\ast H_{\infty,q})_P}
   \nbigo_{\Hhat_{\infty,q}}(\ast H_{\infty,q})_P.
\end{equation}
Similarly,
for a $\nabla^{\lambda}_x$-invariant lattice
$\nbigl$ of $(\nbigv,\nabla^{\lambda})$,
we obtain
\begin{equation}
\label{eq;20.7.22.21}
 \Psi_q^{\ast}(\nbigl,\nabla^{\lambda})_P
\stackrel{c_{\nbigl,1}}{\lrarr}
 \nbigltilde^{\infty}_P
\stackrel{c_{\nbigl,2}}{\lrarr}
  \nbigltilde^{\infty}_P
   \otimes_{\nbigk_{\Hhat_{\infty,q},P}}
   \nbigo_{\Hhat_{\infty,q},P}.
\end{equation}

\begin{lem}
\label{lem;20.7.22.30}
 The morphisms
 $c_{2,\nbigv}\circ c_{1,\nbigv}$
 and
  $c_{2,\nbigl}\circ c_{1,\nbigl}$
 are isomorphisms.
\end{lem}
\pf
Because
 $\nbigvtilde^{\infty}
=
 \Psi_q^{\ast}(\nbigv,\nabla^{\lambda})
 \otimes_{\nbigo_{\Hhat_{\infty,q}}}
 \nbigk_{\Hhat_{\infty,q}}$
and
 $\nbigltilde^{\infty}
=
 \Psi_q^{\ast}(\nbigl,\nabla^{\lambda})
 \otimes_{\nbigo_{\Hhat_{\infty,q}}}
 \nbigk_{\Hhat_{\infty,q}}$,
 the claim is obvious. 
\hfill\qed

\begin{lem}
Let $(\nbigv_i,\nabla^{\lambda})$ $(i=1,2)$
be finite dimensional $\cnum(\!(x_q^{-1})\!)$-vector spaces
with $\lambda$-connection.
Then, there exist natural isomorphisms
\begin{equation}
 \Psi_q^{\ast}
 \bigl(
 (\nbigv_1,\nabla^{\lambda})
 \oplus(\nbigv_2,\nabla^{\lambda})
 \bigr)
 \simeq
 \Psi_q^{\ast}\bigl((\nbigv_1,\nabla^{\lambda})\bigr)
 \oplus
 \Psi_q^{\ast}\bigl((\nbigv_2,\nabla^{\lambda})\bigr)
\end{equation}
\begin{equation}
\label{eq;20.8.1.30}
  \Psi_q^{\ast}
 \bigl(
 (\nbigv_1,\nabla^{\lambda})
 \otimes(\nbigv_2,\nabla^{\lambda})
 \bigr)
 \simeq
 \Psi_q^{\ast}\bigl((\nbigv_1,\nabla^{\lambda})\bigr)
 \otimes
 \Psi_q^{\ast}\bigl((\nbigv_2,\nabla^{\lambda})\bigr)
\end{equation}
\begin{equation}
 \Psi_q^{\ast}\Bigl(
 \Hom\bigl(
 (\nbigv_1,\nabla^{\lambda}),
  (\nbigv_2,\nabla^{\lambda})
 \bigr)
 \Bigr)
 \simeq
 \nhom\Bigl(
 \Psi_q^{\ast}\bigl((\nbigv_1,\nabla^{\lambda})\bigr),
  \Psi_q^{\ast}\bigl((\nbigv_2,\nabla^{\lambda})\bigr)
 \Bigr).
\end{equation}
In particular,
for a finite dimensional $\cnum(\!(x_q^{-1})\!)$-vector space
$\nbigv$ with a $\lambda$-connection $\nabla^{\lambda}$,
we obtain
\[
 \Psi_q^{\ast}\bigl(
 (\nbigv,\nabla^{\lambda})^{\lor}
 \bigr)
 \simeq
 \Psi_q^{\ast}\bigl((\nbigv,\nabla^{\lambda})\bigr)^{\lor}.
\]
\end{lem}
\pf
We explain only the claim for the tensor product.
The others are similar.
We set $\nbigv_3:=\nbigv_1\otimes\nbigv_2$.
By the construction,
there exists a natural isomorphism
\[
 \nbigvtilde_1^{\infty}
  \otimes_{\nbigk_{\Hhat_{\infty,q}}(\ast H_{\infty,q})}
  \nbigvtilde_2^{\infty}
  \simeq
 \nbigvtilde_3^{\infty}.
\]
The natural morphisms
$\Psi_q^{\ast}(\nbigv_i,\nabla^{\lambda})
\lrarr\nbigvtilde^{\infty}_i$ $(i=1,2)$
induce the following:
\[
 \Psi_q^{\ast}(\nbigv_1,\nabla^{\lambda})
 \otimes
 \Psi_q^{\ast}(\nbigv_2,\nabla^{\lambda})
 \lrarr
 \nbigvtilde^{\infty}_3.
\]
It factors through
$\Psi_q^{\ast}(\nbigv_3,\nabla^{\lambda})$,
and it induces the desired isomorphism
(\ref{eq;20.8.1.30}).
\hfill\qed

\begin{rem}
\label{rem;20.8.7.10}
We can compare the constructions
in this section and
{\rm\S\ref{subsection;21.8.12.32}}
by replacing
$(t,y,x)$ with
$(t_1,\beta_1,(1+|\lambda|^2)w)$.
See Remark {\rm\ref{rem;20.8.7.11}}.
\hfill\qed
\end{rem}

\subsection{Equivalence}

The following proposition is fundamental in our study.

\begin{prop}
\label{prop;20.7.18.1}
The functor $\Psi^{\ast}$ induces an equivalence
of the categories of the following objects.
\begin{itemize}
 \item
 finite dimensional
 $\cnum(\!(x_q^{-1})\!)$-modules with $\lambda$-connection
 of Poincar\'{e} rank $<q$.
 \item
      locally free
      $\nbigo_{\Hhat_{\infty,q}}(\ast H_{\infty,q})$-modules
      of finite rank
      with level $\leq 1$.
\end{itemize}
 It also induces an equivalence of the categories of
 the following objects.
\begin{itemize}
 \item finite dimensional $\cnum(\!(x_q^{-1})\!)$-modules with
       regular $\lambda$-connection.
 \item locally free
       $\nbigo_{\Hhat_{\infty,q}}(\ast H_{\infty,q})$-modules
       of finite rank whose level is $0$.
\end{itemize}
\end{prop}
\pf
Let us prove that the functor is essentially surjective.
(See also \S\ref{subsection;21.9.12.1}
for a simplified proof in easy cases.)
Let $\nbige$ be a locally free
$\nbigo_{\Hhat_{\infty,q}}(\ast H_{\infty,q})$-module
of rank $r$
whose level is less than $1$.
Let $\nbige_0$ be an $\nbigo_{\Hhat_{\infty,q}}$-lattice
of $\nbige$ such that
the monodromy of $\Loc_1(\nbige_0)$ is unipotent.
We set
\[
 \nbige^{\infty}:=
 \nbigk_{\Hhat_{\infty,q}}
 \otimes_{\nbigo_{\Hhat_{\infty,q}}} \nbige,
 \quad\quad
  \nbige^{\infty}_0:=
 \nbigk_{\Hhat_{\infty,q}}
 \otimes_{\nbigo_{\Hhat_{\infty,q}}} \nbige_0.
\]
We obtain the naturally defined operators
$\del_{\nbige^{\infty},t}:
\nbige^{\infty}\lrarr\nbige^{\infty}$
and
$\del_{\nbige^{\infty},t}:
\nbige_0^{\infty}\lrarr\nbige^{\infty}_0$.

Note that
$\Psi_q^{\ast}(x_q^{-1})
=y_q^{-1}(1-2\sqrt{-1}\lambda t/y_q^q)^{-1/q}$
is a global section of
$\nbigk_{\Hhat_{\infty,q}}$ on $H_{\infty,q}$,
which is denoted by $x_q^{-1}$ to simplify the description.
Any section $f$ of
$\nbigk_{\Hhat_{\infty,q}}$ on $H_{\infty,q}$
is expanded as
$f=\sum_{j=0}^{\infty} f_j(t)x_q^{-j}$
for $f_j(t)\in C^{\infty}(H_{\infty,q})$.

There exists a global frame $\vecu$ of
$\nbige^{\infty}_0$.
We obtain
$A\in M_r(\nbigk_{\Hhat_{\infty,q}}(H_{\infty,q}))$
determined by
$\del_{\nbige^{\infty},t}\vecu
=\vecu A$.
We have the expansion
$A=\sum_{j=0}^{\infty}A_jx_q^{-j}$.

\begin{lem}
We may assume that
$A_0$ is constant and nilpotent.
\end{lem}
\pf
Let $\nbigc^{\infty}_{H_{\infty,q}}$
denote the sheaf of $C^{\infty}$-functions
on $H_{\infty,q}$.
We obtain the locally free
$\nbigc^{\infty}_{H_{\infty,q}}$-module
$\Loc_1(\nbige_0)\otimes\nbigc^{\infty}_{H_{\infty,q}}$.
It is equipped with the frame $\vecu_0$
induced by $\vecu$,
and the connection $\del_t$
induced by $\del_{\nbige^{\infty},t}$.
We obtain 
$\del_t\vecu_0=\vecu_0A_0$.
Because the monodromy of
$\Loc_1(\nbige_0)$ is unipotent,
there exists a frame $\vecu'_0$
of 
$\Loc_1(\nbige_0)\otimes\nbigc^{\infty}_{H_{\infty,q}}$
such that
$\del_t\vecu'_0=\vecu'_0\cdot N$
for a constant nilpotent matrix $N$.
There exists $G\in \GL_r(C^{\infty}(H_{\infty,q}))$
such that
$\vecu'_0=\vecu_0\cdot G$.
By considering $\vecu\cdot G$
instead of $\vecu$,
we may assume that $A_0$
is a constant nilpotent matrix.
\hfill\qed

\begin{lem}
\label{lem;20.7.17.10}
 There exist matrices
 $G^{(i)}\in M_r(C^{\infty}(H_{\infty,q}))$ $(i=1,2,\ldots)$
 such that the following condition is satisfied.
\begin{itemize}
 \item
      We set $\vecu^{(i)}=
      \vecu(I_r+x_q^{-1}G^{(1)})(I_r+x_q^{-2}G^{(2)})
      \cdots (I_r+x_q^{-i}G^{(i)})$.
      Let $A^{(i)}$ be the matrix determined by
      $\del_{\nbige^{\infty},t}\vecu^{(i)}
      =\vecu^{(i)}A^{(i)}$.
      Then, for the expansion
      $A^{(i)}=\sum_{j=0}^{\infty} A^{(i)}_jx_q^{-j}$,
      $A^{(i)}_j$ $(j\leq i)$ are constant.      
\end{itemize}
\end{lem}
\pf
We shall construct such $G^{(i)}$ inductively.
Suppose that we have already constructed
$G^{(i)}$ $(i\leq \ell-1)$.
By the assumption,
$A^{(\ell-1)}_j$ $(j<\ell)$ are constant.
Moreover, by the construction,
$A^{(\ell-1)}_0=A_0$ is a constant nilpotent matrix.
There exists
$G^{(\ell)}(t)\in M_r(C^{\infty}(H_{\infty,q}))$
such that
\[
 \del_tG^{(\ell)}(t)+\bigl(
 -G^{(\ell)}(t)A_0+A_0G^{(\ell)}(t)
 \bigr)
 +A^{(\ell-1)}_{\ell}(t)
\]
is constant.
We set
$\vecu^{(\ell)}
=\vecu^{(\ell-1)}(I_r+x_q^{-\ell}G^{(\ell)})$.
We obtain
\[
 A^{(\ell)}
 =(I_r+x_q^{-\ell}G^{(\ell)})^{-1}
\Bigl(
\del_tG^{(\ell)}x_q^{-\ell}
+G^{(\ell)}\del_tx_q^{-\ell}
+A^{(\ell-1)}(I_r+x_q^{-\ell}G^{(\ell)})
\Bigr).
\]
Note that 
$\del_tx_q^{-j_0}
=2\sqrt{-1}\lambda j_0q^{-1}x_q^{-j_0-q}
=-2\sqrt{-1}\lambda \del_x(x_q^{-j_0})$.
Because
\[
 A^{(\ell)}_{\ell}
 = \del_tG^{(\ell)}(t)+\bigl(
 -G^{(\ell)}(t)A_0+A_0G^{(\ell)}(t)
 \bigr)
 +A^{(\ell-1)}_{\ell}(t),
\]
$A^{(\ell)}_{\ell}$ is constant.
Hence, the induction can proceed.
\hfill\qed

\vspace{.1in}
By Lemma \ref{lem;20.7.17.10},
there exists a global frame $\vecu^{(\infty)}$
of $\nbige_0^{\infty}$
with the following property.
\begin{itemize}
 \item Let $A^{(\infty)}$ be the matrix determined by
       $\del_{\nbige^{\infty},t}\vecu^{(\infty)}
       =\vecu^{(\infty)}A^{(\infty)}$.
       Then, for the expansion
       $A^{(\infty)}=\sum_{j=0}^{\infty} A_j^{(\infty)}x_q^{-j}$,
       $A^{(\infty)}_j$ are constant.
       Moreover, $A^{(\infty)}_0$ is nilpotent.
\end{itemize}
We set
$\nbigv=\bigoplus_{i=1}^r\cnum(\!(x_q^{-1})\!)e_i$.
We define the $\lambda$-connection by
\[
 \nabla^{\lambda}_x\vece
=\vece (-2\sqrt{-1})^{-1}A^{(\infty)}.
\]
Then, Poincar\'e rank of $(\nbigv,\nabla^{\lambda})$
is strictly smaller than $q$,
and there exists an isomorphism
$\Psi_q^{\ast}(\nbigv,\nabla^{\lambda})
\simeq\nbige$.
If the monodromy of $\Loc_q(\nbige_0)$ is the identity,
we can easily observe that
$\Atilde^{(\infty)}_j=0$ $(0\leq j<q)$,
i.e., $\nabla^{\lambda}_x$ is regular.
Hence, 
$\Psi_q^{\ast}$ is essentially surjective.

\vspace{.1in}

The functor $\Psi_q^{\ast}$ is clearly faithful.
Let us prove that $\Psi_q^{\ast}$ is full.
Let $(\nbigv_i,\nabla^{\lambda})$ $(i=1,2)$ be 
$\cnum(\!(x_q^{-1})\!)$-modules with $\lambda$-connection
whose Poincar\'e rank are strictly smaller than $q$.
Let 
$F:\Psi_q^{\ast}(\nbigv_1,\nabla^{\lambda})
\lrarr
 \Psi_q^{\ast}(\nbigv_2,\nabla^{\lambda})$
 be a morphism of
 $\nbigo_{\Hhat_{\infty,q}}(\ast H_{\infty,q})$-modules.
Let us prove that 
$F$ is induced by a morphism 
$(\nbigv_1,\nabla^{\lambda})\lrarr
 (\nbigv_2,\nabla^{\lambda})$.

There exist $\nabla^{\lambda}_x$-invariant
lattices $\nbigl_i\subset\nbigv_i$ such that 
the induced endomorphisms
$F_{\nabla^{\lambda}}$ of 
 $\nbigl_{i|\infty_{x,q}}$ are nilpotent.
Let $\vecv_i$ be a frame of $\nbigl_i$.
Let $A_i(x_q^{-1})$ be determined by
$\nabla^{\lambda}_x\vecv_i
=\vecv_i A_i(x_q^{-1})$.
Note that $A_i(0)$ are nilpotent.

We obtain frames
$\vecvtilde_i:=\Psi_q^{-1}(\vecv_i)$
of
$\nbigvtilde^{\infty}_i:=
\Psi_q^{-1}(\nbigv_i)\otimes
_{\Psi_q^{-1}(\nbigo_{\inftyhat_{x,q}}(\ast \infty_{x,q}))}
 \nbigk_{\Hhat_{\infty,q}}(\ast H_{\infty,q})$.
The sheaves are equipped with
 the operators
 $\del_{\nbigvtilde^{\infty}_i,t}$. 
We set 
$\Atilde_i:=-2\sqrt{-1}\Psi_q^{\ast}A_i$,
and then 
$\del_{\nbigvtilde_i^{\infty},t}\vecvtilde_i=
\vecvtilde_i\Atilde_i(x_q^{-1})$
holds.

Let $B$ be the matrix determined by 
$F(\vecv_1)=\vecv_2B$
whose entries are sections of
$\nbigk_{\Hhat_{\infty,q}}(\ast H_{\infty,q})$
on $H_{\infty,q}$.
Because $F\circ\del_{\nbigvtilde^{\infty}_1,t}
=\del_{\nbigvtilde^{\infty}_2,t}\circ F$,
we obtain
\begin{equation}
\label{eq;20.7.17.21}
 \del_tB-B\Atilde_1+\Atilde_2B=0. 
\end{equation}
There exists the expansion
$B=\sum_{j\geq -N}^{\infty} B_j(t)x_q^{-j}$,
where $B_j(t)$ are matrices whose entries are
$C^{\infty}$-functions on $H_{\infty,q}$.

\begin{lem}
\label{lem;20.7.17.20}
 $B_j(t)$ are constant.
\end{lem}
\pf
Suppose that there exists $j$ such that
$B_j(t)$ are non-constant.
Let $j_0$ be the minimum of such $j$.
We have the expansions
$\Atilde_i=\sum_{j=0}^{\infty}\Atilde_{i,j}x_q^{-j}$.
We obtain that
$\del_tB_{j_0}(t)+\Atilde_{2,0}B_{j_0}-B_{j,0}\Atilde_{1,0}$
is constant.
Because $\Atilde_{i,0}$ are nilpotent,
we obtain that $B_{j_0}$ is constant,
which contradicts our choice of $j_0$.
Hence, we obtain that $B_j$ are constant
for any $j$.
\hfill\qed

\vspace{.1in}
By Lemma \ref{lem;20.7.17.20},
$F$ is induced by the morphism
$f:\nbigv_1\lrarr\nbigv_2$
defined as
$f(\vecv_1)=\vecv_2B$.
The relation (\ref{eq;20.7.17.21})
implies that $f$ is compatible
with the $\lambda$-connections.
Hence,
$\Psi_q^{\ast}$ is full,
and the proof of Proposition \ref{prop;20.7.18.1}
is completed.
\hfill\qed

\begin{rem}
Proposition {\rm\ref{prop;20.7.18.1}} allows us to understand
the formal structure of
$\nbigo_{\Hhat_{\infty,q}}(\ast H_{\infty,q})$-modules
with level $\leq 1$
in terms of $\lambda$-flat bundles.
It is useful for our understanding of 
the asymptotic behaviour of monopoles.
See {\rm\S\ref{subsection;21.8.13.101}}.
\hfill\qed
\end{rem}

\begin{rem}
Proposition {\rm\ref{prop;17.11.20.1}} and Proposition {\rm\ref{prop;20.7.18.1}}
provides us with an equivalence between some classes of
$\lambda$-connections and $2\sqrt{-1}\lambda T$-difference modules.
See {\rm\S\ref{subsection;17.9.13.1}} and {\rm\S\ref{subsection;17.9.13.2}}
for some explicit examples.
\hfill\qed
\end{rem}

We also obtained the following
from the above proof.

\begin{cor}
\label{cor;21.8.27.1}
Let $(\nbigv,\nabla^{\lambda})$
 be a $\cnum(\!(x_q^{-1})\!)$-module 
 with a $\lambda$-connection
 whose Poincar\'e rank is strictly less than $q$.
 Then, the constructions in
 {\rm\S\ref{subsection;20.7.18.4}} and
 {\rm\S\ref{subsection;17.10.7.2}}
 induce natural bijective correspondences
 between the following objects.
\begin{itemize}
 \item
      $\nabla^{\lambda}_x$-invariant lattices $\nbigl$
      of $(\nbigv,\nabla^{\lambda})$ such that
 the induced endomorphism of $\nbigl_{|\infty}$
 is nilpotent.
 \item $\nbigo_{\Hhat_{\infty,q}}$-lattices
       $\nbige_0$ of
       $\Psi_q^{\ast}(\nbigv,\nabla^{\lambda})$ such that
       the monodromy of
       $\Loc_1(\nbige_0)$ is unipotent.       
 \item
      $\cnum[\![y_q^{-1}]\!]$-lattices $\nbign_0$ of
      the $(2\sqrt{-1}\lambda T)$-difference module
 $\Psi_q^{\ast}(\nbigv,\nabla^{\lambda})_{|\Hhat_{\infty,q}\langle 0\rangle}$
 such that
      $\Phi^{\ast}(\nbign_0)=\nbign_0$
      and that 
      the induced automorphism of
      $\nbign_{0|\infty}$ is unipotent.
\hfill\qed
\end{itemize}
 \end{cor}

\subsubsection{Simpler cases of Proposition \ref{prop;20.7.18.1}}
\label{subsection;21.9.12.1}

Let us explain a simplified proof of the essential surjectivity 
in Proposition \ref{prop;20.7.18.1}
in some special cases, as an explanation of basic ideas of the proof.
Let $\nbige$ and $\nbige_0$ be as in the proof of
Proposition \ref{prop;20.7.18.1}.

Let us consider the case $\lambda=0$.
We obtain the $\cnum[\![y_q^{-1}]\!]$-module
$\nbigl=\nbige_{0|\Hhat_{\infty,q}\langle 0\rangle}$.
If $\lambda=0$,
we obtain the $\cnum[\![y_q^{-1}]\!]$-automorphism $M$ of $\nbigl$.
By the assumption,
the induced automorphism $M_{|y^{-1}_q=0}$ of $\nbigl/y_q^{-1}\nbigl$
is unipotent.
Hence, there exists an 
$\cnum[\![y_q^{-1}]\!]$-endomorphism $N$ of $\nbigl$
such that (i) $\exp(2\sqrt{-1}TN)=M$,
(ii) the induced endomorphism $N_{|y_q^{-1}=0}$ of
$\nbigl/y_q^{-1}\nbigl$ is nilpotent.
We obtain the meromorphic Higgs field $\nabla^0$ of $\nbigl$
by $\nabla^0_x=N$.
Then, it is easy to see that
$\Psi_q^{\ast}(\nbigl,\nabla^0)\simeq\nbige_0$,
and we obtain the essential surjectivity in this case.
We also remark that we can obtain a similar equivalence
even in the analytic case, not only in the formal case,
if $\lambda=0$.

\vspace{.1in}
Let us consider the case where
$\lambda\neq 0$ but $\rank\nbige=1$.
Let $v$ be a frame of $\nbige_0$.
There exists a global section 
$a=\sum_{j=1}^{\infty}a_j(t)x_q^{-j}$
of $\nbigk_{\Hhat_{\infty,q}}$
such that
$\del_{t}v=va$.
We would like to find a section
$b=\sum_{j=1}^{\infty}b_j(t)x_q^{-j}$ of $\nbigk_{\Hhat_{\infty,q}}$
such that 
$\del_{t}(b)+a\in x_q^{-1}\cnum[\![x_q^{-1}]\!]$.
Indeed, if there exists such $b$,
we put $v'=v\exp(b)$ for which $\del_{t}v'=v'(\del_tb+a)$ holds,
and hence we obtain that
$(\nbige_0^{\infty},\del_t)$ comes from
a $\lambda$-connection of rank one.
To find $b$ such that
$\del_tb+a\in x_q^{-1}\cnum[\![x_q^{-1}]\!]$,
we note the following.
\begin{itemize}
 \item For a $C^{\infty}$-function $f$
       on $H_{\infty,q}\simeq S^1_T$,
       there exists a function $g$
       such that $f-\del_tg$ is constant.
 \item We have
       $\del_tx_q^{-j}=2\sqrt{-1}\lambda j q^{-1}x_q^{-j-q}$.
\end{itemize}
Then, by using a standard inductive argument,
we can easily construct $b$ such that
$\del_tb+a\in x_q^{-1}\cnum[\![x_q^{-1}]\!]$.

\subsection{Example 1}
\label{subsection;17.9.13.1}

Let $\gminia\in x_q\cnum[x_q]$
such that $\deg_{x_q}\gminia<q$.
Let us consider
the $\cnum(\!(x_q^{-1})\!)$-module
$\nbigv=\cnum(\!(x_q^{-1})\!)\cdot v$
with the $\lambda$-connection
$\nabla^{\lambda}v
=v\cdot d\gminia$.
We obtain
an $\nbigk_{\Hhat_{\infty,q}}(\ast H_{\infty,q})$-module
$\nbigvtilde^{\infty}$
as in (\ref{eq;20.7.18.2}),
which is equipped
with the induced frame $\vtilde=\Psi_q^{-1}(v)$
and the operator $\del_{\nbigvtilde^{\infty},t}$
defined as in (\ref{eq;20.7.18.3}).
We have
\[
 \del_{\nbigvtilde^{\infty},t}\vtilde
=\vtilde
\Bigl(
 -2\sqrt{-1}
 \bigl(
 \del_x\gminia(x_q)
 \bigr)_{|x_q=y_q(1-2\sqrt{-1}\lambda ty^{-1})^{1/q}}
 \Bigr).
\]
Let us compute the corresponding difference module
explicitly.

We set
$\nbigvtilde^{\infty,\cov}:=\varpi_q^{-1}\nbigvtilde^{\infty}$
which is
an $\nbigk_{\Hhat^{\cov}_{\infty,q}}(\ast H^{\cov}_{\infty,q})$-module.
It is equipped with the induced frame
$\vtilde^{\cov}=\varpi_q^{-1}(\vtilde)$.
It is also equipped with the operator
induced by $\del_{\nbigvtilde^{\infty\cov},t}$,
which we denote by $\del_t$ to simplify the description.

\subsubsection{}
Let us study the case $\lambda\neq 0$.
Note that
\[
 \del_{t}\vtilde^{\cov}
=\vtilde^{\cov}
 \del_t
 \Bigl(
 \lambda^{-1}
 \gminia\bigl(y_q(1-2\sqrt{-1}\lambda ty^{-1})^{1/q}\bigr)
 \Bigr).
\]
Because $\deg_{x_q}\gminia<q$,
\[
 -\lambda^{-1}\gminia\bigl(y_q(1-2\sqrt{-1}\lambda ty_q^{-q})^{1/q}\bigr)
 +\lambda^{-1}\gminia(y_q)
\]
is a section of
$y_q^{-1}\nbigk_{\Hhat^{\cov}_{\infty,q}}(\ast H^{\cov}_{\infty,q})$
on $H^{\cov}_{\infty,q}$.
Hence, 
\[
 \utilde=\vtilde^{\cov}\cdot
 \exp\Bigl(
 -\lambda^{-1}\gminia(y_q(1-2\sqrt{-1}\lambda ty^{-1})^{1/q})
 +\lambda^{-1}\gminia(y_q)
 \Bigr)
\]
is a section of
$\nbigvtilde^{\infty\,\cov}$
such that $\del_{t}\utilde=0$.
Note that $\kappa_1^{\ast}(\vtilde^{\cov})=\vtilde^{\cov}$.
We also note that $\kappa_1^{\ast}(x_q)=x_q$
which implies that
$\kappa_1^{\ast}(y_q(1-2\lambda\sqrt{-1}ty^{-1})^{1/q})
=y_q(1-2\lambda\sqrt{-1}ty^{-1})^{1/q}$.
Hence, we obtain
\[
 \kappa_1^{\ast}(\utilde)
 =\utilde\cdot
 \exp\Bigl(
 \lambda^{-1}\gminia
 \bigl(y_q(1+2\sqrt{-1}\lambda Ty^{-1})^{1/q}\bigr)
-\lambda^{-1}\gminia(y_q)
 \Bigr).
\]
Therefore,
the corresponding difference module
is
the $\cnum(\!(y_q^{-1})\!)$-module
$\cnum(\!(y_q^{-1})\!)\cdot u$
with the difference operator defined as
\[
 \Phi^{\ast}u
=u\cdot
 \exp\Bigl(
 \lambda^{-1}
 \gminia\bigl(y_q(1+2\sqrt{-1}\lambda Ty^{-1})^{1/q}\bigr)
-\lambda^{-1}
 \gminia(y_q)
 \Bigr).
\]

\subsubsection{}
Let us study the case $\lambda=0$.
Note that
\[
 \del_t\vtilde^{\cov}
 =\vtilde^{\cov}
 \bigl(-2\sqrt{-1}(\del_x\gminia)(y_q)\bigr).
\]
We set
$\utilde=
\vtilde^{\cov}
\exp\bigl(
2t\sqrt{-1}(\del_x\gminia)(y_q)\bigr)$,
which is a frame of
$\nbigvtilde^{\infty\,\cov}$
such that
$\del_t\utilde=0$.
Because $\kappa_1^{\ast}\vtilde=\vtilde$,
we obtain
\[
 \kappa_1^{\ast}(\utilde)
 =\utilde
 \cdot
 \exp\bigl(2T\sqrt{-1}(\del_x\gminia)(y_q)\bigr).
\]
Therefore, 
the corresponding difference module
is
the $\cnum(\!(y_q^{-1})\!)$-module
$\cnum(\!(y_q^{-1})\!)\cdot u$
with the difference operator $\Phi^{\ast}$
defined as
\[
 \Phi^{\ast}u
 =u\cdot
 \exp\bigl(2T\sqrt{-1}(\del_x\gminia)(y_q)\bigr).
\]

\subsection{Example 2}
\label{subsection;17.9.13.2}

We set 
$\nbigv=\bigoplus_{i=1}^r
 \cnum(\!(x_q^{-1})\!)v_i$.
Let $A\in M_r(\cnum)$.
 We consider the $\lambda$-connection
 $\nabla^{\lambda}$ on $\nbigv$
 defined by
$\nabla^{\lambda}_x\vecv
=\vecv\cdot\bigl(
 x^{-1}A
 \bigr)$.
 We obtain
the $\nbigk_{\Hhat_{\infty,q}}(\ast H_{\infty,q})$-module
$\nbigvtilde^{\infty}$
as in (\ref{eq;20.7.18.2})
with the operator
$\del_{\nbigvtilde^{\infty},t}$
and the frame $\vecvtilde=\Psi_q^{-1}(\vecv)$.
We obtain
\[
 \del_{\nbigvtilde^{\infty},t}\vecvtilde
=\vecvtilde\Bigl(
 -2\sqrt{-1}(y-2\sqrt{-1}\lambda t)^{-1}A
 \Bigr).
\]

Let us compute the corresponding
difference module explicitly.
As in \S\ref{subsection;17.9.13.1},
we set
$\nbigvtilde^{\infty,\cov}=\varpi_q^{-1}(\nbigvtilde^{\infty})$,
which is equipped with the induced frame
$\vecvtilde^{\cov}=\varpi_q^{-1}(\vecvtilde)$
and the induced differential operator $\del_t$.

\subsubsection{}
Let us study the case $\lambda\neq 0$.
We obtain a frame $\vecutilde$
of $\nbigvtilde^{\infty\,\cov}$
satisfying $\del_t\vecutilde=0$
as follows:
\[
 \vecutilde:=
\vecvtilde\cdot
 \exp\Bigl(
 -\lambda^{-1}A
 \log\bigl((1-2\sqrt{-1}\lambda ty^{-1})\bigr)
 \Bigr).
\]
The following holds:
\[
 \kappa_1^{\ast}(\vecutilde)
 =\vecutilde\cdot
  \exp\Bigl(
 \lambda^{-1}A
 \log\bigl(
 y^{-1}(y+2\sqrt{-1}\lambda T)
 \bigr)
 \Bigr).
\]
Hence, 
the corresponding difference module is
$\bigoplus_{i=1}^r \cnum(\!(y_q^{-1})\!)u_i$
with the difference operator $\Phi^{\ast}$
defined as
\[
 \Phi^{\ast}\vecu
=\vecu
 \exp\Bigl(
 \lambda^{-1}A
 \log\bigl(
 y^{-1}(y+2\sqrt{-1}\lambda T)
 \bigr)
 \Bigr).
\]
\subsubsection{}
Let us study the case $\lambda=0$.
We obtain a frame $\vecutilde$ of
$\nbigvtilde^{\infty\cov}$
such that $\del_t\vecutilde=0$
as follows:
\[
 \vecutilde=\vecvtilde
 \exp\bigl(
  2\sqrt{-1}ty^{-1}A
 \bigr).
\]
Note that
$\kappa_1^{\ast}(\vecutilde)
 =\vecutilde\exp\Bigl(
 2\sqrt{-1}Ty^{-1}A
 \Bigr)$.
Hence, 
the corresponding difference module is
$\bigoplus_{i=1}^r \cnum(\!(y_q^{-1})\!)u_i$
with the difference operator $\Phi^{\ast}$
defined as
$\Phi^{\ast}\vecu
=\vecu
 \exp\Bigl(
 2\sqrt{-1}Ty^{-1}A
 \Bigr)$.

\subsection{Comparison of good filtered bundles}

Let $(\nbigv,\nabla^{\lambda})$
be a $\cnum(\!(x_q^{-1})\!)$-module
of finite rank with a $\lambda$-connection.
Let $\nbigp_{\ast}\nbigv$
be a filtered bundle over $\nbigv$.
Recall that $(\nbigp_{\ast}\nbigv,\nabla^{\lambda})$
is called good
if there exist $p\in q\seisuu_{>0}$ and 
a decomposition
\[
 \varphi_{q,p}^{\ast}(\nbigp_{\ast}\nbigv)
 =\bigoplus_{\gminia\in x_q\cnum[x_q]}
 \nbigp_{\ast}\nbigv_{\gminia}
\]
such that
$(\nabla^{\lambda}-d\gminia\id)
\nbigp_a\nbigv_{\gminia}
\subset
\nbigp_a\nbigv_{\gminia}\cdot dx_q/x_q$
for any $\gminia$ and $a$.

\begin{prop}
\label{prop;20.7.20.131}
 Suppose that Poincar\'{e} rank of
 $(\nbigv,\nabla^{\lambda})$ is strictly less than $q$,
 and suppose that each
 $\nbigp_a\nbigv$ is $\nabla_x^{\lambda}$-invariant.
 Then,
 the induced filtered bundle
 $\Psi_q^{\ast}(\nbigp_{\ast}\nbigv,\nabla^{\lambda})$
 is good
 if and only if
 $(\nbigp_{\ast}\nbigv,\nabla^{\lambda})$
 is good.
\end{prop}
\pf
Let us consider the case
where $(\nbigv,\nabla^{\lambda})$ is unramified,
i.e.,
there exists
\begin{equation}
\label{eq;20.7.20.120}
 (\nbigv,\nabla^{\lambda})=
  \bigoplus_{\gminia\in \nbigi(\nbigv,\nabla^{\lambda})}
  (\nbigv_{\gminia},\nabla^{\lambda}_{\gminia})
\end{equation}
such that
$(\nbigv,\nabla^{\lambda}-d\gminia\id)$
are regular.
We obtain the induced decomposition
\begin{equation}
\label{eq;20.7.20.121}
 \Psi_q^{\ast}
   (\nbigv,\nabla^{\lambda})=
  \bigoplus_{\gminia\in \nbigi(\nbigv,\nabla^{\lambda})}
 \Psi_q^{\ast}(\nbigv_{\gminia},\nabla^{\lambda}_{\gminia}).
\end{equation}

By the equivalence in Proposition \ref{prop;20.7.18.1},
the filtered bundle
$\nbigp_{\ast}\nbigv$ is compatible
with the decomposition (\ref{eq;20.7.20.120})
if and only if
the induced filtered bundle
$\Psi_q^{\ast}(\nbigp_{\ast}\nbigv,\nabla^{\lambda})$ is compatible
with the decomposition (\ref{eq;20.7.20.121}),
according to Corollary \ref{cor;21.8.27.1}.
By the computation in \S\ref{subsection;17.9.13.1},
it is enough to consider the case
where $\nbigi(\nbigv,\nabla^{\lambda})=\{0\}$,
which is reduced to the following lemma.

\begin{lem}
\label{lem;20.7.20.130}
 Let $\nbige_0$ be a locally free
 $\nbigo_{\Hhat_{\infty,q}}$-module
 of rank $r$.
 We set
 $\nbige_0^{\infty}:=\nbigk_{\Hhat_{\infty,q}}
 \otimes_{\nbigo_{\Hhat_{\infty,q}}}\nbige_0$.
 Then,
 the monodromy of
 $\Loc_q(\nbige_0)$ is the identity
 if and only if
 there exists a frame $\vecv$ of
 $\nbige_0^{\infty}$
 such that
 $\del_t\vecv=\vecv\cdot x_q^{-q}A$
 for $A\in M_r(\cnum[\![x_q^{-1}]\!])$.
\end{lem}
\pf
Let $\nbigc^{\infty}_{H_{\infty,q}}$
denote the sheaf of $C^{\infty}$-functions
on $H_{\infty,q}$.
Note that
$\Loc_q(\nbige_0)^{\infty}:=
\Loc_q(\nbige_0)\otimes\nbigc^{\infty}_{H_{\infty,q}}$
is naturally isomorphic to
$\nbige_0^{\infty}/x_q^{-q}\nbige_0^{\infty}$.
We also remark that
the tuple
$y_q^{-j}$ $(j=0,\ldots,q-1)$
is a frame of the local system
$\Loc_q(\nbigo_{\Hhat_{\infty,q}})$,
and  equal to the tuple
induced by
$x_q^{-j}$ $(j=0,\ldots,q-1)$.

Let $\vecv=(v_1,\ldots,v_r)$
be a frame of $\nbige_0^{\infty}$.
Then, the tuple
$\vecv_0=\bigl(y_q^{-j}v_i\,\big|\,
j=0,\ldots,q-1,\,i=1,\ldots,r
\bigr)$
is a frame of $\nbige_0^{\infty}$
over $\nbigc^{\infty}_{H_{\infty,q}}$.
It is called the induced frame.

Assume that there exists a frame $\vecv$
of $\nbige_0^{\infty}$ satisfying the condition
in the lemma.
It induces a frame $\vecv_0$ of
$\Loc_q(\nbige_0)^{\infty}$ as above.
We obtain $\del_t\vecv_0=0$.
Hence, the monodromy of
$\Loc_q(\nbige_0)$ is the identity.

Conversely,
suppose that the monodromy of
$\Loc_q(\nbige_0)$ is the identity.
It implies that
$\Loc_q(\nbige_0)$
is a free module over the sheaf of algebras
$\Loc_q(\nbigo_{\Hhat_{\infty,q}})$.
There exists a frame
$\vecu=(u_1,\ldots,u_r)$ of
$\Loc_q(\nbige_0)$
over
$\Loc_q(\nbigo_{\Hhat_{\infty,q}})$.
 
Let $\vecv$ be any frame of
$\nbige_0^{\infty}$
which induces $\vecu$
as a tuple of sections of
$\Loc_q(\nbige_0)^{\infty}$.
Then, we obtain
$\del_t\vecu=\vecu\cdot x_q^{-q}A$
for a global section $A$
of $M_r(\nbigk_{\Hhat_{\infty,q}})$ on $H_{\infty,q}$.
As in the proof of Proposition \ref{prop;20.7.18.1},
after an appropriate gauge transform of the frame,
we obtain $A\in M_r(\cnum[\![x_q^{-1}]\!])$.
Thus, we obtain Lemma \ref{lem;20.7.20.130}
and the claim of Proposition \ref{prop;20.7.20.131}
in the unramified case.
\hfill\qed

\vspace{.1in}

For $p\in q\seisuu_{>0}$,
there exists the following natural
commutative diagram:
\[
 \begin{CD}
  \Hhat_{p,\infty}
  @>{\nbigr_{q,p}}>>
  \Hhat_{\infty,q}\\
  @V{\Psi_p}VV @V{\Psi_q}VV \\
  \inftyhat_{x,p}
  @>{\varphi_{q,p}}>>
  \inftyhat_{x,q}.
 \end{CD}
\]
Hence,
there exists the following natural isomorphism:
\[
 \Psi_p^{\ast}\varphi_{q,p}^{\ast}(\nbigp_{\ast}\nbigv,\nabla^{\lambda})
 \simeq
 \nbigr_{q,p}^{\ast}\Psi_q^{\ast}(\nbigp_{\ast}\nbigv,\nabla^{\lambda}).
\]
By the consideration in the unramified case,
$\varphi_{q,p}^{\ast}
 (\nbigp_{\ast}\nbigv,\nabla^{\lambda})$ is good
if and only if
$\Psi_p^{\ast}\varphi_{q,p}^{\ast}
 (\nbigp_{\ast}\nbigv,\nabla^{\lambda})$
is good.
Thus, we obtain the claim in the general case.
\hfill\qed

\subsection{Comparison of the associated graded pieces}

Let $(\nbigp_{\ast}\nbigv,\nabla^{\lambda})$
be a good filtered $\lambda$-flat bundle
whose Poincar\'{e} rank is strictly less than $1$.
We obtain the good filtered bundle
$\Psi_q^{\ast}(\nbigp_{\ast}\nbigv,\nabla^{\lambda})$
on $(\Hhat_{\infty,q},H_{\infty,q})$.

Take any $P\in H_{\infty,q}$.
For any $v\in\nbigv$,
we obtain
an element 
$\Psi_q^{-1}(v)_P$
of $\nbigvtilde^{\infty}_P$.
By using Lemma \ref{lem;20.7.22.30},
we obtain
\[
 (c_{2,\nbigv}\circ c_{1,\nbigv})^{-1}
 \Bigl(
  c_{2,\nbigv}\bigl(\Psi_q^{-1}(v)_P\bigr)
  \Bigr)
  \in \Psi_q^{\ast}(\nbigv,\nabla^{\lambda})_P.
\]
Thus, we obtain a $\cnum$-linear map
$\nbigv\lrarr \Psi_q^{\ast}(\nbigv,\nabla^{\lambda})_P$.

\begin{prop}
\label{prop;20.7.22.50}
 The induced morphism
$\Gr^{\nbigp}_a(\nbigv)
\lrarr
 \Gr^{\nbigp}_a(\Psi_q^{\ast}(\nbigv,\nabla^{\lambda}))_P$
is a $\cnum$-linear isomorphism.
 Moreover,  under the isomorphism,
 $q^{-1}(2\sqrt{-1})T\Res(\nabla^{\lambda})$
 is equal to
 $\gbigf_{\Psi_q^{\ast}(\nbigv,\nabla^{\lambda})}$.
 (See {\rm\S\ref{subsection;20.7.22.31}}
 for the endomorphism
 $\gbigf_{\Psi_q^{\ast}(\nbigv,\nabla^{\lambda})}$.)
 In particular,
 the weight filtrations of
 the nilpotent parts of
 $\Res(\nabla^{\lambda})$
 and
 $\gbigf_{\Psi_q^{\ast}(\nbigv,\nabla^{\lambda})}$
 are equal.
\end{prop}
\pf
The first claim is clear.
Let us study the case where
$(\nbigv,\nabla^{\lambda})$ is regular.
Let $\vecv$ be a frame of $\nbigp_a\nbigv$
compatible with the parabolic structure
as in \S\ref{subsection;17.10.5.301},
i.e.,
there is a decomposition
$\vecv=\coprod_{a-1<b\leq a}\vecv_b$
such that
each $\vecv_b$ is a tuple of sections of
$\nbigp_b(\nbigv)$
and induces a base of
$\Gr^{\nbigp}_b(\nbigv)$.
For $j\in\seisuu_{\geq 0}$
and $a-1<b\leq a$,
we obtain tuples of sections
$\Psi_q^{-1}(x_q^{-j}\vecv_b)$
of $\widetilde{\nbigp_a\nbigv}^{\infty}$.
Let $[\Psi_q^{-1}(x_q^{-j}\vecv_b)]$
denote the induced sections of
a locally free $\nbigc^{\infty}_{H_{\infty,q}}$-module
\begin{equation}
\label{eq;20.7.22.40}
 \widetilde{\nbigp_a\nbigv}^{\infty}
 \big/ \widetilde{\nbigp_{<a-q}\nbigv}^{\infty}.
\end{equation}
Then, the union of
$[\Psi_q^{-1}(x_q^{-j}\vecv_b)]$
$(a-1<b<a,\,\,0\leq j\leq q-1)$
and
$[\Psi_q^{-1}(x_q^{-j}\vecv_a)]$
$(0\leq j\leq q)$
is a frame of (\ref{eq;20.7.22.40})
over $\nbigc^{\infty}_{H_{\infty,q}}$.

We set $r(a):=\dim_{\cnum} \Gr^{\nbigp}_a(\nbigv)$.
Let $[\vecv_a]$ denote the base of
$\Gr^{\nbigp}_a(\nbigv)$
induced by $\vecv_a$.
We obtain the matrix $A\in M_{r(a)}(\cnum)$
by
$\Res\nabla^{\lambda}[\vecv_a]=[\vecv_a]$.
Namely,
for $\vecv_a=(v_{a,1},\ldots,v_{a,r(a)})$,
$\Res(\nabla)[v_{a,j}]
 =\sum A_{i,j}[v_{a,i}]$.

For any $a-1<b<a$,
and for any $j=0,\ldots,q-1$,
we obtain
$\del_t[\Psi_q^{-1}(x_q^{-j}\vecv_b)]=0$.
For $j=1,\ldots,q$,
we obtain
$\del_t[\Psi_q^{-1}(x_q^{-j}\vecv_a)]=0$.
Because
$\del_t[\Psi_q^{-1}(\vecv_a)]
 =[\Psi_q^{-1}(x_q^{-q}\vecv_a)]q^{-1}(-2\sqrt{-1}A)$,
 \[
 \varpi_q^{-1}\Psi_q^{-1}[v_{a,j}]
 +q^{-1}2\sqrt{-1}t\sum A_{i,j}\Psi_q^{-1}[x_q^{-q}v_{a,i}]
 \quad
  (i=1,\ldots,r(a))
 \]
are flat sections of
the pull back of (\ref{eq;20.7.22.40})
by $\varpi_q$.
Thus, we obtain the claim of the proposition
in the case where $(\nbigv,\nabla^{\lambda})$ is regular.

As in the proof of Proposition \ref{prop;20.7.20.131},
by using the computation in \S\ref{subsection;17.9.13.1}
and the pull back by $\varphi_{q,p}$
such that $\varphi_{q,p}^{\ast}(\nbigv,\nabla^{\lambda})$
is unramified,
we can deduce the claim in the general case.
\hfill\qed

\subsection{Some functoriality}

Let $(\nbigv_i,\nabla^{\lambda})$ $(i=1,2)$
be finite dimensional $\cnum(\!(x_q^{-1})\!)$-vector spaces
with $\lambda$-connection.
Let $(\nbigp_{\ast}\nbigv_i,\nabla^{\lambda})$
be good filtered $\lambda$-flat bundles over
$(\nbigv_i,\nabla^{\lambda})$.

\begin{prop}
There exist the following natural isomorphisms:
\begin{equation}
 \Psi_q^{\ast}(\nbigp_{\ast}\nbigv_1,\nabla^{\lambda})
  \oplus
 \Psi_q^{\ast}(\nbigp_{\ast}\nbigv_2,\nabla^{\lambda})
 \simeq
 \Psi_q^{\ast}\bigl(
 (\nbigp_{\ast}\nbigv_1,\nabla^{\lambda})
  \oplus
 (\nbigp_{\ast}\nbigv_2,\nabla^{\lambda})
 \bigr)
\end{equation} 
\begin{equation}
 \Psi_q^{\ast}(\nbigp_{\ast}\nbigv_1,\nabla^{\lambda})
  \otimes
 \Psi_q^{\ast}(\nbigp_{\ast}\nbigv_2,\nabla^{\lambda})
 \simeq
 \Psi_q^{\ast}\bigl(
 (\nbigp_{\ast}\nbigv_1,\nabla^{\lambda})
  \otimes
 (\nbigp_{\ast}\nbigv_2,\nabla^{\lambda})
 \bigr)
\end{equation}
\begin{equation}
 \nhom\Bigl(
   \Psi_q^{\ast}(\nbigp_{\ast}\nbigv_1,\nabla^{\lambda}),
 \Psi_q^{\ast}(\nbigp_{\ast}\nbigv_2,\nabla^{\lambda})
  \Bigr)
  \simeq
  \Psi_q^{\ast}\Bigl(
  \Hom\bigl(
  (\nbigp_{\ast}\nbigv_1,\nabla^{\lambda}),
  (\nbigp_{\ast}\nbigv_2,\nabla^{\lambda}) 
  \bigr)
  \Bigr).
\end{equation}
In particular,
 for a good filtered $\lambda$-flat bundle
 $(\nbigp_{\ast}\nbigv,\nabla^{\lambda})$
 over a finite dimensional $\cnum(\!(x_q^{-1})\!)$-vector space
 $\nbigv$ with a $\lambda$-connection $\nabla^{\lambda}$,
 we obtain
 $\Psi_q^{\ast}(\nbigv,\nabla^{\lambda})^{\lor}
 \simeq
 \Psi_q^{\ast}\Bigl(
  (\nbigv,\nabla^{\lambda})^{\lor}
 \Bigr)$. 
\hfill\qed
\end{prop}

\section{Appendix: Pull back and descent in the $\real$-direction}

To emphasize the dependence of $\Hhat_{\infty,q}$ on $T$,
we denote it by $\Hhat_{\infty,q,T}$.
We denote $\pi_q:\Hhat_{\infty,q,T}\lrarr S^1_T$
by $\pi_{q,T}$.
As explained,
we use the covering maps
$\nbigr_{q,p}:\Hhat_{\infty,p,T}\lrarr\Hhat_{\infty,q,T}$
for $p\in q\seisuu_{>0}$
to define the notion of good filtered bundles.
There also exists the naturally defined covering map
$\nbigt_{T,mT}:
\Hhat_{\infty,q,mT}\lrarr\Hhat_{\infty,q,T}$,
which we may use to define the notion of good filtered
bundles in another way.
Although we will not use it in this monograph,
we explain it because it is also a natural way.

For any $m\in\seisuu_{\geq 1}$,
there exists the following natural commutative diagram:
\[
 \begin{CD}
 \Hhat_{\infty,q,mT}
 @>{\nbigt_{T,mT}}>>
 \Hhat_{\infty,q,T} \\
 @V{\pi_{q,mT}}VV @V{\pi_{q,T}}VV \\
 S^1_{mT}
 @>{\varphi}>>
 S^1_T
 \end{CD}
\]
Note that
$\nbigt_{T,mT}^{-1}\nbigo_{\Hhat_{\infty,q,T}}
=\nbigo_{\Hhat_{\infty,q,mT}}$.
There exits the pull back 
$\nbigt_{T,mT}^{\ast}
=\nbigt_{T,mT}^{-1}$
of 
$\nbigo_{\Hhat_{\infty,q,T}}$-modules
and the push-forward $\nbigt_{T,mT\ast}$ of
$\nbigo_{\Hhat_{\infty,q,mT}}$-modules.

We obtain the action of $\seisuu/m\seisuu$
on $\Hhat_{\infty,q,mT}$
induced by $n(t,y_q)=(t+nT,y_q(1+2\sqrt{-1}\lambda Tny^{-1})^{1/q})$.
Let $\nbige_1$ be an $\nbigo_{\Hhat_{\infty,q,mT}}$-module
which is equivariant with respect to 
the $\seisuu/m\seisuu$-action.
Then,
the push-forward 
$\nbigt_{T,mT\ast}\nbige_1$
is equipped with the naturally induced
$\seisuu/m\seisuu$-action.
The invariant part is called 
the descent of $\nbige_1$
with respect to the $\seisuu/m\seisuu$-action.

Let $\bigl(
 \nbigp_{\ast}\nbige_{|\Hhat_{\infty,q,T}\langle t\rangle}\,\big|\,
 t\in S^1_T
 \bigr)$
be a filtered bundle over
a locally free $\nbigo_{\Hhat_{\infty,q,T}}(H_{\infty,q,T})$-module.
Because 
$\nbigt_{T,mT}^{\ast}(\nbige)_{|\Hhat_{\infty,q,mT}\langle t_1\rangle}
=\nbige_{|\Hhat_{\infty,q,T}\langle \varphi(t_1)\rangle}$,
we obtain a filtered bundle over
$\nbigt_{T,m_T}^{\ast}(\nbige)$
as follows:
\[
 \nbigp_{\ast}\bigl(\nbigt^{\ast}_{T,mT}(\nbige)
  _{|\Hhat_{\infty,q,mT}\langle t_1\rangle}\bigr):=
 \nbigp_{\ast}\nbige_{|\Hhat_{\infty,q,T}\langle\varphi(t_1)\rangle}
\quad
 (t_1\in S^1_{mT}).
\]

Let $\bigl(
 \nbigp_{\ast}\nbige_{1|\Hhat_{\infty,q,mT}\langle t_1\rangle}\,
 \big|\,
 t_1\in S^1_{mT}
 \bigr)$
be a filtered bundle over
a locally free $\nbigo_{\Hhat_{\infty,q,mT}}(\ast H_{\infty,q,mT})$-module
$\nbige_1$.
We obtain the filtered bundle
over $\nbigt_{T,mT\ast}\nbige_1$ given as follows:
\begin{equation}
 \label{eq;17.9.12.3}
 \nbigp_{\ast}
 \nbigt_{T,mT\ast}(\nbige_1)
 _{|\Hhat_{\infty,q,T}\langle t\rangle}
=\bigoplus_{t_1\in \varphi^{-1}(t)}
 \nbigp_{\ast}\nbige
 _{1|\Hhat_{\infty,q,mT}\langle t_1\rangle}.
\end{equation}
If $\nbige_1$ and $\nbigp_{\ast}\nbige_1$ are
$\seisuu/m\seisuu$-equivariant,
then we obtain the filtered bundle
over the descent of $\nbige_1$
by taking the invariant part of (\ref{eq;17.9.12.3}).
The following lemma is clear
by the construction.
\begin{lem}
 Let $\nbige$ be a locally free
 $\nbigo_{\Hhat_{\infty,q}}(\ast H_{\infty,q})$-module.
 Let $\nbigp_{\ast}\nbige$ be a filtered bundle 
 over $\nbige$.
 Let $\nbigp_{\ast}\nbigt_{T,mT}^{\ast}\nbige$
 be the induced filtered bundle
 over $\nbigt_{T,mT}^{\ast}\nbige$.
 Then, $\nbigp_{\ast}\nbige$
 is the descent of 
 $\nbigp_{\ast}\nbigt_{T,mT}^{\ast}\nbige$.
\hfill\qed
\end{lem}

\begin{prop}
 Let $\nbige$ be a locally free
 $\nbigo_{\Hhat_{\infty,q,T}}(\ast H_{\infty,q,T})$-module,
and let 
$\bigl(
 \nbigp_{\ast}\nbige_{|\Hhat_{\infty,q,T}\langle t\rangle}
 \,\big|\,t\in S^1_T
 \bigr)$
be a filtered bundle over $\nbige$.
Then, it is good 
if and only if 
the induced filtered bundle over
$\nbigt_{T,mT}^{\ast}\nbige$
is good.
\end{prop}
\pf
We set $\nbige'=\nbigt_{T,mT}^{\ast}\nbige$.
Let us consider the case where
the level of $\nbige$ is less than $1$.
We use the notation $\Psi_{q,T}$ to emphasize
the dependence on $T$.
There exists
a $\cnum(\!(x_q^{-1})\!)$-module with $\lambda$-connection
$(\nbigv,\nabla^{\lambda})$
such that
$\nbige=\Psi_{q,T}^{\ast}(\nbigv,\nabla)$.
By construction, we obtain the natural isomorphism
$\nbige'=
\nbigt_{T,mT}^{\ast}\circ
 \Psi_{q,T}^{\ast}(\nbigv,\nabla_{\lambda})
 \simeq
\Psi_{q,mT}^{\ast}(\nbigv,\nabla^{\lambda})$.
If $\nbigp_{\ast}\nbige$ is good,
there exists a good filtered $\lambda$-flat bundle
$(\nbigp_{\ast}\nbigv,\nabla^{\lambda})$
over $(\nbigv,\nabla^{\lambda})$
such that
$\nbigp_{\ast}\nbige=\Psi_{q,T}^{\ast}(\nbigp_{\ast}\nbigv)$.
Because 
$\nbigp_{\ast}\nbige'
=\Psi_{q,mT}^{\ast}(\nbigp_{\ast}\nbigv,\nabla^{\lambda})$,
it is also good.
The converse can be proved similarly.

We use the notation
$\LLhat^{\lambda}_{q,T}(\ell,\alpha)$
to emphasize the dependence on $T$.
According to the computation in \S\ref{subsection;17.10.4.30},
there exists the natural isomorphism
$\nbigt_{T,mT}^{\ast}\LLhat^{\lambda}_{q,T}(\ell,\alpha)
\simeq
 \LLhat^{\lambda}_{q,mT}(m\ell,\alpha^m)$.
Under the isomorphism,
we obtain
\[
 \nbigp_{\ast}
 \nbigt_{T,mT}^{\ast}\LLhat^{\lambda}_{q,T}(\ell,\alpha)
=\nbigp_{\ast}\LLhat^{\lambda}_{q,mT}(m\ell,\alpha^m)
\]
by the construction of the filtrations.

\vspace{.1in}

Let $\nbige$ be a locally free 
$\nbigo_{\Hhat_{\infty,q,T}}(\ast H_{\infty,q,T})$-module.
We have $p\in q\seisuu_{\geq 1}$
such that
$\nbige_1:=
 \nbigr_{p,q}^{\ast}\nbige$ is unramified.
Let $\nbige':=\nbigt_{T,mT}^{\ast}\nbige$
and $\nbige'_1:=\nbigt_{T,mT}^{\ast}\nbige_1$.
It is easy to observe that
$\nbige'$ is the descent of $\nbige'_1$
with respect to the ramified covering
$\nbigr_{p,q,mT}:
\Hhat_{\infty,p,mT}\lrarr \Hhat_{\infty,q,mT}$.

Let $\nbigp_{\ast}\nbige$ be a good filtered bundle
over $\nbige$.
The induced filtered bundle
$\nbigp_{\ast}\nbige_1$
over $\nbige_1$ is good.
Let $\nbigp_{\ast}\nbige'_1$
be the filtered bundle
obtained as the pull back of 
$\nbigp_{\ast}\nbige_1$
by $\nbigt_{T,mT}$.
Then, by the previous consideration,
$\nbigp_{\ast}\nbige'_1$ is good.
We can easily observe that 
$\nbigp_{\ast}\nbige'_1$
is $\Gal_{q,p}$-equivariant,
and 
$\nbigp_{\ast}\nbige'$
is the descent of $\nbigp_{\ast}\nbige'_1$.
Hence, $\nbigp_{\ast}\nbige'$ is good.

Conversely,
suppose that $\nbigp_{\ast}\nbige'$ is good.
The pull back
$\nbigp_{\ast}\nbige_1'$ over $\nbige'_1$
is good.
We can easily observe that
$\nbigp_{\ast}\nbige_1'$ is
equivariant with respect to
$\seisuu/m\seisuu$-action
and $\Gal_{q,p}$-action.
Let $\nbigp_{\ast}\nbige_1$
be the descent of $\nbigp_{\ast}\nbige'_1$.
It is good, $\Gal_{q,p}$-equivariant,
and the descent is $\nbigp_{\ast}\nbige$.
Hence, we obtain the claim of the proposition.
\hfill\qed

\begin{cor}
 Let $\nbige'$ be a locally free
 $\nbigo_{\Hhat_{\infty,q,mT}}(\ast H_{\infty,q,mT})$-module.
 Let $\nbigp_{\ast}\nbige'$ be
 a $\seisuu/m\seisuu$-equivariant
 good filtered bundle over $\nbige'$.
 Then, the descent of $\nbigp_{\ast}\nbige'$ is also good.
\hfill\qed
\end{cor}

\subsection{Examples}
\label{subsection;17.10.4.30}

Let $q\in\seisuu_{\geq 1}$ and $\ell\in\seisuu$.
We consider the automorphism $\Phi^{\ast}$ of
$\cnum(\!(y_q^{-1})\!)$
induced by $\Phi^{\ast}(y)=y+\varrho$.
We consider the difference $\cnum(\!(y_q^{-1})\!)$-module
$\LLhat^{\lambda}_q(\ell)$ given by
$\LLhat^{\lambda}_q(\ell)=\cnum(\!(y_q^{-1})\!)\,e$
with 
\[
 \Phi^{\ast}(e)
=\Bigl(1+\frac{|\varrho|^2}{4}\Bigr)^{\ell/q}
 y_q^{-\ell}(1+\varrho y^{-1})^{-\ell/q}
 \exp\Bigl(
 \frac{\ell}{q}
 G(\varrho y^{-1})
 \Bigr)\cdot e.
\]
Here, $G(x)=1-x^{-1}\log(1+x)$.
Because this example has the following property,
the corresponding
$\nbigo_{\Hhat_{\infty,q,T}}(\ast H_{\infty,q,T})$-modules
are convenient with respect to the pull back and push-forward
with respect to $\nbigt_{T,mT}$.

\begin{lem}
 For any $m\in\seisuu_{\geq 1}$,
we have
\[
 (\Phi^m)^{\ast}(e)=
 \Bigl(1+\frac{|\varrho|^2}{4}\Bigr)^{m\ell/q}
  y_q^{-m\ell}(1+\varrho my^{-1})^{-m\ell/q}
 \exp\Bigl(
 \frac{m\ell}{q}
 G(\varrho m y^{-1})
 \Bigr)\cdot e.
\]
\end{lem}
\pf
We give only a formal argument.
to simplify the description.
We have
\[
 (\Phi^m)^{\ast}e
=
 \Bigl(1+\frac{|\varrho|^2}{4}\Bigr)^{m\ell/q}
 \prod_{j=1}^m(y+j\varrho)^{-\ell/q}
 \cdot
 \exp\Bigl(
 \frac{\ell}{q}
 \sum_{j=0}^{m-1}
 G\bigl(\varrho(y+j\varrho)^{-1}\bigr)
 \Bigr)\cdot e.
\]
We have 
\begin{multline}
 \frac{\ell}{q}
 \sum_{j=0}^{m-1}
 \Bigl(
 1-\frac{y+\varrho}{\varrho}
 \log\Bigl(
 \frac{y+(j+1)\varrho}{y+j\varrho}
 \Bigr)
 \Bigr)
 = \\
 \frac{\ell m}{q}
-\frac{\ell}{q}
 \sum_{j=0}^{m-1}
 \frac{y}{\varrho}
 \log\Bigl(
 \frac{y+(j+1)\varrho}{y+j\varrho}
 \Bigr)
-\frac{\ell}{q}
 \sum_{j=0}^{m-1}
 j\log\Bigl(
 \frac{y+(j+1)\varrho}{y+j\varrho}
 \Bigr)
 = \\
 \frac{\ell m}{q}
-\frac{\ell m}{q}\frac{y}{\varrho m}
 \log\Bigl(\frac{y+m\varrho}{y}\Bigr)
  \\
-\frac{\ell}{q}
 \Bigl(
 \sum_{j=0}^{m-1}(j+1)\log(y+(j+1)\varrho)
-\sum_{j=0}^{m-1}\log(y+(j+1)\varrho)
-\sum_{j=0}^{m-1}j\log(y+j\varrho)
 \Bigr)
\\
=\frac{\ell m}{q}G(\varrho m y^{-1})
-\frac{\ell}{q}
 \Bigl(
 m\log (y+m\varrho)
-\sum_{j=1}^m
 \log(y+j\varrho)
 \Bigr).
\end{multline}
Then, we obtain the claim of the lemma.
\hfill\qed

\chapter{Filtered bundles}
\label{section;20.8.8.21}

In \S\ref{subsection;20.7.31.20},
we introduce the notion of filtered bundles
on locally free
$\nbigo_{\nbigmbar^{\lambda}\setminus Z}
(\ast H^{\lambda}_{\infty})$-modules
of Dirac type.
We also define the stability condition for such filtered bundles.
Moreover, we explain the comparison with
parabolic difference modules.
In \S\ref{subsection;17.12.16.1},
we introduce the notion of filtered bundles
on a ramified covering space $\nbigb^{\lambda}_q$
of a neighbourhood of $H^{\lambda}_{\infty}$.
In \S\ref{subsection;20.7.31.21},
we explain an extension of
a mini-holomorphic bundle
$\nbigb^{\lambda}_q\setminus H^{\lambda}_{\infty,q}$
with a Hermitian metric
to a filtered object
on $(\nbigb^{\lambda}_q,H^{\lambda}_{\infty,q})$.
We introduce the conditions called
norm estimates and strong adaptedness in this context.
In \S\ref{subsection;20.7.31.22},
we study the analytification of
the construction in \S\ref{subsection;17.10.28.20},
which is also a ramified version of the construction in
\S\ref{subsection;21.8.12.32}.
It is our purpose to compare the norm estimate
for $\lambda$-connections
and
the norm estimate for
$\nbigo_{\nbigb^{\lambda}_q}(\ast H^{\lambda}_{\infty,q})$-modules
(Proposition \ref{prop;17.10.14.20}).

\section{Filtered bundles in the global case}
\label{subsection;20.7.31.20}

We use the notation in \S\ref{subsection;17.10.2.1}.
Let $Z$ be a finite subset of $\nbigm^{\lambda}$.
Let $\pi^{\lambda}:\nbigmbar^{\lambda}\lrarr S^1_T$
denote the projection induced by
$(t_1,\beta_1)\longmapsto t_1$.
\index{projection $\pi^{\lambda}$}
We set $\nbigmbar^{\lambda}\langle t_1\rangle:=
(\pi^{\lambda})^{-1}(t_1)$
and $H^{\lambda}_{\infty}\langle t_1\rangle:=
(\pi^{\lambda})^{-1}(t_1)\cap H^{\lambda}_{\infty}$
for any $t_1\in S^1_T$.
\index{space $\nbigmbar^{\lambda}\langle t_1\rangle$}
\index{space $H^{\lambda}_{\infty}\langle t_1\rangle$}
If there is no risk of confusion,
the point
$H^{\lambda}_{\infty}\langle t_1\rangle
\in \nbigmbar^{\lambda}\langle t_1\rangle$
is also denoted by $\infty$.

Let $\nbige$ be a locally free
 $\nbigo_{\nbigmbar^{\lambda}\setminus Z}(\ast H^{\lambda}_{\infty})$-module
of Dirac type.
Let $\nbige_{|\nbigmbar^{\lambda}\langle t_1\rangle\setminus Z}$
denote the pull back of
$\nbige$ by $\nbigmbar^{\lambda}\langle t_1\rangle\setminus Z
\lrarr \nbigmbar^{\lambda}\setminus Z$,
which is naturally a locally free
$\nbigo_{\nbigmbar^{\lambda}\langle t_1\rangle\setminus Z}(\ast\infty)$-module.

\begin{df}
\index{filtered bundle}
\label{df;21.8.16.2}
A filtered bundle over $\nbige$ is 
a family of filtered bundles
$\nbigp_{\ast}(\nbige_{|\nbigmbar^{\lambda}\langle t_1\rangle\setminus Z})$
$(t_1\in S^1_T)$
over $\nbige_{|\nbigmbar^{\lambda}\langle t_1\rangle \setminus Z}$.
(See {\rm\S\ref{subsection;21.8.15.2}}
for the notion of filtered bundles
on punctured Riemann surfaces.)
Such a family is also denoted by
$\nbigp_{\ast}\nbige$.
\index{filtered bundle $\nbigp_{\ast}\nbige$}
 \hfill\qed
\end{df} 

For any $P\in H^{\lambda}_{\infty}$,
let $\nbigo_{\Hhat^{\lambda}_{\infty},P}$
denote the completion of
the local ring
$\nbigo_{\nbigmbar^{\lambda}_{\infty},P}$
with respect to the maximal ideal.
\index{ring $\nbigo_{\Hhat^{\lambda}_{\infty},P}$}
Thus, we obtain the sheaf of algebras
$\nbigo_{\Hhat^{\lambda}_{\infty}}$
on $H^{\lambda}_{\infty}$.
\index{sheaf $\nbigo_{\Hhat^{\lambda}_{\infty}}$}
Let $\Hhat^{\lambda}_{\infty}$
denote the ringed space
obtained as the topological space
$H^{\lambda}_{\infty}$
with the sheaf of algebras
$\nbigo_{\Hhat^{\lambda}_{\infty}}$.
\index{ringed space $H^{\lambda}_{\infty}$}
We also set
$\nbigo_{\Hhat^{\lambda}_{\infty}}(\ast H^{\lambda}_{\infty})_P:=
 \nbigo_{\nbigmbar^{\lambda}}(\ast H^{\lambda}_{\infty})_P
 \otimes_{\nbigo_{\nbigmbar^{\lambda},P}}
 \nbigo_{\Hhat^{\lambda}_{\infty},P}$,
 and we obtain the sheaf of algebras
 $\nbigo_{\Hhat^{\lambda}_{\infty}}(\ast H^{\lambda}_{\infty})$.
\index{sheaf $\nbigo_{\Hhat^{\lambda}_{\infty}}(\ast H^{\lambda}_{\infty})$}
 
\begin{rem}
The ringed space
$H^{\lambda}_{\infty}$
is naturally isomorphic to
the formal space
$\Hhat_{\infty,1}$ in {\rm\S\ref{subsection;20.7.21.1}}
by the change of the variables
from
$(t_1,\beta_1)$ to
$(t,y)$.
The sheaf of algebras
$\nbigo_{\Hhat^{\lambda}_{\infty}}(\ast H^{\lambda}_{\infty})$
is naturally isomorphic to
$\nbigo_{\Hhat_{\infty,1}}(\ast H_{\infty,1})$. 
\hfill\qed
\end{rem}

Let $k^{\lambda}:
 \Hhat^{\lambda}_{\infty}
 \lrarr \nbigmbar^{\lambda}\setminus Z$
denote the morphism of ringed spaces
obtained as the inclusion
$H^{\lambda}_{\infty}\lrarr
 \nbigmbar^{\lambda}\setminus Z$
with the natural morphism of
sheaves of algebras
$\nbigo_{\nbigmbar^{\lambda}\setminus Z|H^{\lambda}_{\infty}}
\lrarr
\nbigo_{\Hhat^{\lambda}_{\infty}}$.
\index{morphism $k^{\lambda}$}
We set
$\nbige_{|\Hhat^{\lambda}_{\infty}}:=
(k^{\lambda})^{\ast}(\nbige)$,
which is naturally
a locally free
$\nbigo_{\Hhat^{\lambda}_{\infty}}(\ast H^{\lambda}_{\infty})$-module.
A filtered bundle over $\nbige$
naturally induces
a filtered bundle over $\nbige_{|\Hhat^{\lambda}_{\infty}}$
in the sense of Definition \ref{df;20.7.21.2}.

\begin{df}
\label{df;17.10.28.10}
$\nbigp_{\ast}\nbige$ is called good
if the induced filtered bundle over
$\nbige_{|\Hhat^{\lambda}_{\infty}}$ is good
in the sense of Definition {\rm\ref{df;20.7.28.20}}.
\index{good filtered bundle}
\hfill\qed
\end{df}

\begin{df}
A (good) filtered bundle of Dirac type on
$(\nbigmbar^{\lambda};Z,H^{\lambda}_{\infty})$
is a locally free
$\nbigo_{\nbigmbar^{\lambda}\setminus Z}(\ast H^{\lambda}_{\infty})$-module
of Dirac type $\nbige$
equipped with a (good) filtered bundle
$\nbigp_{\ast}\nbige$. 
\index{good filtered bundle of Dirac type on
$(\nbigmbar^{\lambda};Z,H^{\lambda}_{\infty})$}
\hfill\qed
\end{df}

\begin{rem}
 \index{tensor product (filtered bundle)}
 \index{direct sum (filtered bundle)}
 \index{inner homomorphism (filtered bundle)}
 \index{dual (filtered bundle)}
For filtered bundles $\nbigp_{\ast}(\nbige_i)$ $(i=1,2)$
over 
locally free
$\nbigo_{\nbigmbar^{\lambda}\setminus Z}
(\ast H^{\lambda}_{\infty})$-modules
 $\nbige_i$,
we naturally define the filtered bundles
 $\nbigp_{\ast}(\nbige_1)\oplus\nbigp_{\ast}(\nbige_2)$
 over $\nbige_1\oplus\nbige_2$,
 $\nbigp_{\ast}(\nbige_1)\otimes\nbigp_{\ast}(\nbige_2)$
 over $\nbige_1\otimes\nbige_2$
 and
 $\nhom(\nbigp_{\ast}\nbige_1,\nbigp_{\ast}\nbige_2)$
 over $\nhom(\nbige_1,\nbige_2)$
 by applying the constructions in
 {\rm \S\ref{subsection;20.8.1.40}}
 to
 $\nbigp_{\ast}(\nbige_{i|\nbigmbar^{\lambda}\langle t_1\rangle\setminus Z})$.
If $\nbigp_{\ast}\nbige_i$ $(i=1,2)$ are good,
then 
$\nbigp_{\ast}(\nbige_1)\oplus\nbigp_{\ast}(\nbige_2)$,
$\nbigp_{\ast}(\nbige_1)\otimes\nbigp_{\ast}(\nbige_2)$
and 
$\nhom(\nbigp_{\ast}\nbige_1,\nbigp_{\ast}\nbige_2)$
are also good,
which follows from  Lemma {\rm\ref{lem;21.8.26.50}}.
Similarly,
for a filtered bundle $\nbigp_{\ast}\nbige$ over $\nbige$,
we naturally obtain the dual $(\nbigp_{\ast}\nbige)^{\lor}$
over $\nbige^{\lor}$.
If $\nbigp_{\ast}\nbige$ is good,
then $(\nbigp_{\ast}\nbige)^{\lor}$ is also good.
\hfill\qed 
\end{rem}

For any $t_1\in\real$,
we set
$\Mbar^{\lambda}\langle t_1\rangle:=
\{t_1\}\times\proj^1_{\beta_1}
\subset \Mbar^{\lambda}$.
\index{space $\Mbar^{\lambda}\langle t_1\rangle$}
The projection $\varpi^{\lambda}:\Mbar^{\lambda}\lrarr\nbigmbar^{\lambda}$
induces an isomorphism
$\Mbar^{\lambda}\langle \ttilde_1\rangle\simeq
\nbigmbar^{\lambda}\langle t_1\rangle$
for any lift $\ttilde_1\in\real$ of $t_1\in S^1_T$.
We set $Z^{\cov}:=(\varpi^{\lambda})^{-1}(Z)$.
We obtain the locally free
$\nbigo_{\Mbar^{\lambda}\setminus Z^{\cov}}
(\ast H^{\lambda\cov}_{\infty})$-module
$\nbige^{\cov}:=(\varpi^{\lambda})^{-1}(\nbige)$.
\index{sheaf $\nbige^{\cov}$}
As the restriction, we obtain
the locally free
$\nbigo_{\Mbar^{\lambda}\langle t_1\rangle\setminus Z^{\cov}}
(\ast H^{\lambda\cov}_{\infty})$-modules
$\nbige^{\cov}_{|\Mbar^{\lambda}\langle t_1\rangle\setminus Z^{\cov}}$
$(t_1\in\real)$.
There exists the induced filtered bundles
$\nbigp_{\ast}(\nbige^{\cov}_{|\Mbar^{\lambda}\langle t_1\rangle
 \setminus Z^{\cov}})$
$(t_1\in\real)$. 
By choosing $t_1^0,t_1^1\in\real$,
we obtain the isomorphism
\begin{equation}
\label{eq;21.8.16.3}
 \nbige^{\cov}_{|\{t_1^0\}\times\inftyhat}
 \simeq
 \nbige^{\cov}_{|\{t_1^1\}\times\inftyhat}.
\end{equation}
In Definition {\rm\ref{df;21.8.16.2}},
we do not impose any relation between
the filtered bundles
$\nbigp_{\ast}(\nbige^{\cov}
 _{|\{t_1^0\}\times\inftyhat})$
and
$\nbigp_{\ast}(\nbige^{\cov}
 _{|\{t_1^1\}\times\inftyhat})$
under the isomorphism {\rm(\ref{eq;21.8.16.3})}. 
Recall that there exists the slope decomposition
$\nbige_{|\Hhat^{\lambda}_{\infty}}
=\bigoplus_{\omega\in\rnum}
\nbigs_{\omega}(\nbige_{|\Hhat^{\lambda}_{\infty}})$,
which induces the decomposition
$\nbige^{\cov}_{|\Hhat^{\lambda\cov}_{\infty}}
=\bigoplus_{\omega\in\rnum}
\nbigs_{\omega}(\nbige^{\cov}_{|\Hhat^{\lambda}_{\infty}})$.
If the goodness condition in Definition {\rm\ref{df;17.10.28.10}}
is satisfied,
there exist the decompositions
\[
 \nbigp_{\ast}(\nbige^{\cov}_{|\{t_1\}\times\inftyhat})
=\bigoplus_{\omega\in\rnum}
 \nbigp_{\ast}
 \nbigs_{\omega}(\nbige^{\cov}_{|\{t_1\}\times\inftyhat})
 \quad(t_1\in\real),
\]
and there exist the following isomorphisms
as in Lemma {\rm\ref{lem;20.7.20.110}}: 
\begin{equation}
\label{eq;21.9.6.1}
 \nbigp_{a-t_1^0\omega/T}
 \nbigs_{\omega}(\nbige^{\cov}_{|\{t_1^0\}\times\inftyhat})
 \simeq
 \nbigp_{a-t_1^1\omega/T}
 \nbigs_{\omega}(\nbige^{\cov}_{|\{t_1^1\}\times\inftyhat}).
\end{equation}
We shall later discuss more detailed coherence property
in the ramified situation
(See \S\ref{subsection;17.10.10.11}.)

\subsection{Subbundles and quotient bundles}
\label{subsection;18.11.21.3}

Let $\nbige_1\subset\nbige$
be a locally free 
$\nbigo_{\nbigmbar^{\lambda}\setminus Z}(\ast H^{\lambda}_{\infty})$-submodule
of $\nbige$.
We say that $\nbige_1$ is saturated
if $\nbige_2=\nbige/\nbige_1$ is torsion-free.
In that case, we obtain filtered bundles
$\nbigp_{\ast}\nbige_i$ over $\nbige_i$
by applying the constructions in \S\ref{subsection;21.8.26.51}.
Namely, by setting 
\[
 \nbigp_{b}(\nbige_{1|\nbigmbar^{\lambda}\langle t_1\rangle\setminus Z}):=
 \nbigp_b(\nbige_{|\nbigmbar^{\lambda}\langle t_1\rangle\setminus Z })
 \cap
 \nbige_{1|\nbigmbar^{\lambda}\langle t_1\rangle\setminus Z },
\]
we obtain the induced filtered bundle over $\nbige_1$,
which we denote by $\nbigp_{\ast}\nbige_1$.
By setting
\[
 \nbigp_b\bigl(
 \nbige_{2|\nbigmbar^{\lambda}\langle t_1\rangle\setminus Z}
 \bigr)
 :=
 \Image\Bigl(
 \nbigp_b\bigl(
 \nbige_{|\nbigmbar^{\lambda}\langle t_1\rangle\setminus Z}
  \bigr)
\lrarr
 \nbige_{2|\nbigmbar^{\lambda}\langle t_1\rangle\setminus Z}
 \Bigr),
\]
we obtain a filtered bundle 
$\nbigp_{\ast}\nbige_2$ over $\nbige_2$.
We obtain the following lemma from
Lemma \ref{lem;21.8.26.40}.
\begin{lem}
\label{lem;21.9.17.20}
If $\nbigp_{\ast}\nbige$ is good,
$\nbigp_{\ast}\nbige_i$ are also good.
\hfill\qed
\end{lem}
\index{subbundle (filtered bundle)}
\index{quotient bundle (filtered bundle)}

\subsection{Degree and slope}

Let $\nbigp_{\ast}\nbige$
be a good filtered bundle of Dirac type on
$(\nbigmbar^{\lambda};Z,H^{\lambda}_{\infty})$.
Note that
the numbers
$\deg(\nbigp_{\ast}\nbige_{|\nbigmbar^{\lambda}\langle t_1\rangle})$
are well defined for any $t_1\in S^1_T\setminus \pi^{\lambda}(Z)$,
and they induce a function on $S^1_T\setminus\pi^{\lambda}(Z)$.
We have the slope decomposition
$\nbige_{|\Hhat_{\infty}^{\lambda}}
=\bigoplus\nbigs_{\omega}(\nbige_{|\Hhat_{\infty}^{\lambda}})$.
We set
$r(\omega):=\rank\nbigs_{\omega}(\nbige_{|\Hhat_{\infty}^{\lambda}})$.

\begin{lem}
\label{lem;21.9.6.4}
 Let $0\leq \ttilde_1^0<\ttilde_1^1<T$ $(i=0,1)$.
Suppose that the induced points $t_1^i\in S^1_T$ are contained in
a connected component of $S^1_T\setminus\pi^{\lambda}(Z)$.
Then, the following holds.
\[
 \deg(\nbigp_{\ast}\nbige_{|\nbigmbar^{\lambda}\langle t_1^1\rangle})
-\deg(\nbigp_{\ast}\nbige_{|\nbigmbar^{\lambda}\langle t_1^0\rangle})
 =\sum_{\omega\in\rnum} \frac{r(\omega)\omega}{T}
 (\ttilde_1^1-\ttilde_1^0).
\]  
\end{lem}
\pf
It follows from (\ref{eq;21.9.6.1})
and the definition
$\deg(\nbigp_{\ast}\nbige_{|\nbigmbar^{\lambda}\langle t_1\rangle})$
(see the formula (\ref{eq;21.9.6.3})).
\hfill\qed

\vspace{.1in}
By Lemma \ref{lem;21.9.6.4},
on each connected component of
$S^1_T\setminus\pi^{\lambda}(Z)$,
the function is affine with respect to
the local coordinate $t_1$,
and hence it is integrable over $S^1_T$.
The degree of $\nbigp_{\ast}\nbige$
is defined as
\index{degree $\deg(\nbigp_{\ast}\nbige)$}
\[
 \deg(\nbigp_{\ast}\nbige):=
 \frac{1}{T}
 \int_{0}^T
 \deg\bigl(
  \nbigp_{\ast}\nbige_{|\nbigmbar^{\lambda}\langle t_1\rangle\setminus Z}
   \bigr)
  \,dt_1.
\]
We define
$\mu(\nbigp_{\ast}\nbige):=
 \deg(\nbigp_{\ast}\nbige)/\rank\nbige$.
\index{slope $\mu(\nbigp_{\ast}\nbige)$}

Let us rewrite $\deg(\nbigp_{\ast}\nbige)$.
For each $P\in Z$,
we obtain the tuple of integers
$(\ell_1(P),\ldots,\ell_{\rank\nbige}(P))$
as in Lemma {\rm\ref{lem;21.9.6.2}}.
\begin{lem}
\label{lem;21.9.6.5}
The following holds.
\[
 \deg(\nbigp_{\ast}\nbige)=
 \lim_{\substack{\epsilon>0\\ \epsilon\to 0}}
 \deg\bigl(
 \nbigp_{\ast}\nbige_{|\nbigmbar^{\lambda}\langle -\epsilon\rangle}
  \bigr)
 +\sum_{\omega\in\rnum}
 \frac{r(\omega)\cdot\omega}{2}
 +\sum_{P\in Z}
 \Bigl(1-\frac{\pi^{\lambda}(P)}{T}\Bigr)
 \sum_{k=1}^{\rank(\nbige)}\ell_k(P).
\]
Here, we naturally regard $0\leq \pi^{\lambda}(P)<T$. 
\end{lem}
\pf
Let $0\leq \ttilde_1^0<T$ such that
the induced point $t_1^0\in S^1_T$ is contained in $\pi^{\lambda}(Z)$.
Then, the following holds:
\[
 \lim_{\substack{\epsilon>0\\ \epsilon\to 0}}
 \Bigl(
 \deg(\nbigp_{\ast}\nbige_{|\nbigmbar^{\lambda}\langle t_1^0+\epsilon\rangle})
-\deg(\nbigp_{\ast}\nbige_{|\nbigmbar^{\lambda}\langle t_1^0-\epsilon\rangle})
 \Bigr)
=\sum_{P\in \nbigmbar^{\lambda}\langle t_1\rangle}
 \sum_{k=1}^{\rank\nbige} \ell_k(P).
\]
Then, the claim follows from Lemma \ref{lem;21.9.6.4}.
\hfill\qed

\begin{lem}
\label{lem;20.8.7.20}
 Let $\nbigp_{\ast}\nbige$ be a filtered bundle
 over $\nbige$.
 Let $0\lrarr \nbige_1\lrarr \nbige\lrarr\nbige_2\lrarr 0$
be an exact sequence of 
locally free
$\nbigo_{\nbigmbar^{\lambda}\setminus Z}(\ast H^{\lambda}_{\infty})$-modules
of Dirac type.
Let $\nbigp_{\ast}\nbige_i$
be the induced filtered bundles over $\nbige_i$.
Then,
we obtain
$\deg(\nbigp_{\ast}\nbige_1)
+\deg(\nbigp_{\ast}\nbige_2)
=\deg(\nbigp_{\ast}\nbige)$.
As a result,
we obtain
\[
 \mu(\nbigp_{\ast}\nbige)-\mu(\nbigp_{\ast}\nbige_1)
=\frac{r_2}{r_1}
 \bigl(
 \mu(\nbigp_{\ast}\nbige_2)-\mu(\nbigp_{\ast}\nbige)
 \bigr).
\]
\end{lem}
\pf
For each $t_1\in S^1_T\setminus \pi^{\lambda}(Z)$,
we obtain
$\deg(\nbigp_{\ast}\nbige_{1|\nbigmbar^{\lambda}\langle t_1\rangle})
+\deg(\nbigp_{\ast}\nbige_{2|\nbigmbar^{\lambda}\langle t_1\rangle})
=\deg(\nbigp_{\ast}\nbige_{|\nbigmbar^{\lambda}\langle t_1\rangle})$.
Then, the claim of the lemma follows.
\hfill\qed

\begin{cor}
\label{cor;20.8.8.1}
 Let $\nbigp_{\ast}\nbige$
 and $\nbigp_{\ast}\nbige_i$ $(i=1,2)$
 be as in Lemma {\rm\ref{lem;20.8.7.20}}.
 Then, one of the following holds:
 (i)
 $\mu(\nbigp_{\ast}\nbige_1)<
 \mu(\nbigp_{\ast}\nbige)<
 \mu(\nbigp_{\ast}\nbige_2)$,
 (ii)
  $\mu(\nbigp_{\ast}\nbige_1)>
 \mu(\nbigp_{\ast}\nbige)>
 \mu(\nbigp_{\ast}\nbige_2)$,
 (iii)
   $\mu(\nbigp_{\ast}\nbige_1)=
 \mu(\nbigp_{\ast}\nbige)=
 \mu(\nbigp_{\ast}\nbige_2)$.
\hfill\qed
\end{cor}

\subsection{Stability condition}

\begin{df}
\index{stable(filtered bundle)}
\index{semistable(filtered bundle)}
\index{polystable(filtered bundle)}
\label{df;21.8.22.20}
 Let $\nbigp_{\ast}\nbige$
 be a good filtered bundle of Dirac type
 over $(\nbigmbar^{\lambda};Z,H^{\lambda}_{\infty})$.
 It is called stable (semistable)
 if the following holds
 for any saturated
 $\nbigo_{\nbigmbar^{\lambda}\setminus Z}(\ast H^{\lambda}_{\infty})$-submodule
 $\nbige_1$ of Dirac type
 with $0<\rank\nbige_1<\rank\nbige$:
\[
 \mu(\nbigp_{\ast}\nbige_1)
<(\leq)\mu(\nbigp_{\ast}\nbige)
\]
It is called polystable
if it is the direct sum
of stable filtered bundles
$\nbigp_{\ast}(\nbige)=\bigoplus\nbigp_{\ast}\nbige_i$
with $\mu(\nbigp_{\ast}\nbige_i)=\mu(\nbigp_{\ast}\nbige)$.
\hfill\qed
\end{df}

We have some standard lemmas.
\begin{lem}
\label{lem;21.9.17.30}
 Let $f:\nbigp_{\ast}\nbige_1
\lrarr \nbigp_{\ast}\nbige_2$
 be a morphism of stable good filtered bundles of Dirac type
 on $(\nbigmbar^{\lambda};Z,H^{\lambda}_{\infty})$
with $\mu(\nbigp_{\ast}\nbige_1)=\mu(\nbigp_{\ast}\nbige_2)$.
Then,
$f$ is an isomorphism or $0$.
\end{lem}
\pf
Let $\nbige_3$ be the image of $f$.
Suppose that $\nbige_3\neq 0$.
We obtain filtered bundles
$\nbigp^{(i)}_{\ast}\nbige_3$ $(i=1,2)$ over $\nbige_3$
induced by $\nbigp_{\ast}\nbige_i$.
By Corollary \ref{cor;20.8.8.1},
we obtain
\[
\mu(\nbigp_{\ast}\nbige_1)\leq
\mu(\nbigp^{(1)}_{\ast}\nbige_3)
\leq
\mu(\nbigp^{(2)}_{\ast}\nbige_3)
\leq\mu(\nbigp_{\ast}\nbige_2).
\]
Because $\mu(\nbigp_{\ast}\nbige_1)=\mu(\nbigp_{\ast}\nbige_2)$,
the equalities hold.
Hence, we obtain that
$\nbigp_{\ast}\nbige_1\simeq
\nbigp^{(1)}_{\ast}\nbige_3
\simeq\nbigp^{(2)}_{\ast}\nbige_3
\simeq\nbigp_{\ast}\nbige_2$.
\hfill\qed

\begin{cor}
\label{cor;21.9.17.31}
 Let $\nbigp_{\ast}\nbige$ be a stable good filtered
 bundle of Dirac type on
 $(\nbigmbar^{\lambda};Z,H_{\lambda}^{\infty})$.
 Any endomorphism $f$ of $\nbigp_{\ast}\nbige$ is
 the multiplication of a complex number.
\hfill\qed
\end{cor}

As a result,
we obtain the following proposition.
\begin{prop}
\label{prop;21.9.17.32}
 Any polystable good filtered bundle of Dirac type
 $\nbigp_{\ast}\nbige$ on
 $(\nbigmbar^{\lambda};Z,H_{\lambda}^{\infty})$
 has a unique decomposition
\begin{equation}
\label{eq;20.7.20.140}
 \nbigp_{\ast}\nbige
=\bigoplus_{i\in\Lambda}
 \nbigp_{\ast}\nbige_i\otimes_{\cnum} U_i,
\end{equation}
where
$\nbigp_{\ast}\nbige_i$
are stable good filtered bundles of Dirac type
on $(\nbigmbar^{\lambda};Z,H^{\lambda}_{\infty})$
such that 
$\nbigp_{\ast}\nbige_i
\not\simeq
 \nbigp_{\ast}\nbige_j$ $(i\neq j)$,
and $U_i$ denotes $\cnum$-vector spaces.
\hfill\qed 
\end{prop}

\subsection{Good filtered bundles of Dirac type
and parabolic difference modules}

\subsubsection{Polystable parabolic difference modules}

Let $\Phi^{\ast}$ be an automorphism of the field $\cnum(\beta_1)$
defined by $\Phi^{\ast}(f)(\beta_1)=f(\beta_1+2\sqrt{-1}\lambda T)$.
Let $\vecV$ be a $2\sqrt{-1}\lambda T$-difference module
as in Definition \ref{df;20.7.29.1}.
We set
$\vecV_{|\inftyhat}:=
\cnum(\!(y^{-1})\!)\otimes_{\cnum(y)}\vecV$.
We recall that
a filtered bundle $\nbigp_{\ast}(\vecV_{|\inftyhat})$
over $\vecV_{|\inftyhat}$
is called a good parabolic structure of
the difference module $\vecV$ at $\infty$
(Definition \ref{df;20.7.30.1})
if $(\nbigp_{\ast}\vecV_{|\inftyhat},\Phi^{\ast})$
is a good filtered difference module in the sense of 
Definition {\rm\ref{df;20.7.28.11}}.
\index{good parabolic structure at infinity}
We also recall that
a parabolic structure of the difference module $\vecV$
consists of
a parabolic structure at finite place
$(V,m,(\vectau_x,\vecL_x)_{x\in\cnum})$
(see Definition {\rm\ref{df;17.12.1.10}})
and a good parabolic structure at infinity
$\nbigp_{\ast}(\vecV_{|\inftyhat})$,
and it is called a parabolic difference module
(Definition \ref{df;21.9.17.1}).
\index{parabolic difference module}
Such a tuple
$(\vecV,V,m,(\vectau_x,\vecL_x)_{x\in\cnum},
\nbigp_{\ast}(\vecV_{|\inftyhat}))$
is denoted by $\vecV_{\ast}$.

Let $\vecV'\subset\vecV$ be a difference submodule.
By setting
$V'=V\cap\vecV'$,
$L'_{x,i}=L_{x,i}\cap
(\vecV'\otimes_{\cnum(\beta_1)}\cnum(\!(\beta_1-x)\!))$,
and
$\nbigp_a(\vecV'_{|\inftyhat})=
\nbigp_a(\vecV_{|\inftyhat})\cap \vecV'_{|\inftyhat}$ $(a\in\real)$,
we obtain a parabolic structure of $\vecV'$.
(See Lemma \ref{lem;21.8.26.14}.)
The induced parabolic difference module is denoted by $\vecV'_{\ast}$.

Let $\vecV''=\vecV/\vecV'$ be the quotient difference module.
Let $V''$ denote the image of $V\lrarr \vecV''$,
which is isomorphic to $V/V'$.
Let $L''_{x,i}$ denote the image of
$L_{x,i}\lrarr \vecV''\otimes_{\cnum(\beta_1)}\cnum(\!(\beta_1-x)\!)$,
which is isomorphic to $L_{x,i}/L'_{x,i}$.
Let $\nbigp_a(\vecV''_{|\inftyhat})$ $(a\in\real)$
denote the image of
$\nbigp_a(\vecV_{|\inftyhat})\lrarr
\vecV''_{|\inftyhat}$,
which is isomorphic to
$\nbigp_a(\vecV_{|\inftyhat})\big/
\nbigp_a(\vecV'_{|\inftyhat})$.
Note that it is a good parabolic structure of $\vecV''$
at infinity (see Lemma \ref{lem;21.8.26.14}).
Thus, we obtain a parabolic structure of
the difference module $\vecV''$.
The induced parabolic difference module is denoted by
$\vecV''_{\ast}$.

For a parabolic difference module $\vecV_{\ast}$,
we define the degree $\deg(\vecV_{\ast})$
and the slope $\mu(\vecV_{\ast})$
by the formulas (\ref{eq;17.12.4.4})
and (\ref{eq;21.9.17.2}), respectively.
\begin{lem}
\label{lem;21.9.17.21}
Let $\vecV'\subset\vecV$ be a difference submodule,
and we set $\vecV''=\vecV/\vecV'$.
Then, we have
$\deg(\vecV_{\ast})=\deg(\vecV'_{\ast})+\deg(\vecV''_{\ast})$.
\end{lem}
\pf
It is enough to check that
$\deg(L_{i,x},L_{i-1,x})=
\deg(L'_{i,x},L'_{i-1,x})
+\deg(L''_{i,x},L''_{i-1,x})$
for any $x\in\cnum$
and $i=1,\ldots,m(x)$.
Note that for any 
$N\subset L_{i,x}\cap L_{i-1,x}$
such that $(L_{i,x}\cap L_{i-1,x})/N$ is torsion,
we have
\[
 \deg(L_{i,x},L_{i-1,x})=
 \dim_{\cnum}(L_{i,x}/N)
-\dim_{\cnum} (L_{i-1,x}/N).
\]
We set $N':=N\cap (V'\otimes\cnum(\!(\beta_1-x)\!))$
and $N''=N/N'$.
Then, 
$N'\subset L'_{i,x}\cap L'_{i-1,x}$
and
$N''\subset L''_{i,x}\cap L''_{i-1,x}$,
and the quotients are torsion.
Hence, we obtain
\begin{multline}
 \deg(L'_{i,x},L'_{i-1,x})
+\deg(L''_{i,x},L''_{i-1,x})
 =\\
 \dim(L'_{i,x}/N')-\dim(L'_{i-1,x}/N')
 +\dim(L''_{i,x}/N'')-\dim(L''_{i-1,x}/N'')
 \\
=\dim(L_{i,x}/N)-\dim(L_{i-1,x}/N)
=\deg(L_{i,x},L_{i-1,x}).
\end{multline}
Thus, we are done.
\hfill\qed

\vspace{.1in}
As a consequence of Lemma \ref{lem;21.9.17.21},
a similar statement to Corollary \ref{cor;20.8.8.1} holds.
The stability, semistability and polystability
conditions are defined in the standard way
as in Definition \ref{df;21.9.17.3}.
\index{stable (difference module)}
\index{semistable (difference module)}
\index{polystable (difference module)}
We obtain standard statements similar to 
Lemma \ref{lem;21.9.17.30},
Corollary \ref{cor;21.9.17.31}
and Proposition \ref{prop;21.9.17.32}.

\subsubsection{Equivalence}

Let
$\nbigp_{\ast}\nbige$
be a good filtered bundle of Dirac type on
$(\nbigmbar^{\lambda};Z,H_{\lambda}^{\infty})$.
We use the notation in \S\ref{subsection;17.10.5.100}.
From a locally free
$\nbigo_{\nbigmbar^{\lambda}\setminus Z}
(\ast H^{\lambda}_{\infty})$-module
$\nbige$,
we obtained a difference module $\vecV(\nbige)$
equipped with a parabolic structure at finite place
\begin{equation}
\label{eq;21.9.17.13}
(V(\nbige),m_Z,(\vectau_{Z,x},\vecL_{Z,x}(\nbige))_{x\in\cnum}).
\end{equation}
There exist natural isomorphisms
\begin{equation}
\label{eq;21.9.17.4}
 \vecV(\nbige)_{|\inftyhat}
\simeq
 \nbige^{\cov}_{|\{0\}\times\inftyhat}
\simeq
 \nbige_{|\Hhat^{\lambda}_{\infty}\langle 0\rangle}.
\end{equation}
Hence, we obtain the filtered bundle
$\nbigp_{\ast}(\vecV(\nbige)_{|\inftyhat}):=
\nbigp_{\ast}(\nbige_{|\Hhat^{\lambda}_{\infty}\langle 0\rangle})$
under the isomorphism (\ref{eq;21.9.17.4}).
By the definition of good filtered bundle over $\nbige$,
$\nbigp_{\ast}(\vecV(\nbige)_{|\inftyhat})$
is a good parabolic structure of
$\vecV(\nbige)$ at infinity.
Let $\vecV_{\ast}(\nbigp_{\ast}\nbige)$
denote the difference module
$\vecV(\nbige)$
with the parabolic structure at finite place (\ref{eq;21.9.17.13})
and the good parabolic structure at infinity
$\nbigp_{\ast}(\vecV(\nbige)_{|\inftyhat})$.
By Lemma \ref{lem;21.9.6.5},
we obtain
\begin{equation}
\label{eq;21.9.17.14}
 \deg(\nbigp_{\ast}\nbige)
=\deg(\vecV_{\ast}(\nbigp_{\ast}\nbige)).
\end{equation}

The following proposition is obvious
by Lemma \ref{lem;17.12.4.2} and Proposition \ref{prop;21.9.17.41}.
\begin{prop}
\label{prop;21.9.17.50}
The above construction induces
an equivalence between
good filtered bundles $\nbige$ of Dirac type on
$(\nbigmbar^{\lambda};Z,H^{\lambda}_{\infty})$
and parabolic difference modules
$\vecV_{\ast}=(\vecV,V,m_Z,(\vectau_{Z,x},\vecL_{Z,x})_{x\in\cnum},
\nbigp_{\ast}(\vecV_{|\inftyhat}))$,
where $m_Z$ (resp. $\vectau_{Z,x}$ $(x\in\cnum)$)
are determined by $Z$ as in
{\rm(\ref{eq;21.9.17.10})}
(resp. {\rm(\ref{eq;21.9.17.11})} and {\rm(\ref{eq;21.9.17.12})}).
The construction is functorial with respect to
direct sums, tensor products and inner homomorphisms.
\hfill\qed 
\end{prop}

For any locally free
$\nbigo_{\nbigmbar^{\lambda}\setminus Z}(\ast H^{\lambda}_{\infty})$-module
$\nbige'\subset\nbige$
such that $\nbige''=\nbige/\nbige'$ is also locally free,
we obtain the induced difference submodule
$\vecV(\nbige')\subset\vecV(\nbige)$
for which
$\vecV(\nbige)/\vecV(\nbige')\simeq\vecV(\nbige'')$.
By the constructions,
the parabolic difference module
$\vecV(\nbigp_{\ast}\nbige')$
is equal to
the difference module
$\vecV(\nbige')$ equipped with the parabolic structure
induced by $\vecV_{\ast}(\nbigp_{\ast}\nbige)$
and $\vecV(\nbige')\subset\vecV(\nbige)$.
The parabolic difference module
$\vecV_{\ast}(\nbigp_{\ast}\nbige'')$
is equal to
the difference module
$\vecV(\nbige'')$ equipped with the parabolic structure
induced by $\vecV_{\ast}(\nbigp_{\ast}\nbige)$
and $\vecV(\nbige)\lrarr\vecV(\nbige'')$.

Conversely, let $\vecV'\subset\vecV(\nbige)$
be a difference submodule.
We set $\vecV''=\vecV(\nbige)/\vecV'$.
We obtain the induced parabolic difference modules
$\vecV'_{\ast}$ and $\vecV''_{\ast}$.
We have good filtered bundle $\nbigp_{\ast}\nbige'$
(resp. $\nbigp_{\ast}\nbige''$)
on $(\nbigmbar^{\lambda};Z,H^{\lambda}_{\infty})$
such that
$\vecV_{\ast}(\nbigp_{\ast}\nbige')\simeq\vecV'_{\ast}$
(resp. $\vecV_{\ast}(\nbigp_{\ast}\nbige'')\simeq\vecV''_{\ast}$).
By the construction, it is easy to see that
there exists an exact sequence
$0\lrarr\nbige'\lrarr\nbige\lrarr\nbige''\lrarr 0$
which induces the exact sequences
\[
0\lrarr
\nbigp_{\ast}(\nbige'_{|\nbigmbar^{\lambda}\langle t_1\rangle\setminus Z})
\lrarr
\nbigp_{\ast}(\nbige_{|\nbigmbar^{\lambda}\langle t_1\rangle\setminus Z})
\lrarr
\nbigp_{\ast}(\nbige''_{|\nbigmbar^{\lambda}\langle t_1\rangle\setminus Z})
\lrarr 0.
\]

As a result, we also obtain the following proposition.
\begin{prop}
\label{prop;21.9.17.51}
The construction preserves
the stability, semistability and polystability conditions.
\hfill\qed
\end{prop}

\section{Filtered bundles on ramified coverings}
\label{subsection;17.12.16.1}

We generalize the notion of filtered bundles
to the ramified situation.

\subsection{The case $\lambda=0$}
\label{subsection;21.8.19.1}

We use the notation in \S\ref{subsection;17.10.2.1}.
We set
$Y^{0,\cov}:=\Mbar^0\setminus (\real\times\{0\})
=\real_t\times(\proj_w^1\setminus\{0\})$.
\index{space $Y^{0\cov}$}
We also set
$Y^{0,\cov,\ast}:=Y^{0,\cov}\setminus H^0_{\infty}
=\real_t\times\cnum_w^{\ast}$.
\index{space $Y^{0\cov\ast}$}
They are preserved by the $\seisuu$-action $\kappa$.
The quotient spaces are denoted by
$Y^{0}$ and $Y^{0\ast}$.
They are naturally open subsets of
$\nbigmbar^0$.
\index{space $Y^0$}
\index{space $Y^{0\ast}$}

For any positive integer $q$,
we fix a $q$-th root $w_q$ of $w$
such that $w_{qm}^{m}=w_q$ for any $m\in\seisuu_{\geq 1}$.
\index{variable $w_q$}
There exists the ramified covering
$\real_t\times\proj^1_{w_q}\lrarr \Mbar^0=\real_t\times\proj^1_w$.
Let $Y^{0\cov}_q$, $Y^{0\cov\ast}_q$
and $H^{0\cov}_{\infty,q}$
denote the inverse images of $Y^{0\cov}$, $Y^{0\cov\ast}$
and $H^{0\cov}_{\infty}$, respectively.
\index{space $Y^{0\cov}_q$}
\index{space $Y^{0\cov\ast}_q$}
\index{space $H^{0\cov}_{\infty,q}$}
Let $\kappa_q$ be  the $\seisuu$-action
on $\real_t\times\proj^1_{w_q}$
defined by
$\kappa_{q,n}(t,w_q)=(t+nT,w_q)$,
which preserves the subsets
$Y_q^0$, $Y_q^{0\ast}$ and $H^{0}_{\infty,q}$.
\index{action $\kappa_q$}
The quotient spaces are denoted by
$Y_q^0$, $Y_q^{0\ast}$ and $H^{0}_{\infty,q}$,
respectively.
\index{space $Y_q^0$}
\index{space $Y_q^{0\ast}$}
\index{space $H^0_{\infty,q}$}
The projection $Y_q^{0\cov}\lrarr Y_q^0$
is denoted by $\varpi_q^0$.
\index{projection $\varpi_q^0$}
Note that
$Y_q^{0\cov}$ and $Y_q^{0}$
are naturally equipped with mini-complex structures
induced by the mini-complex structure
of $\real_t\times\proj^1_{w_q}$.
Clearly,
$Y_1^{0\cov}$
and
$Y_1^{0}$
are equal to
$Y^{0\cov}$
and $Y^0$,
respectively.

Let $\Psi_q^0:Y^0_q\lrarr \proj^1_{w_q}\setminus\{0\}$
denote the projection.
\index{map $\Psi_q^0$}
The restriction
$Y^{0\ast}_q\lrarr \cnum_{w_q}\setminus\{0\}$
is denoted by
$\Psi_q^0$ or simply $\Psi_q$.
\index{map $\Psi_q$}

\subsection{Ramified coverings for general $\lambda$}
\label{subsection;20.7.22.100}

Let $Y^{\lambda\,\cov}$ denote 
the open subset of $\Mbar^{\lambda}$
defined as follows:
\index{space $Y^{\lambda\cov}$} 
\[
 Y^{\lambda,\cov}:=
 \Mbar^{\lambda}\setminus
 \bigl\{
 (t_1,\beta_1)\,\big|\,
 \beta_1-2\sqrt{-1}\lambda t_1=0
 \bigr\}.
\] 
We also put
$Y^{\lambda,\cov\ast}:=
 Y^{\lambda,\cov}\setminus 
H^{\lambda,\cov}_{\infty}$.
\index{space $Y^{\lambda\cov\ast}$}
The $\seisuu$-action $\kappa$ on $\Mbar^{\lambda}$
preserves
$Y^{\lambda\,\cov}$
and $Y^{\lambda,\cov\ast}$.
The quotient spaces are denoted by
$Y^{\lambda}$ and $Y^{\lambda\ast}$,
respectively,
which are naturally open subsets of
$\nbigmbar^{\lambda}$.
\index{space $Y^{\lambda}$}
\index{space $Y^{\lambda\ast}$}
Clearly,
$Y^{\lambda\ast}=Y^{\lambda}\setminus H^{\lambda}_{\infty}$.
Note that
$Y^{\lambda,\cov,\ast}=Y^{0,\cov,\ast}$
and 
$Y^{\lambda,\ast}=Y^{0,\ast}$
under the $C^{\infty}$-identifications
$M^{\lambda}=M^0$
and $\nbigm^{\lambda}=\nbigm^0$.
Recall the relation
$w=(1+|\lambda|^2)^{-1}(\beta_1-2\sqrt{-1}\lambda t_1)$.
By using the morphism
$\Psi^{\lambda}:\Mbar^{\lambda}\lrarr \proj^1_w$
in \S\ref{subsection;21.8.13.21}
and the projection
$\varpi^{\lambda}:\Mbar^{\lambda}\lrarr\nbigmbar^{\lambda}$,
we have
$Y^{\lambda}=(\Psi^{\lambda})^{-1}(\proj^1_{w}\setminus\{0\})$
and
$Y^{\lambda\cov}=
(\varpi^{\lambda})^{-1}(Y^{\lambda})$.

\vspace{.1in}

Let us construct ramified coverings
of $Y^{\lambda\cov}$. 
Let $\iota_{\lambda}:\real_u\times\proj^1_x
 \simeq
 \real_{t_1}\times\proj^1_{\beta_1}$
 denote the diffeomorphism
defined by
$\iota_{\lambda}(u,x):=
 (u,x+2\sqrt{-1}\lambda u)$.
\index{map $\iota_{\lambda}$}
\index{variable $(u,x)$}
We have
$\iota_{\lambda}^{-1}(t_1,\beta_1)
=(t_1,\beta_1-2\sqrt{-1}\lambda t_1)$.
Let $x_q$ be a $q$-th root of $x$,
and let $\real_{u}\times\proj^1_{x_q}
 \lrarr \real_{t_1}\times\proj^1_{\beta_1}$
be the $C^{\infty}$-map
obtained as the composite of 
the ramified covering
$(u,x_q)\longmapsto (u,x_q^q)$
and the diffeomorphism $\iota_{\lambda}$.
\index{variable $x_{q}$}
Let $Y_q^{\lambda,\cov}$,
$Y_q^{\lambda,\cov\ast}$
and
$H_{\infty,q}^{\lambda,\cov}$
be the inverse image of
$Y^{\lambda,\cov}$,
$Y^{\lambda,\cov\ast}$,
and $H_{\infty}^{\lambda,\cov}$,
respectively.
\index{space $Y_q^{\lambda\cov\ast}$}
\index{space $Y_q^{\lambda\cov}$}
\index{space $H_{\infty}^{\lambda,\cov}$}
Let $\kappa_q$  denote the $\seisuu$-action
on $\real_u\times\proj^1_{x_q}$
defined by 
$\kappa_{q,n}(u,x_q)=(u+nT,x_q)$.
\index{action $\kappa_q$}
It induces the $\seisuu$-action on
$Y_q^{\lambda,\cov}$,
$Y_q^{\lambda,\cov\ast}$,
and $H_{\infty,q}^{\lambda,\cov}$.
The quotient spaces are denoted by
$Y_q^{\lambda}$,
$Y_q^{\lambda\ast}$,
and $H_{\infty,q}^{\lambda}$.
\index{space $Y_q^{\lambda}$}
\index{space $Y_q^{\lambda\ast}$}
\index{space $H_{\infty,q}^{\lambda}$}
The projection $Y_q^{\lambda\cov}\lrarr Y_q^{\lambda}$
is denoted by $\varpi_q^{\lambda}$.
\index{projection $\varpi_q^{\lambda}$}
The function $x_q$ is well defined on
$Y_q^{\lambda\ast}$,
which is equal to
$(1+|\lambda|^2)^{1/q}w_q$.

Because $Y_q^{\lambda\cov\ast}\lrarr Y^{\lambda\cov\ast}$
and $Y_q^{\lambda\ast}\lrarr Y^{\lambda\ast}$
are covering spaces,
we obtain the induced mini-complex structures
on $Y_q^{\lambda\cov\ast}$
and $Y_q^{\lambda\ast}$.

\begin{lem}
\label{lem;21.8.21.11}
The mini-complex structures on
$Y_q^{\lambda\cov\ast}$ and $Y_q^{\lambda\ast}$
uniquely extend to mini-complex structures
on $Y_q^{\lambda\cov}$ and $Y_q^{\lambda}$,
respectively.
\end{lem}
\pf
The uniqueness is clear.
Let $P$ be any point of $H^{\lambda\,\cov}_{\infty,q}$.
There exists a neighbourhood $U_P$ of $P$
in $Y_q^{\lambda\cov\ast}$
such that
\[
 \beta_{1,q}^{-1}:=
 x_q^{-1}(1+2\sqrt{-1}\lambda ux^{-1})^{-1/q}
\]
is a well defined function on $U_P$.
\index{variable $\beta_{1,q}^{-1}$}
It is a $q$-th root of 
$\beta_1^{-1}$.
We obtain the mini-complex coordinate system
$\bigl(
 t_1,
 \beta_{1,q}^{-1}
 \bigr)$
on $U_P$,
which is compatible with the mini-complex structure
of $U_P\cap Y_q^{\lambda\,\cov,\ast}$.
By varying $P\in H^{\lambda\cov}_{\infty,q}$,
we obtain a mini-complex structure of $Y^{\lambda,\cov}_{q}$.
It induces a mini-complex structure on $Y^{\lambda}_q$
which is the extension of the mini-complex structure
on $Y^{\lambda\ast}_q$.
\hfill\qed

\vspace{.1in}

We emphasize
$Y_q^{\lambda\cov\ast}=Y_q^{0\cov\ast}$
and 
$Y_q^{\lambda\ast}=Y_q^{0\ast}$
as Riemannian manifolds.
We also note that
$Y_q^{\lambda}$ (resp. $Y_q^{\lambda\cov}$)
is the fiber product of
$Y^{\lambda}$ (resp. $Y^{\lambda\cov}$)
and $\proj^1_{w_q}\setminus\{0\}$
over $\proj^1_w\setminus \{0\}$.

Let $\Psi^{\lambda}_q:
Y^{\lambda}_q\lrarr \proj^1_{w_q}\setminus\{0\}$
denote the projection
induced by
$(t_1,\beta_{1,q})\longmapsto
(1+|\lambda|^2)^{-1/q}(\beta_{1,q}^q-2\sqrt{-1}\lambda t_1)^{1/q}$.
It is proper.
The restriction 
$Y^{\lambda\ast}_q\lrarr \cnum_{w_q}\setminus\{0\}$
is also denoted by $\Psi_q$
because it is identified with
the projection
$\Psi_q:Y^{0\ast}_q\lrarr \cnum_{w_q}\setminus\{0\}$.
\index{map $\Psi_q^{\lambda}$}
\index{map $\Psi_q$}

\vspace{.1in}

For any $P\in H^{\lambda\cov}_{\infty,q}$,
let $\nbigo_{\Hhat^{\lambda\,\cov}_{\infty,q},P}$
denote the completion
of the local ring
$\nbigo_{Y^{\lambda\cov}_q,P}$
with respect to the maximal ideal.
\index{ring $\nbigo_{\Hhat^{\lambda\,\cov}_{\infty,q},P}$}
Thus, we obtain the sheaf of algebras
$\nbigo_{\Hhat^{\lambda\,\cov}_{\infty,q}}$
on $H^{\lambda\,\cov}_{\infty,q}$.
\index{sheaf $\nbigo_{\Hhat^{\lambda\,\cov}_{\infty,q}}$}
The ringed space
$H^{\lambda\,\cov}_{\infty,q}$
with $\nbigo_{\Hhat^{\lambda\cov}_{\infty,q}}$
is denoted by $\Hhat^{\lambda\cov}_{\infty,q}$.
\index{ringed space $\Hhat^{\lambda\cov}_{\infty,q}$.}
Similarly,
for any $P\in H^{\lambda}_{\infty,q}$,
let $\nbigo_{\Hhat^{\lambda}_{\infty,q},P}$
denote the completion
of the local ring
$\nbigo_{Y^{\lambda}_q,P}$
with respect to the maximal ideal.
\index{ring $\nbigo_{\Hhat^{\lambda}_{\infty,q},P}$}
We obtain the sheaf of algebras
$\nbigo_{\Hhat^{\lambda}_{\infty,q}}$
on $H^{\lambda}_{\infty,q}$,
and the ringed space
$\Hhat^{\lambda}_{\infty,q}$.
\index{sheaf $\nbigo_{\Hhat^{\lambda}_{\infty,q}}$}
\index{ringed space $\Hhat^{\lambda}_{\infty,q}$.}
The ringed space
$\Hhat^{\lambda\cov}_{\infty,q}$
is equipped with the induced $\seisuu$-action $\kappa$,
and $\Hhat^{\lambda}_{\infty,q}$
is the quotient space.

For any $P\in H^{\lambda}_{\infty,q}$,
we set $\nbigo_{\Hhat^{\lambda}_{\infty,q}}
(\ast H^{\lambda}_{\infty,q})_P
:=\nbigo_{Y^{\lambda}_q}(\ast H^{\lambda}_{\infty,q})_P
\otimes_{\nbigo_{Y^{\lambda}_q,P}}
\nbigo_{\Hhat^{\lambda}_{\infty,q},P}$.
Thus, we obtain the sheaf of algebras
$\nbigo_{\Hhat^{\lambda}_{\infty,q}}(\ast H^{\lambda}_{\infty,q})$
on $H^{\lambda}_{\infty,q}$.
\index{sheaf $\nbigo_{\Hhat^{\lambda}_{\infty,q}}(\ast H^{\lambda}_{\infty,q})$}
Similarly,
we obtain the sheaf of algebras
$\nbigo_{\Hhat^{\lambda\cov}_{\infty,q}}
(\ast H^{\lambda\cov}_{\infty,q})$
as the formal completion of
$\nbigo_{Y^{\lambda\cov}_q}(\ast H^{\lambda\cov}_{\infty,q})$
along $H^{\lambda\cov}_{\infty,q}$.
\index{sheaf $\nbigo_{\Hhat^{\lambda\cov}_{\infty,q}}
(\ast H^{\lambda\cov}_{\infty,q})$}

\begin{rem}
\label{rem;21.8.19.3}
The ringed space
$\Hhat^{\lambda\cov}_{\infty,q}$
with the $\seisuu$-action $\kappa$
is isomorphic to
$\Hhat^{\cov}_{\infty,q}$ with the $\seisuu$-action
in {\rm\S\ref{subsection;20.7.21.1}}
 by the change of the variables
 $(t_1,\beta_{1,q}^{-1},x^{-1}_q)=(t,y_q^{-1},x^{-1}_q)$.
The ringed space
$\Hhat^{\lambda}_{\infty,q}$
is isomorphic to
$\Hhat_{\infty,q}$ with the $\seisuu$-action
in {\rm\S\ref{subsection;20.7.21.1}}.
\hfill\qed
\end{rem}

Let $k^{\lambda}_q:
\Hhat^{\lambda}_{\infty,q}\lrarr
Y^{\lambda}_q$
denote the morphism of ringed spaces
obtained as the inclusion
$H^{\lambda}_{\infty,q}
\lrarr Y^{\lambda}_q$
and the natural morphism of sheaves of algebras
$\nbigo_{Y^{\lambda}_q|H^{\lambda}_{\infty,q}}
\lrarr \nbigo_{\Hhat^{\lambda}_{\infty,q}}$.
\index{morphism $k^{\lambda}_q$}
For any $\nbigo_{Y^{\lambda}_q}$-module $\nbigf$,
we set
$\nbigf_{|\Hhat^{\lambda}_{\infty,q}}:=
 (k^{\lambda}_q)^{\ast}(\nbigf)$.
Similarly,
we obtain the morphism of ringed spaces
$\Hhat^{\lambda\cov}_{\infty,q}
\lrarr
Y^{\lambda\cov}_{q}$.
The pull back of
$\nbigo_{Y^{\lambda\cov}_q}$-module $\nbigf$
denoted by
$\nbigf_{|\Hhat^{\lambda\cov}_{\infty,q}}$.
\index{restriction $\nbigf_{|\Hhat^{\lambda\cov}_{\infty,q}}$}

For any $p\in q\seisuu_{\geq 1}$,
let 
$\nbigr_{q,p}:
 Y^{\lambda,\cov}_{p}\lrarr
 Y^{\lambda,\cov}_q$
 denote the $\seisuu$-equivariant map induced by
 $x_p\longmapsto x_p^{p/q}$.
\index{morphism $\nbigr_{q,p}$}
Note that the pull back of functions
induces
$\nbigr_{q,p}^{-1}\nbigo_{Y^{\lambda\cov}_q}
\lrarr\nbigo_{Y^{\lambda\cov}_p}$. 
We define the pull back of
$\nbigo_{Y^{\lambda,\cov}_q}(\ast H^{\lambda,\cov}_{\infty,q})$-modules
and the push-forward of
$\nbigo_{Y^{\lambda,\cov}_{p}}(\ast H^{\lambda,\cov}_{\infty,p})$-modules
in natural ways.
Let $\Gal_{q,p}$ denote the Galois group
of the ramified covering $\nbigr_{q,p}$.
We also define the descent of
$\Gal_{q,p}$-equivariant
$\nbigo_{Y^{\lambda,\cov}_{p}}(\ast H^{\lambda,\cov}_{\infty,p})$-modules
in a natural way.
\index{group $\Gal_{q,p}$}

We obtain the induced ramified covering
$\nbigr_{q,p}:
 Y^{\lambda}_{p}\lrarr
 Y^{\lambda}_q$.
We define the pull back of
$\nbigo_{Y^{\lambda}_q}(\ast H^{\lambda}_{\infty,q})$-modules
and the push-forward of
$\nbigo_{Y^{\lambda}_{p}}(\ast H^{\lambda}_{\infty,p})$-modules.
We also define the descent of
$\Gal_{q,p}$-equivariant
$\nbigo_{Y^{\lambda}_{p}}(\ast H^{\lambda}_{\infty,p})$-modules.

\subsection{Filtered bundles}
\label{subsection;20.8.8.30}

Let $\nbigb^{\lambda}_q$
be a neighbourhood of $H^{\lambda}_{\infty,q}$
in $Y^{\lambda}_{q}$.
\index{space $\nbigb^{\lambda}_q$}
Let $\pi^{\lambda}_q:\nbigb^{\lambda}_q\lrarr S^1_T$
be the map induced by
$\pi^{\lambda}_q(u,x_q)=u$.
In terms of the local coordinate systems
$(t_1,\beta_{1,q}^{-1})$,
it is described as
$\pi^{\lambda}_q(t_1,\beta_{1,q}^{-1})=t_1$.
\index{map $\pi^{\lambda}_q$}
We set
$\nbigb^{\lambda}_q\langle t_1\rangle:=
(\pi^{\lambda}_q)^{-1}(t_1)$.
\index{space $\nbigb^{\lambda}_q(\langle t_1\rangle)$}
We also set
$H^{\lambda}_{\infty,q}\langle t_1\rangle:=
(\pi^{\lambda}_q)^{-1}(t_1)\cap H^{\lambda}_{\infty,q}$.
\index{space $H^{\lambda}_{\infty,q}\langle t_1\rangle$}
The point
$H^{\lambda}_{\infty,q}\langle t_1\rangle
\in \nbigb^{\lambda}_q\langle t_1\rangle$
is also denoted by $\infty$
if there is no risk of confusion.
For $t_1\in S^1_T\simeq H^{\lambda}_{\infty,q}$,
let $\Hhat^{\lambda}_{\infty,q}\langle t_1\rangle$
denote the ringed space
$\bigl(
H^{\lambda}_{\infty,q},
\nbigo_{\Hhat^{\lambda}_{\infty,q},t_1}
\bigr)$.
It is the formal completion of
$\nbigb^{\lambda}_q\langle t_1\rangle$
along $H^{\lambda}_{\infty,q}\langle t_1\rangle$.
\index{ringed space $\Hhat^{\lambda}_{\infty,q}\langle t_1\rangle$}

Let $\nbige$ be a locally free
$\nbigo_{\nbigb^{\lambda}_q}(\ast H^{\lambda}_{\infty,q})$-module
of finite rank.
We obtain the locally free
$\nbigo_{\nbigb^{\lambda}_q\langle t_1\rangle}(\ast \infty)$-module
$\nbige_{|\nbigb^{\lambda}_q\langle t_1\rangle}$
as the pull back of $\nbige$
by the inclusion
$\nbigb^{\lambda}_q\langle t_1\rangle
\lrarr \nbigb^{\lambda}_q$.
\index{restriction $\nbige_{|\nbigb^{\lambda}_q\langle t_1\rangle}$}

\begin{df}
\label{df;20.7.24.1}
A family of filtered bundles
$\nbigp_{\ast}\nbige=
\Bigl(
\nbigp_{\ast}\bigl(\nbige_{|\nbigb^{\lambda}_q\langle t_1\rangle}\bigr)\,
 \Big|\,
t_1\in S^1_T
\Bigr)$
is called a filtered bundle over $\nbige$.
\index{filtered bundle}
\hfill\qed
\end{df}

Note that
a filtered bundle
$\nbigp_{\ast}\nbige
=\bigl(
\nbigp_{\ast}(\nbige_{|\nbigb^{\lambda}_1\langle t_1\rangle})
\,\big|\,t_1\in S^1_T
\bigr)$
over $\nbige$ induces
a filtered bundle
$\nbigp_{\ast}\bigl(
\nbige_{|\Hhat^{\lambda}_{\infty,q}}
\bigr)
=\bigl(
\nbigp_{\ast}(\nbige_{|\Hhat^{\lambda}_{\infty,q}\langle t_1\rangle})
\,\big|\,
t_1\in S^1_T
\bigr)$
over
$\nbige_{|\Hhat^{\lambda}_{\infty,q}}$.

\begin{df}
\label{df;21.8.19.2}
A filtered bundle $\nbigp_{\ast}\nbige$
over $\nbige$ is called good
if the induced filtered bundle over
 $\nbige_{|\Hhat^{\lambda}_{\infty,q}}$
 is good
 in the sense of Definition
 {\rm\ref{df;20.7.28.20}}.
\index{good filtered bundle}
\hfill\qed
\end{df}

\begin{df}
A (good) filtered bundle on $(\nbigb^{\lambda}_q,H^{\lambda}_{\infty,q})$
 is a locally free
 $\nbigo_{\nbigb^{\lambda}_q}(\ast H^{\lambda}_{\infty,q})$-module $\nbige$
 of finite rank
 equipped with a (good) filtered bundle $\nbigp_{\ast}\nbige$.
\index{good filtered bundle on $(\nbigb^{\lambda}_q,H^{\lambda}_{\infty,q})$}
 \hfill\qed
\end{df}

Let $p\in q\seisuu_{>0}$.
We set
$\nbigb^{\lambda}_p:=
\nbigr_{q,p}^{-1}(\nbigb^{\lambda}_q)
\subset Y^{\lambda}_p$,
which is a neighbourhood of
$H^{\lambda}_{\infty,p}$.
Let $\nbigp_{\ast}\nbige$
be a filtered bundle on
$(\nbigb^{\lambda}_{q},H^{\lambda}_{\infty,q})$.
We obtain a locally free
$\nbigo_{\nbigb^{\lambda}_p}(\ast H^{\lambda}_{\infty,p})$-module
$\nbigr^{\ast}_{q,p}(\nbige)$.
We also obtain
a family of filtered bundles
$\nbigr^{\ast}_{q,p}\bigl(
\nbigp_{\ast}(\nbige_{|\nbigb^{\lambda}_q\langle t_1\rangle})\bigr)$
$(t_1\in S^1_T)$
over $\nbigr^{\ast}_{q,p}(\nbige)_{|\nbigb^{\lambda}_q\langle t_1\rangle}$.
Thus, we obtain a filtered bundle
$\nbigr^{\ast}_{q,p}(\nbigp_{\ast}\nbige)$
on $(\nbigb^{\lambda}_p,H^{\lambda}_{\infty,p})$.
\index{pull back (filtered bundle)}

Let $\nbigp_{\ast}\nbige'$ be
a filtered bundle on $(\nbigb^{\lambda}_p,H^{\lambda}_{\infty,p})$.
We obtain a locally free
$\nbigo_{\nbigb^{\lambda}_q}(\ast H^{\lambda}_{\infty,q})$-module
$\nbigr_{q,p\ast}(\nbige')$.
Each $\nbigr_{q,p\ast}(\nbige')_{|\nbigb_q^{\lambda}\langle t_1\rangle}$
is equipped with a filtered bundle
$\nbigr_{q,p\ast}\bigl(
\nbigp_{\ast}(\nbige'_{|\nbigb_p^{\lambda}\langle t_1\rangle})
\bigr)$.
Thus, we obtain a filtered bundle
$\nbigr_{q,p\ast}(\nbigp_{\ast}\nbige')$.
\index{push-forward (filtered bundle)}
If $\nbigp_{\ast}\nbige'$ is
$\Gal_{q,p}$-equivariant,
then
the descent of $\nbigp_{\ast}\nbige'$
is defined to be the $\Gal_{q,p}$-invariant part of
$\nbigr_{q,p\ast}(\nbigp_{\ast}\nbige')$.
\index{descent (filtered bundle)}
The following lemma is obvious by definition.
\begin{lem}
\label{lem;21.7.7.3}
Let $\nbigp_{\ast}\nbige$ be a filtered
bundle on $(\nbigb^{\lambda}_{q},H^{\lambda}_{\infty,q})$.
\begin{itemize}
 \item $\nbigr_{q,p}^{\ast}(\nbigp_{\ast}\nbige)$
 is naturally $\Gal_{q,p}$-equivariant,
       and $\nbigp_{\ast}\nbige$ is the descent of
       $\nbigr_{q,p}^{\ast}(\nbigp_{\ast}\nbige)$.
 \item $\nbigp_{\ast}\nbige$ is good
       if and only if
       $\nbigr_{q,p}^{\ast}(\nbigp_{\ast}\nbige)$ is good.
\hfill\qed
\end{itemize}
\end{lem}

\begin{rem}
 \index{tensor product (filtered bundle)}
 \index{direct sum (filtered bundle)}
 \index{inner homomorphism (filtered bundle)}
 \index{dual (filtered bundle)}
 For filtered bundles $\nbigp_{\ast}(\nbige_i)$ $(i=1,2)$
over 
locally free
$\nbigo_{\nbigb_q^{\lambda}}
(\ast H^{\lambda}_{\infty})$-modules
 $\nbige_i$,
we naturally define the filtered bundles
 $\nbigp_{\ast}(\nbige_1)\oplus\nbigp_{\ast}(\nbige_2)$
 over $\nbige_1\oplus\nbige_2$,
 $\nbigp_{\ast}(\nbige_1)\otimes\nbigp_{\ast}(\nbige_2)$
 over $\nbige_1\otimes\nbige_2$
 and
 $\nhom(\nbigp_{\ast}\nbige_1,\nbigp_{\ast}\nbige_2)$
 over $\nhom(\nbige_1,\nbige_2)$
 by applying the constructions in
 {\rm \S\ref{subsection;20.8.1.40}}
 to
 $\nbigp_{\ast}(\nbige_{i|\nbigb^{\lambda}_q(t_1)})$.
 If $\nbigp_{\ast}\nbige_i$ $(i=1,2)$ are good,
then 
 $\nbigp_{\ast}(\nbige_1)\oplus\nbigp_{\ast}(\nbige_2)$,
 $\nbigp_{\ast}(\nbige_1)\otimes\nbigp_{\ast}(\nbige_2)$
 and 
 $\nhom(\nbigp_{\ast}\nbige_1,\nbigp_{\ast}\nbige_2)$
 are also good.
 Similarly,
 for a filtered bundle $\nbigp_{\ast}\nbige$ over $\nbige$,
 we naturally obtain the dual $(\nbigp_{\ast}\nbige)^{\lor}$
 over $\nbige^{\lor}$.
 If $\nbigp_{\ast}\nbige$ is good,
 then $(\nbigp_{\ast}\nbige)^{\lor}$ is also good.
 \hfill\qed 
\end{rem}

\subsection{Local lattices and the weight filtration
   on the graded pieces}
\label{subsection;17.10.10.11}

Take $t^0_1\in S^1_T$
and $0<\epsilon<T/10$.
Let $\nbigb^{\lambda}_q\langle t^{0}_1,\epsilon\rangle$
denote the inverse image 
of $\openopen{t_1^0-\epsilon}{t_1^0+\epsilon}$
by $\pi^{\lambda}_q:\nbigb^{\lambda}_{q}\lrarr S^1_T$.
\index{space $\nbigb^{\lambda}_q\langle t_1^0,\epsilon\rangle$}
We set
\index{space $H^{\lambda}_{\infty,q}\langle t_1^0,\epsilon\rangle$}
\[
 H^{\lambda}_{\infty,q}\langle t_1^0,\epsilon\rangle:=
(\pi^{\lambda}_q)^{-1}\bigl(
\openopen{t^{0}_1-\epsilon}{t^{0}_1+\epsilon}
\bigr)\cap H^{\lambda}_{\infty,q}.
\]
As in \S\ref{subsection;17.10.10.1},
let $\Hhat^{\lambda}_{\infty,q}\langle t^{0}_1,\epsilon\rangle$
denote the ringed space
$H^{\lambda}_{\infty,q}\langle t^{0}_1,\epsilon\rangle$
with $\nbigo_{\Hhat^{\lambda}_{\infty,q}\langle t^{0}_1,\epsilon\rangle}:=
 \nbigo_{\Hhat^{\lambda}_{\infty,q}|
 H^{\lambda}_{\infty,q}\langle t^0_1,\epsilon\rangle}$.
\index{ringed space $\Hhat^{\lambda}_{\infty,q}\langle t^{0}_1,\epsilon\rangle$}
For any $\nbigo_{\nbigb^{\lambda}_q}$-module $\nbigf$,
let $\nbigf_{|\Hhat^{\lambda}_{\infty,q}\langle t^{0}_1,\epsilon\rangle}$
denote the pull back of
$\nbigf$ by the morphism of the ringed spaces
$\Hhat^{\lambda}_{\infty,q}\langle
t^{0}_1,\epsilon\rangle\lrarr \nbigb^{\lambda}_q$.
Similarly,
for any $\nbigo_{\Hhat^{\lambda}_{\infty,q}}$-module $\nbigf$,
let $\nbigf_{|\Hhat^{\lambda}_{\infty,q}\langle t^{0}_1,\epsilon\rangle}$
denote the pull back of
$\nbigf$ by the morphism of the ringed spaces
$\Hhat^{\lambda}_{\infty,q}\langle t^{0}_1,\epsilon\rangle
\lrarr \Hhat^{\lambda}_{\infty,q}$.

Let $\nbigp_{\ast}\nbige$ be a good filtered bundle
on $(\nbigb^{\lambda}_q,H^{\lambda}_{\infty,q})$.
For each $a\in\real$,
there exists the lattice
$\vecP_a^{(t^{0}_1)}\bigl(
\nbige_{|\Hhat^{\lambda}_{\infty,q}\langle t^{0}_1,\epsilon\rangle}
 \bigr)$
of
$\nbige_{|\Hhat^{\lambda}_{\infty,q}\langle t^{0}_1,\epsilon\rangle}$
as in \S\ref{subsection;17.10.10.1}.
Let $k^{\lambda}_q:
\Hhat^{\lambda}_{\infty,q}\langle t^{0}_1,\epsilon\rangle
\lrarr\nbigb^{\lambda}_q\langle t^{0}_1,\epsilon\rangle$
denote the naturally induced morphism.
We obtain the following morphisms of sheaves:
\[
\nbige_{|\nbigb^{\lambda}_q\langle t^{0}_1,\epsilon\rangle}
\stackrel{b_1}{\lrarr}
k^{\lambda}_{q\ast}\Bigl(
\nbige_{|\Hhat^{\lambda}_{\infty,q}\langle t^{0}_1,\epsilon\rangle}
\Bigr)
\stackrel{b_2}{\lrarr}
k^{\lambda}_{q\ast}\Bigl(
\nbige_{|\Hhat^{\lambda}_{\infty,q}\langle t^{0}_1,\epsilon\rangle}
\Big/
\vecP^{(t^{0}_1)}_a\bigl(
\nbige_{|\Hhat^{\lambda}_{\infty,q}\langle t^{0}_1,\epsilon\rangle}
\bigr)
\Bigr).
\]
We define 
$\vecP^{(t^{0}_1)}_a\bigl(
\nbige_{|\nbigb^{\lambda}_q\langle t^{0}_1,\epsilon\rangle}\bigr)$
as the kernel of the epimorphism $b_2\circ b_1$.
\index{lattice $\vecP^{(t^{0}_1)}_a(\nbige_{|\nbigb^{\lambda}_q\langle t^{0}_1,\epsilon\rangle})$}
Thus, we obtain
the locally free
$\nbigo_{\nbigb^{\lambda}_q\langle t^{0}_1,\epsilon\rangle}$-lattice
$\vecP^{(t^{0}_1)}_a\bigl(
\nbige_{|\nbigb^{\lambda}_q\langle t^{0}_1,\epsilon\rangle}\bigr)$
of $\nbige_{|\nbigb^{\lambda}_q\langle t^{0}_1,\epsilon\rangle}$
satisfying the condition
$\vecP_a^{(t^{0}_1)}\bigl(\nbige_{|\nbigb^{\lambda}_q\langle
t^0_1,\epsilon\rangle}\bigr)_{|
 \Hhat^{\lambda}_{\infty,q}(t^{0}_1,\epsilon)}
 =\vecP_a^{(t^{0}_1)}\bigl(
 \nbige_{|\Hhat^{\lambda}_{\infty,q}(t^{0}_1,\epsilon)}
 \bigr)$.
We obtain the following local system on
$H^{\lambda}_{\infty,q}(t^{0}_1,\epsilon)$:
\[
 \vecG^{(t^{0}_1)}_a\bigl(
  \nbige_{|\nbigb^{\lambda}_q\langle t^{0}_1,\epsilon\rangle}
  \bigr)
 :=
 \vecP^{(t^{0}_1)}_a\bigl(
  \nbige_{|\nbigb^{\lambda}_q\langle t^{0}_1,\epsilon\rangle}
  \bigr)
 \Big/
  \vecP^{(t^{0}_1)}_{<a}\bigl(
  \nbige_{|\nbigb^{\lambda}_q\langle t^{0}_1,\epsilon\rangle}
  \bigr).
\]
\index{local system
$\vecG^{(t^{0}_1)}_a\bigl(
\nbige_{|\nbigb^{\lambda}_q\langle t^{0}_1,\epsilon\rangle}
\bigr)$}
There exists the natural isomorphism.
\begin{multline}
\label{eq;20.7,21.10}
 \vecG^{(t^{0}_1)}_a\bigl(
 \nbige_{|\nbigb^{\lambda}_q\langle t^{0}_1,\epsilon\rangle}\bigr)
 _{|H^{\lambda}_{\infty,q}\langle t^{0}_1,\epsilon\rangle}
 \simeq
 \\
 \vecG^{(t^{0}_1)}_a\bigl(
 \nbige_{|\Hhat^{\lambda}_{\infty,q}\langle t^{0}_1,\epsilon\rangle}
  \bigr)
 =\bigoplus_{\omega}
 \vecG^{(t^{0}_1)}_a\Bigl(
  \nbigs_{\omega}\bigl(\nbige_{|\Hhat^{\lambda}_{\infty,q}}\bigr)
 _{|\Hhat^{\lambda}_{\infty,q}\langle t^{0}_1,\epsilon\rangle}
  \Bigr).
\end{multline}
There exists the filtration $W$
on
$\vecG^{(t^{0}_1)}_a\bigl(
\nbige_{|\Hhat^{\lambda}_{\infty,q}\langle
 t^{0}_1,\epsilon\rangle}\bigr)$
as in \S\ref{subsection;17.10.10.1}.
It induces the filtration $W$ on 
$\vecG^{(t^{0}_1)}_a\bigl(
\nbige_{|\nbigb^{\lambda}_q\langle t^{0}_1,\epsilon\rangle}\bigr)$.
\index{monodromy weight filtration $W$}

\subsection{Convenient frame}
\label{subsection;17.10.11.1}

\begin{df}
\index{convenient frame}
A frame $\vecv$ of
 $\vecP^{(t^{0}_1)}_a\bigl(
 \nbige_{|\nbigb^{\lambda}_q\langle t^{0}_1,\epsilon\rangle}
 \bigr)$
is called a convenient frame around $P$
if the following conditions are satisfied.
\begin{itemize}
 \item There exists a decomposition
       $\vecv=\coprod_{a-1<b\leq a}\vecv_{b}$
       such that
       $\vecv_{b}$ is a tuple of sections of
       $\vecP^{(t^{0}_1)}_b\bigl(
       \nbige_{|\nbigb^{\lambda}_q\langle t^{0}_1,\epsilon\rangle}\bigr)$,
       and that $\vecv_{b}$ induces
       a frame of
       $\vecG^{(t^{0}_1)}_b\bigl(
       \nbige_{|\nbigb^{\lambda}_q\langle t^{0}_1,\epsilon\rangle}\bigr)$.
 \item There exists a decomposition
       $\vecv_{b}=\coprod_{\omega\in\rnum}\vecv_{b,\omega}$
       such that
       $\vecv_{b,\omega}$ induces a base of
       \[
  \vecG^{(t^{0}_1)}_b\Bigl(
  \nbigs_{\omega}\bigl(\nbige_{|\Hhat^{\lambda}_{\infty,q}}\bigr)
 _{|\Hhat^{\lambda}_q\langle t^{0}_1,\epsilon\rangle}
       \Bigr)
       \]
       in the decomposition {\rm(\ref{eq;20.7,21.10})}.
 \item There exists a decomposition
       $\vecv_{b,\omega}=\coprod_{k\in\seisuu}\vecv_{b,\omega,k}$
       such that
       $\coprod_{k\leq\ell}\vecv_{b,\omega,k}$
       induces a base of
       \[
  W_{\ell}\vecG^{(t^{0}_1)}_b\Bigl(
  \nbigs_{\omega}\bigl(\nbige_{|\Hhat^{\lambda}_{\infty,q}}\bigr)
 _{|\Hhat^{\lambda}_{\infty,q}\langle t^{0}_1,\epsilon\rangle}
       \Bigr).
       \]
\end{itemize}
It is easy to see the existence of a convenient frame.
\hfill\qed
\end{df}

We set
$\nu(\nbige):=
 \max\bigl\{
  |\omega_1|+|\omega_2|\,\big|\,
  \nbigs_{\omega_i}\nbige_{|\Hhat^{\lambda}_{\infty,q}}\neq 0
  \bigr\}$.
  Suppose that $\epsilon>0$ satisfies
\begin{equation}
\label{eq;20.7.21.32}
 10q\nu(\nbige)\epsilon/T<
 \min\bigl\{
 |a-b|\,\big|\,
 a,b\in
 \Par\bigl(
 \nbigp_{\ast}\nbige_{|\nbigb^{\lambda}_q\langle t^{0}_1\rangle}
 \bigr),\,
 a\neq b
 \bigr\}.
\end{equation}

\begin{lem}
 Let $\vecv=\coprod_{b,\omega,k} \vecv_{b,\omega,k}$
 be a convenient frame of
 $\vecP^{(t^{0}_1)}_a\bigl(
  \nbige_{|\nbigb^{\lambda}\langle t^{0}_1,\epsilon\rangle}\bigr)$.
 Then, the following holds.
 \begin{itemize}
  \item 
 For each $t_1\in H^{\lambda}_{\infty,q}\langle t^{0}_1,\epsilon\rangle$,
 $\vecv_{b,\omega,k}$
	induces a tuple of sections of
	\[
	\nbigp_{b-(t_1-t_1^{0})q\omega/T}
	\bigl(
	\nbigs_{\omega}(\nbige_{|\Hhat^{\lambda}_{\infty,q}})
	_{|\Hhat^{\lambda}_{\infty,q}\langle t_1\rangle}
	\bigr).
	\]
 \item
	The induced tuple $[\vecv_{b,\omega,k}]$
	of elements of
      $\Gr^{\nbigp}_{b-(t_1-t_1^{0})q\omega/T}\bigl(
      \nbigs_{\omega}(\nbige_{|\Hhat^{\lambda}_{\infty,q}})
      _{|\Hhat^{\lambda}_{\infty,q}\langle t_1\rangle}\bigr)$
	are contained in
	\begin{equation}
	 \label{eq;20.7.21.30}
      W_k\Gr^{\nbigp}_{b-(t_1-t_1^{0})q\omega/T}\bigl(
      \nbigs_{\omega}(\nbige_{|\Hhat^{\lambda}_{\infty,q}})
      _{|\Hhat^{\lambda}_{\infty,q}\langle t_1\rangle}\bigr).
	\end{equation}
 \item	$\coprod_{k\leq \ell}[\vecv_{b,\omega,k}]$
	is a base of {\rm(\ref{eq;20.7.21.30})}.
  \end{itemize}
\end{lem}
\pf
It follows from 
Lemma \ref{lem;20.7.21.31}.
\hfill\qed

\section{Hermitian metrics and filtrations}
\label{subsection;20.7.31.21}

\subsection{Prolongation by growth conditions}
\label{subsection;17.10.25.100}

Let $\nbigb^{\lambda}_q$
denote an open neighbourhood of
$H^{\lambda}_{\infty,q}$
in $Y^{\lambda}_q$
as in \S\ref{subsection;20.8.8.30}.
We set
$\nbigb^{\lambda\ast}_q:=
\nbigb^{\lambda}_q\setminus H^{\lambda}_{\infty,q}$.
\index{space $\nbigb^{\lambda\ast}_q$}
Let $(E,\delbar_{E})$
be a mini-holomorphic bundle
on $\nbigb^{\lambda\ast}_q$.
Let $h$ be a Hermitian metric of $E$.
We explain a general procedure
to construct
an $\nbigo_{\nbigb^{\lambda}_q}(\ast H^{\lambda}_{\infty,q})$-module
$\nbigp^h E$,
although it is not necessarily locally free.
\index{sheaf $\nbigp^hE$}

Let $U$ be any open subset 
in $\nbigb^{\lambda}_q$.
Let $\nbigp^h E(U)$
denote the space of mini-holomorphic sections $s$ of
$E$ on $U\setminus H^{\lambda}_{\infty,q}$
such that the following holds.
\begin{itemize}
\item
For any point $P$ of $U\cap H^{\lambda}_{\infty,q}$,
we take a relatively compact neighbourhood $U_P$ 
of $P$ in $\nbigb^{\lambda}_q$.
Then, there exists $N(P)>0$ such that
\begin{equation}
\label{eq;20.8.8.2}
 \bigl|
 s_{|U_P\setminus H^{\lambda}_{\infty,q}}
 \bigr|_h
=O\bigl(|x_q|^{N(P)}\bigr).
\end{equation}
See \S\ref{subsection;20.7.22.100} for
     the function $x_q$.
     We may replace it with $w_q$ or $\beta_{1,q}$
     in (\ref{eq;20.8.8.2}).
\end{itemize}
Thus, we obtain 
an $\nbigo_{\nbigb^{\lambda}_q}(\ast H^{\lambda}_{\infty,q})$-module
$\nbigp^h E$.

For each $t_1\in S^1_T$,
we set
$\nbigb^{\lambda\ast}_q\langle t_1\rangle:=
\nbigb^{\lambda}_q\langle t_1\rangle
\setminus H^{\lambda}_{\infty,q}$.
\index{space $\nbigb^{\lambda\ast}_q\langle t_1\rangle$}
We obtain the holomorphic vector bundle
$(E,\delbar_E)_{|\nbigb^{\lambda\ast}_1\langle t_1\rangle}$
with the Hermitian metric
$h_{|\nbigb^{\lambda\ast}_q\langle t_1\rangle}$.
By applying the construction in \S\ref{subsection;20.8.8.31}
to
$(E,\delbar_E,h)_{|\nbigb^{\lambda\ast}_1\langle t_1\rangle}$
we obtain the families of sheaves
$\nbigp_{\ast}^h(E_{|\nbigb^{\lambda\ast}_q\langle t_1\rangle})
=\bigl(
 \nbigp_b^h(E_{|\nbigb^{\lambda\ast}_q\langle t_1\rangle})
 \,\big|\,
 b\in\real
 \bigr)$.
The tuple
$\bigl(
 \nbigp_{\ast}^h(E_{|\nbigb^{\lambda\ast}_q\langle t_1\rangle})
 \,\big|\,
 t_1\in S^1_T
 \bigr)$
is denoted by
$\nbigp_{\ast}^h(E)$.
\index{filtered object $\nbigp_{\ast}^h(E)$}
 
\begin{rem}
$\nbigp_{\ast}^h(E_{|\nbigb^{\lambda\ast}_q\langle t_1\rangle})$
 are not necessarily filtered bundles.
 \hfill\qed
\end{rem}

\begin{rem}
Set
$\nbigp^h\bigl(E_{|\nbigb^{\lambda\ast}_q(t_1)}\bigr)
:=
 \bigcup_{a\in\real}
 \nbigp_{a}^h(E_{|\nbigb^{\lambda\ast}_q(t_1)})$.
We have the naturally defined morphism
$(\nbigp^hE)_{|\nbigb^{\lambda}_q(t_1)}
 \lrarr
 \nbigp^h(E_{|\nbigb^{\lambda\ast}_q(t_1)})$,
but it is not necessarily surjective.
\hfill\qed
\end{rem}

\begin{rem}
 We shall often use the notation
 $\nbigp E$
 and $\nbigp_{\ast}E$
 instead of
 $\nbigp^hE$ and $\nbigp^h_{\ast}E$,
 respectively,
 if there is no risk of confusion.
\hfill\qed 
\end{rem}

\subsection{Norm estimate for good filtered bundles}
\label{subsection;17.10.15.1}

Let $\nbigp_{\ast}\nbige$
be a good filtered bundle over a locally free 
$\nbigo_{\nbigb^{\lambda}_q}(\ast H^{\lambda}_{\infty,q})$-module
$\nbige$.
Let $(E,\delbar_E)$ be the mini-holomorphic bundle
obtained as the restriction of $\nbige$
to $\nbigb^{\lambda\ast}_q$.

Let $P\in H^{\lambda}_{\infty,q}$.
We set
$t^{0}_1=\pi^{\lambda}_q(P)\in S^1_T\simeq H^{\lambda}_{\infty,q}$.
Take $\epsilon$ satisfying
(\ref{eq;20.7.21.32}).
Let $\vecv=\coprod \vecv_{b,\omega,k}$
be a convenient frame
on $\nbigb^{\lambda}_q\langle t^{0}_1,\epsilon\rangle$
as in \S\ref{subsection;17.10.11.1}.
We set
$\nbigb^{\lambda\ast}_q\langle t^{0}_1,\epsilon\rangle:=
\nbigb^{\lambda}_q\langle t^{0}_1,\epsilon\rangle
\setminus H^{\lambda}_{\infty,q}$.
For $v_i\in \vecv_{b,\omega,k}$,
we set
$\omega(v_i):=\omega$,
$b(v_i)=b$ and $k(v_i)=k$.
Let $h_P$ be a Hermitian metric of 
$E_{|\nbigb^{\lambda\ast}_q\langle t^{0}_1,\epsilon\rangle}$
determined by the following condition:
\[
 h_P(v_i,v_j):=
 \left\{
 \begin{array}{ll}
 |x_q|^{2b(v_i)-2\omega(v_i)q(t_1-t^{0}_1)/T}
 (\log|x_q|)^{k(v_i)}
 & (i=j)\\
 0 & (i\neq j).\\
 \end{array}
 \right.
\]
The following lemma is clear by the construction.
\begin{lem}
Let $h_P'$  be the Hermitian metric of
$E_{|\nbigb^{\lambda\ast}_q\langle t^{0}_1,\epsilon\rangle}$
 constructed in the same way
 from another convenient frame. 
 Then, for any relatively compact neighbourhood $U_P$
 of $P$ in $\nbigb^{\lambda}_q\langle t^{0}_1,\epsilon\rangle$,
 $h_P$ and $h_P'$ are mutually bounded
 on $U_P\setminus H^{\lambda}_{\infty,q}$.
\hfill\qed
\end{lem}

\begin{df}
\label{df;17.10.24.20}
\index{norm estimate}
Let $h$ be a Hermitian metric of $E$.
We say that $(E,\delbar_E,h)$ satisfies the norm estimate
with respect to
$\nbigp_{\ast}\nbige$ if the following holds.
\begin{itemize}
\item
     For any $P\in H^{\lambda}_{\infty,q}$,
     let $t^{0}_1=\pi^{\lambda}_q(P)$,
     and we take $\epsilon>0$ satisfying
     {\rm(\ref{eq;20.7.21.32})}.
     We construct a Hermitian metric
     $h_P$ of
     $E_{|\nbigb^{\lambda\ast}_q\langle t^{0}_1,\epsilon\rangle}$
     from a convenient frame as above.
     Then, for any relatively compact neighbourhood
     $U_P$ of $P$ in
     $\nbigb^{\lambda}_q\langle t^{0}_1,\epsilon\rangle$,
     $h$ and $h_P$ are mutually bounded
     on $U_P\setminus H^{\lambda}_{\infty,q}$.
\hfill\qed
\end{itemize}
\end{df}
If $(E,\delbar_E,h)$ satisfies the norm estimate
with respect to $\nbigp_{\ast}\nbige$,
then we obtain
$\nbigp^h E=\nbige$
and 
 $\nbigp^h(E_{|\nbigb^{\lambda\ast}_q\langle t_1\rangle})
=\nbigp_{\ast}(\nbige_{|\nbigb^{\lambda}_q\langle t_1\rangle})$
$(t_1\in S^1_T)$.

Take $p\in q\seisuu_{\geq 1}$.
We set $\nbigb^{\lambda}_p:=\nbigr_{q,p}^{-1}(\nbigb^{\lambda}_q)$.
We obtain 
$\nbigr_{q,p}^{-1}(E,\delbar_E,h)$
on $\nbigb^{\lambda\ast}_p=
\nbigb^{\lambda}_p\setminus H^{\lambda}_{\infty,p}$
and a good filtered bundle
 $\nbigr_{q,p}^{\ast}(\nbigp_{\ast}\nbige)$
on $(\nbigb^{\lambda}_p,H^{\lambda}_{\infty,p})$. 

\begin{lem}
 \label{lem;17.10.15.2}
$(E,\delbar_E,h)$ satisfies the norm estimate 
with respect to $\nbigp_{\ast}\nbige$
if and only if 
$\nbigr_{q,p}^{-1}(E,\delbar_E,h)$
satisfies the norm estimate
with respect to
$\nbigr_{q,p}^{\ast}(\nbigp_{\ast}\nbige)$.
\end{lem}
\pf
Set $m:=p/q$.
Let $P$ be any point of $H^{\lambda}_{\infty,q}$,
and set $t^{0}_1:=\pi^{\lambda}_q(P)$.
Take $\epsilon>0$ satisfying (\ref{eq;20.7.21.32})
for
$\nbigp_{\ast}(\nbige_{|\nbigb^{\lambda\ast}_q\langle t^0_1\rangle })$
and
$\nbigr_{q,p}^{\ast}
\nbigp_{\ast}(\nbige_{|\nbigb^{\lambda\ast}_q\langle t^0_1\rangle})$.

Let $\vecv=\bigcup \vecv_{b,\omega,k}$
be a convenient frame of $\vecP^{(t^{0}_1)}_a\nbige$
on $\nbigb^{\lambda}_q(t^{0}_1,\epsilon)$.
We construct a Hermitian metric $h_P$ 
of $E_{|\nbigb^{\lambda\ast}_q(t^{0}_1,\epsilon)}$
from $\vecv_P$ as above.

We set $P'=\nbigr_{q,p}^{-1}(P)\in H^{\lambda}_{\infty,p}$.
For each $a-1<b\leq a$,
let $n(b)\in\seisuu$ be determined by
$ma-1< mb+n(b)\leq ma$.
We set
\[
 \vecv'_{b',\omega,k}:=
 \coprod_{mb+n(b)=b'}
 \beta_{1,p}^{n(b)}\cdot
 \nbigr_{q,p}^{\ast}(\vecv_{b,\omega,k}).
\]
Then, 
$\vecv'=\bigcup_{b',\omega,k} \vecv'_{b',\omega,k}$
is a convenient frame of
$\vecP^{(t^{0}_1)}_{ma}(\nbigr_{q,p}^{\ast}\nbige)$
on $\nbigb^{\lambda}_{p}(t^{0}_1,\epsilon)$.
We construct a Hermitian metric $h_{P'}$
of $\nbigr_{q,p}^{-1}(E)
 _{|\nbigb^{\lambda\ast}_p(t^{0}_1,\epsilon)}$
from $\vecv'$ as above.
We can easily observe that
$h_{P'}$ and $\nbigr_{q,p}^{-1}h_{P}$
are mutually bounded.
Then, the claim of the lemma is clear.
\hfill\qed

\subsection{Strong adaptedness}

We introduce a condition called strong adaptedness,
which is weaker than the norm estimate.
Let $\nbigp_{\ast}\nbige$ 
and $(E,\delbar_E)$ be as in \S\ref{subsection;17.10.15.1}.

\begin{df}
\label{df;20.8.8.3}
\index{strongly adapted}
\index{strong adaptedness}
 Let $h$ be a Hermitian metric of $E$.
We say that $h$ is strongly adapted to $\nbigp_{\ast}\nbige$
if the following condition is satisfied.
\begin{itemize}
\item
 For any $P\in H^{\lambda}_{\infty,q}$,
 we take a neighbourhood $U_P$
 and a Hermitian metric $h_{P}$ of 
 $E_{|U_P\setminus H^{\lambda}_{\infty,q}}$
 as in {\rm\S\ref{subsection;17.10.15.1}}.
 Then, for any $\delta>0,$
there exists a  constant $C_{\delta}\geq 1$
 such that 
 $C_{\delta}^{-1}|w_q|^{-\delta}h_P
 \leq
 h
 \leq 
 C_{\delta}|w_q|^{\delta}h_P$.
\hfill\qed
\end{itemize} 
\end{df}

We obtain the following as in Lemma \ref{lem;17.10.15.2}.
\begin{lem}
Take $p\in q\seisuu_{\geq 1}$.
Then, $h$ is strongly adapted to $\nbigp_{\ast}\nbige$
if and only if
$\nbigr_{q,p}^{-1}(h)$ is strongly adapted to
 $\nbigr_{q,p}^{\ast}(\nbigp_{\ast}\nbige)$.
\hfill\qed
\end{lem}

\begin{rem}
The strong adaptedness in Definition {\rm\ref{df;20.8.8.3}}
is stronger than the adaptedness in the context of
 harmonic bundles,
 i.e.,
 it is a locally uniform estimate in the $S^1_T$-direction.
 \hfill\qed
\end{rem}

\section{Comparison with $\lambda$-connections}
\label{subsection;20.7.31.22}

\subsection{Some sheaves on $Y_q^{\lambda\cov}$
and $Y_q^{\lambda}$}

There exists the complex vector field
$\del_{\betabar_1}$ on
$Y_1^{\lambda}$.
By taking the lift with respect to the covering map
$\nbigr_{1,q}:
Y_q^{\lambda\cov\ast}\lrarr Y_1^{\lambda\cov\ast}$,
we obtain the complex vector field
$\del_{\betabar_1}$ on $Y_q^{\lambda\cov\ast}$
such that
$\nbigr_{1,q\ast}\del_{\betabar_1}=\del_{\betabar_1}$.
Locally around any point of $H^{\lambda\cov}_{\infty,q}$,
$\nbigr_{1,q}$ is described as
$(t_1,\beta_{1,q}^{-1})\longmapsto
(t_1,\beta_{1,q}^{-q})$,
and we have
$\del_{\betabar_1}=
\betabar_{1,q}^{-q}(q^{-1}\betabar_{1,q}\del_{\betabar_{1,q}})$.
Hence
$\del_{\betabar_1}$ extends to a complex vector field
on $Y_q^{\lambda\cov}$,
which is also denoted by $\del_{\betabar_1}$.
Similarly, we obtain the vector field
$\del_{t_1}$ on $Y_q^{\lambda\cov}$.
\index{complex vector field $\del_{\betabar_1}$}
\index{complex vector field $\del_{t_1}$}

Let $U$ be any open subset of $Y^{\lambda\cov}_q$.
Let $\nbigk_{Y^{\lambda\cov}_q}(U)$
denote the space of $C^{\infty}$-functions $f$ on $U$
such that $\del_{\betabar_1}(f)=0$.
We obtain a sheaf of algebras
$\nbigk_{Y^{\lambda\cov}_q}$ on $Y^{\lambda\cov}_q$.
\index{sheaf $\nbigk_{Y^{\lambda\cov}_q}$}
We have the induced operator
$\del_{t_1}:\nbigk_{Y^{\lambda\cov}_q}
\lrarr\nbigk_{Y^{\lambda\cov}_q}$.
It is an epimorphism,
and the kernel is $\nbigo_{Y^{\lambda\cov}_q}$.
We set
$\nbigk_{Y^{\lambda\cov}_q}(\ast H^{\lambda}_{\infty,q})
:=
 \nbigk_{Y^{\lambda\cov}_q}
 \otimes_{\nbigo_{Y^{\lambda\cov}_q}}
 \nbigo_{Y^{\lambda\cov}_q}(\ast H^{\lambda\cov}_q)$.
\index{sheaf $\nbigk_{Y^{\lambda\cov}_q}(\ast H^{\lambda}_{\infty,q})$}
We have the induced operator
$\del_{t_1}:\nbigk_{Y^{\lambda\cov}_q}(\ast H^{\lambda\cov}_{\infty,q})
\lrarr\nbigk_{Y^{\lambda\cov}_q}(\ast H^{\lambda}_{\infty,q})$.
It is an epimorphism,
and the kernel is
$\nbigo_{Y^{\lambda\cov}_q}(\ast H^{\lambda\cov}_{\infty,q})$. 

Because
$\nbigk_{Y^{\lambda\cov}_q}$
is naturally $\seisuu$-equivariant,
we obtain the sheaf of algebras
$\nbigk_{Y^{\lambda}_q}$ on $Y^{\lambda}_q$
as the descent.
\index{sheaf $\nbigk_{Y^{\lambda}_q}$}
It is naturally equipped with the operator
$\del_{t_1}:\nbigk_{Y^{\lambda}_q}\lrarr
\nbigk_{Y^{\lambda}_q}$,
which is an epimorphism,
and the kernel is $\nbigo_{Y^{\lambda}_q}$.
Similarly,
$\nbigk_{Y^{\lambda\cov}_q}(\ast H^{\lambda\cov}_{\infty,q})$
is naturally $\seisuu$-equivariant,
we obtain the sheaf of algebras
$\nbigk_{Y^{\lambda}_q}(\ast H^{\lambda}_{\infty,q})$
on $Y^{\lambda}_q$ as the descent.
\index{sheaf $\nbigk_{Y^{\lambda}_q}(\ast H^{\lambda}_{\infty,q})$}
It is naturally equipped with the operator
$\del_{t_1}:\nbigk_{Y^{\lambda}_q}(\ast H^{\lambda}_{\infty,q})
 \lrarr
\nbigk_{Y^{\lambda}_q}(\ast H^{\lambda}_{\infty,q})$,
which is an epimorphism,
and the kernel is $\nbigo_{Y^{\lambda}_q}(\ast H^{\lambda}_{\infty,q})$.
We have
$\nbigk_{Y^{\lambda}_q}(\ast H^{\lambda}_{\infty,q})
=\nbigo_{Y^{\lambda}_q}(\ast H^{\lambda}_{\infty,q})
\otimes_{\nbigo_{Y^{\lambda}_q}}
\nbigk_{Y^{\lambda}_q}$.

Because the complex vector field
$\del_{\betabar_1}$ on $Y^{\lambda\cov}_q$
is invariant with respect to the $\seisuu$-action,
we obtain $\del_{\betabar_1}$ on $Y^{\lambda}_q$.
We may obtain $\nbigk_{Y^{\lambda}_q}$
as the sheaf of $C^{\infty}$-functions $f$
such that $\del_{\betabar_1}f=0$.

\subsection{The induced $\nbigo_{\nbigb_q^{\lambda}}$-modules
from $\lambda$-connections}
\label{subsection;21.8.12.40}

We explain the analytic version of the construction in
\S\ref{subsection;17.10.28.20},
which is also a ramified version of
the construction in \S\ref{subsection;21.8.12.32}.
Note that the variable $x_q$ in this subsection
is related with the variable $w$ 
in \S\ref{subsection;21.8.12.32}
by $w=(1+|\lambda|^2)^{-1}x_q^q$.

There exists the natural map
$\real_u\times\proj^1_{x_q}
\lrarr \proj^1_{x_q}$
induced by $(u,x_q)\longmapsto x_q$.
It induces the proper map
$Y_q^{\lambda\cov}\lrarr \proj^1_{x_q}\setminus\{0\}$.
\index{map $\Psi^{\lambda}_q$}
In terms of the local coordinate systems
$(t_1,\beta_{1,q}^{-1})$ and $x_q^{-1}$,
it is described as
\[
 (t_1,\beta_{1,q}^{-1})\longmapsto
\beta_{1,q}^{-1}(1-2\sqrt{-1}\lambda t_1\beta_{1,q}^{-q})^{-1/q}.
\]
It induces the map
$\Psi^{\lambda}_q:Y_q^{\lambda}\lrarr \proj^1_{x_q}\setminus\{0\}$.
It is equal to the map
$\Psi^{\lambda}_q$ in \S\ref{subsection;20.7.22.100}
under the identification
$x_q=(1+|\lambda|^2)^{1/q}w_q$.
Clearly,
$(\Psi_q^{\lambda})^{-1}(\infty)=H^{\lambda}_{\infty,q}$.
Note that the restriction of
$\Psi_q^{\lambda}$ to $Y_q^{\lambda\ast}$
is equal to
$\Psi_q^{0}$ to $Y_q^{0\,\ast}$
under the $C^{\infty}$-identification
$Y_q^{\lambda\ast}=Y_q^{0\ast}$.
We also note that we obtain the morphism
$(\Psi_q^{\lambda})^{\ast}:
(\Psi_q^{\lambda})^{-1}
\nbigo_{\proj^1_{x_q}\setminus\{0\}}
\lrarr
\nbigk_{Y^{\lambda}_q}$
by the pull back of functions.

Let $U_{x,q}$ be any open subset of
$\proj^1_{x_q}\setminus\{0\}$.
\index{space $U_{x,q}$}
We set
$\nbigb^{\lambda}_{q}:=
(\Psi_q^{\lambda})^{-1}(U_{x,q})$.
\index{space $\nbigb^{\lambda}_{q}$}
Let $\nbigv$ be a locally free
$\nbigo_{U_{x,q}}(\ast\infty)$-module
with a $\lambda$-connection
$\nabla^{\lambda}:
\nbigv\lrarr\nbigv\otimes\Omega^1_{U_{x,q}}$
whose Poincar\'{e} rank is strictly less than $q$.
By the inner product with
$\del_x=q^{-1}x_q^{-q+1}\del_{x_q}$,
we obtain a morphism of sheaves
$\nabla^{\lambda}_x:
\nbigv\lrarr\nbigv$.
We set \index{sheaf $\nbigvtilde^{\infty}$}
\[
 \nbigvtilde^{\infty}:=
 (\Psi_q^{\lambda})^{-1}(\nbigv)
  \otimes_{(\Psi_q^{\lambda})^{-1}\nbigo_{U_{x,q}}}
  \nbigk_{\nbigb^{\lambda}_{q}}
  =
   (\Psi_q^{\lambda})^{-1}(\nbigv)
  \otimes_{(\Psi_q^{\lambda})^{-1}\nbigo_{U_{x,q}}(\ast\infty)}
  \nbigk_{\nbigb^{\lambda}_{q}}(\ast H^{\lambda}_{\infty,q}).
\]
There uniquely exists
$\del_{\nbigvtilde^{\infty},t_1}:
\nbigvtilde^{\infty}\lrarr
\nbigvtilde^{\infty}$
determined by the following condition
for any local sections
$f$ and $s$ of
$\nbigk_{\nbigb^{\lambda}_{q}}$
and $\nbigv$:
\index{operator $\del_{\nbigvtilde^{\infty},t_1}$}
\begin{equation}
\del_{\nbigvtilde^{\infty},t_1}\bigl(
f\cdot(\Psi_q^{\lambda})^{-1}(s)\bigr)
=\del_{t_1}(f)(\Psi^{\lambda}_q)^{-1}(s)
+f(\Psi^{\lambda}_q)^{-1}(-2\sqrt{-1}\nabla^{\lambda}_xs).
\end{equation}
Let $(\Psi^{\lambda}_q)^{\ast}(\nbigv,\nabla^{\lambda})$
denote the
$\nbigo_{\nbigb^{\lambda}_q}(\ast H^{\lambda}_{\infty,q})$-module
obtained as the kernel of
$\del_{\nbigvtilde^{\infty},t_1}:
\nbigvtilde^{\infty}\lrarr\nbigvtilde^{\infty}$.
\index{sheaf $(\Psi^{\lambda}_q)^{\ast}(\nbigv,\nabla^{\lambda})$}

Similarly,
for any lattice $\nbigl$ of $\nbigv$
such that
$\nabla^{\lambda}_x\nbigl\subset\nbigl$,
we set
\index{sheaf $\nbigltilde^{\infty}$}
\[
 \nbigltilde^{\infty}:=
 (\Psi_q^{\lambda})^{-1}(\nbigl)
 \otimes_{(\Psi_q^{\lambda})^{-1}(\nbigo_{U_{x,q}})}
 \nbigk_{\nbigb^{\lambda}_q}.
\]
It is equipped with the induced operator
$\del_{\nbigltilde^{\infty},t_1}:
 \nbigltilde^{\infty}\lrarr\nbigltilde^{\infty}$.
\index{operator $\del_{\nbigltilde^{\infty},t_1}$}
The kernel is denoted by
$(\Psi_q^{\lambda})^{\ast}(\nbigl,\nabla^{\lambda})$
which is an $\nbigo_{\nbigb^{\lambda}_q}$-module.
\index{sheaf $(\Psi_q^{\lambda})^{\ast}(\nbigl,\nabla^{\lambda})$}

Let
$(\nbigv,\nabla^{\lambda})_{|\inftyhat_{x,q}}$
denote the formal completion of
$(\nbigv,\nabla^{\lambda})$
at $\infty$,
i.e.,
$\nbigv_{|\inftyhat_{x,q}}:=
 \nbigv_{\infty}\otimes_{\nbigo_{U_{x,q}}}\cnum[\![x_q^{-1}]\!]$.
\index{formal completion $(\nbigv,\nabla^{\lambda})_{|\inftyhat_{x,q}}$}
By the construction in \S\ref{subsection;17.10.7.2},
we obtain
a locally free
$\nbigo_{\Hhat^{\lambda}_{\infty,q}}(\ast H^{\lambda}_{\infty,q})$-module
$(\Psi^{\lambda}_q)^{\ast}\bigl(
(\nbigv,\nabla^{\lambda})_{|\inftyhat_{x,q}}
\bigr)$.
Similarly, from the completion
$\nbigl_{|\inftyhat_{x,q}}$
of $\nbigl$
which is a $\nabla^{\lambda}_x$-invariant lattice of
$(\nbigv,\nabla^{\lambda})_{|\inftyhat_{x,q}}$,
we obtain a lattice
$(\Psi^{\lambda}_q)^{\ast}\bigl(
(\nbigl,\nabla^{\lambda})_{|\inftyhat_{x,q}}
\bigr)$
of
$(\Psi^{\lambda}_q)^{\ast}\bigl(
(\nbigv,\nabla^{\lambda})_{|\inftyhat_{x,q}}
\bigr)$.

\begin{prop}
\mbox{{}}\label{prop;20.7.21.50}
\begin{itemize}
\item 
$(\Psi^{\lambda}_q)^{\ast}(\nbigv,\nabla^{\lambda})$
is a locally free 
$\nbigo_{\nbigb^{\lambda}_q}(\ast H^{\lambda}_{\infty,q})$-module.
 Moreover, there exists a natural isomorphism:
\begin{equation}
\label{eq;20.7.21.40}
 (\Psi^{\lambda}_q)^{\ast}(\nbigv,\nabla^{\lambda})
 _{|\Hhat^{\lambda}_{\infty,q}}
 \simeq
 (\Psi^{\lambda}_q)^{\ast}
  \bigl(
  (\nbigv,\nabla^{\lambda})_{|\inftyhat_{x,q}}
 \bigr).
 \end{equation}
  \item For any $\nabla^{\lambda}_x$-invariant lattice
	$\nbigl$ of $\nbigv$,
	$(\Psi_q^{\lambda})^{\ast}(\nbigl,\nabla^{\lambda})$
	is a lattice of
	$(\Psi^{\lambda}_q)^{\ast}(\nbigv,\nabla^{\lambda})$.
 Moreover, there exists a natural isomorphism:
\begin{equation}
\label{eq;20.7.21.41}
 (\Psi^{\lambda}_q)^{\ast}(\nbigl,\nabla^{\lambda})
 _{|\Hhat^{\lambda}_{\infty,q}}
 \simeq
 (\Psi^{\lambda}_q)^{\ast}
  \bigl(
  (\nbigl,\nabla^{\lambda})_{|\inftyhat_{x,q}}
 \bigr).
 \end{equation}
 \end{itemize}
 \end{prop}
\pf
It is easy to see that
$(\Psi^{\lambda}_q)^{\ast}(\nbigv,\nabla^{\lambda})_{
|\nbigb^{\lambda}_{q}\setminus
H^{\lambda}_{\infty,q}}$
is a locally free
$\nbigo_{\nbigb^{\lambda}_q\setminus H^{\lambda}_{\infty,q}}$-module.
Let $\nbigl$ be a $\nabla^{\lambda}_x$-invariant lattice
of $\nbigv$.
Take any frame $\vecv$ of $\nbigl$.
We set $r:=\rank(\nbigv)$.
We obtain a section
$A$ of $M_r(\nbigo_{U_{x,q}})$
determined by
$\nabla^{\lambda}_x\vecv=\vecv\cdot A$.
We set $\vecvtilde=(\Psi_q^{\lambda})^{-1}(\vecv)$
which is frame of $\nbigvtilde^{\infty}$.
Let $P$ be any point of $H^{\lambda}_{\infty,q}$.
We take a mini-complex
local coordinate system $(t_1,\beta_{1,q}^{-1})$
around $P$.
We obtain the following local section of
$M_r\bigl(
 \nbigk_{\nbigb^{\lambda}_q}
 \bigr)$ around $P$:
\[
 \Atilde=
 -2\pi\sqrt{-1}
 A\bigl(\beta_{1,q}^{-1}(1-2\sqrt{-1}\lambda t\beta_{1,q}^{-q})^{-1/q}
 \bigr).
\]
We obtain
$\del_{\nbigvtilde^{\infty},t_1}\vecvtilde
 =\vecvtilde\cdot\Atilde$.
There exists a local section
$G$ of $M_r\bigl(\nbigk_{\nbigb^{\lambda}_q}\bigr)$ around $P$
such that
(i) $G$ is invertible, i.e.,
$G^{-1}\in M_r\bigl(\nbigk_{\nbigb^{\lambda}_q}\bigr)$,
(ii) $G^{-1}\del_{t}G=-\Atilde$.
Then,
we obtain
$\del_{\nbigvtilde^{\infty},t_1}(\vecvtilde\cdot G)=0$.
Because
$\vecvtilde\cdot G$
is a local frame of
$(\Psi^{\lambda}_q)^{\ast}(\nbigv,\nabla^{\lambda})$
around $P$,
we obtain that
$(\Psi^{\lambda}_q)^{\ast}(\nbigl,\nabla^{\lambda})$
is a locally free
$\nbigo_{\nbigb^{\lambda}_q}$-module.
Then, the claims are clear.
\hfill\qed

\vspace{.1in}

Let $\nbigp_{\ast}\nbigv$
be a filtered bundle over $\nbigv$
such that each $\nbigp_a\nbigv$
is $\nabla^{\lambda}_x$-invariant.
The family
$(\Psi^{\lambda}_q)^{\ast}(\nbigp_a\nbigv,\nabla^{\lambda})$
$(a\in\real)$
induces a filtered bundle over
$(\Psi^{\lambda}_q)^{\ast}(\nbigv,\nabla^{\lambda})$,
denoted by
$(\Psi^{\lambda}_q)^{\ast}(\nbigp_{\ast}\nbigv,\nabla^{\lambda})$.
We obtain the following lemma
from Proposition \ref{prop;20.7.20.131}
and Proposition \ref{prop;20.7.21.50}.
\begin{prop}
\label{prop;20.7.24.10}
 $(\nbigp_{\ast}\nbigv,\nabla^{\lambda})$
is a good filtered $\lambda$-flat bundle
if and only if
 $(\Psi^{\lambda}_q)^{\ast}(\nbigp_{\ast}\nbigv,\nabla^{\lambda})$
is a good filtered bundle.
\hfill\qed  
\end{prop}

\begin{prop}
Let $(\nbigp_{\ast}\nbigv_i,\nabla^{\lambda})$
be good filtered $\lambda$-flat bundles
on $(U_{x,q},\infty_{x,q})$
whose Poincar\'{e} rank are strictly smaller than $q$.
 There exist the following natural isomorphisms:
\begin{equation}
 \Psi_q^{\ast}\bigl(
 (\nbigp_{\ast}\nbigv_1,\nabla^{\lambda})
 \oplus
  (\nbigp_{\ast}\nbigv_2,\nabla^{\lambda})
  \bigr)
\simeq
  \Psi_q^{\ast}(\nbigp_{\ast}\nbigv_1,\nabla^{\lambda})
\oplus
  \Psi_q^{\ast}(\nbigp_{\ast}\nbigv_2,\nabla^{\lambda}).
\end{equation}
\begin{equation}
 \Psi_q^{\ast}\bigl(
 (\nbigp_{\ast}\nbigv_1,\nabla^{\lambda})
 \otimes
  (\nbigp_{\ast}\nbigv_2,\nabla^{\lambda})
  \bigr)
\simeq
  \Psi_q^{\ast}(\nbigp_{\ast}\nbigv_1,\nabla^{\lambda})
 \otimes
  \Psi_q^{\ast}(\nbigp_{\ast}\nbigv_2,\nabla^{\lambda}).
\end{equation}
\begin{equation}
 \Psi_q^{\ast}\Bigl(
  \nhom\bigl(
  (\nbigp_{\ast}\nbigv_1,\nabla^{\lambda}),
  (\nbigp_{\ast}\nbigv_2,\nabla^{\lambda})
\bigr)
  \Bigr)
  \simeq
  \nhom\Bigl(
  \Psi_q^{\ast}(\nbigp_{\ast}\nbigv_1,\nabla^{\lambda}),
 \Psi_q^{\ast}(\nbigp_{\ast}\nbigv_2,\nabla^{\lambda})
 \Bigr).
\end{equation}
 In particular,
 for a good filtered $\lambda$-flat bundle
 $(\nbigp_{\ast}\nbigv,\nabla^{\lambda})$
 on $(U_{x,q},\infty_{x,q})$ whose Poincar\'e rank is
 strictly smaller than $q$,
 we obtain
 $\Psi_q^{\ast}(\nbigp_{\ast}\nbigv,\nabla^{\lambda})^{\lor}
 \simeq
 \Psi_q^{\ast}\bigl(
 (\nbigp_{\ast}\nbigv,\nabla^{\lambda})^{\lor}
 \bigr)$.
\hfill\qed
\end{prop}

\begin{rem}
The construction is available
in the case $\infty\not\in U_{x,q}$,
and it is the same as the construction
in {\rm\S\ref{subsection;21.8.13.21}}.
\hfill\qed 
\end{rem}

\subsection{Norm estimates for $\lambda$-connections}

Let us recall the norm estimate condition
for $\lambda$-flat bundles
with a Hermitian metric.
Let $U_{x,q}$ denote a neighbourhood
of $\infty$ in $\proj^1_{x_q}$.
\index{space $U_{x,q}$}
We set $U_{x,q}^{\ast}:=U_{x,q}\setminus\{\infty\}$.
\index{space $U^{\ast}_{x,q}$}
Let $(\nbigp_{\ast}\nbigv,\nabla^{\lambda})$
be a good filtered $\lambda$-flat bundle
on $(U_{x,q},\infty)$.
We set
$\Gr^{\nbigp}_a(\nbigv):=\nbigp_a\nbigv/\nbigp_{<a}\nbigv$.
There exists the residue endomorphism
$\Res(\nabla^{\lambda})$ on each $\Gr^{\nbigp}_b(\nbigv)$.
(See \cite[\S2.5.2]{Mochizuki-wild}.)
Let $W$ be the monodromy weight filtration of 
the nilpotent part of $\Res(\nabla^{\lambda})$.
We set
$(V,\nabla^{\lambda}):=
(\nbigv,\nabla^{\lambda})_{|U_{x,q}^{\ast}}$.

A frame $\vecv$ of $\nbigp_a\nbigv$ is called convenient
for $(\nbigp_{\ast}\nbigv,\nabla^{\lambda})$
if the following conditions are satisfied.
\index{convenient frame}
\begin{itemize}
\item
 There exists a decomposition
 $\vecv=\coprod_{a-1<b\leq a}\vecv_b$
     such that $\vecv_{b}$ is a tuple of
     sections of $\nbigp_b\nbigv$
     and induces a base of
     $\Gr^{\nbigp}_b(\nbigv)$.
\item
 There exists a decomposition
 $\vecv_b=\coprod_{k\in\seisuu}\vecv_{b,k}$
 such that
 $\coprod_{k\leq \ell}\vecv_{b,k}$
 induces a frame of 
 $W_{\ell}\Gr^{\nbigp}_b(\nbigv)$.
\end{itemize}
For a convenient frame $\vecv$,
we define a Hermitian metric $h_0$
of $V$ as follows.
For $v_i\in\vecv_{b,k}$,
we put $b(v_i):=b$ and $k(v_i):=k$,
and we set
\[
 h_0(v_i,v_j):=
 \left\{
 \begin{array}{ll}
 |x_q|^{2b(v_i)}(\log|x_q|)^{k(v_i)}
 & (i=j)\\
 0 & (i\neq j).
 \end{array}
 \right.
\]
The following lemma is clear.
\begin{lem}
 Let $h_0'$ be a Hermitian metric of $V$
 obtained from another convenient frame $\vecv'$
 as above.
 Then, for any relatively compact neighbourhood
 $U'$ of $\infty$ in $U_{x,q}$,
 $h_0$ and $h_0'$ are mutually bounded
 on $U'\setminus\{\infty\}$.
 \hfill\qed
\end{lem}

\begin{df}
\index{norm estimate (good filtered $\lambda$-flat bundle)}
Let $h_V$ be a Hermitian metric of $V$.
We say that 
$(V,\nabla^{\lambda},h)$ satisfies the norm estimate 
with respect to $(\nbigp_{\ast}\nbigv,\nabla^{\lambda})$
 if $h_V$ and $h_0$ are mutually bounded,
 where $h_0$ is a Hermitian metric of $V$
 obtained from a convenient frame $\vecv$
 for $(\nbigp_{\ast}\nbigv,\nabla^{\lambda})$.
\hfill\qed
 \end{df}

\subsection{Comparison of the norm estimates}

Let $U_{x,q}$ be a neighbourhood of $\infty$
in $\proj^1_{x_q}\setminus\{0\}$.
We set
$\nbigb^{\lambda}_q:=
(\Psi^{\lambda}_{q})^{-1}(U_{x,q})$.
Let $(\nbigp_{\ast}\nbigv,\nabla^{\lambda})$ be
a good filtered $\lambda$-flat bundle
on $(U_{x,q},\infty)$.
We obtain a good filtered bundle
$(\Psi^{\lambda}_q)^{\ast}(\nbigp_{\ast}\nbigv,\nabla^{\lambda})$
on $(\nbigb^{\lambda}_q,H^{\lambda}_{\infty,q})$.
We set
$U_{x,q}^{\ast}:=U_{x,q}\setminus\{\infty\}$,
$\nbigb^{\lambda\ast}_q:=(\Psi_q^{\lambda})^{-1}(U_{x,q}^{\ast})$
and $(V,\nabla^{\lambda}):=
(\nbigv,\nabla^{\lambda})_{|U_{x,q}^{\ast}}$.

\begin{prop}
\label{prop;17.10.14.20}
Let $h_V$ be a Hermitian metric of $V$
such that
$(V,\nabla^{\lambda},h_V)$ satisfies
the norm estimate
with respect to $(\nbigp_{\ast}\nbigv,\nabla^{\lambda})$.
Suppose that the Poincar\'e rank of
$(\nbigv,\nabla^{\lambda})$ is strictly less than $q$.
Then,
$\Psi_q^{\ast}(\nbigv,\nabla^{\lambda})_{|\nbigb^{\lambda\ast}_q}$
with $\Psi_q^{-1}(h_V)$
satisfies the norm estimate
with respect to the filtered bundle
$(\Psi_q^{\lambda})^{\ast}(\nbigp_{\ast}\nbigv,\nabla^{\lambda})$.
\end{prop}
\pf
Let us study the case where
there exists $\gminia\in x_q\cnum[x_q]$
such that 
$\nabla^{\lambda}-d\gminia\id$ is logarithmic
with respect to $\nbigp_{\ast}\nbigv$.
Note that $\deg_{x_q}\gminia<q$.
Set $r:=\rank(\nbigv)$.
Let $I_r\in M_r(\cnum)$ denote the identity matrix.
Let $\vecv$ be a convenient frame of $\nbigp_a\nbigv$.
There exists a decomposition
$\vecv=\coprod_{a-1<b\leq a}\coprod_{k\in\seisuu}
 \vecv_{b,k}$
such that 
(i) $\coprod_{k\in\seisuu}\vecv_{b,k}$
is a tuple of sections of
$\nbigp_{b}\nbigv$,
and induces a base of
$\Gr^{\nbigp}_b(\nbigv)$,
(ii) $\coprod_{k\leq \ell}\vecv_{b,k}$
induces a base of
$W_{\ell}\Gr^{\nbigp}_b(\nbigv)$.
For $v_i\in \vecv_{b,k}$,
we set $b(v_i)=b$ and $k(v_i)=k$.
Let $h_{V,0}$ be the Hermitian metric determined by
$h_{V,0}(v_i,v_i)=
 |x_q|^{2b(v_i)}(\log|x_q|)^{k(v_i)}$
 and
 $h_{V,0}(v_i,v_j)=0$ $(i\neq j)$.
We set $\htilde_0:=(\Psi_q^{\lambda})^{-1}(h_{V,0})$.

Let $A$ be the matrix valued function
determined by
$\nabla^{\lambda}_x\vecv=
 \vecv\cdot
 \bigl(\del_{x}\gminia I_r+A(x_q^{-1})x_q^{-q}\bigr)$.
Then, $A(x_q^{-1})$ is holomorphic with respect to $x_q^{-1}$.
Set $\vecvtilde:=(\Psi_q^{\lambda})^{-1}(\vecv)$
which is a frame of $\nbigvtilde^{\infty}$.
We obtain
\begin{equation}
 \label{eq;20.7.22.101}
 \del_{\nbigvtilde^{\infty},t_1}\vecvtilde
=\vecvtilde
 (-2\sqrt{-1})
  (\Psi_q^{\lambda})^{-1}
 \Bigl(
\del_x\gminia I_r
+x_q^{-q}A(x_q^{-1})
 \Bigr).
\end{equation}

Take $P\in H^{\lambda}_{\infty,q}$.
Set $t_1^{0}:=\pi^{\lambda}_q(P)$.
The restriction 
$\vecvtilde_{|\nbigb^{\lambda}_q\langle t_1^{0}\rangle}$
is a holomorphic frame of
$(\Psi_q^{\lambda})^{\ast}(\nbigv,\nabla^{\lambda})
 _{|\nbigb^{\lambda}_q\langle t_1^{0}\rangle}$.
 (See \S\ref{subsection;20.8.8.30}
 for $\nbigb^{\lambda}_q\langle t_1^0\rangle$.)
On a neighbourhood $U_P$ of $P$,
there exists the frame $\vecu_P$ of
$(\Psi_q^{\lambda})^{\ast}(\nbigv,\nabla^{\lambda})$
such that
$\vecu_{P|\nbigb^{\lambda}_q\langle t_1^{0}\rangle}=
\vecvtilde_{|\nbigb^{\lambda}_q\langle t_1^{0}\rangle}$.
By Proposition \ref{prop;20.7.22.50},
$\vecu_P$ is a convenient frame
of the filtered bundle
$(\Psi_q^{\lambda})^{\ast}(\nbigp_{\ast}\nbigv,\nabla^{\lambda})$.
Let $h_P$ be the Hermitian metric of
$(\Psi_q^{\lambda})^{\ast}(\nbigv,\nabla^{\lambda})_{|U_P\setminus
 H^{\lambda}_{\infty,q}}$
 determined by
 $h_P(u_i,u_i)=|x_{q}|^{2b(v_i)}(\log|x_q|)^{k(v_i)}$
 and
 $h_P(u_i,u_j)=0$ $(i\neq j)$.

Let $G$ be the $M_r(\cnum)$-valued
$C^{\infty}$-function on $U_P\setminus H^{\lambda}_{\infty,q}$
determined by
$\vecu_P=\vecv_PG$.
By (\ref{eq;20.7.22.101}),
there exists a function
$\gminic$ such that
(i) $\gminic=O(|x_q|^{-1})$,
(ii)
$|G-(1+\gminic)I_r|=O(|x_q|^{-q})$.
Hence, $\htilde_0$ and $h_P$ are mutually bounded,
and we obtain the claim of the proposition in this case.

\vspace{.1in}
Let us study the case that there exists
a finite subset
$\nbigi\subset S(q)=\{\sum_{j=1}^{q-1}\gminia_jx_q^{j}\}$
and an isomorphism
\begin{equation}
\label{eq;20.7.22.110}
F:
(\nbigp_{\ast}\nbigv,\nabla^{\lambda})_{|\inftyhat}
\simeq
 \bigoplus_{\gminia\in\nbigi}
 (\nbigp_{\ast}\nbigv_{\gminia},\nabla_{\gminia}^{\lambda})_{|\inftyhat},
\end{equation}
where $\nabla^{\lambda}_{\gminia}-d\gminia\id$ are logarithmic
with respect to $\nbigp_{\ast}\nbigv_{\gminia}$.
We obtain the induced isomorphism:
\begin{equation}
 \label{eq;17.10.14.10}
(\Psi_q^{\lambda})^{\ast}(F):
 (\Psi_q^{\lambda})^{\ast}
(\nbigp_{\ast}\nbigv,\nabla^{\lambda})_{|\Hhat^{\lambda}_{\infty,q}}
\simeq
 \bigoplus_{\gminia\in\nbigi}
  (\Psi_q^{\lambda})^{\ast}
 (\nbigp_{\ast}\nbigv_{\gminia},\nabla_{\gminia}^{\lambda})
 _{|\Hhat^{\lambda}_{\infty,q}}.
\end{equation}
Let $\nbigc^{\infty}_{U_{x,q}}$
and 
$\nbigc^{\infty}_{\nbigb^{\lambda}_q}$
denote the sheaf of $C^{\infty}$-functions
on $U_{x,q}$ and $\nbigb^{\lambda}_q$,
respectively.
\index{sheaf $\nbigc^{\infty}_{U_{x,q}}$}
\index{sheaf $\nbigc^{\infty}_{\nbigb^{\lambda}_q}$}
There exists an isomorphism of sheaves
\[
 F_{C^{\infty}}:
 \nbigv\otimes_{\nbigo_{U_{x,q}}} \nbigc^{\infty}_{U_{x,q}}
\simeq
 \bigoplus_{\gminia\in\nbigi}
 \nbigv_{\gminia}
 \otimes_{\nbigo_{U_{x,q}}}
 \nbigc^{\infty}_{U_{x,q}}
\]
which induces (\ref{eq;20.7.22.110}).
It induces an isomorphism of sheaves
\[
(\Psi_q^{\lambda})^{\ast}(F_{C^{\infty}}):
(\Psi_q^{\lambda})^{\ast}
(\nbigp_{\ast}\nbigv,\nabla^{\lambda})
 \otimes_{\nbigo_{\nbigb^{\lambda}_q}}
  \nbigc^{\infty}_{\nbigb^{\lambda}_q}
\simeq
 \bigoplus_{\gminia\in\nbigi}
  (\Psi_q^{\lambda})^{\ast}
 (\nbigp_{\ast}\nbigv_{\gminia},\nabla_{\gminia}^{\lambda})
 \otimes_{\nbigo_{\nbigb^{\lambda}_q}}
 \nbigc^{\infty}_{\nbigb^{\lambda}_q},
\]
which induces (\ref{eq;17.10.14.10}).
It is easy to see that
a Hermitian metric $h_V$
of $\nbigv_{|U_{x,q}^{\ast}}$
satisfies the norm estimate for
$(\nbigp_{\ast}\nbigv,\nabla^{\lambda})$
if and only if
the induced Hermitian metric
$F_{C^{\infty}\ast}(h_V)$
of
$\bigoplus \nbigv_{\gminia|U_{x,q}^{\ast}}$
satisfies the norm estimate for
$\bigoplus(\nbigp_{\ast}\nbigv_{\gminia},
\nabla_{\gminia}^{\lambda})$.
Similarly,
it is easy to see that
$(\Psi_q^{\lambda})^{-1}(h_V)$
satisfies the norm estimate for
$(\Psi_q^{\lambda})^{\ast}(\nbigp_{\ast}\nbigv,\nabla^{\lambda})$
if and only if
$(\Psi_q^{\lambda})^{-1}(F_{C^{\infty}\ast}(h_V))$
satisfies the norm estimate for
$\bigoplus
(\Psi_q^{\lambda})^{\ast}(\nbigp_{\ast}\nbigv_{\gminia},
 \nabla_{\gminia}^{\lambda})$.
Hence, we obtain the claim of the proposition in this case.

In the general case,
there exists $p\in q\seisuu_{>0}$
such that
$\varphi_{q,p}^{\ast}(\nbigp_{\ast}\nbigv,\nabla)$
has a decomposition as in (\ref{eq;20.7.22.110}).
For a Hermitian metric $h_V$
of $\nbigv_{|U_{x,q}^{\ast}}$,
it is easy to see that
$h_V$ satisfies the norm estimate
for $(\nbigp_{\ast}\nbigv,\nabla^{\lambda})$
if and only if
$\varphi_{q,p}^{-1}(h_V)$
satisfies the norm estimate for
$\varphi_{q,p}^{\ast}(\nbigp_{\ast}\nbigv,\nabla^{\lambda})$.
Then, by using Lemma \ref{lem;17.10.15.2},
we obtain the claim of the proposition
in the general case.
\hfill\qed

\subsection{Non-integrable case}
\label{subsection;21.8.13.1}

Let us explain the ramified version of the construction
in \S\ref{subsection;21.8.12.11},
where the integrability of a $\lambda$-connection is not assumed.
Let $U_{x,q}$ be an open subset of $\proj^1_{x_q}$
such that $\infty\in U_{x,q}$.

\subsubsection{A general construction}
\label{subsection;21.8.12.41}

Let $(V,\delbar_V)$ be a holomorphic vector bundle
on $U_{x,q}$.
Let
$\DD^{\lambda\,(1,0)}:C^{\infty}(U_{x,q},V)
\lrarr C^{\infty}(U_{x,q},V\otimes\Omega^{1,0}_{U_{x,q}}((q+1)\infty))$
be the differential operator
$\DD^{\lambda\,(1,0)}(fs)
=\lambda\del_{U_{x,q}}(f)s+f\DD^{\lambda\,(1,0)}(s)$
for any $f\in C^{\infty}(U_{x,q})$
and $s\in C^{\infty}(U_{x,q},V)$.
We do not assume $[\delbar_{V},\DD^{\lambda\,(1,0)}]=0$
at this stage.

We set $\Vtilde=(\Psi^{\lambda}_q)^{-1}(V)$ on
$\nbigb_q^{\lambda}:=(\Psi^{\lambda}_q)^{-1}(U_{x,q})$.
We set
$\nbigb_q^{\lambda\ast}:=
\nbigb_q^{\lambda}\setminus H^{\lambda}_{\infty,q}$.
Let $\Vtilde^{\ast}$ denote the restriction of
$\Vtilde$ to $\nbigb_q^{\lambda\ast}$.
In \S\ref{subsection;21.8.12.3},
we constructed
an $S^1_T$-equivariant differential operator
$\delbar_{\Vtilde^{\ast}}:
 C^{\infty}(\nbigb_q^{\lambda\ast},\Vtilde^{\ast})
 \lrarr
 C^{\infty}(\nbigb_q^{\lambda\ast},
 \Vtilde^{\ast}\otimes\Omega^{0,1}_{\nbigb_q^{\lambda\ast}})$
satisfying the mini-complex Leibniz rule.

\begin{lem}
The differential operator
$\delbar_{\Vtilde^{\ast}}$ uniquely extends to
 a differential operator
$\delbar_{\Vtilde}:
 C^{\infty}(\nbigb_q^{\lambda},\Vtilde)
 \lrarr
 C^{\infty}(\nbigb_q^{\lambda},
 \Vtilde\otimes\Omega^{0,1}_{\nbigb_q^{\lambda}})$
satisfying the mini-complex Leibniz rule.
\end{lem}
\pf
We just repeat the argument in the proof of Lemma \ref{lem;21.8.12.33}.
We may assume that $U_{x,q}$ is an open disc around $\infty$.
We have the complex vector fields
$\del_{\wbar}=(1+|\lambda|^2)q^{-1}\xbar_q^{-q+1}\del_{\xbar_q}$
and
$\del_{w}=(1+|\lambda|^2)q^{-1}\xbar_q^{-q+1}\del_{x_q}$ on $U_{x,q}$.
Let $\del_{V,\wbar}$ (resp. $\DDlambda_w$)
denote the differential operator
on $V$ induced by
$\del_{\wbar}$ and $\delbar_V$
(resp. $\del_w$ and $\DD^{\lambda\,(1,0)}$).

Let $\vecv$ be a holomorphic frame of $(V,\delbar_V)$ on $U_{x,q}$.
We have $\del_{V,\wbar}\vecv=0$.
Let $A$ be the matrix valued $C^{\infty}$-function on
$U_{x,q}$ determined by
$\DDlambda_w\vecv=\vecv A$.
The pull back $\vecvtilde:=(\Psi_q^{\lambda})^{-1}\vecv$
is a $C^{\infty}$-frame of $\Vtilde$ on $\nbigb_{q}^{\lambda}$.

We have the projection
$\varpi^{\lambda}_q:Y^{\lambda\cov}_q\lrarr Y^{\lambda}_q$.
For any $P\in H^{\lambda}_{\infty,q}$,
we choose $\Ptilde\in H^{\lambda\cov}_{\infty,q}$
such that $\varpi_q^{\lambda}(\Ptilde)=P$.
Then, $(t_1,\tau_{1,q})=(t_1,\beta_{1,q}^{-1})$ around $\Ptilde$
induces a local mini-complex coordinate system
on a neighbourhood $U_P$ of $P$.
By the formulas (\ref{eq;21.8.12.30}) and (\ref{eq;21.8.12.31}),
on $U_P\setminus H^{\lambda}_{\infty,q}$,
the operators
$\del_{\Vtilde^{\ast},t_1}$
and
$\del_{\Vtilde^{\ast},\taubar_{1,q}}$
are described as follows
with respect to the frame $\vecvtilde$:
\[
 \del_{\Vtilde^{\ast},t_1}\vecvtilde
=\vecvtilde \frac{-2\sqrt{-1}}{1+|\lambda|^2}
 (\Psi^{\lambda})^{\ast}(A),
\quad
 \del_{\Vtilde^{\ast},\taubar_{1,q}}\vecvtilde=0.
\]
Note that $A$ is $C^{\infty}$ with respect to
$x_q^{-1}$ and $\xbar_q^{-1}$,
and that
$(\Psi^{\lambda})^{\ast}(x_q^{-1})
=\tau_{1,q}^q(1-2\sqrt{-1}\lambda t_1\tau_{1,q}^q)^{-1}$
is $C^{\infty}$ around $P$.
Hence, $\del_{\Vtilde^{\ast},t_1}$
and $\del_{\Vtilde^{\ast},\taubar_1}$
uniquely extend to operators 
$\del_{\Vtilde,t_1}$
and $\del_{\Vtilde,\taubar_{1,q}}$
on $C^{\infty}(\nbigb_q^{\lambda},\Vtilde)$,
which implies the claim of the lemma.
\hfill\qed

\begin{lem}
\label{lem;21.8.12.34}
The above construction induces an equivalence between
the following objects.
\begin{itemize}
 \item Holomorphic vector bundles
       $(V,\delbar_V)$ on $U_{x,q}$ equipped with
       a differential operator
      $\DD^{\lambda\,(1,0)}:C^{\infty}(U_{x,q},V)
\lrarr C^{\infty}(U_{x,q},V\otimes\Omega^{1,0}_{U_{x,q}}((q+1)\infty))$
such that
$\DD^{\lambda\,(1,0)}(fs)
=\lambda\del_{U_{x,q}}(f)s+f\DD^{\lambda\,(1,0)}(s)$
for any $f\in C^{\infty}(U_{x,q})$
and $s\in C^{\infty}(U_{x,q},V)$.
 \item $S^1_T$-equivariant vector bundles $\Vtilde$ on $\nbigb_q^{\lambda}$
       equipped with an $S^1_T$-equivariant linear differential operator
       $\delbar_{\Vtilde}:
       C^{\infty}(\nbigb_q^{\lambda},\Vtilde)
       \lrarr
       C^{\infty}(\nbigb_q^{\lambda},
       \Vtilde\otimes\Omega^{0,1}_{\nbigb_q^{\lambda}})$
       satisfying the mini-complex Leibniz rule.
\end{itemize}
\end{lem}
\pf
By repeating the argument in the proof of Lemma \ref{lem;21.8.12.10},
we indicate the inverse construction.
There exists a $C^{\infty}$-vector bundle $V$ on $U_{x,q}$
equipped with an $S^1_T$-equivariant isomorphism
$\Vtilde\simeq (\Psi_q^{\lambda})^{-1}(V)$.
We set
$U_{x,q}^{\ast}=U_{x,q}\setminus\{\infty\}$ and
$V^{\ast}:=V_{|U_{x,q}^{\ast}}$.
We obtain a $\lambda$-connection
$\DD^{\lambda}_{V^{\ast}}$
corresponding to $\delbar_{\Vtilde|\nbigb_q^{\lambda\ast}}$.

Fix $t_1^0\in S^1_T$.
Note that
$\nbigb^{\lambda}_q\langle t^0_1\rangle$
is naturally equipped with a complex structure,
and $(\Vtilde,\delbar_{\Vtilde})_{|\nbigb^{\lambda}_q\langle t_1^0\rangle}$
is naturally a holomorphic vector bundle on
$\nbigb^{\lambda}_q\langle t_1^0\rangle$.
(See \S\ref{subsection;20.8.8.30}
for $\nbigb^{\lambda}_q\langle t_1^0\rangle$.)
The induced morphism
$\Psi^{\lambda}_{q,t_1^0}:\nbigb^{\lambda}_q\langle t_1^0\rangle
\lrarr U_{x,q}$
is holomorphic, and the image is an open neighbourhood of $\infty$.
We have the isomorphism of $C^{\infty}$-vector bundles
$\Vtilde_{|\nbigb^{\lambda}_1\langle t_1^0\rangle}
\simeq (\Psi_q^{\lambda})^{-1}(V)_{|\nbigb^{\lambda}_q\langle t_1^0\rangle}$.
Note that
$(\Psi_q^{\lambda})^{-1}(V)_{|\nbigb^{\lambda\ast}_q\langle t_1^0\rangle}$
is equipped with the holomorphic structure induced by
$\delbar_{V^{\ast}}$,
and that
$\Vtilde_{|\nbigb^{\lambda\ast}_q\langle t_1^0\rangle}
\simeq
(\Psi_q^{\lambda})^{-1}(V)_{|\nbigb^{\lambda\ast}_q\langle t_1^0\rangle}$
is holomorphic.
(See \S\ref{subsection;17.10.25.100}
for $\nbigb^{\lambda\ast}_1\langle t_1^0\rangle$.)
It implies that $\delbar_{V^{\ast}}$
uniquely extends to a holomorphic structure $\delbar_V$ of $V$.

Let $\vecv$ be a holomorphic frame of $(V,\delbar_V)$ on $U_{x,q}$.
We set $\vecvtilde=(\Psi_q^{\lambda})^{-1}(\vecv)$.
There exists a matrix-valued $C^{\infty}$-function $A$ on $U_{x,q}$
determined by 
\[
 \del_{\Vtilde,t_1}\vecvtilde
 =\vecvtilde\cdot \frac{-2\sqrt{-1}}{1+|\lambda|^2}
 (\Psi_q^{\lambda})^{\ast}(A).
\]
Let $\DDlambda_{V^{\ast},w}$ denote the differential operator
of $V^{\ast}$
induced by $\DD^{\lambda}_{V^{\ast}}$ and $\del_w$.
By the construction of $\DD^{\lambda}_{V^{\ast}}$
(see the proof of Lemma \ref{lem;21.8.12.2}),
we have
$\DD^{\lambda}_{V^{\ast},w}\vecv_{|U_{x,q}^{\ast}}
 =\vecv_{|U_{x,q}^{\ast}}\cdot A$.
It implies that $\DD^{\lambda}_{V^{\ast}}$
uniquely extends to a differential operator $\DD^{\lambda}$
of $V$.

Thus, we obtain $(V,\delbar_V,\DD^{\lambda\,(1,0)})$
from $(\Vtilde,\delbar_{\Vtilde})$.
The two constructions are mutually inverse.
\hfill\qed

\begin{cor}
\label{cor;21.8.13.24}
The above construction induces an equivalence between
the following objects.
\begin{itemize}
 \item Locally free $\nbigo_{U_{x,q}}$-module $\nbigl$
       equipped with a meromorphic $\lambda$-connection
       $\DDlambda:\nbigl\lrarr\nbigl\otimes\Omega^{1}_{U_{x,q}}((q+1)\infty)$.
 \item $S^1_T$-equivariant
       locally free $\nbigo_{\nbigb^{\lambda}_q}$-modules
       $\nbigltilde$.
\hfill\qed
\end{itemize}
\end{cor}

\subsubsection{Comparison with
the construction in \S\ref{subsection;21.8.12.40}}

Let $\nbigl$ be a holomorphic vector bundle
on $U_{x,q}$ equipped with a meromorphic $\lambda$-connection
$\nabla^{\lambda}:
\nbigl\lrarr\nbigl\otimes\Omega^{1}_{U_{x,q}}((q+1)\infty)$.
Let $(V,\delbar_V)$ be the holomorphic vector bundle
corresponding to $\nbigl$.
It is equipped with
$\DD^{\lambda\,(1,0)}:C^{\infty}(U_{x,q},V)
\lrarr
C^{\infty}\bigl(U_{x,q},
V\otimes\Omega^{1,0}_{U_{x,q}}((q+1)\infty)\bigr)$
induced by $\nabla^{\lambda}$.
We obtain $(\Vtilde,\delbar_{\Vtilde})$ on $\nbigb_q^{\lambda}$,
which is mini-holomorphic by (\ref{eq;21.8.12.1}).

\begin{lem}
$(\Psi_q^{\lambda})^{\ast}(\nbigl,\nabla^{\lambda})$
is the sheaf of mini-holomorphic sections of
$(\Vtilde,\delbar_{\Vtilde})$.
\end{lem}
\pf
By the construction,
it is easy to see that
$\nbigltilde^{\infty}$
in (\ref{eq;21.8.20.1})
is the sheaf of $C^{\infty}$-sections $u$ of $\Vtilde$
such that $\del_{\Vtilde,\betabar_{1}}u=0$.
Then, the claim is clear.
\hfill\qed

\subsubsection{Comparison with the construction in the formal case}
\label{subsection;21.8.13.2}

Let $(V,\delbar_V)$ be a holomorphic vector bundle
with $\DD^{\lambda\,(1,0)}$ as in \S\ref{subsection;21.8.12.41}.
Suppose that
$[\delbar_V,\DD^{\lambda\,(1,0)}]_{|\widehat{\infty}_{x,q}}=0$,
i.e., the Taylor series of
$[\delbar_V,\DD^{\lambda\,(1,0)}]$ at $x_q^{-1}=0$ is $0$.

Let $\nbigl$ denote the $\nbigo_{U_{x,q}}$-module
obtained as the sheaf of holomorphic sections of $(V,\delbar_V)$,
and we obtain
the $\nbigo_{\inftyhat_{x,q}}$-module
$\nbiglhat:=\nbigl_{|\widehat{\infty}_{x,q}}$.
It is equipped with the meromorphic $\lambda$-connection
$\nabla^{\lambda}$ induced by $\DD^{\lambda\,(1,0)}$.
In \S\ref{subsection;17.10.7.2},
we obtain
the $\nbigk_{\Hhat^{\lambda}_{\infty,q}}$-module
\[
\widetilde{\nbiglhat}^{\infty}
=\nbigk_{\Hhat^{\lambda}_{\infty,q}}\otimes
_{(\Psi_q^{\lambda})^{-1}\nbigo_{\inftyhat_{x,q}}}
(\Psi_q^{\lambda})^{-1}(\nbiglhat)
\]
equipped with the action of $\del_{t_1}$,
and 
the $\nbigo_{\Hhat^{\lambda}_{\infty,q}}$-module
$(\Psi_q^{\lambda})^{\ast}(\nbiglhat,\nabla^{\lambda})$
as the kernel of $\del_{t_1}$.
Let us explain how to obtain
$(\widetilde{\nbiglhat}^{\infty},\del_{t_1})$
and
$\widetilde{\nbiglhat}$
naturally from $(\Vtilde,\delbar_{\Vtilde})$
induced by $(V,\delbar_V+\DD^{\lambda\,(1,0)})$.

\vspace{.1in}
Let us introduce sheaves of algebras 
$\nbigc^{\infty}_{\Hhat^{\lambda\cov}_{\infty,q}}$
on $H^{\lambda\cov}_{\infty,q}$,
and 
$\nbigc^{\infty}_{\Hhat^{\lambda}_{\infty,q}}$
on $H^{\lambda}_{\infty,q}$.
For any open subset $U\subset H^{\lambda\cov}_{\infty,q}$,
we set
$\nbigc^{\infty}_{\Hhat^{\lambda\cov}_{\infty,q}}(U)
=\bigl\{
\sum_{i,j\geq 0} a_{i,j}(t_1)\beta_{1,q}^{-i}\betabar_{1,q}^{-j}
\bigr\}$.
There exist naturally defined linear differential operators
$\del_{t_1}$ and $\del_{\betabar_1}$ on
$\nbigc^{\infty}_{\Hhat^{\lambda\cov}_{\infty,q}}(U)$.
Thus, we obtain the sheaf of algebras
$\nbigc^{\infty}_{\Hhat^{\lambda\cov}_{\infty,q}}$
equipped with the differential operators
$\del_{t_1}$ and $\del_{\betabar_1}$
on $H^{\lambda\cov}_{\infty,q}$.
We note that the kernel of
$\del_{\betabar_1}:
\nbigc^{\infty}_{\Hhat^{\lambda\cov}_{\infty,q}}
\lrarr
\nbigc^{\infty}_{\Hhat^{\lambda\cov}_{\infty,q}}$
is $\nbigk_{\Hhat^{\lambda\cov}_{\infty,q}}$.
The sheaf 
$\nbigc^{\infty}_{\Hhat^{\lambda\cov}_{\infty,q}}$
is naturally equivariant,
and we obtain
$\nbigc^{\infty}_{\Hhat^{\lambda}_{\infty,q}}$
on $H^{\lambda}_{\infty,q}$ as the descent,
which is equipped with
$\del_{t_1}$ and $\del_{\betabar_1}$.
The kernel of
$\del_{\betabar_1}:
\nbigc^{\infty}_{\Hhat^{\lambda}_{\infty,q}}
\lrarr
\nbigc^{\infty}_{\Hhat^{\lambda}_{\infty,q}}$
is $\nbigk_{\Hhat^{\lambda}_{\infty,q}}$.

Let $\nbigc^{\infty}_{\nbigb_q^{\lambda}}$
be the sheaf of $C^{\infty}$-functions on $\nbigb_q^{\lambda}$.
By taking the Taylor series along $H^{\lambda}_{\infty,q}$,
we obtain the morphism of sheaves of algebras
$\nbigc^{\infty}_{\nbigb_q^{\lambda}|H^{\lambda}_{\infty,q}}
\lrarr \nbigc^{\infty}_{\Hhat^{\lambda}_{\infty,q}}$.
Thus, we obtain a morphism of ringed spaces
$k^{\lambda}_{q,C^{\infty}}:
(H^{\lambda}_{\infty,q},\nbigc^{\infty}_{\Hhat^{\lambda}_{\infty,q}})
\lrarr
 (\nbigb^{\lambda}_q,\nbigc^{\infty}_{\nbigb^{\lambda}_{q}})$.
For a locally free
$\nbigc^{\infty}_{\nbigb_q^{\lambda}}$-module
$\nbigw$,
we obtain the locally free
$\nbigc^{\infty}_{\Hhat^{\lambda}_{\infty,q}}$-module
$(k^{\lambda}_{q,C^{\infty}})^{\ast}(\nbigw)$.
If $\nbigw$ is equipped with
linear differential operators
$\del_{\nbigw,\kappa}$ $(\kappa=t_1,\betabar_1)$
satisfying 
$\del_{\nbigw,\kappa}(fu)
=\del_{\kappa}(f)\cdot s+f\del_{\nbigw,\kappa}(s)$
for any local sections
$f$ and $s$ of
$\nbigc^{\infty}_{\nbigb^{\lambda}_q}$
and $\nbigw$, respectively,
then
$(k^{\lambda}_{q,C^{\infty}})^{\ast}\nbigw$ is equipped with 
the induced linear differential operators
$\del_{(k^{\lambda}_{q,C^{\infty}})^{\ast}\nbigw,\kappa}$
$(\kappa=t_1,\betabar_1)$
satisfying 
$\del_{(k^{\lambda}_{q,C^{\infty}})^{\ast}\nbigw,\kappa}(fu)
=\del_{\kappa}(f)\cdot s
+f\del_{(k^{\lambda}_{q,C^{\infty}})^{\ast}\nbigw,\kappa}(s)$
for any local sections
$f$ and $s$ of
$\nbigc^{\infty}_{\Hhat^{\lambda}_{\infty,q}}$
and $(k^{\lambda}_{q,C^{\infty}})^{\ast}\nbigw$, respectively.

\vspace{.1in}

To simplify the description,
the sheaf of $C^{\infty}$-sections of $V$
is also denoted by $V$.
We obtain
the $\nbigc^{\infty}_{\Hhat^{\lambda}_{\infty,q}}$-module
$(k^{\lambda}_{q,C^{\infty}})^{\ast}\Vtilde_{C^{\infty}}$,
which is equipped with the induced linear differential operators
$\del_{(k^{\lambda}_{q,C^{\infty}})^{\ast}\Vtilde,t_1}$
and $\del_{(k^{\lambda}_{q,C^{\infty}})^{\ast}\Vtilde,\betabar_1}$
By (\ref{eq;21.8.12.1}),
the two operators
$\del_{(k^{\lambda}_{q,C^{\infty}})^{\ast}\Vtilde,t_1}$
and $\del_{(k^{\lambda}_{q,C^{\infty}})^{\ast}\Vtilde,\betabar_1}$
are commutative.
The following lemma is obvious
by the construction of $\widetilde{\nbiglhat}$.
\begin{lem}
\label{lem;21.8.12.60} 
The kernel of $\del_{(k^{\lambda}_{q,C^{\infty}})^{\ast}\Vtilde,\betabar_1}$
on $(k^{\lambda}_{q,C^{\infty}})^{\ast}\Vtilde_{C^{\infty}}$
is naturally isomorphic to
$\widetilde{\nbiglhat}^{\infty}$
in a way compatible with the actions of $\del_{t_1}$.
Hence, the intersection of the kernels of
$\del_{(k^{\lambda}_{q,C^{\infty}})^{\ast}\Vtilde,\betabar_1}$
and
$\del_{(k^{\lambda}_{q,C^{\infty}})^{\ast}\Vtilde,t_1}$
is naturally isomorphic to
$(\Psi_q^{\lambda})^{\ast}(\nbiglhat,\nabla^{\lambda})$.
\hfill\qed
\end{lem}

\chapter{Basic examples of monopoles around infinity}
\label{section;17.10.2.20}

We explain some basic example of monopoles,
which are models in the study of the asymptotic behaviour
of monopoles.
In particular,
we shall later use the examples in \S\ref{subsection;17.9.26.1}
and \S\ref{subsection;17.10.25.122},
which are summarized in \S\ref{subsection;21.8.14.2}.

\section{Examples of monopoles with pure slope $\ell/p$}
\label{subsection;17.9.26.1}

Let $p$ be a positive integer.
Let $Y_p^{0\cov\ast}$ be as in \S\ref{subsection;21.8.19.1}.
The pull back of the function $w$
by $Y_p^{0\cov\ast}\lrarr \nbigm^0$
is also denoted by $w$,
which is equal to $w_p^p$.
Hence, for instance,
$dw$ also means $pw_p^{p-1}dw_p$.
We also use the complex vector fields
$\del_w=p^{-1}w_p^{1-p}\del_{w_p}$
and $\del_{\wbar}=p^{-1}\wbar_p^{1-p}\del_{\wbar_p}$.

Let $\LL^{\cov\ast}_{p}(\ell)$ denote the product line bundle
on $Y^{0\cov\ast}_p$
with a global frame $e$.
Let $h$ be the Hermitian metric
determined by $h(e,e)=1$.
We take a positive number $T$ and
an integer $\ell$,
and let $\nabla$ and $\phi$
be the unitary connection and the Higgs field
given as follows:
\[
 \nabla e=e\Bigl(
 \frac{\ell}{pT}\frac{t}{2}
 \Bigl(
 \frac{d\wbar}{\wbar}
-\frac{dw}{w}
 \Bigr)
 \Bigr),
\quad\quad
 \phi=\sqrt{-1}\frac{\ell}{pT}\log|w|.
\]

\begin{lem}
$(\LL^{\cov\ast}_p(\ell),h,\nabla,\phi)$ satisfies
the Bogomolny equation
with respect to the metric
$dt\,dt+dw\,d\wbar$
on $Y_p^{0\cov\ast}$.
\index{monopole $\LL^{\cov\ast}_p(\ell)$}
\end{lem}
\pf
We obtain the following equalities:
\[
 \nabla\circ\nabla
=\frac{\ell}{pT}
 \frac{dt}{2}
 \Bigl(
 \frac{d\wbar}{\wbar}
-\frac{dw}{w}
 \Bigr),
\quad\quad
 \nabla\phi
=\sqrt{-1}\frac{\ell}{pT}
 \frac{1}{2}
 \Bigl(
 \frac{dw}{w}
+\frac{d\wbar}{\wbar}
 \Bigr).
\]
For the expression
$F=F_{t\wbar}dt\,d\wbar+F_{tw}\,dt\,dw+F_{w\wbar}dw\,d\wbar$,
we obtain
\[
 F_{t\wbar}
=\frac{\ell}{pT}\frac{1}{2}
 \wbar^{-1},
\quad
 F_{tw}
=-\frac{\ell}{pT}\frac{1}{2}w^{-1},
\quad
 F_{w\wbar}=0.
\]
We also obtain the following equalities:
\[
 \nabla_t\phi=0,
\quad
 \nabla_w\phi
=\sqrt{-1}\frac{\ell}{pT}\frac{1}{2}w^{-1},
\quad
 \nabla_{\wbar}\phi
=\sqrt{-1}\frac{\ell}{pT}\frac{1}{2}\wbar^{-1}.
\]
Recall the following relation
for the Hodge star operator:
\[
 \ast (dt\,dw)
=-\sqrt{-1}\,dw,
\quad
 \ast(dt\,d\wbar)
=\sqrt{-1}d\wbar,
\quad
 \ast(dw\,d\wbar)
=-2\sqrt{-1}\,dt.
\]
We have $F_{w\wbar}=\nabla_t\phi=0$.
We have
\[
 \ast
 \bigl(
 F_{tw}\,dt\,dw
 \bigr)
=-\sqrt{-1}F_{tw}\,dw
=\sqrt{-1}\frac{\ell}{2pT}w^{-1}\,dw
=
 \nabla_w\phi\,dw.
\]
We also have
\[
 \ast\bigl(
 F_{t\wbar}\,dt\,d\wbar
 \bigr)
=\sqrt{-1}F_{t\wbar}\,d\wbar
=\sqrt{-1}\frac{\ell}{2pT}
 \wbar^{-1}\,d\wbar
=\nabla_{\wbar}\phi\,d\wbar.
\]
Hence, the Bogomolny equation 
is satisfied.
\hfill\qed

\subsection{Equivariance with respect to the $\seisuu$-action}
\label{subsection;21.8.16.10}

Recall that
$Y_p^{0\cov\ast}$ is equipped with the naturally induced
$\seisuu$-action $\kappa_p$
(see \S\ref{subsection;21.8.19.1}).
In particular,
$\kappa_{p,1}:
Y_p^{0\cov\ast}\lrarr Y_p^{0\cov\ast}$
is given by
$\kappa_{p,1}(t,w_p)=(t+T,w_p)$.
We obtain
\[
 \kappa_{p,1}^{-1}(\nabla)
 \kappa_{p,1}^{-1}(e)
=\kappa_{p,1}^{-1}(e)
 \left(
 \frac{\ell}{pT}\frac{t}{2}
 \Bigl(
 \frac{d\wbar}{\wbar}
 -\frac{dw}{w}\Bigr)
+\frac{\ell}{2p}
 \Bigl(
 \frac{d\wbar}{\wbar}
-\frac{dw}{w}
 \Bigr)
 \right),
\]
\[
 \kappa_{p,1}^{\ast}\phi
=\sqrt{-1}\frac{\ell}{pT}
 \log|w|.
\]
We also have
\[
 \nabla\Bigl(
 e\cdot
 \bigl(|w_p|w_p^{-1}\bigr)^{\ell}
 \Bigr)
=e\bigl(
 |w_p|w_p^{-1}
 \bigr)^{\ell}
 \cdot
 \left(
 \frac{\ell}{pT}\frac{t}{2}
 \Bigl(
 \frac{d\wbar}{\wbar}
 -\frac{dw}{w}\Bigr)
+\frac{\ell}{2p}
 \Bigl(
 \frac{d\wbar}{\wbar}
-\frac{dw}{w}
 \Bigr)
 \right).
\]
By the correspondence
$\kappa_{p,1}^{-1}(e)\longmapsto
 e\cdot (|w_p|w_p^{-1})^{\ell}$,
the monopole 
$(\LL^{\cov\ast}_p(\ell),h,\nabla,\phi)$
is equivariant
with respect to the $\seisuu$-action $\kappa_p$.
Recall that the quotient of $Y_p^{0\ast\cov}$
by the $\seisuu$-action is denoted by
$Y^{0\ast}_p$.
(See \S\ref{subsection;21.8.19.1}.)
As the descent, we obtain
a monopole 
$(\LL^{\ast}_p(\ell),h,\nabla,\phi)$
on 
$Y_p^{0\ast}$.
\index{monopole $\LL^{\ast}_p(\ell)$}

\subsection{The underlying mini-holomorphic bundle
at $\lambda=0$}
\label{subsection;17.10.25.121}

Let $\LL^{0\cov\ast}_p(\ell)$ denote
the $\seisuu$-equivariant mini-holomorphic bundle
on $Y_p^{0\cov\ast}$
underlying the monopole
$(\LL^{\cov\ast}_p(\ell),h,\nabla,\phi)$.
\index{mini-holomorphic bundle $\LL^{0\cov\ast}_p(\ell)$}

\begin{prop}
\label{prop;20.8.8.40}
 There exists a mini-holomorphic frame
$u^{0}_{p,\ell}$ of $\LL^{0\cov\ast}_p(\ell)$
on $Y^{0\cov\ast}_{p}$
such that the following holds:
\begin{itemize}
\item
$|u^0_{p,\ell}|_h
=|w|^{-\ell t/(pT)}=|w_p|^{-\ell t/T}$.
\item
$\kappa_{p,1}^{\ast}(u^0_{p,\ell})
=w_p^{-\ell}u^0_{p,\ell}$
under the identification
$\kappa_{p,1}^{\ast}(e)
=e(|w_p|w_p^{-1})^{\ell}$.
\end{itemize}
\end{prop}
\pf
We have the following equalities:
\[
 \bigl(
 \nabla_t-\sqrt{-1}\phi
 \bigr)e
=e\,\frac{\ell}{pT}\log|w|,
\quad\quad
 \nabla_{\wbar}e
=e\,\frac{\ell}{pT}\frac{t}{2}\frac{1}{\wbar}.
\]
We define the section $u^0_{p,\ell}$ as follows:
\[
 u^0_{p,\ell}:=e\cdot
 \exp\Bigl(
 -\frac{\ell}{pT}
 t\log|w|
 \Bigr).
\]
Then, $u^0_{p,\ell}$ has the desired properties.
\hfill\qed

\vspace{.1in}
Let $\LL^{0\ast}_{p}(\ell)$
denote the mini-holomorphic bundle
on $Y^{0}_p$
underlying the monopole
$(\LL^{\ast}_p(\ell),h,\nabla,\phi)$.
Let $\LL^{0}_p(\ell)$
and $\nbigp_{\ast}\LL^{0}_p(\ell)$
denote 
the $\nbigo_{Y^{0}_p}(\ast H^{0}_{\infty,p})$-module
and the family of filtered sheaves
obtained from
$\LL^{0\ast}_p(\ell)$ with $h$
by the procedure in \S\ref{subsection;17.10.25.100}.

\begin{cor}
\label{cor;21.8.21.25}
 $\LL^{0}_p(\ell)$
is a locally free 
$\nbigo_{Y^{0}_p}(\ast H^{0}_{\infty,p})$-module,
and 
$\nbigp_{\ast}\LL^{0}_{p}(\ell)$
is a good filtered bundle over  $\LL^{0}_p(\ell)$.
\index{$\nbigo_{Y^{0}_p}(\ast H^{0}_{\infty,p})$-module $\LL^0_p(\ell)$}
\index{good filtered bundle $\nbigp_{\ast}\LL^{0}_{p}(\ell)$}
 (See Definition {\rm\ref{df;20.7.24.1}} and
Definition {\rm\ref{df;21.8.19.2}}.) 
The norm estimate holds for
$\nbigp_{\ast}\LL^{0}_{p}(\ell)$
with the metric $h$. 
Moreover,
$\nbigp_{\ast}\LL^{0}_{p}(\ell)_{|\Hhat^{0}_{\infty,p}}$
is isomorphic to
$\nbigp^{(0)}_{\ast}\LLhat^{0}_p(\ell,1)$
in {\rm\S\ref{subsection;17.10.7.1}}.
(See {\rm\S\ref{subsection;20.7.22.100}}
for $\Hhat^0_{\infty,q}$.)
\hfill\qed
\end{cor}

\subsection{The underlying mini-holomorphic bundle at
   general $\lambda$}
\label{subsection;17.10.25.120}

Let $\lambda$ be any complex number.
Let $\LL^{\lambda\cov\ast}_p(\ell)$
be the $\seisuu$-equivariant mini-holomorphic bundle
on $Y^{\lambda\cov\ast}_p$
underlying $(\LL^{\cov\ast}_p(\ell),h,\nabla,\phi)$.
\index{mini-holomorphic bundle $\LL^{\lambda\cov\ast}_p(\ell)$}
(See \S\ref{subsection;20.7.22.100} for $Y^{\lambda\cov\ast}_p$,
and recall that
$Y^{\lambda\cov\ast}_p=Y^{0\ast\cov}_p$ as Riemannian manifolds.)
We shall prove the following propositions
in \S\ref{subsection;17.10.3.20}--\ref{subsection;17.9.17.100}.

\begin{prop}
\label{prop;17.10.3.30}
There exists a neighbourhood $\nbigu_p^{\lambda}$
of $H^{\lambda\cov}_{\infty,p}$ in $Y^{\lambda\cov}_p$
and a mini-holomorphic frame
$u^{\lambda}_{p,\ell}$ of $\LL^{\lambda\cov\ast}_p(\ell)$
on a neighbourhood 
$\nbigu_p^{\lambda}\setminus H^{\lambda\cov}_{\infty,p}$ 
such that the following holds:
\begin{description}
\item[(A1)]
For any $P\in H^{\lambda\cov}_{p}$,
there exist $C(P)>1$
and a neighbourhood $U_P$ of $P$ in $\nbigu_p^{\lambda}$
such that 
the following holds on $U_P\setminus H^{\lambda\cov}_p$:
\[
 C(P)^{-1}|w|^{-\ell t_1/(pT)}
\leq
 |u_{p,\ell}^{\lambda}|_h
\leq
 C(P)|w|^{-\ell t_1/(pT)}.
\]
\item[(A2)]
On 
$\bigl(
 \nbigu^{\lambda}_p\cap\kappa_{p,1}^{-1}(\nbigu^{\lambda}_p)
 \bigr)\setminus H^{\lambda\cov}_p$,
we have the following equality:
\[
  \kappa_{p,1}^{\ast}(u_{p,\ell}^{\lambda})
=u_{p,\ell}^{\lambda}
 \cdot
 (1+|\lambda|^2)^{\ell/p}
 (\beta_1+2\sqrt{-1}\lambda T)^{-\ell/p}
 \times
 \exp\Bigl(
 \frac{\ell}{p}G(2\sqrt{-1}\lambda T/\beta_1)
 \Bigr),
\]
where $G(x)=1-x^{-1}\log(1+x)=x/2+O(x^2)$.
\end{description}
\end{prop}

Let $\LL^{\lambda\ast}_{p}(\ell)$
denote the mini-holomorphic bundle
on $Y^{\lambda\ast}_p$
underlying the monopole
$(\LL^{\ast}_p(\ell),h,\nabla,\phi)$.
(See \S\ref{subsection;20.7.22.100} for $Y^{\lambda\ast}_p$,
and recall $Y^{\lambda\ast}_p=Y^{0\ast}_p$
as Riemannian manifolds.)
Let $\LL^{\lambda}_p(\ell)$
and $\nbigp_{\ast}\LL^{\lambda}_p(\ell)$
denote 
the $\nbigo_{Y^{\lambda}_p}(\ast H^{\lambda}_{\infty,p})$-module
and the family of sheaves obtained from
$\LL^{\lambda\ast}_p(\ell)$ with $h$
by the procedure in \S\ref{subsection;17.10.25.100}.

\begin{cor}
\label{cor;17.10.25.131}
 $\LL^{\lambda}_p(\ell)$
is a locally free 
$\nbigo_{Y^{\lambda}_p}(\ast H^{\lambda}_{\infty,p})$-module,
and 
$\nbigp_{\ast}\LL^{\lambda}_{p}(\ell)$
is a good filtered bundle over  $\LL^{\lambda}_p(\ell)$.
\index{$\nbigo_{Y^{\lambda}_p}(\ast H^{\lambda}_{\infty,p})$-module
$\LL^{\lambda}_p(\ell)$}
\index{good filtered flat bundle $\nbigp_{\ast}\LL^{\lambda}_{p}(\ell)$}
The norm estimate holds for
 $\nbigp_{\ast}\LL^{\lambda}_{p}(\ell)$
with the metric $h$.
Moreover,
$\nbigp_{\ast}\LL^{\lambda}_{p}(\ell)_{|\Hhat^{\lambda}_{\infty,p}}$
is isomorphic to
$\nbigp^{(0)}_{\ast}\LLhat^{\lambda}_p(\ell,1)$
in {\rm\S\ref{subsection;17.10.7.1}}.
(See Remark {\rm\ref{rem;21.8.19.3}} for the comparison of
$\Hhat^{\lambda\cov}_{\infty,p}$
and
$\Hhat^{\cov}_{\infty,p}$ in {\rm\S\ref{subsection;20.7.21.1}}.)
\hfill\qed
\end{cor}

\subsection{Complex structures}
\label{subsection;17.10.3.20}

We set 
$\Ytilde^{0\cov\ast}:=\real_s\times Y^{0\cov\ast}$
as a $C^{\infty}$-manifold.
By the complex coordinate system
$(z,w)=(s+\sqrt{-1}t,w)$,
$\Ytilde^{0\cov\ast}$ is a complex manifold.
Let $\Ytilde^{\lambda\cov\ast}$ denote the complex manifold
induced by  the complex coordinate system
$(\xi,\eta)=(z+\lambda\wbar,w-\lambda\zbar)$.

Let $\Ytilde_p^{0\cov\ast}:=\real_s\times Y_p^{0\cov\ast}$.
It  is naturally a covering space of $\Ytilde^{0\cov\ast}$,
and hence equipped with the complex structure.
The complex coordinate system
$(z,w)$ of $\Ytilde^{0\cov\ast}$
induces local coordinate systems
of $\Ytilde_p^{0\cov\ast}$.
The complex vector fields
$\del_z$, $\del_{\zbar}$, $\del_w$ and $\del_{\wbar}$
are naturally lifted to complex vector fields
on $\Ytilde_p^{0\cov\ast}$,
which are denoted by the same notation.

Let $\Ytilde_p^{\lambda\cov\ast}$
be the complex manifold
induced by the local coordinate system $(\xi,\eta)$.
The complex vector fields
$\del_{\xibar}$, $\del_{\xi}$, $\del_{\eta}$, $\del_{\etabar}$
on $\Ytilde^{\lambda\cov\ast}$
are naturally lifted to
the complex vector fields on $\Ytilde_p^{\lambda\cov\ast}$
which are denoted by the same notation.

By using the local coordinate system
$(\alpha_1,\beta_1)$
given as in \S\ref{subsection;17.10.2.1},
we obtain the natural identification
$\Ytilde_p^{\lambda\cov\ast}
=\real_{\Re\alpha_1}\times Y_p^{\lambda\cov\ast}$.

\subsection{The induced instanton and the underlying holomorphic bundle}

Let $(\LL^{\cov\ast}_p(\ell),h,\nabla,\phi)$ be the monopole
on $Y^{0\cov\ast}_p$ in \S\ref{subsection;17.9.26.1}.
We obtain the induced instanton
$(\Etilde,\htilde,\nablatilde)$
on $\Ytilde^{0\cov\ast}_p$
as explained in \S\ref{subsection;17.10.3.10}.
The connection is described as follows
with respect to the actions of
the complex vector fields:
\[
 \nablatilde_{\wbar}e
=e\cdot
 \Bigl(
 \frac{\ell}{pT}
 \frac{t}{2}\wbar^{-1}
 \Bigr),
\quad
 \nablatilde_we
=e\cdot
 \Bigl(
 -\frac{\ell}{pT}\frac{t}{2}w^{-1}
 \Bigr),
\]
\[
 \nablatilde_ze
=e\cdot
 \Bigl(
 \frac{\sqrt{-1}}{2}
 \frac{\ell}{pT}
 \log|w|
 \Bigr),
\quad
 \nablatilde_{\zbar}e
=e\cdot
 \Bigl(
 \frac{\sqrt{-1}}{2}
 \frac{\ell}{pT}
 \log|w|
 \Bigr).
\]

Let us look at 
the holomorphic bundle
on $\Ytilde_p^{\lambda\cov\ast}$
underlying $(\Etilde,\htilde,\nablatilde)$.
The actions of $\nablatilde_{\xibar}$
and $\nablatilde_{\etabar}$
are given as follows:
\[
 \nablatilde_{\xibar}e=
 e\cdot
 \frac{1}{1+|\lambda|^2}
 \frac{\ell}{pT}
 \Bigl(
 \frac{\sqrt{-1}}{2}\log|w|
-\frac{\lambda t}{2}w^{-1}
 \Bigr),
\]
\[
 \nablatilde_{\etabar}e=
 e\cdot
 \frac{1}{1+|\lambda|^2}
 \frac{\ell}{pT}
 \Bigl(
 \frac{t}{2}\wbar^{-1}
-\frac{1}{2}\lambda\sqrt{-1}
 \log|w|
 \Bigr).
\]

In particular, we obtain{\small
\begin{multline}
 \nablatilde_{\etabar}e=
e\cdot\frac{1}{1+|\lambda|^2}
 \frac{\ell}{pT}  \times
 \\
 \Bigl(
 \frac{1}{4\sqrt{-1}}
 \frac{\xi-\xibar
 +\lambdabar\eta+|\lambda|^2\xi}
 {\etabar+\lambdabar\xi}
-\frac{\sqrt{-1}}{2}
 \lambda
 \log|\eta+\lambda\xibar|
-\frac{\lambda}{4\sqrt{-1}}
+\frac{\sqrt{-1}}{2}\lambda
 \log(1+|\lambda|^2)
 \Bigr). 
\end{multline}
}

\subsection{$C^{\infty}$-frame}

We put as follows:
\begin{multline}
 \vtilde^{\lambda}:=
 e\exp\Bigl(
 -\frac{1}{2\sqrt{-1}}
 \frac{\ell}{pT}
 (\xi-\xibar)
 \log\bigl|\etabar+\lambdabar\xi\bigr|
\Bigr.\\
\Bigl.
+\frac{1}{1+|\lambda|^2}
 \frac{\ell}{pT}
 \Bigl(
 -\frac{\lambdabar}{2\sqrt{-1}}
 (\eta+\lambda\xibar)
 \log|\etabar+\lambdabar\xi|
-\frac{\lambda}{2\sqrt{-1}}
 (\etabar+\lambdabar\xi)
 \log|\etabar+\lambdabar\xi|
 \Bigr.
 \\
 \Bigl.
-
\frac{\lambda\sqrt{-1}}{2}
 \bigl(1+\log(1+|\lambda|^2)\bigr)
 (\etabar+\lambdabar\xi)
-
\frac{\lambdabar\sqrt{-1}}{2}
 \bigl(1+\log(1+|\lambda|^2)\bigr)
 (\eta+\lambda\xibar)
 \Bigr)
 \Bigr).
\end{multline}

Then, the following holds.
\begin{itemize}
\item 
 $\nablatilde_{\etabar}\vtilde^{\lambda}=0$.
\item
 $\vtilde^{\lambda}$ is invariant
 with respect to the $\real$-action
     given by $s\cdot(z,w)=(z+s,w)$ in terms of
     the coordinate system $(z,w)$,
 or equivalently $s\cdot(\xi,\eta)=(\xi+s,\eta-\lambda s)$
 in terms of the coordinate system $(\xi,\eta)$.
\item
 $|\vtilde^{\lambda}|_h=
 \exp\bigl(
 -\frac{1}{2\sqrt{-1}}
 (\xi-\xibar)
 \frac{\ell}{pT}
 \log|\etabar+\lambdabar\xi|
 \bigr)$.
\end{itemize}
Under the identification
$\kappa_{p,1}^{\ast}(e)=|w_p|^{\ell}w_p^{-\ell}e$,
we obtain
\[
 \kappa_{p,1}^{\ast}(\vtilde^{\lambda})
=\vtilde^{\lambda}
 w_p^{-\ell}(1+|\lambda|^2)^{-\ell/p}.
\]
We also obtain the following equality:
\[
 \nablatilde_{\xibar}
 \vtilde^{\lambda}
=\vtilde^{\lambda}
 \frac{\ell}{pT}
 \Bigl(
 \frac{-1}{2\sqrt{-1}}
 \frac{\lambda}{\eta+\lambda\xi-\lambda(\xi-\xibar)}
 (\xi-\xibar)
-\frac{\sqrt{-1}}{2}
 \log(1+|\lambda|^2)
 \Bigr)
\]

In terms of the local coordinate system
$(\alpha_1,\beta_1)$ given in \S\ref{subsection;17.10.3.12},
we obtain the following:
\[
 \nablatilde_{\betabar_1}\vtilde^{\lambda}=0,
\]
\[
 \nablatilde_{\alphabar_1}\vtilde^{\lambda}
=\vtilde^{\lambda}
 \frac{\ell}{pT}
 \Bigl(
 \frac{-1}{2\sqrt{-1}}
 \frac{\lambda(\alpha_1-\alphabar_1)}
 {\beta_1-\lambda(\alpha_1-\alphabar_1)}
-\frac{\sqrt{-1}}{2}
 \log(1+|\lambda|^2)
 \Bigr).
\]
We set $t_1=\Image\alpha_1$
as in \S\ref{subsection;17.10.3.12}.
Because $\vtilde^{\lambda}$ is $\real_s$-invariant,
we obtain
\[
 \nablatilde^{\lambda}_{t_1}
 \vtilde^{\lambda}
=\frac{2}{\sqrt{-1}}
 \nablatilde_{\alphabar_1}\vtilde^{\lambda}
=\vtilde^{\lambda}
 \frac{\ell}{pT}
 \Bigl(
 \frac{\beta_1}{\beta_1-2\sqrt{-1}\lambda t_1}
-1
-\log(1+|\lambda|^2)
 \Bigr).
\]

\subsection{Proof of Proposition \ref{prop;17.10.3.30}}
\label{subsection;17.9.17.100}

We take a large $R>0$.
We consider the open subset
$\nbigu_p^{\lambda}(R)\subset Y^{\lambda\cov\ast}_p$
determined by the condition:
\[
 \nbigu_p^{\lambda}(R):=
 \left\{
 \bigl|
 \beta_1-2\sqrt{-1}\lambda t_1
 \bigr|<R,\,\,\,
 \frac{\bigl|2\lambda(t_1+T)\bigr|}
 {\bigl|\beta_1-2\sqrt{-1}\lambda t_1\bigr|}
 <\frac{1}{2},\,\,\,
 \frac{\bigl|2\lambda t_1\bigr|}
 {\bigl|\beta_1-2\sqrt{-1}\lambda t_1\bigr|}
 <\frac{1}{2}
 \right\}.
\]
Note that we have the well defined function
$w_p=(\beta_1-2\sqrt{-1}t_1)^{1/p}$
on $Y_p^{\lambda\cov\ast}$.
Hence, we have
the well defined function
\[
(\beta_1-2\sqrt{-1}\lambda T)^{1/p}
=
(\beta_1-2\sqrt{-1}t_1)^{1/p}
 \Bigl(
 1+\frac{2\sqrt{-1}\lambda(t_1+T)}{\beta_1-2\sqrt{-1}\lambda t_1}
 \Bigr)^{1/p}
\]
on $\nbigu_p^{\lambda}(R)$.
We also have the well defined branch of
the function
\[
-\log\Bigl(
1+\frac{2\sqrt{-1}\lambda t_1}
 {\beta_1-2\sqrt{-1}\lambda t_1}
 \Bigr),
\]
which we also denote by
\[
 \log\Bigl(
 1-\frac{2\sqrt{-1}t_1}{\beta_1}
 \Bigr).
\]

Because $\vtilde^{\lambda}$ is $\real_s$-invariant,
we obtain a $C^{\infty}$-frame
$v^{\lambda}$
of $\LL^{\lambda\cov\ast}_p(\ell)$.
On $\nbigu^{\lambda}_p(R)$,
we consider the following section:
\[
 u_{p,\ell}^{\lambda}
:=v^{\lambda}
 \cdot
 \exp\Bigl[
 \frac{\ell}{pT}
 \Bigl(
 t_1\log(1+|\lambda|^2)
+\frac{\beta_1}{2\sqrt{-1}\lambda}
 \log\Bigl(
 1-\frac{2\sqrt{-1}\lambda}{\beta_1}t_1
 \Bigr)
+t_1
 \Bigr)
 \Bigr].
\]
Then, we have
$\nablatilde^{\lambda}_{\betabar_1}
 u_{p,\ell}^{\lambda}=0$
and 
$\nablatilde^{\lambda}_{t_1}
 u_{p,\ell}^{\lambda}=0$.
As for the norm,
we have
\[
 \bigl|
 u_{p,\ell}^{\lambda}
 \bigr|_h
=\exp\Bigl[
-t_1\frac{\ell}{pT}\log|w|
+\frac{\ell}{pT}
 \Re\Bigl(
 \frac{\beta_1}{2\sqrt{-1}\lambda}
 \log(1-2\sqrt{-1}\lambda t_1/\beta_1)
+t_1
 \Bigr)
 \Bigr].
\]
Note that
\[
 \frac{\beta_1}{2\sqrt{-1}\lambda}
\log\Bigl(
 1-\frac{2\sqrt{-1}\lambda}{\beta_1}t_1
\Bigr)
+t_1
=\frac{\beta_1}{2\sqrt{-1}\lambda}
\times
 O\Bigl(
 \Bigl(
 t_1
 \frac{2\sqrt{-1}\lambda}{\beta_1}
 \Bigr)^2
 \Bigr).
\]
Hence, we obtain {\bf (A1)}.

Let us check {\bf (A2)}.
We have the following:
\begin{multline}
 \Phi^{\ast}(u_{p,\ell}^{\lambda})
=u_{p,\ell}^{\lambda}w_p^{-\ell}(1+|\lambda^{2}|)^{-\ell/p}
 \times \\
 \exp\Bigl(
 \frac{\ell}{pT}
 \Bigl(
 (t_1+T)\log(1+|\lambda|^2)
+\frac{\beta_1+2\sqrt{-1}\lambda T}{2\sqrt{-1}\lambda}
 \log\Bigl(
 \frac{\beta_1-2\sqrt{-1}\lambda t_1}
 {\beta_1+2\sqrt{-1}\lambda T}
 \Bigr)
 +t_1+T
 \Bigr)
 \Bigr)\times \\
 \exp\Bigl(
 -\frac{\ell}{pT}
 \Bigl(
 t_1\log(1+|\lambda|^2)
+\frac{\beta_1}{2\sqrt{-1}\lambda}
 \log\Bigl(
 \frac{\beta_1-2\sqrt{-1}\lambda t_1}
 {\beta_1}
 \Bigr)
 +t_1
 \Bigr)\Bigr).
\end{multline}

We have
\[
 \frac{\ell}{pT}
 (t_1+T)\log(1+|\lambda|^2)
-\frac{\ell}{pT}t_1\log(1+|\lambda|^2)
=\frac{\ell}{p}\log(1+|\lambda|^2).
\]
We also have
\begin{multline}
 \frac{\ell}{pT}
 \Bigl(
 \frac{\beta_1}{2\sqrt{-1}\lambda}
 \log\Bigl(
 \frac{\beta_1-2\sqrt{-1}\lambda t_1}
 {\beta_1+2\sqrt{-1}\lambda T}
 \Bigr)
+T\log\Bigl(
 \frac{\beta_1-2\sqrt{-1}\lambda t_1}
 {\beta+2\sqrt{-1}\lambda T}
 \Bigr)
  \\
+t_1+T
-\frac{\beta_1}{2\sqrt{-1}\lambda}
 \log\Bigl(
 \frac{\beta_1-2\sqrt{-1}\lambda t_1}{\beta_1}
 \Bigr)
-t_1
 \Bigr)
\\
=\frac{\ell}{pT}
 \Bigl(
 \frac{\beta_1}{2\sqrt{-1}\lambda}
 \log\Bigl(
 \frac{\beta_1}{\beta_1+2\sqrt{-1}\lambda T}
 \Bigr)
+T
+T\log\Bigl(
 \frac{\beta_1-2\sqrt{-1}\lambda t_1}
 {\beta_1+2\sqrt{-1}\lambda T}
 \Bigr)
\Bigr)
\\
=\frac{\ell}{pT}
 \Bigl(
 -\frac{\beta_1}{2\sqrt{-1}\lambda}
 \log\Bigl(
 1+\frac{2\sqrt{-1}\lambda}{\beta_1}T
 \Bigr)
+T
 \Bigr)
+\frac{\ell}{p}\cdot \log
 \Bigl(
 \frac{\beta_1-2\sqrt{-1}
 \lambda t_1}
{\beta_1+2\sqrt{-1}\lambda T}
 \Bigr).
\end{multline}
Hence,
we have
\begin{multline}
 \Phi^{\ast}(u_{p,\ell}^{\lambda})
=u_{p,\ell}^{\lambda}
 (1+|\lambda|^2)^{\ell/p}
 (\beta_1+2\sqrt{-1}\lambda T)^{-\ell/p}
 \times \\
 \exp\Bigl(
 \frac{\ell}{pT}
  \Bigl(
 \frac{-\beta_1}{2\sqrt{-1}\lambda}
 \log\Bigl(
 1+\frac{2\sqrt{-1}\lambda}{\beta_1}T
\Bigr)
+T
 \Bigr)
 \Bigr).
\end{multline}
Then, the claim of the lemma follows.
\hfill\qed

\section{Examples of monopoles induced by
wild harmonic bundles}
\label{subsection;17.10.25.130}

Let $w_p$ be a $p$-th root of the variable $w$
as in \S\ref{subsection;21.8.19.1}.
\index{variable $w_p$}
We take $\gminia(w_p)\in \cnum[w_p]$
such that 
$\deg_{w_p}\gminia\leq p$
and $\gminia(0)=0$.
We obtain the harmonic bundle on $\cnum_{w_p}^{\ast}$
given as 
$\nbigl^{\ast}(\gminia)=\nbigo_{\cnum_{w_p}^{\ast}}\,e_0$
with the Higgs field $\theta_0=d\gminia$
and the metric $h_0(e_0,e_0)=1$.
\index{harmonic bundle $(\nbigl^{\ast}(\gminia),\theta_0,h_0)$}
As explained in \S\ref{subsection;17.9.29.1},
there exists the monopole 
$(\vecnbigl^{\ast}(\gminia),h,\nabla,\phi)$
on $Y_p^{0\ast}$ 
induced by
$(\nbigl^{\ast}(\gminia),\theta_0,h_0)$.
\index{monopole $\vecnbigl^{\ast}(\gminia)$}
Let us describe the monopole explicitly.

Let $\Psi_p:Y_p^{0\ast}\lrarr \cnum_{w_p}^{\ast}$
denote the projection as in \S\ref{subsection;21.8.19.1},
i.e., the map induced by $(t,w_p)\longmapsto w_p$.
We have $\vecnbigl^{\ast}(\gminia)=\Psi_p^{-1}\nbigl^{\ast}(\gminia)$
as $C^{\infty}$-bundle
with the induced metric $h=\Psi_p^{-1}(h_0)$.
We set $e:=\Psi_p^{-1}(e_0)$.
The unitary connection and the Higgs field are given as follows:
\[
 \nabla e
=e\Bigl(
 -\sqrt{-1}\bigl(\del_w\gminia+\overline{\del_w\gminia}\bigr)\,dt
 \Bigr),
\quad\quad
 \phi e
=e\bigl(\del_{w}\gminia-\overline{\del_w\gminia}\bigr).
\]
Here, $\del_w\gminia=p^{-1}w_p^{1-p}\del_{w_p}\gminia$.
Note that $\del_w\gminia$
is of the form
$\sum_{j=0}^{p-1} a_jw_p^{-j}$,
i.e.,
it is an element of the set $S(p)$.
(See \S\ref{subsection;20.8.8.10}.)

\subsection{The underlying mini-holomorphic bundle at $\lambda=0$}
\label{subsection;20.7.23.1}

Let us describe the mini-holomorphic bundle
$\vecnbigl^{0\ast}(\gminia)$ on $Y_p^{0\ast}$
underlying
$(\vecnbigl^{\ast}(\gminia),h,\nabla,\phi)$.
\index{mini-holomorphic bundle $\vecnbigl^{0\ast}(\gminia)$}
Recall that
$\varpi_p^0:Y_p^{0\cov}\lrarr Y_p^0$
denotes the projection
(see \S\ref{subsection;21.8.19.1}).
The restriction
$Y_p^{0\cov\ast}\lrarr Y_p^{0\ast}$
is also denoted by $\varpi_p^0$.
\index{map $\varpi_p^0$}

\begin{prop}
\label{prop;17.10.25.111}
There exists a mini-holomorphic frame $u_{\gminia}^0$
of $(\varpi_p^0)^{\ast}\vecnbigl^{0\ast}(\gminia)$
on $Y_p^{0\cov\ast}$
satisfying the following:
\begin{description}
\item[(A1)]
 For any $P\in H^0_{\infty,p}$,
 there exists a neighbourhood $U_P$ of $P$ in
 $Y^{0\cov}_p$ and a constant $C(P)\geq 1$
such that
\[
 C(P)^{-1}\leq
  |u_{\gminia}^0|_h
 \leq C(P).
\]
\item[(A2)]
$\kappa_{p,1}^{\ast}u^0_{\gminia}
=u_{\gminia}^0\exp\bigl(2\sqrt{-1}T\del_w\gminia\bigr)$.
\end{description}
\end{prop}
\pf
We have
$\del_{t}e=e\bigl(-2\sqrt{-1}\del_{w}\gminia\bigr)$
and
$\del_{\wbar}e=0$.
Hence, 
we have the following mini-holomorphic frame
of $(\varpi_p^0)^{\ast}\vecnbigl^{0\ast}(\gminia)$
on $Y^{0\cov\ast}_p$:
\[
 u_{\gminia}^0:=(\varpi_p^0)^{-1}(e)\cdot
 \exp\bigl(
 2\sqrt{-1}t\del_w\gminia\bigr).
\]
It has the desired properties.
\hfill\qed

\subsection{The associated $\lambda$-connection and 
 the induced mini-holomorphic bundle at $\lambda$}

The $\lambda$-flat  bundle associated with
the harmonic bundle $(\nbigl^{\ast}(\gminia),\theta_0,h_0)$
is described as follows:
\[
 (\delbar+\lambda\theta_0^{\dagger})e_0
=e_0\,\lambda d\gminiabar,
\quad\quad
 (\lambda\del+\theta)e_0=e_0d\gminia.
\]
We set 
$v_0^{\lambda}=
e_0\exp(-\lambda\gminiabar+\lambdabar\gminia)$
on $\cnum_{w_p}^{\ast}$.
By the construction,
it is holomorphic with respect to $\delbar+\lambda\theta^{\dagger}$.
We also have
\[
 \DDlambda v_0^{\lambda}=
 (\lambda\del+\theta)v_0^{\lambda}=
 v_0^{\lambda}(1+|\lambda|^2)d\gminia.
\]
Hence, we obtain
$\DD^{\lambda}_wv_0^{\lambda}
=v_0^{\lambda}(1+|\lambda|^2)\del_w\gminia$.

\vspace{.1in}

Let us describe the mini-holomorphic bundle
$\vecnbigl^{\lambda\ast}(\gminia)$
on the mini-complex manifold $Y_p^{\lambda\ast}$
\index{mini-holomorphic bundle $\vecnbigl^{\lambda\ast}(\gminia)$}
underlying the monopole $(\vecnbigl^{\ast}(\gminia),h,\nabla,\phi)$.
Recall that
$\varpi_p^{\lambda}:Y_p^{\lambda\cov}\lrarr Y_p^{\lambda}$
denote the projection (see \S\ref{subsection;20.7.22.100}).
\index{projection $\varpi_p^{\lambda}$}
The restriction
$Y_p^{\lambda\cov\ast}\lrarr Y_p^{\lambda\ast}$
is also denoted by $\varpi_p^{\lambda}$.
Let us emphasize that
$Y_p^{\lambda\cov\ast}=Y_p^{0\cov\ast}$
and
$Y_p^{\lambda\ast}=Y_p^{0\ast}$
as Riemannian manifolds,
and that the restrictions of
$\varpi_p^{\lambda}$ and
$\varpi_p^{0}$
to $Y_p^{\lambda\cov\ast}=Y_p^{0\cov\ast}$
are the same.

\begin{prop}
\label{prop;17.10.25.112}
There exists
a mini-holomorphic frame 
$u^{\lambda}_{\gminia}$ 
of $(\varpi_p^{\lambda})^{-1}\vecnbigl^{\lambda\ast}(\gminia)$
such that the following holds:
\begin{description}
\item[(A1)]
For any $P\in H^{\lambda\cov}_{\infty,p}$,
there exists a constant $C(P)\geq 1$ and a neighbourhood $U_P$ of $P$
in $Y_p^{\lambda\cov}$
such that the following holds
on $U_P\setminus H^{\lambda\cov}_{\infty,p}$:
\[
 C(P)^{-1} \leq 
 \bigl|u^{\lambda}_{\gminia}\bigr|_h
\leq C(P).
\]
\item[(A2)]
We obtain the following equality:
\[
 \kappa_{p,1}^{\ast}u^{\lambda}_{\gminia}
=u^{\lambda}_{\gminia}
 \cdot
 \exp\Bigl[
(\lambda^{-1}+\lambdabar)
\Bigl(
 \gminia_{|w=(1+|\lambda|^2)^{-1}(\beta_1+2\sqrt{-1}\lambda T)}
-\gminia_{|w=(1+|\lambda|^2)^{-1}\beta_1}
 \Bigr)
\Bigr].
\]
\end{description}
\end{prop}
\pf
Set $v^{\lambda}:=\Psi_p^{-1}(v^{\lambda}_0)$.
By Corollary \ref{cor;21.8.12.6} with
the formulas (\ref{eq;21.8.12.30}) and (\ref{eq;21.8.12.31}),
we have the following:
\[
\del_{\betabar_1}v^{\lambda}=0,
\quad\quad
 \del_{t_1}v^{\lambda}
=v^{\lambda}\bigl(-2\sqrt{-1}\Psi_p^{\ast}(\del_w\gminia)\bigr).
\]
We remark the following:
\begin{multline}
 \del_{t_1}
 \Bigl(
 \gminia_{|w=(1+|\lambda|^2)^{-1}(\beta_1-2\sqrt{-1}\lambda t_1)}
 \Bigr)
 = \\
 \bigl(\del_w\gminia\bigr)
  _{|w=(1+|\lambda|^2)^{-1}(\beta_1-2\sqrt{-1}\lambda t_1)}
 \times \del_{t_1}(\beta_1-2\sqrt{-1}\lambda t_1)
 \\
=\frac{-2\sqrt{-1}\lambda}{1+|\lambda|^2}
 \bigl(\del_w\gminia\bigr)
 _{|w=(1+|\lambda|^2)^{-1}(\beta_1-2\sqrt{-1}\lambda t_1)}.
\end{multline}
Hence, we obtain the following mini-holomorphic frame
$u^{\lambda}_{\gminia}$ on $Y_p^{\lambda\ast}$:
\[
 u_{\gminia}^{\lambda}:=
 (\varpi_p^{\lambda})^{-1}(v^{\lambda})
 \cdot
 \exp\Bigl(
 -(\lambda^{-1}+\lambdabar)
 \bigl(
 \gminia_{|w=(1+|\lambda|^2)^{-1}(\beta_1-2\sqrt{-1}\lambda t_1)}
-\gminia_{|w=(1+|\lambda|^2)^{-1}\beta_1}
 \bigr)
 \Bigr).
\]
We can easily check that 
$u^{\lambda}_{\gminia}$ has the desired properties.
\hfill\qed

\subsection{Special case}
\label{subsection;17.10.25.122}
Let us consider the case
$\gminia=\gamma w_p^p=\gamma w$.
In this case,
the associated unitary connection
and the Higgs fields 
are given as follows:
\[
 \nabla e=e\Bigl(-\sqrt{-1}(\gamma+\gammabar)\,dt\Bigr),
\quad\quad
 \phi e=e\cdot \bigl(\gamma-\gammabar\bigr).
\]
We obtain the mini-holomorphic bundle
$\vecnbigl^{\lambda\ast}(\gamma w)$
on $Y^{\lambda\ast}_p$ 
underlying
the monopole
$(\vecnbigl^{\ast}(\gamma w),h,\nabla,\phi)$.
Proposition \ref{prop;17.10.25.111}
and Proposition \ref{prop;17.10.25.112}
are restated as follows.
(See also \S\ref{subsection;17.10.2.21}.)
\begin{prop}
\label{prop;17.10.25.110}
There exists a mini-holomorphic frame $u_{\gamma w}^{\lambda}$
 of the mini-holomorphic bundle
 $(\varpi_p^{\lambda})^{\ast}\vecnbigl^{\lambda\ast}(\gamma w)$
on $Y_p^{\lambda\cov\ast}$
satisfying the following:
\begin{description}
\item[(A1)]
We have 
$|u_{\gamma w}^{\lambda}|_h=
 \exp\bigl(-2\Image(\gamma)t_1\bigr)$.
 In particular,
 for any $P\in H^{\lambda}_{\infty,p}$,
 there exists a neighbourhood $U_P$ of $P$ in
 $Y^{\lambda\cov}_p$ and a constant $C(P)\geq 1$
such that
$C(P)^{-1}\leq
  |u_{\gamma w}^{\lambda}|_h
 \leq C(P)$.
\item[(A2)]
$\kappa_{p,1}^{\ast}u^{\lambda}_{\gamma w}
=u_{\gamma w}^{\lambda}\exp\bigl(2\sqrt{-1}T\gamma\bigr)$.
\hfill\qed
\end{description}
\end{prop}

Let $\vecnbigl^{\lambda}(\gamma w)$
and $\nbigp_{\ast}\vecnbigl^{\lambda}(\gamma w)$
denote the $\nbigo_{Y^{\lambda}_p}(\ast H^{\lambda}_{\infty,p})$-module
and the family  of sheaves
obtained from $\vecnbigl^{\lambda\ast}(\gamma w)$ with $h$
as in \S\ref{subsection;17.10.25.100}.
\begin{cor}
\label{cor;17.10.25.130}
$\vecnbigl^{\lambda}(\gamma w)$
is a locally free 
$\nbigo_{Y^{\lambda}_p}(\ast H^{\lambda}_{\infty,p})$-module,
and $\nbigp_{\ast}\vecnbigl^{\lambda}(\gamma w)$
is a good filtered bundle over $\vecnbigl^{\lambda}(\gamma w)$
in the sense of Definition {\rm\ref{df;20.7.24.1}}
and Definition {\rm\ref{df;21.8.19.2}}.
The norm estimate holds for
$\nbigp_{\ast}\vecnbigl^{\lambda}(\gamma w)$
with the metric $h$.
Moreover,
$\nbigp_{\ast}\vecnbigl^{\lambda}(\gamma w)_{|\Hhat^{\lambda}_{\infty,p}}$
is isomorphic to
$\nbigp^{(0)}_{\ast}\LLhat^{\lambda}_p(0,\alpha)$
in {\rm\S\ref{subsection;17.10.7.1}},
where $\alpha:=\exp(2\sqrt{-1}T\gamma)$.
(See Remark {\rm\ref{rem;21.8.19.3}} for the comparison of
$\Hhat^{\lambda\cov}_{\infty,p}$
and
$\Hhat^{\cov}_{\infty,p}$ in {\rm\S\ref{subsection;20.7.21.1}}.)
\hfill\qed
\end{cor}

\subsection{Notation}
\label{subsection;21.8.14.2}

For $\ell\in\seisuu$,
let $(\LL_p^{\ast}(\ell),h_{p,\ell},\nabla_{p,\ell},
\phi_{p,\ell})$
be the monopole on $Y_p^{0\ast}$
introduced in \S\ref{subsection;21.8.16.10}.
For $\alpha\in\cnum^{\ast}$,
we choose $\gamma\in\cnum$ such that
$\alpha=\exp(2\sqrt{-1}T\gamma)$,
and
we obtain the monopole
$(\vecnbigl^{\ast}(\gamma w_p^p),
h_{\gamma},\nabla_{\gamma},\phi_{\gamma})$
as above.
We obtain the following monopole:
\index{monopole $\LL^{\ast}_p(\ell,\alpha)$}
\begin{equation}
\label{eq;21.8.16.12}
 (\LL^{\ast}_p(\ell,\alpha),h_{\LL,p,\ell,\alpha},
 \nabla_{\LL,p,\ell,\alpha},\phi_{\LL,p,\ell,\alpha})
=(\LL^{\ast}_p(\ell,0),h,\nabla,\phi)\otimes
 (\vecnbigl^{\ast}(\gamma w_p^p),h_{\gamma},\nabla_{\gamma},\phi_{\gamma}).
\end{equation}
For each $\lambda$,
let $\LL_p^{\lambda\ast}(\ell,\alpha)$
denote the mini-holomorphic bundle
on $Y^{\lambda\ast}_p$
underlying the monopole (\ref{eq;21.8.16.12}).
\index{mini-holomorphic bundle $\LL_p^{\lambda\ast}(\ell,\alpha)$}
By the procedure in \S\ref{subsection;17.10.25.100},
it extends to
an $\nbigo_{Y^{\lambda}_p}(\ast H^{\lambda}_p)$-module
$\LL_p^{\lambda}(\ell,\alpha)$,
and the filtered bundle
$\nbigp_{\ast}\LL_p^{\lambda}(\ell,\alpha)$.
Let 
$\LL^{\lambda}_p(\ell)$ and 
$\nbigp_{\ast}\LL^{\lambda}_p(\ell)$ be as in 
\S\ref{subsection;17.10.25.121}--\S\ref{subsection;17.10.25.120}.
Let $\vecnbigl^{\lambda}(\gamma w_p^p)$
and $\nbigp_{\ast}\vecnbigl^{\lambda}(\gamma w_p^p)$
be as in \S\ref{subsection;17.10.25.122}.
There exist the following isomorphisms:
\index{sheaf $\LL^{\lambda}_p(\ell,\alpha)$}
\index{filtered bundle $ \nbigp_{\ast}\LL^{\lambda}_p(\ell,\alpha)$}
\[
 \LL^{\lambda}_p(\ell,\alpha)\simeq
 \LL^{\lambda}_p(\ell)
\otimes
 \vecnbigl^{\lambda}(\gamma w_p^p),
\quad\quad
 \nbigp_{\ast}\LL^{\lambda}_p(\ell,\alpha)\simeq
 \nbigp_{\ast}\LL^{\lambda}_p(\ell)
\otimes
 \nbigp_{\ast}\vecnbigl^{\lambda}(\gamma w_p^p).
\]
By
Corollary \ref{cor;21.8.21.25},
Corollary \ref{cor;17.10.25.131} and Corollary \ref{cor;17.10.25.130},
there exists an isomorphism
\begin{equation}
\label{eq;21.8.21.21}
 \nbigp_{\ast}\LL^{\lambda}_p(\ell,\alpha)_{|\Hhat^{\lambda}_{\infty,p}}
\simeq
 \nbigp^{(0)}_{\ast}\LLhat^{\lambda}_p(\ell,\alpha).
\end{equation}
(See \S\ref{subsection;17.10.7.1}
for $\nbigp_{\ast}^{(0)}\LLhat^{\lambda}_p(\ell,\alpha)$.)
In particular, $\nbigp_{\ast}\LL^{\lambda}_p(\ell,\alpha)$
is a good filtered bundle.
We also obtain the following lemma
from Corollary \ref{cor;21.8.21.25},
Corollary \ref{cor;17.10.25.131} and Corollary \ref{cor;17.10.25.130}.

\begin{lem}
\label{lem;21.8.23.3}
The norm estimate holds for
$\nbigp_{\ast}\LL^{\lambda}_q(\ell,\alpha)$
with $h_{\LL,q,\ell,\alpha}$.
\hfill\qed
\end{lem}

Recall
$\varpi_p^{\lambda}:Y_p^{\lambda\cov}\lrarr Y_p^{\lambda}$
denote the projection.
We set
\index{sheaf $\LL^{\lambda\cov}_p(\ell,\alpha)$}
\index{filtered bundle $\nbigp_{\ast}\LL^{\lambda\cov}_p(\ell,\alpha)$}
\[
 \LL^{\lambda\cov}_p(\ell,\alpha):=
 (\varpi^{\lambda}_p)^{-1} \LL^{\lambda}_p(\ell,\alpha),
 \quad
 \nbigp_{\ast}\LL^{\lambda\cov}_p(\ell,\alpha):=
 (\varpi^{\lambda}_p)^{-1}\nbigp_{\ast}\LL^{\lambda}_p(\ell,\alpha).
\]
We also set
$\LL^{\lambda\cov\ast}_p(\ell,\alpha):=
(\varpi_p^{\lambda})^{-1}(\LL^{\lambda\ast}_p(\ell,\alpha))$.
As the tensor product of
the frames $u^{\lambda}_{p,\ell}$
in Proposition \ref{prop;20.8.8.40}
and Proposition \ref{prop;17.10.3.30}
and $u^{\lambda}_{\gamma w}$ in Proposition \ref{prop;17.10.25.110},
we obtain the mini-holomorphic frame
$u^{\lambda}_{p,\ell,\alpha}$
of $\LL^{\lambda\cov}_p(\ell,\alpha)$
on a neighbourhood $\nbigu^{\lambda}_{p}$
of $H^{\lambda\cov}_{\infty,p}$.
The $\seisuu$-action is described as
\begin{equation}
\label{eq;21.8.21.30}
 \kappa_{p,1}^{\ast}(u^{\lambda}_{p,\ell,\alpha})
=u^{\lambda}_{p,\ell,\alpha}\cdot
\alpha
(1+|\lambda|^2)^{\ell/p}
(\beta_1+2\sqrt{-1}\lambda T)^{-\ell/p}
\exp\Bigl(
\frac{\ell}{p}G(2\sqrt{-1}\lambda T/\beta_1)
\Bigr),
\end{equation}
on $\kappa_{p,1}^{-1}(\nbigu^{\lambda}_{p})
\cap \nbigu^{\lambda}_p$,
where $G(x)=1-x^{-1}\log(1+x)$.
For any $P\in H^{\lambda\cov}_{\infty,p}$,
and for any relatively compact neighbourhood
$U_P$ of $P$ in $\nbigu^{\lambda}_{p}$,
there exists $C(P)>1$ such that
\begin{equation}
\label{eq;21.8.21.31}
 C(P)^{-1}|w|^{-\ell t_1/(pT)}
 \leq
 |u^{\lambda}_{p,\ell,\alpha}|_{h_{\LL^{\cov},p,\ell,\alpha}}
 \leq
 C(P)|w|^{-\ell t_1/(pT)}.
\end{equation}
If $\lambda=0$, we may choose $\nbigu^{0}_p=Y^{0\cov}_p$,
and 
(\ref{eq;21.8.21.30})
is rewritten as
$\kappa_{p,1}^{\ast}(u^{0}_{p,\ell,\alpha})
=\alpha w_{p}^{-\ell}u_{p,\ell,\alpha}$.

\section{Examples of monopoles induced by
tame harmonic bundles (1)}
\label{subsection;21.8.13.22}

Take $(a,\alpha)\in\real\times\cnum$.
We have the harmonic bundle
given as 
$\nbigl(a,\alpha)=\nbigo_{\cnum_{w_p}^{\ast}}e_0$
with
the Higgs field $\theta_0 e_0=e_0(-\alpha dw/w)$
and the metric $h_0$ determined by
$h_0(e_0,e_0)=|w|^{2a}$.
\index{harmonic bundle $(\nbigl(a,\alpha),\theta_0,h_0)$}
We have $\del_we_0=e_0\cdot adw/w$
and $\theta_0^{\dagger}e_0=e_0(-\alphabar\,d\wbar/\wbar)$.

We obtain the monopole
$(\vecnbigl(a,\alpha),h,\nabla,\phi)$
on $Y_p^{0\ast}$ induced by
$(\nbigl(a,\alpha),\theta_0,h_0)$
as in \S\ref{subsection;17.9.29.1}.
\index{monopole $\vecnbigl(a,\alpha)$}
We describe it explicitly.

We have $\vecnbigl(a,\alpha)=\Psi_p^{-1}\nbigl(a,\alpha)$,
which is equipped with the induced metric $h=\Psi_p^{-1}(h_0)$.
We set
$e:=\Psi_p^{-1}(e_0)$.
The unitary connection $\nabla$ and the Higgs field $\phi$
are given as follows:
\[
 \nabla e=
 e\bigl(aw^{-1}dw
 +\sqrt{-1}(\alpha w^{-1}+\alphabar \wbar^{-1})\,dt\bigr),
\quad\quad
 \phi e=e\bigl(-\alpha w^{-1}+\alphabar \wbar^{-1}\bigr).
\]

\subsection{The mini-holomorphic bundle
at $\lambda=0$}

Let us describe the mini-holomorphic bundle
$\vecnbigl^0(a,\alpha)$
on $Y_p^{0\ast}$
underlying
$(\vecnbigl(a,\alpha),h,\nabla,\phi)$.
\index{mini-holomorphic bundle $\vecnbigl^0(a,\alpha)$}

\begin{prop}
There exists a mini-holomorphic frame 
 $u^0_{a,\alpha}$ of
 the mini-holomorphic bundle
$(\varpi^0_p)^{-1}\vecnbigl^0(a,\alpha)$
on $Y_p^{0\cov\ast}$
such that the following holds.
\begin{description}
\item[(A1)]
 For any $P\in H^0_{\infty,p}$,
 there exists a neighbourhood $U_P$ of $P$
 in $Y_p^{0\cov\ast}$
 and  a constant $C(P)\geq 1$
 such that the following holds:
\[
 C(P)^{-1}|w|^{a}
\leq
 \bigl|u^0_{a,\alpha}\bigr|_h
\leq
 C(P)|w|^{a}.
\]
\item[(A2)]
$\kappa_{p,1}^{\ast}
 u^0_{a,\alpha}
=u^0_{a,\alpha}
 \exp\bigl(
 -2\sqrt{-1}\alpha T w^{-1}
 \bigr)$. 
\end{description}
\end{prop}
\pf
We obtain
$\nabla_{\wbar}e=0$,
$\nabla_{w}e=eaw^{-1}$,
and
$\nabla_te=e\bigl(\sqrt{-1}(\alpha w^{-1}+\alphabar \wbar^{-1})\bigr)$.
We also obtain
\[
 (\nabla_t-\sqrt{-1}\phi)e
=e(2\sqrt{-1}\alpha w^{-1}).
\]
Hence, we obtain the following mini-holomorphic frame
of $(\varpi_p^{0})^{-1}\vecnbigl^0(a,\alpha)$:
\[
 u^0_{a,\alpha}:=
 (\varpi^0_p)^{-1}(e)\exp\bigl(
 2\sqrt{-1}\alpha t w^{-1}
 \bigr).
\]
It has the desired properties.
\hfill\qed

\subsection{The mini-holomorphic bundle at $\lambda$}

The $\lambda$-flat bundle
associated with the harmonic bundle $(\nbigl(a,\alpha),\theta_0,h_0)$
is described as follows:
\[
 (\delbar+\lambda\theta_0^{\dagger})e_0
=e_0\bigl(-\alphabar\lambda\bigr) d\wbar/\wbar,
\quad\quad
 (\lambda\del+\theta_0)e_0
=e_0\bigl(\lambda a-\alpha\bigr)dw/w.
\]
We set
$v_0^{\lambda}:=
 e_0\exp\bigl(\alphabar\lambda \log|w|^2\bigr)$
on $\cnum_{w_p}^{\ast}$,
which is holomorphic with respect to
$\delbar+\lambda\theta_0^{\dagger}$.
We obtain
\[
 \DDlambda v^{\lambda}_0
 =v^{\lambda}_0
 \bigl(-\alpha+\lambda a+\lambda^2\alphabar\bigr)dw/w,
\quad\quad
 \bigl|v^{\lambda}_0\bigr|_{h_0}
=|w|^{a+2\Re(\alphabar\lambda)}.
\]

Let us describe the mini-holomorphic bundle
$\vecnbigl^{\lambda}(a,\alpha)$
on $Y_p^{\lambda\cov\ast}$
underlying 
$(\vecnbigl(a,\alpha),h,\nabla,\phi)$.
\index{mini-holomorphic bundle $\vecnbigl^{\lambda}(a,\alpha)$}
\begin{prop}
There exist a neighbourhood $\nbigu^{\lambda}_p$
of $H^{\lambda\cov}_{\infty,p}$
in $Y^{\lambda\cov}_{p}$
and a mini-holomorphic frame $u^{\lambda}_{a,\alpha}$ of
$(\varpi_p^{\lambda})^{-1}\vecnbigl^{\lambda}(a,\alpha)$
on $Y_p^{\ast\cov}$
such that the following holds.
\begin{description}
\item[(A1)]
For any $P\in H^{\lambda\cov}_{\infty,p}$,
there exist a constant $C(P)\geq 1$
and a neighbourhood $U_P$ of $P$
in $\nbigu^{\lambda}_p$
such that the following holds on
$U_P\setminus H^{\lambda\cov}_{\infty,p}$:
\[
 C(P)^{-1}|w|^{a+2\Re(\lambda\alphabar)}
 \leq
 \bigl|
 u^{\lambda}_{a,\alpha}
 \bigr|_{h}
 \leq 
 C(P)|w|^{a+2\Re(\lambda\alphabar)}.
\]
\item[(A2)]
On $\bigl(\kappa_{p,1}^{-1}(\nbigu^{\lambda}_p)
\cap
\nbigu^{\lambda}_p\bigr)
\setminus
H^{\lambda\cov}_{\infty,p}$,
the following holds:
\[
 \kappa_{p,1}^{\ast}
 u^{\lambda}_{a,\alpha}
=u^{\lambda}_{a,\alpha}
 \exp\Bigl(
 \bigl(-\alpha\lambda^{-1}+a+\lambda\alphabar\bigr)
 \log\bigl(
 \beta_1^{-1}(\beta_1+2\sqrt{-1}\lambda T)
 \bigr)
 \Bigr).
\]
\end{description}
\end{prop}
\pf
We set $v^{\lambda}:=\Psi_p^{-1}(v_0^{\lambda})$.
By Corollary \ref{cor;21.8.12.6}
with the formulas (\ref{eq;21.8.12.30}) and (\ref{eq;21.8.12.31}),
we obtain $\del_{\betabar_1}v^{\lambda}=0$
and 
\[
 \del_{t_1}v^{\lambda}
=v^{\lambda}\Bigl(
 -2\sqrt{-1}
 \frac{-\alpha+a\lambda+\lambda^2\alphabar}
 {\beta_1-2\sqrt{-1}\lambda t_1}
 \Bigr).
\]
On an appropriate open neighbourhood $\nbigu^{\lambda}_p$
of $H^{\lambda}_{\infty,p}$,
we obtain the following mini-holomorphic frame
on $\nbigu^{\lambda}_p\setminus H^{\lambda\cov}_{\infty,p}$:
\[
 u^{\lambda}_{a,\alpha}
 =(\varpi_p^{\lambda})^{-1}
 (v^{\lambda})\exp\Bigl(
-\lambda^{-1} \bigl(
-\alpha+\lambda a
+\lambda^2\alphabar
 \bigr)
 \log\bigl(\beta_1^{-1}(\beta_1-2\sqrt{-1}\lambda t_1)\bigr)
 \Bigr).
\]
It has the desired properties.
\hfill\qed

\section{Examples of monopoles induced by
tame harmonic bundles (2)}

\subsection{Case of $\lambda=0$}
\label{subsection;17.10.4.20}

Let $N$ be the $r\times r$ matrix whose
$(i,j)$-entries are $1$ $(i=j+1)$
or $0$ $(i\neq j+1)$.
Let $U^{\ast}_{w,p}:=\{w_p\in\cnum\,|\,|w_p^p|>R\}$ for some $R>1$.
We set
$V_N:=\bigoplus_{i=1}^r \nbigo_{U^{\ast}_{w,p}}v_{0,i}$
on which we define the Higgs field $\theta_N$
by
$\theta_N\vecv_0=\vecv_0 N\,dw/w$.
Let $h_V$ be a harmonic metric of the Higgs bundle $(V_N,\theta_N)$
such that 
$h_V$ is mutually bounded with
the metric $h_{V,1}$
determined as follows:
\[
 h_{V,1}(v_i,v_j)=
\left\{
 \begin{array}{ll}
 (\log|w|)^{r+1-2i}  &(i=j)\\
 0 &(i\neq j)
 \end{array}
\right.
\]
Note that such a harmonic metric $h_V$ exists.
(See \cite{Simpson90} and \cite[\S6]{mochi2}.)

We obtain the induced monopole $(E_N,h,\nabla,\phi)$
on $\Psi_p^{-1}(U^{\ast}_{w,p})\subset Y^{0\ast}_p$.
Let us describe the underlying mini-holomorphic bundle
$(E^0_N,\delbar)$ on 
$\Psi_p^{-1}(U^{\ast}_{w,p})$.
\index{monopole $(E_N,h,\nabla,\phi)$}

\begin{prop}
There exist a neighbourhood $\nbigu^{0}_p$
of $H^{0\cov}_{\infty,p}$ in $Y^{0\cov}_p$
and a mini-holomorphic frame
$\vecu^0_N=(u^0_{N,1},\ldots,u^0_{N,r})$ of 
$(\varpi_p^0)^{-1}(E^0_N,\delbar)$
on $\nbigu^{0}_p\setminus H^{0\cov}_{\infty,p}$
satisfying the following conditions:
\begin{description}
\item[(A1)]
 Let $h_0$ be the metric determined by
\[
 h_{0}(u^0_{N,i},u^0_{N,j})=
\left\{
 \begin{array}{ll}
 (\log|w|)^{r+1-2i}  &(i=j)\\
 0 &(i\neq j)
 \end{array}
\right.
\]
Then, $h_0$ and $h$ are mutually bounded
locally around any point of $H^{0\cov}_{\infty,p}$.
\item[(A2)]
On $(\kappa_{p,1}^{-1}(\nbigu^0_p)\cap\nbigu^0_p)
 \setminus H^{0\cov}_{\infty,p}$,
we obtain
$\kappa_{p,1}^{\ast}\vecu^0_N
=\vecu^0_N\exp\bigl(2\sqrt{-1}NTw^{-1}\bigr)$.
\end{description}
\end{prop}
\pf
Let $\vecv^{0}$ be the pull back $\Psi_p^{-1}(\vecv_0)$.
We have
$\del_{E^0_N,\wbar}\vecv^0=0$
and 
$\del_{E^0_N,t}\vecv^0
=\vecv^0(-2\sqrt{-1}Nw^{-1})$.
We obtain the following frame:
\[
 \vecu^0_N=
 (\varpi_p^0)^{-1}
 (\vecv^0_N)\exp\Bigl(2\sqrt{-1}Ntw^{-1}\Bigr).
\]
It has the desired property.
\hfill\qed

\subsection{Case of $\lambda\neq 0$}

We consider the $\lambda$-connection
$\DDlambda_N$ on $V_N$
given by
$\DDlambda_N\vecv_0=\vecv_0 N\,dw/w$.
Let $h_V$ be a harmonic metric
of the $\lambda$-flat bundle
$(V_N,\DDlambda_N)$
such that $h_V$ is mutually bounded
with the metric $h_{V,1}$ determined as follows:
\[
 h_{V,1}(v_i,v_j)=
\left\{
 \begin{array}{ll}
 (\log|w|)^{r+1-2i}  &(i=j)\\
 0 &(i\neq j)
 \end{array}
\right.
\]
Note that such a harmonic metric $h_V$ exists
(See \cite{Simpson90} and \cite[\S6]{mochi2}.)
We obtain the induced monopole $(E_N,h,\nabla,\phi)$
on $\Psi_p^{-1}(U_{w,p}^{\ast})$.
Let us describe the underlying mini-holomorphic bundle
$(E^{\lambda}_N,\delbar)$ on $Y^{\lambda\ast}_p$.

\begin{prop}
There exist a neighbourhood $\nbigu^{\lambda}_p$
of $H^{\lambda\cov}_{\infty,p}$ in $Y^{\lambda\cov}_p$
and a mini-holomorphic frame
$\vecu^{\lambda}_N=
 (u^{\lambda}_{N,1},\ldots,u^{\lambda}_{N,r})$ 
of $(\varpi_p^{\lambda})^{-1}(E^{\lambda}_N,\delbar)$
on $\nbigu^{\lambda}_p\setminus H^{\lambda\cov}_{\infty,p}$
with the following frame:
\begin{description}
\item[(A1)]
 Let $h_0$ be the metric determined by
\[
 h_{0}(u^{\lambda}_{N,i},u^{\lambda}_{N,j})=
\left\{
 \begin{array}{ll}
 (\log|w|)^{r+1-2i}  &(i=j)\\
 0 &(i\neq j)
 \end{array}
\right.
\]
Then, $h_0$ and $h$ are mutually bounded
locally around any point of $H^{\lambda\cov}_{\infty,p}$.
\item[(A2)]
On $(\kappa_{p,1}^{-1}(\nbigu^{\lambda}_p)\cap\nbigu^{\lambda}_p)
 \setminus H^{\lambda\cov}_{\infty,p}$,
 we obtain
\[
 \kappa_{p,1}^{\ast}\vecu^{\lambda}_N
=\vecu^{\lambda}_N 
 \exp\bigl(
 -\lambda^{-1}\log(\beta_1(\beta_1+2\sqrt{-1}\lambda T)^{-1})
 \bigr).
\]  
\end{description}
\end{prop}
\pf
Let $\vecv^{\lambda}$ be the pull back of $\vecv_0$.
By the formulas (\ref{eq;21.8.12.30}) and (\ref{eq;21.8.12.31}),
we have
$\del_{E^{\lambda}_N,\betabar_1}\vecv^{\lambda}=0$
and 
$\del_{E^{\lambda}_N,t_1}\vecv^{\lambda}
=\vecv^{\lambda}(-2\sqrt{-1}N(\beta_1-2\sqrt{-1}\lambda t_1)^{-1})$.
We obtain the following frame:
\[
 \vecu^{\lambda}_N=
 (\varpi_p^{\lambda})^{-1}(\vecv^{\lambda}_N)
 \exp\Bigl(-\lambda^{-1}N
 \log\bigl(\beta_1^{-1}(\beta_1-2\sqrt{-1}\lambda t_1)\bigr)
 \Bigr).
\]
Then, it has the desired property.
\hfill\qed

\chapter[Asymptotic behaviour around infinity]{Asymptotic behaviour of periodic monopoles
 around infinity}
\label{section;17.10.5.40}

In \S\ref{subsection;20.7.31.40},
we introduce some notation
for covering spaces
and the decomposition of a mini-holomorphic bundle
induced by the decomposition of the spectral curve.

In \S\ref{subsection;16.9.28.1},
first, we explain the main result of this section
(Theorem \ref{thm;16.9.20.22}),
which is roughly the asymptotic orthogonality
of the spectral decomposition in the slope level
for a periodic monopole satisfying
Condition \ref{condition;21.8.21.1}
and Condition \ref{condition;21.8.21.2}.
Note that if Condition \ref{condition;21.8.21.1}
is satisfied,
and Condition \ref{condition;21.8.21.2} is also satisfied
on an appropriate ramified covering.
Then, with the aid of the estimate
for asymptotic doubly periodic instantons
in \S\ref{subsection;20.8.1.1}--\S\ref{subsection;20.8.1.2}
developed in \cite{Mochizuki-doubly-periodic},
we explain the second main result of this section
(Proposition \ref{prop;17.10.5.3}).
As a result, with the aid of the estimates for
asymptotic harmonic bundles in 
\S\ref{subsection;17.10.21.10}
and \cite[\S5.5]{Mochizuki-doubly-periodic},
we may understand the asymptotic behaviour of
periodic monopoles.
We also note that
Condition \ref{condition;21.8.21.1}
is equivalent to the GCK-condition
(Proposition \ref{prop;20.7.29.20}).

The rest of this subsection is devoted to the proof of
Theorem \ref{thm;16.9.20.22}.
In \S\ref{subsection;17.9.28.10},
we shall prove that
the unitary connection asymptotically preserves
an orthogonal decomposition associated with the Higgs field,
under the assumption on the behaviour of the eigenvalues of
the Higgs field.
After the preliminary
in \S\ref{subsection;20.7.31.41}--\S\ref{subsection;20.7.31.42},
we shall prove that the assumption is satisfied
and that the spectral decomposition is asymptotically orthogonal.

\section{Preliminary}
\label{subsection;20.7.31.40}

\subsection{Notation}
\label{subsection;20.8.8.110}

Set $U_w(R):=\{|w|>R\}\cup\{\infty\}$ in $\proj^1_w$
for any  $R>0$.
\index{space $U_w(R)$}
We set $U_w^{\ast}(R):=U_w(R)\setminus\{\infty\}$.
\index{space $U^{\ast}_w(R)$}
We fix a $q$-th root $w_q$ of $w$ for any $q\in\seisuu_{\geq 1}$
such that $w_{qm}^m=w_q$ for any $m\in\seisuu_{\geq 1}$.
\index{variable $w_q$}
For any $q\in\seisuu_{\geq 1}$,
let $\proj^1_{w_q}\lrarr \proj^1_w$
be the ramified covering map induced by
$w_q\longmapsto w_q^q$.
Let $U_{w,q}(R)$ and $U_{w,q}^{\ast}(R)$
denote the inverse image of
$U_w(R)$ and $U_w^{\ast}(R)$,
respectively.
\index{space $U_{w,q}(R)$}
\index{space $U^{\ast}_{w,q}(R)$}
For any $p\in q\seisuu_{\geq 1}$,
there exists the natural morphism
$\varphi_{q,p}:U_{w,p}(R)\lrarr U_{w,q}(R)$.
\index{map $\varphi_{q,p}$}

Take $T>0$.
Let $\nbigb^0_q(R)$ denote 
$(\real_t/T\seisuu)\times U_{w,q}(R)$
equipped with the natural mini-holomorphic structure.
\index{mini-complex manifold $\nbigb^0_q(R)$}
Similarly, we set
$\nbigb^{0\ast}_q(R):=
 (\real/T\seisuu)\times U_{w,q}^{\ast}(R)
 \subset \nbigb^0_q(R)$,
which is equipped with the metric $dt\,dt+dw\,d\wbar$.
\index{mini-complex manifold $\nbigb^{0\ast}_q(R)$}
We also put
$H^0_{\infty,q}:=S^1_T\times\{\infty\}=
 \nbigb^0_q(R)\setminus \nbigb^{0\ast}_q(R)$.
\index{space $H^0_{\infty,q}$}
 
On $U_{w,q}(R)$ and $\nbigb^0_q(R)$,
the complex vector fields
$q^{-1}w_q^{-q+1}\del_{w_q}$
and 
$q^{-1}\wbar_q^{-q+1}\del_{\wbar_q}$
are denoted by $\del_{w}$ and $\del_{\wbar}$.
\index{vector fields $\del_w$, $\del_{\wbar}$}
Similarly,
the complex $1$-forms
$qw_q^{q-1}dw_q$
and 
$q\wbar_q^{q-1}d\wbar_q$
are denoted by
$dw$ and $d\wbar$.
\index{one forms $dw$, $d\wbar$}
The function $w_q^q$ is denoted by $w$.
\index{function $w$}

We also set
$\nbigb^{0,\cov}_q(R):=\real_t\times U_{w,q}(R)$,
$\nbigb^{0,\cov,\ast}_q(R):=\real_t\times U^{\ast}_{w,q}(R)$
and $H^{0,\cov}_q:=\real_t\times\{\infty\}$.
\index{space $\nbigb^{0\cov}_q(R)$}
\index{space $\nbigb^{0\cov\ast}_q(R)$}
\index{space $H^{0\cov}_q$}
Let $\varpi^0_q:\nbigb^{0\,\cov}_q(R)\lrarr \nbigb^{0}_q(R)$
denote the map induced by $\real\lrarr\real/T\seisuu$.
The restriction
$\nbigb^{0\cov\ast}_q(R)\lrarr\nbigb^{0\ast}_{q}(R)$
is denoted by $\varpi^0_q$ or simply $\varpi_q$.
\index{projection $\varpi^0_q$}
\index{projection $\varpi_q$}

Let $\kappa$ be the free $\seisuu$-action
on $\nbigb^{0,\cov}_{q}(R)$
defined by 
$\kappa_n(t,w_q)=(t+nT,w_q)$ $(n\in\seisuu)$.
The quotient space is naturally identified with
$\nbigb^{0}_q(R)$.
\index{action $\kappa$}

By setting $(x,y)=(\Re(w),\Image(w))=(\Re(w_q^q),\Image(w_q^q))$,
we obtain the local coordinate systems $(t,x,y)$ on 
$\nbigb^{0\ast}_{q}(R)$.
\index{real coordinate system $(t,x,y)$}
In particular,
we obtain the vector fields
$\del_{a}$ $(a=t,x,y)$,
which are globally well defined.

\begin{rem}
The variables $x$ and $y$
are different from the variables $x$ and $y$
in {\rm\S\ref{section;20.8.8.20}--\S\ref{section;20.8.8.21}}.
\hfill\qed
\end{rem}

For $R_1>R$ and $p\in q\seisuu_{\geq 1}$,
we have the naturally defined ramified coverings
$\nbigr_{q,p}:\nbigb^0_{p}(R_1)\lrarr\nbigb^0_q(R)$,
and 
$\nbigr_{q,p}:\nbigb^{0\cov}_{p}(R_1)\lrarr\nbigb^{0\cov}_q(R)$.
\index{morphism $\nbigr_{q,p}$}

For a vector bundle $E$ with a Hermitian metric $h$
on $\nbigb^{0\ast}_q(R)$ (resp. $U_{w,q}^{\ast}$)
and for 
an $\End(E)$-valued differential form $s$ on
$\nbigb^{0\ast}_q(R)$ (resp. $U_{w,q}^{\ast}$),
let $|s|_{h}$ denote the function on $\nbigb^{0\ast}_q(R)$
(resp. $U_{w,q}^{\ast}$)
obtained as the point-wise norm of $s$ with respect to $h$
and the natural Riemannian metric of $\nbigb^{0\ast}_q(R)$
(resp. $U_{w,q}^{\ast}$).
We omit to denote the dependence on the metrics of the base spaces
because they are fixed.
\index{norm \mbox{$|s|_h$}}

\subsection{Decomposition of holomorphic bundles
 with an automorphism}

Let $V$ be a locally free $\nbigo_{U_{w,q}^{\ast}(R)}$-module
of finite rank with an automorphism $F$.
For each $Q\in U_{w,q}^{\ast}(R)$,
let $\Sp(F_{|Q})$ denote the set of the eigenvalues
of $F_{|Q}$.
Let 
$\Sp(F)\subset U_{w,q}^{\ast}(R)\times\cnum^{\ast}$
denote the spectral curve
given by
$\Sp(F)=\coprod_{Q\in U_{w,q}^{\ast}(R)}\Sp(F_{|Q})$.
It is a closed complex analytic curve in $U^{\ast}_{w,q}(R)\times\cnum^{\ast}$.
\index{spectral curve $\Sp(F)$}

\begin{lem}
\label{lem;17.9.27.1}
Suppose that the closure of $\Sp(F)$ 
in $U_{w,q}(R)\times\proj^1$
is complex analytic.
Then, 
there exist $R_1>R$, $p\in q\seisuu_{\geq 1}$,
and a decomposition
\[
 \varphi_{q,p}^{-1}(V,F)
=\bigoplus_{(\ell,\alpha)\in\seisuu\times\cnum^{\ast}}
  (V_{p,\ell,\alpha},F_{p,\ell,\alpha})
\]
on $U_{w,p}(R_1)$,
such that the following holds
for some $C>0$ and $\delta>0$:
\begin{itemize}
\item
For any $Q\in U_{w,p}^{\ast}(R_1)$,
 any eigenvalue $\rho_Q$ of
 $F_{p,\ell,\alpha|Q}$ satisfies
\[
 \Bigl|
  \rho_Q\cdot
 \bigl(
 \alpha w_p(Q)^{-\ell}
 \bigr)^{-1}
-1
\Bigr|
\leq C|w_p(Q)|^{-\delta}.
\]
\end{itemize}
\end{lem}
\pf
By replacing $U_{w,q}(R)$ with
$U_{w,q}(R_2)$ for some large $R_2>0$,
we may assume that the irreducible decomposition
$\overline{\Sp(F)}
=\bigcup_{i\in\Lambda}\overline{\Sp(F)}_i$
satisfies the following condition.
\begin{itemize}
\item
The germs of $\overline{\Sp(F)}_i$
at $\overline{\Sp(F)}_i\cap(\{\infty\}\times\proj^1)$
are also irreducible.
\end{itemize}
There exists the corresponding decomposition
$(V,F)
=\bigoplus_{i\in\Lambda}(V_i,F_i)$.
Let $Z\lrarr \overline{\Sp(F)}$
be the normalization.
Let $Z=\coprod_{i\in\Lambda} Z_i$
denote the decomposition
into the connected components.

The induced maps
$a_i:Z_i\lrarr U_{w,q}(R)$
are coverings ramified at $\infty$.
Let $m_i$ be the ramification indexes.
For each $i$,
there exits a coordinate $\zeta_{i}$ of $Z_i$
such that $a_i^{\ast}w_q=\zeta_i^{m_i}$.
We set $Z_i^{\ast}=a_i^{-1}(U_{w,q}^{\ast}(R))$.
The restriction $Z_i^{\ast}\lrarr U_{w,q}(R)$
is also denoted by $a_i$.
Let $\upsilon_i:Z_i^{\ast}\lrarr \cnum^{\ast}$
be the map obtained as the composition of
$Z_i^{\ast}\lrarr U_{w,q}^{\ast}(R)\times \cnum^{\ast}
 \lrarr \cnum^{\ast}$.

On $Z_i^{\ast}$,
there exists an $\nbigo_{Z_i^{\ast}}$-module
with an automorphism $(\Vtilde_i,\Ftilde_i)$
such that
$(V_i,F_i)=a_{i\ast}(\Vtilde_i,\Ftilde)$,
and that
$\Sp(\Ftilde_i)\subset Z_i^{\ast}\times\cnum^{\ast}$
is equal to the image of
$\id\times\upsilon_i:
 Z_i^{\ast}\lrarr Z_i^{\ast}\times\cnum^{\ast}$.
There exist an integer $\ell_i$,
a non-zero complex number $\gamma_i$
and a holomorphic function $\gminib_i(u)$ on 
a neighbourhood of $u=0$ with $\gminib_i(0)=0$,
such that 
$\upsilon_i(\zeta_i)=
 \gamma_i
 \zeta_i^{-\ell_i}
 \exp\bigl(
 2\sqrt{-1}T\gminib_i(\zeta_i^{-1})
 \bigr)$.
Then, we obtain the claim of the lemma.
\hfill\qed

\subsection{Some basic mini-holomorphic
bundles on $\nbigb^{0\ast}_q(R)$}
\label{subsection;17.9.26.2}

\subsubsection{Basic rank one bundles}
\label{subsection;21.8.16.20}

For $q\in\seisuu_{>0}$ and $(\ell,\alpha)\in\seisuu\times\cnum^{\ast}$,
we introduced the monopole
$(\LL^{\ast}_q(\ell,\alpha),
h_{\LL,q,\ell,\alpha},\nabla_{\LL,q,\ell,\alpha},
\phi_{\LL,q,\ell,\alpha})$ on $Y_q^{0\ast}$
in (\ref{eq;21.8.16.12}).
\index{monopole $\LL^{\ast}_q(\ell,\alpha)$}
Let $\LL^{0\ast}_q(\ell,\alpha)$
denote the mini-holomorphic bundle on $Y_q^{0\ast}$
underlying the monopole.
\index{mini-holomorphic bundle $\LL^{0\ast}_q(\ell,\alpha)$}
As explained in \S\ref{subsection;21.8.14.2},
it extends to
an $\nbigo_{Y_q^0}(\ast H^0_{\infty,q})$-module
$\LL^0_q(\ell,\alpha)$.
\index{$\nbigo_{Y_q^0}(\ast H^0_{\infty,q})$-module
$\LL^0_q(\ell,\alpha)$}

Let $\LL^{0\cov}_q(\ell,\alpha)$ denote the pull back
of $\LL^0_q(\ell,\alpha)$ by
$\varpi^0:Y_q^{0\cov}\lrarr Y_q^0$.
\index{sheaf $\LL^{0\cov}_q(\ell,\alpha)$}
As explained in \S\ref{subsection;21.8.14.2},
it is equipped with a global mini-holomorphic frame
$e_{q,\ell,\alpha}$,
i.e.,
$\LL^{0\cov}_q(\ell,\alpha):=
 \nbigo_{Y^{0,\cov}_{q}}(\ast H^{0,\cov}_{\infty,q})
 \,e_{q,\ell,\alpha}$,
for which $\seisuu$-action is described as follows:
\[
 \kappa_1^{\ast}(e_{q,\ell,\alpha})=
  w_q^{-\ell}\alpha e_{q,\ell,\alpha}.
\]
\index{mini-holomorphic bundle $\LL^{0\cov}_q(\ell,\alpha)$}
Let $\LL^{0\cov\ast}_{q}(\ell,\alpha)$
denote the restriction of
$\LL^{0\cov}_q(\ell,\alpha)$ to $Y^{0\cov\ast}_q$,
which is equal to the pull back of
$\LL^{0\ast}_q(\ell,\alpha)$.
It is equipped with the metric
$h_{\LL^{\cov},q,\ell,\alpha}$
obtained as the pull back of
$h_{\LL,q,\ell,\alpha}$,
for which the following holds:
\index{metric $h_{\LL^{\cov},q,\ell,\alpha}$}
\[
  h_{\LL^{\cov},q,\ell,\alpha}(e_{q,\ell,\alpha},e_{q,\ell,\alpha})
=\bigl|\alpha w_q^{-\ell}\bigr|^{2t/T}.
\]

The restriction of
$(\LL^{0\ast}_q(\ell,\alpha),h_{\LL,q,\ell,\alpha})$
to $\nbigb^{0\ast}_q(R)$ is denoted by the same notation
$(\LL^{0\ast}_q(\ell,\alpha),h_{\LL,q,\ell,\alpha})$.
Similarly,
the restriction of
$(\LL^{0\cov\ast}_q(\ell,\alpha),h_{\LL^{\cov},q,\ell,\alpha})$
to $\nbigb^{0\cov\ast}_q(R)$ is also denoted by
$(\LL^{0\cov\ast}_q(\ell,\alpha),h_{\LL^{\cov},q,\ell,\alpha})$.

\subsubsection{Mini-holomorphic bundles
induced by holomorphic bundles with automorphism}

Let $\Psi_q$ denote the projection
$\nbigb_q^{0\ast}(R)\lrarr U_{w,q}^{\ast}(R)$.
Let $V$ be a holomorphic vector bundle 
with an endomorphism $g$ on $U_{w,q}^{\ast}$.
We obtain the induced mini-holomorphic bundle
$\Psi_q^{\ast}(V,g)$ on $\nbigb_q^{0\ast}(R)$
as explained in \S\ref{subsection;17.9.29.2}.
\index{mini-holomorphic bundle
$\Psi_q^{\ast}(V,g)$}

\subsection{Decomposition of mini-holomorphic bundles}
\label{subsection;17.10.5.20}

Let $\nbigv$ be a locally free $\nbigo_{\nbigb^{0\ast}_q(R)}$-module
of finite rank.
We obtain the automorphism $F$ of
$\nbigv_{|\{0\}\times U_{w,q}^{\ast}(R)}$
obtained as the monodromy 
along each loop $S^1_T\times\{w\}$
(see \S\ref{subsection;17.10.5.1}).
We set $\Sp(\nbigv):=\Sp(F)$.
\index{spectral curve $\Sp(\nbigv)$}

\begin{lem}
\label{lem;17.9.27.5}
Suppose that the closure of $\Sp(\nbigv)$
in $U_{w,q}(R)\times\proj^1$ is a complex analytic curve.
Then, the following holds:
\begin{itemize}
\item
There exist
$R_1>R$, $p\in q\seisuu_{\geq 1}$,
a finite subset
$S\subset\seisuu\times\cnum^{\ast}$,
a tuple of 
locally free $\nbigo_{U_{w,q}^{\ast}}$-module
with an endomorphism
$(V_{p,\ell,\alpha},f_{p,\ell,\alpha})$
$((\ell,\alpha)\in S)$,
and the following isomorphism
on $\nbigb_p^{0\ast}(R_1)$:
\index{sheaf with an endomorphism
$(V_{p,\ell,\alpha},f_{p,\ell,\alpha})$}
\begin{equation}
\label{eq;17.9.27.2}
 \nbigr_{q,p}^{\ast}\nbigv
\simeq
 \bigoplus_{(\ell,\alpha)\in S}
 \LL^{0\ast}_{p}(\ell,\alpha)
 \otimes
 \Psi_p^{\ast}(V_{p,\ell,\alpha},f_{p,\ell,\alpha}).
     \end{equation}
\end{itemize}
\end{lem}
\pf
Applying Lemma \ref{lem;17.9.27.1}
to $V:=\nbigv_{|\{0\}\times U_{w,q}^{\ast}}$
with the automorphism $F$,
we obtain a decomposition as in Lemma \ref{lem;17.9.27.1}.
On each $Z^{\ast}_i$,
there exists the endomorphism 
$\ftilde_i$ of $\Vtilde_i$
such that 
(i) $\exp(2\sqrt{-1}T\ftilde_i)=\Ftilde_i$,
(ii) $\Sp(\ftilde_i)\subset Z_i^{\ast}\times\cnum$ is equal to 
the image of $\gminib_i:Z_i^{\ast}\lrarr \cnum$.
Then, we obtain the claim of the lemma.
\hfill\qed

\vspace{.1in}
We identify the both sides of
(\ref{eq;17.9.27.2}) by the isomorphism.
By setting
\begin{equation}
 \nbigv_{p,\ell}:=
 \bigoplus_{\alpha\in\cnum^{\ast}}
 \LL^{0\ast}_{p}(\ell,\alpha)
 \otimes
 \Psi_p^{\ast}(V_{p,\ell,\alpha},f_{p,\ell,\alpha}),
\end{equation}
we obtain the following decomposition
of mini-holomorphic bundles:
\begin{equation}
\label{eq;17.9.27.3}
 \nbigr_{q,p}^{\ast}\nbigv
=\bigoplus_{\ell\in\seisuu}\nbigv_{p,\ell}.
\end{equation}
As a characterization of the decomposition
(\ref{eq;17.9.27.3}),
there exists $C>1$ such that the following holds.
\begin{itemize}
\item
 For any eigenvalue $\rho_Q$ of the monodromy of
 $\nbigv_{p,\ell}$ along $S^1_T\times\{Q\}$,
we obtain
$C^{-1}<|\rho_Q\big/w_p(Q)^{-\ell}|<C$.
\end{itemize}

\section{Estimate of periodic monopoles around infinity}
\label{subsection;16.9.28.1}

\subsection{Setting}

We consider a monopole
$(E,h,\nabla,\phi)$ on $\nbigb^{0\ast}_{q}(R)$
for some $R>0$.
Let $(E,\delbar_E)$ denote
the underlying mini-holomorphic bundle
on $\nbigb^{0\ast}_q(R)$.
We assume the following condition.
\begin{condition}
\label{condition;21.8.21.1} \mbox{{}}
\begin{description}
\item[(B1)]
 $\bigl|\nabla(\phi)_{|(t,w)}\bigr|_{h}\to 0$
 when $|w|\to\infty$.
\item[(B2)]
 The closure of $\Sp(E,\delbar_E)$ 
 in $U_{w,q}\times\proj^1$
 is a complex analytic curve.
\end{description}
 \end{condition}

Moreover,
we also assume the following condition.
\begin{condition}
\label{condition;21.8.21.2}
There exists a decomposition
of mini-holomorphic bundles
\begin{equation}
 \label{eq;17.9.27.4}
 (E,\delbar_E)
=\bigoplus_{(\ell,\alpha)\in\seisuu\times\cnum^{\ast}}
 (E_{\ell,\alpha},\delbar_{E_{\ell,\alpha}}),
\end{equation}
such that the following holds
for some $C>0$ and $\delta>0$:
\begin{itemize}
\item
 For any eigenvalue $\rho_Q$ of
 the monodromy of $(E_{\ell,\alpha},\delbar_{E_{\ell,\alpha}})$
 along $S^1_T\times\{Q\}$,
 we have
\[
  \bigl|
 \rho_Q\cdot(\alpha w_q(Q)^{-\ell})^{-1}-1
 \bigr|
\leq C|w_q(Q)|^{-\delta}.
\]
\end{itemize}
\end{condition}
 
\begin{rem}
\label{rem;21.7.7.2}
Condition {\rm\ref{condition;21.8.21.2}} implies
{\bf (B2)} in Condition {\rm\ref{condition;21.8.21.1}}.
Conversely, if {\bf (B2)} is satisfied,
after the pull back by an appropriate covering map
$\nbigr_{q,p}:\nbigb_p(R_1)\lrarr\nbigb_q(R)$,
Condition {\rm\ref{condition;21.8.21.2}} is satisfied,
according to Lemma {\rm\ref{lem;17.9.27.5}}.
\hfill\qed
 \end{rem}

We set 
$E_{\ell}:=\bigoplus_{\alpha}E_{\ell,\alpha}$.
We obtain the decomposition
\begin{equation}
\label{eq;17.9.27.6}
(E,\delbar_E)
=\bigoplus_{\ell\in\seisuu}
 (E_{\ell},\delbar_{E_{\ell}}).
\end{equation}

\subsection{Asymptotic orthogonality of
 the decomposition (\ref{eq;17.9.27.6})}

We set $h^{\circ}_{\ell}:=h_{|E_{\ell}}$.
We obtain the metric
$h^{\circ}:=\bigoplus h^{\circ}_{\ell}$
on $E=\bigoplus E_{\ell}$.
\index{metric $h^{\circ}$}
We set $I_q(w_q):=|w_q^q|\cdot\log|w_q^q|$.
\index{function $I_q(w_q)$}
We shall prove the following theorem
in \S\ref{subsection;17.9.28.10}--\ref{subsection;17.10.5.2}.

\begin{thm}
\label{thm;16.9.20.22}
Let $s$ be the unique automorphism of $E$
such that 
(i) $h^{\circ}=h\cdot s$,
(ii) $s$ is self-adjoint with respect to
 both $h$ and $h^{\circ}$.
For any $m\in\seisuu_{\geq 0}$,
there exist positive constants
$A_i(m)$ $(i=1,2)$ such that 
\[
 \Bigl|
 \nabla_{\kappa_1}\circ\cdots\circ
 \nabla_{\kappa_m}(s-\id_E)
 \Bigr|_h
\leq A_1(m)\exp\bigl(-A_2(m)\cdot I_q(w_q)\bigr)
\]
for any 
$(\kappa_1,\ldots,\kappa_{m})
\in\{t,x,y\}^{m}$.
\end{thm}

\begin{rem}
Theorem {\rm\ref{thm;16.9.20.22}}
and consequences below
were originally prepared 
for the study of the Nahm transform
of periodic monopoles
{\rm\cite{Mochizuki-periodic-Nahm}}.
\hfill\qed
\end{rem}

Let $\nabla^{\circ}$ and $\phi^{\circ}$
be the Chern connection 
and the Higgs field of 
$(E,\delbar_E,h^{\circ})$
as in \S\ref{subsection;16.9.20.20}.
\index{connection $\nabla^{\circ}$}
\index{Higgs field $\phi^{\circ}$}
Let $F(\nabla^{\circ})$
denote the curvature of $\nabla^{\circ}$.

\begin{cor}
\label{cor;17.10.5.5}
For any $m\in\seisuu_{\geq 0}$,
there exist positive constants $A_i(m)$ $(i=3,4)$ such that 
\[
 \Bigl|
 \nabla^{\circ}_{\kappa_1}\circ
 \cdots \circ
 \nabla^{\circ}_{\kappa_{m}}
 \bigl(
 \nabla-\nabla^{\circ}
 \bigr)
 \Bigr|_{h^{\circ}}
\leq
 A_{3}(m)\exp\bigl(-A_4(m)\cdot I_q(w_q)\bigr),
\]
\[
  \Bigl|
 \nabla^{\circ}_{\kappa_1}\circ
 \cdots \circ
 \nabla^{\circ}_{\kappa_{m}}
 (\phi-\phi^{\circ})
 \Bigr|_{h^{\circ}}
\leq
 A_{3}(m)\exp\bigl(-A_4(m)\cdot I_q(w_q)\bigr),
\]
\[
  \Bigl|
 \nabla^{\circ}_{\kappa_1}\circ
 \cdots \circ
 \nabla^{\circ}_{\kappa_{m}}
 \bigl(
 F(\nabla)-F(\nabla^{\circ})
 \bigr)
 \Bigr|_{h^{\circ}}
\leq
 A_{3}(m)\exp\bigl(-A_4(m)\cdot I_q(w_q)\bigr),
\]
for any 
$(\kappa_1,\ldots,\kappa_{m})
 \in \{t,x,y\}^{m}$.
\hfill\qed
\end{cor}

We obtain the following from Theorem \ref{thm;16.9.20.22}.
\begin{cor}
\label{cor;20.7.23.10}
 For any $m\in\seisuu_{\geq 0}$,
there exist positive constants $A_i(m)$ $(i=5,6)$ such that 
\[
 \Bigl|
 \nabla^{\circ}_{\kappa_1}\circ
 \cdots \circ
 \nabla^{\circ}_{\kappa_{m}}
 \bigl(F(\nabla^{\circ})
-\ast\nabla^{\circ}\phi^{\circ}\bigr)
 \Bigr|_{h^{\circ}}
\leq
 A_{5}(m)\exp\bigl(-A_6(m)\cdot I_q(w_q)\bigr)
\]
for any 
$(\kappa_1,\ldots,\kappa_{m})
 \in \{t,x,y\}^{m}$.
\hfill\qed
\end{cor}

\subsection{Eigen decomposition in the level $1$}
\label{subsection;17.10.5.50}

Let us look at the decomposition of mini-holomorphic bundles
\begin{equation}
 \label{eq;17.9.27.7}
(E_{\ell},\delbar_{E_{\ell}})
=\bigoplus_{\alpha\in\cnum^{\ast}}
 (E_{\ell,\alpha},\delbar_{E_{\ell,\alpha}}).
\end{equation}
For each $\alpha$,
there exists a locally free $\nbigo_{U_{w,q}^{\ast}}$-module
$V_{\ell,\alpha}$ with an endomorphism $f_{\ell,\alpha}$
and an isomorphism of mini-holomorphic bundles:
\begin{equation}
 \label{eq;17.10.6.300}
 E_{\ell,\alpha}
\simeq
 \LL^{0\ast}_{q}(\ell,\alpha)
\otimes
 \Psi_q^{\ast}(V_{\ell,\alpha},f_{\ell,\alpha}).
\end{equation}
We impose that the eigenvalues of $f_{\ell,\alpha|w_q}$
goes to $0$ as $|w_q|\to\infty$.
Let $h'_{\ell,\alpha}$
denote the restriction of
$h^{\circ}_{\ell}
 \otimes
 (h_{\LL,q,\ell,\alpha})^{-1}$
to 
$\Psi_q^{\ast}(V_{\ell,\alpha},f_{\ell,\alpha})$.
We have the Fourier expansion
of $h'_{\ell,\alpha}$ along $S^1_T$:
\[
 h'_{\ell,\alpha}
=\sum_{j\in\seisuu}
 \Psi_q^{-1}(h'_{V,\ell,\alpha;j})
 e^{2\pi\sqrt{-1}jt/T}.
\]
Here,
$h_{V,\ell,\alpha,j}'$
are sesqui-linear pairings
$V_{\ell,\alpha}
 \otimes
 \overline{V_{\ell,\alpha}}
\lrarr\cnum$.
We set $h_{V,\ell,\alpha}:=h_{V,\ell,\alpha;0}'$.

We set 
$h^{\shikaku}_{\ell}:=
 \bigoplus_{\alpha\in\cnum^{\ast}}
 h_{\LL,q,\ell,\alpha}
 \otimes
 \Psi_q^{-1}(h_{V,\ell,\alpha})$
on $E_{\ell}$.
\index{metric $h^{\shikaku}_{\ell}$}
We obtain the following as in the case of
doubly periodic instantons
\cite{Mochizuki-doubly-periodic}.
\begin{prop}
\label{prop;17.10.5.3}
Let $s_{\ell}$ be the unique endomorphism
of $E_{\ell}$ determined by the condition
$h^{\shikaku}_{\ell}=h^{\circ}_{\ell}s_{\ell}$.
For any $m\in\seisuu_{\geq 0}$,
there exist
$A_{i}(m)>0$ $(i=10,11)$
such that the following holds
for any 
$(\kappa_1,\ldots,\kappa_{m})
 \in \{t,x,y\}^{m}$:
\[
 \Bigl|
 \nabla^{\circ}_{\kappa_1}\circ
 \cdots
 \circ
 \nabla^{\circ}_{\kappa_m}
 \bigl(
 \id_{E_{\ell}}-s_{\ell}
 \bigr)
\Bigr|_{h^{\circ}_{\ell}}
\leq
 A_{10}(m)\cdot
 \exp\bigl(
 -A_{11}(m)|w_q^q|
 \bigr).
\]
\end{prop}
\pf
We have only to compare
$h_{\LL,q,\ell,1}^{-1}\otimes
 h_{\ell}^{\shikaku}$
and $h_{\LL,q,\ell,1}^{-1}\otimes h^{\circ}_{\ell}$.
Hence, it is enough to study the case $\ell=0$.
A similar problem was closely studied in
the case of the doubly periodic instantons
\cite[Theorem 5.11, Theorem 5.17]{Mochizuki-doubly-periodic}.
We shall explain an outline of the proof of
a generalization of \cite[Theorem 5.17]{Mochizuki-doubly-periodic}
in Proposition \ref{prop;17.10.23.31} below.
Here, we explain how the proof of 
Proposition \ref{prop;17.10.5.3}
is reduced to Proposition \ref{prop;17.10.23.31}.

We set $S^1_1=\real/\seisuu$.
Let $(\Etilde_0,\delbar_{\Etilde_0})$
be the holomorphic bundle
with the metric $\htilde_0$
on $S^1_1\times\nbigb^{0\ast}_q(R)$
induced by the mini-holomorphic bundle
$(E_0,\delbar_{E_0})$
with the metric $h^{\circ}_0$
as explained in \S\ref{subsection;13.11.29.2}.
Let $F(\htilde_0)$ denote the curvature 
of the Chern connection of 
$(\Etilde_0,\delbar_{\Etilde_0},\htilde_0)$.
By Corollary \ref{cor;20.7.23.10},
there exists $C>0$ such that
\[
 \Lambda F(\htilde_0)=
 O\bigl(\exp(-CI_q(w_q))\bigr). 
\]
By {\bf (B1)} in Condition \ref{condition;21.8.21.1}
and Corollary \ref{cor;20.7.23.10},
we obtain $F(\htilde_0)\to 0$ as $|w_q|\to\infty$.
Hence, $\Etilde_{0|S^1_1\times S^1_T\times\{w_q\}}$
are semistable of degree $0$
if $|w_q|$ is sufficiently large,
and we obtain the spectral curve
$\Sigma_{\Etilde_0}$ of $\Etilde_0$.
(See \cite[\S2.1]{Mochizuki-doubly-periodic}.)

Let us prove that the closure of 
$\Sigma_{\Etilde_0}$
in $U_{w,q}\times \nbigt^{\lor}$
is a complex analytic curve.
Let $S^1_1\times\nbigb^{0\ast}_q(R)\lrarr \nbigb^{0\ast}_q(R)$
be the projection.
Let $s$ denote the standard local coordinates of 
$S^1_1=\real/\seisuu$.
By considering $z=s+\sqrt{-1}t$,
we may regard
$S^1_1\times\nbigb^{0\ast}_q(R)$
as the quotient of
$\cnum_z\times U_{w,q}^{\ast}(R)$
by the lattice
\[
 \Gamma=\bigl\{
 n_1+T\sqrt{-1}n_2\,\big|\,
 n_1,n_2\in\seisuu
 \bigr\}
\subset\cnum_z.
\]
The dual lattice $\Gamma^{\lor}$ is given as
$\Gamma^{\lor}=
 \bigl\{
 T^{-1}\pi m_1+\pi \sqrt{-1}m_2\,\big|\,
 m_1,m_2\in\seisuu
 \bigr\}$.
Let $\nbigt$ denote $S^1_1\times S^1_T=\cnum_z/\Gamma$,
and $\nbigt^{\lor}$ denote the dual
$\cnum_{\zeta}/\Gamma^{\lor}$.
There exists the isomorphism
$\cnum\big/2\pi\sqrt{-1}\seisuu
\simeq
 \cnum^{\ast}$ 
induced by the exponential.
We obtain the morphism
\[
 \cnum\big/2\pi\sqrt{-1}\seisuu
\lrarr
 \nbigt^{\lor}=\cnum\big/\Gamma^{\lor}
\]
induced by
$\gamma\longmapsto
 (2T)^{-1}\gamma$.
It induces the following holomorphic map:
\[
 \nbigj:
 \cnum^{\ast}\times U_{w,q}^{\ast}
\lrarr
 \nbigt^{\lor}\times U_{w,q}.
\]
We can easily observe
$\Sigma_{\Etilde_0}
=\nbigj(\Sp(E_0))$.
Hence, the closure of 
$\Sigma_{\Etilde_0}$
in $\nbigt^{\lor}\times U_{w,q}$
is a complex analytic curve.
Then, the claim of Proposition \ref{prop;17.10.5.3}
follows from Proposition \ref{prop;17.10.23.31}.
\hfill\qed

\vspace{.1in}

Let $\nabla^{\shikaku}_{\ell}$
and $\phi^{\shikaku}_{\ell}$
denote the Chern connection
and the Higgs field of
$(E_{\ell},\delbar_{E_{\ell}},h^{\shikaku}_{\ell})$.
We obtain a Hermitian metric
$h^{\shikaku}=\bigoplus h^{\shikaku}_{\ell}$,
a unitary connection
$\nabla^{\shikaku}=\bigoplus \nabla^{\shikaku}_{\ell}$
and
an anti-Hermitian endomorphism
$\phi^{\shikaku}=\bigoplus \phi^{\shikaku}_{\ell}$
on $E$.

\begin{cor}
\label{cor;17.10.5.4}
For any $m\in\seisuu_{\geq 0}$,
there exist positive constants
$A_{i}(m)$ $(i=12,13)$
such that 
for any 
$(\kappa_1,\ldots,\kappa_{m})
 \in \{t,x,y\}^{m}$:
\[
 \Bigl|
 \nabla^{\circ}_{\kappa_1}\circ
 \cdots
 \circ
 \nabla^{\circ}_{\kappa_m}
 \bigl(
 \nabla^{\circ}
-\nabla^{\shikaku}
 \bigr)
\Bigr|_{h^{\circ}}
\leq
 A_{12}(m)\cdot
 \exp\bigl(
 -A_{13}(m)|w_q^q|
 \bigr).
\]
\[
 \Bigl|
 \nabla^{\circ}_{\kappa_1}\circ
 \cdots
 \circ
 \nabla^{\circ}_{\kappa_m}
 \bigl(
 \phi^{\circ}
-\phi^{\shikaku}
 \bigr)
\Bigr|_{h^{\circ}}
\leq
 A_{12}(m)\cdot
 \exp\bigl(
 -A_{13}(m)|w_q^q|
 \bigr).
\]
\[
 \Bigl|
 \nabla^{\circ}_{\kappa_1}\circ
 \cdots
 \circ
 \nabla^{\circ}_{\kappa_m}
 \bigl(
 F(\nabla^{\circ})
-F(\nabla^{\shikaku})
 \bigr)
\Bigr|_{h^{\circ}}
\leq
 A_{12}(m)\cdot
 \exp\bigl(
 -A_{13}(m)|w_q^q|
 \bigr).
\]
\hfill\qed
\end{cor}
We obtain the following.
\begin{cor}
\label{cor;17.9.27.20}
For any $m\in\seisuu_{\geq 0}$,
we have positive constants
$A_{i}(m)$ $(i=14,15)$
such that the following holds
for any 
$(\kappa_1,\ldots,\kappa_{m})
 \in \{t,x,y\}^{m}$:
\[
 \Bigl|
 \nabla^{\shikaku}_{\kappa_1}\circ
 \cdots
 \circ
 \nabla^{\shikaku}_{\kappa_m}
 \bigl(
 F(\nabla^{\shikaku})
-\ast \nabla^{\shikaku}\phi^{\shikaku}
 \bigr)
\Bigr|_{h^{\circ}}
\leq
 A_{14}(m)\cdot
 \exp\bigl(
 -A_{15}(m)|w_q^q|
 \bigr).
\]
\hfill\qed
\end{cor}

\subsection{Asymptotic harmonic bundles}
\label{subsection;17.10.6.301}

We have the Higgs field
\[
 \theta_{\ell,\alpha}
=f_{\ell,\alpha}\cdot(qw_q^{q-1}dw_q)
\]
of $(V_{\ell,\alpha},\delbar_{V_{\ell,\alpha}})$.
Let $F(h_{V,\ell,\alpha})$
denote the curvature of the Chern connection
$\nabla^{\ell,\alpha}$
of
$(V_{\ell,\alpha},\delbar_{V_{\ell,\alpha}},h_{V,\ell,\alpha})$.
Let 
$\bigl(
 \theta_{\ell,\alpha}
\bigr)^{\dagger}$
 denote the adjoint of $\theta_{\ell,\alpha}$
 with respect to $h_{V,\ell,\alpha}$.
The following is a direct consequence of
Corollary \ref{cor;17.9.27.20}.

\begin{cor}
\label{cor;20.7.23.11}
 For any $m\in\seisuu_{\geq 1}$,
we have $A_{i}(m)$ $(i=20,21)$
such that
\begin{equation}
\label{eq;17.10.6.302}
\left|
 \nabla^{\ell,\alpha}_{\kappa_1}\circ
 \cdots
 \circ
 \nabla^{\ell,\alpha}_{\kappa_m}
 \Bigl(
  F(h_{V,\ell,\alpha})+
 \bigl[
 \theta_{\ell,\alpha},
 \bigl(
 \theta_{\ell,\alpha}
\bigr)^{\dagger}
 \bigr]
 \Bigr)
\right|_{h_{V,\ell,\alpha}}
\leq
 A_{20}(m)\cdot
 \exp\bigl(-A_{21}(m)|w_q^q|\bigr)
\end{equation}
for any 
$(\kappa_1,\ldots,\kappa_{m})
 \in \{x,y\}^{m}$.
\hfill\qed
\end{cor}

\subsection{Curvature}

We have the expression
\[
 F(\nabla)=
 F(\nabla)_{w,\wbar}dw\,d\wbar
+F(\nabla)_{w,t}dw\,dt
+F(\nabla)_{\wbar,t}d\wbar\,dt.
\]

\begin{cor}
\label{cor;17.10.6.1}
The following estimates hold:
\[
 \bigl|
 F(\nabla)_{w\wbar}
 \bigr|_h
=O\Bigl(
 |w_q^q|^{-2}(\log|w_q|)^{-2}
 \Bigr),
\]
\[
 \bigl|
 F(\nabla)_{wt}
 \bigr|_h
=O(|w_q^q|^{-1}),
\quad
 \bigl|
 F(\nabla)_{\wbar t}
 \bigr|_h
=O(|w_q^q|^{-1}).
\]
\end{cor}
\pf
Let $F(h'_{\ell,\alpha})$
be the curvature of the Chern connection of
$E_{\ell,\alpha}
 \otimes
 \LL^{0\ast}_{q}(\ell,\alpha)^{-1}$
with the metric $h'_{\ell,\alpha}$.
We have the expression
$F(h'_{\ell,\alpha})
=F(h'_{\ell,\alpha})_{w\wbar}dw\,d\wbar
+F(h'_{\ell,\alpha})_{wt}dw\,dt
+F(h'_{\ell,\alpha})_{\wbar t}d\wbar\,dt$.
As in the case of doubly periodic instantons
\cite[Theorem 5.14]{Mochizuki-doubly-periodic},
we obtain the following estimates
for some $\delta>0$
from
Corollary \ref{cor;17.10.5.4},
Corollary \ref{cor;20.7.23.11},
and the estimates for asymptotic harmonic bundles
(see \cite[Theorem 5.14, \S5.5]{Mochizuki-doubly-periodic}):
\[
 \Bigl|
 F(h'_{\ell,\alpha})_{w\wbar}
 \Bigr|_{h'_{\ell,\alpha}}
=O\bigl(
 |w|^{-2}(\log|w|)^{-2}
 \bigr),
\]
\[
 \Bigl|
 F(h'_{\ell,\alpha})_{wt}
 \Bigr|_{h'_{\ell,\alpha}}
=O\bigl(
 |w|^{-1-\delta}
 \bigr),
\quad
  \Bigl|
 F(h'_{\ell,\alpha})_{\wbar t}
 \Bigr|_{h'_{\ell,\alpha}}
=O\bigl(
 |w|^{-1-\delta}
 \bigr).
\]
Let $F(h_{\LL,q,\ell,\alpha})$
denote the curvature of
$\LL^{0\ast}_q(\ell,\alpha)$
with the metric $h_{\LL,q,\ell,\alpha}$.
By an easy computation,
we obtain the following:
\[
 F(h_{\LL,q,\ell,\alpha})_{w\,\wbar}=0,
\quad
 \bigl|
  F(h_{\LL,q,\ell,\alpha})_{w\,t}
 \bigr|
=O(|w|^{-1}),
\quad
 \bigl|
  F(h_{\LL,q,\ell,\alpha})_{\wbar\,t}
 \bigr|
=O(|w|^{-1}).
\]
Then, we obtain the claim of the lemma
from Corollary \ref{cor;17.10.5.5}.
\hfill\qed

\subsection{Another equivalent decay condition}

Let $(E,h,\nabla,\phi)$ be a monopole
on $\nbigb_{q}^{0\ast}(R)$ for some $R>0$.

\begin{prop}
\label{prop;20.7.29.20}
Condition {\rm\ref{condition;21.8.21.1}} is satisfied
for $(E,h,\nabla,\phi)$
if and only if the following holds.
\begin{description}
 \item[(GCK)]
$|F(h)|_h\to 0$ $(|w_q|\to\infty)$
and $|\phi|_h=O\bigl(\log|w_q|\bigr)$.
\end{description}
\index{GCK-condition}
\end{prop}
\pf
Suppose that {\bf (GCK)} is satisfied.
The condition {\bf (B1)} follows from the Bogomolny equation
and $|F(h)|_h\to 0$ $(|w_q|\to \infty)$.
Let $(E,\delbar_E)$ denote the mini-holomorphic bundle
on $\nbigb^{0\ast}_q(R)$
underlying the monopole.
Let $M(t,w_q)$ denote the holomorphic family of monodromy
of $E_{|\{t\}\times U_{w,q}^{\ast}(R)}$.
Because $|\phi|_h=O\bigl(\log|w_q|\bigr)$,
we obtain the following estimate for some $N>0$:
\begin{equation}
\label{eq;17.9.13.30}
 |M(t,w_q)|_h+|M(t,w_q)^{-1}|_h=O(|w_q|^{N}).
\end{equation}
Let 
$P(w_q,x):=\det(x\id-M(t,w_q))=\sum a_j(w_q)x^j$
denote the characteristic polynomials
of the monodromy automorphisms.
Then, $\Sp(E,\delbar_E)$ is the zero set of $P(w_q,x)$.
By (\ref{eq;17.9.13.30}),
$a_j(w_q)$ are meromorphic at $w_q=\infty$.
Hence, the closure of $\Sp(E,\delbar_E)$ 
in $U_{w,q}(R)\times \proj^1$
is also a complex analytic subset,
i.e., {\bf (B2)} is satisfied.

Suppose Condition \ref{condition;21.8.21.1} is satisfied.
We obtain the estimate for the curvature by the Bogomolny equation.
Let us study the growth order of $|\phi|_h$.
By taking the pull back via an appropriate covering
$\nbigr_{q,p}:\nbigb_{p}^{0\ast}(R_1)\lrarr\nbigb_{q}^{0\ast}(R)$,
we may assume that Condition \ref{condition;21.8.21.2}
is also satisfied.
The norm of the Higgs field of the monopole
$\bigl(
 \LL^{0\ast}_q(\ell,\alpha),h_{\LL,q,\ell,\alpha}
\bigr)$ is
$O\bigl(\log|w|\bigr)$
(see \S\ref{subsection;17.9.26.1}).
There exists a constant $C_0>0$
such that any eigenvalue $a$ of $f_{\ell,\alpha|Q}$ 
$(Q\in U_{w,q}^{\ast})$
satisfies $|a|<C_0$.
Then, we obtain that
$\bigl|f_{\ell,\alpha}\bigr|_{h_{V,\ell,\alpha}}$
are bounded
according to Simpson's main estimate
for asymptotic harmonic bundles.
(See Proposition \ref{prop;12.8.7.3}.)
Hence, 
by the formula (\ref{eq;17.10.6.21}),
the Higgs field of
$\bigl(
 \Psi_q^{\ast}(V_{\ell,\alpha},f_{\ell,\alpha}),\Psi_q^{-1}(h_{V,\ell,\alpha})
 \bigr)$
is bounded.
Then, we obtain 
$|\phi|_h=O(\log|w|)$
from Corollary \ref{cor;17.10.5.5}
and Corollary \ref{cor;17.10.5.4},
i.e., {\bf (GCK)} is satisfied.
\hfill\qed

\section{Connections and orthogonal decompositions}
\label{subsection;17.9.28.10}

\subsection{Statement}
\label{subsection;21.8.21.4}

We consider a monopole
$(E,h,\nabla,\phi)$ on $\nbigb^{0\ast}_q(R_0)$
satisfying the following condition
for some $R_0>0$ and $C_0>0$.
\begin{condition}
\label{condition;21.8.21.3}\mbox{{}}
\begin{itemize}
\item
$\bigl|\nabla(\phi)_{|(t,w_q)}\bigr|_h\to 0$
when $|w_q|\to \infty$.
\item
 There exist a finite subset
 $\Lambda\subset\seisuu$
 and an orthogonal decomposition 
 $(E,\phi)=
 \bigoplus_{\ell\in\Lambda} 
 (E^{\bullet}_{\ell},\phi^{\bullet}_{\ell})$
 such that
 any eigenvalue $\alpha$ of
 $\phi^{\bullet}_{\ell|(t,w_q)}$
 satisfies 
\[
 \Bigl|
 \alpha-\sqrt{-1}\ell T^{-1}\log|w_q|
 \Bigr|\leq C_0.
\]
\item
For any $\ell_1,\ell_2\in\Lambda$ $(\ell_1\neq\ell_2)$,
we obtain
$|\ell_1-\ell_2|\cdot T^{-1}\log R_0>10C_0$.
\hfill\qed
\end{itemize}
\end{condition}
 
We obtain the decomposition
$\nabla=\nabla^{\bullet}+\rho$,
where $\nabla^{\bullet}$ is the unitary connection
preserving the decomposition 
$E=\bigoplus E^{\bullet}_{\ell}$,
and $\rho$ is a section of
$\bigoplus_{\ell_1\neq \ell_2}
 \Hom(E^{\bullet}_{\ell_1},E^{\bullet}_{\ell_2})
 \otimes\Omega^1$.
The inner products of $\rho$ and $\del_{\kappa}$ $(\kappa=t,x,y)$
are denoted by $\rho_{\kappa}$.
We set
$\End(E)^{\bullet}:=\bigoplus \End(E^{\bullet}_{\ell})$
and
$\End(E)^{\top}=\bigoplus_{\ell_1\neq\ell_2}
\Hom(E^{\bullet}_{\ell_1},E^{\bullet}_{\ell_2})$.
\index{bundle $\End(E)^{\bullet}$}
\index{bundle $\End(E)^{\top}$}
For any section $s$ of
$\End(E)\otimes\Omega^p$,
let $s=s^{\bullet}+s^{\top}$ be the unique decomposition
such that $s^{\bullet}$ and $s^{\top}$
are sections of
$\End(E)^{\bullet}\otimes\Omega^p$
and
$\End(E)^{\top}\otimes\Omega^p$,
respectively.
Note that
$(\nabla\phi)^{\top}=[\rho,\phi]$.
\index{section $s^{\bullet}$}
\index{section $s^{\top}$}

We shall prove the following proposition
in \S\ref{subsection;16.9.10.61}--\S\ref{subsection;16.9.10.62}.

\begin{prop}
\label{prop;16.9.10.60}
There exist positive constants
$R_1$ and $C_i$ $(i=1,2)$
depending only on $\rank E$, $R_0$, $C_0$ and $\Lambda$
such that 
$\bigl|
 \rho
 \bigr|_h
\leq
 C_1\exp\bigl(-C_2I_q(w_q) \bigr)$
on $\nbigb^{0\ast}_q(R_1)$.
For any positive integer $k$,
there exist positive constants
$C_i(k)$ $(i=1,2)$
such that 
\[
 \bigl|
 \nabla^{\bullet}_{\kappa_1}\circ\cdots
 \circ\nabla^{\bullet}_{\kappa_{k}}
 \rho
 \bigr|_h
\leq 
C_1(k)\exp\bigl(-C_2(k)I_q(w_q) \bigr)
\]
on $\nbigb^{0\ast}_q(R_1)$
for any 
$(\kappa_1,\ldots,\kappa_{k})
 \in\{t,x,y\}^{k}$.
\end{prop}

We introduce a notation.
Let $f$ and $g$ be functions on an open subset
of $\nbigb^{0\ast}_q(R)$ for some $R>R_0$.
We say $f=O(g)$
if there exists a constant $B$
depending only on $\rank E$, $C_0$ and $\Lambda$
such that $|f|\leq B|g|$ on
$\nbigb^{0\ast}_q(R)$.

\subsection{Preliminary}
\label{subsection;16.9.10.61}

For any $\epsilon>0$,
there exists $R_{10}(\epsilon)$ such that 
$|\nabla\phi|_h<\epsilon$
on $\nbigb^{0\ast}_q(R_{10}(\epsilon))$.
If $\epsilon$ is sufficiently small,
as remarked in Lemma \ref{lem;17.10.5.10},
we obtain
\[
 \bigl|
 \nabla_{\kappa_1}\circ\cdots
 \circ\nabla_{\kappa_{k}}\phi
 \bigr|_h\leq B(k)\epsilon
\]
on $\nbigb^{0\ast}_q(2R_{10}(\epsilon))$
for any 
$(\kappa_1,\ldots,\kappa_{k})
 \in\{t,x,y\}^{k}$.
Here, $B(k)$ are positive constants
which are independent of $\epsilon$.

\begin{lem}
\label{lem;16.9.10.10}
We obtain
$|\rho_{\kappa}|_h
=O\bigl(
 \bigl|(\nabla_{\kappa}\phi)^{\top}\bigr|
 \bigr)$ 
for $\kappa=t,x,y$.
We also obtain
the following for any
$\kappa_1,\kappa_2\in\{t,x,y\}$:
\[
 \bigl|\nabla^{\bullet}_{\kappa_1}\rho_{\kappa_2}\bigr|_h
=O\Bigl(
 \epsilon
 \bigl|(\nabla_{\kappa_2}\phi)^{\top}\bigr|_h
+\bigl|
 \nabla_{\kappa_1}^{\bullet}
 (\nabla_{\kappa_2}\phi)^{\top}
 \bigr|
 \Bigr).
\]
\end{lem}
\pf
The first claim follows from
$[\phi,\rho_{\kappa}]
=-\nabla_{\kappa}(\phi)^{\top}$.
We have the following equality:
\[
\nabla_{\kappa_1}^{\bullet}\bigl(
 (\nabla_{\kappa_2}\phi)^{\top}
 \bigr)
=
 \nabla^{\bullet}_{\kappa_1}
 \bigl(
 \bigl[\phi,\rho_{\kappa_2}\bigr]
 \bigr)
=\bigl[
 \nabla_{\kappa_1}^{\bullet}\phi,
 \rho_{\kappa_2}
 \bigr]
+\bigl[
 \phi,\nabla_{\kappa_1}^{\bullet}\rho_{\kappa_2}
 \bigr].
\]
Then, the second claim follows.
\hfill\qed

\vspace{.1in}
We recall some equalities.

\begin{lem}
\label{lem;16.9.10.62}
For $(a,b,c)=(t,x,y),(x,y,t),(y,t,x)$,
we obtain
\[
 (\nabla_x^2+\nabla_y^2+\nabla_t^2)(\nabla_a\phi)
=4\bigl[\nabla_b\phi,\nabla_c\phi\bigr]
-\bigl[\phi,[\phi,\nabla_a\phi]\bigr].
\]
\end{lem}
\pf
It is enough to consider the case
$(a,b,c)=(t,x,y)$.
We obtain
\begin{multline}
 \nabla_x^2(\nabla_t\phi)
=\nabla_x(\nabla_t\nabla_x\phi+[F_{x,t},\phi])
=\nabla_t\nabla_x^2\phi
+[F_{x,t},\nabla_x\phi]
+\nabla_x[F_{x,t},\phi]
\\
=\nabla_t\nabla_x^2\phi
+2[F_{x,t},\nabla_x\phi]
+[\nabla_x(F_{x,t}),\phi].
\end{multline}
Similarly, we also obtain
$\nabla_y^2\nabla_t\phi
=\nabla_t\nabla_y^2\phi
+2\bigl[
 F_{y,t},\nabla_y\phi
 \bigr]
+\bigl[\nabla_y(F_{y,t}),\phi\bigr]$.
We have the following equality,
which follows from the Bogomolny equation:
\begin{equation}
\label{eq;21.8.22.10}
(\nabla_x^2+\nabla_y^2+\nabla_t^2)\phi=0.
\end{equation}
By using
$\ast F=\nabla\phi$,
we obtain the following:
\begin{multline}
(\nabla_x^2+\nabla_y^2+\nabla_t^2)
 (\nabla_t\phi)
=-2\bigl[\nabla_y\phi,\nabla_x\phi\bigr]
+2\bigl[\nabla_x\phi,\nabla_y\phi\bigr]
+\bigl[-\nabla_x\nabla_y\phi+\nabla_y\nabla_x\phi,\phi\bigr]
\\
=4\bigl[\nabla_x\phi,\nabla_y\phi\bigr]
-\bigl[
 [F_{x,y},\phi],\phi
 \bigr]
=4\bigl[
 \nabla_x\phi,\nabla_y\phi
 \bigr]
+\bigl[
 [\phi,\nabla_t\phi],\phi
 \bigr]
 \\
=4\bigl[
 \nabla_x\phi,\nabla_y\phi
 \bigr]
-\bigl[
 \phi,[\phi,\nabla_t\phi]
 \bigr].
\end{multline}
Thus, we obtain Lemma \ref{lem;16.9.10.62}.
\hfill\qed

\subsection{Step 1}

We obtain the following equalities
for any $\kappa_1,\kappa_2\in\{t,x,y\}$:
\begin{multline}
 \label{eq;16.9.10.1}
 \del_{\kappa_1}^2
 \bigl|
 (\nabla_{\kappa_2}\phi)^{\top}
 \bigr|_h^2
=2\bigl|
 \nabla_{\kappa_1}
 (\nabla_{\kappa_2}\phi)^{\top}
 \bigr|_h^2
 \\
+2\Re h\bigl(
 \nabla_{\kappa_1}^2\nabla_{\kappa_2}\phi,
 (\nabla_{\kappa_2}\phi)^{\top}
 \bigr)
-2\Re h\bigl(
 \nabla_{\kappa_1}^2
 (\nabla_{\kappa_2}\phi)^{\bullet},
 (\nabla_{\kappa_2}\phi)^{\top}
 \bigr).
\end{multline}

Let us look at the last term
in the right hand side of (\ref{eq;16.9.10.1}).
\begin{lem}
The following holds on $\nbigb^{0\ast}_q(R_{10}(\epsilon))$:
\[
 h\bigl(
 \nabla_{\kappa_1}^2
 (\nabla_{\kappa_2}\phi)^{\bullet},
 (\nabla_{\kappa_2}\phi)^{\top}
 \bigr)
=O\Bigl(
 \epsilon\cdot
 \Bigl(
 \bigl|
 (\nabla_{\kappa_1}\phi)^{\top}
 \bigr|_h
+ \bigl|
 \nabla^{\bullet}_{\kappa_1}
 (\nabla_{\kappa_1}\phi)^{\top}
 \bigr|_h
 \Bigr)
\cdot 
 \bigl|
 (\nabla_{\kappa_2}\phi)^{\top}
 \bigr|_h
 \Bigr).
\]
\end{lem}
\pf
We have the following equality:
\begin{multline}
 \Bigl(
 \nabla_{\kappa_1}^2
 (\nabla_{\kappa_2}\phi)^{\bullet}
 \Bigr)^{\top}
= \\
 \Bigl(
 \bigl[
 \nabla^{\bullet}_{\kappa_1}\rho_{\kappa_1},
 (\nabla_{\kappa_2}\phi)^{\bullet}
 \bigr]
+2\bigl[ \rho_{\kappa_1},
 \nabla^{\bullet}_{\kappa_1}
 (\nabla_{\kappa_2}\phi)^{\bullet}
 \bigr]
+\Bigl[
 \rho_{\kappa_1},
 \bigl[\rho_{\kappa_1},
 (\nabla_{\kappa_2}\phi)^{\bullet}\bigr]
 \Bigr]
 \Bigr)^{\top}.
\end{multline}
Note that
$|(\nabla_{\kappa_2}\phi)^{\bullet}|_h=O(\epsilon)$.
We also have
$\nabla_{\kappa_1}^{\bullet}
(\nabla_{\kappa_2}\phi)^{\bullet}
=(\nabla_{\kappa_1}\nabla_{\kappa_2}\phi)^{\bullet}
-[\rho_{\kappa_1},[\rho_{\kappa_2},\phi]]^{\bullet}$.
Then, the claim follows from
Lemma \ref{lem;16.9.10.10}.
\hfill\qed

\vspace{.1in}
Let us look at the sum of the second terms 
in the right hand side of (\ref{eq;16.9.10.1})
for $\kappa_1=t,x,y$.

\begin{lem}
We have the following estimate on
$\nbigb^{0\ast}_q(R_{10}(\epsilon))$:
\begin{equation}
\sum_{\kappa_1=t,x,y}
 h\bigl(
 \nabla_{\kappa_1}^2
 (\nabla_{\kappa_2}\phi),
 (\nabla_{\kappa_2}\phi)^{\top}
 \bigr)
=
 \Bigl|
 \bigl[\phi,(\nabla_{\kappa_2}\phi)^{\top}\bigr]
 \Bigr|^2_h
+
O\Bigl(
 \epsilon\cdot
 \bigl|
 (\nabla\phi)^{\top}
 \bigr|_h
\cdot
 \bigl|(\nabla_{\kappa_2}\phi)^{\top}\bigr|
 \Bigr).
\end{equation}
\end{lem}
\pf
We explain the case $\kappa_2=t$.
By Lemma \ref{lem;16.9.10.62},
we have the following equality:
\begin{multline}
\Bigl(
 (\nabla_x^2+\nabla_y^2+\nabla_t^2)
 \nabla_{t}\phi
\Bigr)^{\top}
=\Bigl(
 4[\nabla_x\phi,\nabla_y\phi]
-\bigl[
 \phi,[\phi,\nabla_t\phi]
 \bigr]
\Bigr)^{\top}
 \\
=
 O\Bigl(
 \epsilon\cdot
 \Bigl(
 \bigl|(\nabla_x\phi)^{\top}\bigr|_h
+\bigl|(\nabla_y\phi)^{\top}
 \bigr|_h
 \Bigr)
 \Bigr)
-\bigl[
 \phi,[\phi,\nabla_t\phi]
 \bigr]^{\top}.
\end{multline}
We have the following equality:
\begin{multline}
 -h\bigl(
 [\phi,[\phi,\nabla_t\phi]]^{\top},
 (\nabla_t\phi)^{\top}
 \bigr)
=-h\bigl(
 [\phi,[\phi,\nabla_t\phi]],
 (\nabla_t\phi)^{\top}
 \bigr)
 \\
=h\bigl(
 [\phi,\nabla_t\phi],\,
 [\phi,(\nabla_t\phi)^{\top}]
 \bigr)
=\bigl|
 [\phi,(\nabla_t\phi)^{\top}]
 \bigr|_h^2.
\end{multline}
Then, we obtain the claim of the lemma.
\hfill\qed

\begin{lem}
There exist $R_{11}>0$ and $C_{11}>0$
such that the following inequality holds
on $\nbigb^{0\ast}_q(R_{11})$:
\[
 -(\del_x^2+\del_y^2+\del_t^2)
 \bigl|(\nabla\phi)^{\top}\bigr|_h^2
\leq
 -C_{11}
 \bigl|
 (\nabla\phi)^{\top}
 \bigr|_h^2
 \cdot (\log|w|)^2.
\]
\end{lem}
\pf By the previous lemmas,
we obtain the following estimate
on $\nbigb^{0\ast}_q(R_{10}(\epsilon))$:
\begin{multline}
 -(\del_x^2+\del_y^2+\del_t^2)
 \bigl|
 (\nabla_{\kappa_2}\phi)^{\top}
 \bigr|_h^2
=-2\bigl|
 \nabla\bigl(
 (\nabla_{\kappa_2}\phi)^{\top}\bigr)
 \bigr|_h^2
-2\bigl|
 \bigl[\phi,(\nabla_{\kappa_2}\phi)^{\top}\bigr]
 \bigr|_h^2
 \\
+O\Bigl(
 \epsilon\cdot
 \bigl|(\nabla\phi)^{\top}\bigr|_h
\cdot
 \bigl|(\nabla_{\kappa_2}\phi)^{\top}\bigr|_h
 \Bigr)
 \\
+O\Bigl(
 \epsilon\cdot
 \Bigl(
 \bigl|\nabla^{\bullet}_{t}(\nabla_t\phi)^{\top}\bigr|_h
+\bigl|\nabla^{\bullet}_{x}(\nabla_x\phi)^{\top}\bigr|_h
+\bigl|\nabla^{\bullet}_{y}(\nabla_y\phi)^{\top}\bigr|_h
 \Bigr)
\cdot \bigl|(\nabla_{\kappa_2}\phi)^{\top}\bigr|_h
 \Bigr).
\end{multline}
Because $|\rho|_h=O(\epsilon)$ 
on $\nbigb^{0\ast}_q(R_{10}(\epsilon))$,
it is equivalent to the following estimate:
\begin{multline}
 -(\del_x^2+\del_y^2+\del_t^2)
 \bigl|
 (\nabla_{\kappa_2}\phi)^{\top}
 \bigr|_h^2
=
-2\bigl|
 \nabla\bigl(
 (\nabla_{\kappa_2}\phi)^{\top}\bigr)
 \bigr|_h^2
-2\bigl|
 \bigl[\phi,(\nabla_{\kappa_2}\phi)^{\top}\bigr]
 \bigr|_h^2
 \\
+O\Bigl(
 \epsilon\cdot
 \bigl|(\nabla\phi)^{\top}\bigr|_h
\cdot
 \bigl|(\nabla_{\kappa_2}\phi)^{\top}\bigr|_h
 \Bigr)
 \\
+O\Bigl(
 \epsilon\cdot
 \Bigl(
 \bigl|\nabla_{t}(\nabla_t\phi)^{\top}\bigr|_h
+\bigl|\nabla_{x}(\nabla_x\phi)^{\top}\bigr|_h
+\bigl|\nabla_{y}(\nabla_y\phi)^{\top}\bigr|_h
 \Bigr)
\cdot \bigl|(\nabla_{\kappa_2}\phi)^{\top}\bigr|_h
 \Bigr).
\end{multline}
By taking the sum for $\kappa_2=t,x,y$,
we obtain the following estimate
on $\nbigb^{0\ast}_q(R_{10}(\epsilon))$:
\begin{multline}
-(\del_x^2+\del_y^2+\del_t^2)
 \bigl|\bigl(
 \nabla\phi
 \bigr)^{\top}\bigr|_h^2
=
-\sum_{\kappa_2=t,x,y}
 2\bigl|
 \nabla\bigl(
 (\nabla_{\kappa_2}\phi)^{\top}\bigr)
 \bigr|_h^2
-2\bigl|
 \bigl[\phi,(\nabla\phi)^{\top}\bigr]
 \bigr|_h^2
\\
+O\Bigl(
 \epsilon\cdot
 \bigl|(\nabla\phi)^{\top}\bigr|^2_h
 \Bigr)
 \\
+O\Bigl(
 \epsilon\cdot
 \Bigl(
 \bigl|\nabla_{t}(\nabla_t\phi)^{\top}\bigr|_h
+\bigl|\nabla_{x}(\nabla_x\phi)^{\top}\bigr|_h
+\bigl|\nabla_{y}(\nabla_y\phi)^{\top}\bigr|_h
 \Bigr)
\cdot \bigl|(\nabla\phi)^{\top}\bigr|_h
 \Bigr).
\end{multline}
Note that there exists a constant $C_{10}>0$
such that 
$\bigl|
 \bigl[
 \phi,(\nabla\phi)^{\top}
 \bigr]
 \bigr|_h
\geq 
 C_{10}(\log|w|)\cdot
 \bigl|(\nabla\phi)^{\top}\bigr|_h$
on $\nbigb^{0\ast}_q(R_0)$.
Then, we obtain the claim of the lemma.
\hfill\qed

\vspace{.1in}

Set $U^{\ast}_{w,q}(R):=\{|w_q|>R\}$.
For any function $f$ on $\nbigb^{0\ast}_q(R)$,
let $\int_{S^1_T}f$
denote the function on $U^{\ast}_{w,q}(R)$
given by
$w_q\longmapsto \int_{S^1_T\times\{w_q\}}f\,dt$.
By the previous lemma,
we obtain the following inequality
on $U^{\ast}_{w,q}(R_{11})$:
\begin{equation}
\label{eq;16.9.10.20}
-(\del_x^2+\del_y^2)
 \int_{S^1_T}
 \bigl|(\nabla\phi)^{\bot}\bigr|_h^2
\leq
 -
 C_{11}
 \bigl(\log|w|
 \bigr)^2\cdot
 \int_{S^1_T}
 \bigl|(\nabla\phi)^{\bot}\bigr|_h^2.
\end{equation}

\begin{lem}
There exist $R_{12}>R_{11}$ and $C_{12}>0$
such that the following holds
on $U^{\ast}_{w,q}(R_{12})$:
\[
 \int_{S^1_T}
 \bigl|(\nabla\phi)^{\bot}\bigr|_h^2
=O\Bigl(
 \exp
 \bigl(
 -C_{12}I_q(w_q)
 \bigr) 
\Bigr).
\]
\end{lem}
\pf
Let $w=w_q^q$.
We have the following equality,
which can be checked by a direct computation:
\begin{multline}
\frac{\del^2}{\del w\del\wbar}
\Bigl(
 \exp\bigl(-C|w|\log|w|\bigr)
\Bigr)
=
\exp\bigl(
 -C|w|\log|w|
 \bigr)
 \cdot
 \frac{C^2}{4}(\log|w|+1)^2
\\
-\exp\bigl(
 -C|w|\log|w|
 \bigr)
 \cdot
 \frac{C}{4|w|}
 (\log|w|+2)
\end{multline}
We may assume $\log|w|\geq 1$ on 
$U^{\ast}_{w,q}(R_{11})$.
Hence, we have the following inequality:
\begin{equation}
\label{eq;16.9.10.30}
 -\frac{\del}{\del w}
 \frac{\del}{\del\wbar}
 \exp\bigl(-C|w|\log|w|\bigr)
\geq
 -C^2\exp\bigl(-C|w|\log|w|\bigr)
 \bigl(\log|w|\bigr)^2.
\end{equation}
Then, by the standard argument for Ahlfors lemma \cite{a}
using the inequalities 
(\ref{eq;16.9.10.20})
and (\ref{eq;16.9.10.30}),
we obtain the desired estimate.
(See \cite{Simpson90} and \cite{Mochizuki-wild}.
See also
the proof of Lemma \ref{lem;17.10.7.300}.)
\hfill\qed

\vspace{.1in}
As a corollary,
we obtain the following estimate on $\nbigb^{0\ast}_q(R_{12})$:
\begin{equation}
\label{eq;16.9.10.50}
 \int_{S^1_T}
 \bigl|
 \rho_{\kappa}
 \bigr|_h^2
=O\Bigl(
 \exp\bigl(-C_{12}I_q(w_q)\bigr)
 \Bigr).
\end{equation}

\subsection{Step 2}
\label{subsection;16.9.10.62}

Let $F_{\phi}$  denote the endomorphism
on $\bigoplus_{\ell_1\neq\ell_2}
 \Hom(E_{\ell_1}^{\bullet},E_{\ell_2}^{\bullet})$
obtained as the inverse of
the adjoint of $\phi$.

\begin{lem}
For $(a,b,c)=(x,y,t),(y,t,x),(t,x,y)$,
the following equalities hold:
\begin{equation}
\label{eq;16.9.10.51}
 \nabla_a\rho_b
-\nabla_b\rho_a
-[\rho_a,\rho_b]^{\top}
-2[\rho_a,\rho_b]^{\bullet}
-[\rho_c,\phi]=0.
\end{equation}
We also have the following equality:
\begin{equation}
\label{eq;16.9.10.52}
 \sum_{\kappa=x,t,y}
 \nabla_{\kappa}\rho_{\kappa}
+F_{\phi}
\left(
 \sum_{\kappa=x,y,t}
 \Bigl(
 2\bigl[\rho_{\kappa},\nabla_{\kappa}\phi\bigr]
-\bigl[\rho_{\kappa},[\rho_{\kappa},\phi]\bigr]
\Bigr)^{\top}
\right)=0.
\end{equation}
\end{lem}
\pf
For $(a,b,c)=(x,y,t),(y,t,x),(t,x,y)$,
we have the equalities
\[
[\nabla_a,\nabla_b]=\nabla_c\phi.
\]
We obtain the following:
\[
 \nabla^{\bullet}_a\rho_b
-\nabla_b^{\bullet}\rho_a
+[\rho_a,\rho_b]^{\top}
-[\rho_c,\phi]=0.
\]
Then, we obtain (\ref{eq;16.9.10.51}).

We have the equality
$\sum_{\kappa=x,y,t}\nabla_{\kappa}\nabla_{\kappa}\phi=0$.
Hence,
we obtain the following:
\[
\sum_{\kappa}
 \bigl[
 \nabla^{\bullet}_{\kappa}\rho_{\kappa},\phi
 \bigr]
+
 \Bigl(
 \sum_{\kappa}
 2[\rho_{\kappa},\nabla_{\kappa}^{\bullet}\phi]
 \Bigr)^{\top}
+
 \Bigl(
 \sum_{\kappa}
 \bigl[\rho_{\kappa},[\rho_{\kappa},\phi]\bigr]
 \Bigr)^{\top}=0
\]
Note that
$\nabla_{\kappa}^{\bullet}\rho_{\kappa}
=\nabla_{\kappa}\rho_{\kappa}$
and
$[\rho_{\kappa},\nabla_{\kappa}\phi]
=[\rho_{\kappa},\nabla^{\bullet}_{\kappa}\phi]
+[\rho_{\kappa},[\rho_{\kappa},\phi]]$.
Then, we obtain
(\ref{eq;16.9.10.52}).
\hfill\qed

\vspace{.1in}

For any $w_{q,0}\in U^{\ast}_{w,q}(R_{12})$
and for any $r_0>0$,
let
$B_1(w_{q,0},r_0)$
denote the connected component of  
$\bigl\{w_q\in U^{\ast}_{w,q}\,\big|\,
|w_q^q-w_{q,0}^q|<r_0
\bigr\}$
such that $w_{q,0}\in B_1(w_{q,0},r_0)$.
We set
$B(w_{q,0},r_0):=S^1_T\times B_1(w_{q,0},r_0)$.
By (\ref{eq;16.9.10.50}),
we have $C_{13}>0$
such that the following holds
if $B(w_{q,0},10)\subset \nbigb_q^{0\ast}(2R_{12})$:
\[
 \bigl\|
 \rho_{\kappa}
 \bigr\|_{h,L^2(B(w_{q,0},10))}
=O\Bigl(
 \exp\bigl(
 -C_{13}I_q(w_{q,0})
 \bigr)
 \Bigr).
\]

Let us apply a standard bootstrapping argument.
Let $\nu:\real_{\geq 0}\lrarr[0,1]$ be a $C^{\infty}$-function
such that
(i) $\nu(u)=1$ $(u\leq 8)$,
(ii) $\nu(u)=0$ $(u\geq 9)$,
(iii) $\nu^{\alpha}$ is $C^{\infty}$ for any $\alpha>0$.
Note that
the $C^{\infty}$-functions
$2\del_u\nu^{\alpha/2}$ $(\alpha>0)$
are equal to $\nu^{-\alpha/2}\del_u\nu^{\alpha}$
on $\{\nu(u)\neq 0\}$.
We set $\chi(w_q)=\rho(|w_q^q-w_{q,0}^q|)$
on $B_1(w_{q,0},10)$
and $\chi(w_q)=0$ for $w_q\not\in B_1(w_{q,0},10)$.
For any $\alpha>0$,
we set
$\rho^{(\alpha)}_{\kappa}:=\chi^{\alpha}\rho_{\kappa}$.
We obtain
\begin{multline}
\label{eq;21.9.14.1}
 \nabla_a(\rho_b^{(\alpha)})
-\nabla_b(\rho_a^{(\alpha)})  
=[\rho_a^{(\alpha/2)},\rho_b^{(\alpha/2)}]^{\top}
+2[\rho_a^{(\alpha/2)},\rho_b^{(\alpha/2)}]^{\bullet} \\
+\chi^{-\alpha/2}(\del_a\chi^{\alpha})\rho_b^{(\alpha/2)}
-\chi^{-\alpha/2}\del_b\chi^{\alpha}\rho_a^{(\alpha/2)},
\end{multline}
\begin{multline}
\label{eq;21.9.14.2}
 \sum_{\kappa}\nabla_{\kappa}(\rho_{\kappa}^{(\alpha/2)})
 =-F_{\phi}\left(
 \sum_{\kappa}\Bigl(
 2\chi^{\alpha/2}[\rho^{(\alpha/2)}_{\kappa},\nabla_{\kappa}\phi]
 -\bigl[\rho_{\kappa}^{(\alpha/2)},
 [\rho^{(\alpha/2)}_{\kappa},\phi]
  \bigr]
 \Bigr)^{\top}
 \right)
 \\
 +\sum_{\kappa} (\chi^{-\alpha/2}\del_{\kappa}\chi^{\alpha})
 \rho^{(\alpha/2)}_{\kappa}.
\end{multline}
There exists $C_{20}>0$, which is independent of $w_{q,0}$,
such that
\[
 \|\rho_{\kappa}^{(\alpha)}\|_{h,L^2(B(w_{q,0},10))}
 \leq
 C_{20}\exp\bigl(
 -C_{13}I_q(w_{q,0})
 \bigr).
\]
Note that the system of differential equations
(\ref{eq;21.9.14.1}, \ref{eq;21.9.14.2})
is elliptic.
Moreover, because the curvature of $\nabla$ is dominated
independently from $w_{q,0}$,
there exists a unitary frame on $B(w_{q,0},10)$
with respect to which the connection form of $\nabla$
is dominated independently from $w_{q,0}$.
Hence, for any $\alpha>0$,
there exist $C_{i}(1,\alpha)>0$  $(i=20,21)$
which are independent of $w_{q,0}$,
such that
\[
  \|\rho_{\kappa}^{(\alpha)}\|_{h,L^2_1(B(w_{q,0},10))}
 \leq
 C_{20}(1,\alpha)\exp\bigl(
 -C_{21}(1,\alpha)I_q(w_{q,0})
 \bigr),
\]
where the Sobolev norms are considered
for the connection $\nabla$.
By an induction of $k\in\seisuu_{>0}$,
we can prove that there exist
$C_{i}(k,\alpha)>0$ $(i=20,21)$
which are independent of $w_{q,0}$,
such that
\[
  \|\rho_{\kappa}^{(\alpha)}\|_{h,L^2_k(B(w_{q,0},10))}
 \leq
 C_{20}(k,\alpha)\exp\bigl(
 -C_{21}(k,\alpha)I_q(w_{q,0})
 \bigr).
\]
Because $\nabla^{\bullet}=\nabla-\rho$,
we obtain the same estimate 
with respect to the Sobolev norms 
for $\nabla^{\bullet}$.
Thus, the proof of Proposition \ref{prop;16.9.10.60}
is completed.
\hfill\qed

\section{Some lemmas from linear algebra}
\label{subsection;20.7.31.41}

\subsection{Eigenvalues}
\label{subsection;16.9.8.3}

Let $d_{\cnum^r}$ denote the standard
Euclidean distance on $\cnum^r$.
Let $d_{\Sym^r\cnum}$ denote the distance on
$\Sym^r\cnum=\cnum^r/\gbigs_r$
given by
\[
 d_{\Sym^r\cnum}(x,y)
=\min_{\substack{x'\in \pi^{-1}(x)\\ y'\in\pi^{-1}(y)}}
 d_{\cnum^r}(x',y'),
\]
where
$\gbigs_r$ denotes the $r$-th symmetric group,
and $\pi:\cnum^r\lrarr\Sym^r\cnum$ denotes the projection.

For any $Y\in M_r(\cnum)$,
let $\Sp(Y)$ denote the set of the eigenvalues of $Y$,
and let $\EE_{\alpha}(Y)$
denote the generalized eigen space
corresponding to $\alpha\in\Sp(Y)$.
Note that
$\Sp(Y)$ induces a point in
$\Sym^r\cnum$.

\begin{lem}
 \label{lem;16.9.8.1}
 There exists $C>0$ depending only on $r$
 such that the following holds.
\begin{itemize}
 \item Let $A\in M_r(\cnum)$ be normal, i.e.,
$\lefttop{t}\Abar\,A=A\,\lefttop{t}\Abar$.
       For any $\epsilon>0$,
       and for any $B\in M_r(\cnum)$ satisfying
       $|B|\leq \epsilon$,
       we obtain
       $d_{\Sym^r\cnum}\bigl(\Sp(A),\Sp(A+B)\bigr)
       \leq C\epsilon$.
\end{itemize}
\end{lem}
\pf
Let $A$, $\epsilon$ and $B$ be as in the statement of
the lemma.
Let us begin with a preliminary.

\begin{lem}
\label{lem;13.11.11.2}
 If there exist $\alpha\in\cnum$
 and a non-zero vector $\vecv\in\cnum^r$
 such that
 $\bigl|
 A\vecv-\alpha\vecv
 \bigr|
\leq
 \epsilon|\vecv|$,
then
there exists an eigenvalue $\beta$
of $A$
such that
$|\beta-\alpha|\leq\epsilon r$.
\end{lem}
\pf
We may assume that
$A$ is a diagonal matrix
whose $(i,i)$-entries are $\gamma_i$.
The assumption can be reworded as follows:
\[
 \sum_{i=1}^r|\gamma_i-\alpha|^2|v_i|^2
\leq
 \epsilon^2|\vecv|^2.
\]
There exists $i_0$ such that
$|v_{i_0}|^2\geq r^{-1}|\vecv|^2$.
Because
\[
 \frac{|\vecv|^2}{r^2}
 |\gamma_{i_0}-\alpha|^2
\leq
 |\gamma_{i_0}-\alpha|^2
 |v_{i_0}|^2
\leq
 \epsilon^2|\vecv|^2,
\]
we obtain
$|\gamma_{i_0}-\alpha|\leq \epsilon r$.
\hfill\qed

\vspace{.1in}

For any $\alpha\in\cnum$ and $\rho>0$,
we set
$D(\alpha,\rho):=
 \bigl\{
 w\in\cnum\,\big|\,
 |w-\alpha|<\rho
 \bigr\}$.
We set
\[
 D(\Sp(A),r\epsilon)
:=\bigcup_{\alpha\in\Sp(A)}
 D(\alpha,r\epsilon).
\]
Lemma \ref{lem;13.11.11.2}
implies
$\Sp(A+B)\subset
 B(\Sp(A),r\epsilon)$.
Let $B(\Sp(A),r\epsilon)
=\coprod_{i\in S(A)}U_i$
denote the decomposition
into the connected components.

\begin{lem}
\label{lem;13.12.5.1}
We obtain the following equality
for each $i\in S(A)$:
\[
 \sum_{\alpha\in U_i\cap\Sp(A)}
 \dim \EE_{\alpha}(A)
=
 \sum_{\alpha\in U_i\cap\Sp(A+B)}
 \dim \EE_{\alpha}(A+B)
\]
\end{lem}
\pf
For the continuous family of matrices
$A_t:=A+tB$ $(0\leq t\leq 1)$,
each $\Sp(A_t)$ is contained in 
$B(\Sp(A),r\epsilon)=\coprod_{i\in S(A)} U_i$.
Because $\sum_{\alpha\in U_i\cap\Sp(A_t)}
 \EE_{\alpha}(A_t)$ is invariant,
the claim of the lemma follows.
\hfill\qed

\vspace{.1in}

Note that
$\max_{\alpha,\beta\in U_i}|\alpha-\beta|<2r^2\epsilon$.
Hence, we obtain Lemma \ref{lem;16.9.8.1}
from Lemma \ref{lem;13.12.5.1}.
\hfill\qed

\subsection{Almost commuting Hermitian matrix and anti-Hermitian
matrix}
\label{subsection;13.11.11.20}

Let 
$\Re:\cnum\lrarr\real$
and $\sqrt{-1}\Image:\cnum\lrarr \sqrt{-1}\real$
denote the maps defined by
$\alpha=\Re(\alpha)+\sqrt{-1}\Image(\alpha)$.
They induce the map
$\Re:\Sym^r\cnum\lrarr\Sym^r\real$
and 
$\sqrt{-1}\Image:\Sym^r\cnum\lrarr\Sym^r(\sqrt{-1}\real)$.

\begin{lem}
 \label{lem;16.9.8.2}
There exists a constant $C>0$
depending only on $r$
such that the following holds:
\begin{itemize}
 \item
      For any $\epsilon>0$,
      a Hermitian matrix $A\in M_r(\cnum)$,
      and an anti-Hermitian matrix $B\in M_r(\cnum)$
      such that
      $\bigl|[A,B]\bigr|\leq\epsilon$,
      we obtain
\[
 d_{\Sym^r\cnum}\bigl(
 \Sp(A),\Re(\Sp(A+B))
 \bigr)
\leq C\epsilon^{1/2},
\]
\[
 d_{\Sym^r\cnum}\bigl(
 \Sp(B),
 \sqrt{-1}\Image(\Sp(A+B))
 \bigr)
\leq
 C\epsilon^{1/2}.
\]
 \item There exists a normal matrix $H\in M_r(\cnum)$
       such that
       (i) the characteristic polynomial of $H$ is
       equal to the characteristic polynomial of $A+B$,
       (ii) $|H-(A+B)|\leq C\epsilon^{1/2}$.      
\end{itemize}
\end{lem}
\pf
Applying Lemma \ref{lem;12.6.20.11}
below to $\Sp(A)$
with $L_1=100$
and $\epsilon_0=\epsilon^{1/2}$,
there exist constants $C_0,C_1>0$
with $C_1/C_0>20$
and $4\leq C_0\leq 4(100 r)^r$
and a decomposition
$\Sp(A)
=\coprod_{i\in \Lambda}
 \nbigs_i$
such that the following holds.
\begin{itemize}
 \item For $\alpha,\beta\in\nbigs_i$,
       we obtain $|\alpha-\beta|\leq C_0\epsilon^{1/2}$.
 \item For $\alpha\in\nbigs_i$
       and $\beta\in\nbigs_j$ with $i\neq j$,
       we obtain
       $|\alpha-\beta|\geq C_1\epsilon^{1/2}$.
\end{itemize}
In the following of this proof,
$C^{(i)}_{10}$ will denote positive constants
depending only on $r$.

We may assume that $A$ is diagonal,
i.e.,
$A=\bigoplus_{i\in \Lambda}\Gamma_i$,
such that $\Gamma_i$ are diagonal
whose $(k,k)$-entries are 
elements of $\nbigs_i$.
Because $|[A,B]|\leq\epsilon$,
for the block decomposition
$B=\sum_{i,j\in\Lambda} B_{ij}$,
we obtain $|B_{ij}|\leq C_{10}^{(1)}\epsilon^{1/2}$ $(i\neq j)$.

We choose $a_i\in\nbigs_i$ for each $i\in\Lambda$,
and we set
$A':=\bigoplus_{i\in\Lambda}a_i I_{r_i}$.
We also put $B':=\bigoplus B_{i,i}$.
 Because $[A',B']=0$,
 we obtain
\[
\Re\bigl(
 \Sp(A'+B')
\bigr)
=\Sp(A'),
\quad
\sqrt{-1}\Image\bigl(
 \Sp(A'+B')
 \bigr)
=\Sp(B').
\]
Because $\bigl|
 (A'+B')-(A+B)
 \bigr|\leq C_{10}^{(2)}\epsilon^{1/2}$,
we obtain
$d_{\Sym^r\cnum}\bigl(
 \Sp(A'+B'),
 \Sp(A+B)
 \bigr) 
\leq C_{10}^{(3)}\epsilon^{1/2}$
by Lemma \ref{lem;16.9.8.1}.
We also obtain
$d_{\Sym^r\cnum}\bigl(
 \Sp(A'),
 \Sp(A)
 \bigr) 
\leq C_{10}^{(4)}\epsilon^{1/2}$
and
$d_{\Sym^r\cnum}\bigl(
 \Sp(B'),
 \Sp(B)
 \bigr) 
\leq C_{10}^{(4)}\epsilon^{1/2}$
by Lemma \ref{lem;16.9.8.1}.
Because
$d_{\Sym^r\cnum}\bigl(
 \Sp(A'+B'),
 \Sp(A+B)
 \bigr) 
\leq C_{10}^{(3)}\epsilon^{1/2}$,
and because $A'+B'$ is normal,
there exists a normal matrix $H$
such that
$\Sp(H)=\Sp(A+B)$
and that
$|H-(A'+B')|\leq C_{10}^{(5)}\epsilon^{1/2}$.
Thus, we are done.
\hfill\qed

\subsection{Decomposition 
of finite tuples in metric spaces (Appendix)}

Let $(X,d)$ be a metric space.
We take
$(x_1,\ldots,x_n)\in X^n$.
Take any $\epsilon_0>0$
and $L_1>8$.

\begin{lem}
\label{lem;12.6.20.11}
There exist 
a non-negative integer $N\leq n$
and a decomposition 
$\{1,\ldots,n\}=
 \coprod_{j\in\Lambda}S_j$
with the following properties:
\begin{itemize}
 \item  For any $a\in S_j$
	and $b\in S_k$
	with $j\neq k$,
	we obtain
	$d\bigl(x_{a},x_b\bigr)>
	(L_1-4)\bigl(L_1n\bigr)^{N}\epsilon_0$.
 \item  For any $a,b\in S_j$,
	we obtain
	$d(x_a,x_b)\leq 4\bigl(L_1n\bigr)^{N}\epsilon_0$.
\end{itemize}
\end{lem}
\pf
We make general preparations.
Let $\Gamma$ be any finite graph.
Let $V(\Gamma)$ denote the set of vertexes of
$\Gamma$.
We obtain the decomposition of the graph
$\Gamma=\coprod_{j\in C(\Gamma)}\Gamma_j$
corresponding to the decomposition of
the geometric realization into the connected components.
Let $|V(\Gamma_j)|$ denote the number of
the vertices of $V(\Gamma_j)$,
and we put $m(\Gamma):=
 \max
 \bigl\{
 |V(\Gamma_j)|\,\big|\,
 j\in C(\Gamma)
 \bigr\}$.
For any positive number $\delta$,
a finite set $I$
and $\vecy=(y_i\,|\,i\in I)\in X^I$,
let $\Gamma(\vecy,\delta)$
denote the unique graph
determined by the conditions.
 \begin{itemize}
\item $V\bigl(\Gamma(\vecy,\delta)\bigr)=I$.
\item
$a,b\in I$ are connected by an edge
 if and only if $d\bigl(y_{a},y_{b}\bigr)\leq \delta$
 and $a\neq b$.
 \end{itemize}

Let us construct a decomposition
as in the claim of the proposition.
We set $S^{(0)}:=\{1,\ldots,m\}$,
$\vecx^{(0)}:=\vecx$,
and 
$\Gamma^{(0)}:=\Gamma(\vecx^{(0)},L_1\epsilon_0)$.
We shall inductively construct a decreasing sequence of subsets
$S^{(0)}\supset S^{(1)}\supset \cdots\supset S^{(N)}$
and graphs $\Gamma^{(0)},\Gamma^{(1)},\ldots,\Gamma^{(N)}$
with $V(\Gamma^{(j)})=S^{(j)}$
such that 
$m(\Gamma^{(j)})>1$ $(j<N)$
and $m(\Gamma^{(N)})=1$.
Suppose that we have already constructed
$(S^{(j)},\Gamma^{(j)})$ 
for $j=0,\ldots,\ell$ with $m(\Gamma^{(j)})>1$ $(j<\ell)$.
If $m(\Gamma^{(\ell)})=1$,
we stop here.
Let us consider the case $m(\Gamma^{(\ell)})>1$.
We have the decomposition
$S^{(\ell)}=\coprod_{j\in C(\Gamma^{(\ell)})}S^{(\ell)}_j$
according to the decomposition of the graph $\Gamma^{(\ell)}$
into the connected components of the geometric realization.
For each $j\in C(\Gamma^{(\ell)})$,
we choose an element $a^{(\ell)}_j\in S^{(\ell)}_j$.
Then, we define
$S^{(\ell+1)}:=\bigl\{
 a^{(\ell)}_j\,\big|\,
 j\in C(\Gamma^{(\ell)})
 \bigr\}$,
$\vecx^{(\ell+1)}:=
 \bigl(
 x_a\,\big|\,
 a\in S^{(\ell+1)}
 \bigr)$ and
\[
 \Gamma^{(\ell+1)}:=
 \Gamma\bigl(\vecx^{(\ell+1)},
 L_1(L_1n)^{\ell+1}\epsilon_0\bigr).
\]

The inductive procedure finishes at some $\ell=N$.
By the construction,
there exist the maps
$\pi_i:S^{(i)}\lrarr S^{(i+1)}$
determined by
$\pi_i(a)=a_j^{(i)}\in S^{(i+1)}$
for $a\in S^{(i)}_j$.
Let $\pi:S\lrarr S^{(N)}$
denote the induced map.
For $c\in S^{(N)}$,
we set $S_c:=\pi^{-1}(c)$.
By the construction,
if $a$ is contained in 
$S_c$,
we obtain
\[
 d(a,c)
\leq
 (L_1n)^{N}\epsilon_0
+(L_1n)^{N-1}\epsilon_0
+\cdots
+(L_1n)\epsilon_0
\leq
 2(L_1n)^{N}\epsilon_0.
\]
Hence,
for $a,b\in S_c$,
we obtain
$d(a,b)\leq 4(L_1n)^{N}\epsilon_0$.

If $a_i\in S_{c_i}$ $(i=1,2)$
with $c_1\neq c_2$,
we obtain
\begin{multline}
 d(a_1,a_2)
\geq
 d(c_1,c_2)
-d(a_1,c_1)-d(a_2,c_2)
\geq
 L_1(L_1n)^{N}\epsilon_0
 -4(L_1n)^{N}\epsilon_0
 \\
\geq
 (L_1-4)(L_1n)^{N}\epsilon_0.
\end{multline}
Hence, the decomposition
$\{1,\ldots,n\}=\coprod_{c\in S^{(\ell)}}S_c$
has the desired property.
\hfill\qed

\section{Vector bundles with a connection on a circle (I)}

\subsection{Statement}
\label{subsection;16.9.11.12}

Let $E$ be a vector bundle on $S^1_T:=\real/T\seisuu$
of rank $r$
with a Hermitian metric $h$,
a unitary connection $\nabla$,
and a self-adjoint section $\psi$.
Let $t$ be the standard coordinate of $\real$,
which induces local  coordinates on $S^1_T$.
For an $\End(E)$-valued differential form $s$,
let $|s|_h$ denote the norm of $s$
with respect to $h$ and the standard metric $dt\,dt$ of $S^1_T$.

We set 
$\nablatilde:=\nabla+\psi\,dt$.
Let $M$ denote the monodromy of
the connection $\nablatilde$
along the loop
$\gamma:[0,1]\lrarr S^1_T$
defined by $\gamma(s)=Ts$.
For each eigenvalue $\alpha$ of $M$,
we have the well defined number
$T^{-1}\log|\alpha|$.
The tuple of 
the numbers $T^{-1}\log|\alpha|$ with the multiplicity
induces an element
in $\Sym^r\real$,
denoted by
$T^{-1}\log|\Sp(M)|$.
For each $Q\in S^1_T$,
the tuple of the eigenvalues with the multiplicity of 
$-\psi_{|Q}\in\End(E_{|Q})$
induces an element of 
$\Sym^r\real$,
denoted by
$\Sp(-\psi_{|Q})$.

We shall prove the following proposition in
\S\ref{subsection;20.7.23.100}--\S\ref{subsection;20.7.23.101}.
\begin{prop}
\label{prop;16.9.11.10}
There exist positive constants
$\epsilon_0>0$ and $C>0$
depending only on $r$
such that the following holds for any $0<\epsilon<\epsilon_0$
\begin{itemize}
 \item 
 If $|\nabla(\psi)|_h\leq \epsilon$,
 we obtain
$d_{\Sym^r\real}\Bigl(
 \Sp(-\psi_{|Q}),\,
 T^{-1}\log|\Sp(M)|
 \Bigr)
\leq
 C\epsilon^{1/2}$.
\end{itemize}
\end{prop}

In \S\ref{subsection;20.7.23.100}--\S\ref{subsection;20.7.23.101},
we suppose $|\nabla\psi|_h\leq \epsilon$
for a positive constant $\epsilon>0$.

\subsection{Preliminary}
\label{subsection;20.7.23.100}

Let $U\in\GL(E_{|0})$ 
denote the monodromy of 
the unitary connection $\nabla$ along
the loop $\gamma$.
There exists $e^{\sqrt{-1}\theta_0}\in S^1$
such that
(i) $0\leq \theta_0\leq 2\pi$,
(ii) $|\theta_0-\arg(\alpha)|> \pi/2r$
for any $\alpha\in\Sp(U)$,
where $\arg(\alpha)$ denotes
any real number such that
$\exp(\sqrt{-1}\arg(\alpha))=\alpha$.
There exists an unitary frame
$\vece$ of $E$
such that
$\nabla\vece=\vece\,(-T^{-1}A)\,dt$,
such that the following holds.
\begin{itemize}
\item
     $A$ is a constant diagonal anti-Hermitian matrix
     such that
\[
     \Sp(A)\subset
     \bigl\{
     \sqrt{-1}a\,\big|\,
     \theta_0-2\pi+\pi/2r<a<\theta_0-\pi/2r
     \bigr\}.
\]          
\end{itemize}
Note that
$U\vece_{|0}
=\vece_{|0}\exp(A)$.

For any $L^2$-section $\rho$ of $\End(E)$,
let $B_{\rho}$ be the matrix valued function
determined by
$\rho\vece=\vece\,B_{\rho}$.
There exists the decomposition
$B_{\rho}=B_{\rho,0}+B_{\rho,1}$,
where $B_{\rho,0}=T^{-1}\int B_{\rho}\,dt$
and $\int B_{\rho,1}\,dt=0$.
It induces a decomposition
$\rho=\rho_0+\rho_1$.
In particular,
we obtain the decomposition
$\psi=\psi_0+\psi_1$.
We obtain a similar decomposition 
of any $L^2$-section of $E$.

\begin{lem}
 \label{lem;17.9.27.30}
There exists a positive constant $C_0>0$
depending only on $r$
such that
$|\nabla\psi_0|_{h}\leq C_0\epsilon$
and
$|\nabla\psi_1|_{h}\leq C_0\epsilon$.
\end{lem}
\pf
Because $A$ is constant,
we obtain
$\nabla(\psi_i)
=(\nabla\psi)_i$ $(i=0,1)$.
Let $B_{\nabla_t\psi}$
be the matrix expressing
$\nabla_t\psi$ with respect to
the frame $\vece$,
i.e.,
$\nabla_t(\psi)\vece=\vece\cdot B_{\nabla_t\psi}$,
then $T^{-1}\int B_{\nabla_t\psi}\,dt$
expresses
$\nabla_t\psi_0$ with respect to
the frame $\vece$.
Hence, there exists $C'_0>0$
depending only on $r$
such that
$|(\nabla_t\psi)_0|_h
\leq
 C'_0\int |\nabla_t\psi|_h\,dt$.
Then, the claim easily follows from
the assumption
$|\nabla\psi|_{h}\leq \epsilon$.
\hfill\qed

\vspace{.1in}
For a section $s$ of $\End(E)$,
we set
\[
|s|_{L^2}:=\left(\int_{0}^T|s|_h^2\right)^{1/2},
\quad
|s|_{L_k^2}=\left(
\sum_{0\leq j\leq k}|\nabla_t^js|^2_{L^2}
\right)^{1/2}.
\]

\begin{cor}
There exist $C_1>0$,
depending only on $r$,
such that 
$|\psi_1|_{L_1^2}\leq C_1\epsilon$
and
$\sup|\psi_1|_h\leq C_1\epsilon$.
\end{cor}
\pf
Note that 
the eigenvalues $\sqrt{-1}\beta$ of $\ad(A)$
satisfy
$|\beta|<2\pi-(\pi/r)$.
Hence, 
we obtain
$|\psi_1|_{L_1^2}\leq C_1'\epsilon$
from the estimate of $|\nabla\psi_1|$
in the previous lemma.
Because $\dim S^1_T=1$,
we obtain 
$\sup|\psi_1|_h\leq C_1''\epsilon$.
\hfill\qed

\subsection{A decomposition of function spaces}

By the frame $\vece$,
we obtain a $C^{\infty}$-isometry of vector bundles
$\End(E)\simeq S^1_T\times M_r(\cnum)$
on $S^1_T$.
For $k\geq 0$,
let $L_k^2(M_r(\cnum))$
denote the space of
$L_k^2$-maps from $S^1_T$ to $M_r(\cnum)$.
Let $L_k^2(M_r(\cnum))_1$
denote the subspace of
$F\in L_k^2(M_r(\cnum))$
such that
$\int_{S^1_T}F\,dt=0$.
Let $L_k^2(M_r(\cnum))_0$
denote the space of the constant maps
from $S^1_T$ to $M_r(\cnum)$.
We obtain the decomposition
$L_k^2(M_r(\cnum))
=L_k^2(M_r(\cnum))_0\oplus
 L_k^2(M_r(\cnum))_1$.
We obtain the corresponding decomposition
\[
 L_k^2\bigl(\End(E)\bigr)
=L_k^2\bigl(\End(E)\bigr)_0
\oplus
 L_k^2\bigl(\End(E)\bigr)_1.
\]
Namely,
$L_k^2\bigl(\End(E)\bigr)_0$
denotes the space of
sections of $\End(E)$,
which are constant with respect to $\vece$,
and 
$L_k^2\bigl(\End(E)\bigr)_1$
denotes the space of 
$L_k^2$-sections of $\End(E)$
which are represented by matrices $B$
with respect to $\vece$
such that $\int_{S^1_T} B\,dt=0$.
There exist similar decompositions
for $C^{\infty}(\End(E))$
and $L_k^2\bigl(\End(E)\otimes\Omega^1\bigr)$,
etc.

\begin{lem}
\label{lem;20.7.23.110}
 There exist $C_2>0$ and $\epsilon_1>0$,
depending only on $r$,
with the following property.
\begin{itemize}
\item
If $\epsilon\leq \epsilon_1$,
we obtain
$|a|_{L_1^2}\leq
 C_2\bigl|
 (\nabla_{t}+\psi_0)a
 \bigr|_{L^2}$
for any 
$a\in L_1^2(\End(E))_1$.
\end{itemize}
\end{lem}
\pf
In this proof,
$C_2^{(i)}$ denote positive constants
depending only on $r$.
Let $B\in M_r(\cnum)$ be determined by
$\psi_0\vece=\vece\,B$,
which satisfies $\lefttop{t}\Bbar=B$.
By the assumption $|\nabla(\psi)|\leq \epsilon$
and Lemma \ref{lem;17.9.27.30},
we obtain
$\bigl|
 [A,B]\bigr|\leq C_2^{(1)}\epsilon$.
We set $L=-T^{-1}A+B$.
Because
$\lefttop{t}\Lbar=T^{-1}A+B$,
we obtain
$\bigl|
 \bigl[
 L,\lefttop{t}\Lbar
 \bigr]
\bigr|
\leq C_2^{(2)}\epsilon$.
By Lemma \ref{lem;16.9.8.2},
there exists a decomposition
$L=L_1+L_2$
such that
(i) $L_1$ is normal,
and the characteristic polynomials of $L_1$ and $L$
 are the same,
(ii) $|L_2|\leq C_2^{(3)}\epsilon^{1/2}$ holds.
Note that $L_1$ is diagonalizable
by a unitary matrix,
and that any eigenvalue $\gamma$ of $\ad L_1$
satisfies $|\Image(\gamma)|<2\pi-(\pi/2r)$.
Let $\nabla^{(0)}$ be determined by
$\nabla_{t}^{(0)}\vece=\vece\,L_1$,
and let $g^{(0)}$ be determined by
$g^{(0)}\vece=\vece\,L_2$.
We can easily check that
there exists $C_2^{(4)}>0$ depending only on $r$
such that
$|a|_{L_1^2}\leq
 C_2^{(4)}\bigl|\nabla_{t}^{(0)}a\bigr|_{L^2}$
for any $a\in L_1^2(\End(E))_1$.
Because $|g^{(0)}|\leq C_2^{(3)}\epsilon^{1/2}$
and
$\bigl|
 \nabla_{t}^{(0)}a
 \bigr|_{L^2}
 \leq
 \bigl|(\nabla_{t}+\psi_0)a\bigr|_{L^2}
 +\bigl|[g^{(0)},a]
 \bigr|_{L^2}$,
there exists $C_2^{(5)}>0$ such that
$\bigl|
 \nabla_{t}^{(0)}a
 \bigr|_{L^2}
\leq
 C_2^{(5)}\bigl|(\nabla_{t}+\psi_0)a\bigr|_{L^2}$
if $\epsilon_1$ is sufficiently small.
Then, the claim of the lemma follows.
\hfill\qed

\subsection{Gauge transformation}

Let $\epsilon_1$ be as in Lemma \ref{lem;20.7.23.110}.
We assume that $\epsilon<\epsilon_1$.
We take 
a neighbourhood $\nbigu_0\subset L_1^2(\End(E))_1$
of $0$
such that
$1+a$ is invertible for any
$a\in\nbigu_0$.

\begin{prop}
\label{prop;16.9.11.11}
There exist positive constants
$\epsilon_2$ and $C_3$,
depending only on $r$,
such that the following holds:
\begin{itemize}
\item
If $\epsilon\leq \epsilon_2$,
there exists
$(a,b)\in\nbigu_0\times L^2(\End(E))_0$
satisfying
$\bigl|(\nabla_{t}+\psi_0)a\bigr|_{L^2}
\leq  C_3\epsilon$,
$|b|\leq C_3\epsilon$,
and 
\[
(1+a)^{-1}\circ(\nabla_{t}+\psi_0+b)\circ(1+a)
=\nabla_{t}+\psi.
\]
\end{itemize}
\end{prop}
\pf
We define the map
$\Psi:
 \nbigu_0\times
 L^2(\End(E))_0
\lrarr
 L^2(\End(E))$
by
\[
  \Psi(a,b)
=(1+a)^{-1}\circ(\nabla_{t}+\psi_0+b)\circ(1+a)
-(\nabla_{t}+\psi_0).
\]
The derivative 
$T_{(a,b)}\Psi:
 L_1^2(\End(E))_1\oplus
 L^2(\End(E))_0
\lrarr
 L^2(\End(E))$
of $\Psi$ at $(a,b)$
is as follows:
\begin{multline}
 T_{(a,b)}\Psi(u,v)
\equiv
 (1+a+u)^{-1}\circ(\nabla_{t}+\psi_0+b+v)\circ(1+a+u)
-(1+a)^{-1}\circ(\nabla_{t}+\psi_0+b)\circ(1+a)
 \\
\equiv
-(1+a)^{-1}\circ u\circ (1+a)^{-1}\circ
 (\nabla_{t}+\psi_0+b)\circ(1+a)
 \\
+(1+a)^{-1}\circ(\nabla_{t}+\psi_0+b)
 \circ(1+a)\circ(1+a)^{-1}\circ u
+(1+a)^{-1}\circ v\circ (1+a).
\end{multline}
Let us obtain a bound of 
$T_{(a,b)}\Psi(u,v)
-\bigl((\nabla_{t}+\psi_0)u+v\bigr)$.
Note that
$(\nabla_{t}+\psi_0)u
=(\nabla_{t}+\psi_0)\circ u
-u\circ(\nabla_{t}+\psi_0)$.
We clearly have
$\Ad(1+a)^{-1}(v)-v
=O\bigl(|a|_{L_{1}^2}\cdot |v|\bigr)$.
We set
\begin{equation}
 \gbiga:=
 (1+a)^{-1}\circ b\circ (1+a)\circ (1+a)^{-1}\circ u
-(1+a)^{-1}\circ u\circ (1+a)^{-1}\circ b\circ(1+a),
\end{equation}
\begin{multline}
 \gbigb:=
 (1+a)^{-1}\circ(\nabla_{t}+\psi_0)\circ(1+a)
 \circ(1+a)^{-1}\circ u
-(1+a)^{-1}\circ u\circ (1+a)^{-1}\circ
 (\nabla_{t}+\psi_0)\circ(1+a)
 \\
-(\nabla_{t}+\psi_0)\bigl((1+a)^{-1}u\bigr),
\end{multline}
\begin{equation}
 \gbigc:=(\nabla_{t}+\psi_0)
 \bigl(
 (1+a)^{-1}-1
\bigr)u.
\end{equation}
We obtain
\begin{multline}
 (1+a)^{-1}\circ\bigl(\nabla_{t}+\psi_0+b\bigr)
 \circ(1+a)\circ(1+a)^{-1}\circ u
 \\
 -(1+a)^{-1}\circ u\circ(1+a)^{-1}
 \circ\bigl(\nabla_{t}+\psi_0+b\bigr)\circ(1+a)
-(\nabla_{t}+\psi_0)u
=\gbiga+\gbigb+\gbigc.
\end{multline}
We have
$|\gbiga|_{L^2}=O\bigl(
 |b|\cdot|u|_{L_1^2}
 (1+|a|_{L_1^2})
 \bigr)$.
Because
$\gbigb=\bigl[
 (1+a)^{-1}\circ(\nabla_{t}+\psi_0)(a),
 (1+a)^{-1}\circ u
 \bigr]$,
we obtain
\[
 |\gbigb|_{L^2}
 =O\Bigl(
 C^{(0)}_3(1+|a|_{L_1^2})^2
 \bigl|
 (\nabla_{t}+\psi_0)a
 \bigr|_{L^2}
 \,|u|_{L_{1}^2}
 \Bigr).
\]
For the term $\gbigc$,
we have the following equality:
\[
\gbigc=(\nabla_{t}+\psi_0)
 \bigl(
 ((1+a)^{-1}-1)u
 \bigr)
=\bigl((1+a)^{-1}-1\bigr)\circ
 (\nabla_{t}+\psi_0)(u)
+(\nabla_{t}+\psi_0)\bigl((1+a)^{-1}\bigr)\circ u.
\]
There exist $C^{(1)}_{3},C^{(2)}_{3}>0$
such that 
\begin{multline}
 |\gbigc|_{L^2}
\leq
 C_{3}^{(1)}|a|_{L_{1}^2}
 \bigl|(\nabla_{t}+\psi_0)u\bigr|_{L^2}
+C_3^{(1)}\bigl|
 (\nabla_{t}+\psi_0)\bigl((1+a)^{-1}\bigr)
 \bigr|_{L^2}
 |u|_{L_{1}^2}
 \\
\leq
  C_3^{(2)}|a|_{L_{1}^2}
 \bigl|(\nabla_{t}+\psi_0)u\bigr|_{L^2}
+C_3^{(2)}
 (1+|a|_{L_1^2})
 \bigl|
 (\nabla_{t}+\psi_0)a
 \bigr|_{L^2}
 |u|_{L_{1}^2}.
\end{multline}
Therefore, there exists $C_3^{(3)}>0$
such that the following holds:
\begin{multline}
 \Bigl|
 T_{(a,b)}\Psi(u,v)-
 \bigl(
 (\nabla_{t}+\psi_0)u+v
 \bigr)
\Bigr|_{L^2}
 \\
\leq
 C_3^{(3)}\Bigl(
 |a|_{L_1^2}|v|
+ |b|\,|u|_{L_{1}^2}(1+|a|_{L_1^2})
+(1+|a|_{L_{1}^2})^2
 \bigl|(\nabla_{t}+\psi_0)a\bigr|_{L^2}
 |u|_{L_{1}^2}
\Bigr)
 \\
+C_3^{(3)}\Bigl(
 |a|_{L_1^2}\bigl|(\nabla_{t}+\psi_0)u\bigr|_{L^2}
+
 (1+|a|_{L_1^2})
 \bigl|
 (\nabla_{t}+\psi_0)(a)
 \bigr|_{L^2}
 |u|_{L_{1}^2}
\Bigr).
\end{multline}

Let $\Gamma$ be the composite of the following maps:
{\small
\[
 \begin{CD}
\nbigu_0\times L^2(\End(E))_0
@>{\Psi}>>
 L^2(\End(E))
@>{(\nabla_{t}+\psi_0)^{-1}\oplus\id}>>
 L_1^2(\End(E))_1
\oplus
 L^2(\End(E))_0.
 \end{CD}
\]
}
We regard
$L_1^2(\End(E))_1$
as a Banach space
by the norm
$a\longmapsto
 \bigl|(\nabla_{t}+\psi_0)a\bigr|_{L^2}$.
Hence, the operator norm of 
$T_{a,b}\Gamma-\id$
is 
$O\Bigl(
\bigl(
 |(\nabla_t+\psi_0)a|_{L_1^2}
 +\bigl|b\bigr|
 \bigr)\cdot (1+|a|_{L_1^2})^2
 \Bigr)$.
There exists $C_{4}>0$
depending only on $r$
such that
if $\bigl|(\nabla_{t}+\psi_0)a\bigr|_{L^2}+|b|\leq C_4$
then 
$|T_{a,b}\Gamma-\id|<1/2$.
As in the proof of 
the inverse mapping theorem (see \cite{Lang}),
the image of $\Gamma$ contains 
\[
 \Bigl\{
c+d\in L_1^2(\End(E))_1\oplus L^2(\End(E))_0\,\Big|\,
 \bigl|(\nabla_{t}+\psi_0)c\bigr|_{L^2}+|d|
\leq
 C_4/2
\Bigr\}.
\]
It means that the image of 
$\Psi$ contains
$\Bigl\{
c+d\in L^2(\End(E))\oplus L^2(\End(E))_0\,\Big|\,
 \bigl|c\bigr|_{L^2}+|d|
\leq
 C_4/2
\Bigr\}$.
Because we have
 $|\psi_1|_{L^2}\leq C_1\epsilon$,
we obtain the claim of the proposition.
\hfill\qed

\subsection{Proof of Proposition \ref{prop;16.9.11.10}}
\label{subsection;20.7.23.101}

Let us prove Proposition \ref{prop;16.9.11.10}.
We use the notation in Proposition \ref{prop;16.9.11.11}.
There exist
$B,B'\in M_r(\cnum)$
determined by
$\psi_0\vece=\vece B$
and $b\vece=\vece B'$.
Let $\nablatilde':=\nabla+(\psi_0+b)\,d\theta$.
We have
$\nablatilde'\vece
=\vece\,(-T^{-1}A+B+B')\,d\theta$.
Let $\Re\Sp(T^{-1}A-B-B')$
denote the point of $\Sym^r\real$
induced by the real part of the eigenvalues of $T^{-1}A-B-B'$
with the multiplicity.
Let $M'$ denote the monodromy of
$\nablatilde'$,
which is conjugate to
$\exp\bigl(A-TB-TB')$.
Because $M$ and $M'$ are conjugate,
we have
$T^{-1}\log|\Sp(M)|
=T^{-1}\log|\Sp(M')|
=\Re\Sp(T^{-1}A-B-B')$.
Then, the claim follows from 
Lemma \ref{lem;16.9.8.2}.
\hfill\qed

\section{Vector bundles with a connection on a circle (II)}
\label{subsection;20.7.31.42}

\subsection{Additional assumption on the eigenvalues
 of the monodromy}
\label{subsection;13.11.12.4}

We continue to use the notation in 
\S\ref{subsection;16.9.11.12}.
We assume that 
 $|\nabla(\psi)|_{h}\leq \epsilon<\epsilon_0$
 as in Proposition \ref{prop;16.9.11.10}.
We impose the following additional assumptions
on the eigenvalues of $\psi$.
\begin{condition}
\label{condition}\mbox{{}}
\begin{description}
 \item[(A1)]
There exist positive constants
 $C_{i}$ $(i=10,11)$ with
 $1<C_{10}$ and  $100C_{10}<C_{11}$,
and an element
$\nbigs\in\Sym^r\real$
for which the following holds
for any $Q\in S^1_T$:
\[
 d_{\Sym^r\cnum}\bigl(
 \Sp(-\psi_{|Q}), \nbigs
 \bigr)
<C_{10}/2,
\quad\quad
 C_{11}<\min\bigl\{
 |\gamma_1-\gamma_2|\,\big|\,
 \gamma_i\in\nbigs,
 \gamma_1\neq\gamma_2
 \bigr\}.
\]
We also assume that $\epsilon_0$ is so small
that 
$d_{\Sym^r\cnum}\bigl(
 T^{-1}\log\bigl|\Sp(M(\nablatilde))\bigr|,\nbigs
\bigr)
<C_{10}$.
(See Proposition {\rm\ref{prop;16.9.11.10}}.)
\item[(A2)]
There exist $k_0\geq 1$ and $C_{12}>0$
such that
$|\nabla^j(\psi)|_h\leq C_{12}\epsilon$
for any $1\leq j\leq k_0$.
\hfill\qed
\end{description}
\end{condition}
 
We naturally regard $\nbigs$ as a subset of $\real$.
We obtain the orthogonal decomposition
$(E,\psi)
=\bigoplus_{\alpha\in\nbigs}
 (E_{\alpha}^{\bullet},\psi_{\alpha})$
such that 
the eigenvalues $\gamma$ of 
$-\psi_{\alpha|Q}$ 
satisfy $|\gamma-\alpha|<C_{10}/2$.
It induces the decomposition
$\nabla_{t}=\nabla_{t}^{\bullet}
+\rho$,
where
$\nabla_{t}^{\bullet}$
is a unitary connection preserving the decomposition,
and
$\rho$ is a section of
$\bigoplus_{\alpha\neq\beta}
 \Hom(E_{\alpha}^{\bullet},E_{\beta}^{\bullet})$.
We obtain the decomposition
$\nabla_t^{\bullet}=
 \bigoplus_{\alpha\in\nbigs}
 \nabla^{\bullet}_{\alpha,t}$.
In general,
for any vector bundle $E'$
and for any section $s$ of 
$\End(E)\otimes E'$,
we obtain the decomposition
$s=s^{\bullet}+s^{\top}$
into the sections of
$\bigoplus_{\alpha\in\cnum}
 \End(E^{\bullet}_{\alpha})\otimes E'$
and 
$\bigoplus_{\alpha\neq\beta}
 \Hom(E^{\bullet}_{\alpha},E^{\bullet}_{\beta})
\otimes E'$.

\begin{lem}
There exists $C_{13}>0$,
depending only on
$r$, $k_0$ and $C_i$ $(i=10,11,12)$,
such that
the following holds
for any $1\leq j\leq k_0$:
\begin{equation}
\label{eq;16.9.14.10}
\bigl|
 (\nabla_{t}^{\bullet})^j\psi
 \bigr|_{h}
\leq C_{13}\epsilon,
\quad
 \bigl|
 \bigl[
 \psi,(\nabla_{t}^{\bullet})^{j-1}\rho
 \bigr]
 \bigr|_{h}
\leq 
C_{13}\epsilon,
\quad
\bigl|
 (\nabla_{t}^{\bullet})^{j-1}\rho
 \bigr|_{h}
\leq
 C_{13}\epsilon.
\end{equation}
\end{lem}
\pf
There exists the following description:
\[
 \nabla_{t}^j\psi
=(\nabla_{t}^{\bullet})^j\psi
+\ad\bigl(
 (\nabla_{t}^{\bullet})^{j-1}\rho
\bigr)\psi
+\sum_{i\leq j-1}
 \nbigp_{j,i}\cdot(\nabla_{t}^{\bullet})^i(\psi).
\]
Here, $\nbigp_{j,i}$
are given as a sum of 
some compositions of
$\ad(\nabla^{m}_{t}\rho)$ $(m<j-1)$.
Note the equalities
$(\nabla_{t}^{\bullet})^j\psi
=\bigl(
 (\nabla_{t}^{\bullet})^j\psi
 \bigr)^{\bullet}$
and 
$\ad\bigl(
 (\nabla_{t}^{\bullet})^{j-1}\rho
\bigr)\psi
=\Bigl(
 \ad\bigl(
 (\nabla_{t}^{\bullet})^{j-1}\rho
\bigr)\psi
 \Bigr)^{\top}$.
We also remark that
an estimate for 
$\bigl|
 \bigl[
 \psi,(\nabla_{t}^{\bullet})^{j-1}\rho
 \bigr]
 \bigr|_{h}$
 implies an estimate for
$\bigl|
 (\nabla_{t}^{\bullet})^{j-1}\rho
 \bigr|_h$.
Then, we obtain (\ref{eq;16.9.14.10})
by an easy induction.
\hfill\qed

\begin{lem}
\label{lem;16.9.9.1}
There exists $C_{14}>0$ and $\epsilon_{10}>0$
depending only on $r$, $k_0$ and $C_i$ $(i=10,11,12)$,
such that the following holds
for any $k\leq k_0-1$
and for any $a\in
 \bigoplus_{\alpha\neq\beta}
 L_{k+1}^2(\Hom(E_{\alpha}^{\bullet},E_{\beta}^{\bullet}))$
if $\epsilon\leq \epsilon_{10}$:
\[
 |a|_{L_{k+1}^2}
\leq
 C_{14}\bigl|
 (\nabla^{\bullet}_{t}+\psi)(a)
 \bigr|_{L_k^2}.
\]
\end{lem}
\pf
We set 
$\Gamma_0:=
 \bigoplus \alpha\cdot \id_{E_{\alpha}^{\bullet}}$.
Because $\Gamma_0$ is flat with respect to
$\nabla^{\bullet}$, and self-adjoint with respect to the metric,
we obtain
$|(\nabla^{\bullet}_t+\Gamma_0)a|_{L^2}\geq |\nabla^{\bullet}_ta|_{L^2}$
and 
$|(\nabla^{\bullet}_t+\Gamma_0)a|_{L^2}\geq |a|_{L^2}$.
We also obtain 
\[
 \bigl|
 (\nabla^{\bullet}_{t})^j
 (\nabla^{\bullet}_{t}+\Gamma_0)(a)
 \bigr|_{L^2}
\geq 
 \bigl|(\nabla^{\bullet}_{t})^{j+1} a\bigr|_{L^2},
\]
\[
  \bigl|
 (\nabla^{\bullet}_{t})^j
 (\nabla^{\bullet}_{t}+\Gamma_0)(a)
 \bigr|_{L^2}
\geq
 \bigl|(\nabla^{\bullet}_{t})^{j}a\bigr|_{L^2}.
\]
Then, the claim of the lemma follows.
\hfill\qed

\subsection{Gauge transformations}
\label{subsection;13.11.12.101}
Let $0\leq k\leq k_0$.
We use the norm on
$\bigoplus_{\alpha\neq \beta}
 L_{k+1}^2(\Hom(E_{\alpha}^{\bullet},E_{\beta}^{\bullet}))$
given by
$\bigl|
 (\nabla^{\bullet}_{t}+\psi)(a)
 \bigr|_{L_k^2}$.
Let
\[
 \nbigu_0\subset
 \bigoplus_{\alpha\neq\beta}
 L_{k+1}^2(\Hom(E_{\alpha}^{\bullet},E_{\beta}^{\bullet}))
\]
be a small neighbourhood of $0$
such that
$1+a$ is invertible
for any $a\in \nbigu_0$.
Let
\[
 \nbigu_1\subset
 \bigoplus_{\alpha} L_k^2(\End(E_{\alpha}^{\bullet}))
\]
be a small neighbourhood of $0$.
We define the map
$\Phi:\nbigu_0\times\nbigu_1\lrarr
 L_k^2\bigl( \End(E) \bigr)$
given by
\[
 \Phi(a,b)
:=(1+a)^{-1}(\nabla_{t}^{\bullet}+\psi)(a)
+(1+a)^{-1}b(1+a).
\]

\begin{prop}
\label{prop;16.9.14.11}
There exist positive constants
$C_{15}$ and $\epsilon_{11}$
depending only on $r$, $k_0$ and $C_i$ $(i=10,11,12)$
such that the following holds:
\begin{itemize}
\item
If $\epsilon<\epsilon_{11}$,
there uniquely exists
$(a,b)\in\nbigu_0\times\nbigu_1$
satisfying the following conditions for any $k\leq k_0$:
\begin{equation}
\label{eq;16.9.14.20}
 \Phi(a,b)=\rho,
\quad
 \bigl|
 (\nabla_{t}^{\bullet}+\psi)(a)
 \bigr|_{L^2_k}
\leq
 C_{15}|\rho|_{L^2_k},
 \quad
 |b|_{L^2_k}\leq C_{15}|\rho|_{L^2_k}.
\end{equation}
Moreover, $a$ and $b$ are $C^{\infty}$.
\end{itemize}
\end{prop}
\pf
The derivative of $\Phi$
at $(a,b)$ is as follows:
\begin{multline}
 T_{(a,b)}\Phi(u,v)
=(1+a)^{-1}(\nabla_{t}^{\bullet}+\psi)(u)
-(1+a)^{-1}u(1+a)(\nabla_{t}^{\bullet}+\psi)(a)
\\
+(1+a)^{-1}v(1+a)-(1+a)^{-1}u(1+a)^{-1}b(1+a)
+(1+a)^{-1}bu.
\end{multline}
Hence, 
there exists a positive constant $C_{16}$
depending only on 
$r$, $k_0$ and $C_i$ $(i=10,11,12)$
such that the following holds
for $(u,v)\in T_a\nbigu_0\oplus v\in T_b\nbigu_1$:
\begin{multline}
\bigl|
 T_{a,b}\Phi(u,v)
 -(\nabla_t^{\bullet}+\psi)u-v
 \bigr|_{L_k^2}
 \leq
 C_{16}
(1+|a|_{L_{k+1}^2})^2
 |(\nabla^{\bullet}_t+\psi)a|_{L_k^{2}}
 \bigl|
 (\nabla^{\bullet}_t+\psi)u
 \bigr|_{L_k^2}
 \\
+C_{16}(1+|a|_{L_{k+1}^2})^3|b|_{L_{k}^2}|u|_{L_k^2} 
 +C_{16}
 (1+|a|_{L_{k+1}^2})|a|_{L_{k+1}^2}|v|_{L_{k}^2}.
\end{multline}
Then, we obtain $(a,b)$ satisfying
(\ref{eq;16.9.14.20})
as in the case of Proposition \ref{prop;16.9.11.11}.

Let us observe that $(a,b)$ is $C^{\infty}$.
Note that
$(\nabla^{\bullet}_{t}+\psi)(a)
+ba
=\bigl((1+a)\rho\bigr)^{\top}$
and 
$b=\bigl((1+a)\rho\bigr)^{\bullet}$.
By the second equality,
$b$ is $L_{k+1}^2$.
By the first, we obtain that $a$ is $L_{k+2}^2$.
By an easy inductive argument,
we obtain that $a$ and $b$ are $L_{\ell}^2$ for any $\ell$.
\hfill\qed

\subsection{Comparison of the decompositions}

Let $\nablatilde:=\nabla+\psi\,dt$.
By Proposition \ref{prop;16.9.11.10},
there exists the $\nablatilde$-invariant decomposition 
$E=\bigoplus_{\alpha\in\nbigs} E_{\alpha}$
with the following property:
\begin{itemize}
\item
 For any eigenvalue $\gamma$
 of $M(\nablatilde)$ on $E_{\alpha}$,
 we obtain
 $\bigl|T^{-1}\log|\gamma|-\alpha\bigr|\leq C_{10}$. 
\end{itemize}

We assume that 
$\epsilon<\epsilon_{11}$.
Let $(a,b)$ be as in Proposition {\rm\ref{prop;16.9.14.11}}.
\begin{prop}
We obtain
$E_{\alpha}^{\bullet}=(1+a)E_{\alpha}$.
\end{prop}
\pf
We set 
$\Etilde_{\alpha}:=(1+a)^{-1}E_{\alpha}^{\bullet}$.
The connection
$\nabla_{t}^{\bullet}+\psi+b$ preserves
the decomposition
$E=\bigoplus E_{\alpha}^{\bullet}$.
Because
$\nabla_{t}+\psi
=(1+a)^{-1}(\nabla^{\bullet}_{t}+\psi+b)(1+a)$,
the connection
$\nabla_{t}+\psi$ preserves the decomposition
$\bigoplus \Etilde_{\alpha}$.
By considering the eigenvalues of the monodromy
on $\Etilde_{\alpha}$,
we obtain that
$\Etilde_{\alpha}=E_{\alpha}$.
\hfill\qed

\begin{cor}
There exists a constant $C_{20}$
depending only on $r$, $C_i$ $(i=10,11,12)$
such that the following holds:
\begin{itemize}
 \item
For any $Q\in S^1_T$,
and 
for any $u_{\alpha}\in E_{\alpha|Q}$
and $u_{\beta}\in E_{\beta|Q}$
with $\alpha\neq\beta$,
we obtain
\[
 \bigl|
h(u_{\alpha},u_{\beta})
\bigr|
\leq
 C_{20}
 \bigl|
 \rho
 \bigr|_{L^2}\cdot
|u_{\alpha}|_h\cdot|u_{\beta}|_h.
\]
\end{itemize}
\end{cor}
\pf
Note that
$\sup|a|_h\leq C_{21}|a|_{L_k^1}
\leq |\rho|_{L^2}$.
Then, the claim follows from
$h\bigl(
 (1+a_{|Q})u_{\alpha},
 (1+a_{|Q})u_{\beta}
 \bigr)=0$.
\hfill\qed

\vspace{.1in}

We set
$\End(E)^{\circ}:=
 \bigoplus_{\alpha\in\nbigs}
 \End(E_{\alpha})$
and 
$\End(E)^{\bot}:=
 \bigoplus_{\alpha\neq \beta}
 \Hom(E_{\alpha},E_{\beta})$.
 For any vector bundle $E'$,
 and for any section $g$ of $\End(E)\otimes E'$,
there exists the corresponding decomposition
$g=g^{\circ}+g^{\bot}$.

\begin{cor}
There exists a positive constant
$C_{21}>0$
such that the following holds:
\begin{itemize}
\item
Let $g$ be any section of 
$\End(E)^{\bullet}\otimes\Omega^p$.
Then,
$\bigl|
 g^{\bot}
 \bigr|_{L_k^2}
\leq
 C_{21}\cdot
 \bigl|
 \rho
 \bigr|_{L_k^2}\cdot
 \bigl|g\bigr|_{L_k^2}$.
\end{itemize}
\end{cor}
\pf
Let $\pi_{\alpha}$ denote the projection 
of $E$ onto $E_{\alpha}\subset E$
with respect to the decomposition
$E=\bigoplus E_{\alpha}$.
Let $\pi_{\alpha}^{\bullet}$ denote 
the orthogonal projection of $E$
onto $E_{\alpha}^{\bullet}\subset E$.
We have
$\pi_{\alpha}=(1+a)\circ\pi_{\alpha}^{\bullet}\circ(1+a)^{-1}$.

We have
$g^{\bot}=
 \sum_{\alpha\neq\beta}
 \pi_{\alpha}\circ g\circ \pi_{\beta}$.
By using
$\pi_{\alpha}^{\bullet}\circ g\circ \pi_{\beta}^{\bullet}=0$
for any $\alpha\neq\beta$,
we can easily obtain the desired estimate.
\hfill\qed

\section{Proof of Theorem \ref{thm;16.9.20.22}}
\label{subsection;17.10.5.2}

Let $(E,h,\nabla,\phi)$ be a monopole
on $\nbigb^{0\ast}_q(R)$
satisfying Condition \ref{condition;21.8.21.1}
and Condition \ref{condition;21.8.21.2}
as in \S\ref{subsection;16.9.28.1}.
There exists the decomposition (\ref{eq;17.9.27.6}).
We put $S(E):=\{\ell\,|\,E_{\ell}\neq 0\}$.
We set $\psi=-\sqrt{-1}\phi$.

\begin{lem}
There exist positive numbers $C_0$, $R_0$
and an orthogonal decomposition
$(E,\phi)=\bigoplus_{\ell\in S(E)}
  (E_{\ell}^{\bullet},\phi_{\ell}^{\bullet})$
on $\nbigb_q^{0\ast}(R_0)$
such that 
any eigenvalue $\alpha$ of
$\phi^{\bullet}_{\ell|(t,w_q)}$
satisfies
$\Bigl|
 \alpha-\sqrt{-1}\ell T^{-1}\log|w_q|
 \Bigr|
\leq C_0$.
In other words,
Condition {\rm\ref{condition;21.8.21.3}}
is satisfied.
\end{lem}
\pf
It follows from Proposition \ref{prop;16.9.11.10}.
\hfill\qed

\begin{prop}
\label{prop;16.9.14.101}
There exists a section $a$
of $\bigoplus_{\ell_1\neq\ell_2}
 \Hom(E_{\ell_1}^{\bullet},E_{\ell_2}^{\bullet})$
such that the following holds:
\begin{itemize}
\item
$E_{\ell}^{\bullet}=(1+a)E_{\ell}$
for any $\ell\in S(E)$.
\item
For any $k\in\seisuu_{\geq 0}$,
there exist positive constants $C_{i}(k)>0$ $(i=10,11)$
such that 
\[
 \bigl|
 \nabla_{\kappa_1}\circ
 \cdots \circ\nabla_{\kappa_{k}}(a)
\bigr|_h
\leq 
 C_{10}(k)
 \cdot
 \exp\bigl(
 -C_{11}(k)I_q(w_q)
 \bigr)
\]
for any $(\kappa_1,\ldots,\kappa_{k})\in\{t,x,y\}^{k}$:
\end{itemize}
\end{prop}
\pf
From the decomposition $E=\bigoplus E^{\bullet}_{\ell}$,
we obtain the decomposition 
$\End(E)=\End(E)^{\bullet}\oplus \End(E)^{\top}$.
For any section $s$ of
$\End(E)\otimes\Omega^p$,
we use the decomposition
$s=s^{\bullet}+s^{\top}$
into the section $s^{\bullet}$ of
$\End(E)^{\bullet}\otimes\Omega^p$
and the section $s^{\top}$
of $\End(E)^{\top}\otimes\Omega^p$.
We obtain
the decomposition 
$\nabla=\nabla^{\bullet}+\rho$
as in \S\ref{subsection;21.8.21.4}.
The norm of $\rho$ is dominated by
Proposition \ref{prop;16.9.10.60}.

By Proposition \ref{prop;16.9.14.11},
there exist $R_{1}>R$
and a section 
$(a,b)$ of  
$\End(E)^{\top}\oplus 
 \End(E)^{\bullet}$
on $\nbigb^{0\ast}_q(R_{1})$
such that 
\[
 (1+a)^{-1}
 \circ
 (\nabla^{\bullet}_{t}+\psi+b)
 \circ
 (1+a)
=\nabla^{\bullet}_t+\rho_t+\psi.
\]
Here, 
for each $w_q\in U^{\ast}_{w,q}(R_{1})$,
the restriction
$(a,b)_{|S^1_T\times \{w_q\}}$
is $C^{\infty}$,
and we obtain
\[
 \bigl|
 (\nabla^{\bullet}_{t}+\psi)
 a_{|S^1_T\times \{w_q\}}\bigr|_{L_{k}^2}
+\bigl|
 b_{|S^1_T\times\{w_q\}}
 \bigr|_{L_k^2}
\leq C(k)
 \bigl|
 \rho_{t|S^1_T\times\{w_q\}}
 \bigr|_{L^2_k}
\]
for any $k$.
It is easy to see that
$(a,b)$ are $C^{\infty}$-sections of
$\End(E)^{\top}\oplus \End(E)^{\bullet}$
by using the inverse function theorem.

Because
$(\nabla^{\bullet}_{t}+\psi)(a)
+ba=\bigl((1+a)\rho_t\bigr)^{\top}$
and 
$b=\bigl((1+a)\rho_t\bigr)^{\bullet}$,
we obtain
\[
 (\nabla^{\bullet}_{t}+\psi)(a)
+\bigl((1+a)\rho_t\bigr)^{\bullet}a
-\bigl((1+a)\rho_t\bigr)^{\top}=0.
\]
For any 
$\veckappa=(\kappa_1,\ldots,\kappa_{k})\in \{t,x,y\}^k$,
we set
$k:=|\veckappa|$
and 
$\nabla^{\bullet}_{\veckappa}:=
 \nabla^{\bullet}_{\kappa_1}\circ\cdots\circ
 \nabla^{\bullet}_{\kappa_{k}}$.
By an inductive argument,
we obtain equalities
\[
 \bigl(
 \nabla^{\bullet}_{t}+\psi+\nbigb_{\veckappa}
 \bigr)\nabla^{\bullet}_{\veckappa}a
+\nbigc_{\veckappa}=0,
\]
where $\nbigb_{\veckappa}$
and $\nbigc_{\veckappa}$
are constructed by some linear algebraic operations
from 
$\nabla^{\bullet}_{\veckappa_1}\psi$,
$\nabla^{\bullet}_{\veckappa_2}a$
and 
$\nabla^{\bullet}_{\veckappa_3}\rho$,
where $|\veckappa_i|<|\veckappa|$.
We obtain the estimates for
$\bigl|
 \nabla^{\bullet}_{\veckappa}(a)_{|S^1_T\times \{w_q\}}
 \bigr|_{L^2}$
by an easy induction.
Therefore, we obtain the claim of
Proposition \ref{prop;16.9.14.101}.
\hfill\qed

\vspace{.1in}

Let us complete the proof of Theorem \ref{thm;16.9.20.22}.
For any endomorphism $f$ of $E$,
let $f_h^{\dagger}$ denote the adjoint
with respect to $h$.
Let $\Pi_{\ell}^{\bullet}$ denote the projection
$E\lrarr E_{\ell}^{\bullet}\subset E$
with respect to the orthogonal decomposition
$E=\bigoplus E_{\ell}^{\bullet}$.
Note that
$(\Pi_{\ell}^{\bullet})_h^{\dagger}=\Pi_{\ell}^{\bullet}$
and 
$\sum_{\ell\in S(E)}
 (\Pi_{\ell}^{\bullet})_h^{\dagger}
\circ
 \Pi_{\ell}^{\bullet}=\id_E$.
Let $\Pi_{\ell}$ denote the projection
$E\lrarr E_{\ell}$
with respect to the decomposition
$E=\bigoplus_{\ell\in S(E)} E_{\ell}$.
Note that
\[
 \Pi_{\ell}
=(\id_E+a)^{-1}\circ
 \Pi_{\ell}^{\bullet}\circ
 (\id_E+a).
\]

Let $s$ be the automorphism of $E$
determined by $h^{\circ}=h\cdot s$.
Then, the following holds.
\[
 s=\sum_{\ell\in S(E)}
 (\Pi_{\ell})^{\dagger}_h\circ
 \Pi_{\ell}.
\]
We obtain the desired estimate for $s$
from Proposition \ref{prop;16.9.14.101}.
Thus, we obtain Theorem \ref{thm;16.9.20.22}.
\hfill\qed

\chapter{The filtered bundles associated with periodic monopoles}
\label{section;20.7.25.10}

In \S\ref{section;21.8.23.1},
we summarize some notation, which will be used in this chapter.
In \S\ref{section;21.7.7.1}--\S\ref{subsection;17.10.10.10},
we shall prove that a periodic monopole
satisfying the GCK-condition
induces a good filtered bundle.
Moreover, the norm estimate is satisfied.
In \S\ref{subsection;20.8.1.3},
we shall study the converse.
Namely, we shall prove that
if a monopole $(E,\delbar_E,h)$ 
is strongly adapted to
a good filtered bundle,
then the monopole satisfies the GCK-condition.

\section{Notation}
\label{section;21.8.23.1}

\subsection{Some spaces and morphisms}
\label{subsection;21.8.22.1}

We summarize some notation.
For any $q\in\seisuu_{\geq 1}$,
let $w_q$ be a $q$-th root of the variable $w$
such that $w_{q\ell}^{\ell}=w_q$ for any $q,\ell\in\seisuu_{>0}$.
\index{variable $w_q$}
Let $\proj^1_{w_q}\lrarr \proj^1_w$
be the ramified covering induced by
$w_q\longmapsto w_q^q$.
The pull back of $w$ is also denoted by $w$,
i.e., $w=w_q^q$ on $\proj^1_{w_q}$.

Let $\lambda\in\cnum$. 
We have the map
$\Psi^{\lambda}:\nbigmbar^{\lambda}\lrarr \cnum_{w}$
induced by
$(t_1,\beta_1)\longmapsto
(1+|\lambda|^2)^{-1}(\beta_1-2\sqrt{-1}\lambda t_1)$.
\index{map $\Psi^{\lambda}$}
(See \S\ref{subsection;17.10.2.1}
for the mini-complex manifold $\nbigmbar^{\lambda}$
and mini-complex local coordinate systems $(t_1,\beta_1)$.)
For $q\in\seisuu_{\geq 1}$,
let $Y^{\lambda}_{q}$
be the fiber product of
$\nbigmbar^{\lambda}$ and $\proj^1_{w_q}\setminus\{0\}$
over $\proj^1_{w_q}$.
\index{space $Y^{\lambda}_q$}
Let $\Psi^{\lambda}_q:
Y^{\lambda}_q\lrarr \proj^1_{w_q}\setminus\{0\}$
denote the induced proper map.
\index{map $\Psi^{\lambda}_q$}
We set
$H^{\lambda}_{\infty,q}=(\Psi^{\lambda}_q)^{-1}(\infty)$
and 
$Y^{\lambda\ast}_q=
Y^{\lambda}_q\setminus H^{\lambda}_{\infty,q}
=(\Psi^{\lambda}_q)^{-1}(\proj^1_{w_q}\setminus\{0,\infty\})$.
\index{space $H^{\lambda}_{\infty,q}$}
\index{space $Y^{\lambda\ast}_q$}
Let $\pi^{\lambda}_q:Y^{\lambda}_q\lrarr S^1_T$
denote the morphism obtained as the composite
$Y^{\lambda}_q\lrarr \nbigmbar^{\lambda} \lrarr S^1_T$.
\index{projection $\pi^{\lambda}_q$}

Let $\varpi^{\lambda}:\Mbar^{\lambda}\lrarr \nbigmbar^{\lambda}$
denote the projection in \S\ref{subsection;21.8.21.10}.
\index{projection $\varpi^{\lambda}$}
Let $Y^{\lambda\cov}_q$ denote the fiber product of
$Y^{\lambda}_q$ and $\Mbar^{\lambda}$ over $\nbigmbar^{\lambda}$.
\index{space $Y^{\lambda\cov}_q$}
Let $\varpi^{\lambda}_q:Y^{\lambda\cov}_q\lrarr Y^{\lambda}_q$
denote the projection.
\index{projection $\varpi^{\lambda}_q$}
Note that $Y^{\lambda\cov}_q$ is equipped with
the naturally induced free $\seisuu$-action $\kappa_q$,
and $Y^{\lambda}_q$ is identified with the quotient space.
We set
$H^{\lambda\cov}_{\infty,q}:=
(\varpi_q^{\lambda})^{-1}(H^{\lambda}_{\infty,q})$
and
$Y^{\lambda\cov\ast}_q:=
Y^{\lambda\cov}\setminus H^{\lambda\cov}_{\infty,q}$.
\index{space $H^{\lambda\cov}_{\infty,q}$}
\index{space $Y^{\lambda\cov\ast}_q$}

We recall that
$Y^{\lambda}_q$ and $Y^{\lambda\cov}_q$
are naturally equipped with
the mini-complex structures.
(See Lemma \ref{lem;21.8.21.11}.)
We have the natural local mini-complex coordinate systems
given by $(t_1,\beta_{1,q}^{-1})$,
where $\beta_{1,q}$ denotes a $q$-th root of
the variable $\beta_1$.
\index{variable $\beta_{1,q}$}

Let $\Hhat^{\lambda}_{\infty,q}$
denote the formal space obtained
as the completion of $Y^{\lambda}_q$
along $H^{\lambda}_{\infty,q}$.
\index{ringed space $\Hhat^{\lambda}_{\infty,q}$}
(See \S\ref{subsection;20.7.22.100} for
more details on $\Hhat^{\lambda}_{\infty,q}$.)
Let $\Hhat^{\lambda\cov}_{\infty,q}$
denote the ringed space obtained
as the completion of $Y^{\lambda\cov}_q$
along $H^{\lambda\cov}_{\infty,q}$.
\index{ringed space $\Hhat^{\lambda\cov}_{\infty,q}$}

From the ramified covering map
$\varphi_{q,p}:\proj^1_{w_p}\lrarr \proj^1_{w_q}$
for $p\in q\seisuu_{\geq 1}$,
we obtain the induced ramified covering maps
$\nbigr_{q,p}:
Y^{\lambda}_p\lrarr Y^{\lambda}_q$.
\index{map $\varphi_{q,p}$}
\index{map $\nbigr_{q,p}$}
The induced map
$H^{\lambda}_{\infty,p}\lrarr H^{\lambda}_{\infty,q}$
is an isomorphism,
and
the induced map
$Y^{\lambda\ast}_p\lrarr Y^{\lambda\ast}_q$
is a covering map.

\vspace{.1in}

Note that 
$Y^{\lambda\ast}_q$ and $Y^{\lambda\cov\ast}_q$
are equipped with the Riemannian metric
induced by the Riemannian metric of $\nbigm^{\lambda}$,
which are locally Euclidean
compatible with the mini-complex structures.
Because $\nbigm^0=\nbigm^{\lambda}$
as Riemannian manifolds,
we have
$Y^{\lambda\ast}_q=Y^{0\ast}_q$
and
$Y^{\lambda\cov\ast}_q=Y^{0\cov\ast}_q$
as Riemannian manifolds.
Because the maps
$\varpi^{\lambda}_q:Y^{\lambda\cov\ast}_q\lrarr Y^{\lambda\ast}_q$
and
$\varpi^{0}_q:Y^{0\ast\cov}_q\lrarr Y^{0\ast}_q$
are the same,
i.e., independent of $\lambda$,
the restriction of
$\varpi^{\lambda}_q$ to $Y^{\lambda\cov\ast}_q$
is also denoted simply by $\varpi_q$.
Similarly,
because the maps
$\Psi^{\lambda}_q:Y^{\lambda\ast}_q\lrarr
 \proj^1_{w_q}\setminus\{0,\infty\}$
and 
$\Psi^{0}_q:Y^{0\ast}_q\lrarr
 \proj^1_{w_q}\setminus\{0,\infty\}$
are the same,
the restriction of 
$\Psi^{\lambda}_q$
to $Y^{\lambda\ast}_q$ is also denoted simply by
$\Psi_q$.

The composite
$\Psi^{\lambda}_q\circ
\varpi^{\lambda}_q:Y_q^{\lambda\cov}\lrarr \proj^1_{w_q}\setminus\{0\}$
is denoted by $\Psi_q^{\lambda\cov}$.
\index{morphism $\Psi_q^{\lambda\cov}$}
The restriction
$Y_q^{\lambda\cov\ast}\lrarr
\proj^1_{w_q}\setminus\{0,\infty\}$
is independent of $\lambda$,
which is also denoted by $\Psi_q^{\cov}$.

\subsection{Neighbourhoods of infinity}
\label{subsection;21.8.21.20}

Let $U_w$ be a neighbourhood of $\infty$ in $\proj^1_w$.
\index{space $U_w$}
We set $U_w^{\ast}=U_w\setminus\{\infty\}$.
\index{space $U_w^{\ast}$}
Let $U_{w,q}$ and $U_{w,q}^{\ast}$
denote the pull back of $U_w$ and $U_{w}^{\ast}$
by the ramified covering
$\proj^1_{w_q}\lrarr \proj^1_{w}$, respectively.
\index{space $U_{w,q}$}
\index{space $U_{w,q}^{\ast}$}
We set
$\nbigb_q^{\lambda}:=(\Psi_q^{\lambda})^{-1}(U_{w,q})$
and 
$\nbigb_q^{\lambda\ast}:=\nbigb_q^{\lambda}\setminus H^{\lambda}_{\infty,q}$.
\index{space $\nbigb_q^{\lambda}$}
\index{space $\nbigb_q^{\lambda\ast}$}
The induced map
$\nbigb_q^{\lambda}\lrarr U_{w,q}$
is also denoted by $\Psi^{\lambda}_q$,
and the restriction
$\nbigb_q^{\lambda\ast}\lrarr U_{w,q}^{\ast}$
is also denoted simply by $\Psi_q$.

For $t_1\in S^1_T$,
we set
$\nbigb^{\lambda}_q\langle t_1\rangle
:=(\pi_q^{\lambda})^{-1}(t_1)\cap\nbigb^{\lambda}_q$.
\index{space $\nbigb^{\lambda}_q\langle t_1\rangle$}
Similarly,
we set
$H^{\lambda}_{\infty,q}\langle t_1\rangle
=(\pi_q^{\lambda})^{-1}(t_1)\cap H^{\lambda}_{\infty,q}$
and 
$\nbigb^{\lambda\ast}_q\langle t_1\rangle
=(\pi^{\lambda}_q)^{-1}(t_1)\cap
\nbigb^{\lambda\ast}_q$.
\index{space $H^{\lambda}_{\infty,q}\langle t_1\rangle$}
\index{space $\nbigb^{\lambda\ast}_q\langle t_1\rangle$}
The point
$H^{\lambda}_q\langle t_1\rangle
\in\nbigb^{\lambda}_q\langle t_1\rangle$
is also denoted by $\infty$
if there is no risk of confusion.
The formal completion of
$\nbigb^{\lambda}_q\langle t_1\rangle$
along $H^{\lambda}_{\infty,q}\langle t_1\rangle$
is denoted by
$\Hhat^{\lambda}_{\infty,q}\langle t_1\rangle$.
\index{ringed space $\Hhat^{\lambda}_{\infty,q}\langle t_1\rangle$}

We set
$\nbigb^{\lambda\cov}_q:=
(\varpi^{\lambda}_q)^{-1}(\nbigb^{\lambda}_q)$
and
$\nbigb^{\lambda\cov\ast}_q:=
\nbigb^{\lambda\cov}\setminus H^{\lambda\cov}_{\infty,q}$.
\index{space $\nbigb^{\lambda\cov}_q$}
\index{space $\nbigb^{\lambda\cov\ast}_q$}
The induced maps
$\nbigb_q^{\lambda\cov}\lrarr \nbigb_q^{\lambda}$
and 
$\nbigb_q^{\lambda\cov}\lrarr U_{w,q}$
are also denoted by
$\varpi^{\lambda}_q$
and
$\Psi^{\lambda\cov}_q$,
and the restrictions to $\nbigb^{\lambda\cov\ast}$
are also denoted simply by
$\varpi_q$ and $\Psi_q^{\cov}$.
\index{projection $\varpi^{\lambda}_q$}
\index{map $\Psi^{\lambda\cov}_q$}
\index{projection $\varpi_q$}
\index{map $\Psi_q^{\cov}$}

We emphasize
$\nbigb^{\lambda\ast}_q=\nbigb^{0\ast}_q$
and
$\nbigb^{\lambda\cov\ast}_q=\nbigb^{0\cov\ast}_q$
as Riemannian manifolds.

\subsection{Norm of differential forms}

For a vector bundle $E$ with a Hermitian metric $h$
on $\nbigb^{\lambda\ast}_q$,
and for an $\End(E)$-valued differential form $s$,
let $|s|_{h}$ denote the induced function on $\nbigb^{\lambda\ast}_q$
obtained as the norm of $s$ with respect to $h$ and
the natural Riemannian metric
$dt\,dt+dw\,d\wbar$ of $\nbigb^{\lambda\ast}_q$.
Similarly,
for a vector bundle $E$ with a Hermitian metric $h$
on $U_{w,q}^{\ast}$,
and for an $\End(E)$-valued differential form $s$,
let $|s|_{h}$ denote the induced function on $U_{w,q}^{\ast}$
obtained as the norm of $s$ with respect to $h$ and
the natural Riemannian metric $dw\,d\wbar$ of $U_{w,q}^{\ast}$.
We omit to denote the dependence on the metrics
of the base spaces, because they are fixed.
\index{norm \mbox{$|s|_h$}}

\section{Meromorphic prolongation}
\label{section;21.7.7.1}
 
\subsection{Statements}
\label{subsection;17.10.6.100}

Let $U_w$ be a neighbourhood of $\infty$ in $\proj^1_w$.
Let $q\in\seisuu_{\geq 1}$.
For each $\lambda\in\cnum$,
we obtain
the induced neighbourhood
$\nbigb^{\lambda}_q=(\Psi^{\lambda}_q)^{-1}(U_{w,q})$
of $H^{\lambda}_{\infty,q}$ in $Y^{\lambda}_q$
as in \S\ref{subsection;21.8.21.20},
and we set
$\nbigb^{\lambda\ast}_q:=
\nbigb^{\lambda}_q\setminus H^{\lambda}_{\infty,q}$.

Let $(E,h,\nabla,\phi)$ be a monopole
on $\nbigb_{q}^{0\ast}$
satisfying the GCK-condition.
(See Proposition \ref{prop;20.7.29.20}.)
Fix $\lambda\in\cnum$.
Let $(E^{\lambda},\delbar_{E^{\lambda}})$ 
denote the mini-holomorphic bundle
on $\nbigb^{\lambda\ast}_{q}$
underlying the monopole $(E,h,\nabla,\phi)$.
(See \S\ref{subsection;21.8.23.2}.)
 We obtain the
 $\nbigo_{\nbigb^{\lambda}_q}(\ast H^{\lambda}_{\infty,q})$-module
 $\nbigp^hE^{\lambda}$
 by the procedure in \S\ref{subsection;17.10.25.100}.
We shall prove the following proposition 
in \S\ref{subsection;17.10.6.2}.
\begin{prop}
\label{prop;17.9.17.50}
$\nbigp^h(E^{\lambda})$ is a locally free
$\nbigo_{\nbigb^{\lambda}_{q}}(\ast H^{\lambda}_{\infty,q})$-module.
\end{prop}

Let $\del_{E^{\lambda},h,\beta_1}$
be obtained from $\del_{E^{\lambda},\betabar_1}$ and $h$
as in \S\ref{subsection;17.10.5.120}.
According to Proposition \ref{prop;17.9.29.101},
there exists the following relation:
\begin{equation}
 \label{eq;17.10.26.11}
 \bigl[
 \del_{E^{\lambda},\betabar_1},
 \del_{E^{\lambda},h,\beta_1}
 \bigr]
=(1+|\lambda|^2)^{-1}\cdot F_{\wbar,w}.
\end{equation}
We obtain the following lemma from 
Corollary \ref{cor;17.12.16.2}
and Corollary \ref{cor;17.10.6.1}.
\begin{lem}
\label{lem;21.7.9.3}
 $\bigl(E^{\lambda},\delbar_{E^{\lambda}},h
 \bigr)_{|\nbigb^{\lambda\ast}_q\langle t_1\rangle}$
 is acceptable for any $t_1\in S^1_T$.
 (See {\rm\S\ref{subsection;20.8.8.31}}
 for the acceptability condition.)
\hfill\qed
\end{lem}

As in \S\ref{subsection;17.10.5.122},
for each $t_1\in S^1_T$,
we obtain
the locally free
$\nbigo_{\nbigb^{\lambda}_q\langle t_1\rangle}(\ast \infty)$-module
$\nbigp^h\bigl(
E^{\lambda}_{|
 \nbigb^{\lambda\ast}_q\langle t_1\rangle}
\bigr)$,
and the filtered bundle
$\nbigp^h_{\ast}\bigl(
 E^{\lambda}_{|\nbigb^{\lambda\ast}_q\langle t_1\rangle}
\bigr)$ 
over
$\nbigp^h\bigl(
E^{\lambda}_{|\nbigb^{\lambda\ast}_q\langle t_1\rangle}\bigr)$.
We shall prove the following lemma 
in \S\ref{subsection;17.10.6.2}.
\begin{lem}
\label{lem;17.10.5.300}
 $\nbigp^h(E^{\lambda})
 _{|\nbigb^{\lambda}_q\langle t_1\rangle }$
is naturally isomorphic to
 $\nbigp^h\bigl(
 E^{\lambda}
 _{|\nbigb^{\lambda\ast}_q\langle t_1\rangle}
 \bigr)$.
\end{lem}

Let $P$ be any point of $H^{\lambda}_{\infty,q}$.
We take a mini-holomorphic frame $\vecu=(u_1,\ldots,u_r)$ of
$\nbigp^h(E^{\lambda})$
on a neighbourhood $U_P$ of $P$.
Let $H(h,\vecu)$ be the Hermitian matrix valued function
on $U_P\setminus H^{\lambda}_{\infty,q}$
determined by $H(h,\vecu)_{ij}:=h(u_i,u_j)$.
\index{function $H(h,\vecu)$}
We shall prove the following lemma in \S\ref{subsection;17.10.6.2}.
\begin{lem}
\label{lem;17.10.7.9}
There exist $C_1\geq 1$ and $N>0$
such that 
$C_1^{-1}|x_q|^{-N}
\leq
 |H(h,\vecu)|
\leq
 C_1|x_q|^N$.
\end{lem}

\subsection{Proof}
\label{subsection;17.10.6.2}

Take $P\in H^{\lambda}_{\infty,q}$.
We put $t_1^{0}=\pi^{\lambda}_q(P)\in S^1_T$.
Let $U_P$ be a neighbourhood of $P$
in $\nbigb_q^{\lambda}$
of the form
$\{|t_1-t_1^{0}|<\epsilon_1\}
\times
\{|\beta_{1,q}^{-1}|<\epsilon_2\}$.
Let $s$ be a holomorphic section of
$\nbigp^h\bigl(E^{\lambda}
_{|\nbigb_q^{\lambda\ast}\langle t_1\rangle}\bigr)$.
The restriction
$s_{|U_P\cap\nbigb^{\lambda\ast}_q\langle t_1\rangle}$
uniquely extends to
a mini-holomorphic section
$\stilde$ of $E$
on $U_P\setminus H^{\lambda}_{\infty,q}$.

\begin{lem}
\label{lem;17.10.6.200}
$\stilde$ is a section of
$\nbigp^h(E^{\lambda})$
on $U_P$.
More precisely,
there exist $C_1\geq 1$ and $N_1>0$  such that 
\[
 C_1^{-1}|\beta_1|^{-N_1}\cdot |s|_h(t_1^{0},\beta_{1,q}^{-1})
\leq
 \bigl|\stilde\bigr|_h(t_1,\beta_{1,q}^{-1})
\leq
 C_1|\beta_1|^{N_1}\cdot |s|_h(t_1^{0},\beta_{1,q}^{-1}).
\]
\end{lem}
\pf
Because
$\del_{E^{\lambda},t_1}=\nabla_{t_1}-\sqrt{-1}\phi$.
we obtain
$\nabla_{t_1}\stilde=\sqrt{-1}\phi \stilde$.
We obtain
\[
 \left|
 \del_{t_1}\bigl|\stilde\bigr|_h^2
\right|
=\Bigl|
 2\Re h(\sqrt{-1}\phi \stilde,\stilde)
 \Bigr|
\leq
 C_1\log|\beta_1-2\sqrt{-1}\lambda t_1|\cdot |\stilde|_h^2
\leq
 \bigl(
 C_2\log|\beta_1|+C_3
 \bigr)\cdot|\stilde|_h^2.
\]
Hence, we obtain the claim of the lemma.
\hfill\qed

\vspace{.1in}
By Lemma \ref{lem;17.10.6.200},
the following morphism is an isomorphism:
\begin{equation}
 \label{eq;17.10.6.101}
  H^0\Bigl(
  U_P,
 \nbigp^h(E^{\lambda})
 \Bigr)
\lrarr
 H^0\Bigl(
U_P\cap \nbigb^{\lambda}_q\langle t_1^{0}\rangle,
 \nbigp^h
 (E^{\lambda}_{|\nbigb^{\lambda\ast}_q\langle t_1^{0}\rangle})
 \Bigr).
\end{equation}
Then, we obtain Lemma \ref{lem;17.10.5.300}.

The following natural morphism is an isomorphism:
\[
 A:=
 H^0\Bigl(U_P,
 \nbigo_{\nbigb^{\lambda}_{q}}(\ast H^{\lambda}_{\infty,q})
 \Bigr)
\lrarr
 H^0\Bigl(\nbigb^{\lambda}_q\langle t_1^{0}\rangle,
 \nbigo_{\nbigb^{\lambda}_q\langle t_1^{0}\rangle}(\ast \{\infty\})
 \Bigr).
\]
As remarked in \S\ref{subsection;17.10.6.100},
$\nbigp^h(
E^{\lambda}_{|\nbigb^{\lambda\ast}_q\langle t_1^{0}\rangle})$
is a locally free sheaf,
and hence
the both side of (\ref{eq;17.10.6.101})
are free $A$-modules.
Hence, there exists a local frame of
$\nbigp^h(E^{\lambda})$,
i.e., we obtain Proposition \ref{prop;17.9.17.50}.
We also obtain Lemma \ref{lem;17.10.7.9}
from Lemma \ref{lem;17.10.6.200}.
\hfill\qed

\section{Filtered prolongation}

\label{subsection;17.10.10.10}

\subsection{Statement}

We continue to use the notation in \S\ref{subsection;17.10.6.100}.
As a result of Lemma \ref{lem;17.10.5.300},
the family
\index{Filtered bundle $\nbigp^h_{\ast}(E^{\lambda})$}
\[
 \nbigp^h_{\ast}\bigl(
 E^{\lambda}
 \bigr)
:=
 \Bigl(
 \nbigp^h_{\ast}\bigl(
 E^{\lambda}_{|\nbigb^{\lambda\ast}_q\langle t_1\rangle}
\bigr)
 \,\Big|\,
 t_1\in S^1_T
 \Bigr)
\]
is a filtered bundle over
$\nbigp^h\bigl(E^{\lambda} \bigr)$.
 \index{filtered bundle $\nbigp^h_{\ast}\bigl(
 \nbigr_{q,p}^{-1}E^{\lambda}
 \bigr)$}
 (See Definition \ref{df;20.7.24.1}.)
We shall prove the following proposition
in this subsection.

\begin{thm}
\label{thm;17.10.5.130}
The filtered bundle
$\nbigp^h_{\ast}\bigl(
E^{\lambda}
\bigr)$
is good.
(See Definition {\rm\ref{df;21.8.19.2}}.)
\end{thm}

\subsection{Refined statement}
\label{subsection;20.8.8.100}

Suppose that Condition \ref{condition;21.8.21.2}
is satisfied for $(E,h,\nabla,\phi)$,
i.e., there exists a decomposition of
the mini-holomorphic bundle
$(E^0,\delbar_{E^0})$ on $\nbigb^{0\ast}_q$
as in (\ref{eq;17.9.27.4}):
\begin{equation}
\label{eq;21.8.21.40}
 (E^0,\delbar_{E^0})=
  \bigoplus_{(\ell,\alpha)\in\seisuu\times\cnum^{\ast}}
 (E^0_{\ell,\alpha},\delbar_{E^0_{\ell,\alpha}}).
\end{equation}
(See Remark \ref{rem;21.7.7.2}.)
We shall explain a more precise claim
in this situation.

There exist the holomorphic bundles
$(V_{\ell,\alpha},\delbar_{\ell,\alpha})$ on $U_{w,q}^{\ast}$
with a holomorphic automorphism $f_{\ell,\alpha}$,
and isomorphisms of mini-holomorphic bundles
as in (\ref{eq;17.10.6.300}):
\begin{equation}
\label{eq;21.8.21.41}
 (E^0_{\ell,\alpha},\delbar_{E^0_{\ell,\alpha}})
  \simeq
  \LL^{0\ast}_q(\ell,\alpha)
  \otimes
  \Psi_q^{\ast}(V_{\ell,\alpha},f_{\ell,\alpha}).
\end{equation}

\subsubsection{The filtered bundles associated with
the basic monopoles of rank one}

For $(\ell,\alpha)\in\seisuu\times\cnum^{\ast}$,
there exists the monopole
$(\LL^{\ast}_q(\ell,\alpha),h_{\LL,q,\ell,\alpha},
\nabla_{\LL,q,\ell,\alpha},\phi_{\LL,q,\ell,\alpha})$ on
$Y_q^{\lambda\ast}=Y_q^{0\ast}$
as in (\ref{eq;21.8.16.12}).
Let $\LL^{\lambda\ast}_q(\ell,\alpha)$
denote the underlying mini-holomorphic bundle
on $Y_q^{\lambda\ast}$.
\index{mini-holomorphic bundle $\LL^{\lambda\ast}_q(\ell,\alpha)$}
By the procedure in \S\ref{subsection;17.10.25.100},
it extends to
a locally free $\nbigo_{Y^{\lambda}_q}(\ast H^{\lambda}_q)$-module
$\LL_q^{\lambda}(\ell,\alpha)$,
and the good filtered bundle
$\nbigp_{\ast}\LL_q^{\lambda}(\ell,\alpha)$.
\index{good filtered bundle $\nbigp_{\ast}\LL_q^{\lambda}(\ell,\alpha)$}
There exists an isomorphism
\begin{equation}
\label{eq;21.8.21.22}
 \nbigp_{\ast}\LL^{\lambda}_q(\ell,\alpha)_{|\Hhat^{\lambda}_{\infty,p}}
\simeq
 \nbigp^{(0)}_{\ast}\LLhat^{\lambda}_q(\ell,\alpha).
\end{equation}
(See \S\ref{subsection;17.10.7.1}
for $\nbigp_{\ast}^{(0)}\LLhat^{\lambda}_q(\ell,\alpha)$.)

\subsubsection{The filtered bundles associated with
the asymptotic harmonic bundles}
\label{subsection;21.8.21.100}

We set $\theta_{\ell,\alpha}=f_{\ell,\alpha}\,dw$.
As explained in 
\S\ref{subsection;17.10.5.50}--\ref{subsection;17.10.6.301},
we obtain
Hermitian metrics $h_{V,\ell,\alpha}$
of
$(V_{\ell,\alpha},\delbar_{\ell,\alpha},\theta_{\ell,\alpha})$
for which the estimate (\ref{eq;17.10.6.302}) holds.
Let $\delbar_{\ell,\alpha}+\del_{\ell,\alpha}$ be the Chern connection
determined by $\delbar_{\ell,\alpha}$ and $h_{V,\ell,\alpha}$.
Let $\theta_{\ell,\alpha}^{\dagger}$
denote the adjoint of $\theta_{\ell,\alpha}$.
According to the estimate for asymptotic harmonic bundles
in \cite[\S5.5]{Mochizuki-doubly-periodic}
(see also Proposition \ref{prop;12.8.7.3} below),
we obtain
\[
 \bigl[\delbar_{\ell,\alpha},\del_{\ell,\alpha}\bigr]
=O\bigl(|w_q|^{-2}(\log|w_q|)^{-2}dw_q\,d\wbar_q\bigr),
\]
\[
 \bigl[\theta_{\ell,\alpha},\theta^{\dagger}_{\ell,\alpha}\bigr]
=O\bigl(|w_q|^{-2}\log|w_q|^{-2}dw_q\,d\wbar_q\bigr).
\]

Let $V_{\ell,\alpha}^{\lambda}$ be the holomorphic bundle
$(V_{\ell,\alpha},\delbar_{\ell,\alpha}+\lambda\theta^{\dagger}_{\ell,\alpha})$
with the metric $h_{V,\ell,\alpha}$.
The Chern connection is
$\delbar_{\ell,\alpha}+\lambda\theta^{\dagger}_{\ell,\alpha}
+\del_{\ell,\alpha}-\lambdabar\theta_{\ell,\alpha}$,
and the curvature is
\begin{equation}
\label{eq;17.10.6.10}
 \bigl[
 \delbar_{\ell,\alpha}+\lambda\theta^{\dagger}_{\ell,\alpha},
 \del_{\ell,\alpha}-\lambdabar\theta_{\ell,\alpha}
 \bigr] 
=O\bigl(|w_q|^{-2}(\log|w_q|)^{-2}dw_q\,d\wbar_q\bigr).
\end{equation}
Hence, 
we obtain the locally free
$\nbigo_{U_{w,p}}(\ast\infty)$-module
$\nbigp^{h_{V,\ell,\alpha}}V^{\lambda}_{\ell,\alpha}$
and the filtered bundle
$\nbigp^{h_{V,\ell,\alpha}}_{\ast}V^{\lambda}_{\ell,\alpha}$
over $\nbigp^{h_{V,\ell,\alpha}}V^{\lambda}_{\ell,\alpha}$.

Set $\DD^{\lambda\,(1,0)}:=
 \lambda\del_{\ell,\alpha}+\theta_{\ell,\alpha}$
and 
$\DD^{\lambda}=\DD^{\lambda\,(1,0)}+
 \delbar_{\ell,\alpha}+\lambda\theta_{\ell,\alpha}^{\dagger}$.
Let $\vecv_{\ell,\alpha}$ be a holomorphic frame 
of $\nbigp^{h_{V,\ell,\alpha}}_aV^{\lambda}_{\ell,\alpha}$.
Let $A$  be the matrix valued function 
on $U_{w,q}^{\ast}$
determined by
$\DD^{\lambda(1,0)}\vecv_{\ell,\alpha}
=\vecv_{\ell,\alpha}A\,dw$.
We obtain the following lemma.
(See \S\ref{subsection;17.12.16.4} below.)
\begin{lem}
$A$ is a $C^{\infty}$-function on $U_{w,q}$,
and the Taylor series of $\delbar A$ at $\infty$ is $0$.
In particular,
the Taylor series of $A$ at $\infty$ is an element of 
 $M_r(\cnum[\![w_q^{-1}]\!])$.
Moreover, $A_{|\infty}$ is nilpotent.
\hfill\qed
\end{lem}

We use the variable $x_q=(1+|\lambda|^2)^{1/q}w_q$
instead of $w_q$.
We obtain
a filtered bundle 
$\nbigp_{\ast}\nbigvhat^{\lambda}_{\ell,\alpha}
:=\nbigp^{h_{V,\ell,\alpha}}_{\ast}
 V^{\lambda}_{\ell,\alpha|\inftyhat_{x,q}}$
on $(\inftyhat_{x,q},\infty)$.
(See \S\ref{subsection;20.7.21.1} for the ringed space
$\inftyhat_{x,q}$,
and Definition \ref{df;20.7.21.2}
for filtered bundle in this situation.)
By $\DDhat_{\ell,\alpha}^{\lambda}\vecv_{|\inftyhat_{x,q}}
=\vecv_{|\inftyhat_{x,q}}A_{|\inftyhat_{x,q}}$,
we obtain the formal $\lambda$-connection
$\DDhat_{\ell,\alpha}^{\lambda}$.
It is independent of the choice of $\vecv$.
By Proposition \ref{prop;17.12.16.3},
we obtain the following lemma.
\begin{lem}
$(\nbigp_{\ast}\nbigvhat^{\lambda}_{\ell,\alpha},
 \DDhat^{\lambda}_{\ell,\alpha})$
are good filtered $\lambda$-flat bundles.
\index{good filtered $\lambda$-flat bundle
$(\nbigp_{\ast}\nbigvhat^{\lambda}_{\ell,\alpha},
\DDhat^{\lambda}_{\ell,\alpha})$}
\hfill\qed
\end{lem}
 
Applying the procedure in \S\ref{subsection;17.10.7.2},
we obtain filtered bundles
$(\Psi^{\lambda}_q)^{\ast}\bigl(
 \nbigp_{\ast}\nbigvhat^{\lambda}_{\ell,\alpha},
 \DDhat^{\lambda}_{\ell,\alpha}\bigr)$
over locally free
$\nbigo_{\Hhat^{\lambda}_{\infty,q}}(\ast
H^{\lambda}_{\infty,q})$-modules
$(\Psi_q^{\lambda})^{\ast}
 \bigl(\nbigp\nbigvhat^{\lambda}_{\ell,\alpha},\DDhat^{\lambda}\bigr)$.
According to Proposition \ref{prop;20.7.20.131},
they are good filtered bundles.
\index{good filtered bundle $(\Psi^{\lambda}_q)^{\ast}\bigl(
 \nbigp_{\ast}\nbigvhat^{\lambda}_{\ell,\alpha},
 \DDhat^{\lambda}_{\ell,\alpha}\bigr)$}

\subsubsection{Refined statement}

\begin{thm}
\label{thm;17.10.7.10}
There exists an isomorphism
\begin{equation}
\label{eq;17.10.7.131}
 \nbigp^h E^{\lambda}_{|\Hhat^{\lambda}_{\infty,q}}
\simeq
 \bigoplus_{(\ell,\alpha)}
 \LLhat^{\lambda}_q(\ell,\alpha)
\otimes
 (\Psi_q^{\lambda})^{\ast}
 \bigl(\nbigp\nbigvhat^{\lambda}_{\ell,\alpha},
 \DDhat^{\lambda}_{\ell,\alpha}\bigr)
\end{equation}
such that
it induces an isomorphism of filtered bundles,
i.e., the following holds for any $t_1\in S^1_T$:
\[
 \nbigp^h_{\ast} E^{\lambda}_{|\Hhat^{\lambda}_{\infty,q}\langle t_1\rangle}
\simeq
 \bigoplus_{(\ell,\alpha)}
 \nbigp^{(0)}_{\ast}
 \LLhat^{\lambda}_q(\ell,\alpha)
 _{|\Hhat^{\lambda}_{\infty,q}\langle t_1\rangle}
\otimes
 (\Psi_q^{\lambda})^{\ast}
 \bigl(\nbigp_{\ast}\nbigvhat^{\lambda}_{\ell,\alpha},
 \DDhat^{\lambda}_{\ell,\alpha}\bigr)
 _{|\Hhat^{\lambda}_{\infty,q}\langle t_1\rangle}.
\]
\end{thm}

Theorem \ref{thm;17.10.5.130}
immediately follows from
Theorem \ref{thm;17.10.7.10},
Proposition \ref{prop;20.7.20.131} and
Lemma \ref{lem;21.7.7.3}.
We shall prove Theorem \ref{thm;17.10.7.10}
in \S\ref{subsection;17.10.14.1}--\ref{subsection;17.10.14.2}.

\subsection{Norm estimate}

Before going to the proof of Theorem \ref{thm;17.10.7.10},
we state the norm estimate for monopoles satisfying 
the GCK-condition,
which will be proved in
\S\ref{subsection;20.7.24.50}.

\begin{prop}
\label{prop;17.10.14.21}
The norm estimate holds for
 $(E^{\lambda},\delbar_{E^{\lambda}},h)$
with respect to 
the good filtered bundle $\nbigp^h_{\ast}E^{\lambda}$.
 (See Definition {\rm\ref{df;17.10.24.20}}
 for the norm estimates
 for good filtered bundles
 with a Hermitian metric.)
\end{prop}

\begin{rem}
According to Proposition {\rm\ref{prop;17.10.14.21}},
the GCK-condition for monopoles 
implies the norm estimate
 with respect to the associated filtered bundles
 $\nbigp^h_{\ast}E^{\lambda}$.
In particular, the metric $h$ is strongly adapted to
 $\nbigp^h_{\ast}E^{\lambda}$.
We shall study the converse
in Proposition {\rm\ref{prop;17.10.12.5}}.
\hfill\qed
\end{rem}

\subsection{Step 0}
\label{subsection;17.10.14.1}

We use the notation in \S\ref{subsection;16.9.28.1}.
Note that
$\nbigb_q^{\lambda\ast}
=\nbigb_q^{0\ast}$
as Riemannian manifolds.
We may assume that
$U_w=\{w\in\cnum\,|\,|w|>R\}\cup\{\infty\}$.

Let $(E,h,\nabla,\phi)$ be a monopole
as in \S\ref{subsection;20.8.8.100}.
There exist the decomposition (\ref{eq;21.8.21.40})
on $\nbigb^{0\ast}_{q}$,
holomorphic vector bundles with an automorphism
$(V_{\ell,\alpha},f_{\ell,\alpha})$
on $U_{w,q}^{\ast}$,
and an isomorphism (\ref{eq;21.8.21.41}).
As in \S\ref{subsection;21.8.21.100}
(see also \S\ref{subsection;17.10.5.50}),
we have the Hermitian metrics $h_{V,\ell,\alpha}$
of $V_{\ell,\alpha}$,
and the Hermitian metric
$h^{\shikaku}$ of $E$
as
\index{metric $h^{\shikaku}$}
\[
 h^{\shikaku}
 =\bigoplus_{\ell,\alpha}
 h_{\LL,q,\ell,\alpha}
 \otimes
 \Psi_q^{-1}(h_{V,\ell,\alpha}).
\]
Let $\nabla^{\shikaku}$
and $\phi^{\shikaku}$
denote the Chern connection
and the anti-Hermitian endomorphism
associated with
the mini-holomorphic bundle
$(E^0,\delbar_{E^0})$
with $h^{\shikaku}$
as in \S\ref{subsection;17.10.5.50}.
\index{connection $\nabla^{\shikaku}$}
\index{endomorphism $\phi^{\shikaku}$}
They are compatible with the decomposition
(\ref{eq;21.8.21.40}),
i.e.,
we have the decomposition
$\nabla^{\shikaku}=\bigoplus \nabla^{\shikaku}_{\ell,\alpha}$
and 
$\phi^{\shikaku}=\bigoplus \phi^{\shikaku}_{\ell,\alpha}$.

Let $s$ be the automorphism of $E$
such that (i) $h=h^{\shikaku}s$,
(ii) $s$ is self-adjoint with respect to both
 $h$ and $h^{\shikaku}$.
\begin{prop}
\label{prop;21.8.21.50}
 For any $m\in\seisuu_{\geq 0}$,
there exist positive constants
$A_i(m)$ $(i=1,2)$ such that 
the following holds for any 
$(\kappa_1,\ldots,\kappa_{m})
\in\{t,w,\wbar\}^{m}$.
\[
 \Bigl|
 \nabla^{\shikaku}_{\kappa_1}\circ\cdots\circ
 \nabla^{\shikaku}_{\kappa_m}(s-\id_E)
 \Bigr|_{h^{\shikaku}}
\leq A_1(m)\exp\bigl(-A_2(m)\cdot |w_q^q|\bigr),
\]
\[
 \Bigl|
 \nabla^{\shikaku}_{\kappa_1}\circ\cdots\circ
 \nabla^{\shikaku}_{\kappa_m}(\nabla-\nabla^{\shikaku})
 \Bigr|_{h^{\shikaku}}
\leq A_1(m)\exp\bigl(-A_2(m)\cdot |w_q^q|\bigr),
\]
\[
 \Bigl|
 \nabla^{\shikaku}_{\kappa_1}\circ\cdots\circ
 \nabla^{\shikaku}_{\kappa_m}(\phi-\phi^{\shikaku})
 \Bigr|_{h^{\shikaku}}
\leq A_1(m)\exp\bigl(-A_2(m)\cdot |w_q^q|\bigr),
\]
\[
 \Bigl|
 \nabla^{\shikaku}_{\kappa_1}\circ\cdots\circ
 \nabla^{\shikaku}_{\kappa_m}(F(\nabla)-F(\nabla^{\shikaku}))
 \Bigr|_{h^{\shikaku}}
\leq A_1(m)\exp\bigl(-A_2(m)\cdot |w_q^q|\bigr).
\]
\end{prop}
\pf
It follows from
Theorem \ref{thm;16.9.20.22},
Corollary \ref{cor;17.10.5.5},
Proposition \ref{prop;17.10.5.3}
and Corollary \ref{cor;17.10.5.4}.
\hfill\qed

\subsection{Step 1}

We define the differential operators
$\del^{\shikaku}_{E^{\lambda},\betabar_1}$
and 
$\del^{\shikaku}_{E^{\lambda},t_1}$
on $C^{\infty}(\nbigb^{\lambda\ast}_q,E)$
as follows:
\[
 \del^{\shikaku}_{E^{\lambda},\betabar_1}:=
 \frac{1}{1+|\lambda|^2}
 \left(
 \frac{\lambda\sqrt{-1}}{2}\nabla_t^{\shikaku}
+\nabla^{\shikaku}_{\wbar}
-\frac{\lambda}{2}\phi^{\shikaku}
 \right),
\]
\[
 \del^{\shikaku}_{E^{\lambda},t_1}:=
 \frac{1-|\lambda|^2}{1+|\lambda|^2}
 \nabla_t^{\shikaku}
-\frac{2\lambda\sqrt{-1}}{1+|\lambda|^2}
 \nabla^{\shikaku}_w
+\frac{2\lambdabar\sqrt{-1}}{1+|\lambda|^2}
 \nabla^{\shikaku}_{\wbar}
-\sqrt{-1}\phi^{\shikaku}.
\]
We set
$\nu_{\betabar_1}:=
 \del^{\shikaku}_{E^{\lambda},\betabar_1}
-\del_{E^{\lambda},\betabar_1}$
and 
$\nu_{t_1}:=
 \del^{\shikaku}_{E^{\lambda},t_1}
-\del_{E^{\lambda},t_1}$.
By the construction,
we obtain the following equalities:
\[
 \nu_{\betabar_1}=
\frac{1}{1+|\lambda|^2}
 \left(
 \frac{\lambda\sqrt{-1}}{2}
 (\nabla_t^{\shikaku}-\nabla_t)
+(\nabla^{\shikaku}_{\wbar}-\nabla_{\wbar})
-\frac{\lambda}{2}(\phi^{\shikaku}-\phi)
 \right),
\]
\[
 \nu_{t_1}=
\frac{1-|\lambda|^2}{1+|\lambda|^2}
 (\nabla_t^{\shikaku}-\nabla_t)
-\frac{2\lambda\sqrt{-1}}{1+|\lambda|^2}
 (\nabla^{\shikaku}_w-\nabla_w)
+\frac{2\lambdabar\sqrt{-1}}{1+|\lambda|^2}
 (\nabla^{\shikaku}_{\wbar}-\nabla_{\wbar})
-\sqrt{-1}(\phi^{\shikaku}-\phi).
\]
We obtain the following estimates from
Proposition \ref{prop;21.8.21.50}.
\begin{lem}
\label{lem;17.10.7.101}
For any $m\in\seisuu_{\geq 0}$,
there exists $\epsilon_m>0$
such that
\[
 \bigl(
 \del^{\shikaku}_{E^{\lambda},\kappa_1}
 \circ
 \cdots
 \circ
  \del^{\shikaku}_{E^{\lambda},\kappa_m}
 \bigr)
 \nu_{\betabar_1}
=O\Bigl(
 \exp\bigl(-\epsilon_{m}|w_q|^q\bigr)
 \Bigr),
\]
\[
 \bigl(
 \del^{\shikaku}_{E^{\lambda},\kappa_1}
 \circ
 \cdots
 \circ
  \del^{\shikaku}_{E^{\lambda},\kappa_m}
 \bigr)
 \nu_{t_1}
=O\Bigl(
 \exp\bigl(-\epsilon_{m}|w_q|^q\bigr)
 \Bigr),
\]
for any $(\kappa_1,\ldots,\kappa_m)\in\{t_1,\betabar_1\}^m$.
\hfill\qed
\end{lem}

\subsection{Step 2}
\label{subsection;21.8.13.23}

Let $\LL^{\lambda\cov}_q(\ell,\alpha)$
denote
the $\nbigo_{Y^{\lambda\cov}_q}(\ast H^{\lambda\cov}_{\infty,q})$-module
obtained as the pull back of
$\LL^{\lambda}_q(\ell,\alpha)$
by the projection $\varpi^{\lambda}_q$.
The restriction of $\LL^{\lambda\cov}_q(\ell,\alpha)$
to $Y^{\lambda\cov\ast}_q$
is denoted by
$\LL^{\lambda\cov\ast}_q(\ell,\alpha)$
which is equal to the pull back of
$\LL^{\lambda\ast}_q(\ell,\alpha)$ by $\varpi_q$.
It is equipped with the metric
$h_{\LL^{\cov},q,\ell,\alpha}$
obtained as the pull back of
$h_{\LL,q,\ell,\alpha}$.

There exists the induced decomposition
of $C^{\infty}$-vector bundles
on $\nbigb_q^{\lambda\cov\ast}$:
\begin{equation}
\label{eq;20.7.24.20}
 E^{\lambda\cov}:=
\varpi_q^{-1}(E^{\lambda})
=\bigoplus_{\ell,\alpha}
 \LL^{\lambda\cov\ast}_{q}(\ell,\alpha)
 \otimes
(\Psi^{\cov}_q)^{-1}(V_{\ell,\alpha}).
\end{equation}
The bundle $E^{\lambda\cov}$
is equipped with the operators
$\del^{\shikaku}_{E^{\lambda\cov},\betabar_1}$
and 
$\del^{\shikaku}_{E^{\lambda\cov},t_1}$
induced by 
$\del^{\shikaku}_{E^{\lambda},\betabar_1}$
and 
$\del^{\shikaku}_{E^{\lambda},t_1}$,
respectively.
They preserve the decomposition (\ref{eq;20.7.24.20}).

Let $b\in\real$.
Let $\nbigp^{h_{V,\ell,\alpha}}V^{\lambda}_{\ell,\alpha}$
be the filtered bundle on
$(U_{w,q},\infty)$
as in \S\ref{subsection;21.8.21.100}.
Let $\vecv_{\ell,\alpha}$ be 
holomorphic frames of
$\nbigp^{h_{V,\ell,\alpha}}_bV^{\lambda}_{\ell,\alpha}$
on $U_{w,q}$.
We obtain the $C^{\infty}$-frame 
\begin{equation}
\label{eq;17.10.5.310}
 \vecvtilde_{\ell,\alpha}=(\vtilde_{\ell,\alpha,i}\,|\,i=1,\ldots,\rank V_{\ell,\alpha}):=
 (\Psi^{\cov}_q)^{-1}(\vecv_{\ell,\alpha|U^{\ast}_{w,q}})
\end{equation}
of $(\Psi^{\cov}_q)^{-1}(V_{\ell,\alpha})$
on $\nbigb^{\lambda\cov\ast}_q$.

As explained in \S\ref{subsection;21.8.14.2},
there exist a neighbourhood
$\nbigu^{\lambda}_q$ of $H^{\lambda\cov}_{\infty,q}$
in $\nbigb^{\lambda\cov}_{q}$
and mini-holomorphic frames $u^{\lambda}_{q,\ell,\alpha}$
of $\LL^{\lambda\cov}_{q}(\ell,\alpha)$
on $\nbigu^{\lambda}_q$
for which (\ref{eq;21.8.21.30})
and (\ref{eq;21.8.21.31}) hold.

We set $\nbigu^{\lambda\ast}_q:=
\nbigu^{\lambda}_q\setminus H^{\lambda}_{\infty,q}$.
We obtain the $C^{\infty}$-frame
$u^{\lambda}_{q,\ell,\alpha}
 \otimes\vecvtilde_{\ell,\alpha}$
of the $C^{\infty}$-vector bundle
$\LL^{\lambda\cov\ast}_{q}(\ell,\alpha)
 \otimes
 (\Psi^{\cov}_q)^{-1}(V_{\ell,\alpha})$
on $\nbigu^{\lambda\ast}_q$.

We define the matrix valued $C^{\infty}$-function $A_{\ell,\alpha}$
on $U_{w,q}$
by $\DD^{\lambda(1,0)}\vecv_{\ell,\alpha}
 =\vecv_{\ell,\alpha}A_{\ell,\alpha}\,dw$.
By Proposition \ref{prop;21.8.12.5}
with the formulas (\ref{eq;21.8.12.30}) and (\ref{eq;21.8.12.31}),
we obtain
\begin{equation}
\label{eq;17.10.7.30}
 \del^{\shikaku}_{E^{\lambda\cov},\betabar_1}
 \bigl(
  u^{\lambda}_{q,\ell,\alpha}
 \otimes\vecvtilde_{\ell,\alpha}
 \bigr)=0, 
\end{equation}
\begin{equation}
\label{eq;21.8.3.11}
  \del^{\shikaku}_{E^{\lambda\cov},t_1}
 \bigl(
  u^{\lambda}_{q,\ell,\alpha}
 \otimes\vecvtilde_{\ell,\alpha}
 \bigr)
=
 \bigl(
  u^{\lambda}_{q,\ell,\alpha}
 \otimes\vecvtilde_{\ell,\alpha}
 \bigr)
\frac{1}{1+|\lambda|^2}
 \Bigl(
 -2\sqrt{-1}\bigl( 
 (\Psi_q^{\cov})^{-1}(A_{\ell,\alpha|U_{w,q}^{\ast}})
 \bigr)
 \Bigr).  
\end{equation}

Let
$H\bigl(h^{\shikaku},
 u^{\lambda}_{\ell,\alpha}\otimes\vecvtilde_{\ell,\alpha}\bigr)$
denote the Hermitian-matrix valued function 
on $\nbigb^{\lambda\cov\ast}_q$
determined by
\[
H\bigl(h^{\shikaku},
 u^{\lambda}_{\ell,\alpha}\otimes\vecvtilde_{\ell,\alpha}\bigr)_{i,j}:=
 (\varpi_q^{-1}h^{\shikaku})\Bigl(
 u^{\lambda}_{\ell,\alpha}\otimes \vtilde_{\ell,\alpha,i},
 u^{\lambda}_{\ell,\alpha}\otimes \vtilde_{\ell,\alpha,j}
 \Bigr).
\]

\begin{lem}
Let $P$ be any point of $H^{\lambda\cov}_{\infty,q}$.
Let $U_P$ be a relatively compact neighbourhood of $P$
in $\nbigb^{\lambda\cov}_q=(\varpi_q^{\lambda})^{-1}(\nbigb^{\lambda}_q)$.
Then, there exist $C_P\geq 1$ and $N_P>0$
such that 
\[
 C_P^{-1}|\beta_1|^{-N_P}
\leq
\bigl|
 H\bigl(h^{\shikaku},
 u^{\lambda}_{\ell,\alpha}\otimes\vecvtilde_{\ell,\alpha}\bigr)
\bigr|
 \leq
 C_P|\beta_1|^{N_P}.
\]
\end{lem}
\pf
It follows from Lemma \ref{lem;17.10.7.100},
Proposition \ref{prop;20.8.8.40},
Proposition \ref{prop;17.10.3.30}
and Proposition \ref{prop;17.10.25.110}.
\hfill\qed

\subsection{Step 3}
\label{subsection;21.8.12.55}

By using the frame 
$u^{\lambda}_{q,\ell,\alpha}
 \otimes\vecvtilde_{\ell,\alpha}$,
we extend the bundle
$\LL^{\lambda\cov\ast}_{q}(\ell,\alpha)
 \otimes
(\Psi^{\cov}_q)^{-1}(V_{\ell,\alpha})$
to a $C^{\infty}$-bundle
$\gbigp_b\bigl(
\LL^{\lambda\cov\ast}_{q}(\ell,\alpha),
(V_{\ell,\alpha},\DD^{\lambda}_{\ell,\alpha})
\bigr)$
on $\nbigb^{\lambda\cov}_q$.
Namely,
we construct 
$\gbigp_b\bigl(
\LL^{\lambda\cov\ast}_{q}(\ell,\alpha),
(V_{\ell,\alpha},\DD^{\lambda}_{\ell,\alpha})
\bigr)$
from the product bundle
$\cnum^{\rank V_{\ell,\alpha}}\times \nbigu^{\lambda}_q$
and 
$\LL^{\lambda\cov\ast}_{q}(\ell,\alpha)
 \otimes
 (\Psi^{\cov}_q)^{-1}(V_{\ell,\alpha})$
by the gluing isomorphism
\[
\Bigl(
\cnum^{\rank V_{\ell,\alpha}}\times \nbigu^{\lambda}_q
\Bigr)_{|\nbigu^{\lambda\ast}_q}
\simeq
 \Bigl(
 \LL^{\lambda\cov\ast}_{q}(\ell,\alpha)
 \otimes
 (\Psi^{\cov}_q)^{-1}(V_{\ell,\alpha})
 \Bigr)_{|\nbigu^{\lambda\ast}_q}
\]
induced by the identification  of the canonical frame
and $u^{\lambda}_{q,\ell,\alpha}
\otimes\vecvtilde_{\ell,\alpha}$.
We naturally regard
$u^{\lambda}_{q,\ell,\alpha}\otimes
 \vecvtilde_{\ell,\alpha}$
as a $C^{\infty}$-frame of 
$\gbigp_b\bigl(
 \LL^{\lambda\cov\ast}_{q}(\ell,\alpha),
 (V_{\ell,\alpha},\DD^{\lambda}_{\ell,\alpha})
 \bigr)$
on $\nbigu_q^{\lambda}$.
We define 
\[
 \gbigp_b\bigl(
 E^{\lambda\cov}
\bigr)
:=\bigoplus_{\ell,\alpha}
 \gbigp_b\bigl(
 \LL^{\lambda\cov\ast}_{q}(\ell,\alpha),
 (V_{\ell,\alpha},\DD^{\lambda}_{\ell,\alpha})
 \bigr).
\]
In this way,
$E^{\lambda\cov}$ extends to
a $C^{\infty}$-vector bundle
$\gbigp_b\bigl(
 E^{\lambda\cov}
\bigr)$
on $\nbigb^{\lambda\cov}_q$.

We can naturally regard
$\del^{\shikaku}_{E^{\lambda\cov},\betabar_1}$
and
$\del^{\shikaku}_{E^{\lambda\cov},t_1}$
as $C^{\infty}$-differential operators
on 
$\gbigp_b\bigl(
 E^{\lambda\cov}
\bigr)$
by the expressions (\ref{eq;17.10.7.30})
and (\ref{eq;21.8.3.11}).
We may assume that there exists
a $q$-th root $\tau_{1,q}$ of $\beta_{1}^{-1}$
on $\nbigu_q^{\lambda}$.
We set
$\del^{\shikaku}_{E^{\lambda},\taubar_{1,q}}
=-q\taubar_{1,q}^{-q-1}\del^{\shikaku}_{E^{\lambda},\betabar_1}$.
Because
\begin{equation}
\label{eq;17.10.7.102}
 \del^{\shikaku}_{E^{\lambda},\taubar_{1,q}}
 u^{\lambda}_{q,\ell,\alpha}
 \otimes
 \vecvtilde_{\ell,\alpha}
=0,
\end{equation}
it also induces a differential operator 
on
$\gbigp_b\bigl(
 E^{\lambda\cov}
\bigr)$.

\subsection{Step 4}
\label{subsection;21.8.12.52}

We have the $C^{\infty}$-frame
$\veca=\bigcup_{\ell,\alpha}
\bigl(
 u^{\lambda}_{q,\ell,\alpha}
 \otimes
 \vecv_{\ell,\alpha}
 \bigr)$
 of 
$\gbigp_b\bigl(
E^{\lambda\cov}
\bigr)$
on $\nbigu^{\lambda}_q$.
Let $\nbiga_{\betabar_1}$ and $\nbiga_{t_1}$ be the matrix valued functions
on $\nbigu^{\lambda\ast}_{q}$
determined by
\[
 \nu_{\betabar_1}
 \veca
=\veca\cdot \nbiga_{\betabar_1},
\quad\quad
 \nu_{t_1}
 \veca
=\veca\cdot \nbiga_{t_1}.
\]

\begin{lem}
\label{lem;17.10.7.110}
$\nbiga_{t_1}$ and $\nbiga_{\betabar_1}$
are $C^{\infty}$-functions on $\nbigu^{\lambda}_q$
whose Taylor series along $H^{\lambda\cov}_{\infty,q}$
are $0$.
\end{lem}
\pf
From Lemma \ref{lem;17.10.7.101}
and the expressions
(\ref{eq;17.10.7.30})
(\ref{eq;21.8.3.11})
and (\ref{eq;17.10.7.102}),
we obtain
\[
 \del_{t_1}^{\ell_1}
 \del_{\taubar_{1,p}}^{\ell_2}
 \nbiga_{t_1}
=O\Bigl(
 \exp\bigl(-\epsilon_{\ell_1,\ell_2}|w_q|^q\bigr)
 \Bigr),
\quad
  \del_{t_1}^{\ell_1}
 \del_{\taubar_{1,p}}^{\ell_2}
 \nbiga_{\betabar_1}
=O\Bigl(
 \exp\bigl(-\epsilon_{\ell_1,\ell_2}|w_q|^q\bigr)
 \Bigr)
\]
for any $(\ell_1,\ell_2)\in\seisuu^2_{\geq 0}$
for some $\epsilon_{\ell_1,\ell_2}>0$.
Then, 
by the elliptic regularity
(see Lemma \ref{lem;21.9.13.1} below),
for any $(\ell_1,\ell_2,\ell_3)\in\seisuu^3_{\geq 0}$
there exists $\epsilon'_{\ell_1,\ell_2,\ell_3}>0$
such that
\[
 \del_{t_1}^{\ell_1}
 \del_{\taubar_{1,p}}^{\ell_2}
 \del_{\tau_{1,p}}^{\ell_3}
 \nbiga_{t_1}
=O\Bigl(
 \exp\bigl(-\epsilon_{\ell_1,\ell_2,\ell_3}|w_q|^q\bigr)
 \Bigr),
\]
\[
  \del_{t_1}^{\ell_1}
 \del_{\taubar_{1,p}}^{\ell_2}
 \del_{\tau_{1,p}}^{\ell_3}
 \nbiga_{\betabar_1}
=O\Bigl(
 \exp\bigl(-\epsilon_{\ell_1,\ell_2,\ell_3}|w_q|^q\bigr)
 \Bigr).
\]
Then, the claim of the lemma follows.
\hfill\qed

\begin{lem}
\label{lem;21.9.13.1}
For $r_1,r_2>0$,
we set
$B(r_1,r_2)=\bigl\{
(\sigma,\zeta)\in\real\times\cnum\,\big|\,
|\sigma|<r_1,\,|\zeta|<r_2
 \bigr\}$.
Let $f$ be a $C^{\infty}$-function on $B(r_1,r_2)$.
Suppose that there exists
$\ell\in\seisuu_{\geq 0}$ and $A>0$
such that 
$|\del_{\sigma}^{j_1}\del_{\zetabar}^{j_2}f|<A$
for any $j_1+j_2\leq \ell+1$ on $B(r_1,r_2)$.
Then, for any $0<r_i'<r_i$ $(i=1,2)$,
there exist $A'$, depending only on
$r_i$, $r_i'$, $\ell$ and $A$
such that 
$|\del_{\sigma}^{k_1}\del_{\zeta}^{k_2}\del_{\zetabar}^{k_3}f|<A'$
on $B(r_1',r_2')$
for any $k_1+k_2+k_3\leq \ell$.
\end{lem}
\pf
Let $r_i'<r_i''<r_i$.
Let $\chi:\real\times\cnum\lrarr [0,1]$ be a $C^{\infty}$-function
such that $\chi(\sigma,\zeta)=1$ on $B(r_1'',r_2'')$
and that the support is contained in $B(r_1,r_2)$.
We set $\ftilde=\chi\cdot f$, which we can naturally regard
as a $C^{\infty}$-function on $\real\times\cnum$
with compact support.
There exists $A_1>0$,
depending only on $r_i$, $r_i''$ and $A$,
such that 
$|\del_{\sigma}^{j_1}\del_{\zetabar}^{j_2}\ftilde|<A_1$
for any $j_1+j_2\leq \ell+1$ on $\real\times\cnum$.
Choose any $p>3$.
There exists $A_2>0$,
depending only on $r_i$, $r_i''$ and $A$
such that the $L_{\ell+1}^{p}$-norm of
$\ftilde$ is smaller than $A_2$.
By the Sobolev embedding theorem
(for example, see \cite[Theorem 7.10]{Gilbarg-Trudinger}),
there exists $A_3>0$,
depending only on $r_i$, $r_i''$ and $A$
such that the $C^{\ell}$-norm 
$\ftilde$ is smaller than $A_3$.
\hfill\qed

\subsection{Step 5}

By Lemma \ref{lem;17.10.7.110},
the mini-holomorphic structure $\delbar_{E^{\lambda\cov}}$ of 
$E^{\lambda\cov}$
extends to a mini-holomorphic structure
on $\gbigp_b(E^{\lambda\cov})$.
We obtain the locally free
$\nbigo_{\nbigb^{\lambda\cov}_q}$-module
as the associated sheaf of mini-holomorphic sections,
which is also denoted by $\gbigp_b(E^{\lambda\cov})$.
We obtain the following locally free
$\nbigo_{\nbigb^{\lambda\cov}_q}(\ast H^{\lambda\cov}_{\infty,q})$-module:
\[
\nbigp'E^{\lambda\cov}:=
\gbigp_b(E^{\lambda\cov})
\otimes
\nbigo_{\nbigb^{\lambda\cov}_q}(\ast H^{\lambda\cov}_{\infty,q}).
\]
We can easily observe that it is independent of the choice of $b$
up to canonical isomorphisms.
The $\seisuu$-action on
$E^{\lambda}$
also naturally extends to
a $\seisuu$-action
on $\nbigp'E^{\lambda\cov}$.

By the construction,
there exists the naturally defined
$\seisuu$-equivariant morphism
\begin{equation}
\label{eq;17.10.7.120}
 \nbigp'E^{\lambda\cov}
\lrarr
 \nbigp^h E^{\lambda\cov}.
\end{equation}
\begin{lem}
\label{lem;17.10.7.121}
The morphism
{\rm(\ref{eq;17.10.7.120})}
is an isomorphism.
\end{lem}
\pf
Both sides are locally free
$\nbigo_{\nbigb^{\lambda\cov}_q}(\ast H^{\lambda\cov}_q)$-modules,
and the restriction of (\ref{eq;17.10.7.120}) to
$\nbigb^{\lambda\cov\ast}_q$ is an isomorphism.
Then, the claim follows.
\hfill\qed

\vspace{.1in}
As a result of Lemma \ref{lem;17.10.7.121},
the $\seisuu$-equivariant morphism
\begin{equation}
\label{eq;21.8.12.50}
 \nbigp'(E^{\lambda\cov})_{|\Hhat^{\lambda\cov}_{\infty,q}}
\lrarr
 \nbigp^h(E^{\lambda\cov})_{|\Hhat^{\lambda\cov}_{\infty,q}}
\end{equation}
is an isomorphism.
\begin{lem}
There exists a natural $\seisuu$-equivariant isomorphism
\begin{equation}
\label{eq;21.8.13.1}
 \gbigp_b(E^{\lambda\cov})_{|\Hhat^{\lambda\cov}_{\infty,q}}
\simeq
 \bigoplus_{(\ell,\alpha)}
\Bigl(
\bigl(
 \nbigo_{\Hhat^{\lambda\cov}_{\infty,q}}u^{\lambda}_{q,\ell,\alpha}
\bigr)
\otimes
 (\varpi_q)^{-1}
 (\Psi^{\lambda})^{\ast}
 \bigl(\nbigp_b\nbigvhat_{\ell,\alpha},\DDhat^{\lambda}_{\ell,\alpha}
 \bigr)
\Bigr).
\end{equation}
As a result,
there exists a natural $\seisuu$-equivariant isomorphism
\begin{equation}
\label{eq;21.8.12.51}
  \nbigp'(E^{\lambda\cov})_{|\Hhat^{\lambda\cov}_{\infty,q}}
\simeq
\bigoplus_{(\ell,\alpha)}
\Bigl(
\bigl(
\nbigo_{\Hhat^{\lambda\cov}_{\infty,q}}(\ast H^{\lambda\cov}_{\infty,q})
 u^{\lambda}_{q,\ell,\alpha}
\bigr)
\otimes
 (\varpi_q)^{-1}
 (\Psi^{\lambda})^{\ast}
 \bigl(\nbigp\nbigvhat_{\ell,\alpha},\DDhat^{\lambda}_{\ell,\alpha}
 \bigr)
\Bigr).
 \end{equation}
\end{lem}
\pf
Let $\nbigc^{\infty}_{\Hhat^{\lambda\cov}_{\infty,q}}$
denote the sheaf of algebras on $H^{\lambda\cov}_{\infty,q}$
as in \S\ref{subsection;21.8.13.2}.
By taking the Taylor series along
$H^{\lambda\cov}_{\infty,q}$,
we obtain the morphism of the sheaf of algebras
$\nbigc^{\infty}_{\nbigb^{\lambda}_q|H^{\lambda\cov}_{\infty,q}}
\lrarr
\nbigc^{\infty}_{\Hhat^{\lambda\cov}_{\infty,q}}$.
Hence, we obtain the morphism of the ringed spaces
$k^{\lambda\cov}_{q,C^{\infty}}:
(H^{\lambda\cov}_{\infty,q},
\nbigc^{\infty}_{\Hhat^{\lambda\cov}_{\infty,q}})
\lrarr
(\nbigb^{\lambda\cov}_{\infty,q},
\nbigc^{\infty}_{\nbigb^{\lambda\cov}_{\infty,q}})$.
For any $\seisuu$-equivariant
$C^{\infty}$-bundle $W$ on $\nbigb^{\lambda\cov}_{\infty,q}$,
we naturally obtain
the $\nbigc^{\infty}_{\Hhat^{\lambda\cov}_{\infty,q}}$-module
$(k^{\lambda\cov}_{q,C^{\infty}})^{\ast}W$.

By the construction of $\gbigp_b(E^{\lambda\cov})$,
there exists the following $C^{\infty}$-isomorphism:
\begin{equation}
\label{eq;21.8.13.3}
 \gbigp_b(E^{\lambda})
 \otimes
 \nbigc^{\infty}_{\nbigb^{\lambda}_q}
 \simeq  
 \bigoplus_{\ell,\alpha}
 \Bigl(
 \bigl(
  \nbigc^{\infty}_{\nbigb^{\lambda}_q}\cdot
 u^{\lambda}_{q,\ell,\alpha}
    \bigr)
 \otimes_{\nbigc^{\infty}_{\nbigb^{\lambda}_q}}
 (\Psi_q^{\lambda\cov})^{-1}
 \bigl(
 \nbigp_bV^{\lambda}_{\ell,\alpha}
 \bigr)
   \otimes
  \nbigc^{\infty}_{\nbigb^{\lambda}_q}
 \Bigr).
\end{equation}
Let $W_1$ and $W_2$  denote the left hand side and the right hand side
of (\ref{eq;21.8.13.3}), respectively.
Both $W_i$ are equipped with
the linear differential operator 
$\delbar_{W_i}:C^{\infty}(\nbigb^{\lambda\cov}_q,W_i)
\lrarr
C^{\infty}(\nbigb^{\lambda\cov}_q,W_i
\otimes\Omega^{0,1}_{\nbigb^{\lambda\cov}_q})$
satisfying the mini-complex Leibniz rule.
(See \S\ref{subsection;21.8.12.3} for the mini-complex Leibniz rule.)
Let $\del_{W_i,t_1}$ (resp. $\del_{W,\betabar_1}$)
denote the differential operators on $W_i$
induced by $\delbar_{W_i}$ and $\del_{t_1}$
(resp. $\del_{\betabar_1}$).
Under the induced isomorphism
of $\nbigc^{\infty}_{\Hhat^{\lambda\cov}_{\infty,q}}$-modules
\begin{equation}
\label{eq;21.8.13.11}
 (k^{\lambda\cov}_{q,C^{\infty}})^{\ast}W_1
 \simeq
 (k^{\lambda\cov}_{q,C^{\infty}})^{\ast}W_2,
\end{equation}
we have
$\del_{(k^{\lambda\cov}_{q,C^{\infty}})^{\ast}W_1,t_1}
=\del_{(k^{\lambda\cov}_{q,C^{\infty}})^{\ast}W_2,t_1}$
and
$\del_{(k^{\lambda\cov}_{q,C^{\infty}})^{\ast}W_1,\betabar_1}
=\del_{(k^{\lambda\cov}_{q,C^{\infty}})^{\ast}W_2,\betabar_1}$.
The intersection of the kernels of
$\del_{(k^{\lambda\cov}_{q,C^{\infty}})^{\ast}W_1,t_1}$
and 
$\del_{(k^{\lambda\cov}_{q,C^{\infty}})^{\ast}W_1,\betabar_1}$
is identified with the left hand side of
(\ref{eq;21.8.13.1}).
By using Lemma \ref{lem;21.8.12.60},
we obtain that 
the intersection of the kernel of 
$\del_{(k^{\lambda\cov}_{q,C^{\infty}})^{\ast}W_2,t_1}$
and 
$\del_{(k^{\lambda\cov}_{q,C^{\infty}})^{\ast}W_2,\betabar_1}$
is identified with the right hand side of
(\ref{eq;21.8.13.1}).
Hence, we obtain the isomorphism
(\ref{eq;21.8.13.1}) from (\ref{eq;21.8.13.11}).
\hfill\qed

\vspace{.1in}

We obtain the isomorphism (\ref{eq;17.10.7.131})
from (\ref{eq;21.8.12.50}) and (\ref{eq;21.8.12.51}).

\subsection{Step 6}
\label{subsection;17.10.14.2}

We take $t_1^{0}\in\real_{t_1}$ and $c\in\real$.
For $\ell\in\seisuu$,
we set 
\[
b(t^{0}_1,c,\ell):=c+\frac{\ell t^{0}_1}{T}.
\]
Let
$\vecv_{\ell,\alpha}=
 (v_{\ell,\alpha,i}\,|\,i=1,\ldots,\rank V_{\ell,\alpha})$
 be a holomorphic frame
of $\nbigp^{h_{V,\ell,\alpha}}_{b(t_1^0,c,\ell)}V^{\lambda}_{\ell,\alpha}$
such that it is compatible with the parabolic structure.
(See \S\ref{subsection;17.10.5.301}.)
Let $a(\ell,\alpha,i)$ denote the parabolic degree of
$v_{\ell,\alpha,i}$
(see (\ref{eq;20.8.8.50})).
Recall the following
(see Lemma \ref{lem;17.10.7.100}).
\begin{lem}
\label{lem;20.7.24.30}
 Let $H^{\ell,\alpha}$ be the Hermitian matrix valued function
 on $U_{w,q}^{\ast}$ determined by
\[
 H^{\ell,\alpha}_{i,j}:=
 h_{V,\ell,\alpha}\bigl(
 v_{\ell,\alpha,i},
 v_{\ell,\alpha,j}\bigr)|w_q|^{-a(\ell,\alpha,i)-a(\ell,\alpha,j)}.
\]
Then, there exist $C_1>1$ and $N>0$ such that
$C_1^{-1}
 \Bigl(
 \log |w_q|
 \Bigr)^{-N}
\leq
|H^{\ell,\alpha}|
 \leq 
 C_1\Bigl(
 \log|w_q|
 \Bigr)^N$
and 
 $C_1^{-1}
 \Bigl(
 \log |w_q|
 \Bigr)^{-N}
\leq
|(H^{\ell,\alpha})^{-1}|
 \leq 
 C_1\Bigl(
 \log|w_q|
 \Bigr)^N$
\hfill\qed
\end{lem}

As before,
we obtain the $C^{\infty}$-frame
$u^{\lambda}_{q,\ell,\alpha}
 \otimes
\vecvtilde_{\ell,\alpha}$
of
$\LL_q^{\lambda\cov\ast}(\ell,\alpha)
\otimes
(\Psi_q^{\cov})^{-1}(V_{\ell,\alpha})$
on $\nbigu^{\lambda}_q$.
Let $\pi^{\lambda\cov}_q:
\nbigb^{\lambda\cov}_q\lrarr\real$
be the map induced by
$(t_1,\beta_1)\longmapsto t_1$.
Let $\nbigb^{\lambda\cov}_q\langle t_1\rangle$
denote the fiber.
We set
$\nbigu^{\lambda}\langle t_1\rangle:=
\nbigb^{\lambda\cov}_q\langle t_1\rangle
\cap\nbigu^{\lambda}$
and 
$\nbigu^{\lambda\ast}\langle t_1\rangle:=
\nbigb^{\lambda\cov}_q\langle t_1\rangle
\cap\nbigu^{\lambda\ast}$.
The restriction of
$u^{\lambda}_{q,\ell,\alpha}
\otimes
\vecvtilde_{\ell,\alpha}$
to
$\nbigu^{\lambda\ast}_q\langle t_1\rangle$
induces a $C^{\infty}$-frame
\[
\bigl(
 s_{\ell,\alpha,i}\,\big|\,
 i=1,\ldots,\rank V_{\ell,\alpha}
\bigr)
\]
of the restriction of 
$ \LL_q^{\lambda\cov\ast}(\ell,\alpha) \otimes
(\Psi^{\cov}_q)^{-1}(V_{\ell,\alpha})$
to
$\nbigu^{\lambda\ast}_q\langle t_1\rangle$.
We obtain a $C^{\infty}$-frame
\[
\bigcup_{(\ell,\alpha)}
 \bigl(
 s_{\ell,\alpha,i}\,\big|\,
 i=1,\ldots,\rank V_{\ell,\alpha}
\bigr)
\]
of $E^{\lambda\cov}
 _{|\nbigu^{\lambda\ast}_q\langle t_1\rangle}$.
By Proposition \ref{prop;21.8.21.50},
$h$ and $h^{\shikaku}$ are mutually bounded.
By Lemma \ref{lem;20.7.24.30},
we obtain the following lemma.
\begin{lem}
Let $H$  be the Hermitian-matrix valued function
on $\nbigu_q^{\lambda\ast}\langle t_1\rangle$
determined by
\[
 H_{(\ell_1,\alpha_1,i_1),(\ell_2,\alpha_2,i_2)}:=
 h\bigl(
 s_{\ell_1,\alpha_1,i_1},
 s_{\ell_2,\alpha_2,i_2}
 \bigr)\cdot
 |\beta_{1,q}|^{-a(\ell_1,\alpha_1,i_1)-a(\ell_2,\alpha_2,i_2)
 +(\ell_1+\ell_2)t^0_1/T}.
\]
Then, there exist $C_1>1$ and $N_1>0$
such that
$C_1^{-1}(\log|\beta_1|)^{-N_1}<
 |H|<
 C_1(\log|\beta_1|)^{N_1}$.
\hfill\qed
\end{lem}

The tuple
$(s_{\ell,\alpha,i})$ 
induces a holomorphic frame of
\begin{equation}
 \label{eq;17.9.18.10}
 \bigoplus_{(\ell,\alpha)}
\Bigl(
 \nbigp_{-\ell t_1/T}^{(0)}\LL^{\lambda\cov}_q(\ell,\alpha)
 _{|\Hhat^{\lambda\cov}_{\infty,q}}
\otimes
 (\varpi_q^{\lambda})^{\ast}\Psi^{\ast}_{q}\bigl(
 \nbigp_{b(t_1^{0},c,\ell)}\nbigvhat^{\lambda}_{\ell,\alpha},
 \DDhatlambda_{\ell,\alpha}
 \bigr)
\Bigr)_{|\{t_1\}\times\inftyhat}.
\end{equation}
Note that
$\del_{E^{\lambda},\betabar_1}s_{\ell,\alpha,i}
=O\bigl(
 \exp(-\epsilon|\beta_1|)
 \bigr)$ for some $\epsilon>0$.
For any large $M>0$,
there exist holomorphic sections
$\sbar_{\ell,\alpha,i}$
of
$E^{\lambda\cov}_{|\nbigu^{\lambda\ast}_q\langle t_1\rangle}$
such that
\[
 |\sbar_{\ell,\alpha,i}
-s_{\ell,\alpha,i}|_h
=O\Bigl(
 |\beta_1|^{-M}
 \Bigr).
\]
If $M$ is sufficiently large,
$\bigcup_{(\ell,\alpha)}
\bigl(
\sbar_{\ell,\alpha,i}\,\big|\,
 i=1,\ldots,\rank V_{\ell,\alpha}
\bigr)$
is a holomorphic frame of
$\nbigp^h_{c}
 \Bigl(
 E^{\lambda\cov}_{|\nbigu^{\lambda\ast}_q\langle t_1\rangle}
 \Bigr)$.
Hence, we obtain that
the completion of 
$\nbigp^h_{c}
 \Bigl(
 E^{\lambda\cov}_{|\nbigu^{\lambda\ast}_q\langle t_1\rangle}
 \Bigr)$
at $\infty$
is equal to (\ref{eq;17.9.18.10}).
Therefore, the proof of Theorem \ref{thm;17.10.7.10}
is completed.
\hfill\qed

\subsection{Proof of Proposition \ref{prop;17.10.14.21}}
\label{subsection;20.7.24.50}

According to Lemma \ref{lem;17.10.15.2},
it is enough to study the case
where Condition \ref{condition;21.8.21.2}
is satisfied for $(E,h,\nabla,\phi)$,
i.e.,
there exists a decomposition (\ref{eq;21.8.21.40}).
We continue to use the notation in
\S\ref{subsection;17.10.14.1}--\S\ref{subsection;17.10.14.2}.

\begin{lem}
\label{lem;21.8.21.110}
 For each $(\ell,\alpha)$,
there exists a good filtered $\lambda$-flat bundle
$(\nbigp_{\ast}\nbigv^{\lambda}_{\ell,\alpha},\nabla^{\lambda})$
on $(U_{x,q},\infty)$
with an isomorphism
\begin{equation}
\label{eq;20.7.24.51} 
 (\nbigp_{\ast}\nbigv^{\lambda}_{\ell,\alpha},\nabla^{\lambda})
 _{|\inftyhat_{x,q}}
 \simeq
 (\nbigp_{\ast}\nbigvhat^{\lambda}_{\ell,\alpha},
  \DDhat^{\lambda}_{\ell,\alpha}).
\end{equation}
If $\lambda\neq 0$,
we may assume that the Stokes structure of
$(\nbigv_{\ell,\alpha},\DD_{\ell,\alpha}^{\lambda})$
is trivial.
\end{lem}
\pf
If $\lambda=0$,
we may choose $(\nbigp_{\ast}\nbigv^0_{\ell,\alpha},\nabla^0)$
as
$(\nbigp^{h_{V,\ell,\alpha}}_{\ast}V^0_{\ell,\alpha},f_{\ell,\alpha}dw)$.
Let us consider the case $\lambda\neq 0$.
For each $\gminia\in w_{q'}\cnum[w_{q'}]$,
we set $L(\gminia):=(\cnum(\!(w_{q'}^{-1})\!),\lambda d+d\gminia)$.
We obtain the $\cnum(\!(w_q)\!)$-module with the $\lambda$-connection
$\varphi_{q,q'\ast}L(\gminia)$.
According to the Hukuhara-Levelt-Turrittin theorem,
for an appropriate $q'$,
there exist a finite subset
$W(\ell,\alpha)\subset w_q'\cnum[w_q']$,
a tuple of
$\cnum(\!(w_q)\!)$-modules with regular singular $\lambda$-connection
$R_{\gminia}$ $(\gminia\in W(\ell,\alpha))$,
and an isomorphism
\[
(\nbigvhat_{\ell,\alpha},\DDhat^{\lambda}_{\ell,\alpha})
\simeq
\bigoplus_{\gminia\in W(\ell,\alpha)}
\varphi_{q,q'\ast}L(\gminia)
\otimes R_{\gminia}.
\]
Because $\varphi_{q',q\ast}L(\gminia)$ and
$R_{\gminia}$ are obtained as the formal completions of
meromorphic flat bundles on $(U_{w,q},\infty)$,
there exist
meromorphic $\lambda$-flat bundles
$(\nbigv_{\ell,\alpha},\DDlambda_{\ell,\alpha})$
on $(U_{w,q},\infty)$
such that
$(\nbigv_{\ell,\alpha},\DDlambda_{\ell,\alpha})_{|\inftyhat_{w,q}}
\simeq
(\nbigvhat_{\ell,\alpha},\DDhat^{\lambda}_{\ell,\alpha})$.
Hence, we obtain
good filtered $\lambda$-flat bundles
$(\nbigp_{\ast}\nbigv_{\ell,\alpha},\DDlambda_{\ell,\alpha})$
over $(\nbigv_{\ell,\alpha},\DDlambda_{\ell,\alpha})$
such that
$(\nbigp_{\ast}\nbigv_{\ell,\alpha},\DDlambda_{\ell,\alpha})_{|\inftyhat_{w,q}}
\simeq
(\nbigp_{\ast}\nbigvhat_{\ell,\alpha},\DDhat^{\lambda}_{\ell,\alpha})$.
\hfill\qed

\vspace{.1in}

Let $\nbigc^{\infty}_{U_{x,q}}$
denote the sheaf of $C^{\infty}$-functions
on $U_{x,q}$.
There exists an isomorphism
\[
 F_{\ell,\alpha}:
 \nbigp_{\ast}\nbigv^{\lambda}_{\ell,\alpha}
  \otimes
 \nbigc^{\infty}_{U_{x,q}}
 \simeq
 \nbigp^{h_{V,\ell,\alpha}}_{\ast}V^{\lambda}_{\ell,\alpha}
 \otimes
 \nbigc^{\infty}_{U_{x,q}}
\]
which induces (\ref{eq;20.7.24.51}).
By applying Proposition \ref{prop;20.7.27.2}
to $(V^{\lambda}_{\ell,\alpha},\DD^{\lambda}_{\ell,\alpha})$
with $h_{V,\ell,\alpha}$,
we obtain that
the norm estimate holds for
$(\nbigp_{\ast}\nbigv^{\lambda}_{\ell,\alpha},
\nabla^{\lambda})$
with $h_{0,\ell,\alpha}:=F_{\ell,\alpha}^{\ast}(h_{V,\ell,\alpha})$.

We obtain the following filtered bundle
\[
 \nbigp_{\ast}\nbige_0:= 
 \bigoplus_{\ell,\alpha}
 \nbigp_{\ast}\LL^{\lambda}_q(\ell,\alpha)
 \otimes
 \Psi_q^{\ast}\bigl(
  \nbigp_{\ast}\nbigv^{\lambda}_{\ell,\alpha},
  \nabla^{\lambda}
 \bigr).
\]
By Proposition \ref{prop;17.10.14.20}
and Lemma \ref{lem;21.8.23.3},
the norm estimate holds for
$\nbigp_{\ast}\nbige_0$
with
$h_{\nbige_0}=\bigoplus_{\ell,\alpha}
 h_{\LL,q,\ell,\alpha}
 \otimes
 \Psi_q^{-1}(h_{0,\ell,\alpha})$.

Let $\nbigc^{\infty}_{\nbigb^{\lambda}_q}$
denote the sheaf of $C^{\infty}$-functions
on $\nbigb^{\lambda}_q$.
By the proof of Theorem \ref{thm;17.10.7.10},
the morphisms $F_{\ell,\alpha}$ induce
the following isomorphism:
\[
 F:
  \nbige_0\otimes\nbigc^{\infty}_{\nbigb^{\lambda}_q}
  \simeq
  \nbigp^h(E^{\lambda})
  \otimes\nbigc^{\infty}_{\nbigb^{\lambda}_q}.
\]
Moreover, it induces
the following isomorphism
of the filtered bundles:
\[
\bigl(
\nbigp_{\ast}
  \nbige_0
  \bigr)_{|\Hhat^{\lambda}_{\infty,q}}
  \simeq
  \nbigp^h_{\ast}(E^{\lambda})_{|\Hhat^{\lambda}_{\infty,q}}.
\]
Note that
$h_{\nbige_0}
=F^{\ast}(h^{\shikaku})$,
and hence
$h_{\nbige_0}$
and 
$F^{\ast}(h)$ are mutually bounded.
Then, we can check that 
the norm estimate holds for
$(\nbigp^h_{\ast}E^{\lambda},h)$
by using the argument in
\S\ref{subsection;17.10.14.2}.
\hfill\qed

\section{Strong adaptedness and the GCK-condition}
\label{subsection;20.8.1.3}

\subsection{Statements}
\label{subsection;17.10.13.500}

Let $\nbigb_q^{\lambda}$
be a neighbourhood of
$H^{\lambda}_{\infty,q}$
in $Y^{\lambda}_q$
as in \S\ref{subsection;21.8.21.20}.
We set
$\nbigb_q^{\lambda\ast}:=
\nbigb_q^{\lambda}\setminus H^{\lambda}_{\infty,q}$.
\index{space $\nbigb_q^{\lambda}$}
\index{space $\nbigb_q^{\lambda\ast}$}
Let $\nbigp_{\ast}\nbige$ be a good filtered bundle
over a locally free 
$\nbigo_{\nbigb^{\lambda}_q}(\ast H^{\lambda}_{\infty,q})$-module
$\nbige$.
We obtain the mini-holomorphic bundle
$(E,\delbar_{E})$ on $\nbigb^{\lambda\ast}_q$
as the restriction of $\nbige$.
Let $h$ be a Hermitian metric of $E$
such that 
$(E,\delbar_{E},h)$ is a monopole.

\begin{prop}
\label{prop;17.10.12.5}
Suppose that $h$ is strongly adapted to $\nbigp_{\ast}\nbige$
in the sense of Definition {\rm\ref{df;20.8.8.3}}.
Then, the monopole satisfies the GCK-condition.
Moreover, $(E,\delbar_E,h)$ satisfies the norm estimate
with respect to $\nbigp_{\ast}\nbige$.
\end{prop}

We obtain the following corollary immediately.
\begin{cor}
Suppose that $(E,\delbar_E,h)$ satisfies the norm estimate
with respect to $\nbigp_{\ast}\nbige$.
Then, the monopole satisfies the GCK-condition.
\end{cor}
\pf
If $(E,\delbar_E,h)$ satisfies the norm estimate
with respect to $\nbigp_{\ast}\nbige$,
then $h$ is strongly adapted to $\nbigp_{\ast}\nbige$.
Hence, the claim follows from
Proposition \ref{prop;17.10.12.5}.
\hfill\qed

\vspace{.1in}

For the proof of Proposition \ref{prop;17.10.12.5},
we shall use the following auxiliary proposition,
which is also useful for the proof of
Theorem \ref{thm;17.9.30.20}.

\begin{prop}
\label{prop;17.10.12.6}
There exists a Hermitian metric $h_1$ of $E$
such that the following conditions are satisfied.
\begin{itemize}
 \item[(i)]
      The norm estimate holds for
      $(\nbigp_{\ast}\nbige,h_1)$.
\item[(ii)]
$|G(h_1)|_{h_1}=O(|w|^{-2-\epsilon})=O(|w_q|^{-q(2+\epsilon)})$ 
for some $\epsilon>0$,
and 
	   $|\delbar_EG(h_1)|_{h_1}\to 0$ as $|w_q|\to \infty$.
	   (See {\rm\S\ref{subsection;21.8.21.20}}
	   for $|\cdot|_{h_1}$.)
\item[(iii)]
$|F(h_1)|_{h_1}\to 0$ and
$|\nabla_{h_1} \phi_{h_1}|_{h_1}\to 0$ as $|w_q|\to \infty$.
\item[(iv)]
$|\phi_{h_1}|_{h_1}=O\bigl(\log|w_q|\bigr)$.
\item[(v)]
Let $\del_{E,h,\beta_1}$ denote the operator
induced by $\del_{E,\betabar_1}$ and $h$
as in {\rm\S\ref{subsection;17.10.5.120}}.
Then, we obtain
\[
\Bigl|
 \bigl[
 \del_{E,\betabar_1},\del_{E,h_1,\beta_1}
 \bigr]\Bigr|_{h_1}
=O\bigl(|w|^{-2}(\log|w|)^{-2}\bigr).
\]
\end{itemize}
If $\rank E=1$,
we may also assume that
$\del_{t_1}^{\ell_1}
 \del_{\beta_1}^{\ell_2}
 \del_{\betabar_1}^{\ell_3}
 G(h_1)=O(|w|^{-k})$
for any $(\ell_1,\ell_2,\ell_3)\in\seisuu_{\geq 0}^3$
and for any $k>0$.
\end{prop}

\subsection{Some estimates for tame harmonic bundles (1)}
\label{subsection;17.10.13.20}

Let $U_w$ be a neighbourhood of $\infty$ in $\proj^1_w$.
Let $U_{w,q}$ be the pull back of $U_w$
by the ramified covering $\proj^1_{w_q}\lrarr \proj^1_w$,
as in \S\ref{subsection;21.8.21.20}.
We set $U_{w,q}^{\ast}=U_{w,q}\setminus\{\infty\}$.
The pull back of $w$ by the ramified covering is also denoted by $w$,
i.e., $w=w_q^q$.

Let $(V,\delbar_V,\theta,h)$ be a tame harmonic bundle
on $U_{w,q}^{\ast}$.
For simplicity,
we assume that there exists a neighbourhood $U'_{w,q}$
of $U_{w,q}$ in $\proj^1_{w,q}$
such that $(V,\delbar_V,\theta,h)$
is the restriction of a harmonic bundle
on $U'_{w,q}\setminus\{\infty\}$.
We obtain the associated filtered Higgs bundle
$(\nbigp^h_{\ast}V,\theta)$ on $(U_{w,q},\infty)$.
Let $\delbar_V+\del_V$ denote the Chern connection
determined by $\delbar_V$ and $h$.
Let $F(h)$ denote the curvature of
$\delbar_V+\del_V$.
Let $\theta^{\dagger}$ denote the adjoint of $\theta$
with respect to $h$.
We have the expressions $\theta=f\,dw/w$
and $\theta^{\dagger}=f^{\dagger}d\wbar/\wbar$.

\begin{lem}
\mbox{{}}\label{lem;17.10.13.1}
\begin{itemize}
\item
 $f$ is bounded with respect to $h$.
\item
 $\delbar_V\del_Vf=O\bigl(|w|^{-2}(\log|w|)^{-2}dw\,d\wbar\bigr)$
 and 
 $\del_V\delbar_Vf^{\dagger}=O\bigl(|w|^{-2}(\log|w|)^{-2}dw\,d\wbar\bigr)$
 with respect to $h$.
 \item
 $\del_Vf=O(|w|^{-1}dw)$
 and $\delbar_Vf^{\dagger}=O(|w|^{-1}d\wbar)$
 with respect to $h$.
 \item
 $\del_{V,w}\del_{V,w}f=O(|w|^{-2})$
 and 
 $\del_{V,\wbar}\del_{V,\wbar}f^{\dagger}=O(|w|^{-2})$
 with respect to $h$,
where $\del_{V,w}$ and $\del_{V,\wbar}$
denote the inner product of $\del_w$ and $\del_{\wbar}$
with the Chern connection,
respectively.
\end{itemize}
\end{lem}
\pf
The first claim is proved in \cite{Simpson90}.
It is also proved
\[
 F(h)=O\bigl(|w|^{-2}(\log|w|)^{-2}dw\,d\wbar\bigr)
\]
in \cite{Simpson90}.
Because
$\delbar_V\del_Vf=
 (\delbar_V\del_V+\del_V\delbar_V)f
=[F(h),f]$, we obtain 
$\delbar_V\del_Vf=O\bigl(|w|^{-2}(\log|w|)^{-2}dw\,d\wbar\bigr)$.
By taking the adjoint,
we obtain the estimate for
$\del_V\delbar_Vf^{\dagger}$.

We have the holomorphic map
$g_q:\cnum_a\lrarr \cnum^{\ast}_{w_q}$ given by
$g_q(a)=e^{a/q}$,
i.e., $e^a=w$.
Let $\Utilde$ denote the pull back of
$U_{w,q}$ by $g_q$.
There exists $L\in\real$ such that
$\{a\,|\,\Re(a)>L\}\subset \Utilde$.
Let $(\Vtilde,\del_{\Vtilde},\thetatilde,\htilde)$
denote the pull back of
$(V,\delbar_V,\theta,h)$ by $g_q$.
Let $\ftilde$ and $\ftilde^{\dagger}$
denote the pull back of $f$ and $f^{\dagger}$ by $g_q$.
We have $\thetatilde=f\,da$.
For the adjoint $\thetatilde^{\dagger}$ of $\thetatilde$
with respect to $\htilde$, we have
$\thetatilde^{\dagger}=f^{\dagger}d\abar$.
Let $\del_{\Vtilde,a}$ and $\del_{\Vtilde,\overline{a}}$
denote the inner product of $\del_a$ and $\del_{\overline{a}}$
with the pull back
$\del_{\Vtilde}+\delbar_{\Vtilde}$ of $\del_V+\delbar_V$ by $g_q$.
Note that the curvature of the Chern connection
$\del_{\Vtilde}+\delbar_{\Vtilde}$ is
$O\bigl(\Re(a)^{-1}\bigr)$.

\vspace{.1in}
There exists $C_0>0$ such that
$|\ftilde|_h<C_0$ on $\Utilde$.
There exists $C_1>0$ such that the following holds for any
$a_0$ with $\Re(a_0)>L+10$.
\begin{itemize}
 \item
      $|\del_{\Vtilde,\overline{a}}\del_{\Vtilde,a}\ftilde|_{\htilde}\leq
      C_{1}|a_0|^{-2}$
      on the disc $\{|a-a_0|<2\}$.
\end{itemize}
By Proposition \ref{prop;21.9.13.11} below,
there exists $C_2>0$ such that the following holds
for any $a_0$ satisfying $\Re(a_0)>L+10$.
\begin{itemize}
 \item $|\del_{\Vtilde,a}\ftilde|_{\htilde}\leq C_2$
       on the disc $\{|a-a_0|<1\}$.
\end{itemize}
It implies that
$\del_Vf=O(|w|^{-1}dw)$
and 
$\delbar_Vf^{\dagger}=O(|w|^{-1}d\wbar)$
as $|w|\to\infty$.

We have the relation
$F(\htilde)+[\thetatilde,\thetatilde^{\dagger}]=0$,
i.e.,
$F(\htilde)_{a,\overline{a}}+[\ftilde,\ftilde^{\dagger}]=0$.
We obtain
$\del_{\Vtilde,a}F(\htilde)_{a,\overline{a}}=O(1)$
and 
$\del_{\Vtilde,\overline{a}}F(\htilde)_{a,\overline{a}}=O(1)$.
Because
$\del_{\Vtilde,\overline{a}}\del_{\Vtilde,a}\del_{\Vtilde,a}\ftilde
=F(\htilde)_{\overline{a},a}\del_{\Vtilde,a}\ftilde
+\del_{\Vtilde,a}(F(\htilde)_{\overline{a},a})\ftilde
=O(1)$,
we obtain $\del_{\Vtilde,a}\del_{\Vtilde,a}\ftilde=O(1)$
by Proposition \ref{prop;21.9.13.11} below.
We also obtain 
$\del_{\Vtilde,\overline{a}}
\del_{\Vtilde,\overline{a}}\ftilde^{\dagger}=O(1)$.
Thus, we obtain the fourth claim.
\hfill\qed

\subsubsection{Appendix}
\label{subsection;21.9.13.2}

For $r>0$, we set
$B(r)=\{\zeta\in\cnum\,|\,|\zeta|<r\}$.
Let $E$ be a vector bundle on $B(r)$
with a Hermitian metric $h_E$
with a unitary connection $\nabla$.
Let $F(\nabla)=F\cdot d\zeta\,d\zetabar$ denote
the curvature of $\nabla$.
Suppose that $|F|_h\leq C_1$ on $B(r)$
for a positive constant $C_1>0$.

\begin{prop}
\label{prop;21.9.13.11}
Let $s$ be any section of $E$ on $B(r)$
such that $|s|_h+|\nabla_{\zetabar}\nabla_{\zeta}s|_h\leq C_2$
for a positive constant $C_2>0$.
Let $0<r'<r$.
Then, there exists $C_3>0$,
depending only on $C_i$ $(i=1,2)$, $r$, $r'$ and $\rank(E)$,
such that
$|\nabla_{\zeta}s|_h+|\nabla_{\zetabar}s|_h\leq C_3$
on $B(r')$.
\end{prop}
\pf
Let $(x,y)$ denote the real coordinate system
determined by $\zeta=x+\sqrt{-1}y$.
Let $\vecv$ be a unitary $C^{\infty}$-frame of $E$
such that
$\nabla_x(\vecv_{|y=0})=0$ and $\nabla_y\vecv=0$.
Let $A_x$ be the matrix valued function determined by
$\nabla_x\vecv=\vecv A_x$.
Because $\nabla_yA_x$ represents $F$ with respect to
the frame $\vecv$,
and because $A_{x|y=0}=0$,
there exists $C_{10}$, depending only on
$r$, $C_1$ and $\rank(E)$
such that $|A_{x}|\leq C_{10}$.

Let $s'$ be a $C^{\infty}$-section of $E$
with compact support.
For any $p\geq 1$,
we set
$\|s'\|_{L^p(B(r))}=\left(
\int_{B(r)}|s'|_h^p\right)^{1/p}$.

\begin{lem}
\label{lem;21.9.13.5}
Suppose that
$\|s'\|_{L^{p}(B(r))}
+\|\nabla_{\zetabar}s'\|_{L^{p}(B(r))}\leq C_{11}$
for a constant $C_{11}>0$.
\begin{description}
 \item[The case $p=2$]
       For any $q>2$ and $0<r_1<r$,
       there exists  $C_{12}>0$,
       depending only on
       $q$, $r$, $r_1$, $C_{i}$ $(i=10,11)$ and $\rank(E)$
       such that
       $\|s'\|_{L^q(B(r_1))}\leq C_{12}$.
 \item[The case $p>2$]
For any $0<r_1<r$,
       there exists $C_{13}>0$,
       depending only on
       $p$, $r$, $r_1$, $C_{i}$ $(i=10,11)$ and $\rank(E)$
       such that
       $\sup_{B(r_1)}|s'|_h\leq C_{13}$.
\end{description}
\end{lem}
\pf
For the expression
$s'=\sum s'_iv_i$,
we obtain
\[
 \nabla_{\zetabar}(s')=
\sum_i \del_{\zetabar}(s'_i)v_i
+\sum_{i,j}\frac{1}{2}A_{x,ji}s'_iv_j.
\]
By the assumption,
there exists $C_{14}>0$,
depending only on 
$p$, $r$, $C_{i}$ $(i=10,11)$ and $\rank(E)$
such that 
$\|s'_i\|_{L^p(B(r))}
+\|\del_{\zetabar}s'_i\|_{L^p(B(r))}\leq C_{14}$.
Then, the claim follows from the elliptic regularity and
the Sobolev embedding theorem.
(For example, see \cite[Theorem 10.1]{ADN2} and
\cite[Theorem 7.10]{Gilbarg-Trudinger}.)
\hfill\qed

\vspace{.1in}
We set $r''=(r+r')/2$.
There exists a decreasing $C^{\infty}$-function
$\rho:\real_{\geq 0}\lrarr [0,1]$
such that
(i) $\rho(u)=1$ $(u\leq r')$,
(ii) $\rho(u)>0$ $(u<r'')$ and $\rho(u)=0$ $(u\geq r'')$,
(iii) $\rho^a$ is $C^{\infty}$ on $\real_{\geq 0}$ for any $a>0$.
We note the conditions imply that
the $C^{\infty}$-function
$(\del_u\rho^a)\rho^{-a/2}=2\del_u\rho^{a/2}$
on $u<r''$ extends to a $C^{\infty}$-function on $\real_{\geq 0}$
which are constantly $0$ on $u\geq r''$.
We obtain the $C^{\infty}$-function $\chi:\cnum\lrarr [0,1]$
by setting $\chi(\zeta)=\rho(|\zeta|)$.

We set $s'_a:=\chi^{a/2}\cdot \nabla_{\zeta}(s)$ for $a>0$.
We have
\begin{multline}
\label{eq;21.9.13.3}
 \bigl\|s'_a\bigr\|_{L^2(B(r))}^2
=\int_{B(r)}
 \chi^a h(\nabla_{\zeta}s,\nabla_{\zeta}s) \\
=-\int_{B(r)}\del_{\zeta}(\chi^a)\cdot h(s,\nabla_{\zeta}s)
-\int_{B(r)}\chi^a h(s,\nabla_{\zetabar}\nabla_{\zeta}s).
\end{multline}
There exists $C_{20}(a)>0$ depending only on
$a$, $r$, $r'$, $C_i$ $(i=1,2)$ and $\rank(E)$
such that
\[
 \left|
 \int_{B(r)}\chi^a h(s,\nabla_{\zetabar}\nabla_{\zeta}(s))
 \right|
 \leq C_{20}(a).
\]
There exists $C_{21}(a)>0$ depending only on
$a$, $r$, $r'$, $C_i$ $(i=1,2)$ and $\rank(E)$
such that
the first term in the right hand side of (\ref{eq;21.9.13.3})
is dominated as follows:
\begin{multline}
\label{eq;21.9.13.4}
 \left|
 \int_{B(r)}\del_{\zeta}(\chi^a)\cdot h(s,\nabla_{\zeta}s)
 \right|
\leq 
 \left(
 \int_{B(r)}
 \bigl(
  \chi^{-a/2}\del_{\zeta}(\chi^a)
  \bigr)^2|s|_h^2
  \right)^{1/2}
  \cdot
 \|s'_a\|_{L^2(B(r))}
 \\
\leq 
  C_{21}(a)\|s'_a\|_{L^2(B(r))}.
\end{multline}
We obtain
$\|s'_a\|_{L^2(B(r))}^2
\leq C_{20}(a)+C_{21}(a)\|s'_a\|_{L^2(B(r))}$.
Hence, there exists 
$C_{22}(a)>0$ depending only on
$a$, $r$, $r'$, $C_i$ $(i=1,2)$ and $\rank(E)$
such that
$\|s'_a\|_{L^2(B(r))}\leq C_{22}(a)$.

Note the following equality:
\begin{equation}
\label{eq;21.9.13.10}
 \nabla_{\zetabar}(s'_{a})
=(\chi^{-a/4}\del_{\zetabar}(\chi^{a/2}))s'_{a/2}
+\chi^{a/2}\nabla_{\zetabar}\nabla_{\zeta}(s).
\end{equation}
Hence, there exists $C_{23}(a)>0$ depending only on
$a$, $r$, $r'$, $C_i$ $(i=1,2)$ and $\rank(E)$
such that
$\|\nabla_{\zetabar}(s'_a)\|_{L^2(B(r))}\leq C_{23}(a)$.

Let $r''<r_1<r$.
By Lemma \ref{lem;21.9.13.5},
for any $p>2$,
there exists
$C_{24}(a,p)>0$ depending only on
$p$, $a$, $r$, $r'$, $r_1$, $C_i$ $(i=1,2)$ and $\rank(E)$
such that
$|s'_a|_{L^p(B(r_1))}\leq C_{24}(a,p)$.
By (\ref{eq;21.9.13.10}),
there exists
$C_{25}(a,p)>0$ depending only on
$p$, $a$, $r$, $r'$, $r_1$, $C_i$ $(i=1,2)$ and $\rank(E)$
such that
$|\nabla_{\zetabar}(s'_a)|_{L^p(B(r_1))}
\leq C_{25}(a,p)$.
By Lemma \ref{lem;21.9.13.5},
there exists
$C_{26}(a,p)>0$ depending only on
$p$, $a$, $r$, $r'$, $r_1$, $C_i$ $(i=1,2)$ and $\rank(E)$
such that
$\sup_{B(r'')}|s'_a|_h\leq C_{26}(a,p)$.
Thus, we obtain the estimate for
$|\nabla_{\zeta}(s)|_h$  on $B(r')$.
Similarly, we obtain the estimate for
$|\nabla_{\zetabar}(s)|_h$  on $B(r')$.
Thus, the proof of Proposition \ref{prop;21.9.13.11}
is completed.
\hfill\qed

\subsection{Some estimates for tame harmonic bundles (2)}

We continue to use the notation in \S\ref{subsection;17.10.13.20}.
For any $\lambda$
we obtain the holomorphic structure
$\delbar_V+\lambda\theta^{\dagger}$ on $V$.
Let $V^{\lambda}$ denote the holomorphic bundle
$(V,\delbar_V+\lambda\theta^{\dagger})$.
Let $\nabla$ denote the Chern connection
determined by $\delbar_V+\lambda\theta^{\dagger}$
and $h$,
i.e.,
$\nabla=\delbar_V+\lambda\theta^{\dagger}
+\del_V-\lambdabar\theta$.
Let $F(\nabla)$ be the curvature of $\nabla$.
We have $F(\nabla)=-(1+|\lambda|^2)[\theta,\theta^{\dagger}]$.
Because $V^{\lambda}$ with $h$ is acceptable,
we obtain the filtered bundle
$\nbigp^h_{\ast} V^{\lambda}$
on $(U_{w,q},\infty)$.

\begin{lem}
We obtain
$\nabla_w F(\nabla)_{w,\wbar}=O(|w|^{-3})$
and 
$\nabla_{\wbar} F(\nabla)_{w,\wbar}=O(|w|^{-3})$
with respect to $h$ as $|w_q|\to\infty$.
\end{lem}
\pf
We have
$F(\nabla)_{w,\wbar}=-(1+|\lambda|^2)|w|^{-2}\cdot[f,f^{\dagger}]$.
We also have
$\nabla_{w}(f)=\del_{V,w}f-\lambdabar[w^{-1}f,f]
=O(|w|^{-1})$
and 
$\nabla_w(f^{\dagger})=
 \del_{V,w}f^{\dagger}
-\lambdabar[w^{-1}f,f^{\dagger}]
 =O(|w|^{-1})$
by Lemma \ref{lem;17.10.13.1}. 
Then, we can easily deduce the claim of the lemma.
\hfill\qed

\vspace{.1in}

Let $s$ be a holomorphic section of $\nbigp^h_{<0}V^{\lambda}$.
Note that $\nabla_{\wbar}s=0$.

\begin{lem}
\label{lem;17.10.13.10}
There exists $\epsilon>0$
such that
$\nabla_{w}s=O(|w|^{-1-\epsilon})$
and 
$\nabla_{w}^2s=O(|w|^{-2-\epsilon})$
 with respect to $h$
as $|w_q|\to\infty$.
\end{lem}
\pf
We consider the map $g_q:\cnum_a\lrarr \cnum_{w_q}^{\ast}$
given by $g_q(a)=e^{a/q}$
as in the proof of Lemma \ref{lem;17.10.13.1}.
Let $(\Vtilde,\htilde,\nablatilde)$ denote the pull back of
$(V,h,\nabla)$ by $g_q$.
Let $\stilde$ denote the pull back of $s$ by $g_q$.
Because
$\nablatilde_{\abar}\nablatilde_a \stilde
=F(\nablatilde)_{\abar,a}\stilde$,
there exists $\epsilon_1>0$
such that the following holds
for any $a_0$ with $\Re(a_0)>L+10$.
\begin{itemize}
 \item On the disc $\{|a-a_0|<3\}$,
we obtain $|\stilde|_{\htilde}=O(e^{-\epsilon_1|a_0|})$,
       $|\nablatilde_{\abar}\nablatilde_a\stilde|
       =O(e^{-\epsilon_1|a_0|})$.
\end{itemize}
We obtain $\nabla_as=O(e^{-\epsilon_1|a_0|})$
on $\{|a-a_0|<2\}$ by Proposition \ref{prop;21.9.13.11}
for any $a_0$ with $\Re(a_0)>L+10$.
Moreover, because 
$\nabla_{\abar}\nabla_a\nabla_as
=2F(\nabla)_{\abar,a}\nabla_as
+\nabla_{a}(F(\nabla)_{\abar,a})s$,
there exists $\epsilon_2>0$
with the following property:
\begin{itemize}
\item
On discs $\{|a-a_0|<2\}$,
we obtain
$|\nabla_as|=O(e^{-\epsilon_2|a_0|})$
and 
$|\nabla_{\abar}\nabla_a\nabla_as|=O(e^{-\epsilon_2|a_0|})$.
\end{itemize}
We obtain
$\nabla_a\nabla_as=O(e^{-\epsilon_2|a_0|})$
on $\{|a-a_0|<1\}$
by Proposition \ref{prop;21.9.13.11},
which implies the claim of the lemma.
\hfill\qed

\subsection{$\lambda$-connections}

Let $(\nbigp_{\ast}\nbigv,\DDlambda)$
be a good filtered $\lambda$-flat bundle
on $(U_{w,q},\infty)$
whose Poincar\'e rank is strictly smaller than $q$.
For simplicity,
we assume that the Stokes structure is trivial,
i.e.,
for an appropriate covering
$\varphi_{q,p}:U_{w,p}\lrarr U_{w,q}$,
there exists a decomposition
\[
\varphi_{q,p}^{\ast}(\nbigp_{\ast}\nbigv,\DDlambda)
=\bigoplus_{\gminia\in w_p\cnum[w_p]}
 \varphi_{q,p}^{\ast}(\nbigp_{\ast}\nbigv_{\gminia},\DDlambda_{\gminia}),
\]
such that $\DDlambda_{\gminia}-d\gminia\id$
are logarithmic with respect to
$\nbigp_{\ast}\nbigv_{\gminia}$.

Let $(V,\DDlambda)$ be a $\lambda$-flat bundle
obtained as the restriction
of $(\nbigv,\DDlambda)$ to $U_{w,q}^{\ast}$.
We obtain the decomposition $\DDlambda=d''_V+d'_V$
into the $(0,1)$-part and the $(1,0)$-part.
For any Hermitian metric $h$ of $V$,
we obtain the operators
$\delta'_h$, $\delta''_h$
and $\DD^{\lambda\star}_h:=\delta'_h-\delta''_h$
as in \S\ref{subsection;17.10.12.1}.

\begin{lem}
\label{lem;17.10.13.100}
There exists a Hermitian metric $h$
satisfying the following conditions,
where we consider the norms with respect to $h$.
 \begin{itemize}
 \item $h$ satisfies the norm estimate for
       $(\nbigp_{\ast}\nbigv,\DDlambda)$.
\item
$|\Lambda\bigl[\DD^{\lambda},\DD^{\lambda\star}_h\bigr]|
=O(|w_q|^{-2q-2\epsilon})$,
where $\Lambda$ is determined by
$\Lambda(dw\,d\wbar)=-2\sqrt{-1}$.
 \item
$\bigl|[d''_{\wbar},\delta'_{h,w}]\bigr|$,
$\bigl|[d''_{\wbar},\delta''_{h,\wbar}]\bigr|$,
$\bigl|[d'_{w},\delta'_{h,w}]\bigr|$
and
$\bigl|[d'_{w},\delta''_{h,\wbar}]\bigr|$
go to $0$
as $w_q\to\infty$.
We also obtain
\begin{equation}
 \label{eq;17.10.24.41}
\bigl|
[d''_{\wbar},\delta'_{h,w}]\bigr|=
 O\bigl(|w|^{-2}(\log|w|)^{-2}\bigr).
\end{equation}
 \item $|d'_w-\lambda \delta'_w|$ and
       $|\delta''_{\wbar}-\lambdabar d'_{\wbar}|$
       go to $0$ as $w_q\to\infty$. 
 \item
$\Bigl|d''_{\wbar}\Lambda\bigl[\DD^{\lambda},\DD^{\lambda\star}_h\bigr]\Bigr|$,
$\Bigl|d'_{w}\Lambda\bigl[\DD^{\lambda},\DD^{\lambda\star}_h\bigr]\Bigr|$,
$\Bigl|\delta'_{w}\Lambda\bigl[\DD^{\lambda},\DD^{\lambda\star}_h\bigr]\Bigr|$
and 
$\Bigl|
\delta''_{\wbar}\Lambda\bigl[\DD^{\lambda},\DD^{\lambda\star}_h\bigr]\Bigr|$
go to $0$
as $w_q\to\infty$.
 \end{itemize}
 If $\rank V=1$,
we can impose
$\Lambda[\DDlambda,\DD^{\lambda\star}_h]=0$.
\end{lem}
\pf
First, let us study the case $\rank V=1$.
Let $v$ be a frame of
$\nbigp_a\nbigv$.
Let $h$ be the Hermitian metric of $V$
determined by
$h(v,v)=|w_q|^{2a}$.
Then,
the Chern connection
associated with $(V,d''_V)$ with $h$
is flat,
and hence it is a harmonic metric,
i.e.,
$[\DDlambda,\DD^{\lambda\star}]=0$.
(See \cite[Lemma 2.31]{Mochizuki-KHII}.)
The norm estimate is satisfied for
$(\nbigp_{\ast}\nbigv,\DDlambda)$
with $h$.
We obtain the meromorphic function $\gminia(w_q)$
determined by
$\DDlambda v=v\,d\gminia$.
Because the Poincar\'e rank of $(\nbigv,\DDlambda)$
is strictly smaller than $q$,
we obtain $|\gminia|=O(|w_q|^{q-1})$ as $|w_q|\to\infty$.
Note that $\del_w(\gminia)=O(|w_q|^{-1})$.
We obtain
\[
d''_{\wbar}v=0,\quad
d'_wv=\del_w(\gminia)v,\quad
\delta'_{w}v=(aq^{-1}w^{-1})v,\quad
\delta''_{\wbar}v=
\Bigl(
 -\overline{\del_w\gminia}
 -\lambdabar a \wbar^{-1}
\Bigr)v.
\]
Then, the claims are easily proved.

To study the general case,
we use the constructions
in \cite{mochi2, Mochizuki-KHII, Mochizuki-wild, Simpson90}.
We explain only an indication.
By pulling back via a ramified covering,
we may assume that there exists
the following decomposition from the beginning:
\[
 (\nbigp_{\ast}\nbigv,\DD^{\lambda})
=\bigoplus_{\gminia\in w_q\cnum[w_q]}
  (\nbigp_{\ast}\nbigv_{\gminia},
 \DD^{\lambda}_{\gminia}).
\]
Here,
$\DD^{\lambda}_{\gminia} -d\gminia\id$
are logarithmic with respect to
$\nbigp_{\ast}\nbigv_{\gminia}$.
It is enough to construct a Hermitian metric
for $(\nbigp_{\ast}\nbigv_{\gminia},
 \DD^{\lambda}_{\gminia})$
for each $\gminia$.
By considering the $\lambda$-connection
associated with the harmonic bundle
$(\nbigo_{U_{w,q}^{\ast}}\cdot e,(1+|\lambda|^2)^{-1}d\gminia)$
with $h(e,e)=1$,
it is enough to consider the case $\gminia=0$,
i.e.,
$\DDlambda$ is logarithmic
with respect to $\nbigp_{\ast}\nbigv$.

\vspace{.1in}

By using the model harmonic bundles
in \cite[\S6]{mochi2},
we can construct an endomorphism
$\gbigf$ of $\nbigv$
such that the following holds
(see \cite{mochi2}):
\begin{itemize}
\item
$\gbigf\nbigp_a\nbigv\subset \nbigp_{<a}\nbigv$ 
for any $a\in\real$.
\item
We set
$\DDlambda_1:=\DDlambda-\gbigf\,dw/w$.
Then, there exists a harmonic metric $h$
of $(V,\DDlambda_1)$
such that
$h$ satisfies the norm estimate
for $(\nbigp_{\ast}\nbigv,\DDlambda_1)$.
\end{itemize}
We obtain the operators
$\delta'_{h}$,
$\delta''_{1,h}$
and $\DD^{\lambda\star}_{1,h}$
from $\DDlambda_1$ and $h$.
Because $h$ is a harmonic metric for
$(V,\DDlambda_1)$,
$[\DDlambda_1,\DD^{\lambda\star}_{1,h}]=0$ holds.

There exists the decomposition
$\DD^{\lambda}_1=d''_V+d'_{1,V}$
into the $(0,1)$-part and the $(1,0)$-part.
There exist the operators
$\delbar_{1,V}$,
$\del_{1,V}$,
$\theta_1$
and $\theta_1^{\dagger}$
such that
$d''_V=\delbar_{1,V}+\lambda\theta_1^{\dagger}$,
$\delta''_{1,h}=\lambdabar\delbar_{1,V}-\theta_1^{\dagger}$,
$d'_{1,V}=\lambda \del_{1,V}+\theta_1$,
and $\delta'_{h}=\del_{1,V}-\lambdabar\theta_1$.
We have the expressions $\theta_1=g_1\,dw$
and $\theta_1^{\dagger}=g_1^{\dagger}d\wbar$.
Let $\del_{1,V,\wbar}$ denote the inner product of
$\delbar_{1,V}$ and $\del_{\wbar}$.
Similarly, $\del_{1,V,w}$ denote the inner product of
$\del_{1,V}$ and $\del_w$.
Then, we obtain
$d''_{\wbar}=\del_{1,V,\wbar}+\lambda g_1^{\dagger}$,
$d'_{1,w}=\lambda\del_{1,V,w}+g_1$,
$\delta''_{1,h,\wbar}=\lambdabar\del_{1,V,\wbar}-g_1^{\dagger}$
and 
$\delta'_{h,w}=\del_{1,V,w}-\lambdabar g_1$.
Because $h$ is a harmonic metric of $(V,\DDlambda)$,
we obtain
$\del_{1,V,\wbar}g_1=0$,
$\del_{1,V,w}g_1^{\dagger}=0$,
and
$[\del_{1,V,\wbar},\del_{1,V,w}]
+[g_1^{\dagger},g_1]=0$.
We obtain the following equalities:
\[
 [d''_{\wbar},\delta'_{h,w}]=
 -(1+|\lambda|^2)[g_1^{\dagger},g_1],
\quad\quad
 [d''_{\wbar},\delta''_{1,h,\wbar}]=
 -(1+|\lambda|^2)\del_{1,V,\wbar}g_1^{\dagger},
\]
\[
 [d'_{1,w},\delta'_{h,w}]=-(1+|\lambda|^2)\del_{1,V,w}g_1,
\quad\quad
 [d'_{1,w},\delta''_{1,h,\wbar}]=-(1+|\lambda|^2)[g_1,g_1^{\dagger}].
\]
By Lemma \ref{lem;17.10.13.1},
they go to $0$ as $|w_q|\to\infty$.
By the estimate for the curvature of wild harmonic bundles
(see \cite[Corollary 7.2.10]{Mochizuki-wild}),
we also obtain
$[d''_{\wbar},\delta'_{h,w}]=O\bigl(|w|^{-2}(\log|w|)^{-2}\bigr)$.

Note that
$[d'_{w},\delta'_{h,w}]
=[d'_{1,w},\delta'_{h,w}]
-\delta'_{h,w}(\gbigf w^{-1})$.
By Lemma \ref{lem;17.10.13.10},
we obtain
$\delta'_{h,w}(\gbigf w^{-1})
=O(|w|^{-2-\epsilon_1})$ for some $\epsilon_1>0$.
Hence, 
$[d'_{w},\delta'_{h,w}]$ goes to $0$.
By taking the adjoint,
we also obtain that
$[d''_{\wbar},\delta''_{h,\wbar}]$ goes to $0$.

Note that
$[d'_w,\delta''_{h,\wbar}]
=[d'_{1,w},\delta''_{1,h,\wbar}]
-d_{1,w}'(\gbigf^{\dagger}\wbar^{-1})
-\delta_{1,h,\wbar}''(\gbigf w^{-1})
-|w|^{-2}[\gbigf,\gbigf^{\dagger}]$.
We have 
$[\gbigf,\gbigf^{\dagger}]=O(|w|^{-2\epsilon})$.
We have
\[
 \delta_{1,h,\wbar}''(\gbigf w^{-1})
=\lambdabar d''_{\wbar}(\gbigf)w^{-1}
-(1+|\lambda|^2)[g_1^{\dagger},\gbigf]w^{-1}
=O(|w|^{-2-\epsilon}).
\]
By taking the adjoint,
we also obtain
$d_{1,w}'(\gbigf^{\dagger}\wbar^{-1})
=O(|w|^{-2-\epsilon})$.
Therefore,
$[d'_w,\delta''_{1,\wbar}]$ goes to $0$
as $|w_q|\to\infty$.
As the adjoint,
$[d''_{\wbar},\delta'_{1,w}]$ goes to $0$
as $|w_q|\to\infty$.

We note that
$d'_{1,w}-\lambda\delta'_{h,w}=(1+|\lambda|^2)g_1=O(|w|^{-1})$
and
$\delta''_{1,h,\wbar}-\lambdabar d''_{\wbar}
=-(1+|\lambda|^2)g_1^{\dagger}=O(|w|^{-1})$.
Because
$d'_w-\lambda\delta'_{h,w}
=(d'_{1,w}-\lambda\delta'_{h,w})+\gbigf w^{-1}$
and
$\delta''_{\wbar}-\lambdabar d'_{\wbar}
=(\delta''_{1,\wbar}-\lambdabar d'_{\wbar})
-\gbigf^{\dagger}_h \wbar^{-1}$,
we obtain the desired estimate for
$d'_w-\lambda\delta'_{h,w}$
and 
$\delta''_{\wbar}-\lambdabar d'_{\wbar}$.

Note that
$\frac{\sqrt{-1}}{2}
 \Lambda\bigl[\DDlambda,\DD^{\lambda\star}_h\bigr]
=[\delta'_{h,w},d''_{\wbar}]
+[d'_w,\delta''_{h,\wbar}]$
and 
$0=\frac{\sqrt{-1}}{2}
 \Lambda\bigl[\DDlambda_1,\DD^{\lambda\star}_{1,h}\bigr]
=[\delta'_{h,w},d''_{\wbar}]
+[d'_{1,w},\delta''_{1,h,\wbar}]$.
Hence, we obtain
\begin{multline}
 \label{eq;17.10.13.30}
 \frac{\sqrt{-1}}{2}
 \Lambda\bigl[\DDlambda,\DD^{\lambda\star}_h\bigr]
=[d'_w,\delta''_{h,\wbar}]
 -[d'_{1,w},\delta''_{1,h,\wbar}]
 \\
=-d'_{1,w}\bigl(\gbigf^{\dagger}\wbar^{-1}\bigr)
-\delta_{1,h,\wbar}''(\gbigf w^{-1})
-\bigl[
 \gbigf,\gbigf^{\dagger}
 \bigr]\cdot|w|^{-2}.
\end{multline}
By the previous consideration,
we obtain
$\frac{\sqrt{-1}}{2}
 \Lambda\bigl[\DDlambda,\DD^{\lambda\star}_h\bigr]
 =O(|w|^{-2-\epsilon})$ for some $\epsilon>0$.

We obtain
\[
\delta'_{h,w}\delta''_{1,h,\wbar}(\gbigf)
=\delta''_{1,h,\wbar}\delta'_{h,w}\gbigf
=\lambdabar (d''_{\wbar}\delta'_{h,w}-\delta'_{h,w} d''_{\wbar})\gbigf
-(1+|\lambda|^2)g_1^{\dagger}\delta_{h,w}'(\gbigf)
=O(|w|^{-2-\epsilon}).
\]
We obtain
\[
 d'_w\delta''_{1,h,\wbar}(\gbigf)
=\bigl(\lambda \delta'_{h,w}+(1+|\lambda|^2)g_1+\gbigf w^{-1}\bigr)
 \delta''_{1,h,\wbar}(\gbigf)
=O(|w|^{-2-\epsilon}).
\]
We have
\begin{multline}
 d''_{\wbar}\delta''_{1,h,\wbar}(\gbigf)
=d''_{\wbar}\bigl(
 \lambdabar d_{\wbar}''-(1+|\lambda|^2g_1^{\dagger})\gbigf
 \bigr) \\
=-(1+|\lambda|^2)d''_{\wbar}(g_1^{\dagger})\gbigf
=-(1+|\lambda|^2)\del_{V,\wbar}(g_1^{\dagger})\cdot \gbigf
=O(|w|^{-2-\epsilon}).
\end{multline}
We obtain
\[
 \delta''_{h,\wbar}
 \delta''_{1,h,\wbar}(\gbigf)
=\Bigl(
 \lambdabar d''_{\wbar}
-(1+|\lambda|^2)g_1^{\dagger}
-\gbigf \wbar^{-1}
 \Bigr)
 \delta''_{1,h,\wbar}(\gbigf)
=O(|w|^{-2-2\epsilon}).
\]
By considering the adjoint,
we obtain the estimate for
$d''_{\wbar}d'_{1,w}(\gbigf^{\dagger})$,
$d'_{w}d'_{1,w}(\gbigf^{\dagger})$,
$\delta''_{h,\wbar}d'_{1,w}(\gbigf^{\dagger})$
and $\delta'_{h,w}d'_{1,w}(\gbigf^{\dagger})$.
We have
$d''_{\wbar}[\gbigf,\gbigf^{\dagger}]
=[\gbigf,d''_{\wbar}\gbigf^{\dagger}]=O(|w|^{-1-\epsilon})$.
We also have
\[
 \delta''_{h,\wbar}[\gbigf,\gbigf^{\dagger}]
=\bigl(
 \lambdabar d''_{\wbar}-(1+|\lambda|^2)g_1^{\dagger}
 -\gbigf^{\dagger}\wbar^{-1}
 \bigr)[\gbigf,\gbigf^{\dagger}]
=O(|w|^{-1-\epsilon}).
\]
By taking the adjoint,
we obtain the estimate for 
 $\delta'_{h,w}[\gbigf,\gbigf^{\dagger}]$
and 
 $d'_{w}[\gbigf,\gbigf^{\dagger}]$.
Thus, we obtain the claim of the lemma.
\hfill\qed

\subsection{Proof of Proposition \ref{prop;17.10.12.6}}

It is enough to consider the case
where $\nbige_{|\Hhat^{\lambda}_{\infty,q}}$
is unramified modulo level $<1$.
There exist 
a finite subset $S\subset\seisuu\times\cnum^{\ast}$,
a tuple of
 $\cnum(\!(w_q^{-1})\!)$-modules 
with $\lamda$-connection
$(\nbigvhat_{\ell,\alpha},\DDhat_{\ell,\alpha}^{\lambda})$
$((\ell,\alpha)\in S)$
and an isomorphism
\begin{equation}
\label{eq;21.8.14.1}
 \nbige_{|\Hhat^{\lambda}_{\infty,q}}
\simeq
 \bigoplus_{(\ell,\alpha)\in S}
 \LLhat^{\lambda}_q(\ell,\alpha)
 \otimes
 (\Psi_q^{\lambda})^{\ast}
 (\nbigvhat_{\ell,\alpha},\DDhat_{\ell,\alpha}^{\lambda}).
\end{equation}
Moreover, there exist good filtered $\lambda$-flat bundles
$(\nbigp_{\ast}\nbigvhat_{\ell,\alpha},\DDhat_{\ell,\alpha}^{\lambda})$
over $(\nbigvhat_{\ell,\alpha},\DDhat_{\ell,\alpha}^{\lambda})$
such that the following holds for any $t_1\in S^1_T$:
\begin{equation}
\label{eq;17.10.12.3}
 \nbigp_{\ast}(\nbige_{|\Hhat^{\lambda}_{\infty,q}\langle t_1\rangle})
\simeq
 \bigoplus
 \Bigl(
 \nbigp^{(0)}_{\ast}\LLhat_q^{\lambda}(\ell,\alpha)
 _{|\Hhat^{\lambda}_{\infty,q}\langle t_1\rangle}
 \Bigr)
 \otimes
 \Bigl(
 \Psi_q^{\ast}\bigl(\nbigp_{\ast}\nbigvhat_{\ell,\alpha},
 \DDhat_{\ell,\alpha}^{\lambda}\bigr)
 _{|\Hhat^{\lambda}_{\infty,q}\langle t_1\rangle}
 \Bigr).
\end{equation}
\begin{lem}
For each $(\ell,\alpha)$,
there exists a good filtered $\lambda$-flat bundle
$\bigl(\nbigp_{\ast}\nbigv_{\ell,\alpha},\DD_{\ell,\alpha}^{\lambda}\bigr)$
on $(U_{w,q},\infty)$
with an isomorphism
\[
 \bigl(\nbigp_{\ast}\nbigv_{\ell,\alpha},\DD_{\ell,\alpha}^{\lambda}\bigr)
_{|\inftyhat_{w,q}}
\simeq
 (\nbigp_{\ast}\nbigvhat_{\ell,\alpha},\DDhat_{\ell,\alpha}^{\lambda}).
\]
For simplicity, if $\lambda\neq 0$,
we assume that the Stokes structure of
$(\nbigv_{\ell,\alpha},\DD_{\ell,\alpha}^{\lambda})$
is trivial.
\end{lem}
\pf
If $\lambda\neq 0$,
the claim is the same as Lemma \ref{lem;21.8.21.110}.
If $\lambda=0$,
by shrinking $U_{w,q}$,
we obtain the decomposition
$\nbige=\bigoplus_{(\ell,\alpha)\in S}
\nbige_{\ell,\alpha}$
which induces the decomposition (\ref{eq;21.8.14.1}),
i.e.,
$\nbige_{\ell,\alpha|\Hhat^{0}_{\infty,q}}
\simeq
\LLhat^{0}_q(\ell,\alpha)\otimes
(\Psi_q^{0})^{\ast}(\nbigvhat_{\ell,\alpha},
\DDhat^0_{\ell,\alpha})$.
There exist
meromorphic Higgs bundles
$(\nbigv_{\ell,\alpha},\DD^0_{\ell,\alpha})$
on $(U_{w,q},\infty)$
which induce
$\nbige_{\ell,\alpha}\otimes
\LL^0_q(\ell,\alpha)^{-1}$.
There exist isomorphisms
$(\nbigv_{\ell,\alpha},\DD^0_{\ell,\alpha})_{|\inftyhat_{w,q}}
\simeq
 (\nbigvhat_{\ell,\alpha},\DD^0_{\ell,\alpha})$.
Hence, we obtain good filtered Higgs bundles
$(\nbigp_{\ast}\nbigv_{\ell,\alpha},\DD^0_{\ell,\alpha})$
over $(\nbigp_{\ell,\alpha},\DD^0_{\ell,\alpha})$
for which 
$(\nbigp_{\ast}\nbigv_{\ell,\alpha},\DD^0_{\ell,\alpha})_{|\inftyhat_{w,q}}
\simeq
 (\nbigp_{\ast}\nbigvhat_{\ell,\alpha},\DD^0_{\ell,\alpha})$
holds.
\hfill\qed

\vspace{.1in}

Recall that there exists a Hermitian metric
$h_{\LL,q,\ell,\alpha}$
of
$\LL^{\lambda\ast}_q(\ell,\alpha)=
 \LL^{\lambda}_q(\ell,\alpha)_{|\nbigb^{\lambda\ast}_q}$
such that
$\bigl(
 \LL^{\lambda\ast}_q(\ell,\alpha),
 h_{\LL,q,\ell,\alpha}
 \bigr)$
is a monopole
and satisfies the norm estimate with respect to
$\nbigp_{\ast}\LL^{\lambda}_q(\ell,\alpha)$.
(See \S\ref{subsection;21.8.14.2}.)

Let $(V_{\ell,\alpha},\DD_{\ell,\alpha}^{\lambda})$
be the $\lambda$-flat bundle on $U_{w,q}^{\ast}$
obtained as the restriction
of $(\nbigp_{\ast}\nbigv_{\ell,\alpha},\DD_{\ell,\alpha}^{\lambda})$.
Let $h_{1,V,\ell,\alpha}$ be a Hermitian metric 
of $V_{\ell,\alpha}$
as in Lemma \ref{lem;17.10.13.100}.
Let $h_{1,\ell,\alpha}$ be the Hermitian metric of
$E_{1,\ell,\alpha}:=
 \LL^{\lambda\ast}_q(\ell,\alpha)
 \otimes
 \Psi_q^{\ast}(V_{\ell,\alpha},\DD^{\lambda}_{\ell,\alpha})$
induced by
$h_{\LL,q,\ell,\alpha}$
and $\Psi_q^{-1}(h_{1,V,\ell,\alpha})$.

\begin{lem}
\label{lem;20.7.25.1}
 The following holds, where we consider the norms
 with respect to $h_{1,\ell,\alpha}$.
\begin{itemize}
 \item $h_{1,\ell,\alpha}$ satisfies the norm estimate
       for the filtered bundle
       $\nbigp_{\ast}\LL^{\lambda}_q(\ell,\alpha)
       \otimes
       \Psi_q^{\ast}(\nbigp_{\ast}\nbigv_{\ell,\alpha},
        \DD^{\lambda}_{\ell,\alpha})$.
 \item
      $\bigl|G(h_{1,\ell,\alpha})\bigr|
     =O\bigl(|w_q|^{-2q-2\epsilon}\bigr)$
     for some $\epsilon>0$.
\item
$\bigl|\nabla_{\betabar_0}G(h_{1,\ell,\alpha})\bigr|$,
$\bigl|\nabla_{\beta_0}G(h_{1,\ell,\alpha})\bigr|$,
$\bigl|\del_{E_{1,\ell,\alpha},t_0}G(h_{1,\ell,\alpha})\bigr|$
and 
$\bigl|\del'_{E_{1,\ell,\alpha},t_0}G(h_{1,\ell,\alpha})\bigr|$
     go to $0$ as $|w_q|\to\infty$.
     (See {\rm\S\ref{subsection;17.10.4.1}} for
     mini-complex coordinate systems $(t_0,\beta_0)$.)
\item
Let $\nabla_{1,\ell,\alpha}$ and $\phi_{1,\ell,\alpha}$
be the Chern connection and the Higgs field
associated with
the mini-holomorphic bundle
$(E_{1,\ell,\alpha},\delbar_{E_{1,\ell,\alpha}})$
and the Hermitian metric $h_{1,\ell,\alpha}$.
Then,
$\Bigl|
\bigl[
 \nabla_{1,\ell,\alpha},
 \nabla_{1,\ell,\alpha}
 \bigr]\Bigr|\to 0$
and 
$\Bigl|\nabla_{1,\ell,\alpha}(\phi_{1,\ell,\alpha})\Bigr|\to 0$
as $|w_q|\to\infty$.
We also have
     $|\phi_{1,h,\ell}|=O\bigl(\log|w|\bigr)$
     as $|w_q|\to\infty$.
 \item
Let $\del_{E_{1,\ell,\alpha},h_{1,\ell,\alpha},\beta_1}$
be the operator induced by
$\del_{E_{1,\ell,\alpha},\betabar_1}$
and $h_{1,\ell,\alpha}$
as in {\rm\S\ref{subsection;17.10.5.120}}.
Then, 
\begin{equation}
\Bigl|
 \bigl[
 \del_{E_{1,\ell,\alpha,\betabar_1}},
 \del_{E_{1,\ell,\alpha},h_{1,\ell,\alpha},\beta_1}
 \bigr]\Bigr|
=O\bigl(|w|^{-2}(\log|w|)^{-2}\bigr).
\end{equation}
\end{itemize}
If $\rank V_{\ell,\alpha}=1$,
 $G(h_{1,\ell,\alpha})=0$ holds.
\end{lem}
\pf
The claim follows from 
our construction of $h_{1,\ell,\alpha}$,
Proposition \ref{prop;17.10.14.20},
Lemma \ref{lem;17.10.13.100}
and the formulas in 
Lemma \ref{lem;17.10.13.101},
Lemma \ref{lem;17.10.12.2}
and Lemma \ref{lem;17.10.24.50}.
See also Remark \ref{rem;21.8.22.11}.
\hfill\qed

\vspace{.1in}

We set
$\nbige_1:=
 \bigoplus_{\ell,\alpha}
 \LL_q^{\lambda}(\ell,\alpha)
\otimes
 (\Psi^{\lambda}_q)^{\ast}
 (\nbigp V_{\ell,\alpha},\DD^{\lambda}_{\ell,\alpha})$.
We obtain the naturally induced good filtered bundle
$\nbigp_{\ast}\nbige_1$ over $\nbige_1$.
By the construction and the isomorphism (\ref{eq;17.10.12.3}),
there exists the isomorphism
$\widehat{f}:\nbige_{|\Hhat^{\lambda}_{\infty,q}}
\simeq
 \nbige_{1|\Hhat^{\lambda}_{\infty,q}}$
which induces an isomorphism of filtered bundles
$\nbigp_{\ast}\nbige_{|\Hhat^{\lambda}_{\infty,q}}
\simeq
 \nbigp_{\ast}\nbige_{1|\Hhat^{\lambda}_{\infty,q}}$.

Let $\nbigc^{\infty}_{\nbigb^{\lambda}_q}$
denote the sheaf of $C^{\infty}$-functions
on $\nbigb^{\lambda}_q$.
There exists an isomorphism
\[
f_{C^{\infty}}:
 \nbige\otimes_{\nbigo}
  \nbigc^{\infty}_{\nbigb^{\lambda}_q}
  \simeq\nbige_1\otimes_{\nbigo}
   \nbigc^{\infty}_{\nbigb^{\lambda}_q}
\]
which induces 
$\widehat{f}$ at $\Hhat^{\lambda}_{\infty,q}$.
Let $h_1$ be the Hermitian metric on $E$
induced by 
the isomorphism $f_{C^{\infty}}$
and the metric
$\bigoplus 
 h_{\LL,q,\ell,\alpha}
 \otimes
 \Psi_q^{-1}(h_{1,V,\ell,\alpha})$.
Then, $h_1$ has the desired property.
\hfill\qed

\subsection{The strong adaptedness and the norm estimate}

Let $\nbigp_{\ast}\nbige$
and $(E,\delbar_E)$ be as in \S\ref{subsection;17.10.13.500}.
Let $h$ be a Hermitian metric of $E$
which is strongly adapted to $\nbigp_{\ast}\nbige$,
such that $(E,\delbar_E,h)$ is a monopole.
Let $h_1$ be a Hermitian metric of $E$
as in Proposition \ref{prop;17.10.12.6}.
By the assumption of the strong adaptedness of $h$
(see Definition \ref{df;20.8.8.3})
for any $\delta>0$, there exist
$C_{\delta}\geq 1$
such that 
$C_{\delta}^{-1}|w_q|^{-\delta}h_1
\leq
 h
\leq
C_{\delta}|w_q|^{\delta}h_1$.

\begin{lem}
\label{lem;17.10.15.10}
$h$ and $h_1$ are mutually bounded,
 i.e.,
$(\nbigp_{\ast}\nbige,h)$ satisfies the norm estimate.
\end{lem}
\pf
Let $s_1$ be the automorphism of $E$
determined by
$h=h_1s_1$.
By Corollary \ref{cor;17.10.15.12},
there exist positive constants $C_2$ and $\epsilon$
such that
\[
 -\Bigl(
 \del_{E,\betabar_0}\del_{E,\beta_0}
+\frac{1}{4}\del_{E,t_0}\del'_{E,t_0}
 \Bigr)
 \log\Tr(s_1)
\leq
\frac{1}{2}
 |G(h_1)|_{h_1}
\leq C_2|w|^{-2-\epsilon}.
\]
Note that
$\del_{\beta_0}\del_{\betabar_0}+\frac{1}{4}\del_{t_0}\del_{t_0}
=\del_{w}\del_{\wbar}+\frac{1}{4}\del_t\del_t$.
There exists $C_3>0$
such that
\[
  -\Bigl(
 \del_{E,\betabar_0}\del_{E,\beta_0}
+\frac{1}{4}\del_{E,t_0}\del'_{E,t_0}
 \Bigr)
\Bigl(
 \log\Tr(s_1)
+C_3|w|^{-\epsilon}
\Bigr)
\leq
 0.
\]
We set $g:=\log\Tr(s_1)+C_3|w|^{-\epsilon}$.
We take $R_0$ such that
$\{|w_q|\geq R_0\}\subset U_{w,q}^{\ast}$.
There exists $C_5>0$
such that 
$g<C_5$ on $\Psi_q^{-1}(\{|w_q|=R_0\})$.
For any $\rho>0$,
we consider the function
\[
G_{\rho}:=g-\bigl(
 C_5+\rho\log|w|
 \bigr).
\]
Note that
(i)
$G_{\rho}$ is subharmonic,
(ii) $G_{\rho}<0$ on $\Psi_q^{-1}(\{|w_q|=R\})$,
(iii) $G_{\rho}\to-\infty$ as $|w_q|\to\infty$.
Hence, we obtain $G_{\rho}\leq 0$
on $\Psi_q^{-1}(\{|w_q|\geq R\})$.
By taking the limit as $\rho\to 0$,
we obtain
\[
 \Tr(s_1)\leq 
 \exp\bigl(
 C_5-C_3|w|^{-\epsilon}
 \bigr)
 \leq \exp(C_5).
\]
i.e., $\Tr(s_1)$ is bounded.

Let $s_2$ be the automorphism of $E$ determined by
$h_1=h\cdot s_2$.
By a similar argument,
we obtain that
$\Tr(s_2)$ is bounded.
We obtain that $h$ and $h_1$ are
mutually bounded.
\hfill\qed

\subsection{Proof of Proposition \ref{prop;17.10.12.5}}

Let $(E,\delbar_E,h)$ be as in Proposition \ref{prop;17.10.12.5}.
We take a Hermitian metric $h_1$ 
as in Proposition \ref{prop;17.10.12.6}.
We obtain the automorphism $s$ of $E$
determined by
$h=h_1\cdot s$.
According to Lemma \ref{lem;17.10.15.10},
$h$ and $h_1$ are mutually bounded,
and hence $s$ and $s^{-1}$ are bounded.
Let $\dvol$ denote the volume form
of $\nbigb^{\lambda\ast}_q$ with the metric $g$.

Let $U_{w,1}$ be a relatively compact neighbourhood of $\infty$
in $U_w$.
Let $\nbigb_{q,1}^{\lambda}$
be the induced relatively compact neighbourhood of
$H^{\lambda}_{\infty,q}$
in $\nbigb_q^{\lambda}$.
(See \S\ref{subsection;21.8.21.20}.)
We set
$\nbigb_{q,1}^{\lambda\ast}:=
\nbigb_{q,1}^{\lambda}
\setminus
H^{\lambda}_{\infty,q}$.

\begin{lem}
\label{lem;17.10.13.201}
We obtain
$\int_{\nbigb^{\lambda\ast}_{q,1}}
 \Bigl(
 \bigl|
 \del_{E,h_1,\beta_0}s
 \bigr|^2_{h_1}
+\bigl|
 \del_{E,h_1,t_0}s
 \bigr|^2_{h_1}
\Bigr)\dvol<\infty$.
\end{lem}
\pf
By Lemma \ref{lem;17.10.13.110},
we have
\[
-\Bigl(
 \del_{\betabar_0}\del_{\beta_0}
+\frac{1}{4}\del_{t_0}^2
 \Bigr)\Tr(s)
=-\frac{1}{2}\Tr\bigl(sG(h_1)\bigr)
-\bigl|
 s^{-1/2}\del_{E,h_1,\beta_0}s
 \bigr|_{h_1}^2
-\frac{1}{4}
 \bigl|
 s^{-1/2}\del'_{E,h_1,t_0}s
 \bigr|^2_{h_1}.
\]
Recall that
$\nbigb^{\lambda\ast}_q
=\nbigb^{0\ast}_q
=S^1_T\times U_{w,q}^{\ast}$
as Riemannian manifolds,
and that 
\[
-\del_{\betabar_0}\del_{\beta_0}
-\frac{1}{4}\del_{t_0}^2
=-\del_{\wbar}\del_w
-\frac{1}{4}\del_t^2.
\]
For any function $F$ on $\nbigb^{0\ast}_q$,
let $\int_{S^1_T}F$ denote the function 
on $U_{w,q}^{\ast}$
obtained as the integration of $F$ along the fiber direction.
We set
\[
 b_1:= \int_{S^1_T}\Tr(s),
\quad
 b_2:=\frac{1}{2}\int_{S^1_T}\Tr(sG(h_1)),
\]
\[
 b_3:=\int_{S^1_T}\bigl|
 s^{-1/2}\del_{E,h_1,\beta_0}s
 \bigr|_{h_1}^2
+\frac{1}{4}
\int_{S^1_T} \bigl|
 s^{-1/2}\del'_{E,h_1,t_0}s
 \bigr|^2_{h_1}.
\]
We obtain 
$-\del_{\wbar}\del_wb_1
=-b_2-b_3$.
Let $\tau_q:=w_q^{-1}$.
Then, we obtain the following:
\[
 -\del_{\tau_q}\del_{\taubar_q}b_1
=-q^2\bigl|\tau_q\bigr|^{-2(q+1)}\cdot b_2
-q^2\bigl|\tau_q\bigr|^{-2(q+1)}\cdot b_3.
\]
Because 
$G(h_1)=O(|w|^{-2-\epsilon})=
O\bigl(|\tau_q|^{2q+q\epsilon}\bigr)$
as in Lemma \ref{lem;20.7.25.1},
we obtain
\[
|\tau_q|^{-2(q+1)}\cdot b_2
=O(|\tau_q|^{-2+q\epsilon}).
\]
Hence, there exists a bounded function $c$
such that
$-\del_{\tau_q}\del_{\taubar_q}c
=-q^2\bigl|\tau_q\bigr|^{-2(q+1)}\cdot b_2$.
We obtain
\[
 -\del_{\tau_q}\del_{\taubar_q}(b_1-c)
=-q^2\bigl|\tau_q\bigr|^{-2(q+1)}\cdot b_3.
\]
Note that $b_1-c$ is bounded,
and $b_3\geq 0$.
Hence, according to \cite[Lemma 2.2]{Simpson90},
we obtain 
\[
\int_{U_{w,1,q}} \bigl|\tau_q\bigr|^{-2(q+1)}\cdot b_3\,
 |d\tau_q\,d\taubar_q|<\infty.
\]
Here, $U_{w,1,q}$ denotes the pull back of
$U_{w,1}$ by the ramified covering $\proj^1_{w_q}\lrarr \proj^1_w$.
It implies that
$\int_{U_{w,1,q}} b_3\,|dw\,d\wbar|<\infty$.
Because $s$ and $s^{-1}$ are bounded,
we obtain the claim of Lemma \ref{lem;17.10.13.201}.
\hfill\qed

\vspace{.1in}

By the construction of the Higgs field,
we have
$\phi_h=\phi_{h_1}-\frac{\sqrt{-1}}{2}s^{-1}\del'_{h_1,t_0}s$.
(See \S\ref{subsection;16.9.20.20}.)
Note that $|\phi_{h_1}|_{h_1}=O\bigl(\log|w|\bigr)$.
There exists $R_0>0$ such that
$\{|w_q^q|>R_0\}\subset U_{w,q,1}^{\ast}$.
Let $B_{(w_q,t)}(2)$ be the ball with radius $2$
centered at $(w_q,t)$,
where we use the distance induced by
the Euclidean metric $dt\,dt+dw\,d\wbar$.
By Lemma \ref{lem;17.10.13.201},
there exists a constant $C>0$
such that the following holds for any $(w,t)$
with $|w|>2R_0+10$:
\[
 \int_{B_{(w,t)}(2)}
 \bigl|
 \phi_h
 \bigr|_h^2
\leq
 C\bigl(
 \log|w|
 \bigr)^2.
\]
Let $(x,y)$ be the local real coordinate system on $\cnum_{w_q}$
defined by $w=x+\sqrt{-1}y$.
By the equality (\ref{eq;21.8.22.10}),
we obtain
\[
-(\del_t^2+\del_x^2+\del_y^2)\bigl|\phi\bigr|^2_{h}
=-2\bigl(
  \bigl|\nabla_x\phi\bigr|^2_h
+\bigl|\nabla_y\phi\bigr|^2_h
+\bigl|\nabla_t\phi\bigr|^2_h
 \bigr)
 \leq 0,
\]
and hence we obtain
$|\phi_h|_h=O(\log|w|)$
by the mean-value theorem for subharmonic functions.

\vspace{.1in}
For the proof of 
the condition $F(\nabla)\to 0$,
we shall apply an argument in 
\cite{Mochizuki-KH-infinite}.
We indicate only an outline.
We set
$\nbigbtilde^{\lambda\ast}_{q}:=
 S^1\times\nbigb^{\lambda\ast}_q$.
It is equipped with the complex structure
induced by the local complex coordinate systems
$(\alpha_0,\beta_0)=(s_0+\sqrt{-1}t_0,\beta_0)$
as in \S\ref{subsection;17.10.4.1}.
It is equipped with the K\"ahler metric
$\gtilde=d\alpha_0d\alphabar_0+d\beta_0d\betabar_0$.
Let $(\Etilde,\delbar_{\Etilde})$ denote the holomorphic bundle
on $\nbigbtilde^{\lambda\ast}_q$
induced by $(E,\delbar_E)$
as in \S\ref{subsection;13.11.29.2}.
The bundle is equipped with
the induced metrics $\htilde$ and $\htilde_1$.
As recalled in \S\ref{subsection;17.10.3.10},
$(\Etilde,\delbar_{\Etilde},\htilde)$
is an instanton.
\begin{lem}
\label{lem;17.10.13.200}
The following conditions are satisfied,
where we consider the norms with respect to $\htilde_1$ 
and the Euclidean metric $\gtilde$.
\begin{itemize}
\item
$|F(\htilde_1)|\to 0$ as $|w_q|\to\infty$.
\item
$|\Lambda F(\htilde_1)|=O(|w_q|^{-2q-\epsilon})$.
\item
$|\nabla_{\htilde_1}\Lambda F(\htilde_1)|\to 0$ as $|w_q|\to\infty$.
 \end{itemize}
\end{lem}
\pf
Let $p$ denote the projection
$\nbigbtilde^{\lambda\ast}_q\lrarr 
 \nbigb^{\lambda\ast}_q$.
We have
$F(\htilde_1)=
 p^{\ast}F(\nabla_{h_1})
+p^{\ast}\nabla_{h_1}\phi_{h_1}\cdot ds_0$.
We also have
$\sqrt{-1}\Lambda F(\htilde_1)
=p^{\ast}G(h_1)$.
Hence, we obtain the claims of the lemma
from the assumption for $h_1$.
\hfill\qed

\vspace{.1in}

Let $\stilde$ be determined by
$\htilde=\htilde_1\stilde$.
Because $\stilde$ is the pull back of $s$,
$\stilde$ and $\stilde^{-1}$
are bounded
with respect to $\htilde_1$.
By Lemma \ref{lem;17.10.13.201},
we have 
$\int \bigl|\del_{\htilde_1}\stilde\bigr|^2_{\htilde_1,\gtilde}<\infty$.
Let $\stilde_1$ be determined by
$\htilde_1=\htilde\cdot \stilde_1$.
Because
$\stilde_1^{-1}\del_{\htilde}\stilde_1
=-\stilde\del_{\htilde_1}\stilde$,
we obtain
$\int \bigl|\del_{\htilde}\stilde_1\bigr|^2_{\htilde,\gtilde}<\infty$.
For any $\epsilon>0$,
there exists a compact subset 
$K_{\epsilon}\subset \nbigbtilde^{\lambda\ast}_q$
such that the following holds.
\begin{itemize}
\item
$\int_{\nbigbtilde^{\lambda\ast}_q\setminus K_{\epsilon}}
 |\stilde_1^{-1}\del_{E,\htilde}\stilde_1|_{\htilde,\gtilde}^2<\epsilon$.
\item
 $\bigl|\del_{\Etilde,\htilde_1} 
     \Lambda F(\htilde_1)\bigr|_{\htilde,\gtilde}<\epsilon$
     on $\nbigbtilde^{\lambda\ast}_q\setminus K_{\epsilon}$.
\item
 $\bigl| F(\htilde_1)\bigr|_{\htilde,\gtilde}<\epsilon$
 on $\nbigbtilde^{\lambda\ast}_q\setminus K_{\epsilon}$.
\end{itemize}
As proved in \cite[\S2.9.1]{Mochizuki-KH-infinite},
there exist $C_i>0$ $(i=1,2)$
which are independent of $(K_{\epsilon},\epsilon)$,
such that the following holds
on $\nbigbtilde^{\lambda\ast}_q\setminus K_{\epsilon}$:
\[
-\Bigl(
 \del_{\beta_0}\del_{\betabar_0}
+\frac{1}{4}\del_{t_0}^2
 -C_2\epsilon
 \Bigr)
 \bigl|
 \stilde_1^{-1}\del_{\htilde}\stilde_1
 \bigr|_{\htilde,\gtilde}^2
\leq
 C_1\epsilon.
\]
By using \cite[Theorem 9.20]{Gilbarg-Trudinger},
we obtain the following.
\begin{lem}
For any $\epsilon>0$,
there exists a compact subset
$K_{\epsilon}'\subset \nbigbtilde^{\lambda\ast}_q$
such that the following holds
on $\nbigbtilde^{\lambda\ast}_q\setminus K_{\epsilon}'$:
\[
 \sup_{\nbigbtilde^{\lambda\ast}_q\setminus K_{\epsilon}'}
 \bigl|
 \stilde_1^{-1}\del_{\htilde}\stilde_1
 \bigr|_{\htilde,\gtilde}^2
\leq
  \epsilon.
\]
\hfill\qed
\end{lem}

According to \cite[Lemma 3.1]{Simpson88},
we have the following relation:
\begin{equation}
\label{eq;17.10.13.400}
 \Delta_{\Etilde,\htilde_1}(\stilde)
=-\stilde\sqrt{-1}\Lambda F(\htilde_1)
+\sqrt{-1}\Lambda
 \delbar(\stilde)\stilde^{-1}\del_{\htilde_1}(\stilde).
\end{equation}
Hence, for any $\epsilon>0$,
there exists a compact subset 
$K^{(3)}_{\epsilon}\subset \nbigbtilde^{\lambda\ast}_q$
such that
$\bigl|
 \Delta_{\Etilde,\htilde_1}(\stilde)
 \bigr|_{\htilde_1}\leq \epsilon$
on $\nbigbtilde^{\lambda\ast}_q\setminus K^{(3)}_{\epsilon}$.

Take a large $p$.
We obtain that
the $L_2^p$-norm  of $\stilde-\id$ on
the disc with radius $r_0$ centered
at $P\in \nbigbtilde^{\lambda\ast}_q$
goes to $0$ as $P$ goes to $\infty$.
By using (\ref{eq;17.10.13.400}),
we obtain that
the $L_3^p$-norm of $\stilde_1-\id$ on 
the disc with radius $r_0$ centered
at $P\in \nbigbtilde^{\lambda\ast}_q$
goes to $0$ as $P$ goes to $\infty$.
It implies that
$|F(\htilde)|_{\htilde,\gtilde}\to 0$ as $P$ goes to $\infty$.
Therefore,
$|F(h)|_{h}\to 0$.
Thus, we obtain 
the claim of Proposition \ref{prop;17.10.12.5}.
\hfill\qed

\section{Some functoriality}
\label{subsection;20.8.1.50}

Let $\nbigb_q^{\lambda\ast}$
and $\nbigb_q^{\lambda}$ be as in \S\ref{subsection;17.10.6.100}.
Let $(E_i,h_{i},\nabla_{i},\phi_{i})$ $(i=1,2)$
be monopoles on $\nbigb^{\lambda\ast}_q$
satisfying GCK-condition.
(See \S\ref{subsection;20.8.8.110}.)
We obtain the locally free 
$\nbigo_{\nbigb^{\lambda}_q}(\ast H^{\lambda}_{\infty,q})$-module
$\nbigp E^{\lambda}_i$
and the associated good filtered bundles
$\nbigp_{\ast}E^{\lambda}_i$
over $\nbigp E^{\lambda}_i$
(Proposition \ref{prop;17.9.17.50} and
Theorem \ref{thm;17.10.5.130}).
Clearly, the monopoles
\[
(E_1,h_1,\nabla_1,\phi_1)\oplus
(E_2,h_2,\nabla_2,\phi_2),
\]
\[
(E_1,h_1,\nabla_1,\phi_1)\otimes
(E_2,h_2,\nabla_2,\phi_2),
\]
\[
\Hom((E_1,h_1,\nabla_1,\phi_1),
(E_2,h_2,\nabla_2,\phi_2))
\]
also satisfy GCK-condition.

\begin{prop}
 There exists the following natural isomorphisms
 of good filtered bundles:
\begin{equation}
 \nbigp_{\ast}(E_1^{\lambda}\oplus E_2^{\lambda})
\simeq
\nbigp_{\ast}(E_1^{\lambda})
\oplus
\nbigp_{\ast}(E_2^{\lambda})
\end{equation}
\begin{equation}
 \nbigp_{\ast}(E_1^{\lambda}\otimes E_2^{\lambda})
\simeq
\nbigp_{\ast}(E_1^{\lambda})
\otimes
\nbigp_{\ast}(E_2^{\lambda})
\end{equation}
\begin{equation}
 \nbigp_{\ast}(\Hom(E_1,E_2)^{\lambda})
\simeq
\nhom\Bigl(
\nbigp_{\ast}(E_1^{\lambda}),
\nbigp_{\ast}(E_2^{\lambda})
\Bigr).
\end{equation}
\end{prop}
\pf
It follows from Lemma \ref{lem;17.10.7.200}.
\hfill\qed

\chapter{Global periodic monopoles of rank one }

We shall establish an equivalence between
monopoles of GCK type and
good filtered bundles in the rank one case.

\section{Preliminary}

\subsection{Ahlfors type lemma}
\label{subsection;21.7.7.10}

Take $R>0$ and $C_0>0$.
Set $U^{\ast}_w(R):=\{w\in\cnum\,|\,|w|>R\}$.
Let $g:U^{\ast}_w(R)\lrarr \real_{\geq 0}$ be 
a $C^{\infty}$-function such that
\[
 -\del_{w}\del_{\wbar}g\leq -C_0g.
\]
We assume that $g=O(|w|^N)$ for some $N>0$.
\begin{lem}
\label{lem;17.10.7.300}
There exists $\epsilon_1>0$, depending only on $C_0$,
such that
$g=O\bigl(
 \exp(-\epsilon_1|w|)
 \bigr)$.
\end{lem}
\pf
There exist $\epsilon_1>0$ and $R_1\geq R$
such that the following holds
on $\{w\in\cnum\,|\,|w|>R_1\}$:
\[
 -\del_{w}\del_{\wbar}\exp\bigl(-\epsilon_1|w|\bigr)
\geq
 -C_0\exp\bigl(-\epsilon_1|w|\bigr),
\]
\[
 -\del_w\del_{\wbar}\exp\bigl(\epsilon_1|w|\bigr)
\geq
 -C_0\exp\bigl(\epsilon_1|w|\bigr).
\]
There exists $C_2>0$ such that
$g<C_2\exp\bigl(-\epsilon_1|w|\bigr)$
on $\{|w|=R_1\}$.
For any $\delta>0$,
we set
\[
 F_{\delta}:=
 C_2\exp\bigl(-\epsilon_1|w|\bigr)
+\delta\exp\bigl(\epsilon_1|w|\bigr).
\]
We obtain $g<F_{\delta}$ on $\{|w|=R_1\}$
and
$-\del_w\del_{\wbar}(g-F_{\delta})
\leq
 -C_0(g-F_{\delta})$.
We set
\[
 Z(\delta):=\Bigl\{
 w\in\cnum\,\Big|\,
 |w|\geq R_1,\,\,
 g(w)>F_{\delta}(w)
 \Bigr\}.
\]
Because we have $g(w)=O(|w|^N)$ around $\infty$
and $g<F_{\delta}$ on $\{|w|=R_1\}$,
$Z(\delta)$ is relatively compact in $\{|w|>R_1\}$.
Hence, 
$g-F_{\delta}=0$ holds
on $\del Z(\delta)$.
On $Z(\delta)$,
we have
\[
 -\del_w\del_{\wbar}(g-F_{\delta})<0.
\]
Hence, if $Z(\delta)\neq\emptyset$,
we obtain
$g-F_{\delta}\leq 0$ on $Z(\delta)$,
which contradicts the construction of $Z(\delta)$.
Hence, $Z(\delta)=\emptyset$ holds.
Namely, we obtain
$g\leq F_{\delta}$ on $\{|w|>R_1\}$
for any $\delta>0$.
Therefore, we obtain
$g\leq C_2\exp(-\epsilon_1|w|)$.
\hfill\qed

\vspace{.1in}

We give a variant of the estimate.
Let $g:U^{\ast}_w(R)\lrarr \real_{\geq 0}$
be a $C^{\infty}$-function such that
the following holds for some $c>0$ and $N>0$.
\begin{itemize}
\item
 $g=O(|w|^N)$.
\item
For any $k>0$, there exists $b_k>0$ such that 
$-\del_w\del_{\wbar}g\leq b_k|w|^{-k}+(-c+|w|^{-2})g$.
\end{itemize}
\begin{lem}
\label{lem;17.10.25.50}
$g=O(|w|^{-k})$ for any $k>0$.
\end{lem}
\pf
There exists $R_0>R$ such that
$g(w)<2^{-1}|w|^{2N}$ on $|w|\geq R_0$.
We fix $k>0$,
and we take $R_k\geq R_0$ such that 
\[
 cR_k^2\geq \frac{b_k}{R_k^{k-2+2N}}+\frac{(k-2)^2}{4}+1.
\]
We set $e_k=R_k^{k-2+2N}$.
We obtain  $g(w)< e_k|w|^{-k+2}$
on $\{|w|=R_k\}$.
We also have the following inequality
on $\{|w|\geq R_k\}$:
\[
 c-|w|^{-2}
\geq \Bigl(
 \frac{b_k}{e_k}+\frac{(k-2)^2}{4}
 \Bigr)|w|^{-2}.
\]
For any small $\epsilon>0$,
we set $\rho_{k,\epsilon}:=e_k|w|^{-k+2}+\epsilon|w|^{2N}$.
For any large $k$, we obtain the following:
\begin{multline}
 -\del_w\del_{\wbar}\rho_{k,\epsilon}
=-e_k\frac{(k-2)^2}{4}|w|^{-k}
 -\epsilon N^2|w|^{2N-2}
 \\
=b_k|w|^{-k}
 -\Bigl(
 \frac{b_k}{e_k}+\frac{(k-2)^2}{4}
 \Bigr) |w|^{-2}e_k|w|^{-k+2}
-\epsilon N^2|w|^{2N-2}
 \\
\geq
 b_k|w|^{-k}
  -\Bigl(
 \frac{b_k}{e_k}+\frac{(k-2)^2}{4}
 \Bigr) |w|^{-2}
 \rho_{k,\epsilon}.
\end{multline}
Let us consider the set
$Z_{\epsilon}:=\{w\in\cnum\,|\,|w|\geq R_k,\,g(w)>\rho_{k,\epsilon}(w)\}$.
Because $Z_{\epsilon}$ is relatively compact 
in $\{|w|\geq R_k\}$,
we have
$g(w)-\rho_{k,\epsilon}(w)=0$ on 
$\del Z_{\epsilon}$.
On $Z_{\epsilon}$,
we have the following:
\[
 -\del_w\del_{\wbar}(g-\rho_{k,\epsilon})
\leq
  -\Bigl(
 \frac{b_k}{e_k}+\frac{(k-2)^2}{4}
 \Bigr) |w|^{-2}
 (g-\rho_{k,\epsilon}) 
 <0.
\]
We obtain $g-\rho_{k,\epsilon}\leq 0$
on $Z_{\epsilon}$,
which contradicts the choice of $Z_{\epsilon}$.
Thus, we obtain $Z_{\epsilon}=\emptyset$,
i.e.,
$g\leq \rho_{k,\epsilon}$ on
$\{|w|\geq R_k\}$ for any $\epsilon$.
By taking the limit $\epsilon\to 0$,
we obtain $g\leq e_k|w|^{-k+2}$
on $\{|w|\geq R_k\}$.
\hfill\qed

\subsection{Poisson equation (1)}

Let $a$ be a $C^{\infty}$-function on 
$S^1_T\times \{|w|>R\}$ 
satisfying the following conditions.
\begin{itemize}
\item
 For any $(\ell_1,\ell_2,\ell_3)\in\seisuu_{\geq 0}^2$,
 we have
 $|\del_t^{\ell_1}\del_w^{\ell_2}\del_{\wbar}^{\ell_3}a(t,w)|
 =O(|w|^{-k})$ for any $k$ as $w\to \infty$.
\item
 $\int_{S^1_T}a(t,w)\,dt=0$ for any $w$.
\end{itemize}

Let $f$ be an $\real$-valued $C^{\infty}$-function on
$S^1_T\times \{|w|>R\}$
such that
(i) $\Delta f=a$,
(ii) $\int_{S^1_T}f(t,w)\,dt=0$ for any $w$,
(iii) $f=O(|w|^N)$ for some $N>0$.
Here, $\Delta=-\del_t^2-\del_x^2-\del_y^2$
for the real coordinate system $(x,y)$ on $\cnum_w$
determined by $w=x+\sqrt{-1}y$.
\begin{lem}
\label{lem;17.10.7.301}
For any $\vecell=(\ell_1,\ell_2,\ell_3)\in\seisuu_{\geq 0}^3$,
we obtain
$\bigl|
 \del_t^{\ell_1}\del_{w}^{\ell_2}\del_{\wbar}^{\ell_3}f
 \bigr|
=O\bigl(|w|^{-k}\bigr)$
for any $k$.
If the support of $a$ is compact,
there exists $\epsilon_{\vecell}>0$
such that
$\bigl|
 \del_t^{\ell_1}\del_w^{\ell_2}\del_{\wbar}^{\ell_3}f
 \bigr|
=O\Bigl(
 \exp\bigl(-\epsilon_{\vecell}|w|\bigr)
 \Bigr)$.
\end{lem}
\pf
The following holds:
\[
 \Delta |f|^2
=-2\Bigl(
 \bigl| \del_tf \bigr|^2
+\bigl| \del_xf \bigr|^2
+\bigl| \del_yf \bigr|^2
 \Bigr)
+2af.
\]
We set 
$g(x,y):=\int_{S^1_T}|f|^2(t,x,y)dt$.
Because $\int_{S^1_T} f(t,x,y)\,dt=0$,
there exists $C>0$ such that
$\int_{S^1_T}
 \bigl|
 \del_tf
 \bigr|^2\,dt
\geq
 Cg$.
We obtain
\[
 \Delta g
\leq
 -Cg+2\Bigl(
 \int_{S^1_T}|a|^2
 \Bigr)^{1/2}
 g^{1/2}.
\]
For any $k>0$,
there exists $b_k>0$ such that the following holds:
\[
 \Delta g\leq
 -Cg+b_k|w|^{-k}g^{1/2}
\leq
 b_k^2|w|^{-2k+2}
+(-C+|w|^{-2})g.
\]
We also have $g=O(|w|^N)$ for some $N>0$.
By Lemma \ref{lem;17.10.25.50},
we obtain 
$g=O\bigl(|w|^{-k}\bigr)$ for any $k$.
Then, we obtain the first claim
by using Lemma \ref{lem;21.9.13.20} below.
We obtain the second claim
by using a similar argument
with Lemma \ref{lem;17.10.7.300}.
\hfill\qed

\subsubsection{Appendix}

For $r>0$,
let $B(r)
=\bigl\{(x_1,x_2,x_3)\in\real^3\,\big|\,
\sum x_i^2< r^2\bigr\}$.
For any $\vecell=(\ell_1,\ell_2,\ell_3)\in\seisuu_{\geq 0}^3$,
we set $|\vecell|=\sum\ell_i$
and 
$\del^{\vecell}_{\vecx}:=
\del_{x_1}^{\ell_1}\del_{x_2}^{\ell_2}\del_{x_3}^{\ell_3}$.
For any $\ell\in\seisuu_{\geq }$0,
we set $\nbigs(\ell)=\bigl\{
\vecell\in\seisuu_{\geq 0}^3\,\big|,
|\vecell|\leq \ell
\bigr\}$.
We set $\Delta=-\sum\del_{x_i}^2$.

Let $a$ be a $C^{\infty}$-function on $B(r)$
such that the following holds.
\begin{itemize}
 \item For any $\ell\in\seisuu_{\geq 0}$,
       there exists $C_{\ell}>0$ such that
       $|\del_{\vecx}^{\vecell}a|
       <C_{\ell}$ for any $\vecell\in\nbigs(\ell)$.
\end{itemize}
Let $f$ be a $C^{\infty}$-function on $B(r)$
such that (i) $\Delta f=a$
and that (ii) $\int_{B(r)} |f|^2<A_0$ for a constant $A_0>0$.

\begin{lem}
\label{lem;21.9.13.20}
 Let $0<r'<r$.
For any $\ell\in\seisuu_{\geq 0}$,
there exists $\Ctilde_{1,r',\ell}>0$,
depending only on $r$, $r'$, $A_0$
and $C_{\ell'}$ $(\ell'\leq \ell)$
such that
$|\del_{\vecx}^{\vecell}f|
<\Ctilde_{1,r',\ell}$ for any $\vecell\in\nbigs(\ell)$
on $B(r')$.
\end{lem}
\pf
Let $0<r'<r_0<r$.
By \cite[Theorem 9.20]{Gilbarg-Trudinger},
there exists $\Ctilde_{1,r_0,1}>0$,
depending only on $r$, $r_0$, $A_0$ and $C_{0}$ 
such that
$|f|<\Ctilde_{1,r_0,0}$ on $B(r_0)$.
Let $r'<r_1<r_0$.
We apply \cite[Theorem 4.6]{Gilbarg-Trudinger}
to the Poisson equation $\Delta f=a$.
Then, there exists $\Ctilde_{1,r_1,1}>0$,
depending only on 
$r$, $r_i$ $(i=0,1)$, $A_0$ and $C_{\ell'}$ $(\ell'\leq 1)$
such that
$|\del_{x_j}f|<\Ctilde_{1,r_1,1}$ $(j=1,2,3)$ on $B(r_1)$.
Let $r'<r_2<r_1$.
We apply \cite[Theorem 4.6]{Gilbarg-Trudinger}
to the Poisson equations
$\Delta (\del_{x_i}f)=\del_{x_i}a$ $(i=1,2,3)$.
Then, 
there exists $\Ctilde_{1,r_2,2}>0$,
depending only on 
$r$, $r_i$ $(i=0,1,2)$, $A_0$ and $C_{\ell'}$ $(\ell'\leq 2)$
such that
$|\del_{x_j}\del_{x_i}f|<\Ctilde_{1,r_2,2}$ $(i,j=1,2,3)$ on $B(r_2)$.
Then, by an easy induction, we obtain the claim of Lemma \ref{lem;21.9.13.20}.
\hfill\qed

\subsection{Poisson equation (2)}

Let $b$ be a $C^{\infty}$-function
on $S^1_T\times \cnum_w$
such that
(i) $\int_{S^1_T\times\cnum_w}
 b\dvol_{S^1_T\times\cnum_w}=0$,
(ii) $|\del_t^{\ell_1}\del_w^{\ell_2}\del_{\wbar}^{\ell_3}b|
=O\bigl(|w|^{-k}\bigr)$ as $|w|\to\infty$
for any $(\ell_1,\ell_2,\ell_3)\in\seisuu_{\geq 0}^3$
and for any $k$.

\begin{lem}
\label{lem;17.9.19.20}
There exists a $C^{\infty}$-function $f$ on $S^1_T\times\cnum_w$
such that $\Delta_{S^1_T\times\cnum_w}f=b$
and $|f|=O(|w|^{-1})$ as $|w|\to\infty$.
\end{lem}
\pf
Set $c:=\int_{S^1_T}b\,dt$.
We obtain $\int_{\real^2} c\,dx\,dy=0$,
and $\del_x^{\ell_1}\del_y^{\ell_2}c=O(|w|^{-k})$
for any $(\ell_1,\ell_2)\in\seisuu_{\geq 0}^2$
and for any $k$.
We have the natural compactification
$\cnum_w\subset \proj^1$.
We may regard $c$ as a $C^{\infty}$-function on $\proj^1$.
There exists an $\real_{>0}$-valued $C^{\infty}$-function $A$
on $\cnum_w$
such that 
$A\Delta_{\real^2}=\Delta_{\proj^1|\cnum}$
and $\dvol_{\real^2}=A\dvol_{\proj^1|\cnum}$.
We obtain
$\int_{\proj^1}Ac\,\dvol_{\proj^1}
=\int_{\real^2}c\,d\vol_{\real^2}=0$.
By the standard theory of the Poisson equations
on compact manifolds (for example, see \cite[Theorem 4.7]{Aubin-book}),
there exists $\gamma$
such that 
$\gamma(\infty)=0$
and
$\Delta_{\real^2}\gamma=c$.
We have $\gamma=O(|w|^{-1})$.

By considering $b-\frac{1}{T}\int_{S^1_T} b$,
we may assume that
$\int_{S^1_T} b=0$ from the beginning.
For any $P\in S^1_T\times\cnum_w$
and $s>0$,
let $B_P(s)$ denote the set of
$Q\in S^1_T\times\cnum_w$
whose distance from $P$ is less than $s$.
Let $V_P(s)$ denote the volume of $B_P(s)$.
Set $k_b(P,s):=V_P(s)^{-1}\int_{B_P(s)}|b|$.
Because
$k_b(P,s)=O(s^{-2})$ as $s\to\infty$,
we obtain $\int_0^rsk_b(P,s)\,ds=O\bigl(\log(r+2)\bigr)$.
We set
$b_+(P):=\max\{b(P),0\}$
and
$b_-(P):=\max\{-b(P),0\}$.
According to \cite[341 page]{Ni-Shi-Tam},
there exist functions $f_{\pm}$
such that $\Delta f_{\pm}=b_{\pm}$
and $|f_{\pm}(w)|=O\bigl(\log(2+|w|)\bigr)$.
We set $f=f_+-f_-$.
Then,
$\Delta(f)=b$
and $|f|=O\bigl(\log(2+|w|)\bigr)$
are satisfied.
By the elliptic regularity
(for example, see \cite[Theorem 3.23]{Mizohata-book}),
$f$ is $C^{\infty}$.
By Lemma \ref{lem;17.10.7.301},
we obtain $|f|=O\bigl(|w|^{-k}\bigr)$ 
for any $k>0$ as $|w|\to\infty$.
Thus, we obtain the claim of the lemma.
\hfill\qed

\begin{lem}
$f$ induces a $C^{\infty}$-function on
$S^1_T\times\proj^1_w$.
In particular,
the $\ell$-th derivative of $f$ is $O(|w|^{-1-\ell})$.
\end{lem}
\pf
There exists the decomposition $f=f^{\circ}+f^{\bot}$,
where $f^{\circ}$ is constant along $S^1_T$,
and $\int_{S^1_T} f^{\bot}=0$.
By Lemma \ref{lem;17.10.7.301},
$f^{\bot}$ and its higher derivatives
are $O(|w|^{-k})$ for any $k$.
Because $f^{\circ}$ is a $C^{\infty}$-function on $\proj^1_w$
such that $f^{\circ}(\infty)=0$,
we obtain the claim of the lemma.
\hfill\qed

\subsection{Subharmonic functions}

Let $f$ be a bounded function
$S^1_T\times\cnum_w\lrarr \real_{\geq 0}$
such that $\Delta f\leq 0$ in the sense of distributions.

\begin{lem}
$f$ is constant.
\end{lem}
\pf
It is enough to consider the case where
$f$ is $C^{\infty}$.
We set $F:=\int_{S^1_T}f$.
We obtain
$-\del_{w}\del_{\wbar}F\leq 0$ on $\cnum$,
and $F$ is bounded.
Then, $F$ induces a bounded subharmonic function on $\proj^1$,
and hence $F$ is constant.
(See \cite[Proposition 2.2]{Simpson88}.)

There exists the decomposition $f=f_0+f_1$
such  that
(i) $f_0$ is constant in the $S^1_T$-direction,
(ii) $\int_{S^1_T}f_1=0$.
Because $F=\int_{S^1_T}f=\int_{S^1_T}f_0$,
we obtain that $f_0$ is constant.

We have the following inequality:
\[
 \Delta|f|^2
=-2\bigl|\del_tf\bigr|^2
-2\bigl|\del_xf\bigr|^2
-2\bigl|\del_yf\bigr|^2
+2\Delta(f)\cdot f
\leq
-\bigl|\del_tf\bigr|^2.
\]
We obtain the following inequality:
\[
 -4\del_w\del_{\wbar}\int_{S^1_T} |f|^2
\leq
 -\int_{S^1_T}\bigl|\del_tf\bigr|^2.
\]
Note that $\int_{S^1_T}|f|^2=\int_{S^1_T}|f_0|^2+\int_{S^1_T}|f_1|^2$,
and that $\int_{S^1_T}|f_0|^2$ is constant.
We also have
\[
 \int_{S^1_T}\bigl|\del_tf\bigr|^2
=\int_{S^1_T}\bigl|\del_tf_1\bigr|^2
\geq
 C_1\int_{S^1_T}\bigl|f_1\bigr|^2.
\]
We obtain the following:
\[
 -4\del_w\del_{\wbar}\int_{S^1}\bigl|f_1\bigr|^2
\leq
 -C_1\int_{S^1_T}\bigl|f_1\bigr|^2.
\]
Because $\int_{S^1_T}\bigl|f_1\bigr|^2$ is bounded,
we obtain that
$\int_{S^1_T}\bigl|f_1\bigr|^2
=O\bigl(\exp(-\epsilon|w|)\bigr)$ for some $\epsilon>0$
by Lemma \ref{lem;17.10.7.300}.
Because $\int_{S^1_T}\bigl|f_1\bigr|^2$
is non-negative and subharmonic,
we obtain that
$\int_{S^1_T}\bigl|f_1\bigr|^2=0$,
and hence $f_1=0$.
\hfill\qed

\begin{cor}
\label{cor;17.10.16.13}
Let $g$ be a bounded function on $S^1_T\times\cnum_w$
such that $\Delta g=0$.
Then, $g$ is constant.
\hfill\qed
\end{cor}

\section{Global periodic monopoles of rank one (1)}
\label{subsection;17.10.2.21}

We use the notation in \S\ref{subsection;17.10.2.1}.
Take a complex number $\gamma$.
Let us consider the product line bundle
$L(\gamma)$ 
on $S^1_T\times \cnum_w$
with a global frame $\vtilde$
and the metric $h_{L(\gamma)}$
given by $h_{L(\gamma)}(\vtilde,\vtilde)=1$.
We consider the unitary connection $\nabla$
and the Higgs field $\phi$
given by
\[
 \nabla\vtilde
=\vtilde\bigl(-\sqrt{-1}(\gamma+\gammabar)\bigr)\,dt,
\quad\quad
 \phi=\gamma-\gammabar.
\]
Because $F(\nabla)=0$ and $\nabla\phi=0$,
$(L(\gamma),h_{L(\gamma)},\nabla,\phi)$
is a monopole on $S^1_T\times\cnum_w$.
\index{monopole $(L(\gamma),h_{L(\gamma)},\nabla,\phi)$}

\begin{rem}
The restriction of this monopole around infinity
on a ramified covering
appeared in {\rm\S\ref{subsection;17.10.25.122}}.
\hfill\qed
\end{rem}

By (\ref{eq;17.10.7.310}) and (\ref{eq;17.10.7.311}),
the underlying mini-holomorphic bundle
$L^{\lambda}(\gamma)$ on $\nbigm^{\lambda}$
is described as follows:
\index{mini-holomorphic bundle $L^{\lambda}(\gamma)$}
\[
 \del_{L^{\lambda}(\gamma),\betabar_1}\vtilde
=\vtilde\frac{\lambda\gammabar}{1+|\lambda|^2},
\quad
 \del_{L^{\lambda}(\gamma),t_1}\vtilde
=\vtilde\frac{-2\sqrt{-1}(\gamma-|\lambda|^2\gammabar)}{1+|\lambda|^2}.
\]

Let $L^{\lambda\cov}(\gamma)$
denote the pull back of 
$L^{\lambda}(\gamma)$
by the projection
$\varpi^{\lambda}:
M^{\lambda}
\lrarr
\nbigm^{\lambda}$.
\index{mini-holomorphic bundle $L^{\lambda\cov}(\gamma)$}
\index{projection $\varpi^{\lambda}$}
On $M^{\lambda}$, we define
\[
 \utilde=
 \exp\Bigl(
 \frac{-\lambda \gammabar\betabar_1}{1+|\lambda|^2}
+\frac{\lambdabar\gamma\beta_1}{1+|\lambda|^2}
+\frac{2\sqrt{-1}}{1+|\lambda|^2}
 (\gamma-|\lambda|^2\gammabar)t_1
 \Bigr)
\cdot
(\varpi^{\lambda})^{-1}(\vtilde).
\]
Then,
$\utilde$ is a mini-holomorphic frame,
i.e.,
$\del_{L^{\lambda\cov}(\gamma),t_1}\utilde=0$
and 
$\del_{L^{\lambda\cov}(\gamma),\betabar_1}\utilde=0$.
Because
\[
 \bigl|\utilde\bigr|_h
=\exp\bigl(
 -2\Image(\gamma)t_1
 \bigr),
\]
$\utilde$ is a global frame of
the locally free $\nbigo_{\Mbar^{\lambda}}(\ast H^{\lambda\cov}_{\infty})$-module
$\nbigp L^{\lambda\cov}(\gamma)$.
Moreover, the good filtered bundle
$\nbigp_{\ast} L^{\lambda\cov}(\gamma)_{|\{t_1\}\times\proj^1_{\beta_1}}$
is described as
\[
 \nbigp_{b}L^{\lambda\cov}(\gamma)_{|\{t_1\}\times\proj^1_{\beta_1}}
=\nbigo_{\proj^1}([b])\utilde_{|\{t_1\}\times\proj^1_{\beta_1}}.
\]
Here, $[b]:=\max\{n\in\seisuu\,|\,n\leq b\}$.

The locally free $\nbigo_{\nbigmbar^{\lambda}}(\ast H^{\lambda}_{\infty})$-module
$\nbigp L^{\lambda}(\gamma)$
is described as the descent of 
$\nbigp L^{\lambda\cov}(\gamma)$
by the action
\[
 \kappa_1^{-1}(\utilde)
=\utilde\exp(2\sqrt{-1}\gamma T).
\]
We also obtain the description of the good filtered bundle
$\nbigp_{\ast} L^{\lambda}(\gamma)$
as the descent.

\begin{lem}
Let $h_1$ be another Hermitian metric of
$L^{\lambda}(\gamma)$
such that 
(i) $(L^{\lambda}(\gamma),h)$ is a monopole,
(ii) $h_1$ is strongly adapted to 
the good filtered bundle $\nbigp_{\ast}L^{\lambda}(\gamma)$.
Then, there exists a positive constant $a$
such that $h_1=a\cdot h$.
\end{lem}
\pf
According to Proposition \ref{prop;17.10.12.5},
$h$ and $h_1$ are mutually bounded.
There exists the function $s:\nbigm^{\lambda}\lrarr\real_{>0}$ 
determined by
$h_1=h\cdot s$.
Note that $\Delta s=0$ on $\nbigm^{\lambda}$
according to Corollary \ref{cor;17.10.15.12}.
Hence, we obtain that $s$ is constant
by Corollary \ref{cor;17.10.16.13}.
\hfill\qed

\subsection{Reformulation}

Let $\alpha$ be a non-zero complex number.
Let $\nbigl^{\cov}(\alpha)$
be the $\nbigo_{\Mbar^{\lambda}}(\ast H^{\lambda\cov}_{\infty})$-module
$\nbigo_{\Mbar^{\lambda}}(\ast H^{\lambda\cov}_{\infty})\,e$.
\index{$\nbigo_{\Mbar^{\lambda}}(\ast H^{\lambda\cov}_{\infty})$-module
$\nbigl^{\cov}(\alpha)$}
Let $\pi^{\lambda}:\Mbar^{\lambda}\lrarr \real_{t_1}$
denote the projection,
i.e.,
$(t_1,\beta_1)\longmapsto t_1$.
\index{projection $\pi^{\lambda}$}
We obtain the good filtered bundles
$\nbigp^{(0)}_{\ast}(\nbigl^{\cov}(\alpha))$
by setting
\index{good filtered bundle $\nbigp^{(0)}_{\ast}(\nbigl^{\cov}(\alpha))$}
\[
 \deg^{\nbigp^{(0)}}(e_{|\pi^{-1}(t_1)})=0.
\]
We consider the $\seisuu$-action
by $\kappa_1^{\ast}e=\alpha \cdot e$.
Then, we obtain a locally free
$\nbigo_{\nbigmbar^{\lambda}}(\ast H^{\lambda}_{\infty})$-module
$\nbigl(\alpha)$
and the good filtered bundle
$\nbigp^{(0)}_{\ast}\nbigl(\alpha)$ over $\nbigl(\alpha)$.
\index{$\nbigo_{\nbigmbar^{\lambda}}(\ast H^{\lambda}_{\infty})$-module
$\nbigl(\alpha)$}
\index{good filtered bundle $\nbigp^{(0)}_{\ast}\nbigl(\alpha)$}

\begin{lem}
There exists a Hermitian metric $h$ of
$\nbigl(\alpha)_{|\nbigm^{\lambda}}$
such that
(i) $(\nbigl(\alpha)_{|\nbigm^{\lambda}},h)$
is a monopole,
(ii) $h$ is strongly adapted to $\nbigp^{(0)}_{\ast}\nbigl(\alpha)$.
Such a metric is unique up to the multiplication of
positive constants.
\hfill\qed
\end{lem}

\section{Global periodic monopoles of rank one (2)}

\subsection{Construction of mini-holomorphic bundles}

Let $P\in\nbigm^{\lambda}$.
Let $\Ptilde=(t^0_1,\beta^0_1)\in (\varpi^{\lambda})^{-1}(P)$
such that $0\leq t^0_1<T$.
We choose $0<\epsilon<T-t^0_1$.
We set
$U:=\openopen{-\epsilon}{T-\epsilon/2}\times \proj^1$
and 
$H^{\lambda}_{\infty,\epsilon}:=
 \openopen{-\epsilon}{T-\epsilon/2}\times \{\infty\}$.
Let $p_i$ $(i=1,2)$
denote the projections
$p_1:U\setminus\{\Ptilde\}\lrarr \openopen{-\epsilon}{T-\epsilon/2}$
and 
$p_2:U\setminus\{\Ptilde\}\lrarr \proj^1$.
Let $\nbigl^{\cov}_{-\epsilon,T-\epsilon/2}(P,\ell)$ 
be the $\nbigo_{U\setminus \{\Ptilde\}}(\ast H^{\lambda}_{\infty,\epsilon})$-module
determined by the following conditions:
\begin{itemize}
\item
 $\nbigl^{\cov}_{-\epsilon,T-\epsilon/2}(P,\ell)
   (\ast p_2^{-1}(\beta^0_1))$
 is isomorphic to the pull back of 
 $\nbigo_{\proj^1}(\ast \{\beta_1^0,\infty\})\,e$.
\item
We have
\[
  \nbigl^{\cov}_{-\epsilon,T-\epsilon/2}(P,\ell)
  _{|p_1^{-1}(a)}
 =\left\{
 \begin{array}{ll}
  \nbigo_{\proj^1}(\ast\infty) e
 & (-\epsilon<a<t_1^0),\\
 \mbox{{}}\\
  \nbigo_{\proj^1}(-\ell \beta_1^0)(\ast\infty) e
 & (t_1^0<a<T-\epsilon/2),
 \end{array}
 \right.
\]
under the above isomorphism.
\end{itemize}
Let $\kappa_1\colon
 \openopen{-\epsilon}{-\epsilon/2}
\lrarr 
 \openopen{T-\epsilon}{T-\epsilon/2}$
be the isomorphism given by
$\kappa_1(t_1,\beta_1)=
 (t_1+T,\beta_1+2\sqrt{-1}\lambda T)$.
We have the isomorphism
\[
 \kappa_1^{\ast}\Bigl(
 \nbigl^{\cov}_{-\epsilon,T-\epsilon/2}(P,\ell)
 _{|\openopen{T-\epsilon}{T-\epsilon/2}\times\proj^1}
 \Bigr)
\simeq
  \nbigl^{\cov}_{-\epsilon,T-\epsilon/2}(P,\ell)
 _{|\openopen{-\epsilon}{-\epsilon/2}\times\proj^1}
\]
given by
$\kappa_1^{\ast}
 \Bigl(
 (\beta_1-\beta_1^0)^{\ell}e
 \Bigr)
\longmapsto
 e$,
or equivalently
$\kappa_1^{\ast}(e)\longmapsto
 (\beta_1-\beta_1^0+2\sqrt{-1}\lambda T)^{-\ell}e$.
Hence, by gluing,
we obtain 
an $\nbigo_{\nbigmbar^{\lambda}\setminus\{P\}}
 (\ast H^{\lambda}_{\infty})$-module
 denoted by $\nbigl(P,\ell)$.
\index{sheaf $\nbigl(P,\ell)$}

\vspace{.1in}

Let $\nbigl^{\cov}(P,\ell)$
denote the pull back of $\nbigl(P,\ell)$
by $\varpi^{\lambda}$.
\index{sheaf $\nbigl^{\cov}(P,\ell)$}
We shall describe $\nbigl^{\cov}(P,\ell)$ explicitly.
Let $\Ptilde_n=(t^{(n)}_1,\beta_1^{(n)})$ denote 
the point of $(\varpi^{\lambda})^{-1}(P)$
such that $nT\leq t^{(n)}_1<(n+1)T$.
We set
\[
 H(\Ptilde_n,+):=
 \bigl\{(t_1,\beta_1^{(n)})\,\big|\, 
 t_1>t_1^{(n)}
\bigr\},
\quad\quad
 H(\Ptilde_n,-):=
 \bigl\{(t_1,\beta_1^{(n)})\,\big|\, 
 t_1<t_1^{(n)}
\bigr\}.
\]
Then, we have the following isomorphism:
\[
\nbigl^{\cov}(P,\ell)
\simeq
 \nbigo_{\nbigmbar^{\lambda}\setminus (\varpi^{\lambda})^{-1}(P)}
 \Bigl(
 \sum_{n<0} \ell H(\Ptilde_n,-)
-\sum_{n\geq 0}\ell H(\Ptilde_n,+)
 \Bigr).
\]

Let us describe the isomorphism
$\kappa_1^{\ast}
 \nbigl^{\cov}(P,\ell)
\simeq
 \nbigl^{\cov}(P,\ell)$.
If $n\geq 0$, 
the isomorphism
\[
 \kappa_1^{\ast}
 \Bigl(
 \nbigl^{\cov}(P,\ell)_{|\openopen{t_1^0+nT}{t_1^0+(n+1)T}\times\proj^1}
 \Bigr)
\simeq
 \nbigl^{\cov}(P,\ell)_{|\openopen{t_1^0+(n-1)T}{t_1^0+nT}\times\proj^1}
\]
is given by
\[
 \kappa_1^{\ast}\Bigl(
 \prod_{m=0}^n
 \bigl(\beta_1-(\beta_1^0+2\sqrt{-1}\lambda mT)^{\ell}\bigr)e
 \Bigr)
\longmapsto
 \prod_{m=0}^{n-1}
 \bigl(\beta_1-(\beta_1^0+2\sqrt{-1}\lambda mT)^{\ell}\bigr)e.
\]
If $n<0$,
the isomorphism
\[
 \kappa_1^{\ast}
 \Bigl(
 \nbigl^{\cov}(P,\ell)_{|\openopen{t_1^0+nT}{t_1^0+(n+1)T}\times\proj^1}
 \Bigr)
\simeq
 \nbigl^{\cov}(P,\ell)_{|\openopen{t_1^0+(n-1)T}{t_1^0+nT}\times\proj^1}
\]
is given by
\[
 \kappa_1^{\ast}
 \Bigl(
 \prod_{m=0}^{-n-1}
 \bigl(\beta_1-(\beta^0_1-2\sqrt{-1}\lambda mT)^{-\ell}\bigr)e
 \Bigr)
\longmapsto
  \prod_{m=0}^{-n}
 \bigl(\beta_1-(\beta^0_1-2\sqrt{-1}\lambda mT)^{-\ell}\bigr) e.
\] 

\subsection{Good filtered bundles}

Let $e^{t_1}$ denote the restriction of $e$
to $p_1^{-1}(t_1)$.
We obtain the filtered bundle
$\nbigp^{(a)}_{\ast}\nbigl^{\cov}(P,\ell)$
over $\nbigl^{\cov}(P,\ell)$
defined by the following condition:
\index{filtered bundle $\nbigp^{(a)}_{\ast}\nbigl^{\cov}(P,\ell)$}
\[
 \deg^{\nbigp^{(a)}}(e^{t_1})
=-\ell\frac{t_1-t_1^0}{T}+a.
\]
Because $\nbigp_{\ast}\nbigl^{\cov}(P,\ell)$
is $\seisuu$-equivariant,
we obtain the filtered bundle
$\nbigp_{\ast}\nbigl(P,\ell)$
over $\nbigl(P,\ell)$.
\index{filtered bundle $\nbigl(P,\ell)$}

\begin{lem}
We have
$\deg\bigl(\nbigp^{(a)}_{\ast}\nbigl(P,\ell)\bigr)
=-(a+\ell/2)$.
\end{lem}
\pf
Indeed,
\begin{multline}
 T\deg\bigl(\nbigp^{(a)}_{\ast}\nbigl(P,\ell)\bigr)
=\int_0^T
 \deg\bigl(\nbigp_{\ast}^{(a)}\nbigl(P,\ell)_{|p_1^{-1}(t_1)}\bigr)\,dt_1
\\
 =\int_{0}^{t_1^0}
 \Bigl(
 \ell\frac{t_1-t_1^0}{T}-a
 \Bigr)\,dt_1
+\int_{t_1^0}^T
 \Bigl(
 \ell\frac{t_1-t_1^0}{T}-\ell-a
 \Bigr)\,dt_1
\\
=\int_0^T\ell\frac{t_1}{T}\,dt_1
-aT-\int_0^T\ell\frac{t_1^0}{T}\,dt_1
-\ell(T-t_1^0)
=-aT-\frac{\ell T}{2}.
\end{multline}
Thus, we obtain the desired equality.
\hfill\qed

\subsection{Monopoles}

Let $\nbigl^{\ast}(P,\ell)$
denote the mini-holomorphic bundle
on $\nbigm^{\lambda}\setminus \{P\}$
obtained as the restriction of $\nbigl(P,\ell)$.
\index{sheaf $\nbigl^{\ast}(P,\ell)$}

\begin{prop}
There exists a Hermitian metric $h$ of $\nbigl^{\ast}(P,\ell)$
such that the following holds.
\begin{itemize}
\item
$(\nbigl^{\ast}(P,\ell),h)$ is a monopole of Dirac type 
on $\nbigmlambda\setminus \{P\}$.
\item
$h$ is strongly adapted to $\nbigp^{(-\ell/2)}_{\ast}\nbigl(P,\ell)$
in the sense of Definition {\rm\ref{df;20.8.8.3}}.
\end{itemize}
Such $h$ is unique up to the multiplication
of positive constants.
\end{prop}
\pf
Let $\Psi^{\lambda}:\nbigmbar^{\lambda}\lrarr \proj^1_w$
be given by
$\Psi^{\lambda}(t_1,\beta_1)=
 (1+|\lambda|^2)^{-1}(\beta_1-2\sqrt{-1}\lambda t_1)$.
Let $U_w:=\{|w|>R\}\cup\{\infty\}$.
Set $U:=(\Psi^{\lambda})^{-1}(U_w)$.

By using Proposition \ref{prop;17.10.12.6}
in the rank one case,
we can construct a Hermitian metric $h_1$
of $\nbigl(P,\ell)$ such that 
(i) $(\nbigl(P,\ell),h_1)$ satisfies the norm estimate
with respect to $\nbigp^{(-\ell/2)}_{\ast}\nbigl(P,\ell)$,
(ii) 
 $\del_{t_1}^{\ell_1}
 \del_{\beta_1}^{\ell_2}
 \del_{\betabar_1}^{\ell_3}
 G(h_1)=O\bigl(|w|^{-k}\bigr)$
for any $(\ell_1,\ell_2,\ell_3)\in\seisuu_{\geq 0}^3$
and for any $k$ as $|w|\to\infty$,
(iii) $G(h_1)=0$ around $P$.
By Corollary \ref{cor;17.10.24.31},
we have 
$G(h_1e^{\varphi})-G(h_1)=2^{-1}\Delta \varphi$.
We also have
$\int G(h_1)=2\pi T\deg(\nbigp^{(-\ell/2)}_{\ast}\nbigl)=0$.
We can take a bounded $C^{\infty}$-function $\varphi$
such that
$-2^{-1}\Delta \varphi=-G(h_1)$
by Lemma \ref{lem;17.9.19.20}.
Then, $h=h_1e^{\varphi}$
has the desired property.

Suppose that $h'$ is another metric satisfying the conditions.
We obtain the $C^{\infty}$ function $\varphi$
on $\nbigm^{\lambda}\setminus P$
such that $h'=he^{\varphi}$,
which is bounded.
We have $\Delta\varphi=0$
on $\nbigm^{\lambda}\setminus P$,
i.e., $\varphi$ is a harmonic function on
$\nbigm^{\lambda}\setminus P$.
Recall that isolated singularities of bounded harmonic functions
are removable.
(Fore example,
see \cite[Theorem 2.3]{Harmonic-Function}.)
Then, we obtain that $\varphi$ is a constant
by Corollary \ref{cor;17.10.16.13}.
\hfill\qed

\section{Global periodic monopoles of rank one (3)}

Let $Z$ be a finite set.
Let $\nbigl$ be a locally free
$\nbigo_{\nbigmbar^{\lambda}\setminus Z}
 (\ast H^{\lambda}_{\infty})$-module
of rank one with Dirac type singularity at $Z$.
For each $Q\in \nbigmbar^{\lambda}\setminus Z$,
let $\nbigl_Q$ denote the stalk of the sheaf $\nbigl$ at $Q$.
For each $P\in Z$,
let $(\tau,\zeta)$ be a mini-complex coordinate system
around $P$ such that $(\tau(P),\zeta(P))=(0,0)$.
For a small positive number,
we set
$P_{+}=(\epsilon,0)$
and $P_-=(-\epsilon,0)$.
By the scattering map,
we obtain the isomorphism
$\nbigl_{P_-}(\ast P_-)
\simeq
\nbigl_{P_+}(\ast P_+)$
by which we identify
$\nbigl_{P_{\pm}}(\ast P_{\pm})$.
We obtain the integer
\[
 \ell(P):=
 \length\Bigl(
 \nbigl_{P_-}\bigl/(\nbigl_{P_-}\cap\nbigl_{P_+})
 \Bigr)
- \length\Bigl(
 \nbigl_{P_+}\bigl/(\nbigl_{P_-}\cap\nbigl_{P_+})
 \Bigr).
\]
Then, there exists a non-zero complex number $\alpha$
and an isomorphism
\[
 \nbigl\simeq
 \nbigl(\alpha)\otimes
 \bigotimes_{P\in Z} \nbigl(P,\ell(P)).
\]
There exists the unique good filtered bundle
$\nbigp_{\ast}\nbigl$ over $\nbigl$
such that
$\deg(\nbigp_{\ast}\nbigl)=0$,
for which
\[
 \nbigp_{\ast}\nbigl
\simeq
 \nbigp^{(0)}_{\ast}
 \nbigl(\alpha)
\otimes
 \bigotimes_{P\in Z}
 \nbigp^{(-\ell(P)/2)}_{\ast}
 \nbigl(P,\ell(P)).
\]
We obtain the following.
\begin{prop}
\label{prop;17.10.13.150}
There exists a Hermitian metric $h$
of $\nbigl_{|\nbigm^{\lambda}\setminus Z}$
such that
(i) $(\nbigl_{|\nbigm^{\lambda}\setminus Z},h)$ is a monopole
of Dirac type,
(ii) $h$ is strongly adapted to $\nbigp_{\ast}\nbigl$.
Such a metric is unique up to the multiplication of
positive constants.
\hfill\qed
\end{prop}

\chapter[Global periodic monopoles]{Global periodic monopoles and filtered difference modules}

We shall prove the main theorem (Theorem \ref{thm;17.9.30.20})
of this monograph.
In \S\ref{subsection;20.8.1.10},
we shall prove that
for a given periodic monopole of GCK-type,
the associated good filtered bundle
is polystable of degree $0$.
In \S\ref{subsection;17.10.26.2},
we shall prove the existence of a monopole of GCK type
which induces a given polystable good filtered bundle with degree $0$.
Then, Theorem \ref{thm;17.9.30.20} follows.

As one of the consequences of Theorem \ref{thm;17.9.30.20},
the classification of singular monopoles of GCK-type
is reduced to the classification of
polystable parabolic difference modules of degree $0$.
In \S\ref{section;21.9.17.110},
we explain {\em smooth} parabolic $0$-difference modules of rank $2$
are equivalent to filtered torsion-free sheaves of rank one
on spectral curves,
which is a variant of
equivalences between Higgs bundles
and sheaves on the spectral curves
in \cite{Beauville-Narasimhan-Ramanan, Hitchin-self-duality},
and we revisit the classification of $\SU(2)$-monopoles of GCK-type
without Dirac type singularity in \cite{Harland}.

\section{Statements}

Let $Z$ be a finite subset in $\nbigm$.
(See \S\ref{subsection;17.10.2.1} for $\nbigm$.)
Let $g$ denote the Euclidean metric
$dt\,dt+dw\,d\wbar$ of $\nbigm$.
\begin{df}
\index{monopole of GCK type}
A monopole $(E,h,\nabla,\phi)$ on $\nbigm\setminus Z$
is called of GCK type if the following holds.
\begin{itemize}
\item
 $|F(\nabla)|_{h}\to 0$ and $|\phi|_h=O(\log|w|)$
 as $w\to\infty$.
\item
 Each point $P\in Z$ is a Dirac type singularity of $(E,h,\nabla,\phi)$.
\hfill\qed
\end{itemize}
\end{df}

Let $(E,h,\nabla,\phi)$ be a monopole on 
$\nbigm\setminus Z$ of GCK-type.
Let $\lambda$ be any complex number.
Let $(E^{\lambda},\delbar_{E^{\lambda}})$
be the underlying mini-holomorphic bundle
on the mini-complex manifold $\nbigm^{\lambda}\setminus Z$.
It is prolonged to the associated good filtered bundle 
$\nbigp^h_{\ast}E^{\lambda}$
on $(\nbigmbar^{\lambda};Z,H^{\lambda}_{\infty})$
as explained in
Proposition \ref{prop;17.9.17.50} and Theorem \ref{thm;17.10.5.130}.
This construction is compatible with
direct sum,
tensor product,
and inner homomorphism
as in \S\ref{subsection;20.8.1.50}.
We shall prove the following theorem
in \S\ref{subsection;17.10.26.1}--\S\ref{subsection;17.10.26.2}.

\begin{thm}
\label{thm;17.9.30.20}
The above procedure induces
a bijection between the equivalence classes
of the following objects:
\begin{itemize}
\item
Monopoles of GCK-type $(E,h,\nabla,\phi)$
on $\nbigm\setminus Z$.
\item
Polystable good filtered bundles of Dirac type 
$\nbigp_{\ast}E^{\lambda}$
on $(\nbigmbar^{\lambda};Z,H^{\lambda}_{\infty})$
with $\deg(\nbigp_{\ast}E^{\lambda})=0$.
(See Definition {\rm\ref{df;21.8.22.20}}
for the polystablity condition.)
\end{itemize}
\end{thm}

We obtain the following corollary,
which is an analogue of the Corlette-Simpson correspondence
between stable good filtered flat bundles
and stable good filtered Higgs bundles.

\begin{cor}
 \label{cor;17.10.28.32}
The above procedure induces
the bijective correspondence of
the equivalence classes of
the following objects through 
monopoles of GCK-type
on $\nbigm\setminus Z$.
\begin{itemize}
\item
Polystable good filtered bundles of Dirac type
$\nbigp_{\ast}\nbige^0$
on $(\nbigmbar^0;Z,H^{0}_{\infty})$
with $\deg(\nbigp_{\ast}\nbige^{0})=0$.
\item
Polystable good filtered bundles of Dirac type
$\nbigp_{\ast}\nbige^{\lambda}$
on $(\nbigmbar^{\lambda};Z,H^{\lambda}_{\infty})$
with $\deg(\nbigp_{\ast}\nbige^{\lambda})=0$.
\hfill\qed
\end{itemize}
\end{cor}

As for the comparison with parabolic difference modules,
together with
(\ref{eq;21.9.17.4}),
Proposition \ref{prop;21.9.17.50}
and Proposition \ref{prop;21.9.17.51},
we obtain the following.
\begin{cor}
\label{cor;21.9.17.60}
There exists the natural bijection between the isomorphism classes
of the following objects.
\begin{itemize}
 \item Monopoles of GCK-type $(E,h,\nabla,\phi)$ on
       $\nbigm\setminus Z$.
 \item Polystable parabolic $2\sqrt{-1}\lambda T$-difference
       modules of degree $0$:
\[
       (\vecV,V,m_Z,(\vectau_{Z,x},\vecL_{Z,x})_{x\in \cnum},
       \nbigp_{\ast}(\vecV_{|\inftyhat})).
\]
\end{itemize} 
Here, $Z$ and $(m_Z,\vectau_{Z,x})$
are related as in 
{\rm(\ref{eq;21.9.17.10})},
{\rm(\ref{eq;21.9.17.11})} and {\rm(\ref{eq;21.9.17.12})}.
\hfill\qed
\end{cor}

\section{Preliminary}
\label{subsection;17.10.26.1}

\subsection{Ambient good filtered bundles with appropriate metric}
\label{subsection;17.10.26.12}

Let $Z$ be a finite subset in $\nbigm^{\lambda}$.
Let $\nbigp_{\ast}\nbige^{\lambda}$
be a good filtered bundle with Dirac type singularity
on $(\nbigmbar^{\lambda};Z,H_{\infty}^{\lambda})$.
Let $(E,\delbar_E)$ denote the mini-holomorphic bundle
with Dirac type singularity on $\nbigmlambda\setminus Z$
obtained as the restriction of $\nbigp\nbige^{\lambda}$.

Let $h_1$ be a Hermitian metric
of $E$ strongly adapted to $\nbigp_{\ast}\nbige^{\lambda}$
such that the following holds.
\begin{condition}
\label{condition;21.8.22.30}\mbox{{}}
\begin{description}
\item[(A1)]
 Around $H^{\lambda}_{\infty}$,
we have $|G(h_1)|_{h_1}=O(|w|^{-\epsilon-2})$
for some $\epsilon>0$,
and 
$(E,\delbar_E,h_1)$ satisfies the norm estimate
with respect to $\nbigp_{\ast}\nbige$
in the sense of Definition {\rm\ref{df;17.10.24.20}}.
Moreover,
we have
\begin{equation}
\label{eq;17.10.24.51}
\bigl|
 \bigl[
 \del_{E,\betabar_1},\del_{E,h_1,\beta_1}
 \bigr]
 \bigr|_{h_1}=O\bigl(|w|^{-2}(\log|w|)^{-2}\bigr).
\end{equation}
(See Proposition {\rm\ref{prop;17.10.12.6}}.)
\item[(A2)]
Around each point of $Z$,
$(E,\delbar_E,h_1)$ is a monopole with Dirac type singularity.
In particular,
it induces a $C^{\infty}$-metric 
of the Kronheimer resolution of $E$.
(See {\rm\S\ref{subsection;17.9.19.2}} for
the Kronheimer resolution.)
\hfill\qed
 \end{description}
\end{condition}
 
\subsection{Degree of filtered subbundles}

Let $\nbigp_{\ast}\nbige_1\subset
\nbigp_{\ast}\nbige$
be a filtered subbundle on 
$\bigl(
 \nbigmbar^{\lambda};Z,H^{\lambda}_{\infty}
 \bigr)$.
Let $(E_1,\delbar_{E_1})$ be the mini-holomorphic bundle of Dirac type
on $(\nbigm^{\lambda},Z)$.
Let $h_{1,E_1}$ denote the metric of $E_1$
induced by $h_1$.
By the Chern-Weil formula in Lemma \ref{lem;17.10.24.41},
the analytic degree
$\deg(E_1,h_{1,E_1})\in \real\cup\{-\infty\}$ makes sense.

\begin{prop}
\label{prop;17.9.10.10}
$2\pi T\deg(\nbigp_{\ast}\nbige_1)
=\deg(E_1,h_{1,E_1})$ holds.
\end{prop}
\pf
We take a metric $h_{0,E_1}$ of $E_1$
which satisfies
Condition \ref{condition;21.8.22.30}
for $\nbigp_{\ast}\nbige_1$.
Because $G(h_{0,E_1})=O(|w|^{-2-\epsilon})$ $(\epsilon>0)$
with respect to $h_{0,E_1}$
around $H^{\lambda}_{\infty}$,
and because $G(h_{0,E_1})=0$ around each point of $Z$,
$G(h_{0,E_1})$ is $L^1$ with respect to $h_{0,E_1}$
and the natural volume form of $\nbigm^{\lambda}$.
Let $\nabla_{0}$ and $\phi_0$ be 
the Chern connection and the Higgs field
associated with $(E_1,\delbar_{E_1})$ with $h_{0,E_1}$.
Because $(E_1,\delbar_{E_1},h_{0,E_1})$ is a monopole
with Dirac type singularity around each point $P$ of $Z$,
we have
$(\nabla_0\phi_0)_{|Q}=O\bigl(d(Q,P)^{-2}\bigr)$
around $P$, and hence $\nabla_0\phi_0$ is $L^1$ around $P$,
with respect to $h_{0,E_1}$ and
the Riemannian metric of $\nbigm^{\lambda}$.
Let $\del_{E_1,\beta_1}$ denote the operator
induced by $\del_{E_1,\betabar_1}$ and $h_{0,E_1}$
as in \S\ref{subsection;17.10.5.120}.
Because 
$\bigl|
[\del_{E_1,\beta_1},\del_{E_1,\betabar_1}]
\bigr|_{h_{0,E_1}}=
O\bigl(|w|^{-2}(\log|w|)^{-2}\bigr)$
around $H_{\infty}^{\lambda}$,
$\bigl|
[\del_{E_1,\beta_1},\del_{E_1,\betabar_1}]
\bigr|_{h_{0,E_1}}$
is $L^1$ around $H_{\infty}^{\lambda}$.
Hence, we may apply the formula (\ref{eq;21.8.19.11}),
and we obtain the following equality:
\[
\int
 \Tr G(h_{0,E_1})\dvol
=\int_{0}^T
 2\pi\pardeg\bigl(
 \nbigp_{\ast}\nbige_{1|\nbigmbar^{\lambda}\langle t_1\rangle
 \setminus Z}\bigr)
 \,dt_1
=2\pi T\deg(\nbigp_{\ast}\nbige_1).
\]
(See \S\ref{subsection;20.7.31.20}
for $\nbigmbar^{\lambda}\langle t_1\rangle$.)

It remains to prove the following equality:
\begin{equation}
\label{eq;17.10.24.40}
 \int
 \Tr G(h_{1,E_1})\dvol
=\int
 \Tr G(h_{0,E_1})\dvol.
\end{equation}
By considering
$\det E_1\subset \bigwedge^{\rank E_1}E$,
it is enough to consider the case $\rank E_1=1$.
According to Proposition \ref{prop;17.10.13.150},
there exists a Hermitian metric $h'_{E_1}$ of $E_1$
such that 
(i) $(E_1,\delbar_{E_1},h'_{E_1})$ is a monopole of GCK-type,
(ii) the meromorphic extension
$\nbigp^{h'_{E_1}}E_1$ is equal to
$\nbigp\nbige_1$.
We obtain $\deg(\nbigp^{h'_{E_1}}_{\ast}E_1)=0$.
By considering
$(E,\delbar_E,h_1)\otimes(E_1,\delbar_{E_1},h'_{E_1})^{-1}$
and 
$\nbigp_{\ast}\nbige\otimes
 (\nbigp_{\ast}^{h'_{E_1}}E_1)^{\lor}$,
we may reduce the issue
to the case where
$\nbigp\nbige_1$ is 
isomorphic to 
$\nbigo_{\nbigmbar^{\lambda}}(\ast H^{\lambda}_{\infty})$.
There exists the section $f$ of $\nbigp\nbige_1$
corresponding to
$1\in\nbigo_{\nbigmbar^{\lambda}}(\ast H^{\lambda}_{\infty})$.
Let $a\in\real$ be determined by
$a:=\inf\{b\in\real\,|\,f\in\nbigp_b\nbige_1\}$.
According to the following lemma,
we may and will assume $a=0$.

\begin{lem}
\label{lem;20.7.25.20}
It is enough to study the case $a=0$.
\end{lem}
\pf
Let $\psi:\cnum_{w}\lrarr\real_{>0}$
be a $C^{\infty}$-function
such that
$\psi(w)=|w|^{-2a}$
on $|w|>R$ for large $R$,
and that $\psi$ is constant around $Z$.
The metric of $E_1$ induced by
$\psi\cdot  h_1$
is equal to
$\psi\cdot h_{1,E_1}$.

Let $\nbigp'_{\ast}\nbige$ be the filtered bundle
over $\nbige$ determined by
\[
 \nbigp'_{c}\nbige_{|\nbigmbar^{\lambda}\langle t_1\rangle\setminus Z}:=
 \nbigp_{c+a}\nbige_{|\nbigmbar^{\lambda}\langle t_1\rangle\setminus Z}.
\]
Let $\nbigp'_{\ast}\nbige_1$ be the filtered bundle
over $\nbige_1$ induced by $\nbigp'_{\ast}\nbige$.
Then, $\psi\cdot h_{0,E_1}$ is a Hermitian metric of
$E_1$
which satisfies Condition \ref{condition;21.8.22.30}
for $\nbigp'_{\ast}\nbige_1$.

Note that the support of $\Delta\log\psi$ is compact.
According to Corollary \ref{cor;17.10.24.31},
the equality
$\int \Tr G(\psi\cdot h_{0,E_1})\dvol
=\int \Tr G(\psi\cdot h_{1,E_1})\dvol$
implies
$\int \Tr G(h_{0,E_1})\dvol
=\int \Tr G(h_{1,E_1})\dvol$.
Thus, Lemma \ref{lem;20.7.25.20} is proved.
\hfill\qed

\begin{lem}
\label{lem;17.10.24.30}
 Let $\nbigb^{\lambda}$ be a neighbourhood of
 $H^{\lambda}_{\infty}$ in $\nbigmbar^{\lambda}$.
 Let $E$ be a mini-holomorphic bundle
 on $\nbigb^{\lambda\ast}:=\nbigb^{\lambda}\setminus H^{\lambda}_{\infty}$
 with a metric $h$ such that 
 $G(h)$ is $L^1$.
Let $f$ be a mini-holomorphic section of $E$
such that 
\[
 C_1^{-1}
\leq
 |f|_h(\log|w|)^{-k}
\leq
 C_1
\]
for some $C_1>1$ and $k\in\real$.
Then,
$\bigl|\nabla_{\beta_0}f\bigr|_h\cdot
 |f|_h^{-1}$
and
$\bigl|
 (\nabla_{t_0}+\sqrt{-1}\ad\phi)f
 \bigr|_h\cdot
 |f|_h^{-1}$
 are $L^2$,
where $\nabla$ and $\phi$ denote the Chern connection
and the Higgs field associated with
$(E,\delbar_E,h)$.
\end{lem}
\pf
It is enough to prove that
\[
\bigl|\nabla_{\beta_0}f\bigr|_h(\log|w|)^{-k},
\quad
\bigl|(\nabla_{t_0}+\sqrt{-1}\ad\phi)f
 \bigr|_h(\log|w|)^{-k}
\]
are $L^2$.
Because $f$ is mini-holomorphic,
we have
$\nabla_{\betabar_0}f=0$
and $(\nabla_{t_0}-\sqrt{-1}\ad\phi)f=0$.
We may assume that 
$\nbigb^{\lambda\ast}=(\Psi^{\lambda})^{-1}(U^{\ast}_{w}(R))$,
where $U_w^{\ast}(R):=\{w\in\cnum\,|\,|w|>R\}$.

Let $\rho:\real\lrarr\{0\leq a\leq 1\}\subset\real_{\geq 0}$
be a $C^{\infty}$-function
such that,
(i) $\rho(t)=0$ $(t\geq 1)$,
(ii) $\rho(t)=1$ $(t\leq 1/2)$,
(iii) $\rho(t)^{1/2}$
and $\del_t\rho(t)\big/\rho(t)^{1/2}$
induce $C^{\infty}$-functions.

For any large positive integer $N$,
by setting
$\chi_N(w):=\rho\bigl(N^{-1}\log|w|\bigr)$,
we obtain $C^{\infty}$-functions
$\chi_N:U_{w}^{\ast}(R)
\lrarr \real_{\geq 0}$
such that
$\chi_N(w)=0$ $(|w|\geq e^{N})$
and 
$\chi_N(w)=1$ $(|w|\leq e^{N/2})$.
Let $\mu:U_{w}^{\ast}(R)\lrarr \real_{\geq 0}$
be a $C^{\infty}$-function
such that
$\mu(w)=1-\rho\bigl(\log(|w|/R)\bigr)$.
We set $\chitilde_N:=\mu\cdot \chi_N$.
We have 
\[
 \del_w\chitilde_N(w)
=\del_w\mu(w)\chi_N(w)
+\mu(w)\rho'(N^{-1}\log|w|^2)N^{-1}w^{-1}.
\]
By the assumption on $\rho$,
$\del_w\chitilde_N(w)\big/\chitilde_N(w)^{1/2}$
naturally give $C^{\infty}$-functions
on $U^{\ast}_{w}(R)$,
and there exists $C_2>0$, which is independent of $N$,
such that the following holds:
\[
 \bigl|
 \del_w\chitilde_N(w)\big/\chitilde_N(w)^{1/2}
 \bigr|
\leq
 C_2|w|^{-1}(\log|w|^2)^{-1}.
\]
Because
$\del_{\beta_0}w=(1+|\lambda|^2)^{-1}$
and 
$\del_{\beta_0}\wbar=-\lambdabar^2(1+|\lambda|^2)^{-1}$,
there exists $C_3>0$, which is independent of $N$,
such that the following holds:
\[
 \bigl|
 \del_{\beta_0}\bigl(
 \chitilde_N(w)\bigr)\big/\chitilde_N(w)^{1/2}
 \bigr|
\leq
 C_3|w|^{-1}(\log|w|^2)^{-1}.
\]

We consider the following integral:
\begin{multline}
\int_{\nbigb^{\lambda\ast}}
 \chitilde_N(w)\cdot h(\nabla_{\beta_0}f,\nabla_{\beta_0}f)
 (\log|w|^2)^{-2k}\,\dvol
= \\
-\int_{\nbigb^{\lambda\ast}}
 \del_{\beta_0}\bigl(\chitilde_N(w)\bigr)\cdot
 h(f,\nabla_{\beta_0}f)
 (\log|w|^2)^{-2k}\,\dvol
 \\
-\int_{\nbigb^{\lambda\ast}}
 \chitilde_N(w)\cdot
 h(f,\nabla_{\betabar_0}\nabla_{\beta_0}f)
 (\log|w|^2)^{-2k}\,\dvol
 \\
+\int_{\nbigb^{\lambda\ast}}
 \chitilde_N(w)\cdot
 h(f,\nabla_{\beta_0}f)
\cdot(-2k)(\log|w|^2)^{-2k-1}
 (w^{-1}\del_{\beta_0}w+\wbar^{-1}\del_{\beta_0}\wbar)
 \,\dvol
\end{multline}
We have the following inequality:
\begin{multline}
 \Bigl|
 \del_{\beta_0}\chitilde_N\cdot
 h(f,\nabla_{\beta_0}f)
 (\log|w|^2)^{-2k}
 \Bigr|
\leq \\
 \Bigl(
 C_3C_1|w|^{-1}(\log|w|^2)^{-1}
 \Bigr)
 \cdot
 \Bigl(
 \chitilde_N^{1/2}(w)\cdot
 \bigl|\nabla_{\beta_0}f\bigr|_h
 (\log|w|^2)^{-k}
 \Bigr).
\end{multline}
We also have the following inequality:
\begin{multline}
 \Bigl|
 \chitilde_N\cdot
 h(f,\nabla_{\beta_0}f)
\cdot(\log|w|^2)^{-2k-1}
 (w^{-1}\del_{\beta_0}w+\wbar^{-1}\del_{\beta_0}\wbar)
 \Bigr|
\leq \\
 2\Bigl(
 C_1
 \chitilde_N^{1/2}\cdot
 |w|^{-1}(\log|w|^2)^{-1}
 \Bigr)
 \cdot
 \Bigl(
 \chitilde_N^{1/2}
 \bigl|
 \nabla_{\beta_0}f
 \bigr|_h 
 (\log|w|^2)^{-k}
 \Bigr).
\end{multline}
Note that
$\nabla_{\betabar_0}\nabla_{\beta_0}f
=(\nabla_{\betabar_0}\nabla_{\beta_0}
 -\nabla_{\beta_0}\nabla_{\betabar_0})f
=-[F_{\beta_0,\betabar_0}(h),f]$.
There exists $C_4>0$
which are independent of $N$,
such that the following holds:
\begin{multline}
\label{eq;20.8.8.200}
 \int_{\nbigb^{\lambda\ast}}
 \chitilde_N\cdot
 \bigl|\nabla_{\beta_0}f\bigr|_h^2(\log|w|^2)^{-2k}\dvol
 \leq
 \\
 C_4+
 C_4\Bigl(
 \int_{\nbigb^{\lambda\ast}}
 \chitilde_N\cdot
 \bigl|\nabla_{\beta_0}f\bigr|_h^2(\log|w|^2)^{-2k}\dvol
 \Bigr)^{1/2}
\\
+\int_{\nbigb^{\lambda\ast}}
 \chitilde_N\cdot
 h\bigl(
 f,[F_{\beta_0,\betabar_0},f]
 \bigr)(\log|w|^2)^{-2k}\dvol.
\end{multline}
Similarly,
there exists $C_5>0$
which is independent of $N$,
such that the following holds:
\begin{multline}
\label{eq;20.8.8.201}
 \int_{\nbigb^{\lambda\ast}}
 \chitilde_N\cdot
 \bigl|(\nabla_{t_0}+\sqrt{-1}\ad\phi) f\bigr|_h^2(\log|w|^2)^{-2k}\dvol
\leq \\
 C_5+
 C_5\Bigl(
 \int_{\nbigb^{\lambda\ast}}
 \chitilde_N\cdot
 \bigl|(\nabla_{t_0}+\sqrt{-1}\ad\phi)f\bigr|_h^2(\log|w|^2)^{-2k}\dvol
 \Bigr)^{1/2}
\\
+\int_{\nbigb^{\lambda\ast}}
 \chitilde_N\cdot
 h\bigl(
 f,[-2\sqrt{-1}\nabla_{t_0}\phi,f]
 \bigr)(\log|w|^2)^{-2k}\dvol.
\end{multline}
Because $G(h)$ is $L^1$,
there exists a constant $C_6>0$,
which is independent of $N$,
such that the following holds:
\begin{multline}
\int_{\nbigb^{\lambda\ast}}
 \chitilde_N\cdot
 h\bigl(
 f,[F_{\beta_0,\betabar_0},f]
 \bigr)(\log|w|^2)^{-2k}\dvol
 \\
+\frac{1}{4}
 \int_{\nbigb^{\lambda\ast}}
 \chitilde_N\cdot
 h\bigl(
 f,[-2\sqrt{-1}\nabla_{t_0}\phi,f]
 \bigr)(\log|w|^2)^{-2k}\dvol
\leq C_6.
\end{multline}
We put
\begin{multline}
 A_N:=
\int_{\nbigb^{\lambda\ast}}
 \chitilde_N\cdot
 \bigl|\nabla_{\beta_0}f\bigr|_h^2(\log|w|^2)^{-2k}\dvol
\\
 +\frac{1}{4}\int_{\nbigb^{\lambda\ast}}
 \chitilde_N\cdot
 \bigl|(\nabla_{t_0}+\sqrt{-1}\ad\phi) f\bigr|_h^2(\log|w|^2)^{-2k}\dvol.
\end{multline}
By (\ref{eq;20.8.8.200}) and (\ref{eq;20.8.8.201}),
there exists $C_7>0$,
which are independent of $N$,
such that the following holds:
\[
A_N\leq C_7+C_7A_N^{1/2}.
\]
Hence, there exists $C_8>0$ such that
$A_N\leq C_8$ for any large $N$.
By taking $N\to\infty$,
we obtain the claim of Lemma \ref{lem;17.10.24.30}.
\hfill\qed

\vspace{.1in}

Let $h_{2,E_1}$ be a Hermitian metric of $E_1$
such that the following holds.
\begin{itemize}
\item
There exists a neighbourhood $N_1$ of $Z$
such that 
$h_{2,E_1}=h_{0,E_1}$ on $\nbigmlambda\setminus N_1$.
\item
There exists a neighbourhood $N_2$ of $Z$
contained in $N_1$
such that
$h_{2,E_1}=h_{1,E_1}$ on $N_2\setminus Z$.
\end{itemize}

We obtain the function $s$ determined by
$h_{1,E_1}=h_{2,E_1}\cdot s$.
According to Corollary \ref{cor;17.10.24.31},
we have the relation
$G(h_{1,E_1})-G(h_{2,E_1})=2^{-1}\Delta \log s$.
The support of $\log s$ is contained in 
$\nbigm^{\lambda}\setminus N_2$.
By using Lemma \ref{lem;17.10.24.30},
we obtain 
$\int\Delta\log s=0$.
Hence, we obtain
$\int G(h_{1,E_1})=\int G(h_{2,E_1})$.

To compare 
$\int G(h_{0,E_1})$ and $\int G(h_{2,E_1})$,
it is enough to compare the integrals
over a neighbourhood for each $P\in Z$.
We apply a variant of the argument 
in \cite{Li-Narasimhan}.
For simplicity of the description,
we consider the case $P=(0,0)$
for the coordinate system $(t_0,\beta_0)$.
Let $U$ denote the connected component of $N_1$
which contains $P$.
We may regard $U$ as an open subset of
$\real\times\cnum$.
We take the Kronheimer resolution
as in \S\ref{subsection;17.10.24.32}.
We have the map
$\varphi:\cnum^2\lrarr \real\times\cnum$
in (\ref{eq;20.7.30.40}).
Set $\Utilde=\varphi^{-1}(U)$.
We have the holomorphic vector bundle
$\Etilde$ and $\Etilde_1$
on $\Utilde$,
induced by $E$ and $E_1$,
respectively.
We may regard $\Etilde_1$
as a saturated subsheaf of $\Etilde$.
Note that $\Etilde/\Etilde_1$ is not necessarily locally free.
The metric $\htilde_0:=\varphi^{-1}(h_{0,E_1})$
induces a $C^{\infty}$-metric of $\Etilde_1$.
The metric $\htilde_2:=\varphi^{-1}(h_{2,E_1})$
may have singularity.

We take a projective morphism
$\psi:\Utilde'\lrarr \Utilde$
such that
(i) $D:=\psi^{-1}(0,0)$ is simple normal crossing,
(ii) $\Utilde'\setminus D\simeq \Utilde\setminus\{(0,0)\}$,
(iii) the saturation 
 $(\psi^{\ast}\Etilde_1)^{\sim}$
 of $\psi^{\ast}\Etilde_1$ in $\psi^{\ast}\Etilde$
 is a subbundle,
i.e.,
 $\psi^{\ast}\Etilde/(\psi^{\ast}\Etilde_1)^{\sim}$
 is locally free.
We have the Hermitian metric 
$\htilde_0':=\psi^{\ast}(\htilde_0)$
of $\psi^{\ast}\Etilde_1$.
The metric $\psi^{\ast}(\htilde_2)$
induces a $C^{\infty}$-metric $\htilde'_2$
of $(\psi^{\ast}\Etilde_1)^{\sim}$.
Let $s$ be the function on $\Utilde'\setminus D$
determined by
$\htilde_0'=\htilde_2'\cdot s$.
There exists a neighbourhood $N'$ of $D$
such that $s=1$ on $\Utilde'\setminus N'$.
Note that $\log s$ is an $L^2$-function on $\Utilde'$,
and we obtain the following equality
of $(1,1)$-currents on $\Utilde'$:
\[
 \delbar\del \log s
=F(\htilde_0')-F(\htilde_2')
+\sum a_i[D_i].
\]
Here, $D_i$ denote the irreducible components of $D$,
$[D_i]$ denote the $(1,1)$-current obtained as
the integrations over $D_i$,
and $a_i$ are constants.
Because $\psi_{\ast}[D_i]=0$,
we obtain
$\int_{\Utilde}\Lambda F(\htilde_{0})
=\int_{\Utilde}\Lambda F(\htilde_{2})$.
By Lemma \ref{lem;17.10.26.10},
we obtain the desired equality (\ref{eq;17.10.24.40}).
\hfill\qed

\subsection{Analytic degree of subbundles}

Let $E_2\subset E$ be a mini-holomorphic subbundle.
Let $h_{1,E_2}$ denote the metric of $E_2$
induced by $h_1$.
By the Chern-Weil formula in Lemma \ref{lem;17.10.24.41},
$\deg(E_2,h_{1,E_2})\in\real\cup\{-\infty\}$ 
makes sense.

\begin{prop}
If $\deg(E_2,h_{1,E_2})\neq-\infty$,
then there exists a good filtered subbundle
$\nbigp_{\ast}\nbige_2\subset
 \nbigp_{\ast}\nbige$ 
such that
$\nbige_{2|\nbigmlambda\setminus Z}
=E_2$.
Moreover, we obtain
$\deg(E_2,h_{1,E_2})=2\pi T\deg(\nbigp_{\ast}\nbige_2)$.
\end{prop}
\pf
We set
$\nbigm^{\lambda}\langle t_1\rangle:=
\nbigmbar^{\lambda}\langle t_1\rangle
\cap
\nbigm^{\lambda}$ for any $t_1\in S^1_T$.
By (\ref{eq;17.10.24.51}) and \cite[Lemma 10.6]{Simpson88},
$E_{2|\nbigm^{\lambda}\langle t_1\rangle\setminus Z}$
extend to locally free
$\nbigo_{\nbigmbar^{\lambda}\langle t_1\rangle\setminus Z}
 (\ast \infty)$-submodules
of $\nbigp\nbige_{|\nbigmbar^{\lambda}\langle t_1\rangle\setminus Z}$.

Let $P$ be any point of $H^{\lambda}_{\infty}$.
Let $U$ be a neighbourhood of $P$ in $\nbigmbar^{\lambda}$.
On $U$,
we use a local mini-complex coordinate system $(t_1,\beta_1^{-1})$.
On $\Utilde:=\real_{s_1}\times U$,
we use the complex coordinate system
$(\alpha_1,\beta_1^{-1})=(s_1+\sqrt{-1}t_1,\beta_1^{-1})$
as in \S\ref{subsection;17.10.3.12}.
We set $D:=\real_{s_1}\times (U\cap H^{\lambda}_{\infty})$.
As in \S\ref{subsection;13.11.29.2},
we obtain the locally free 
$\nbigo_{\Utilde}(\ast D)$-module
$\widetilde{\nbigp\nbige}$
induced by $\nbigp\nbige$.
We also have the holomorphic vector subbundle 
$\Etilde_2$ of
$\widetilde{\nbigp\nbige}_{|\Utilde\setminus D}$
induced by $E_2$.
Let $p:\Utilde\lrarr D$ be the projection
given by $p(\alpha_1,\beta_1^{-1})=\alpha_1$.
By the above consideration,
$\Etilde_{2|p^{-1}(\alpha_1)\setminus D}$
extends to
$\nbigo_{p^{-1}(\alpha_1)}(\ast\infty)$-submodule of
$\widetilde{\nbigp\nbige}_{|p^{-1}(\alpha_1)}$.
By \cite[Theorem 4.5]{Siu-extension},
$\Etilde_2$ extends to 
$\nbigp_{\Utilde}(\ast D)$-submodule
$\widetilde{\nbigp\nbige_2}$
of $\widetilde{\nbigp\nbige}$.
By the construction,
$\widetilde{\nbigp\nbige_2}$
is naturally $\real$-equivariant,
and hence
$E_{2|U\setminus H^{\lambda}_{\infty}}$ extends
to a locally free
$\nbigo_{U}
 \bigl(\ast (H^{\lambda}_{\infty}\cap U)\bigr)$-submodule
of $\nbigp\nbige_{|U}$.
Therefore, $E_{2}$ extends to a locally free 
$\nbigo_{\nbigmbar^{\lambda}\setminus Z}(\ast H^{\lambda}_{\infty})$-module
$\nbigp\nbige_2$.
We obtain the good filtered bundle
$\nbigp_{\ast}\nbige_2$
over $\nbigp\nbige_2$
as in \S\ref{subsection;17.10.24.100}
and \S\ref{subsection;18.11.21.3}.
The claim for the degree follows from the previous proposition.
\hfill\qed

\vspace{.1in}

As a consequence,
we obtain the following.
\begin{cor}
\label{cor;17.9.20.1}
$\nbigp_{\ast}\nbige$ is stable if and only if
$(E,h_1)$ is analytically stable.
\hfill\qed
\end{cor}

\section{Good filtered bundles associated with monopoles of GCK-type}
\label{subsection;20.8.1.10}

Let $Z$ be a finite subset of $\nbigm^{\lambda}$.
Let $(E,\delbar_E,h)$ be a monopole on
$\nbigm^{\lambda}\setminus Z$ of GCK-type.
Let $\nbigp^h_{\ast}E$ be the associated filtered bundle
of Dirac type on $(\nbigmbar^{\lambda};Z,H^{\lambda}_{\infty})$.

\begin{prop}
\label{prop;17.9.30.21}
The good filtered bundle $\nbigp^h_{\ast}E$ is 
polystable with $\deg(\nbigp^h_{\ast}E)=0$.
Moreover, the monopole $(E,\delbar_E,h)$ is irreducible,
if and only if $\nbigp^h_{\ast}E$ is stable.
\end{prop}
\pf
By Corollary \ref{cor;17.10.6.1},
the relation (\ref{eq;17.10.26.11}),
and Proposition \ref{prop;17.10.14.21},
$(E,\delbar_E,h)$ satisfies the condition 
in \S\ref{subsection;17.10.26.12}.
Applying Proposition \ref{prop;17.9.10.10} to $\nbigp^h_{\ast}E$,
we obtain $2\pi T\deg(\nbigp^h_{\ast}E)=\deg(E,\delbar_E,h)=0$.
Let $\nbigp_{\ast}\nbige_1$ be a good filtered subbundle of
$\nbigp^h_{\ast}E$.
Let $(E_1,\delbar_{E_1})$ be the restriction of
$\nbige_1$ to $\nbigm^{\lambda}\setminus Z$.
By Proposition \ref{prop;17.9.10.10},
we obtain
$2\pi T\deg(\nbigp_{\ast}\nbige_1)=
\deg(E_1,h_{E_1})\leq 0$.
Moreover, if $\deg(\nbigp_{\ast}\nbige_1)=0$,
the Chern-Weil formula (\ref{eq;21.8.22.40}) implies that
the orthogonal projection onto
$E_1$ is flat with respect to the Chern connection
and the Higgs fields,
and the orthogonal decomposition
$E=E_1\oplus E_1^{\bot}$ is mini-holomorphic.
Hence, we obtain the decomposition
$\nbigp^h_{\ast}E=
 \nbigp^{h_{E_1}}_{\ast}E_1\oplus
 \nbigp^{h_{E_1^{\bot}}}_{\ast}E_1^{\bot}$.
We also obtain that
$E_1$ and $E_1^{\bot}$ with the induced metrics
are monopoles.
Then, by an easy induction,
we can prove that
$\nbigp^h_{\ast}E$ is polystable.
\hfill\qed

\begin{prop}
\label{prop;17.9.30.22}
Let $h'$ be another metric of $E$ such that
(i) $(E,\delbar_E,h')$ is a monopole,
(ii) any points of $Z$ are Dirac type singularity,
(iii) $h'$ is strongly adapted to $\nbigp_{\ast}E$
in the sense of Definition {\rm\ref{df;20.8.8.3}}.
Then, the following holds.
\begin{itemize}
\item
There exists a mini-holomorphic decomposition
$(E,\delbar_E)
=\bigoplus_{i=1}^{m}
 (E_i,\delbar_{E_i})$,
which is orthogonal with respect to both
$h$ and $h'$.
\item
There exists positive numbers $a_i$ $(i=1,\ldots,m)$
such that $h_{E_i}=a_ih'_{E_i}$.
\end{itemize}
\end{prop}
\pf
By Proposition \ref{prop;17.10.14.21},
$h$ and $h'$ are mutually bounded.
Hence, we obtain the claim
from \cite[Proposition 1.3]{Mochizuki-KH-infinite}.
\hfill\qed

\section{Construction of monopoles}
\label{subsection;17.10.26.2}

Let $Z$ be a finite subset.
Let $\nbigp_{\ast}\nbige$ be a  stable good filtered bundle
of Dirac type on $(\nbigmbar^{\lambda};Z,H_{\infty}^{\lambda})$
with $\deg(\nbigp_{\ast}\nbige)=0$.
Set $E:=\nbigp_a\nbige_{|\nbigmlambda\setminus Z}$.

\begin{prop}
\label{prop;17.9.30.23}
There exists a Hermitian metric $h$ such that
(i) $(E,\delbar_E,h)$ is a monopole of GCK-type,
(ii) $(E,\delbar_E,h)$ satisfies the norm estimate
 with respect to $\nbigp_{\ast}\nbige$.
\end{prop}
\pf
By using Proposition \ref{prop;17.10.12.6}
and Proposition \ref{prop;17.10.13.150},
we can construct a metric $h_0$ of $E$
satisfying
Condition {\rm\ref{condition;21.8.22.30}}
and the following conditions.
\begin{itemize}
\item
 Let $\nabla_{h_0}$ and $\phi_{h_0}$ be
 the Chern connection and the Higgs field
 associated with $(E,\delbar_E,h_0)$.
 Then, we have the decay
 $|F(h_0)|_{h_0}\to 0$ 
 and the estimate  $|\phi_{h_0}|_{h_0}=O(\log|w|)$
 when $|w|\to\infty$.
\item
 $\det(E,\delbar_E,h_0)$ is a monopole.
\end{itemize}

Set $A^{\lambda}:=S^1\times \nbigm^{\lambda}$
and $\Abar^{\lambda}:=S^1\times\nbigmbar^{\lambda}$.
We have the complex structure
on $A^{\lambda}$ given by 
the local complex coordinate system
$(\alpha_0,\beta_0)=(s_0+\sqrt{-1}t_0,\beta_0)$
in \S\ref{subsection;17.10.4.1}.
Let $p:A^{\lambda}\lrarr \nbigm^{\lambda}$
be the projection.
Let $(\Etilde,\delbar_{\Etilde})$ be
the holomorphic bundle
on $A^{\lambda}\setminus p^{-1}(Z)$
induced by $(E,\delbar_E)$
as in \S\ref{subsection;13.11.29.2}.
It is equipped with the induced metric $\htilde_0$.
We have the natural $S^1$-action on $A^{\lambda}$,
for which $(\Etilde,\delbar_{\Etilde},\htilde_0)$
is equivariant.
We obtain that 
the $S^1$-equivariant bundle
$(\Etilde,\delbar_{\Etilde})$ with the metric $\htilde_0$
is analytically stable
with respect to the $S^1$-action
in the sense of \cite{Mochizuki-KH-infinite}.
Applying the Kobayashi-Hitchin correspondence
for analytically stable bundles in \cite{Mochizuki-KH-infinite},
we obtain a metric $h$
such that the following holds.
\begin{itemize}
\item
$(E,\delbar_E,h)$ is a monopole.
\item
$h$ and $h_0$ are mutually bounded.
In particular,
$(E,\delbar_E,h)$ satisfies the norm estimate
with respect to $\nbigp_{\ast}\nbige$.
\item
$\det(h)=\det(h_0)$.
\end{itemize}
By Proposition \ref{prop;17.10.12.5}
and Proposition \ref{prop;17.10.13.151},
we obtain that 
$(E,\delbar_E,h)$
is of GCK-type.
\hfill\qed

\vspace{.1in}
Theorem \ref{thm;17.9.30.20}
follows from 
Proposition \ref{prop;17.9.30.21},
Proposition \ref{prop;17.9.30.22}
and Proposition \ref{prop;17.9.30.23}.
\hfill\qed

\section{Smooth parabolic $0$-difference modules of rank $2$}
\label{section;21.9.17.110}

In \cite{Harland},
by using the Nahm transform
and the Kobayashi-Hitchin correspondence
between wild harmonic bundles and filtered Higgs bundles
on $\proj^1$,
Harland classified irreducible $\SU(2)$-monopoles
$(E,h,\nabla,\phi)$ on $\nbigm$
of GCK-type without Dirac type singularity
in terms of torsion-free sheaves of rank one
with parabolic structure on the spectral curves.
We revisit his result 
directly by using Corollary \ref{cor;21.9.17.60} with $\lambda=0$.

\subsection{Smooth parabolic $0$-difference modules}

\begin{df}
\index{smooth parabolic difference module}
A parabolic difference module
\[
(\vecV,V,m,
(\vectau_{x},\vecL_{x})_{x\in\cnum},
\nbigp_{\ast}(\vecV_{|\inftyhat}))
\]
is called smooth
if $m(x)=0$ for any $x\in\cnum$.
\hfill\qed
\end{df}

Let
$(\vecV,V,m,
(\vectau_{x},\vecL_{x})_{x\in\cnum},
\nbigp_{\ast}(\vecV_{|\inftyhat}))$
be a smooth parabolic $0$-difference module.
Then, the data $m$ and $(\vectau_{x},\vecL_{x})_{x\in\cnum}$ are trivial.
Because $\lambda=0$,
the difference operator of $\vecV$
is a $\cnum(\beta_1)$-linear automorphism $F$.
Indeed, it is a $\cnum[\beta_1]$-automorphism of $V$.
By setting
$V_{|\inftyhat}:=\cnum(\!(\beta_1^{-1})\!)\otimes_{\cnum[\beta_1]}V$,
we obtain a natural isomorphism
$V_{|\inftyhat}\simeq\vecV_{|\inftyhat}$.
Hence, instead of 
$(\vecV,V,m,
(\vectau_{x},\vecL_{x})_{x\in\cnum},
\nbigp_{\ast}(\vecV_{|\inftyhat}))$,
we prefer to consider
$(V,F,\nbigp_{\ast}(V_{|\inftyhat}))$.

\begin{lem}
For a smooth parabolic $0$-difference module
$(V,F,\nbigp_{\ast}(V_{|\inftyhat}))$,
$\tr(F)$ is a polynomial of $\beta_1$,
and $\det(F)$ is a non-zero constant.
\end{lem}
\pf
The first claim is obvious.
Because $\det(F)$ is the induced $\cnum[\beta_1]$-automorphism of
$\det(V)\simeq \cnum[\beta_1]$,
it is a non-zero constant.
\hfill\qed

\vspace{.1in}

For a smooth parabolic $0$-difference module
$(V,F,\nbigp_{\ast}(V_{|\inftyhat}))$,
let $\nbigf_V$ denote the $\nbigo_{\proj^1}(\ast\infty)$-module
induced by $V$.
Let $\nbigp_{\ast}\nbigf_V$ denote the filtered bundle
over $\nbigf_V$
induced by $\nbigp_{\ast}(V_{|\inftyhat})$.
\begin{lem}
We have
 $\deg(V,F,\nbigp_{\ast}(V_{|\inftyhat}))
 =\deg(\nbigp_{\ast}\nbigf_V)$.
\end{lem}
\pf
Recall the formula (\ref{eq;17.12.4.4}).
There is no contribution of
$\deg(L_{i,x},L_{i-1,x})$.
It is easy to see that
$\sum r(\omega)\cdot\omega=0$
because $\det(F)$ is a non-zero constant.
Then, we obtain
 $\deg(V,F,\nbigp_{\ast}(V_{|\inftyhat}))
 =\deg(\nbigp_{\ast}\nbigf_V)$.
\hfill\qed

\vspace{.1in}
As a consequence of
Corollary \ref{cor;21.9.17.60},
we obtain the following.
\begin{prop}
There exists the natural bijection between
the isomorphism classes of the following objects.
\begin{itemize}
 \item Monopoles of GCK-type $(E,h,\nabla,\phi)$
       defined on $\nbigm$
       without Dirac type singularity.
 \item Polystable smooth parabolic $0$-difference
       modules
       $(V,F,\nbigp_{\ast}(V_{|\inftyhat}))$
       such that 
       $\deg(\nbigp_{\ast}\nbigf_V)=0$.
\end{itemize}
Under this equivalence,
$(E,h,\nabla,\phi)$ is an $\SU(n)$-monopole
if and only if $\det(F)=1$ and $\rank(V)=n$ hold
for the corresponding $(V,F,\nbigp_{\ast}(V_{|\inftyhat}))$.
\hfill\qed
\end{prop}

We shall describe
smooth polystable parabolic $0$-difference modules
$(V,F,\nbigp_{\ast}(V_{|\inftyhat}))$ of degree $0$
with $\rank(V)\leq 2$.

\subsection{Rank one case}
\label{subsection;21.9.18.10}

For any $\gamma\in\cnum^{\ast}$,
we set $N(\gamma)=\cnum[\beta_1]$.
It is equipped with the $\cnum[\beta_1]$-automorphism
determined by the multiplication of $\gamma$.
Under the natural identification
$N(\gamma)_{|\inftyhat}
=\cnum(\!(\beta_1^{-1})\!)\otimes_{\cnum[\beta_1]}N(\gamma)
\simeq
 \cnum(\!(\beta_1^{-1})\!)$,
let $\nbigp_{\ast}(N(\gamma)_{|\inftyhat})$
denote the filtered bundle over
$N(\gamma)_{|\inftyhat}$
defined by
$\nbigp_a(N(\gamma)_{|\inftyhat})
=\beta_1^{[a]}\cnum[\![\beta_1^{-1}]\!]$,
where $[a]:=\max\{n\in\seisuu\,|\,n\leq a\}$.
Let $\nbigp_{\ast}N(\gamma)$ denote the parabolic
$0$-difference module of rank $1$.
Note that $\deg(\nbigp_{\ast}N(\gamma))=0$.
The following lemma is clear.
\begin{lem}
Any smooth parabolic difference module of degree $0$
of rank $1$ 
is isomorphic to $\nbigp_{\ast}N(\gamma)$
for some $\gamma$, 
\hfill\qed
\end{lem}

\subsection{Filtered torsion-free sheaves of rank one on an integral scheme}

Let $Q\in\cnum[\beta_1]$ be a non-constant polynomial,
i.e., $\deg(Q)>0$.
Let $\gamma\in\cnum^{\ast}$.
Let $\Sigma^{\circ}(Q,\gamma)\subset \cnum_x^{\ast}\times\cnum_{\beta_1}$
denote the scheme defined by $x^2-Q(\beta_1)x+\gamma$.

\begin{lem}
If $Q(\beta_1)$ is not constant,
then the polynomial $x^2-Q(\beta_1)x+\gamma\in\cnum[\beta_1][x]$
is irreducible.
As a result, $\Sigma^{\circ}(Q,\gamma)$ is integral.
\end{lem}
\pf
As remarked in \cite{Harland}, it
follows from Eisenstein's criterion.
We can also obtain the claim
by the observation that the discriminant $Q(\beta_1)^2-4\gamma$ of
$x^2-Q(\beta_1)x+\gamma$ is not a square of a polynomial.
\hfill\qed

\vspace{.1in}

We set $\vecP=\proj\bigl(\nbigo_{\proj_{\beta_1}}\oplus
\nbigo_{\proj^1_{\beta_1}}(-k\{\infty\})
\bigr)$,
where $k=\deg(Q)$.
We regard $Q(\beta_1)$ and $\gamma$
as sections of
$\nbigo_{\proj^1_{\beta_1}}(k\{\infty\})$
and
$\nbigo_{\proj^1_{\beta_1}}(2k\{\infty\})$,
respectively.
We obtain the subscheme
$\Sigma(Q,\gamma)\subset\vecP$
defined by $x^2-Qx+\gamma$.
(See \cite[\S3]{Beauville-Narasimhan-Ramanan}.)
By naturally regarding
$\cnum^{\ast}_x\times\cnum_{\beta_1}$
as an open subset of $\vecP$,
we may regard $\Sigma(Q,\gamma)$ as the closure of
$\Sigma^{\circ}(Q,\gamma)$ in $\vecP$.
Let $p_{\Sigma(Q,\gamma)}:\Sigma(Q,\gamma)\lrarr\proj^1_{\beta_1}$
denote the projection.
We set
$\Sigma_{\infty}(Q,\gamma):=p_{\Sigma(Q,\gamma)}^{-1}(\infty)$.
In the following, 
$\Sigma(Q,\gamma)$
and $\Sigma_{\infty}(Q,\gamma)$
are also denoted by $\Sigma$
and $\Sigma_{\infty}$,
respectively, if there is no risk of confusion.

\begin{lem}
$\Sigma_{\infty}(Q,\gamma)$ consists of two distinct points,
and they are smooth point of $\Sigma(Q,\gamma)$.
As a result, $\Sigma(Q,\gamma)$ is integral.
On a neighbourhood of each point of $\Sigma_{\infty}(Q,\gamma)$,
$p_{\Sigma(Q,\gamma)}$ is etale.
\end{lem}
\pf
Let $(1-4\gamma Q^{-2})^{1/2} \in\cnum[\![y^{-1}]\!]$
denote the convergent power series
such that
(i) the constant term is $1$,
(ii) $\bigl(
(1-4\gamma Q^{-2})^{1/2}\bigr)^2=1-4\gamma Q^{-2}$.
One of the roots of $x^2-Qx+1$ is
\[
 \gminia_1=\frac{Q+Q(1-4\gamma Q^{-2})^{1/2}}{2}
 \in y^{\deg(Q)}\cnum[\![y^{-1}]\!].
\]
The other is
$\gminia_2=\gamma\gminia^{-1}_1\in
 y^{-\deg(Q)}\cnum[\![y^{-1}]\!]$.
Note 
$\lim_{y^{-1}\to 0}y^{-\deg(Q)}\gminia_1\neq 0$
and
$\lim_{y^{-1}\to 0}y^{\deg(Q)}\gminia_2\neq 0$.
Then, the claim of the lemma is obvious.
\hfill\qed

\subsubsection{Degree of torsion-free sheaves on $\Sigma(Q,\gamma)$}

Let us recall the degree of coherent sheaves on
the integral scheme $\Sigma(Q,\gamma)$.
(See \cite{Huybrechts-Lehn} for a more precise explanation
in the more general situation.)
Let $\gbigp$ be any ample line bundle on $\Sigma(Q,\gamma)$.
For any coherent sheaf $\nbige$ on $\Sigma(Q,\gamma)$,
there exist rational numbers
$\alpha_i(\nbige)$ $(i=0,1)$
such that
$\chi(\nbige\otimes\gbigp^m)=m\alpha_1(\nbige)+\alpha_0(\nbige)$
for any $m\in\seisuu$,
where $\chi(\nbige\otimes\gbigp^m)$
denotes the Euler number of $\nbige\otimes\gbigp^m$
(see \cite[\S1.2]{Huybrechts-Lehn}).
According to
\cite[Definition 1.2.2, Definition 1.2.11]{Huybrechts-Lehn},
we set
\[
 \rank(\nbige):=\frac{\alpha_1(\nbige)}{\alpha_1(\nbigo_{\Sigma(Q,\gamma)})},
 \quad
 \deg(\nbige):=\alpha_0(\nbige)-\rank(\nbige)\alpha_0(\nbigo_{\Sigma(Q,\gamma)}).
\]
Because $\Sigma(Q,\gamma)$ is integral,
there exists a Zariski open dense connected subset $U$ of $\Sigma(Q,\gamma)$
such that $\nbige_{|U}$ is a locally free $\nbigo_U$-module,
and $\rank(\nbige)$ is equal to the rank of $\nbige_{|U}$.
Because $\dim \Sigma(Q,\gamma)=1$,
it is easy to see that $\deg(\nbige)$
is independent of the choice of
an ample line bundle $\gbigp$
though $\alpha_i(\nbige)$ depend on $\gbigp$.
If $\Sigma(Q,\gamma)$ is smooth,
we have $\deg(\nbige)=\int_{\Sigma(Q,\gamma)}c_1(\nbige)$
for the first Chern class $c_1(\nbige)$ in the ordinary sense.

\begin{lem}
\label{lem;21.9.17.120}
Let $\nbige$ be a torsion-free sheaf of rank $1$ on $\Sigma(Q,\gamma)$.
Then, we have
$\deg(p_{\Sigma(Q,\gamma)\ast}\nbige)=\deg(\nbige)-\deg(Q)$.
\end{lem}
\pf
We consider the ample line bundles
$\nbigo_{\proj^1_{\beta_1}}(1)$ on $\proj^1_{\beta_1}$
and 
$\gbigp:=p_{\Sigma}^{\ast}(\nbigo_{\proj^1_{\beta}}(1))$
on $\Sigma$.
We set $\nbigf=p_{\Sigma\ast}\nbige$,
which is a locally free $\nbigo_{\proj_{\beta_1}^1}$-module.
Because
$\chi(\nbigf\otimes\nbigo(m))
=
\chi(\nbige\otimes\gbigp^m)$,
we have
$\alpha_i(\nbigf)
=\alpha_i(\nbige)$.
Hence, we obtain
\[
 \deg(\nbigf)
=\alpha_0(\nbigf)-\rank(\nbigf)\alpha_0(\nbigo_{\proj^1})
=\alpha_0(\nbige)
-2\alpha_0(\nbigo_{\proj^1}).
\]
As explained in \cite{Beauville-Narasimhan-Ramanan},
we have
$p_{\Sigma\ast}(\nbigo_{\Sigma})
\simeq
\nbigo_{\proj^1}\oplus\nbigo_{\proj^1}(-k)$,
where $k=\deg(Q)$.
Hence, we obtain
\[
 \chi(\nbigo_{\Sigma}\otimes\gbigp^m)
 =\chi(\nbigo_{\proj^1}(m))
 +\chi(\nbigo_{\proj^1}(m-k)).
\]
It implies that
$\alpha_0(\nbigo_{\Sigma})
=2\alpha_0(\nbigo_{\proj^1})-k$.
Hence, we obtain
\[
 \deg(\nbige)
=\alpha_0(\nbige)
-\alpha_0(\nbigo_{\Sigma})
=\alpha_0(\nbige)
-2\alpha_0(\nbigo_{\proj^1})+k \\
=\deg(\nbigf)+k. 
\]
Thus, we are done.
\hfill\qed

\subsubsection{Filtered torsion-free sheaf of rank one
on $(\Sigma(Q,\gamma),\Sigma_{\infty}(Q,\gamma))$}

Let $\nbigl$ be a torsion-free
$\nbigo_{\Sigma(Q,\gamma)}(\ast\Sigma_{\infty}(Q,\gamma))$-module
of rank $1$.
\begin{df}
\label{filtered torsion-free sheaf of rank one}
A filtered sheaf $\nbigp_{\ast}\nbigl$
over $\nbigl$ is 
a tuple of $\nbigo_{\Sigma(Q,\gamma)}$-submodules
$\nbigp_{\vecb}(\nbigl)\subset\nbigl$
 $(\vecb=(b_P\,|\,P\in \Sigma_{\infty}(Q,\gamma))
 \in\real^{\Sigma_{\infty}(Q,\gamma)})$
such that the following holds.
\begin{itemize}
 \item For $P\in \Sigma_{\infty}(Q,\gamma)$,
       let $U_P$ be a neighbourhood of $P$
       in $\Sigma(Q,\gamma)$
       such that
       (i) $U_P\cap\Sigma_{\infty}(Q,\gamma)=\{P\}$,
       (ii) $U_P$ is smooth.
       Then, $\nbigp_{\vecb}(\nbigl)_{|U_P}$
       depends only on $b_P$,
       which we denote by
       $\nbigp_{b_P}(\nbigl_{|U_P})$.
       Moreover,
       $\nbigp_{\ast}(\nbigl_{|U_P})$
       is a filtered bundle
       over $\nbigl_{|U_P}$.
       (See {\rm\S\ref{subsection;21.8.15.2}}.)
\end{itemize}
Such $\nbigp_{\ast}\nbigl$ is called a filtered torsion-free
sheaf of rank one on $(\Sigma(Q,\gamma),\Sigma_{\infty}(Q,\gamma))$.
\hfill\qed
\end{df}

Let $\nbigp_{\ast}\nbigl$ be a filtered torsion-free sheaf
of rank one on $(\Sigma(Q,\gamma),\Sigma_{\infty}(Q,\gamma))$.
For $P\in\Sigma_{\infty}(Q,\gamma)$,
let $\nbigl_P$ denote the stalk of $\nbigl$ at $P$.
We have the induced filtration
$\nbigp_{\ast}(\nbigl_P)$.
We obtain the $\cnum$-vector spaces
$\Gr^{\nbigp}_c(\nbigl_P)=
\nbigp_c(\nbigl_P)\big/
\sum_{d<c}\nbigp_{d}(\nbigl_P)$.
We set
\[
 \deg(\nbigp_{\ast}\nbigl)
=\deg(\nbigp_{\vecb}\nbigl)
-\sum_{P\in \Sigma_{\infty}}
 \sum_{b_P-1<c\leq b_P}
 c\dim \Gr^{\nbigp}_c(\nbigl_P).
\]
It is independent of the choice of $\vecb$.

\subsubsection{The induced smooth parabolic difference modules
of rank $2$}
\label{subsection;21.9.18.1}

In this subsection,
$(\Sigma(Q,\gamma),\Sigma_{\infty}(Q,\gamma))$
is denoted by
$(\Sigma,\Sigma_{\infty})$.
Let $\nbigp_{\ast}\nbigl$ be a filtered
torsion-free sheaf of rank one on
$(\Sigma,\Sigma_{\infty})$.
Let us construct
a parabolic $0$-difference module
$\Upsilon(\nbigp_{\ast}\nbigl)$.

We obtain
a locally free $\nbigo_{\proj^1_{\beta_1}}(\ast\infty)$-module
$p_{\Sigma\ast}(\nbigl)$ of rank $2$.
There exists the following isomorphism of the stalks:
\[
 p_{\Sigma\ast}(\nbigl)_{\infty}
 =\bigoplus_{P\in\Sigma_{\infty}}
 \nbigl_P.
\]
We define the filtration
$\nbigp_{\ast}\bigl(p_{\Sigma\ast}(\nbigl)_{\infty}
\bigr)$
by
\[
\nbigp_{a}\bigl(p_{\Sigma\ast}(\nbigl)_{\infty}\bigr)
=\bigoplus_{P\in\Sigma_{\infty}}
 \nbigp_a(\nbigl_P).
\]
Let $\nbigp_a\bigl(p_{\Sigma\ast}\nbigl\bigr)\subset
p_{\Sigma\ast}(\nbigl)$
be determined by 
$\nbigp_a\bigl(p_{\Sigma\ast}\nbigl\bigr)_{\infty}
=\nbigp_a\bigl(p_{\Sigma\ast}(\nbigl)_{\infty}\bigr)$.
Thus, we obtain a filtered bundle
$\nbigp_{\ast}\bigl(p_{\Sigma\ast}(\nbigl)\bigr)$
over $p_{\Sigma\ast}(\nbigl)$.

\begin{lem}
We have
 $\deg(\nbigp_{\ast}(p_{\Sigma\ast}\nbigl))
 =\deg(\nbigp_{\ast}\nbigl)-\deg(Q)$.
\end{lem}
\pf
By the construction,
we have
$\nbigp_a(p_{\Sigma\ast}\nbigl)
=p_{\Sigma\ast}(\nbigp_{a,a}\nbigl)$ for any $a\in\real$.
Then, the claim follows from
Lemma \ref{lem;21.9.17.120}
and the definitions of 
 $\deg(\nbigp_{\ast}(p_{\Sigma\ast}\nbigl))$
and $\deg(\nbigp_{\ast}\nbigl)$.
\hfill\qed

\vspace{.1in}
As explained in \cite[\S3]{Beauville-Narasimhan-Ramanan},
there exists a natural morphism
$F:\nbigp_a(p_{\Sigma\ast}(\nbigl))\lrarr
\nbigp_a(p_{\Sigma\ast}(\nbigl))\otimes\nbigo_{\proj^1}(\deg(Q))$
for any $a\in\real$.
In particular,
we obtain an automorphism
$F$ of $p_{\Sigma\ast}(\nbigl)$.
It induces the automorphism
$F$ of the $\cnum(\!(\beta_1^{-1})\!)$-vector space
$p_{\Sigma\ast}(\nbigl)_{|\inftyhat}$.
It is easy to see that
the decomposition
\begin{equation}
\label{eq;21.9.17.140}
 p_{\Sigma\ast}(\nbigl)_{|\inftyhat}
=\bigoplus_{P\in\Sigma_{\infty}}
\nbigl_{|\widehat{P}} 
\end{equation}
is the eigen decomposition of $F$.

We set
$V:=H^0(\proj^1,p_{\Sigma\ast}(\nbigl))$.
Because
$V_{|\inftyhat}\simeq
p_{\Sigma\ast}(\nbigl)_{|\widehat{\infty}}$,
we obtain a filtered bundle
$\nbigp_{\ast}(V_{|\inftyhat})$
over $V_{|\inftyhat}$.
Because the filtration is compatible
with the eigen decomposition (\ref{eq;21.9.17.140}) of $F$,
$\nbigp_{\ast}(V_{|\inftyhat})$ is a good parabolic structure
of $(V,F)$ at $\infty$.
Thus, we obtain a smooth parabolic difference module
$\Upsilon(\nbigp_{\ast}\nbigl)=
(V,F,\nbigp_{\ast}(V_{|\inftyhat}))$.
If $\deg(\nbigp_{\ast}\nbigl)=\deg(Q)$,
then we obtain
$\deg(\Upsilon(\nbigp_{\ast}\nbigl))=0$.
The following lemma is obvious by the construction.
\begin{lem}
\label{lem;21.9.18.2}
Let $\nbigp_{\ast}\nbigl_i$ $(i=1,2)$
be filtered torsion-free sheaf of rank one
on $(\Sigma,\Sigma_{\infty})$.
Isomorphisms
of smooth parabolic difference modules
$\Upsilon(\nbigp_{\ast}\nbigl_1)
\simeq
\Upsilon(\nbigp_{\ast}\nbigl_2)$ 
bijectively correspond to
isomorphisms
of filtered torsion-free sheaves
$\nbigp_{\ast}\nbigl_1\simeq\nbigp_{\ast}\nbigl_2$.
\hfill\qed
\end{lem}

We obtain the following lemma
because $\Sigma$ is irreducible.
\begin{lem}
\label{lem;21.9.18.3}
For any filtered torsion-free sheaf of rank one $\nbigp_{\ast}\nbigl$,
the induced smooth parabolic $0$-difference module
$\Upsilon(\nbigp_{\ast}\nbigl)$ is stable. 
\hfill\qed
\end{lem}

\subsection{Description of parabolic difference modules of rank $2$}

Let $Q\in\cnum[\beta_1]$ and $\gamma\in\cnum^{\ast}$.
Let $\gbigs(Q,\gamma,c)$
denote the set of the isomorphism classes of
smooth parabolic difference module
$(V,F,\nbigp_{\ast}(V_{|\inftyhat}))$
such that $\rank(V)=2$ and $\deg(\nbigp_{\ast}\nbigf_V)=c$.
Let $\gbigs^{s}(Q,\gamma,c)\subset\gbigs(Q,\gamma,c)$
be the subset of the isomorphism classes of
stable objects.
Let $\gbigs^{ps}(Q,\gamma,c)\subset\gbigs(Q,\gamma,c)$
be the subset of the isomorphism classes of
polystable objects.

\subsubsection{} 

Suppose that $Q$ is non-constant.
Let $\gbigl(Q,\gamma,c)$ denote the set of
the isomorphism classes of torsion-free sheaf of rank one
$\nbigp_{\ast}\nbigl$ on $(\Sigma(Q,\gamma),\Sigma_{\infty}(Q,\gamma))$
such that $\deg(\nbigp_{\ast}\nbigl)=c$.
In \S\ref{subsection;21.9.18.1},
we constructed
$\Upsilon:\gbigl(Q,\gamma,c+\deg(Q))\lrarr \gbigs(Q,\gamma,c)$.
According to Lemma \ref{lem;21.9.18.3},
the image of $\Upsilon$ is contained in $\gbigs^{s}(Q,\gamma,c)$.
According to Lemma \ref{lem;21.9.18.2},
$\Upsilon$ is injective.
The following proposition is a variant of
\cite[Proposition 3.6]{Beauville-Narasimhan-Ramanan}.
The idea goes back to \cite{Hitchin-self-duality}.

\begin{prop}
\label{prop;21.9.18.10}
We have $\gbigs^s(Q,\gamma,c)=\gbigs(Q,\gamma,c)$,
and $\Upsilon$ is a bijection. 
\end{prop}
\pf
Because $\Sigma(Q,\gamma)$ is irreducible,
it is easy to see
$\gbigs^s(Q,\gamma,c)=\gbigs(Q,\gamma)$.

Let $(V,F,\nbigp_{\ast}(V_{|\inftyhat}))$
be a smooth parabolic $0$-difference module.
According to \cite[Proposition 3.6]{Beauville-Narasimhan-Ramanan},
the $\nbigo_{\proj^1_{\beta_1}}(\ast\infty)$-module $\nbigf_V$
equipped with an automorphism $F$
corresponds to
a torsion-free
$\nbigo_{\Sigma(Q,\gamma)}(\ast \Sigma_{\infty}(Q,\gamma))$-module
$\nbigl$ of rank $1$,
where $\nbigl$ and $\nbigf_V$ are related as
$p_{\Sigma\ast}(\nbigl)=\nbigf_V$.
Note that there exists an eigen decomposition
$V_{|\inftyhat}=\EE_{\gminia_1}\oplus\EE_{\gminia_2}$,
which is identified with the following isomorphism
induced by $p_{\Sigma\ast}(\nbigl)=\nbigf_V$:
\begin{equation}
\label{eq;21.9.17.72}
 V_{|\inftyhat}=\bigoplus_{P\in \Sigma_{\infty}(Q,\gamma)}
  \nbigl_{|\widehat{P}}.
\end{equation}

\begin{lem}
A filtered bundle $\nbigp_{\ast}(V_{|\inftyhat})$
over $V_{|\inftyhat}$
is a good parabolic structure of
the $0$-difference module $(V,F)$
if and only if
$\nbigp_{\ast}(V_{|\inftyhat})$
is induced by 
filtered bundles
$\nbigp_{\ast}(\nbigl_{|\widehat{P}})$
over $\nbigl_{|\widehat{P}}$ 
under the isomorphism {\rm(\ref{eq;21.9.17.72})}. 
\end{lem}
\pf
Suppose that $\nbigp_{\ast}(V_{|\inftyhat})$ is
a good parabolic structure at infinity.
Because good filtration is compatible
with the slope decomposition (see Lemma \ref{lem;21.9.17.70}),
we have the following for any $a\in\real$:
\[
 \nbigp_{a}(V_{|\inftyhat})
=\bigl(
 \nbigp_{a}(V_{|\inftyhat})\cap \EE_{\gminia_1}
 \bigr)
 \oplus
\bigl(
 \nbigp_{a}(V_{|\inftyhat})\cap \EE_{\gminia_2}
 \bigr).
\]
It means that $\nbigp_{\ast}(V_{|\inftyhat})$
is induced by
the filtrations
$\nbigl_{|\widehat{P}}$ $(P\in \Sigma_{\infty}(Q,\gamma))$.
We can check the converse similarly.
Indeed, we have already observed it
in \S\ref{subsection;21.9.18.1}.
\hfill\qed

\vspace{.1in}
Let $\nbigp_{\ast}(\nbigl_{|\widehat{P}})$ $(P\in \Sigma_{\infty}(Q,\gamma))$
be filtered bundles over $\nbigl_{|\widehat{P}}$
corresponding to $\nbigp_{\ast}(V_{|\inftyhat})$.
For $\vecb=(b_P)\in\real^{\Sigma_{\infty}(Q,\gamma)}$,
there exist the torsion-free $\nbigo_{\Sigma(Q,\gamma)}$-submodule
 $\nbigp_{\vecb}\nbigl\subset\nbigl$
such that
$\nbigp_{\vecb}(\nbigl)_{|\widehat{P}}
=\nbigp_{b_P}(\nbigl_{|\widehat{P}})$.
Thus, we obtain a filtered torsion-free sheaf of rank one
$\nbigp_{\ast}\nbigl$.
By the construction, we have
$\Upsilon(\nbigp_{\ast}\nbigl)=(V,F,\nbigp_{\ast}(V_{|\inftyhat}))$.
Hence, we obtain that $\Upsilon$ is surjective.
Because $\Upsilon$ is injective,
we obtain Proposition \ref{prop;21.9.18.10}.
\hfill\qed

\subsubsection{}

Suppose that $Q$ is constant.

\begin{prop}
$\gbigs^{s}(Q,\gamma,0)$ is empty,
and $\gbigs^{ps}(Q,\gamma,0)$
consists of
the isomorphism classes of
$\nbigp_{\ast}N(\alpha_1)\oplus\nbigp_{\ast}N(\alpha_2)$,
where
$(\alpha_1,\alpha_2)$ is determined by
$x^2-Qx+\gamma=(x-\alpha_1)(x-\alpha_2)$.
(See {\rm\S\ref{subsection;21.9.18.10}}
for $\nbigp_{\ast}N(\alpha)$.)
\end{prop}
\pf
Let $(V,F,\nbigp_{\ast}(V_{|\inftyhat}))\in \gbigs^{ps}(Q,\gamma,0)$.
If $\alpha_1\neq\alpha_2$,
there exists a decomposition
$(V,F)=(V_1,F_1)\oplus (V_2,F_2)$
such that $F_i$ are the multiplication of $\alpha_i$.
By Lemma \ref{lem;21.9.17.100},
the decomposition is compatible with
the parabolic structure at infinity,
i.e.,
$\nbigp_{\ast}(V_{|\infty})
=\nbigp_{\ast}(V_{1|\infty})
\oplus
\nbigp_{\ast}(V_{2|\infty})$.
Because $(V,F,\nbigp_{\ast}(V_{|\inftyhat}))$
is assumed to be polystable of degree $0$,
we obtain that
$\deg(V_{i},F_i,\nbigp_{\ast}(V_{i|\inftyhat}))=0$.
Thus, we obtain the claim of the lemma
in the case $\alpha_1\neq\alpha_2$.

Let us study the case $\alpha_1=\alpha_2=:\alpha$.
Note that $F$ preserves $\nbigp_{\ast}(V_{|\inftyhat})$.
Hence, $F-\alpha\id$ also preserves
$\nbigp_{\ast}(V_{|\inftyhat})$.
Suppose that $F-\alpha\id\neq 0$.
Because $\rank V=2$,
we have $(F-\alpha\id)^2=0$.
Hence, there exists
a free $\cnum[\beta_1]$-submodule $V_0\subset V$
such that
(i) $\Image(F-\alpha_0)\subset V_0=\Ker(F-\alpha\id)$,
(ii) $V_1=V/V_0$ is also a free $\cnum[\beta_1]$-module.
We obtain the induced parabolic $0$-difference modules
$(V_i,F_i,\nbigp_{\ast}(V_{i|\inftyhat}))$.
Because $(V,F,\nbigp_{\ast}(V_{|\inftyhat}))$
is polystable of degree $0$,
we obtain
\[
\mu(V_0,F_0,\nbigp_{\ast}(V_{0|\inftyhat}))
\leq
\mu(V,F,\nbigp_{\ast}(V_{|\inftyhat}))=0
\leq
\mu(V_1,F_1,\nbigp_{\ast}(V_{1|\inftyhat})).
\]
Because $F-\alpha\id$ induces a non-trivial morphism
\[
 (V_1,F_1,\nbigp_{\ast}(V_{1|\inftyhat}))
\lrarr
(V_0,F_0,\nbigp_{\ast}(V_{0|\inftyhat})),
\]
and because $\rank V_i=1$,
we obtain
$\mu(V_1,F_1,\nbigp_{\ast}(V_{1|\inftyhat}))
\leq
\mu(V_0,F_0,\nbigp_{\ast}(V_{0|\inftyhat}))$.
Thus, we obtain
$\mu(V_0,F_0,\nbigp_{\ast}(V_{0|\inftyhat}))
=
\mu(V_1,F_1,\nbigp_{\ast}(V_{1|\inftyhat}))=0$.
Because $(V,F,\nbigp_{\ast}(V_{|\inftyhat}))$
is assumed to be polystable,
we obtain
$(V,F,\nbigp_{\ast}(V_{|\inftyhat}))
\simeq
(V_0,F_0,\nbigp_{\ast}(V_{0|\inftyhat}))
\oplus
(V_1,F_1,\nbigp_{\ast}(V_{1|\inftyhat}))$,
which implies $F-\alpha\id=0$.
But, it contradicts with the assumption.
Thus, we obtain $F-\alpha\id=0$.
Then, the filtered bundle
$\nbigp_{\ast}\nbigf_V$
induced by $\nbigf_V$ and $\nbigp_{\ast}(V_{|\inftyhat})$
is a polystable bundle of degree $0$
on $(\proj^1,\infty)$ in the ordinary sense.
It is well known such a filtered bundle is isomorphic to
a direct sum of obvious filtered bundles of rank one
on $(\proj^1,\infty)$.
\hfill\qed

\chapter{Appendix}

In \S\ref{section;21.7.7.20},
we shall explain a formal $\lambda$-connection
induced by a Higgs bundle with a Hermitian metric
which asymptotically satisfies the Hitchin equation.
It is a complement to \cite[\S5.5]{Mochizuki-doubly-periodic}.
In \S\ref{subsection;20.8.1.1}--\ref{subsection;20.8.1.2},
the estimates in \cite{Mochizuki-doubly-periodic}
for doubly-periodic instantons
are generalized to the context of
asymptotically doubly-periodic instantons.

\section[Formal $\lambda$-connections]{Formal $\lambda$-connections
 associated with asymptotic harmonic bundles}
\label{section;21.7.7.20}
\subsection{Asymptotic harmonic bundles}
\label{subsection;17.10.21.10}
\index{asymptotic harmonic bundle}

Let $U$ be a neighbourhood of $0$ in $\cnum$.
Let $(V,\delbar_V,\theta)$ be a Higgs bundle
on $U\setminus\{0\}$.
We assume that 
$(V,\delbar_V,\theta)$ is unramifiedly good wild,
\index{unramifiedly good wild}
i.e., 
after shrinking $U$,
there exists a finite subset
$\nbigi\subset z^{-1}\cnum[z^{-1}]$
and a decomposition
\begin{equation}
 \label{eq;17.10.20.30}
 (V,\delbar_V,\theta)
=\bigoplus_{\gminia\in\nbigi}
 (V_{\gminia},\delbar_{V_{\gminia}},\theta_{\gminia})
\end{equation}
such that the following holds for each $\gminia\in\nbigi$.
\begin{itemize}
\item
Let $f_{\gminia}$ be the endomorphism of $V_{\gminia}$
determined by
$\theta_{\gminia}-d\gminia\id_{V_{\gminia}}=f_{\gminia}dz/z$.
Let $\det(t\id_{V_{\gminia}}-f_{\gminia})=\sum a_{\gminia,j}(z)t^j$ 
be the characteristic polynomial of $f_{\gminia}$.
Then, $a_{\gminia,j}$ are holomorphic on $U$.
\end{itemize}

For any $\gminia\in z^{-1}\cnum[z^{-1}]$,
set $\ord(\gminia):=-\deg_{z^{-1}}\gminia$.
Let $p$ be a positive integer such that
\begin{equation}
\label{eq;17.10.21.1}
-p\leq\min\{\ord(\gminia-\gminib)\,|
 \,\gminia,\gminib\in\nbigi,\gminia\neq\gminib\}. 
\end{equation}
(Note that ``$p$'' is denoted as ``$-p$''
in \cite[\S5.5]{Mochizuki-doubly-periodic}.)

Let $h$ be a Hermitian metric of $V$.
Let $F(h)$ denote the curvature of the Chern connection 
$\nabla$ of $(V,\delbar_V,h)$.
Let $\theta^{\dagger}=f^{\dagger}(d\zbar/\zbar)$
denote the adjoint of $\theta$
with respect to $h$.
We assume the following condition.
\begin{itemize}
\item
Let $P(s_1,s_2,s_3,s_4)$ be a polynomial 
of non-commutative variables $s_i$ $(i=1,2,3,4)$.
Then, there exists $\epsilon(P)>0$ such that 
the following estimate holds:
\begin{equation}
 \label{eq;17.10.21.2}
 P\bigl(\ad(f),\ad(f^{\dagger}),\nabla_z,\nabla_{\zbar}\bigr)
 \bigl(
 F(h)+[\theta,\theta^{\dagger}]
 \bigr)=
 O\Bigl(
 \exp\bigl(-\epsilon(P)|z|^{-p}\bigr)
 \Bigr)\,dz\,d\zbar.
\end{equation}
\end{itemize}

\begin{rem}
The above condition is not the same as
that in {\rm\cite[\S5.5]{Mochizuki-doubly-periodic}}.
Although a stronger assumption 
$-p<\min\{\ord(\gminia-\gminib)\,|
 \,\gminia,\gminib\in\nbigi,\gminia\neq\gminib\}$
is imposed in {\rm \cite{Mochizuki-doubly-periodic}},
the condition {\rm(\ref{eq;17.10.21.1})} is enough
for our purpose.
While we considered only the estimate
for $F(h)+[\theta,\theta^{\dagger}]$ 
in  {\rm\cite{Mochizuki-doubly-periodic}},
we impose the decay condition for the higher derivatives in this paper,
which will be used in
Proposition {\rm\ref{prop;17.10.6.201}}.
\hfill\qed
\end{rem}

By shrinking $U$,
we may assume that
there exists the refined decomposition
\begin{equation}
\label{eq;20.7.27.1}
 (V_{\gminia},\delbar_{V_{\gminia}},\theta_{\gminia})
=\bigoplus_{\alpha\in\cnum}
 \bigl(V_{\gminia,\alpha},
  \delbar_{V_{\gminia,\alpha}},\theta_{\gminia,\alpha}\bigr)
\end{equation}
such that the following holds.
\begin{itemize}
\item
Let $f_{\gminia,\alpha}$ be the endomorphism of $V_{\gminia,\alpha}$
determined by
$\theta_{\gminia,\alpha}-(d\gminia+\alpha dz/z)\id_{V_{\gminia,\alpha}}
=f_{\gminia,\alpha}dz/z$.
We consider the characteristic polynomial
$\det(t\id_{V_{\gminia,\alpha}}-f_{\gminia,\alpha})=
 t^{\rank V_{\gminia,\alpha}}
+\sum_{j=0}^{\rank V_{\gminia,\alpha}-1} a_{\gminia,\alpha,j}(z)t^j$.
Then, $a_{\gminia,\alpha,,j}$ are holomorphic on $U$,
and $a_{\gminia,\alpha,j}(0)=0$
for any $0\leq j\leq \rank V_{\gminia,\alpha}-1$.
\end{itemize}

\subsection{Simpson's main estimate}
\index{Simpson's main estimate}

The following proposition is
\cite[Proposition 5.18]{Mochizuki-doubly-periodic},
which is a variant of Simpson's main estimate
in \cite{Simpson88, Mochizuki-wild}.
\begin{prop}
\label{prop;12.8.7.3}
\mbox{{}}\index{asymptotically orthogonal}
\begin{itemize}
\item
If $\gminia\neq\gminib$,
there exists $\epsilon>0$ such that
$V_{\gminia,\alpha}$ and $V_{\gminib,\beta}$
are $O\bigl(
 \exp(-\epsilon|z|^{\ord(\gminia-\gminib)})\bigr)$-asymptotically
orthogonal,
i.e., there exists $C_1>0$ such that
the following holds
for any $u\in V_{\gminia,\alpha|Q}$
and $v\in V_{\gminib,\beta|Q}$:
\[
 \bigl|h(u,v)\bigr|
 \leq C_1\exp\bigl(-\epsilon|z(Q)|^{\ord(\gminia-\gminib)}\bigr)
 \cdot |u|_h\cdot|v|_h.
\]
\item
If $\alpha\neq\beta$,
there exists $\epsilon>0$ such that 
$V_{\gminia,\alpha}$
and $V_{\gminia,\beta}$ 
are $O(|z|^{\epsilon})$-asymptotically orthogonal. 
\item
 $\theta_{\gminia,\alpha}
 -(d\gminia+\alpha\,dz/z)\,\id_{E_{\gminia,\alpha}}$
 is bounded with respect to
 $h$ and the Poincar\'e metric $g_{\Poin}$.
\end{itemize}
\end{prop}
\pf
Because the condition is slightly changed,
we repeat the argument in 
the proof of \cite[Proposition 5.18]{Mochizuki-doubly-periodic}
with minor modifications.
By considering the tensor product
with a harmonic bundle of a rank one,
we may assume 
$-p\leq \min\bigl\{
 \ord(\gminia)\,\big|\,
 \gminia\in\nbigi
 \bigr\}$.
For any $\ell\leq p$,
we define the map
$\eta_{\ell}:
 z^{-1}\cnum[z^{-1}]\lrarr 
 z^{-\ell}\cnum[z^{-1}]$
 by
$ \eta_{\ell}\Bigl(
 \sum \gminia_jz^{j}
  \Bigr)
  =\sum_{j\leq -\ell}\gminia_jz^j$.
Let $\nbigi_{\ell}$ denote the image of $\nbigi$.
For each $\gminib\in \nbigi_{\ell}$,
we set
$V^{(\ell)}_{\gminib}:=
 \bigoplus_{\eta_{\ell}(\gminia)=\gminib}
 \bigoplus_{\alpha\in\cnum}V_{\gminia,\alpha}$.
Let $\pi^{(\ell)}_{\gminib}$ denote the projection of $V$
onto $V^{(\ell)}_{\gminib}$
with respect to the decomposition
$V=\bigoplus V^{(\ell)}_{\gminib}$.

We take total orders $\leq'$ on $\nbigi_{\ell}$
for any $\ell$
such that the induced maps
$\nbigi_{\ell_1}\lrarr\nbigi_{\ell_2}$
are order-preserving for each $\ell_1\leq\ell_2$.
Let $V^{\prime(\ell)}_{\gminib}$
be the orthogonal complement of
$\bigoplus_{\gminic<'\gminib}V_{\gminic}^{(\ell)}$
in 
$\bigoplus_{\gminic\leq'\gminib}V_{\gminic}^{(\ell)}$.
Let $\pi^{\prime(\ell)}_{\gminib}$ be the orthogonal projection
onto $V^{\prime(\ell)}_{\gminib}$.

We put $\zeta_{\ell}:=\eta_{\ell}-\eta_{\ell+1}$.
Let $f$ be the endomorphism of $V$
determined by $\theta=f\,dz$.
We put
$f^{(\ell)}:=
 f-\sum_{\gminia}
 \del_z\bigl(\eta_{\ell+1}(\gminia)\bigr)\pi_{\gminia}$,
$\mu^{(\ell)}:=
 f^{(\ell)}-\sum_{\gminia}
 \del_z\zeta_{\ell}(\gminia)\pi^{\prime}_{\gminia}$
and 
$\nbigr_{\gminib}^{(\ell)}:=
 \pi_{\gminib}^{(\ell)}-\pi_{\gminib}^{\prime(\ell)}$.
We consider the following claims.
\begin{description}
\item[($P_\ell$)]
 $|f^{(\ell')}|_h=O(|z|^{-\ell'-1})$
 for any $\ell'\geq\ell$.
\item[($Q_\ell$)]
 $|\mu^{(\ell')}|_h=O(|z|^{-\ell'})$
 for any $\ell'\geq\ell$.
 \item[($R_\ell$)]
 There exists $C>0$ such that
 $|\nbigr_{\gminib}^{(\ell')}|_h=O\bigl(
 \exp(-C|z|^{-\ell'})
 \bigr)$
for any $\ell'\geq \ell$
and for any $\gminib\in\nbigi_{\ell'}$.
\end{description}
The asymptotic orthogonality of
$V_{\gminia,\alpha}$ and $V_{\gminib,\beta}$ 
$(\gminia\neq\gminib)$
follows from $(R_{1})$.

We have the expression
$\theta^{\dagger}=f^{\dagger}d\zbar$.
We set $\Delta:=-\del_z\del_{\zbar}$.
If a holomorphic section $s$ of $\End(E)$
satisfies $[s,f]=0$,
we obtain the following inequality
for some $C_0>0$ and $\epsilon_0>0$,
which follows from (\ref{eq;17.10.21.2}) with $P=1$:
\begin{equation}
 \label{eq;13.1.11.2}
\Delta\log|s|^2_h
\leq
-\frac{\bigl|[f^{\dagger},s]\bigr|_h^2}{|s|_h^2}
+C_0\exp(-\epsilon_0|z|^{-p}).
\end{equation}
Let $f^{(\ell)\dagger}$ denote the adjoint of $f^{(\ell)}$
with respect to $h$.
Suppose that
the claims $P_{\ell+1}$, $Q_{\ell+1}$ and $R_{\ell+1}$ hold.
Because
$[f-f^{(\ell)},f^{(\ell)}]=0$,
we obtain
$\bigl[
 (f-f^{(\ell)})^{\dagger},
 f^{(\ell)}
 \bigr]=O\bigl(\exp(-\epsilon|z|^{-\ell-1})\bigr)$
by $R_{\ell+1}$.
By applying (\ref{eq;13.1.11.2}) to $f^{(\ell)}$,
we obtain the following for some $C_1>0$
 as in (99) of \cite{Mochizuki-wild}:
\[
 \Delta\log|f^{(\ell)}|^2_h
\leq
-\frac{\bigl|[f^{(\ell)\dagger},f^{(\ell)}]\bigr|_h^2}{|f^{(\ell)}|_h^2}
+C_{1}.
\]
Then, by the same argument
as that in \S7.3.2--\S7.3.3 of \cite{Mochizuki-wild},
we obtain $P_{\ell}$ and $Q_{\ell}$.
We put
\[
 k^{(\ell)}_{\gminib}:=
\log\bigl(|\pi^{(\ell)}_{\gminib}|_h^2\big/
 |\pi^{\prime(\ell)}_{\gminib}|_h^2\bigr)
=\log\bigl(1+|\nbigr^{(\ell)}_{\gminib}|_h^2\big/
 |\pi^{\prime(\ell)}_{\gminib}|_h^2\bigr).
\]
By applying (\ref{eq;13.1.11.2})
to $\pi^{(\ell)}_{\gminib}$,
we obtain
\[
 \Delta\log k^{(\ell)}_{\gminib}
\leq
 -\frac{\bigl|[f^{\dagger},\pi_{\gminib}^{(\ell)}]\bigr|_h^2}
 {|\pi_{\gminib}^{(\ell)}|_h^2}
+C_0\exp\bigl(
 -\epsilon_0|z|^{-p}
 \bigr).
\]
For any $A>0$,
there exists $C_2>0$ and $r_2>0$
such that the following holds for any $|z|<r_2$:
\begin{multline}
 \Delta \exp(-A|z|^{-\ell})
\geq
 -\exp(-A|z|^{-\ell})\,
\Bigl(
 \frac{\ell^2}{4}A^2|z|^{-2(\ell+1)}
 \Bigr)
 \\
\geq
 -\exp(-A|z|^{-\ell})\,
 \frac{\ell^2}{4}A^2C_2|z|^{-2(\ell+1)}
+C_0\exp(-\epsilon_0|z|^{-p}).
\end{multline}
Hence, 
we obtain $R_{\ell}$
by using the argument in \S7.3.4 of \cite{Mochizuki-wild}.
Similarly,
we obtain the asymptotic orthogonality
of $V_{\gminia,\alpha}$ and $V_{\gminia,\beta}$ $(\alpha\neq\beta)$,
and the boundedness of
$\theta_{\gminia,\alpha}-(d\gminia+\alpha dz/z)\id_{V_{\gminia,\alpha}}$
by using the argument in \S7.3.5--\S7.3.7 of \cite{Mochizuki-wild}
with (\ref{eq;13.1.11.2}).
\hfill\qed

\vspace{.1in}
By applying the argument in
\cite[\S7.2.5]{Mochizuki-wild},
we obtain the following corollary.
\begin{cor}
\label{cor;13.1.11.22}
$(V,\delbar_V,h)$ is acceptable,
i.e., the curvature
$F(h)$ is bounded with respect to $h$ and $g_{\Poin}$.
\index{acceptable bundle}
\hfill\qed
\end{cor}

\subsection{The associated filtered bundles
and formal $\lambda$-connections}
\label{subsection;17.12.16.4}

Let $\lambda\in\cnum$.
We obtain the holomorphic bundle 
$V^{\lambda}=(V,\delbar_V+\lambda\theta^{\dagger})$
on $U\setminus\{0\}$.
As a consequence of Proposition \ref{prop;12.8.7.3}
and (\ref{eq;17.10.21.2}) with $P=1$,
we obtain the following estimate with respect to $h$:
\[
 \bigl[\delbar_V+\lambda\theta^{\dagger},
  \del_V-\lambdabar\theta\bigr]
=[\delbar_V,\del_V]
-|\lambda|^2[\theta,\theta^{\dagger}]
=O\Bigl(
 |z|^{-2}(-\log|z|)^{-2}
 \Bigr)\,dz\,d\zbar.
\]
Namely,
$(V^{\lambda},\delbar_V+\lambda\theta^{\dagger},h)$
is acceptable.
Hence, we obtain a filtered bundle
$\nbigp^h_{\ast}V^{\lambda}$ on $(U,0)$
as in Proposition \ref{prop;20.7.26.1}.

We set $\DD^{\lambda}:=
 \delbar_V+\lambda\theta^{\dagger}
 +\lambda\del_V+\theta$.
Let $\vecv$ be a holomorphic frame of
$\nbigp^h_aV^{\lambda}$
which is compatible with the parabolic structure,
i.e.,
there exists a decomposition
$\vecv=\coprod_{a-1<b\leq a}\vecv_b$
such that
(i) $\vecv_b$ is a tuple of holomorphic sections of
$\nbigp_bV^{\lambda}$,
(ii) $\vecv_b$ induces a frame of
$\nbigp^h_bV^{\lambda}\bigl/\sum_{c<b}\nbigp^h_cV^{\lambda}$.
Let $A$ be the matrix valued function determined by
$\DDlambda\vecv=\vecv A\,dz$.

\begin{prop}
\label{prop;17.10.6.201}
There exists $N\in\seisuu_{>0}$ such that
$z^NA$ is $C^{\infty}$,
and that
the following holds:
\begin{itemize}
\item
For any $(\ell_1,\ell_2)\in\seisuu_{\geq 1}\times\seisuu_{\geq 0}$,
we obtain
$\del_{\zbar}^{\ell_1}\del_{z}^{\ell_2}
 (z^NA)
=O\bigl(\exp(-\epsilon(\ell_1,\ell_2)|z|^{-p})\bigr)$
for some $\epsilon(\ell_1,\ell_2)>0$.
\end{itemize}
As a result,
the Taylor series of
$z^NA$ at $0$
is of the form
$\sum_{j=0}^{\infty} B_jz^j\,dz$ $(B_j\in M_{\rank V}(\cnum))$.
\end{prop}
\pf
Let $C$ be the matrix valued function on $U\setminus\{0\}$
determined by
$(\del_V-\lambdabar\theta)\vecv
=\vecv C\,dz$.
Let $F_C$ be the endomorphism
determined by
$F_C\vecv=\vecv\,C$.
By the estimate for the connection form
of acceptable bundles
\cite[Lemma 21.9.3]{Mochizuki-wild},
we obtain $|F_C|_h=O\bigl(|z|^{-1}(-\log|z|)^{N_0}\bigr)$
for some $N_0$.
According to Proposition \ref{prop;12.8.7.3},
we have $|\theta|_h=O(|z|^{-N_1})$ for some $N_1$.
Because $\lambda\del_V+\theta-\lambda(\del_V-\lambdabar\theta)
=(1+|\lambda|^2)\theta$,
we obtain
$A=O\bigl(|z|^{-N}\bigr)$ for some $N\in\seisuu$.
Hence, $z^NA$ is bounded.

Set
$G:=\bigl[
 \delbar_V+\lambda\theta^{\dagger},
\lambda\del_V+\theta
 \bigr]
=\lambda\bigl(
 [\delbar_V,\del_V]+[\theta,\theta^{\dagger}]
 \bigr)$.
We have
$G\vecv=\vecv\delbar A$.
For each $\ell\in\seisuu_{\geq 0}$,
there exists $\epsilon_{\ell}>0$
such that 
$(\del_{\zbar}+\lambda\theta^{\dagger}_{\zbar})^{\ell}G
=O\bigl(
 \exp(-\epsilon_{\ell}|z|^{-p})
 \bigr)$,
which implies
$\del_{\zbar}^{\ell}(\delbar A)=
 O\bigl(
 \exp(-\epsilon_{\ell}|z|^{-p})
 \bigr)$.
We obtain the claim of the proposition
from the following lemma.

\begin{lem}
\label{lem;17.10.6.210}
Let $f$ be a $C^{\infty}$-function on
$U\setminus\{0\}$.
Suppose that
for any $\ell\in\seisuu_{\geq 0}$
there exists
$\epsilon_{\ell}>0$ such that
$|\del_{\zbar}^{\ell}f|=
 O\bigl(
 \exp(-\epsilon_{\ell}|z|^{-p})\bigr)$.
Then, for any 
$(\ell_1,\ell_2)\in\seisuu_{\geq 0}\times\seisuu_{\geq 0}$,
we obtain
$\bigl|
  \del_{z}^{\ell_1}\del_{\zbar}^{\ell_2}f
 \bigr|=
 O\bigl(\exp(-\epsilon_{\ell_1,\ell_2}|z|^{-p})\bigr)$
for some $\epsilon_{\ell_1,\ell_2}>0$.
\end{lem}
\pf
We may assume that $U=\{|z|<1\}$.
Let 
$\Psi:\{\zeta\in\cnum\,|\,\Re(\zeta)<0\}
\lrarr U\setminus\{0\}$
be the covering map
defined by $\Psi(\zeta)=\exp(\zeta)$.
Set $f_1:=\Psi^{\ast}(f)$.
By the condition,
for any $\ell\in\seisuu_{\geq 0}$,
there exists $\epsilon_{1,\ell}>0$ such that
$\bigl|
  \del_{\zetabar}^{\ell}f_1
 \bigr|=
 O\bigl(\exp(-\epsilon_{1,\ell}|\Psi(\zeta)|^{-p})\bigr)$.

Take $\zeta_0\in\cnum$ with
$\Re(\zeta_0)<-10$.
Let $B(\zeta_0,r)$ be the disc with radius $r$
with center $\zeta_0$.
There exist $C_{\ell}>0$ and $\epsilon_{2,\ell}>0$,
which are independent of $\zeta_0$,
such that the following holds
on $B(\zeta_0,2r)$:
\[
 |\del_{\zetabar}^{\ell+i}f|\leq 
 C_{\ell}\exp\bigl(-\epsilon_{2,\ell}|\Psi(\zeta_0)|^{-p}\bigr)
\quad
 (i=0,1,2).
\]
For any $q>1$, there exists $C_{q,\ell}>0$,
independently from $\zeta_0$,
such that
the $L^q_1$-norm of 
$\del_{\zetabar}^{\ell+1}f$ on $B(\zeta_0,7r/4)$
is dominated by
$C_{q,\ell}\exp\bigl(
 -\epsilon_{2,\ell}|\Phi(\zeta_0)|^{-p}
 \bigr)$.
Then,
there exists $C'_{q,\ell}>0$,
independently from $\zeta_0$,
such that the $L^q$-norm of 
$\del_{\zeta}\del_{\zetabar}^{\ell+1}f$ on 
$B(\zeta_0,5r/3)$
is dominated by 
$C'_{q,\ell}
 \exp\bigl(
 -\epsilon_{2,\ell}|\Phi(\zeta_0)|^{-p}
 \bigr)$.
Then, there exists $C''_{\ell}>0$,
independently from $\zeta_0$,
such that the sup norm of 
$\del_{\zeta}\del_{\zetabar}^{\ell}f$ on 
$B(\zeta_0,3r/2)$
is dominated by 
$C''_{\ell}
 \exp\bigl(
 -\epsilon_{2,\ell}|\Phi(\zeta_0)|^{-p}
 \bigr)$.
By an easy induction,
we obtain the claim of Lemma \ref{lem;17.10.6.210}.
The proof of Proposition \ref{prop;17.10.6.201}
is also completed.
\hfill\qed

\subsection{Formal good filtered $\lambda$-flat bundles}
\label{subsection;20.7.27.3}
\index{formal good filtered $\lambda$-flat bundle}

Let $\nbigvhat^{\lambda}$
denote
the free $\cnum(\!(z)\!)$-module
obtained as the formal completion of
$\nbigp^h V^{\lambda}_0$,
i.e.,
$\nbigvhat^{\lambda}:=
 \cnum[\![z]\!]\otimes_{\nbigo_{U,0}}
 \nbigp^h\nbigv^{\lambda}_0$.
As a result of Proposition \ref{prop;17.10.6.201},
$\DDlambda$ induces
a $\lambda$-connection
$\widehat{\DD}^{\lambda}$
of $\nbigvhat^{\lambda}$.
As the formal completion of
$\nbigp^h_{\ast}\nbigv^{\lambda}$,
we obtain a filtered bundle
$\nbigp_{\ast}\nbigvhat^{\lambda}$
over $\nbigvhat^{\lambda}$.

\begin{prop}
\label{prop;17.12.16.3}
$(\nbigp_{\ast}\nbigvhat^{\lambda},\DDhat^{\lambda})$
is a unramifiedly good filtered $\lambda$-flat bundle.
\index{good filtered $\lambda$-flat bundle
$(\nbigp_{\ast}\nbigvhat^{\lambda},\DDhat^{\lambda})$}
More precisely,
there exists the decomposition
\begin{equation}
\label{eq;20.7.26.21}
 (\nbigp_{\ast}\nbigvhat^{\lambda},\DDhat^{\lambda})
=\bigoplus_{\gminia\in\nbigi}
 \bigl(
 \nbigp_{\ast}\nbigvhat_{\gminia},\DDhat_{\gminia}^{\lambda}
 \bigr),
\end{equation}
where $\DDhat_{\gminia}^{\lambda}-(1+|\lambda|^2)d\gminia$
are logarithmic,
and $\nbigi$ is the index set of the decomposition
{\rm(\ref{eq;17.10.20.30})}.
\end{prop}
\pf
Let $\pi_{\gminia}$ be the projection onto $V_{\gminia}$
with respect to the decomposition {\rm(\ref{eq;17.10.20.30})}.
According to Proposition \ref{prop;12.8.7.3},
we obtain the following estimate:
\[
 \bigl[
 \delbar_V+\lambda\theta^{\dagger},\pi_{\gminia}
 \bigr]
=O\bigl(
 \exp(-\epsilon|z|^{-p})
 \bigr).
\]
Hence, for any large $N>0$,
by using the argument in
\cite[Lemma 7.4.7]{Mochizuki-wild},
we can construct holomorphic endomorphisms
$p_{\gminia}$ $(\gminia\in\nbigi)$
of $\nbigp^h_{\ast}V^{\lambda}$
satisfying the following conditions:
\begin{equation}
p_{\gminia}-\pi_{\gminia}=O(|z|^N),
\quad
[p_{\gminia},p_{\gminib}]=0,
\quad
p_{\gminia}\circ p_{\gminia}=p_{\gminia},
\quad
\sum p_{\gminia}=\id. 
\end{equation}

\begin{lem}
\label{lem;20.7.26.10}
 $\widehat{\DD}^{\lambda}
 -\sum_{\gminia\in\nbigi}
 (1+|\lambda|^2)d\gminia\cdot p_{\gminia}$
is logarithmic with respect to
$\nbigp_{\ast}\nbigvhat^{\lambda}$.
\end{lem}
\pf
We consider
\[
\DDlambda_0:=
\DDlambda-\sum (1+|\lambda|^2)d\gminia\cdot p_{\gminia}
=\delbar_V+\lambda\theta^{\dagger}
+\lambda(\del_V-\lambdabar\theta)
+(1+|\lambda|^2)
 \Bigl(
 \theta-\sum d\gminia\cdot p_{\gminia}
 \Bigr).
\]
Note that
$\bigl|
\theta-\sum d\gminia\cdot p_{\gminia}\bigr|_h
 =O\bigl(|z|^{-1}dz\bigr)$.
Hence, we obtain that
$(\DD_0^{\lambda})_{|\widehat{0}}$ is logarithmic.
\hfill\qed

\vspace{.1in}

The claim of Proposition \ref{prop;17.12.16.3}
follows from Lemma \ref{lem;20.7.26.10} and 
\cite[Proposition 2.3.6]{Mochizuki-wild}.
\hfill\qed

\begin{lem}
For each $\gminia$,
there exists a good filtered $\lambda$-flat bundle
$\bigl(
\nbigp_{\ast}\nbigv_{\gminia},\DDlambda_{\gminia}
\bigr)$
on $(U,0)$ with an isomorphism
$(\nbigp_{\ast}\nbigv_{\gminia},\DDlambda_{\gminia})
\otimes_{\nbigo_{U}}\cnum[\![z]\!]
 \simeq
 (\nbigp_{\ast}\nbigvhat_{\gminia},\DDlambda_{\gminia})$. 
\end{lem}
\pf
The claim is obvious in the case $\lambda=0$.
Let us study the case $\lambda\neq 0$.
It is enough to study the case $\gminia=0$.
Therefore,
we may assume that $\nbigi=\{0\}$ from the beginning,
i.e.,
$\DDhat^{\lambda}$ is logarithmic with respect to
$\nbigp_{\ast}\nbigvhat$.

Let $\Sp_{-1<b\leq 0}(\Res(\DDhat^{\lambda}))$
denote the eigenvalues of
$\Res(\DDhat^{\lambda})$
on $\bigoplus_{-1<b\leq 0}\Gr^{\nbigp}_b(\nbigvhat)$.
Let $N$ be a large integer such that
\[
 N>100
   \max
   \Bigl\{
   \bigl|\lambda^{-1}(\alpha-\beta)\bigr|\,\Big|\,
   \alpha,\beta\in \Sp_{-1<b\leq 0}(\DDhat^{\lambda})
    \Bigr\}.
\]
We set $r:=\rank\nbigvhat$.
Let $\vecv$ be a frame of
$\nbigp_0\nbigvhat$ compatible with a parabolic structure.
We obtain $A\in M_r(\cnum[\![z]\!])$
determined by
$\DDhat^{\lambda}\vecv=\vecv\cdot A\,dz/z$.
For the expansion $A=\sum_{j=0}^{\infty} A_jz^j$,
the set of the eigenvalues of $A_0$
is equal to
$\Sp_{-1<b\leq 0}(\Res(\DDhat^{\lambda}))$.

We inductively construct
$G_i\in M_r(\cnum)$ $(i=N,N+1,\ldots,M)$
such that the following condition is satisfied.
\begin{itemize}
 \item Let $I_r\in M_r(\cnum)$ denote the identity matrix.
       We set
       $\vecv^{(M)}:=
       \vecv\cdot (I_r+G_Nz^N)(I_r+G_{N+1}z^{N+1})\cdots (I_r+G_mz^M)$.
       Let $A^{(M)}\in M_r(\cnum[\![z]\!])$ be determined by
       $\DDhat^{\lambda}\vecv^{(M)}=\vecv^{(M)}\cdot A^{(M)}dz/z$.
       Then, for the expansion
       $A^{(M)}=\sum_{j=0}^{\infty} A^{(M)}_jz^j$,
       we obtain $A^{(M)}_i=A_i$ $(i<N)$
       and $A^{(M)}_i=0$ $(N\leq i\leq M)$.
\end{itemize}
Suppose that we have already constructed
$G_i$ $(i=N,\ldots,M-1)$.
There uniquely exists $G_M\in M_r(\cnum)$
such that
$MG_M+[A^{(M-1)}_0,G_M]+A^{(M-1)}_M=0$.
Then, it is easy to check that
the condition is satisfied for $M$.

As the limit of $\vecv^{(M)}$
with respect to the $z$-adic topology,
we obtain a frame $\vecv^{(\infty)}$
of $\nbigp_0\nbigvhat$.
Let $A^{(\infty)}\in M_r(\cnum[\![z]\!])$
be the matrix determined by
$\DDhat^{\lambda}\vecv^{(\infty)}
=\vecv^{(\infty)}A^{(\infty)}dz/z$.
By the construction,
$A^{(\infty)}\in M_r(\cnum[z])$.
Thus, we obtain the claim of the lemma.
\hfill\qed

\subsection{Residues and KMS-structure}

For any $a\in\real$,
we set
$\nbigp_{<a}\nbigvhat^{\lambda}:=
\sum_{b<a}\nbigp_b\nbigvhat^{\lambda}$
and
$\Gr^{\nbigp}_a(\nbigvhat^{\lambda}):=
\nbigp_a\nbigvhat^{\lambda}\big/
\nbigp_{<a}\nbigvhat^{\lambda}$.
\index{vector space $\Gr^{\nbigp}_a(\nbigvhat^{\lambda})$}
We also obtain
$\nbigp_{<a}\nbigvhat^{\lambda}_{\gminia}$
and 
$\Gr^{\nbigp}_a(\nbigvhat_{\gminia}^{\lambda})$
from
$(\nbigp_{\ast}\nbigvhat_{\gminia}^{\lambda},
\DDhat^{\lambda}_{\gminia})$.
There exists the natural decomposition
$\Gr^{\nbigp}_a(\nbigvhat^{\lambda})
=\bigoplus_{\gminia\in\nbigi}
 \Gr^{\nbigp}_a(\nbigvhat_{\gminia}^{\lambda})$.

We obtain the residue endomorphisms
$\Res(\DDhatlambda_{\gminia})$
of $\Gr^{\nbigp}_a(\nbigvhat_{\gminia}^{\lambda})$
for any $a\in\real$.
(See \cite[\S2.5.2]{Mochizuki-wild}.)
There exists the generalized eigen decomposition
\[
  \Gr^{\nbigp}_a(\nbigvhat_{\gminia}^{\lambda})
=\bigoplus_{\alpha\in\cnum}
 \EE_{\alpha}
 \Gr^{\nbigp}_a(\nbigvhat_{\gminia}^{\lambda}),
\]
where 
$\EE_{\alpha}
 \Gr^{\nbigp}_a(\nbigvhat_{\gminia}^{\lambda})$
denote the generalized eigen spaces
according to $\alpha$.
Let $\KMS(\nbigp_{\ast}\nbigvhat_{\gminia}^{\lambda})$
denote the set of 
$(a,\alpha)$ such that
$\EE_{\alpha}
 \Gr^{\nbigp}_a(\nbigvhat_{\gminia}^{\lambda})\neq 0$.
\index{set $\KMS(\nbigp_{\ast}\nbigvhat_{\gminia}^{\lambda})$}
We put
$\gminim^{\lambda}(a,\alpha,\gminia):=
 \dim
 \EE_{\alpha}
 \Gr^{\nbigp}_a(\nbigvhat_{\gminia}^{\lambda})$.
Let $N^{\lambda}_{a,\alpha,\gminia}$
denote the nilpotent part of
$\Res(\DDhat_{\gminia}^{\lambda})$
on 
$\EE_{\alpha}
 \Gr^{\nbigp}_a(\nbigvhat_{\gminia}^{\lambda})\neq 0$.
We obtain the monodromy weight filtration $W$
of $N^{\lambda}_{a,\alpha,\gminia}$.
\index{number $\gminim^{\lambda}(a,\alpha,\gminia)$}
\index{endomorphism $N^{\lambda}_{a,\alpha,\gminia}$}

\subsection{Norm estimate}
\label{subsection;20.7.27.4}
\index{norm estimate}

Let $\vecv$ be a frame of
$\nbigp_{a}\nbigvlambda$
which is compatible with the parabolic filtration
and the weight filtration,
i.e.,
there exists a decomposition
$\vecv=\coprod_{a-1<b\leq a}\coprod_{k\in\seisuu}\vecv_{b,k}$
such that
(i) $\coprod_{k\in\seisuu}\vecv_{b,k}$ is a tuple of sections of
$\nbigp_a\nbigvlambda$,
(ii) $\coprod_{k\leq\ell}\vecv_{b,k}$
induces a frame of
$W_{\ell}\Gr^{\nbigp}(\nbigvhat^{\lambda})$.
For $v_i\in \vecv_{b,k}$,
we set $b(v_i):=b$ and $k(v_i):=k$.

Let $h_0$ be a Hermitian metric of $V$
for which 
\[
 h_0(v_i,v_j):=
 \left\{
 \begin{array}{ll}
 |z|^{-2b(v_i)}(-\log|z|^2)^{k(v_i)}
 & (i=j),\\
 0 & (\mbox{otherwise}).
 \end{array}
 \right.
\]
We can prove the following proposition
by using the same argument as the proof of
\cite[Proposition 8.1.1]{Mochizuki-wild},
which originally goes back to \cite{Simpson90}.
\begin{prop}
\label{prop;20.7.27.2}
 $h_0$ and $h$ are mutually bounded.
\end{prop}
\pf
We explain only an outline.
We decompose
$\DD^{\lambda}=d''+d'=(\delbar_V+\lambda\theta)
+(\lambda\del_V+\theta^{\dagger})$.
Note that
$\bigl|
[d'',d']
\bigr|_h
=O\bigl(\exp(-\epsilon|z|^{-p})\bigr)$.
For any large $N$,
there exists a $C^{\infty}$-section $r_N$ of
$\End(V^{\lambda})$
such that
(i) $|r_N|_h+|d'r_N|_h+|d''r_N|_h=O\bigl(|z|^{2N}\bigr)$,
(ii) $[d'',d'-r_N\,dz]=0$.
We set
$\DD_1^{\lambda}:=
\DD^{\lambda}-r_N\,dz$.
Let $\DD_1^{\lambda}=d''_1+d'_1$
be the decomposition into
the $(0,1)$-part and the $(1,0)$-part.
We obtain the differential operators
$\delta'_{1,h}$ and $\delta''_{1,h}$
as in \S\ref{subsection;17.10.12.1}.
We set
$\DD_{1,h}^{\lambda\star}=\delta'_{1,h}-\delta''_{1,h}$.
Then, we obtain
$\bigl|
[\DD^{\lambda}_1,\DD_{1,h}^{\lambda\star}]
\bigr|_h
=O\bigl(|z|^N\bigr)$.
If $N$ is sufficiently large,
$(\nbigp^h_{\ast}V^{\lambda},\DDlambda_1)$
is a good filtered $\lambda$-flat bundle
by \cite[Proposition 2.3.6]{Mochizuki-wild}.
The formal completion at $0$ has the decomposition
\begin{equation}
\label{eq;20.7.26.20}
 (\nbigp_{\ast}(\nbigvhat^{\lambda}),\DDhat^{\lambda}_1)
 =\bigoplus_{\gminia\in\nbigi}
 (\nbigp_{\ast}\nbigvhat^{\lambda}_1,\DDhat^{\lambda}_{1,\gminia})
\end{equation}
such that
$\DDhat^{\lambda}_{1,\gminia}
-d\gminia\id$ are logarithmic.
Note that the decompositions
(\ref{eq;20.7.26.21})
and (\ref{eq;20.7.26.20}) are not necessarily the same.
But, because
$(\DDhat^{\lambda}-\DDhat^{\lambda}_1)
\nbigp_b\nbigvhat^{\lambda}
\subset
\nbigp_{b-N}\nbigvhat^{\lambda}$
for any $b\in\real$,
we obtain
$\Res(\DDhat^{\lambda})$
and
$\Res(\DDhat_1^{\lambda})$
and
$\Gr^{\nbigp}_a(\nbigvhat^{\lambda}_{\gminia})
=\Gr^{\nbigp}_a(\nbigvhat^{\lambda}_{1,\gminia})$.

As in the proof of
\cite[Proposition 2.3.6]{Mochizuki-wild},
we can construct a Hermitian metric $h_1$ of $V^{\lambda}$
such that
(i) $h_1$ is mutually bounded with $h_0$,
(ii)
$\bigl[\DD^{\lambda}_1,\DD^{\lambda\star}_{1,h_1}\bigr]
=O(|z|^{-2+\epsilon})dz\,d\zbar$ for some $\epsilon>0$.
Then, by \cite[Corollary 4.3]{Simpson90}
(see also \cite[\S4.3]{Mochizuki-JDG}
for the case of general $\lambda$),
we obtain that $h$ and $h_1$ are mutually bounded.
\hfill\qed

\subsection{Comparison of KMS-structures}

We define the map
$\kmsmap(\lambda):
 \real\times\cnum
\lrarr
 \real\times\cnum$
by 
\index{map $\kmsmap(\lambda)$}
\[
 \kmsmap(\lambda,a,\alpha):=
 \bigl(
 a+2\Re(\lambda\alphabar),
 \alpha-a\lambda-\alphabar\lambda^2
 \bigr).
\]
We can prove the following proposition
by using the argument in \cite[\S7]{Simpson90}
(see also \cite[Proposition 8.2.1]{Mochizuki-wild}).
\begin{prop}
\label{prop;20.7.27.200}
 The map $\kmsmap(\lambda)$ induces
a bijection
$\KMS(\nbigp_{\ast}\nbigvhat^{0}_{\gminia})
\simeq
\KMS(\nbigp_{\ast}\nbigvhat^{\lambda}_{\gminia})$.
 It preserves the multiplicities,
 i.e.,
 $\gminim^0(a,\alpha,\gminia)
 =\gminim^{\lambda}\bigl(\kmsmap(\lambda,a,\alpha),\gminia\bigr)$.
The conjugacy classes of
$N^0_{a,\alpha,\gminia}$
and 
$N^{\lambda}_{\kmsmap(\lambda,a,\alpha),\gminia}$
are the same.
\end{prop}
\pf
We explain only an outline.
Let $(\nbigp^h_{\ast}V^{\lambda},\DDlambda_1)$
denote the good filtered $\lambda$-flat bundle
as in the proof of Proposition \ref{prop;20.7.27.2}.

Let $\pi_{\gminia,\alpha}$ denote the projection
onto $V_{\gminia,\alpha}$
with respect to the decomposition (\ref{eq;20.7.27.1}).
Let $p_{\gminia}$ be as in the proof of
Proposition \ref{prop;17.12.16.3}.
As in \cite[Lemma 8.2.4]{Mochizuki-wild},
we can construct endomorphisms $p_{\gminia,\alpha}$  of
$\nbigp^h_{\ast}V^{\lambda}$ $(\gminia\in\nbigi)$
such that
\[
\bigl|
 p_{\gminia,\alpha}-\pi_{\gminia,\alpha}
 \bigr|_h=O\bigl(|z|^{\epsilon}\bigr),
 \quad
 p_{\gminia,\alpha}^2=p_{\gminia,\alpha},
\]
\[
 p_{\gminia,\alpha}\circ p_{\gminib,\beta}=0
 \,\,\,((\gminia,\alpha)\neq (\gminib,\beta)),
 \quad
 \sum_{\alpha} p_{\gminia,\alpha}=p_{\gminia}.
\]
Here, $\epsilon$ denotes a positive number.
We set
$\Phi:=\sum_{\gminia,\alpha}
(d\gminia+\alpha dz/z)\cdot p_{\gminia,\alpha}$.
The following lemma is an analogue of
\cite[Lemma 8.2.5]{Mochizuki-wild}.
\begin{lem}
\label{lem;20.7.27.5}
 Let $v$ be a holomorphic section of
$V^{\lambda}$ such that
$|v|_h\sim |z|^{-b}(-\log |z|)^k$
for $b\in\real$ and $k\in\seisuu$.
Then, we obtain
\[
 \bigl|
 \DDlambda_1-(1+|\lambda|^2)\Phi
 +\lambda \cdot a \cdot v\cdot dz/z
 \bigr|_h
 =O\bigl(|z|^{-b}(-\log|z|)^{k-1}\bigr)\,|dz/z|.
\]
\end{lem}
\pf
It is enough to study the case $b=0$.
We obtain
\[
\DDlambda_1 v-(1+|\lambda|^2)\Phi(v)
=
\lambda \delta'_{1,h}v
+(1+|\lambda|^2)(\theta-\Phi)v
-r_Nv.
\]
Because
$|(\theta-\Phi)v|=O\bigl((-\log|z|)^{k-1}\bigr)$
(see Proposition \ref{prop;12.8.7.3})
and
$|r_Nv|=O(|z|^{N-1})$,
it is enough to prove that
\[
 \int \bigl|
  \delta'_{1,h}v
  \bigr|^2_h
  (-\log |z|)^{1-\delta-2k}
  \bigl|dz\,d\zbar\bigr|
  <\infty
\]
for any $\delta>0$.
We can prove it by using the argument
in \cite[Page 761--762]{Simpson90}.
\hfill\qed

\vspace{.1in}

We obtain the decomposition
$\nbigp^h_{\ast}V^{\lambda}=\bigoplus
\Image_{\gminia,\alpha}
 \nbigp_{\ast}\bigl(
p_{\gminia,\alpha}\bigr)$.
Let
\[
 \vecv_{\gminia,\alpha}=
 \coprod_{-1<a\leq 0}
 \coprod_{k\in\seisuu}
 \vecv_{\gminia,\alpha,a,k}
\]
be
a frame of $\nbigp_0^h(\Image p_{\gminia,\alpha})$
compatible with the parabolic structure
and the weight filtration.
(See \S\ref{subsection;20.7.27.4}.)
We obtain the induced frame
$\vecv=(v_1,\ldots,v_r)
=\coprod_{\gminia,\alpha}\vecv_{\gminia,\alpha}$
of
$\nbigp_0^hV^{\lambda}$.
For each $v_i\in \vecv_{\gminia,\alpha,a,k}$,
we set $\gminia(v_i)=\gminia$,
$\alpha(v_i)=\alpha$,
$a(v_i)=a$,
and $k(v_i)=k$.

We obtain $A\in M_r(\cnum)$
by $\Res(\DDlambda_1)\vecv_{|0}=\vecv_{|0}A$.
We obtain the following lemma by
Proposition \ref{prop;20.7.27.2},
Lemma \ref{lem;20.7.27.5}
and the construction of $\vecv$.
\begin{lem}
\label{lem;20.7.27.121}
 $A_{i,j}=0$ holds unless
 $\gminia(v_i)=\gminia(v_j)$,
 $\alpha(v_i)=\alpha(v_j)$,
 and
 $(a(v_i),k(v_i))\leq ''(a(v_j),k(v_j))$,
where $\leq''$ denotes the lexicographic order
on $\real\times\seisuu$.
Moreover,
$A_{i,i}=
(1+|\lambda|^2)\alpha(v_i)-\lambda a(v_i)$ holds. 
\hfill\qed
\end{lem}

As in \cite[\S8.2.5]{Mochizuki-wild},
there exists
a wild harmonic bundle
\begin{equation}
\label{eq;20.7.27.110}
(\Vtilde,\delbar_{\Vtilde},\thetatilde,\htilde)
=\bigoplus
(\Vtilde_{\gminia,\alpha},\delbar_{\Vtilde_{\gminia,\alpha}},
\thetatilde_{\gminia,\alpha},\htilde_{\gminia,\alpha})
\end{equation}
on $U\setminus\{0\}$
with an isomorphism of 
filtered bundles
$\Psi:\nbigp^h_{\ast}\Vtilde^0
\simeq\nbigp^h_{\ast}V^0$
such that
$(\theta\Psi-\Psi\thetatilde)\nbigp^h_a\Vtilde^0
 \subset
  \nbigp^h_{<a}V^0$
  for any $a\in\real$.
We may assume that
the decompositions
(\ref{eq;20.7.27.1})
and (\ref{eq;20.7.27.110})
are the same
under the isomorphism $\Psi$.
By the norm estimate (Proposition \ref{prop;20.7.27.2}),
$\Psi$ and $\Psi^{-1}$ are bounded.
We identify the bundles
$V$ and $\Vtilde$
via $\Psi$
as $C^{\infty}$-bundles.
The metrics $h$ and $\htilde$
are mutually bounded.
By Proposition \ref{prop;12.8.7.3} and the construction,
we obtain
$|\theta-\thetatilde|_h
 =O\bigl(|z|^{\epsilon}\bigr)\,dz/|z|$
for some $\epsilon>0$.
We also obtain
$|\theta_h^{\dagger}-\thetatilde^{\dagger}_{\htilde}|_{h}
 =O\bigl((-\log|z|)^{-1}\bigr)\,d\zbar/|z|$
by Proposition \ref{prop;12.8.7.3}.

\begin{lem}
\label{lem;20.7.27.120}
 Let $M$ be the $C^{\infty}$-endomorphism of $V$
determined by
 $\theta^{\dagger}_h-\thetatilde^{\dagger}_{\htilde}
 =M\,d\zbar/\zbar$.
Let $U'$ be any relatively compact neighbourhood of
$0$ in $U\cap \{|z|<1\}$.
Then, we obtain
\begin{equation}
\label{eq;20.7.27.101}
 \int_{U'} |M|_{h}^2\cdot
 (-\log|z|)\cdot
 |dz\,d\zbar|<\infty.
\end{equation}
\end{lem}
\pf
To clarify the argument,
we regard $M$ as a homomorphism
$\Vtilde\lrarr V$
obtained as
$\theta\circ\Psi-\Psi\circ\thetatilde=M\,dz/z$,
which is a section of
the bundle $\Hom(\Vtilde,V)$.
Note that
$\Hom(\Vtilde,V)$
with the induced Higgs field $\thetatilde^{(1)}$
and the induced Hermitian metric $\htilde^{(1)}$
satisfies (\ref{eq;17.10.21.2}).
By the proof of Weitzenb\"{o}ck formula
in \cite[Lemma 4.1, P743]{Simpson90},
there exist $C_i>0$ $(i=1,2)$ and $\epsilon_1>0$
such that
\[
-\del_{z}\del_{\zbar}\bigl|\Psi\bigr|^2_{\htilde^{(1)}}
 \leq
 -C_1\bigl|M\bigr|^2_{\htilde^{(1)}}
+C_2|z|^{2\epsilon_1-2}.
\]
We obtain
$-\del_z\del_{\zbar}\bigl(
 |\Psi|^2_{\htilde^{(1)}}
 +C_2\epsilon_1^{-2}|z|^{\epsilon_1}
 \bigr)
 \leq -C_1|M|^2_{\htilde^{(1)}}$.
 Because
$|\Psi|^2_{\htilde^{(1)}}
 +C_2\epsilon_1^{-2}|z|^{\epsilon_1}$
is bounded,
we obtain (\ref{eq;20.7.27.101})
by the argument in \cite[Lemma 7.7]{Simpson90}
(see also  Lemma \ref{lem;20.7.27.100} below).
\hfill\qed

\vspace{.1in}

From 
$(\Vtilde,\delbar_{\Vtilde},\thetatilde,\htilde)$,
we obtain
a good filtered $\lambda$-flat bundle
$(\nbigp^h_{\ast}\Vtilde^{\lambda},
\DDtilde^{\lambda})
=\bigoplus_{\gminia\in\nbigi}
(\nbigp^h_{\ast}\Vtilde^{\lambda}_{\gminia,\alpha},
\DDtilde^{\lambda}_{\gminia,\alpha})$.
Let $\ptilde_{\gminia,\alpha}$
denote the projection of
$\nbigp^h_{\ast}\Vtilde^{\lambda}$
onto
$\nbigp^h_{\ast}\Vtilde^{\lambda}_{\gminia,\alpha}$
with respect to the decomposition.
Note that
$\ptilde_{\gminia,\alpha}=\pi_{\gminia,\alpha}$
under the isomorphism $\Psi$.
Hence, we obtain
$\ptilde_{\gminia,\alpha}
-p_{\gminia,\alpha}=O(|z|^{\epsilon})$.

Let $\vecvtilde_{\gminia,\alpha}
=\coprod_{-1<a\leq 0}
 \coprod_{k\in\seisuu}\vecvtilde_{\gminia,\alpha,a,k}$
be a frame of
$\nbigp_{0}\Vtilde^{\lambda}_{\gminia,\alpha}$
compatible with the parabolic structure
and the weight filtration.
We obtain a frame
$\vecvtilde=\coprod_{\gminia,\alpha}
\vecvtilde_{\gminia,\alpha}$
of $\nbigp_{0}\Vtilde^{\lambda}$.
For $\vtilde_i\in\vecvtilde_{\gminia,\alpha,a,k}$,
we set
$\gminia(\vtilde_i)=\gminia$,
$\alpha(\vtilde_i)=\alpha$,
$a(\vtilde_i)=a$,
and $k(\vtilde_i)=k$.
Let $I=(I_{j,i})$ be determined by
$\vtilde_i=\sum I_{j,i}v_j$.
We put
$\nbigb_{j,i}
=I_{j,i}|z|^{a(\vtilde_i)-a(v_j)-(k(\vtilde_i)-k(v_j))/2}$.
By using the argument in
\cite{Simpson90}
(see \cite[Proposition 8.2.2]{Mochizuki-wild}
for a summary)
with Lemma \ref{lem;20.7.27.120}
and a complement for the wild case
in \cite[Lemma 8.2.8]{Mochizuki-wild},
we can prove the following.
\begin{itemize}
 \item $\nbigb_{j,i}=O(|z|^{\epsilon})$
       unless
       $(\gminia(\vtilde_i),\alpha(\vtilde_i))=(\gminia(v_j),\alpha(v_j))$.
 \item $|\nbigb_{j,i}|=O\bigl((-\log|z|)^{-1}\bigr)$
       unless
       $a(\vtilde_i)=a(v_j)$.
 \item We obtain
       \[
       \int |\nbigb_{j,i}|_h^2
       \frac{|dz\,d\zbar|}{|z|^2(-\log|z|)^2\log(-\log|z|)}
       <\infty
        \]
        unless $k(\vtilde_i)=k(v_j)$.
\end{itemize}
Then, we obtain
$\dim\Gr^{W}_k\EE_{\beta}\Gr^{\nbigp}_a(\nbigvhat_{\gminia}^{\lambda})
=\dim\Gr^W_k\EE_{\beta}\Gr^{\nbigp}_a
(\nbigp\nbigvtilde^{\lambda}_{\gminia})$
for any $\gminia\in\nbigi$,
$-1<a\leq 0$,
$\beta\in\cnum$
and $k\in\seisuu$,
by applying the argument
in the proof of
\cite[Proposition 7.6, Theorem 7]{Simpson90}
with Lemma \ref{lem;20.7.27.121}.
Because the claim of Proposition \ref{prop;20.7.27.200}
holds for
$(\Vtilde,\delbar_{\Vtilde},\thetatilde,\htilde)$,
we obtain the claim of Proposition \ref{prop;20.7.27.200}
for $(V,\delbar_V,\theta,h)$.
\hfill\qed

\subsection{Appendix}

Let $\phi$ be a bounded $\real_{\geq 0}$-valued 
$C^{\infty}$-function on $\cnum^{\ast}$
such that
there exists $R_0>0$ such that
$\phi(w)=0$ for $|w|>R_0$.
We may naturally regard $\phi(w)$
as an $L^1$-function on $\cnum$,
and we obtain the distribution
$\del_{w}\del_{\wbar}\phi$
on $\cnum$.

Let $\rho:\cnum\lrarr\real_{\geq 0}$
be a $C^{\infty}$-function
such that
(i) $0\leq \rho(w)\leq 1$ for any $w$,
(ii) $\rho(w)=1$ for $|w|\leq 1$,
(iii) $\rho(w)=0$ for $|w|\geq 2$.

For any $z\in\cnum^{\ast}$,
we take
$0<\epsilon<|z|/100$,
and set
$\chi_{z,\epsilon}(w):=
\rho\bigl((z-w)/\epsilon\bigr)$.
The following integral is well defined.
\begin{multline}
\label{eq;20.7.27.130}
 F_{\phi}(z):=
 \frac{2}{\pi}
 \int_{\cnum} \del_w\del_{\wbar}\phi(w)
 \cdot\Bigl(
 \chi_{z,\epsilon}(w)
 \log|z-w|
 \Bigr)
 \\
 +\frac{2}{\pi}
 \int_{\cnum}
 \phi(w)
 \cdot\del_{w}\del_{\wbar}\Bigl(
 \bigl(1-\chi_{z,\epsilon}(w)\bigr)
 \cdot\log|z-w|
  \Bigr).
\end{multline}
It is easy to check that $F(z)$
is independent of the choice of $\epsilon$.
\begin{lem}
 $F_{\phi}(z)=\phi(z)$.
\end{lem}
\pf
For $z$ and $\epsilon$ as above,
we obtain
\begin{multline}
 \frac{2}{\pi}\int_{\cnum}
 \del_w\del_{\wbar}
\bigl(
 (1-\chi_{z,\epsilon}(w))\phi(w)
 \bigr)
 \cdot
 \Bigl(
 \chi_{z,\epsilon}(w)\log|z-w|
 \Bigr)
 \\
 +\frac{2}{\pi}\int_{\cnum}
 (1-\chi_{z,\epsilon}(w))\phi(w)
 \cdot\del_{w}\del_{\wbar}\Bigl(
 \bigl(1-\chi_{z,\epsilon}(w)\bigr)
 \cdot\log|z-w|
 \Bigr)
 = \\
  \frac{2}{\pi}
 \int_{|z-w|\geq \epsilon/4}
 \bigl(1-\chi_{z,\epsilon}(w)\bigr)
 \cdot\phi(w)\cdot
\del_w\del_{\wbar}\log|z-w|
 =0.
\end{multline}
Hence, we obtain
$F_{\phi}(z)
=F_{\chi_{z,\epsilon}\phi}(z)$.
Note that
$\chi_{z,\epsilon}\phi$ is naturally a test function
on $\cnum$.
Because
$\frac{1}{2\pi}\log|w|$
is the Green function for
the operator $4\del_{w}\del_{\wbar}$,
we obtain
$F_{\chi_{z,\epsilon}\phi}(z)
=(\chi_{z,\epsilon}\phi)(z)=\phi(z)$.
\hfill\qed

\vspace{.1in}

Let $U$ be a neighbourhood of $0$ in $\cnum$.
We set $U^{\ast}:=U\setminus\{0\}$.
Let $f_i:U^{\ast}\lrarr\real_{\geq 0}$ $(i=1,2)$
be $C^{\infty}$-functions
such that
(i) $f_1$ is bounded,
(ii) $f_2\leq \del_w\del_{\wbar}f_1$ on $U^{\ast}$.
The following lemma is contained in
\cite[Lemma 7.7]{Simpson90}.
\begin{lem}
\label{lem;20.7.27.100}
 For any relatively compact neighbourhood
 $U'$ of $0$ in $U$,
 we obtain
 $\int_{U'}f_2\cdot (-\log|w|)<\infty$.
\end{lem}
\pf
We may assume
$\bigl\{|w|\leq 10\bigr\}\subset U$.
It is enough to consider the case $U'=\{|w|<1/10\}$.
We may naturally regard
$\ftilde_1:=\rho(w/5)\cdot f_1$
as a $C^{\infty}$-function on $\cnum^{\ast}$.
The integral $F_{\ftilde_1}$ in (\ref{eq;20.7.27.130})
is well defined for $z\in \cnum^{\ast}$.
By Lemma \ref{lem;20.7.27.100},
we obtain
$F_{\ftilde_1}(z)=\ftilde_1(z)$.

Suppose $|z|<1/4$.
We obtain
\begin{multline}
 F_{\ftilde_1}(z)=
 \frac{2}{\pi}
 \int_{\cnum} \del_w\del_{\wbar}\ftilde_1(w)
 \cdot\Bigl(
 \chi_{z,\epsilon}(w)
 \log|z-w|
 \Bigr)
 \\
 +\frac{2}{\pi}
 \int_{\cnum}
 \ftilde_1(w)
 \cdot\del_{w}\del_{\wbar}\Bigl(
 \rho(3(z-w))\cdot
 \bigl(1-\chi_{z,\epsilon}(w)\bigr)
 \cdot\log|z-w|
 \Bigr)
 \\
 +\frac{2}{\pi}
 \int_{\cnum}
 \ftilde_1(w)
 \cdot\del_{w}\del_{\wbar}\Bigl(
 \bigl(
 1-\rho(3(z-w))
 \bigr)
 \cdot\log|z-w|
 \Bigr).
\end{multline}
It is rewritten as follows.
\begin{multline}
 \frac{2}{\pi}
 \int_{\cnum} \del_w\del_{\wbar}\ftilde_1(w)
 \cdot\Bigl(
 \chi_{z,\epsilon}(w)
 \log|z-w|
 \Bigr)
\\
 +\frac{2}{\pi}
 \int_{\cnum}
 \del_w\del_{\wbar}
 \ftilde_1(w)
 \cdot\Bigl(
 \rho(3(z-w))\cdot
 \bigl(1-\chi_{z,\epsilon}(w)\bigr)
 \cdot\log|z-w|
 \Bigr)
 \\
 +\frac{2}{\pi}
 \int_{\cnum}
 \ftilde_1(w)
 \cdot\del_{w}\del_{\wbar}\Bigl(
 \bigl(
 1-\rho(3(z-w))
 \bigr)
 \cdot\log|z-w|
 \Bigr).
\end{multline}
Here,
$\del_w\del_{\wbar}\ftilde_1$
in the second integral is the current.
According to \cite[Lemma 2.2]{Simpson90},
$f_2\leq \del_w\del_{\wbar}f_1$
weakly holds on $U$.
Note that
$\bigcup_{|z|<1/4}
\bigl\{
|w-z|<2/3
\bigr\}\subset \{|w|<1\}$.
Hence, we obtain
the following inequality:
\begin{multline}
\frac{2}{\pi}
 \int_{|z-w|\leq 1/3}
  f_2(w)\cdot
  \Bigl(
 -\rho(3(z-w))\log|z-w|
 \Bigr) \leq \\
 -F_{\ftilde_1}(z)
  +\frac{2}{\pi}
 \int_{\cnum}
 \ftilde_1(w)
 \cdot\del_{w}\del_{\wbar}\Bigl(
 \bigl(
 1-\rho(3(z-w))
 \bigr)
 \cdot\log|z-w|
 \Bigr).
\end{multline}
Note that $-F_{\ftilde_1}(z)=-f_1(z)$
is bounded as $|z|\to 0$.
We also note that
$\{|w|<1/10\}\subset
\{|z-w|<1/3\}$
as $|z|\to 0$.
Hence, there exists $C>0$ such that
the following holds for any $|z|<1/100$:
\[
 \int_{U'}
 f_2(w)\cdot\bigl(
 -\log|w-z|
 \bigr)
\leq C. 
\]
By Fatou's lemma,
we obtain
$\int_{U'}f_2(w)(-\log|w|)<\infty$.
\hfill\qed

\section{Family of vector bundles on torus with small curvature}
\label{subsection;20.8.1.1}

\subsection{Preliminary}

Let $T_0:=\real^2/\seisuu^2$.
Let $(x_1,x_2)$ denote the local coordinate systems of $T_0$
induced by the standard coordinate system of $\real^2$.
Set $U:=\bigl\{(\xi_1,\ldots,\xi_{n-2})\,\big|\,|\xi_i|\leq 1\bigr\}
\subset\real^{n-2}$.
For any non-negative integer $k$,
we set 
\[
 S_1(k):=
  \bigl\{(m_1,m_2)\in\seisuu_{\geq 0}^2\,\big|\,m_1+m_2=k\bigr\},
\]
\[
 S_2(k):=
 \bigl\{(m_1,\ldots,m_{n-2})\in\seisuu_{\geq 0}^{n-2}
  \,\big|\,\sum m_i=k\bigr\}.
\]
We set $S(k_1,k_2):=S(k_1)\times S(k_2)$.
We put $\del_{\vecx}^{\vecm}:=\prod\del_{x_i}^{m_i}$
and $\del_{\vecxi}^{\vecm}:=\prod\del_{\xi_i}^{m_i}$.
We put $N(k_i):=|S(k_i)|$
and $N(k_1,k_2):=N(k_1)\cdot N(k_2)$.

In general,
for a vector space $V$
and
for$f\in C^{\infty}(T_0\times U,V)$,
we set
\[
 D^{k_1}_{\vecx}D^{k_2}_{\vecxi}(f):=
 \bigl(\del_{\vecx}^{\vecm_1}\del_{\vecxi}^{\vecm_2}f\,|\,
 (\vecm_1,\vecm_2)\in S(k_1,k_2)\bigr)
 \in C^{\infty}(T_0\times U,V^{N(k_1,k_2)}).
\]
We formally set $D^0f:=f$.

Let $E$ be a topologically trivial $C^{\infty}$ vector bundle
$T_0\times U$
with a Hermitian metric $h$ and a unitary connection $\nabla$.
We set $r:=\rank E$.
Let $F(\nabla)$ denote the curvature of $\nabla$.
Fix a positive integer $M$.
Let $\epsilon$ denote a small positive number.
We suppose the following condition.
\begin{condition}
 $\bigl|
 \nabla_{x_1}^{m_1}\nabla_{x_2}^{m_2}
 \nabla_{\xi_1}^{\ell_1}\cdots\nabla_{\xi_{n-2}}^{\ell_{n-2}}
 F(\nabla)\bigr|\leq \epsilon$ 
if $\sum m_i\leq M$ and $\sum \ell_i\leq M$.
\hfill\qed
\end{condition}

\subsection{Partially almost holomorphic frames}

Let $\tau\in\cnum$ with $\Image(\tau)>0$.
We set $T:=\cnum/(\seisuu+\tau\seisuu)$.
Let $z$ denote the standard coordinate of $\cnum$,
which induces local coordinates on $T$.
We identify the $C^{\infty}$-manifolds
$T_0:=\real^2/\seisuu^2$ and $T$
by the diffeomorphism induced by 
$(x,y)\longmapsto z=x+\tau y$.

For any frame $\vecw$ of $E$,
let $A^{\vecw}_z$ and $A^{\vecw}_{\zbar}$
be the matrix valued functions determined by
$\nabla_{\zbar}\vecw=\vecw A^{\vecw}_{\zbar}$
and 
$\nabla_{z}\vecw=\vecw A^{\vecw}_{z}$.
Similarly,
let $A^{\vecw}_{\xi_i}$
be determined by 
$\nabla_{\xi_i}\vecw=\vecw A^{\vecw}_{\xi_i}$,
and we set
$\lefttop{2}A^{\vecw}:=\sum A^{\vecw}_{\xi_i}d\xi_i$.
Let $H(h,\vecw)$ denote the matrix valued function
whose $(i,j)$-th entries are  $h(w_i,w_j)$.

Take $p\geq 2$.
We fix a norm $\|\cdot\|_{L_k^p(T)}$
on the Banach space of $L_k^p$-functions on $T$.
For a function $f$ on $T\times U$,
let $\|f\|_{L_k^p}$ denote the function
$U\lrarr \real_{\geq 0}$
determined by
$\|f\|_{L_k^p}(\vecxi)=\|f_{|T\times\{\vecxi\}}\|_{L_k^p(T)}$.
The following lemma is 
a reword of \cite[Proposition 4.7]{Mochizuki-doubly-periodic}.

\begin{lem}
\label{lem;17.10.23.10}
If $\epsilon>0$ is sufficiently small,
there exist $\Gamma\in M_r(\cnum)$
 and a frame $\vecu$ of $E$ on $T\times U$
 satisfying the following conditions
for a positive constant $C_1$
depending only on $n$, $r$, $M$ and $p$.
\begin{itemize}
\item
$\bigl|[\Gamma,\lefttop{t}\Gammabar]\bigr|\leq C_1\epsilon$.
\item
 $A^{\vecu}_{\zbar}$ is constant along 
 the $T$-direction,
 and $\bigl|D^k_{\vecx}(A^{\vecu}_{\zbar}-\Gamma)\bigr|
 \leq C_1\epsilon$
for $0\leq k\leq M$.
\item
 $\bigl\|
 D^k_{\vecxi}(A^{\vecu}_z+\lefttop{t}\Gammabar)
 \bigr\|_{L_M^p}
 \leq C_1\epsilon$
 for any $0\leq k\leq M$.
\item
 $\bigl\| 
 D^k_{\vecxi}
 \bigl(\lefttop{2}A^{\vecu}\bigr)
 \bigr\|_{L_M^p}\leq C_1\epsilon$
for $0\leq k\leq M$.
\end{itemize}
Moreover,
$\|D^k_{\vecxi}(H(h,\vecu)-I_r)\|_{L_{M+1}^p}\leq C_1\epsilon$
for $0\leq k\leq M$,
where $I_r\in M_r(\cnum)$ denotes the identity matrix.
\hfill\qed
\end{lem}

\subsection{Spectra}

For each $\vecxi\in U_2$,
we regard
$E_{\vecxi}:=E_{|T\times\{\vecxi\}}$
with $\nabla_{\zbar}$
as a holomorphic vector bundle on $T$.
If $\epsilon$ is sufficiently small,
$E_{\vecxi}$ are semistable of degree $0$
for any $\vecxi\in U_2$.
Let $T^{\lor}$ denote the dual of $T$.
Let
$\RFM:D^b(T)\lrarr D^b(T^{\lor})$
denote the Fourier-Mukai transform.
Because $E_{\vecxi}$ is semistable of degree $0$,
$\RFM(E_{\vecxi})$ is a torsion $\nbigo_{T^{\lor}}$-module.
The support and the length of $\RFM(E_{\vecxi})$
determines a point in $\Sym^r(T^{\lor})$,
which is denoted by $[\Sp(E_{\vecxi})]$.

The eigenvalues of $\Gamma$
determines the point 
$[\Sp(\Gamma)]\in \Sym^r(\cnum)$.
The quotient map
$\Phi:\cnum\lrarr T^{\lor}$
induces
$\Sym^r\cnum\lrarr \Sym^r T^{\lor}$,
which is also denoted by $\Phi$.
Let $d_{\Sym^rT^{\lor}}$
denote the distance on $T^{\lor}$
induced by a $C^{\infty}$-Riemannian metric
of $\Sym^rT^{\lor}$.
The following is proved in
\cite[Corollary 4.10]{Mochizuki-doubly-periodic}.
\begin{cor}
\label{cor;17.10.24.1}
Let $\Gamma$ be as in Lemma {\rm\ref{lem;17.10.23.10}}.
There exist $\epsilon_0>0$ and $C_0>0$
depending only on $r$
such that the following holds
if $\epsilon\leq \epsilon_0$:
\[
 d_{\Sym^r(T^{\lor})}\bigl(
 [\Sp(E_{\vecxi})],
 \Phi([\Sp(\Gamma)])
 \bigr)
\leq
 C_0\epsilon.
\]
\hfill\qed
\end{cor}

\subsection{Additional assumption on spectra}

Let $\epsilon_0$ be as in Corollary \ref{cor;17.10.24.1}.
We assume that $\epsilon<\epsilon_0$.
Moreover, we assume the following.
\begin{condition}
There exist
$0<\rho$,
$0<\delta<1/10$,
$\zeta_0\in\cnum$,
and a finite subset
\[
 Z\subset
 \bigl\{
 \zeta\in\cnum\,\big|\,
 0\leq 
 \Image(\zeta-\zeta_0)\leq 
 (1-\delta)\pi,\,\,
 0\leq 
 \Image\bigl(\taubar(\zeta-\zeta_0)\bigr)
 \leq 
 (1-\delta)\pi
 \bigr\}
\]
such that the following holds.
\begin{itemize}
\item
 For any distinct points $\nu_1,\nu_2\in Z$,
 we have $d_{\cnum}(\nu_1,\nu_2)>100 r^2\rho$,
 where $d_{\cnum}$ denotes the Euclidean distance on $\cnum$.
\item
 For any $\kappa\in \Sp(E_{\vecxi})$,
 there exists $\nu\in Z$
 such that 
 $d_{T^{\lor}}(\Phi(\nu),\kappa)<\rho$,
     where $d_{T^{\lor}}$ denotes the natural distance on $T^{\lor}$.
\end{itemize}
 We also assume that $\epsilon$ is sufficiently smaller than $\rho^2$.
\hfill\qed
\end{condition}

We obtain the $C^{\infty}$-decomposition
$E=\bigoplus_{\nu\in Z} E_{\nu}$
such that
(i) $\nabla_{\zbar}$ preserves the decomposition,
(ii) any $\kappa\in\Sp(E_{\nu|T\times\{\vecxi\}})$
 satisfies $d_{T^{\lor}}(\Phi(\nu),\kappa)<\rho$.
The following lemma is implicitly used
in \cite[\S4.3, \S5.2]{Mochizuki-doubly-periodic}.

\begin{lem}
\label{lem;17.10.23.12}
If $\epsilon$ is sufficiently small,
there exist $\Gamma\in M_r(\cnum)$
and a frame $\vecu$ of $E$ such that the following holds.
\begin{itemize}
\item
$\Gamma$ and $\vecu$ satisfy
the conditions in Lemma {\rm\ref{lem;17.10.23.10}}.
\item
There exists a decomposition
$\vecu=\bigcup_{\nu\in Z}\vecu_{\nu}$
such that 
$\vecu_{\nu}=(u_{\nu,1},\ldots,u_{\nu,\rank E_{\nu}})$
are frames of the subbundles $E_{\nu}$.
\item
We may impose the following
for any $\vecxi\in U_2$:
\begin{equation}
\label{eq;17.10.23.11}
\frac{1}{|T|} \int_{T\times\{\vecxi\}} h(u_{\nu,i},u_{\nu,j})=\left\{
 \begin{array}{ll}
 1 & (i=j),\\
 \mbox{{}}\\
 0 & (\mbox{\rm otherwise}).
 \end{array}
 \right.
\end{equation}
\end{itemize}
\end{lem}
\pf
In the proof of this lemma,
$C_i$ denote positive constants
depending only on $n$, $r$, $M$, $p$ and $\rho$.
We take $\Gamma$ and $\vecu$ as in Lemma \ref{lem;17.10.23.10}.
There exists a unitary matrix $B^{(0)}$
such that the following holds.
\begin{itemize}
\item
 $\Gamma':=(B^{(0)})^{-1}\Gamma B^{(0)}$ is 
 an upper triangular matrix.
\item
 There exist a total ordering $\{\nu_1,\ldots,\nu_m\}$ on $Z$
 and numbers $r_i$ $(i=1,\ldots,m)$ with $\sum r_i=r$
 such that 
 $(j,j)$-th entries of $\Gamma'$ 
 are contained in the $2\rho$-ball centered at $\nu_i$,
where
$\sum_{k\leq i-1} r_k+1\leq j\leq \sum_{k\leq i}r_k$.
\end{itemize}
Set
$\Lambda_i:=\{j\in\seisuu\,|\,
\sum_{k\leq i-1}r_k+1\leq j\leq \sum_{k\leq i}r_k\}$.
If $j_1\in \Lambda_1$, $j_2\in \Lambda_2$
with $\Lambda_1\neq \Lambda_2$,
the $(j_1,j_2)$-entries of $\Gamma'$ is 
smaller than $C_3\epsilon$
because $\bigl|[\Gamma,\lefttop{t}\Gammabar]\bigr|\leq C_1\epsilon$.
Let $\Gamma'_{i}$ denote the $r_i$-square matrix
given by the $(j_1,j_2)$-entries $(j_1,j_2\in \Lambda_i)$
of $\Gamma'$.
Then, there exists
an upper triangular matrix $B^{(1)}\in M_r(\cnum)$
such that the following holds.
\begin{itemize}
\item
 $B^{(1)}_{j_1,j_2}=0$ if $j_1$ and $j_2$ are contained 
 in the same $\Lambda_i$.
\item
 $|B^{(1)}|\leq C_4\epsilon$.
\item
 $(I_r+B^{(1)})^{-1}\Gamma' (I_r+B^{(1)})=\bigoplus \Gamma'_i$.
\end{itemize}
Then, 
$\bigoplus\Gamma'_i$
and 
$\vecu':=\vecu\cdot B^{(0)}(I+B^{(1)})$ 
satisfy the conditions in Lemma \ref{lem;17.10.23.10},
and $\vecu'$ has a decomposition
with the desired property.
Hence, we may assume the existence of the decomposition
$\vecu=\bigcup \vecu_{\nu}$ from the beginning.

Set $r_{\nu}:=\rank E_{\nu}$.
Let $H(h,\vecu_{\nu})$ be the $r_{\nu}$-th square Hermitian matrices
valued function determined by
\[
 H(h,\vecu_{\nu})_{k,\ell}:=
 h(u_{\nu,k},u_{\nu,\ell}).
\]
We obtain
$\bigl\|
 D^{k}_{\vecxi}\bigl(
 H(h,\vecu_{\nu})_{k,\ell}-I_{r_{\nu}}
 \bigr)
 \bigr\|_{L^k_p}
\leq C_5\epsilon$
for $0\leq k\leq M$,
where $I_{r_{\nu}}\in M_{r_{\nu}}(\cnum)$ denotes
the identity matrix.

Let $H_{1,\nu}$ be the function from $U$ to 
the space of $r_{\nu}$-th positive definite Hermitian matrices
determined by
$\bigl(\overline{H_{1,\nu}}\bigr)^{-2}
=\int_T H(h,\vecu_{\nu})$.
Then, we have
$\bigl\|
 D^{k}_{\vecxi}\bigl(
 H_{1,\nu}-I_{r_{\nu}}
 \bigr)
 \bigr\|_{L_k^p}
\leq C_6\epsilon$
for $0\leq k\leq M$.
Then,
$\Gamma$ and $\bigcup_{\nu\in Z}
 \bigl(\vecu_{\nu}\cdot H_{1,\nu}\bigr)$
satisfies (\ref{eq;17.10.23.11}).
\hfill\qed

\vspace{.1in}

There exists the decomposition
$\nabla_{\zbar}=\nabla_{\zbar,0}+f$
such that 
(i) $(E_{\vecxi},\nabla_{\zbar,0})$ is holomorphically isomorphic to
the product vector bundle $\cnum^r\times T$ on $T$
for any $\vecxi\in U$,
(ii) $\nabla_{\zbar,0}f=0$,
(iii) $\Sp(f_{|T\times\vecxi})$ is contained in the $\rho$-ball 
 of $Z$.
Let $\nbigv_{\vecxi}$ denote the space of
holomorphic sections of
$(E_{\vecxi},\nabla_{\zbar,0})$.
Then,
$\nbigv_{\vecxi}$ $(\vecxi\in U)$
naturally give a $C^{\infty}$-bundle $\nbigv$ on $U$,
and $\vecu$ induces a frame of $\nbigv$.
Let $p:T\times U\lrarr U$ denote the projection.
Then, there exists a naturally defined $C^{\infty}$-isomorphism
$E\simeq p^{\ast}(\nbigv)$ on $T\times U$.
We obtain a decomposition
$\nbigv=\bigoplus \nbigv_{\nu}$,
which is compatible with the above isomorphism.

\subsection{Spaces of functions}
\label{subsection;20.7.28.1}

Let $C^M_{\vecxi}L_{M,\vecx}^p$
denote the space of $C^M$-functions
$U\lrarr L_{M}^p(T)$.
Let $C^M_{\vecxi}L_{M,\vecx}^p(E)$ denote the space of
sections $f=\sum f_iu_i$ of $E$
such that $f_i\in C^M_{\vecxi}L_{M,\vecx}^p$,
where $\vecu$ is as in Lemma \ref{lem;17.10.23.12}.
The space is independent of the choice of $\vecu$.
Because $E=p^{\ast}\nbigv$,
we naturally regard
$C^{M}(U,\nbigv)$
as a subspace of
$C^{M}(T\times U,E)$.
There also exists the naturally defined morphism
$C^M_{\vecxi}L_{M,\vecx}^p(E)
\lrarr C^M(U,\nbigv)$
induced by the integral along the fibers.
The kernel is denoted by
$C^M_{\vecxi}L_{M,\vecx}^p(E)_0$.
Similar spaces are defined for 
$\End(E_{\nu})$ and $\Hom(E_{\nu},E_{\mu})$.
We set
\[
 C^M_{\vecxi}L^p_{M,\vecx}(\End(E))^{\circ}:=
 \bigoplus_{\nu\in Z} C^M(U,\End(\nbigv_{\nu})),
\]
\[
 C^M_{\vecxi}L^p_{M,\vecx}(\End(E))^{\bot}:=
 \bigoplus_{\nu}
 C^M_{\vecxi}L_{M,\vecx}^p(\End(E_{\nu}))_{0}
\oplus
 \bigoplus_{\nu\neq \mu}
 C^M_{\vecxi}L_{M,\vecx}^p(\Hom(E_{\nu},E_{\mu})).
\]
We obtain the decomposition
\[
C^M_{\vecxi}L^p_{M,\vecx}(\End(E))=
 C^M_{\vecxi}L^p_{M,\vecx}(\End(E))^{\circ}\oplus
 C^M_{\vecxi}L^p_{M,\vecx}(\End(E))^{\bot}.
\]
For any $s\in C^M_{\vecxi}L^p_{M,\vecx}(\End(E))$,
we obtain the corresponding decomposition
$s=s^{\circ}+s^{\bot}$.

\subsection{Some estimates}
\label{subsection;17.10.24.2}

Let $\gtilde:=\overline{H(h,\vecu)}$.
We obtain
$\bigl\|
 D^{k}_{\vecxi}\bigl(\gtilde-I_r\bigr)
 \bigr\|_{L^p_k}\leq C\epsilon$
for $0\leq k\leq M$,
where $C$ is a positive constant
depending only on 
$n$, $r$, $M$, $p$ and $\rho$.
Let $F_{z\zbar}$ denote
the $dz\,d\zbar$-component of
$F(\nabla)$.
We obtain the decomposition
$F_{z\zbar}=F_{z\zbar}^{\circ}+F_{z\zbar}^{\bot}$
as in \S\ref{subsection;20.7.28.1}.
The following lemma is proved
by the argument in the proof of 
\cite[Lemma 4.12, Lemma 4.15]{Mochizuki-doubly-periodic}.

\begin{lem}
\label{lem;17.10.23.30}
There exist $C_{10}>0$ and $\epsilon_{10}>0$,
depending only on $r$, $M$, $n$, $p$ and $\rho$,
such that if $\epsilon<\epsilon_{10}$,
we obtain
\[
 \bigl\|
 D^{k}_{\vecxi}
 \bigl(
 \gtilde-I
 \bigr)
\bigr\|_{L^p_{M+2}}
\leq
 C_{10}
 \sum_{j=0}^k
 \bigl\|
 D^j_{\vecxi}F^{\bot}_{z\zbar}
 \bigr\|_{L^p_{M}}.
\]
In particular, we have
$\sup
 \bigl|
 \gtilde-I
 \bigr|\leq
 C_{11}
 \bigl\|
 F^{\bot}_{z\zbar}
 \bigr\|_{L^2}$
for a positive constant $C_{11}$
depending only on $n$, $r$ and $\rho$.
\hfill\qed
\end{lem}

\begin{rem}
The frame $\vecu$ in {\rm\cite[\S4.3]{Mochizuki-doubly-periodic}}
should be chosen as in Lemma {\rm\ref{lem;17.10.23.12}}.
\hfill\qed
\end{rem}

\section{Estimates for asymptotic doubly periodic instantons}
\label{subsection;20.8.1.2}

\subsection{Setting}

Let $L$ be a lattice in $\cnum$.
Let $T:=\cnum/L$.
The standard coordinate $z$ of $\cnum$
induces local coordinates on $T$.
For any $R>0$,
let $U_w^{\ast}(R):=\bigl\{w\in\cnum\,\big|\,|w|\geq R\bigr\}$
and $U_w(R):=U_w^{\ast}(R)\cup\{\infty\}$ in $\proj^1_w$.
Take $q\in\seisuu_{\geq 1}$.
Let $\proj^1_{w_q}\lrarr \proj^1_w$ be the ramified covering
given by $w_q\longmapsto w_q^q$.
Let $U_{w,q}(R)$ and $U_{w,q}^{\ast}(R)$
denote the pull back of
$U_{w}(R)$ and $U_{w}^{\ast}(R)$,
respectively.
We use the metric
$dz\,d\zbar+dw\,d\wbar
=dz\,d\zbar+q^2|w_q|^{2(q-1)}dw_q\,d\wbar_w$
on $T\times U_{w,q}^{\ast}(R)$.

Take $R_0>0$.
Let $(E,\delbar_E)$ be a holomorphic vector bundle
on $T\times U_{w,q}^{\ast}(R_0)$ with a Hermitian metric $h$.
Let $F(h)$ denote the curvature of the Chern connection $\nabla_h$
determined by $h$ and $\delbar_E$.
For any $\End(E)$-valued form $s$,
let $|s|_h$ denote the function obtained as the norm of
$s$ with respect to $h$ and the Euclidean metric. 

\begin{condition}
Let $\nbiga$
denote the set of polynomials of non-commutative variables
$t_1,t_2,t_3,t_4$.
For any $P\in\nbiga$,
there exist $C(P)>0$ and $\epsilon(P)>0$ such that
\begin{equation}
 \label{eq;17.10.23.20}
 \Bigl|
 P(\nabla_z,\nabla_{\zbar},\nabla_w,\nabla_{\wbar})
 \bigl(
 \Lambda F(h)
 \bigr)
 \Bigr|_h
\leq C(P) \exp\bigl(-\epsilon(P)I_q(w_q)\bigr).
\end{equation}
Here, $I_q(w_q):=|w_q^q|(\log|w_q|)^{1/2}$.
Namely,
$(E,\delbar_E,h)$ asymptotically satisfies
the condition of instantons.
We also assume
$|F(h)|_h\to 0$ as $|w_q|\to\infty$.
\hfill\qed     
\end{condition}

The following lemma is standard
due to \cite{Uhlenbeck1,Uhlenbeck0}.
\begin{lem}
\label{lem;21.9.13.40}
For any $P\in\nbiga$, we obtain
$\bigl|
P(\nabla_{z},\nabla_{\zbar},\nabla_w,\nabla_{\wbar})
F(h)\bigr|_h\to 0$
as $|w_q|\to\infty$.
\end{lem}
\pf
For $Q\in T\times U_{w,q}^{\ast}$
and $\rho>0$,
we set
\[
B_{Q}(\rho)
=\bigl\{
Q'\in T\times U_{w,q}^{\ast}
 \,\big|\,
d_{T\times U_{w,q}^{\ast}}(Q,Q')<\rho
\bigr\},
\]
where $d_{T\times U_{w,q}^{\ast}}$
denotes the distance induced by
$dz\,d\zbar+dw\,d\wbar$.
Let $w(Q)$ denote the image of
$Q$ via the natural map
$T\times U_{w,q}^{\ast}\lrarr U_{w}^{\ast}$.
Let $p$ be any large real number.
For any unitary frame $\vecu$ of
$E_{|B_Q(2)}$,
let $A^{\vecu}$ denote the Chern connection form
with respect to $\vecu$.
Similarly, let $F^{\vecu}$ denote the curvature form
representing $F(h)$ with respect to $\vecu$.
We have the elliptic system of partial differential equations
\[
\Lambda(dA^{\vecu}+A^{\vecu}\wedge A^{\vecu})=\Lambda F^{\vecu},
\quad
d^{\ast}A^{\vecu}=0.
\]

For any $\delta>0$ and $\ell\in\seisuu_{\geq 0}$,
there exists $R(\delta,\ell)>R$ such that
the following holds
$(\ell_1,\ell_2,\ell_3,\ell_4)\in\seisuu_{\geq 0}^4$
with $\sum\ell_i\leq \ell$
if $|w(Q)|>R(\delta,\ell)$:
\begin{equation}
\label{eq;21.9.16.3}
\bigl|
 \nabla_z^{\ell_1}\nabla_{\zbar}^{\ell_2}
 \nabla_w^{\ell_3}\nabla_{\wbar}^{\ell_4}
 (\Lambda F^{\vecu})
\bigr|\leq \delta.
\end{equation}
We also assume the following
if $|w(Q)|>R(\delta,\ell)$.
\begin{equation}
\label{eq;21.9.16.4}
|F^{\vecu}|\leq \delta.
\end{equation}

If $\delta$ is sufficiently small,
and if $|w(Q)|>R(\delta,\ell)$
which implies (\ref{eq;21.9.16.4}),
according to \cite{Uhlenbeck1,Uhlenbeck0},
there exists a unitary frame
$\vecu$ of $E_{|B_Q(2)}$
such that
(i) $d^{\ast}A^{\vecu}=0$,
(ii) $\|A^{\vecu}_{|B_Q(2)}\|_{L_1^p}\leq
C_1\|F(h)_{|B_Q(2)}\|_{L^p}$
for a positive constant $C_1$
which is independent of $Q$.
Let $1<r_1<r_2<2$.
Let $\kappa:\real_{\geq 0}\lrarr [0,1]$
be a $C^{\infty}$-function such that
(i) $\kappa(u)=1$ $(u\leq r_1)$,
(ii) $\kappa(u)=0$ $(u\geq r_2)$,
(iii) $\kappa^a$ is $C^{\infty}$ for any $a>0$.
We set $\chi(Q'):=\kappa(d(Q,Q'))$
and $A^{\vecu}_a:=\chi^a A^{\vecu}_{|B_Q(2)}$.
There exists $C'_1(a)>0$, which is independent of $Q$,
such that
$\|A^{\vecu}_a\|_{L_1^p}\leq C'_1(a)\delta$.
We have
\begin{equation}
\label{eq;21.9.13.30}
 \Lambda d(A^{\vecu}_a)=
 \Lambda\Bigl(
 \bigl(\chi^{-a/2}d(\chi^a)\bigr)A^{\vecu}_{a/2}
 -A^{\vecu}_{a/2}\wedge A^{\vecu}_{a/2} \Bigr)
+\chi^a\Lambda F^{\vecu},
\end{equation}
\begin{equation}
\label{eq;21.9.13.31}
  d^{\ast}(A^{\vecu}_a)
 =\nbigp\bigl(
 \chi^{-a/2}d\chi^a,
 A^{\vecu}_{a/2}
 \bigr),
\end{equation}
where $\nbigp(\chi^{-a/2}d\chi^a,A^{\vecu}_{a/2})$
is bilinear with respect to 
$\chi^{-a/2}d\chi^a$ and $A^{\vecu}_{a/2}$.
We apply \cite[Theorem 17.2]{ADN2}
to the system (\ref{eq;21.9.13.30}, \ref{eq;21.9.13.31}).
Then, there exists $C_2(a)>0$,
which is independent of $Q$,
such that
$|A^{\vecu}_a|_{L_2^p}\leq C_2(a)\delta$.
By an induction, we can prove that
there exist $C_{j}(a)>0$ $(j=2,\ldots,\ell+1)$,
which are independent of $Q$,
such that 
$|A^{\vecu}_a|_{L_{j}^p}\leq C_j(a)\delta$.
Thus, we obtain Lemma \ref{lem;21.9.13.40}.
\hfill\qed

\vspace{.1in}

By replacing $R_0$ with a larger number,
we may assume that
$(E,\delbar_E)_{|T\times\{w_q\}}$
are semistable of degree $0$ for any $w_q\in U_{w,q}^{\ast}(R_0)$
from the beginning.
We obtain the spectral curve
$\Sigma_E\subset T^{\lor}\times U_{w,q}^{\ast}(R_0)$.
We assume the following.
\begin{condition}
 The closure $\overline{\Sigma}_E$
 of $\Sigma_E$
 in $T^{\lor}\times U_{w,q}(R_0)$ is 
 a complex analytic curve.
We set
$\Sp_{\infty}(E):=
 \overline{\Sigma}_{E}\cap
 \bigl(
 T^{\lor}\times\{\infty\}
 \bigr)$. 
\hfill\qed
\end{condition}

\subsection{Decomposition}

We fix $\Sptilde_{\infty}(E)\subset\cnum$
such that
the projection $\cnum\lrarr T^{\lor}$
induces a bijection
$\Sptilde_{\infty}(E)\simeq\Sp_{\infty}(E)$.
Let $\Psi:T\times U_{w,q}^{\ast}(R_0)\lrarr U_{w,q}^{\ast}(R_0)$
denote the projection.
We obtain the following lemma
(see \cite[\S2.1]{Mochizuki-doubly-periodic}).
\begin{lem}
There exists a holomorphic vector bundle
$\nbigv$ with an endomorphism $f$
on $U_{w,q}^{\ast}(R_0)$
such that the following holds.
\begin{itemize}
\item
 For any $\rho>0$,
 there exists $R_1\geq R_0$ such that 
 if $|w_q^q|\geq R_1$
 and if $\alpha$ is an eigenvalue of
 $f_{|w_q}$,
 then there exists
 $\beta\in \Sptilde_{\infty}(E)$
 such that
 $|\beta-\alpha|\leq \rho$.
\item
Let $(\Psi^{-1}(\nbigv),\delbar_{\Psi^{-1}(\nbigv)})$
denote the holomorphic vector bundle
on $T\times U_{w,q}^{\ast}(R_0)$
obtained as the pull back of $\nbigv$.
As a twist,
we obtain the holomorphic vector bundle
$\Psi^{\ast}(\nbigv,f)=
 \bigl(
 \Psi^{-1}(\nbigv),
 \delbar_{\Psi^{-1}(\nbigv)}+f\,d\zbar
 \bigr)$.
Then, there exists an isomorphism
$\Psi^{\ast}(\nbigv,f)\simeq E$.
\hfill\qed
\end{itemize}
\end{lem}

There exists the decomposition
$(\nbigv,f)=\bigoplus_{\beta\in\Sptilde_{\infty}(E)}
 (\nbigv_{\beta},f_{\beta})$,
where the eigenvalues of $f_{\beta|w_q}$ are convergent to $\beta$
as $|w_q|\to\infty$.
Under the identification
$E=\Psi^{\ast}(\nbigv,f)$,
we obtain the decomposition
$E=\bigoplus_{\alpha\in\Sptilde_{\infty}(E)}
 \Psi^{\ast}(\nbigv_{\alpha},f_{\alpha})$.

\subsection{Estimates}

Let $h_{\alpha}$ be the restriction of $h$
to $\Psi^{\ast}(\nbigv_{\alpha},f_{\alpha})$.
We obtain the Fourier expansion of $h_{\alpha}$
along the $T$-direction.
Let $h_{V,\alpha}$ be the invariant part of $h_{\alpha}$.
In other words,
for any sections $u_1,u_2$ of $\nbigv_{\alpha}$,
we set 
\[
h_{V,\alpha}(u_1,u_2):=
 |T|^{-1}\int_Th\bigl(\Psi^{-1}(u_1),\Psi^{-1}(u_2)\bigr).
\]

We set
$h^{\circ}:=\bigoplus_{\alpha\in\Sptilde_{\infty}(E)}
 \Psi^{-1}(h_{V,\alpha})$.
Let $b$ be the automorphism of $E$
determined by
$h=h^{\circ}\cdot b$.
The following proposition is
essentially proved in \cite[\S5]{Mochizuki-doubly-periodic}.
\begin{prop}
\label{prop;17.10.23.31}
For any $P\in\nbiga$,
there exist $C(P)>0$ and $\epsilon(P)>0$
such that
\[
  P(\nabla_z,\nabla_{\zbar},\nabla_w,\nabla_{\wbar})(b-\id)
\leq C(P)\exp\bigl(-\epsilon(P)|w|\bigr).
\]
\end{prop}
\pf
We shall explain an outline of the proof.
We set 
$E_{\alpha}:=\Psi^{\ast}(\nbigv_{\alpha},f_{\alpha})$.
We obtain the decomposition
$\End(E)=\bigoplus \Hom(E_{\alpha},E_{\beta})$.

Let $U$ be any open subset in $U_w^{\ast}(R_0)$.
We have the map
$\int_T:C^{\infty}(T\times U,\End(E_{\alpha}))
\lrarr
 C^{\infty}(U,\End(\nbigv_{\alpha}))$
induced by the integration along the fibers.
Let 
$C^{\infty}(T\times U,\End(E_{\alpha}))_0$
denote the kernel of $\int_T$.
We also have the natural inclusion
$C^{\infty}(U,\End(\nbigv_{\alpha}))
\lrarr
 C^{\infty}(T\times U,\End(E_{\alpha}))$
induced by the pull back.
Thus, we obtain the decomposition
$C^{\infty}(T\times U,\End(E_{\alpha}))
=C^{\infty}(U,\End(\nbigv_{\alpha}))
\oplus
 C^{\infty}(T\times U,\End(E_{\alpha}))_0$.
We set
\[
 C^{\infty}(T\times U,\End(E))^{\circ}:=
 \bigoplus_{\alpha\in \Sptilde_{\infty}(E)}
 C^{\infty}(U,\End(\nbigv_{\alpha})),
\]
\[
 C^{\infty}(T\times U,\End(E))^{\bot}:=
 \!\!\!
 \bigoplus_{\alpha\in \Sptilde_{\infty}(E)}
 \!\!\!
 C^{\infty}(T\times U,\End(E_{\alpha}))_0
\oplus
 \bigoplus_{\alpha\neq\beta}
 C^{\infty}(T\times U,\Hom(E_{\alpha},E_{\beta})).
\]
We obtain the decomposition
$C^{\infty}(T\times U,\End(E))=
 C^{\infty}(T\times U,\End(E))^{\circ}\oplus
 C^{\infty}(T\times U,\End(E))^{\bot}$.
Any section $s\in C^{\infty}(T\times U,\End(E))$ 
has the corresponding decomposition
$s=s^{\circ}+s^{\bot}$.
We also obtain the function
$\|s\|$ on $U$
determined by
$\|s\|(w_q)^2=|T|^{-1}\int_{T\times\{w_q\}}|s|^2$.
Note that there exists a constant $C>0$
such that the following holds
if $s=s^{\bot}$:
\begin{equation}
\label{eq;17.10.24.10}
 \|\nabla_{\zbar}s\|\geq C\|s\|.
\end{equation}

Let us consider the case
$U=\bigl\{w_q\,\big|\,|w_q^q-w_{q,0}^q|\leq 1
\bigr\}$ for some $w_{q,0}\in U_{w,q}$.
We have the expression
$F(h)=F_{z\zbar}dz\,d\zbar+F_{z\wbar}dz\,d\wbar
+F_{w\zbar}dw\,d\zbar+F_{w\wbar}dw\,d\wbar$.
For any $\epsilon>0$,
if $|w_{q,0}|$ is sufficiently large,
we obtain $|F(h)|\leq \epsilon$ on $T\times U$.
If $\epsilon$ is sufficiently small,
there exists a frame $\vecu$ as in Lemma \ref{lem;17.10.23.12}.
Note that $\vecu$ is a unitary frame
with respect to the metric $h^{\circ}$,
and that $\gtilde$ in \S\ref{subsection;17.10.24.2}
corresponds to the restriction of $b$ to $U$.
Then, we apply the argument for local estimates
explained in \cite[\S5.2]{Mochizuki-doubly-periodic}
by replacing the condition
$\Lambda F(h)=0$
with (\ref{eq;17.10.23.20}).
We use the notation
$g_1=O(g_2)$
if we have a positive constant $C$,
which is independent of $w_{q,0}$,
such that $g_1\leq C g_2$.
In (\ref{eq;17.10.24.11}--\ref{eq;17.10.24.14}) below,
$C_1$ will denote a positive constant
which is independent of $w_{q,0}$.
We obtain the following 
by the argument in the proof of 
\cite[Proposition 5.5]{Mochizuki-doubly-periodic}
and the condition (\ref{eq;17.10.23.20}):
\begin{multline}
\label{eq;17.10.24.11}
 -\del_w\del_{\wbar}
 \bigl\|F_{z\zbar}^{\bot}\bigr\|^2
\leq
 -\bigl\|\nabla_{\zbar}F_{z\zbar}^{\bot}\bigr\|^2
 -\bigl\|\nabla_{\zbar}F_{z\zbar}^{\bot}\bigr\|^2
 -\bigl\|\nabla_{\wbar}F_{z\zbar}^{\bot}\bigr\|^2
 -\bigl\|\nabla_{\wbar}F_{z\zbar}^{\bot}\bigr\|^2
\\
+
O\Bigl(
 \epsilon\|F_{z\zbar}^{\bot}\|^2
+\epsilon\|F_{z\zbar}^{\bot}\| \|F_{w\zbar}^{\bot}\|
+\epsilon\|\nabla_{\wbar}F_{w\zbar}^{\bot}\|\|F_{z\zbar}\|
+\epsilon\|\nabla_{\zbar}F_{z\zbar}^{\bot}\|\|F_{z\zbar}^{\bot}\|
 \Bigr)
\\
+O\Bigl(
 \epsilon\|\nabla_{\wbar}F_{z\zbar}^{\bot}\| \|F_{z\zbar}^{\bot}\|
+\epsilon \|F_{w\zbar}^{\bot}\|^2
+\epsilon\|F_{w\zbar}^{\bot}\| \|\nabla_{\wbar}F_{z\zbar}^{\bot} \|
 \Bigr)
+O\Bigl(
 \exp\bigl(-C_1I_q(w_{q,0})\bigr)
 \Bigr).
\end{multline}
We obtain the following by using  the condition (\ref{eq;17.10.23.20})
and the argument in the proof of 
\cite[Proposition 5.6]{Mochizuki-doubly-periodic}:
\begin{multline}
\label{eq;17.10.24.12}
-\del_w\del_{\wbar}\|F^{\bot}_{z\wbar}\|^2
\leq
-\|\nabla_{\zbar}F^{\bot}_{z\wbar}\|^2
-\|\nabla_{z}F^{\bot}_{z\wbar}\|^2
-\|\nabla_{\wbar}F^{\bot}_{z\wbar}\|^2
-\|\nabla_{w}F^{\bot}_{z\wbar}\|^2
 \\
+O\Bigl( 
 \epsilon\|F_{z\wbar}^{\bot}\|^2\|F_{z\zbar}^{\bot}\|
+\epsilon\|\nabla_{\wbar}F_{w\zbar}^{\bot}\| \|F_{z\wbar}^{\bot}\|
+\epsilon\|F_{z\wbar}^{\bot}\|\|F_{w\zbar}^{\bot}\|
+\epsilon\|\nabla_{\wbar}F_{z\zbar}^{\bot}\| \|F_{z\wbar}^{\bot} \|
 \Bigr)
\\
+O\Bigl(
 \epsilon \|\nabla_{\zbar}F_{z\wbar}^{\bot}\| \|F_{z\zbar}^{\bot}\|
+\epsilon \|F_{z\wbar}^{\bot}\|^2
 \Bigr)
+O\Bigl(
 \exp(-C_1I_q(w_{q,0}))
 \Bigr).
\end{multline}
By a similar argument,
we also have the following:
\begin{multline}
\label{eq;17.10.24.13}
 -\del_w\del_{\wbar}\|F_{w\zbar}^{\bot}\|^2
\leq
-\bigl\|
 \nabla_{z}F_{w\zbar}^{\bot}
 \bigr\|^2
-\bigl\|
 \nabla_{\zbar}F_{w\zbar}^{\bot}
 \bigr\|^2
-\bigl\|
 \nabla_{w}F_{w\zbar}^{\bot}
 \bigr\|^2
-\bigl\|
 \nabla_{\wbar}F_{w\zbar}^{\bot}
 \bigr\|^2
\\
+O\Bigl(
 \epsilon\|F_{w\zbar}^{\bot}\|\|F_{z\zbar}^{\bot}\|
+\epsilon\|\nabla_{\wbar}F_{w\zbar}^{\bot}\| \|F_{w\zbar}^{\bot}\|
+\epsilon\|F_{w\zbar}^{\bot}\|^2
+\epsilon\|\nabla_{\wbar}F_{z\zbar}^{\bot}\|\|F_{w\zbar}^{\bot}\|
 \Bigr)
 \\
+O\Bigl(
 \epsilon\|\nabla_{\zbar}F_{w\zbar}^{\bot}\| \|F_{z\zbar}\|
+\epsilon\|F_{w\zbar}^{\bot}\| \|F_{w\wbar}\|
\Bigr)
+O\Bigl(
 \exp\bigl(-C_1I_q(w_{q,0})\bigr)
 \Bigr).
\end{multline}
By the condition (\ref{eq;17.10.23.20}),
we have
\begin{equation}
\label{eq;17.10.24.14}
 -\del_w\del_{\wbar}\|F_{w\wbar}^{\bot}\|^2
=-\del_w\del_{\wbar}\|F_{z\zbar}^{\bot}\|^2
+O\Bigl(
 \exp\bigl(-C_1I_q(w_{q,0})\bigr)
 \Bigr).
\end{equation}

From local estimates (\ref{eq;17.10.24.11}--\ref{eq;17.10.24.14})
with (\ref{eq;17.10.24.10}),
we obtain the following inequality
on $U_w^{\ast}(R_2)$ for some $R_2>R_0$:
\begin{multline}
 -\del_w\del_{\wbar}
 \Bigl(
 \|F_{z\zbar}^{\bot}\|^2
+\|F_{w\wbar}^{\bot}\|^2
+\|F_{z\wbar}^{\bot}\|^2
+\|F_{w\zbar}^{\bot}\|^2
 \Bigr)
\leq \\
 -C_2
 \Bigl(
 \|F_{z\zbar}^{\bot}\|^2
+\|F_{w\wbar}^{\bot}\|^2
+\|F_{z\wbar}^{\bot}\|^2
+\|F_{w\zbar}^{\bot}\|^2
 \Bigr)
+C_3\exp(-C_1I_q(w_{q})).
\end{multline}
Here, $C_i$ $(i=2,3)$ are positive constants.
By a standard argument of Ahlfors lemma
(see \cite[Lemma 5.2]{Mochizuki-doubly-periodic}),
there exists $C_4>0$ such that
\[
 \|F_{z\zbar}^{\bot}\|^2
+\|F_{w\wbar}^{\bot}\|^2
+\|F_{z\wbar}^{\bot}\|^2
+\|F_{w\zbar}^{\bot}\|^2
=O\Bigl(
 \exp\bigl(-C_4|w_q^q|\bigr)
 \Bigr).
\]

Let $F(h)^{\bot}:=
 F_{z\zbar}^{\bot}dz\,d\zbar+F_{z\wbar}^{\bot}dz\,d\wbar
+F_{w\zbar}^{\bot}dw\,d\zbar+F_{w\wbar}^{\bot}dw\,d\wbar$.
By using a standard boot-strapping argument as in the proof of 
\cite[Proposition 5.8]{Mochizuki-doubly-periodic},
we obtain the following.
\begin{itemize}
\item
For any $P\in\nbiga$,
there exist $C(P)>0$ and $\epsilon(P)>0$ such that
\[
 \bigl|
 P(\nabla_z,\nabla_{\zbar},\nabla_w,\nabla_{\wbar})
 F(h)^{\bot}
 \bigr|
\leq C(P)\exp\bigl(-\epsilon(P)|w|\bigr).
\]
\end{itemize}
By Lemma \ref{lem;17.10.23.30},
we obtain the desired estimate for $b-\id$.
\hfill\qed

\begin{cor}
\label{cor;17.10.21.24}
The Higgs bundles
$(\nbigv_{\alpha},f_{\alpha}dw)$
with the metric $h_{\alpha}$
satisfy the condition in {\rm\S\ref{subsection;17.10.21.10}}.
\hfill\qed
\end{cor}

\backmatter

\printindex

\end{document}